\documentclass[11pt]{article}

\usepackage[utf8]{inputenc}
\usepackage[T1]{fontenc}

\DeclareFontShape{T1}{cmr}{m}{scit}{<->ssub*cmr/m/sc}{}

\usepackage[margin=1in]{geometry}
\usepackage{amsmath,amsthm,amssymb,amsfonts}
\usepackage{mathrsfs}
\usepackage{mathtools}
\usepackage{esint}
\usepackage{xfrac}
\usepackage{nicefrac}
\usepackage{dsfont}
\usepackage{upgreek}
\usepackage{accents}

\let\olddiamond\diamond
\let\oldsquare\square
\usepackage{mathabx}
\renewcommand{\square}{\oldsquare}
\renewcommand{\diamond}{\olddiamond}

\usepackage{booktabs}
\usepackage{array}
\usepackage{paralist}
\usepackage{enumitem}
\usepackage{verbatim}
\usepackage[normalem]{ulem}
\usepackage{cancel}

\usepackage[dvipsnames]{xcolor}
\usepackage{graphicx}
\usepackage{subcaption}

\usepackage{tikz}
\usepackage{pgf}
\usetikzlibrary{calc}
\usetikzlibrary{external}
\usepackage[
colorlinks=true,
pdfstartview=FitV,
linkcolor=blue,
citecolor=blue,
urlcolor=blue
]{hyperref}

\usepackage{titlesec}

\newcommand{\addperiod}[1]{#1.}
\titleformat{\section}{\centering\normalfont\Large}{\thesection.}{0.5em}{}
\titleformat*{\subsection}{\bfseries}
\titleformat{\subsubsection}[runin]{\normalfont\bfseries}{\thesubsubsection.}{0.5em}{\addperiod}
\titleformat*{\paragraph}{\bfseries}
\titleformat*{\subparagraph}{\large\bfseries}

\usepackage[nottoc,notlot,notlof]{tocbibind}
\usepackage{fancyhdr}
\usepackage{extramarks}

\numberwithin{equation}{section}
\numberwithin{figure}{section}

\newtheorem{theorem}{Theorem}[section]

\newtheorem{proposition}[theorem]{Proposition}
\newtheorem{lemma}[theorem]{Lemma}

\newtheorem{theoremA}{Theorem}

\theoremstyle{definition}
\newtheorem{definition}[theorem]{Definition}

\newtheorem{remark}[theorem]{Remark}

\DeclareMathOperator{\dist}{dist}
\DeclareMathOperator{\diam}{diam}

\DeclareMathOperator*{\argmin}{argmin}
\DeclareMathOperator{\var}{var}
\DeclareMathOperator{\cov}{cov}
\DeclareMathOperator{\size}{size}

\newcommand*{\supp}{\ensuremath{\mathrm{supp\,}}}

\DeclareMathSymbol{\shortminus}{\mathbin}{AMSa}{"39}

\let\originalleft\left
\let\originalright\right
\renewcommand{\left}{\mathopen{}\mathclose\bgroup\originalleft}
\renewcommand{\right}{\aftergroup\egroup\originalright}

\newcommand\thickbar[1]{\accentset{\rule{.45em}{.6pt}}{#1}}
\renewcommand{\bar}{\thickbar}
\renewcommand{\hat}{\widehat}
\renewcommand*{\tilde}{\widetilde}

\makeatletter
\let\oldoint\oint
\renewcommand{\oint}{\oldoint\nolimits}
\makeatother

\newcommand*{\N}{\ensuremath{\mathbb{N}}}
\newcommand*{\Z}{\ensuremath{\mathbb{Z}}}

\newcommand*{\R}{\ensuremath{\mathbb{R}}}

\newcommand*{\Zd}{\ensuremath{\mathbb{Z}^d}}
\newcommand*{\Rd}{\ensuremath{\mathbb{R}^d}}
\newcommand{\E}{\mathbb{E}}
\renewcommand{\P}{\ensuremath{\mathbb{P}}}

\renewcommand{\O}{\ensuremath{\mathcal{O}}}

\newcommand{\X}{\ensuremath{\mathcal{X}}}

\newcommand{\A}{\mathcal{A}}
\renewcommand{\S}{\mathcal{S}}

\newcommand{\W}{\mathcal{W}}
\newcommand{\SW}{\mathcal{SW}}

\renewcommand{\a}{\mathbf{a}}
\renewcommand{\b}{\ensuremath{\mathbf{b}}}
\newcommand{\f}{\mathbf{f}}
\newcommand{\g}{\mathbf{g}}
\newcommand{\h}{\mathbf{h}}
\renewcommand{\k}{\mathbf{k}}

\newcommand{\s}{\mathbf{s}}

\newcommand{\G}{\mathbf{G}}

\newcommand{\bfA}{\mathbf{A}}
\newcommand{\bfB}{\mathbf{B}}

\newcommand{\bfG}{\mathbf{G}}
\newcommand{\bfJ}{\mathbf{J}}
\newcommand{\bfR}{\mathbf{R}}

\newcommand{\linear}{\boldsymbol{\ell}}

\newcommand{\bhom}{\bar{\mathbf{b}}}
\newcommand{\shom}{\bar{\mathbf{s}}}
\newcommand{\khom}{\bar{\mathbf{k}}}

\newcommand{\bfAhom}{\overline{\mathbf{A}}}

\newcommand{\pot}{\mathrm{pot}}
\newcommand{\sol}{\mathrm{sol}}
\newcommand{\sym}{\mathrm{sym}}
\renewcommand{\skew}{\mathrm{skew}}

\newcommand{\loc}{\mathrm{loc}}

\newcommand*{\Id}{\ensuremath{\mathrm{I}_d}}

\newcommand*{\Itwod}{\ensuremath{\mathrm{I}_{2d}}}

\newcommand{\Rsymp}{\R^{d\times d}_{\mathrm{sym},+}}
\newcommand{\Rsym}{\R^{d\times d}_{\mathrm{sym}}}
\newcommand{\Rskew}{\R^{d\times d}_{\mathrm{skew}}}

\newcommand{\Lsol}{L^2_{\mathrm{sol}}}
\newcommand{\Lsolo}{L^2_{\mathrm{sol,0}}}
\newcommand{\Lpot}{L^2_{\mathrm{pot}}}
\newcommand{\Lpoto}{L^2_{\mathrm{pot,0}}}

\newcommand{\eps}{\varepsilon}
\newcommand{\ep}{\eps}

\newcommand{\nf}{\nicefrac}
\newcommand{\qand}{\quad \mbox{and} \quad}
\newcommand{\qqand}{\qquad \mbox{and} \qquad}

\renewcommand{\subset}{\subseteq}

\newcommand{\cs}{\mathfrak{c}_s}
\newcommand{\css}[1]{\mathfrak{c}_{#1}}
\newcommand{\cstar}{c_*}
\newcommand{\mstar}{m_*}
\newcommand{\mstarstar}{m_{**}}
\newcommand{\cgamma}{\upgamma}

\newcommand{\hsep}{\textsc{h}}

\makeatletter
\newcommand{\negphantom}{\v@true\h@true\negph@nt}
\newcommand{\neghphantom}{\v@false\h@true\negph@nt}
\newcommand{\negph@nt}{\ifmmode\expandafter\mathpalette
	\expandafter\mathnegph@nt\else\expandafter\makenegph@nt\fi}
\newcommand{\makenegph@nt}[1]{%
	\setbox\z@\hbox{\color@begingroup#1\color@endgroup}\finnegph@nt}
\newcommand{\finnegph@nt}{%
	\setbox\tw@\null
	\ifv@ \ht\tw@\ht\z@\dp\tw@\dp\z@\fi
	\ifh@\wd\tw@-\wd\z@\fi\box\tw@}
\newcommand{\mathnegph@nt}[2]{%
	\setbox\z@\hbox{$\m@th #1{#2}$}\finnegph@nt}
\makeatother

\newcommand{\Hminusul}{\hat{\phantom{H}}\negphantom{H}\underline{H}^{-1}}
\newcommand{\Hminusuls}[1]{\hat{\phantom{H}}\negphantom{H}\underline{H}^{#1}}

\newcommand{\Wminusul}[2]{\hat{\phantom{W}}\negphantom{W}\underline{W}^{#1,#2}}

\newcommand{\Besov}[3]{\mathring{\phantom{B}}\negphantom{B}\underline{B}^{#1}_{#2,#3}}

\def\Xint#1{\mathchoice
	{\XXint\displaystyle\textstyle{#1}}%
	{\XXint\textstyle\scriptstyle{#1}}%
	{\XXint\scriptstyle\scriptscriptstyle{#1}}%
	{\XXint\scriptscriptstyle\scriptscriptstyle{#1}}%
	\!\int}
\def\XXint#1#2#3{{\setbox0=\hbox{$#1{#2#3}{\int}$}
		\vcenter{\hbox{$#2#3$}}\kern-.5\wd0}}

\def\fint{\Xint-}

\newcommand{\avsum}{\mathop{\mathpalette\avsuminner\relax}\displaylimits}
\makeatletter
\newcommand\avsuminner[2]{%
	{\sbox0{$\m@th#1\sum$}%
		\vphantom{\usebox0}%
		\ooalign{%
			\hidewidth
			\smash{\,\rule[.23em]{8.8pt}{1.1pt} \relax}%
			\hidewidth\cr
			$\m@th#1\sum$\cr
		}%
	}%
}
\makeatother

\makeatletter
\newcommand\avsuminnerr[2]{%
	{\sbox0{$\m@th#1\sum$}%
		\vphantom{\usebox0}%
		\ooalign{%
			\hidewidth
			\smash{\,\rule[.23em]{6pt}{0.7pt} \relax}%
			\hidewidth\cr
			$\m@th#1\sum$\cr
		}%
	}%
}
\makeatother

\newcommand{\cu}{\square}

\newcommand{\avsumcube}[3]{\avsum_{#1 \in 3^{#2}\Zd\cap \cu_{#3}}}

\newcommand{\spx}{\triangle}

\DeclareMathAlphabet{\mathmybb}{U}{bbold}{m}{n}
\newcommand{\indc}{{\mathmybb{1}}}

\newcommand{\Tr}{\mathsf{T}}
\newcommand{\Rt}{\mathsf{R}}

\newcommand{\EthmB}{\mathsf{E}}
\newcommand{\nueff}{\nu_{\mathrm{eff}}}

\title{Superdiffusion and anomalous regularization in self-similar random incompressible flows}

\author{Scott Armstrong
\thanks{CNRS \& Laboratoire Jacques-Louis Lions, Sorbonne Universit\'e
and
Courant Institute of Mathematical Sciences, New York University.
{\footnotesize \href{mailto:scottnarmstrong@gmail.com}{scottnarmstrong@gmail.com}.}
}
\and
Ahmed Bou-Rabee
\thanks{Department of Mathematics, University of Pennsylvania.
{\footnotesize \href{mailto:ahmedmb@sas.upenn.edu}{ahmedmb@sas.upenn.edu}.}
}
\and
Tuomo Kuusi
\thanks{Department of Mathematics and Statistics, University of Helsinki.
{\footnotesize \href{mailto:tuomo.kuusi@helsinki.fi}{tuomo.kuusi@helsinki.fi}.}
}
}
\date{September 6, 2026}

\begin{document}

\maketitle
\begingroup
\renewcommand{\thefootnote}{\ifcase\value{footnote}\or*\or**\fi}

\endgroup

\begin{abstract}
We study the long-time behavior of a particle in~$\mathbb{R}^d$,~$d \geq 2$, subject to molecular diffusion and advection by a random incompressible flow. The velocity field is the divergence of a stationary random stream matrix~$\mathbf{k} $ with positive Hurst exponent~$\cgamma > 0$, so the resulting random environment is multiscale and self-similar. In the perturbative regime~$\cgamma \ll 1$, we prove quenched power-law superdiffusion: for a typical realization of the environment, the displacement variance at time~$t$ grows like~$t^{2/(2-\cgamma)}$, the scaling predicted by renormalization group heuristics. We also identify the leading prefactor up to a random (quenched) relative error of order~$\cgamma^{\nf12}\left| \log \cgamma \right|^{\nf92}$. The proof implements a Wilsonian renormalization group scheme at the level of the infinitesimal generator~$\nabla \cdot (\nu \Id + \mathbf{k} ) \nabla$, based on a self-similar induction across scales. We demonstrate that the coarse-grained generator is well-approximated, at each scale~$r$, by a constant-coefficient Laplacian with effective diffusivity growing like~$r^\cgamma$. This approximation is inherently scale-local: reflecting the multifractal nature of the environment, the relative error does not decay with the scale, but remains of order~$\cgamma^{\nf12}\left| \log \cgamma \right|^{\nf72}$. We also prove anomalous regularization under the quenched law: for almost every realization of the drift, solutions of the associated elliptic equation are H\"older continuous with exponent~$1 - C\cgamma^{\nf12}$ and satisfy estimates which are uniform in the molecular diffusivity~$\nu$ and the scale.
\end{abstract}

\setcounter{tocdepth}{2}
\tableofcontents

\section{Introduction}

\subsection{Diffusion and superdiffusion in a quenched, incompressible drift}

We study the long-time behavior of a Brownian particle advected by a random incompressible drift in~$\Rd$ in dimension~$d \geq 2$. The particle position~$X_t$ is described by the stochastic differential equation
\begin{equation}
	dX_t = \f(X_t)\,dt +\sqrt{2\nu}\,dW_t
	\,,
	\label{e.SDE}
\end{equation}
where~$\nu>0$ is the molecular diffusivity,~$\{W_t\}_{t\geq 0}$ is a standard Brownian motion on~$\Rd$ and~$\f$ is a stationary random vector field which~(i)~is incompressible,~i.e.~$\nabla\cdot\f = 0$,~(ii)~is isotropic in law  (with respect to the symmetries of the hypercube), and~(iii) possesses a correlation structure like that of a field with Hurst exponent~$-1+\cgamma$, where~$\cgamma >0$ is a small parameter. The last condition is satisfied, for instance, in the special case that~$\f$ is Gaussian and its covariances behave at large distances like
\begin{equation}
	\bigl| \cov \bigl[ \f(x) , \f(y) \bigr] \bigr|
	\asymp |x-y|^{-2(1-\cgamma)}\,, \quad |x-y| \gg 1\,.
	\label{e.our.cov.ass}
\end{equation}

It is well-understood that the presence of an incompressible vector field~$\f$ in~\eqref{e.SDE} enhances the diffusivity of the particle: the variance of~$X_t$ at time~$t$ will be at least that of the driving Brownian motion~$\sqrt{2\nu} W_t$, namely~$2 \nu t \Id$. Depending on~$\f$, it may be considerably larger, possibly even growing faster than linearly in~$t$ as~$t\to \infty$, a phenomenon known as \emph{superdiffusivity}.

\smallskip

For such a random vector field~$\f$, the primary determinant of whether we see superdiffusivity or merely enhanced diffusivity is the correlation structure of the drift. If~$\f$ has fast decaying correlations, then the effective diffusivity saturates to a finite constant in the large-scale limit and we observe enhanced diffusion. In contrast, if~$\f$ has sufficiently strong long-range correlations, the enhancements due to these larger length scales may accumulate, leading to the divergence of the effective diffusivity in the large-scale limit and, consequently, genuine superdiffusivity of the process~$\{ X_t \}$.

\smallskip

This problem was studied in the physics literature in the 1980s in the case that the incompressible drift~$\f$ is an isotropic Gaussian field with a covariance structure given by~\eqref{e.our.cov.ass}. It was predicted in~\cite{BCGLD,BG} that, in every dimension~$d\geq 2$, the long-range correlation exponent
\begin{equation*}
	\xi \coloneqq 2(1-\cgamma)
\end{equation*}
determines whether we see superdiffusivity or merely enhanced diffusivity, with the critical parameter being~$\xi = 2$ (that is,~$\cgamma=0$). It was predicted in~\cite[Section 4.3.2.2]{BG} by heuristic renormalization group~(RG) computations that, for large times~$t \gg 1$,
\begin{equation}
	\mathbf{E} \bigl[ \bigl| X_{t} \bigr|^2 \bigr]
	\,,\
	\mathbf{E} \bigl[ \bigl| X_{t}  - \mathbf{E} [ X_t ] \bigr|^2 \bigr]
	\asymp
	\left\{
	\begin{aligned}
		& t & \text{if} \ \xi > 2\,,
		&& \text{(enhanced diffusivity)} \\
		& t ( \log t )^{\nf12}   & \text{if} \ \xi = 2 \,,
		&& \text{(borderline superdiffusivity)} \\
		& t^{\frac{4}{2+\xi} }  & \text{if} \ \xi < 2 \,.
		&& \text{(superdiffusivity)} \\
	\end{aligned}
	\right.
	\label{e.physics.predictions}
\end{equation}
There are two sources of randomness in this model: here~$\mathbf{E}$ denotes, for a fixed realization of the vector field~$\f$, the expectation with respect to trajectories with initial condition~$X_0=0$, while the expectation with respect to the law of the vector field is denoted by~$\E$. The prediction~\eqref{e.physics.predictions} is therefore \emph{quenched}: the averaging is taken with respect to particle trajectories, but it should be valid for typical realizations of the random flow~$\f$.

\smallskip

The earlier paper~\cite{BCGLD} made a different prediction for the exponent in the superdiffusive case, namely~$1+\frac 14(2-\xi)$; our understanding is that this was a ``one-loop'' RG computation meant to hold at leading order in~$(2-\xi)$ for~$0< 2-\xi \ll 1$, and in this sense it indeed agrees with~\eqref{e.physics.predictions}. The guess was improved to the exact exponent appearing in~\eqref{e.physics.predictions} by the observation of~\cite{HK} that, due to constraints imposed by the isotropy assumption, no other exponent could be consistent with a nontrivial RG fixed point. The heuristic RG arguments in~\cite{BCGLD,BG} are inherently perturbative in the parameter~$\cgamma = \frac12 (2-\xi)$, with the predictions for~$\xi$ significantly less than~$2$ simply extrapolated.

\smallskip

One of the main physical motivations for studying this model is to give an explanation for Richardson's~$\nf43$ law in fluid turbulence, which states that the expected squared distance at time~$t$ for a pair of particles released at the same location in a turbulent fluid scales like~$t^3$. In this analogy, one should think of~$X_t$ as the relative displacement between the pair of particles and the drift vector field~$\f$ as the difference in the fluid velocity at the two locations (see~\cite[Section 4.1.2]{BG}).\footnote{As explained in~\cite{BG}, the ``frozen'' (time-independent) velocity field~$\f$ is a good approximation of a turbulent fluid only up to a scale-dependent \emph{eddy correlation time}. Beyond that time scale, the velocity field of a real turbulent fluid decorrelates, and the dynamics follows a different scaling law. So the frozen-field RG prediction only describes the regime where eddies have not yet changed significantly.}
To model pair displacement in a turbulent fluid according to Kolmogorov's phenomenological theory of turbulence, one should take~$\xi = -\nf 23$, which corresponds to the case that~$\f$ is a Gaussian field with Hurst exponent~$\nf13$ (and thus~$\f$ is a H\"older continuous field with regularity exponent nearly~$\nf13$). Although this value of~$\xi$ is admittedly well outside the perturbative range, plugging it into~\eqref{e.physics.predictions} does yield~$t^3$, in agreement with Richardson's~$\nf43$ law. More generally, to model pair displacement in a turbulent fluid according to the Kolmogorov-Onsager theory with corrections for intermittency, one should take
\begin{equation*}
	\xi = -\frac23 - \frac13 \mu\,,
\end{equation*}
where~$\mu>0$ is the intermittency parameter, which is the fractal codimension of the active region of the fluid. Inserting this into~\eqref{e.physics.predictions} one obtains the exponent
\begin{equation}
	\frac{4}{2+\xi}
	=
	\frac{12}{4-\mu}\,.
	\label{e.HP}
\end{equation}
This is in exact agreement with the predicted corrections to Richardson's~$\nf43$ law due to the effects of intermittency, derived in~\cite{HP1,HP2}. It is noted in~\cite{BG} that the derivation of Richardson's~$\nf43$ law from Kolmogorov's phenomenological theory, as well as the intermittent corrections to it, are essentially based on dimensional analysis, and thus ``it is remarkable in this respect that~\eqref{e.HP} [...] turns out to be \emph{exact} for this model'' (emphasis in original).

\subsection{Power-law superdiffusivity in the perturbative regime}

In this paper, we prove the physics prediction~\eqref{e.physics.predictions} in the perturbative regime~$0< \cgamma = \frac12(2 -\xi) \ll 1$ by developing a rigorous version of the renormalization group argument.

\smallskip

If~$\cgamma$ is sufficiently small, then, for a typical realization of the vector field~$\f$, the quenched second moment and variance of the displacement of the particle~$X_t$ grow superlinearly in~$t$ with exponent~$\frac{2}{2-\cgamma} = \frac{4}{2+\xi}$, as predicted in~\eqref{e.physics.predictions}, and we are even able to describe the prefactor constant with considerable precision. What we show is roughly that
\begin{equation}
	\mathbf{E}\bigl[ \bigl| X_{t} \bigr|^2 \bigr]
	\,,\
	\mathbf{E}\bigl[ \bigl| X_{t}  - \mathbf{E}[ X_t ] \bigr|^2 \bigr]
	=
	2d \cstar^{\nf12}\cgamma^{-\nf12} t^{\frac{2}{2-\cgamma}}\bigl(
	1 + O(\cgamma^{\nf12}\left| \log \cgamma \right|^{\nf92})
	\bigr)
	\,,
	\qquad
	t \gtrsim  \exp(C \nu^2 \cgamma^{-1})  \,.
	\label{e.our.main.result}
\end{equation}
Here~$\cstar$ is a positive constant which can be explicitly computed in most examples (it denotes the strength of the disorder), the~$O\bigl(\cgamma^{\nf12}\left| \log \cgamma \right|^{\nf92} \bigr)$ term is random (it depends on~$t$ and the realization of~$\f$), and both are uniform in the molecular diffusivity~$\nu$.

\smallskip

Before presenting the precise statement of~\eqref{e.our.main.result}, we briefly describe the structure we assume on the random environment. Rather than working directly with the drift field~$\f$, it is more convenient to parametrize the model in terms of an antisymmetric \emph{stream matrix} $\k(x)$, which determines~$\f$ by~$\f = \nabla\cdot\k$.
We consider~$\k$ to be the primary random field (rather than~$\f$) and for this reason we henceforth let~$\mathbf{E}^{\k,x_0}$ denote the expectation with respect to the law of the trajectories, rather than~$\mathbf{E}^{\f,x_0}$.

\smallskip

We do not assume that~$\k$ is Gaussian, and none of our arguments uses Gaussian calculus; the hypotheses are formulated in terms of a \emph{scale decomposition} of~$\k$ into independent pieces indexed by triadic scales. Roughly speaking, we require that~$\k$ can be written as a sum of antisymmetric matrix-valued random fields~$\mathbf{j}_k$ indexed by~$k\in\Z$ with range of dependence and regularity on the scale of~$3^k$ and typical size of order~$3^{k\cgamma}$, in a way that is invariant (in law) under dilations. An antisymmetric matrix-valued Gaussian field with positive Hurst exponent~$\cgamma>0$ fits into this framework, and in that case our assumptions are essentially equivalent to the covariance decay~\eqref{e.our.cov.ass}. However, the class of admissible environments we can treat is much larger. We refer to Subsection~\ref{ss.assumptions} below for the precise formulation of the assumptions and further discussion.

\begin{theoremA}[Quenched power-law superdiffusivity]
\label{t.superdiffusivity}
Let~$\nu,\cgamma \in (0,\infty)$ and~$\cstar\in (0,\infty)$ and define, for each time~$t\in (0,\infty)$, the length scale
\begin{equation}
	\Rt(t)
	\coloneqq
	\bigl(
	(\nu t)^{2-\cgamma}
	+
	\cstar \cgamma^{-1} t^2
	\bigr)^{\frac1{2(2-\cgamma)}}
	\,.
	\label{e.Rt.def}
\end{equation}
Assume that~$\P$ is the law of a shell sequence satisfying the standing scale-decomposition assumptions in Subsection~\ref{ss.assumptions}, including~\ref{a.j.frd},~\ref{a.j.reg},~\ref{a.j.iso} and~\ref{a.j.nondeg}, and let~$\mathbf{P}^{\k,x_0}$ denote, given a realization of~$\k$, the law of the trajectories of the stochastic process starting at~$x_0$. Then there exist constants~$\cgamma_0(d,\cstar)>0$ and~$C(d,\cstar)  < \infty$ such that, if~$\cgamma \leq \cgamma_0$, then we have the following estimates: for every moment~$p \in [1,C^{-1} \cgamma^{-1} \left| \log\cgamma\right|^{-6} ]$,
\begin{equation}
	\left\{
	\begin{aligned}
		&
		\E \Bigl[
		\bigl| \mathbf{E}^{\k,0}\bigl[ \bigl| X_{t} \bigr|^2 \bigr]- 2d \Rt(t)^2\bigr|^p
		\Bigr]^{\nf1p}
		\leq
		C \bigl(p^{\nf12} + \left| \log \cgamma \right|^{\nf12}\bigr)
		\cgamma^{\nf12}\left| \log \cgamma \right|^4 \Rt(t)^2
		\,, \quad \mbox{and}
		\\ &
		\E \Bigl[
		\bigl| \mathbf{E}^{\k,0} \bigl[ X_{t}\bigr]  \bigr|^{2p}
		\Bigr]^{\nf1p}
		\leq
		C \bigl(p + \left| \log \cgamma \right|\bigr) \cgamma \left| \log \cgamma \right|^7
		\Rt(t)^2
		\,.
	\end{aligned}
	\right.
	\label{e.thm.superdiffusivity}
\end{equation}
\end{theoremA}

The first estimate~\eqref{e.thm.superdiffusivity} is the precise version of~\eqref{e.our.main.result}. The powers of~$\left| \log \cgamma \right|$ appearing in~\eqref{e.thm.superdiffusivity}, and the range of~$p$, are those obtained in the machine-checked proof; see Section~\ref{ss.formalization}.
Observe that the length scale~$\Rt(t)$ defined by~\eqref{e.Rt.def} behaves asymptotically like
\begin{equation*}
	\Rt(t) \approx (\nu t)^{\nf12}
	\quad \text{for } t \ll t_*\,,
	\quad \text{and} \quad
	\Rt(t) \approx \bigl( \cstar^{\nf12}\cgamma^{-\nf12} t \bigr)^{\frac{1}{2-\cgamma}}
	\quad \text{for } t\gg t_*
	\,,
	\qquad
	t_* \coloneqq
	( \cstar^{-1} \cgamma \nu^{2-\cgamma})^{\nf1\cgamma}
	\,.
\end{equation*}
This reflects the fact that, on very small scales, the molecular diffusivity will be stronger than the drift and the effective diffusivity of the particle will thus be essentially the same as the driving Brownian motion; while, on large scales, it has the behavior indicated in~\eqref{e.our.main.result}.
The crossover between the two regimes occurs at the time scale~$t \asymp t_* \coloneqq ( \cstar^{-1} \cgamma \nu^{2-\cgamma})^{\nf1\cgamma}$. The point of introducing~$\Rt(t)$ is to allow us to write one formula valid for all times. The estimate~\eqref{e.thm.superdiffusivity} indicates that the~$O(\cgamma^{\nf12}\left| \log \cgamma \right|^{\nf92})$ term in~\eqref{e.our.main.result} is valid at the level of finite moments. As the moment becomes large---beyond~$p \asymp \left| \log \cgamma \right|$---the estimate however degrades. This is a consequence of intermittency, which is to be expected for a multifractal model like ours.

\smallskip

We expect that the estimates in~\eqref{e.thm.superdiffusivity} are sharp in their dependence on~$\cgamma$, up to the factors of~$\left|\log \cgamma\right|$. Heuristically, this is the expected sensitivity in these quantities to resampling a single scale in the scale decomposition of~$\k$, on the order of~$\Rt(t)$.

\smallskip

As part of the proof of Theorem~\ref{t.superdiffusivity}, we also obtain bounds on higher moments of the displacement: there exist~$c(d,\cstar)>0$ and~$C(d,\cstar)<\infty$ such that, for every~$t\in (0,\infty)$,
\begin{equation}
	\label{e.moment.bounds}
	c\Rt(t)
	\leq
	\E\bigl[\mathbf{E}^{\k,0}[|X_t|^p]\bigr]^{\nf1p}
	\leq C p^{\nf12}\Rt(t)
	\qquad \text{for } p \in [2, c\cgamma^{-1} |\log\cgamma|^{-6}]\,.
\end{equation}
See Remark~\ref{r.higher.moments.of.displacement}. As opposed to~\eqref{e.thm.superdiffusivity}, these bounds come without precise information on the prefactor constants.

\subsection{Renormalization of the generator}

The proof of Theorem~\ref{t.superdiffusivity} does not proceed by direct manipulation of the trajectories of the solution of the SDE~\eqref{e.SDE}. Conceptually, our approach is closer to that of Euclidean field theory.  The primary object is the infinitesimal generator of the process~$\{X_t\}$, which is the random divergence-form operator
\begin{equation}
	L\coloneqq \nabla \cdot \a(x)\nabla, \qquad
	\mbox{where} \quad  \a(x) \coloneqq \nu \Id  + \k(x)
	\,.
	\label{e.generator.intro}
\end{equation}
Throughout, $\{X_t\}$ denotes the Markov process generated by~$L$; this is the meaning we give to~\eqref{e.SDE}.
The renormalization group acts on the operator~$L$ by coarse-graining its coefficients on larger and larger cubes and tracking their induced flow. By a detailed quantitative analysis, we are able to control the random fluctuations of these coarse-grained coefficients and show that, at every scale, they stay within a factor of~$1 + O(\cgamma^{\nf12} \left|\log \cgamma\right|^{\nf72})$ of a deterministic, scale-dependent diffusivity. We are then able to analyze the evolution of this single parameter, which we denote by~$\shom_m$, where~$m\in\Z$ is the scale parameter corresponding to length scale~$3^m$. Using the perturbative assumption~$\cgamma\ll1$, we can control the RG iteration and obtain the following estimate for the effective diffusivity at scale~$r$:
\begin{equation}
	\Bigl|
	\frac{\shom_m}{\nueff(3^m)} - 1
	\Bigr|
	\leq
	C \cgamma^{\nf12}  \left|\log \cgamma\right|
	\,,
	\qquad \mbox{where} \quad
	\nueff(r) \coloneqq
	\bigl(\nu^2 +\cstar\cgamma^{-1}r^{2\cgamma} \bigr)^{\nf12}
	\,.
\end{equation}
Our renormalization scheme is Wilsonian in the literal sense: we integrate out the stream matrix scale-by-scale and track the induced flow of the scale-local coarse-grained coefficients; unlike in translation-invariant Euclidean field theory, these coarse-grained coefficients remain genuinely random (quenched) at every scale, and it is only their fluctuation about the deterministic~$\shom_m$ that we control.

\smallskip

This analysis leads to the following theorem, which states that~$L \approx \nueff(r)\Delta$ at scale~$r \in (0,\infty)$ in the sense that the relative difference between solutions of the Dirichlet problems for the corresponding Poisson equations is at most~$O(\cgamma^{\nf12}\left| \log\cgamma\right|^{\nf72})$. Here and throughout the paper we denote triadic cubes by
\begin{equation*}
	\cu_m \coloneqq \Bigl( -\frac12 3^m , \frac12 3^m \Bigr)^{\!d}
	\,, \qquad m\in\Z\,.
\end{equation*}

\begin{theoremA}[Renormalization of the generator]\label{t.homogenization}
Let~$\nu,\cgamma \in (0,\infty)$ and~$\cstar\in (0,\infty)$. Assume that~$\P$ is the law of a shell sequence satisfying the standing scale-decomposition assumptions in Subsection~\ref{ss.assumptions}, including~\ref{a.j.frd},~\ref{a.j.reg},~\ref{a.j.iso} and~\ref{a.j.nondeg}. Then there exist~$\cgamma_0(d,\cstar)>0$ and~$C(d,\cstar) < \infty$ such that, if~$\cgamma \leq \cgamma_0$, then the following statement is valid.
For each~$m\in\Z$, there exist an effective diffusivity~$\shom_m \in (0,\infty)$ and a nonnegative random variable~$\EthmB(m)$ which satisfy
\begin{equation}
	\bigl|
	\shom_m
	-
	\bigl(\nu^2 +\cstar\cgamma^{-1}3^{2\cgamma m} \bigr)^{\nf12}
	\bigr|
	\leq
	C \cgamma^{\nf12}
	\left|\log \cgamma\right|
	\shom_m
	\label{e.shom.m.bounds.thmB}
\end{equation}
and, for every~$p \in[1, C^{-1} \cgamma^{-1} \left| \log\cgamma \right|^{-1} ]$,
\begin{equation}
	\E \bigl[ \EthmB(m)^{p} \bigr]^{\nf1p}
	\leq
	C \bigl(p^{\nf12} + \left| \log \cgamma \right|^{\nf12}\bigr) \cgamma^{\nf12}
	\left| \log \cgamma \right|^3
	\,,
	\label{e.mathcal.E.estimate}
\end{equation}
such that, for every~$\g \in C^{0,\nf12}(\cu_m;\Rd)$,~$h\in C^{1,\nf12}( \cu_m)$, and pair~$u,v\in H^1(\cu_m)$ satisfying
\begin{equation}
	\left\{
	\begin{aligned}
		& -\nabla\cdot (\nu \Id + \k)  \nabla u = \nabla \cdot \g
		& \mbox{in} & \ \cu_m\,, \\
		& u = h & \mbox{on} & \ \partial \cu_m\,,
	\end{aligned}
	\right.
	\qquad
	\mbox{and}
	\qquad
	\left\{
	\begin{aligned}
		& -\shom_m \Delta v = \nabla \cdot \g
		& \mbox{in} & \ \cu_m\,, \\
		& v = h & \mbox{on} & \ \partial \cu_m\,,
	\end{aligned}
	\right.
	\label{e.intro.dirichlet.pairs}
\end{equation}
we have the estimate
\begin{equation}
	\label{e.intro.homogenization.estimate}
	3^{-m} \| u - v \|_{L^\infty(\cu_m)}
	\leq
	\EthmB(m)
	\Bigl(
	\shom_m^{-1} 3^{\frac12 m}
	[ \g ]_{\underline{W}^{\nf12,\infty}(\cu_m)}
	+
	3^{\frac12 m} \| \nabla h \|_{\underline{W}^{\nf12,\infty}(\cu_m)}
	\Bigr)
	\,.
\end{equation}
Moreover,
\begin{equation}
	\biggl|
	\fint_{\cu_m} \nu | \nabla u |^2
	-
	\fint_{\cu_m} \shom_m | \nabla v |^2
	\biggr|
	\leq
	\EthmB(m)
	\Bigl(
	\shom_m^{-\nf12} 3^{\frac12 m}
	[ \g ]_{\underline{W}^{\nf12,\infty}(\cu_m)}
	+
	\shom_m^{\nf12} 3^{\frac12 m} \| \nabla h \|_{\underline{W}^{\nf12,\infty}(\cu_m)}
	\Bigr)^2
	\,.
	\label{e.intro.energies}
\end{equation}
\end{theoremA}

Theorem~\ref{t.homogenization} is the central generator-level result of the paper and Theorem~\ref{t.superdiffusivity} can be viewed as a probabilistic corollary of it.
It says that on a cube of size~$3^m$,
the diffusion behaves like Brownian motion with diffusivity~$\shom_m$ up to a random relative error of~$\EthmB(m)$ that is typically of order~$\cgamma^{\nf12} \left| \log \cgamma \right|^{\nf72}$ and does not improve with scale. The connection between the two statements is best understood through the probabilistic interpretation of the solution of the Dirichlet problem. Fix a scale~$m\in\Z$ and a cube~$\cu_m$ of side length~$3^m$, and consider the Dirichlet problems~\eqref{e.intro.dirichlet.pairs} with data~$\g$ and~$h$.  Then the solution~$u$ of the heterogeneous problem can be written as
\begin{equation*}
	u(x) =
	\mathbf{E}^{\k,x}\biggl[ h\bigl(X_{\tau(\cu_m)}\bigr)  + \int_0^{\tau(\cu_m)} (\nabla \cdot \g)\bigl(X_s\bigr)\,ds \biggr]
	\,,
\end{equation*}
where~$\tau(\cu_m)$ denotes the first exit time of~$X_t$ from~$\cu_m$. By taking~$\g=0$ and varying~$h$, we can learn about the exit measure of the process starting from any point in the domain. Likewise, by taking~$h=0$ and~$\g(x) = x_1 e_1$ so that~$\nabla \cdot \g =1$, the value of~$u(x)$ is precisely the expected exit time for the process starting at~$x$. Theorem~\ref{t.homogenization} says that these values will be the same, up to a relative error of order~$O(\cgamma^{\nf12}\left| \log \cgamma \right|^{\nf72})$ with random prefactor, if we replace the process~$X_t$ by~$\sqrt{2 \nueff(3^m)} W_t$, where~$W_t$ is a standard Brownian motion. In this sense, we are able to show that
\begin{equation*}
	X_t \approx \sqrt{2 \nueff(3^m)} W_t
	\quad \text{in $\cu_m$, \ up to a relative error of $O(\cgamma^{\nf12}\left| \log \cgamma \right|^{\nf72})$\,.}
\end{equation*}
Time scales~$t$ should be related to length scales~$r$ by~$r^2 = \nueff(r) t$, and inverting this formula (up to a small relative error) gives~$r = \Rt(t)$, which motivated the definition of~$\Rt(t)$ in~\eqref{e.Rt.def}. Heuristically, this suggests that, for each fixed~$t \in (0,\infty)$,
\begin{equation*}
	X_t \approx \sqrt{2 \nueff(\Rt(t))} W_t
	\quad \text{up to a relative error of $O(\cgamma^{\nf12}\left| \log \cgamma \right|^{\nf72})$}
\end{equation*}
and, on this heuristic level, we obtain
\begin{equation*}
	\mathbf{E}^{\k,0} \bigl[ \bigl| X_t \bigr|^2 \bigr]
	\,,\
	\mathbf{E}^{\k,0} \bigl[ \bigl| X_t - \mathbf{E}^{\k,0}[ X_t] \bigr|^2 \bigr]
	\approx
	2d\nueff(\Rt(t)) t
	\approx
	2d\Rt(t)^2\,.
\end{equation*}
The rigorous implementation of this argument, with the stated tails, is carried out in Section~\ref{s.theprocess}. The proof of Theorem~\ref{t.homogenization} is presented in Section~\ref{ss.proof.theorem.B}.

\smallskip

A key point is that the error in the harmonic approximation is genuinely \emph{scale-local}.
While the random error~$\EthmB(m)$ has tail bounds which are uniform in~$m$, there can be no \emph{single} random variable that controls all scales simultaneously. Indeed, as we show in Section~\ref{ss.scale.local}, if~$k \gg 1$, then the random variables~$\EthmB(m)$ and~$\EthmB(m+k)$ are  \emph{almost independent} in a quantitative sense. This is a reflection of the genuine multifractal nature of the environment and the corresponding scale-invariance of the renormalization group flow. Since the size of the molecular diffusivity~$\nu$ is irrelevant to the RG picture, as we zoom out from scale~$3^m$ to a larger scale~$3^{m+k}$, we see the same picture except that we have essentially resampled the field at new frequencies.

\subsection{Anomalous regularization}

In addition to superdiffusive transport, the renormalization picture has a second, more analytic consequence: it produces \emph{anomalous regularization} for solutions of the elliptic equation associated with the generator.
Even as~$\nu \to 0$, solutions do not lose all regularity; instead, the emergent diffusivity generated by the drift enforces near-Lipschitz behavior on all scales. Since this regularity is \emph{uniform in the molecular diffusivity}~$\nu$, it cannot be attributed to microscopic smoothing by~$\nu\Delta$. Indeed, in the genuinely multiscale setting considered here, one instead expects solutions to exhibit intermittent, multifractal behavior, with energy concentrating on sets of strictly positive codimension.

\smallskip

The renormalization group analysis reveals that this multifractality is nevertheless constrained: as one passes to larger scales the effective diffusivity grows, and this emergent ellipticity suppresses oscillations and prevents excessively sharp energy concentration. The next theorem quantifies this effect in a quenched, scale-local form.  For each triadic cube~$\cu_m$, we consider solutions of the equation~$-\nabla\cdot \a\nabla u=\nabla\cdot\g$ in~$\cu_m$. We show that, with high probability, such solutions are H\"older continuous on~$\cu_{m-1}$ with any exponent~$\alpha<1$ up to a loss of order~$\cgamma^{\nf12}$, i.e.\ they are \emph{almost Lipschitz} in the perturbative regime~$\cgamma\ll 1$.  The corresponding estimate involves a random minimal scale~$X_m(\alpha)$, which measures how many scales must be discarded before one enters the regular regime;~$X_m(\alpha)$ has exponentially decaying tails, reflecting that irregular behavior arises only from very rare collections of unfavorable scales.

\begin{theoremA}[Anomalous H\"older regularity]\label{t.regularity}
Let~$\nu,\cgamma \in (0,\infty)$ and~$\cstar\in (0,\infty)$. Assume that~$\P$ is the law of a shell sequence satisfying the standing scale-decomposition assumptions in Subsection~\ref{ss.assumptions}, including~\ref{a.j.frd},~\ref{a.j.reg},~\ref{a.j.iso} and~\ref{a.j.nondeg}. Then there exist~$\cgamma_0(d,\cstar)>0$ and~$C(d,\cstar) < \infty$ such that, if~$\cgamma \leq \cgamma_0$ and~$\alpha \in (0,1-C \cgamma^{\nf12}]$, then the following statement is valid.
For every~$m\in\Z$, there exists an~$\N$-valued  random variable~$X_m(\alpha)$ satisfying, for every~$N\in\N$,
\begin{equation}
	\P \bigl[ X_m(\alpha) \geq N \bigr]
	\leq C\exp \Bigl( - \frac{(1-\alpha)^2(N-C)}{C \cgamma} \Bigr)\,,
	\label{e.Xm.alpha.integrability}
\end{equation}
and such that the following holds: for every~$\g \in C^{0,\nf12}(\cu_m;\Rd)$, every~$h \in C^{1,\nf12}(\cu_m)$, and every solution~$u \in H^1(\cu_m)$ of
\begin{equation*}
	\left\{
	\begin{aligned}
		& -\nabla \cdot (\nu \Id  + \k)  \nabla u = \nabla \cdot \g & \mbox{in} & \ \cu_m \,, \\
		& u = h & \mbox{on} & \ \partial \cu_m \,,
	\end{aligned}
	\right.
\end{equation*}
we have, for every~$x \in \cu_m$ and every~$n \in \Z$ with~$n\leq m$ and~$X_m(\alpha)\leq m-n$,
\begin{align}
	\label{e.energy.density.estimate}
	\lefteqn{
	\nu^{\nf12} \| \nabla u \|_{\underline{L}^2((x+\cu_{n}) \cap \cu_m)}
	} \quad &
	\notag \\ &
	\leq
	C 3^{(1-\alpha)(m-n)}
	\Bigl(
	\nu^{\nf12} \| \nabla u \|_{\underline{L}^2(\cu_{m})}
	+
	\shom_m^{-\nf12} 3^{\frac12 m} [ \g ]_{\underline{W}^{\nf12,\infty}(\cu_m)}
	+
	\shom_m^{\nf12} 3^{\frac12 m} \| \nabla h \|_{\underline{W}^{\nf12,\infty}(\cu_m)} \indc_{\{ x \notin \cu_{m-1} \}}
	\Bigr)
	\,.
\end{align}
Here the full normalized boundary-data norm is
\begin{equation*}
	\|\nabla h\|_{\underline{W}^{\nf12,\infty}(\cu_m)}
	=3^{-\frac12m}\|\nabla h\|_{L^\infty(\cu_m)}
	+[\nabla h]_{C^{0,\nf12}(\cu_m)}\,.
\end{equation*}
\end{theoremA}

The estimate~\eqref{e.energy.density.estimate} is a \emph{non-concentration} bound for the energy density of the solution~$u$.
By the Morrey characterization of H\"older space, it implies~$C^{0,\alpha}$ regularity of~$u$. Indeed, as we will demonstrate, under the hypotheses of the theorem, we have the estimate
\begin{equation}
	\label{e.local.Holder.estimate}
	3^{\alpha m}[ u ]_{C^{0,\alpha}(\cu_{m-1}) }
	\leq
	C 3^{(1-\alpha)X_m(\alpha)}
	\bigl(
	\| u - (u)_{\cu_{m}}  \|_{\underline{L}^2(\cu_{m})}
	+
	\shom_m^{-1} 3^{\frac32 m} [ \g ]_{\underline{W}^{\nf12,\infty}(\cu_m)}
	\bigr)
	\,.
\end{equation}
The bound~\eqref{e.energy.density.estimate} asserts that the energy density behaves as though it had effective dimension at least~$d-2(1-\alpha)$. Equivalently, any set carrying a nontrivial fraction of energy must have fractal codimension at most~$2(1-\alpha)$.
Since~$\alpha$ may be taken up to~$1-C\cgamma^{\nf12}$, this yields the quantitative restriction on the severity of energy intermittency that the fractal codimension is at most~$C\cgamma^{\nf12}$.
As multifractal heuristics suggest that energy should indeed concentrate on sets of codimension~$\asymp c\cgamma^{\nf12}$ (in an~$L^\infty$ sense), this indicates that the near-Lipschitz H\"older exponent of~$1-C\cgamma^{\nf12}$ as well as the corresponding non-concentration estimate are sharp, up to prefactor constants.

\subsection{Statement of the assumptions}
\label{ss.assumptions}

We next present the precise assumptions governing our main results, Theorems~\ref{t.superdiffusivity},~\ref{t.homogenization} and~\ref{t.regularity}, stated above.

\smallskip

We let~$\nu>0$ denote a small parameter called the \emph{molecular diffusivity}. We fix a small parameter
\begin{equation*}
	\cgamma \in (0,\nf14]
\end{equation*}
and consider a sequence~$\{ \ldots, \mathbf{j}_{-2}, \mathbf{j}_{-1}, \mathbf{j}_0,\mathbf{j}_1,\mathbf{j}_2, \ldots \}$ of independent, antisymmetric matrix-valued random fields satisfying, for every~$k \in \Z$,
\begin{equation}
	\label{e.diff.law.shift}
	\mathbf{j}_k \quad \mbox{and} \quad
	3^{\cgamma k} \mathbf{j}_0 (3^{-k}\cdot)
	\quad \mbox{have the same law}
\end{equation}
and satisfying the following properties:
\begin{enumerate}[label=(\textbf{J\arabic*})]
\setcounter{enumi}{0}

\item
\label{a.j.frd} The field~$\mathbf{j}_0$ is a mean-zero,~$\Rd$-stationary random field with range of dependence~$\sqrt{d}$.

\item
\label{a.j.reg}
Local regularity: with probability one,~$\mathbf{j}_0$ belongs to the space~$C^{1,1}_{\mathrm{loc}} (\Rd; \R^{d\times d}_{\mathrm{skew}})$ of antisymmetric matrices with entries which belong to~$C^{1,1}(B_r)$ for every~$r>0$, and
\begin{equation}
	\label{e.kn.reg.ass}
	\P \Bigl[
	\| \mathbf{j}_0 \|_{L^{\infty}(\cu_0)}
	+
	\sqrt{d} \| \nabla \mathbf{j}_0 \|_{L^{\infty}(\cu_0)}
	+
	d \| \nabla^2 \mathbf{j}_0  \|_{L^{\infty}(\cu_0)}
	> t
	\Bigr]
	\leq \exp(-t^2)
	\quad \forall t \in [1, \infty) \,\,.
\end{equation}

\item
\label{a.j.iso} Hypercube (hyperoctahedral) symmetry: the joint law of~$\{ \mathbf{j}_n\}_{n\in\Z}$
is invariant under negation and the symmetries of a~$d$-dimensional hypercube: reflections and permutations across the coordinate planes. Specifically, for every matrix~$R$ with exactly one~$\pm 1$ in each row and column and~$0$s elsewhere, the law of the conjugated sequence~$\{ R^t \mathbf{j}_n(R \cdot) R\}_{n\in\Z}$ is the same as that of~$\{ \mathbf{j}_n \}_{n\in\Z}$, and the law of the sequence~$\{ -\mathbf{j}_n\}_{n\in\Z}$ is the same as that of~$\{ \mathbf{j}_n \}_{n\in\Z}$.

\item
\label{a.j.nondeg}
Non-degeneracy: there exists a positive constant~$\cstar > 0$ such that, for every~$e\in\partial B_1$,
\begin{equation}
	\label{e.nondeg}
	\E \bigl[ \bigl| \bigl(\nabla \Delta^{-1} \nabla \cdot (\mathbf{j}_0 e)\bigr)(0) \bigr|^2  \bigr] = \cstar (\log 3)  \,\,.
\end{equation}
\end{enumerate}

Throughout the paper, we denote the law of the sequence~$\{ \mathbf{j}_n \}_{n\in\Z}$ by~$\P$ and the corresponding expectation by~$\E$.
We also distinguish the normalized Frobenius norm constant
\begin{equation*}
	\cstar^{+}
	\coloneqq
	(\log 3)^{-1}
	\E\biggl[
	\sum_{i,j=1}^d
	|( \mathbf j_0)_{ij} (0) |^{2}
	\biggr]
	\,.
\end{equation*}
The directional nondegeneracy in~\eqref{e.nondeg} and the local regularity assumption imply~$d\cstar\leq \cstar^{+}\leq C(d)$. Note that~\eqref{e.kn.reg.ass} and~\eqref{e.nondeg} together imply~$\cstar \leq \frac32$.

\smallskip

The stream matrix~$\k(\cdot)$ is given by the formal sum
\begin{equation}
	\label{e.k.sum.def}
	\k(x) = \sum_{n=-\infty}^\infty
	\mathbf{j}_n (x)\,.
\end{equation}
Strictly speaking, this sum does not converge; it is the differences~$\k(x) - \k(y)$ which converge.
Since only~$\nabla \cdot \k$ enters, it is naturally defined modulo constants; our scale decomposition produces such a field, and the generator is well-defined even though the full multiscale sum does not converge pointwise. This makes no difference from the point of view of the process~$\{ X_t \}$, since the addition of a constant antisymmetric matrix to the coefficient field has no effect on the generator~$\nabla \cdot \a(x)\nabla$ nor on the law of the Markov process~$\{ X_t\}$.

\smallskip

To compare with~\eqref{e.our.cov.ass}, observe that these assumptions imply the decay upper bound
\begin{equation*}
	\bigl| \cov\bigl[ \nabla \k(x) , \nabla \k(y) \bigr] \bigr| \lesssim (1-\cgamma)^{-1} |x-y|^{-2 (1-\cgamma)}   \,.
\end{equation*}
Note that we do not assume full rotational invariance; the discrete hypercube symmetries suffice for our conclusions. The main role of~\ref{a.j.iso} is to ensure that the effective diffusivity is a scalar.
The negation symmetry forces~$\a$ and~$\a^t$ to have the same law, which turns out to be very convenient.

\subsection{Overview of the proofs}

We next present an outline of the arguments leading to the proofs of Theorems~\ref{t.superdiffusivity},~\ref{t.homogenization} and~\ref{t.regularity}. The paper is divided into four main steps. The central pillar is the second step, which is a Wilsonian RG analysis of the generator.

\begin{enumerate}
\item \textbf{Deterministic coarse-graining theory (Section~\ref{s.coarse.graining.theory}).} Following~\cite{AK.HC}, we develop quantitative, scale-local estimates for general (possibly nonsymmetric) divergence-form operators, including coarse-grained ellipticity/sensitivity bounds and an inhomogeneous theory that converts coarse-scale information into quantitative harmonic approximation on a given cube.

\item \textbf{Renormalization group flow (Section~\ref{s.RG.flow}).}
Using the multiscale structure of the environment, we define infrared truncations of the operator and a running effective diffusivity~$\shom_m$ at scale~$3^m$. We prove by induction in the scale that, on cubes of size~$3^m$, the truncated operator is well-approximated by~$-\shom_m \Delta$, in the sense that the coarse-grained coefficient fields at scale~$3^m$ are close to~$\shom_m \Id$. The difference is measured by a scale-local defect~$\EthmB(m)$, a random variable that is typically of order~$\cgamma^{\nf12}\left| \log \cgamma \right|^{\nf72}$. Using perturbative arguments in the parameter~$\cgamma$, we show that the sequence~$(\shom_m)_{m\in\Z}$ of effective diffusivities satisfies an approximate recurrence. Integrating this, we are able to propagate the asymptotics of~$\shom_m$ to larger scales.

\item \textbf{Anomalous regularization (Section~\ref{s.regularity}).}
We combine the generator approximation with an excess-decay iteration to obtain large-scale~$C^{0,\alpha}$ estimates with~$\alpha = 1 - C\cgamma^{\nf12}$, uniform in the microscopic diffusivity, together with quantitative control of the (random) minimal scale beyond which the regularity holds (Theorem~\ref{t.regularity}). Combining this estimate with the results of the previous two steps, we obtain~$L^\infty$ homogenization estimates at every scale (Theorem~\ref{t.homogenization}).

\item \textbf{Superdiffusion of particle trajectories (Section~\ref{s.theprocess}).}
We translate the generator-level renormalization into estimates for the stochastic process via probabilistic representation formulas for Dirichlet problems for the generator. This yields sharp control of exit times, displacement tails and first and second moments of the stopped process. Choosing a large confinement scale and using the displacement tail bounds, we can remove the stopping time from the representation formulas, which gives us the superdiffusive estimates stated in Theorem~\ref{t.superdiffusivity}.
\end{enumerate}

\noindent
The appendices collect probabilistic and analytic inputs (concentration inequalities, percolation estimates and functional inequalities for Besov spaces) needed to implement the induction and to control atypical environments. Below we give more details, still at a heuristic level, of each of the four main steps listed above.

We emphasize that, while techniques from our previous works~\cite{AKMBook,AK.HC,ABK.SD} inform the ideas here, the arguments in this paper are self-contained. There is no appeal to any of the main results from these works.

\subsubsection{Deterministic coarse-graining theory (Section~\ref{s.coarse.graining.theory})}

Section~\ref{s.coarse.graining.theory} develops a \emph{deterministic} coarse-graining theory for general uniformly elliptic coefficient fields~$\a = \s + \k$. The key objects are \emph{coarse-grained matrices}~$\s(U)$,~$\s_*(U)$,~$\k(U)$ associated with each bounded domain~$U$, which encode the effective linear response of the operator~$-\nabla \cdot \a \nabla$ when restricted to~$U$. These matrices, and the associated \emph{coarse-grained ellipticity constants}~$\lambda_{s,q}(\cu_m)$ and~$\Lambda_{s,q}(\cu_m)$, measure the effective diffusive behavior at each spatial scale without reference to microscopic parameters.

\smallskip

The central insight is that while the microscopic diffusivity~$\nu$ may be arbitrarily small---and eventually disappears in the critical limit~$\cgamma \to 0$---the coarse-grained ellipticity at scale~$3^m$ can be much larger: it captures the emergent diffusivity produced by the interaction of the flow with the random drift at that scale. We can prove a full slate of \emph{coarse-grained functional inequalities}, which generalize key elliptic estimates by replacing the constants of uniform ellipticity with~$\lambda_{s,q}(\cu_m)$ and~$\Lambda_{s,q}(\cu_m)$. This allows all estimates to be expressed in terms of scale-local quantities rather than the bare parameter~$\nu$. In the language of renormalization, the coarse-grained ellipticity constants are the ``state variables'' that we track along the RG flow; the microscopic~$\nu$ is integrated out and does not appear in any of the scale-to-scale estimates.

\smallskip

The definitions of the coarse-grained matrices and their basic properties originate from the theory of quantitative homogenization~\cite{AKMBook}.  The recent work~\cite{AK.HC} introduced the notion of \emph{coarse-grained ellipticity} and extended this framework to high-contrast elliptic homogenization---the regime where the ellipticity ratio can degenerate. The present paper requires two extensions of the theory developed in~\cite{AK.HC} which are developed here and are of independent interest.

\paragraph{Coarse-grained sensitivity estimates.}
We quantify how the coarse-grained matrices and ellipticity constants respond to antisymmetric perturbations~$\h$ of the field (the notation of Section~\ref{ss.sensitivity}). When we ``inject'' a fresh layer~$\mathbf{j}_m$ of the stream matrix, we measure its effect in units of the running diffusivity~$\shom_m$, not~$\nu$. Since~$|\mathbf{j}_m| \simeq 3^{\cgamma m}$ while~$\shom_m \simeq \cgamma^{-\nf12} 3^{\cgamma m}$, the relative perturbation is of order~$\cgamma^{\nf12}$---small enough to be perturbative. If we measured instead in units of~$\nu$, the perturbation would be catastrophically large and our induction would not close.

\paragraph{Equations with nonzero right-hand side.} We extend the coarse-graining machinery to handle inhomogeneous equations~$-\nabla \cdot \a \nabla u = \nabla \cdot \g$ with forcing~$\g \neq 0$. This requires working with the~$2d \times 2d$ block matrix formulation (controlling gradients and fluxes simultaneously) and leads to quantitative homogenization estimates: given a solution~$u$ of the heterogeneous equation and a solution~$v$ of the constant-coefficient problem~$-\s_0 \Delta v = \nabla \cdot \g$ with the same boundary data, we bound~$\| \nabla u - \nabla v \|$ in negative Sobolev norms purely in terms of the coarse-grained ellipticity constants and the regularity of~$\g$. These ``black-box'' homogenization estimates are the bridge between the RG analysis of Section~\ref{s.RG.flow} and the regularity theory of Section~\ref{s.regularity}.

\smallskip

We emphasize that Section~\ref{s.coarse.graining.theory} is entirely deterministic: the coefficient field~$\a$ is arbitrary and no probabilistic structure is used. The randomness enters only in Section~\ref{s.RG.flow}, where the deterministic estimates are applied to the specific random field~$\a = \nu \Id + \k$ given in the assumptions stated in Section~\ref{ss.assumptions}. The philosophy is to separate the analytic machinery (which is robust and reusable) from the probabilistic inputs (which are model-specific).

\subsubsection{The RG flow (Section~\ref{s.RG.flow})}

In the supercritical regime ($\cgamma>0$) the RG flow is genuinely multifractal: as we zoom out we do not approach a single limiting homogenized operator, but rather see the same mechanism repeat across scales, with the effective diffusivity growing like a power so that no microscopic ellipticity scale is privileged. Consequently, the iteration must be closed purely in coarse-grained units---sensitivity estimates, functional inequalities, and infrared-cutoff errors are all normalized by the \emph{effective} diffusivity at the current scale, which we denote by~$\shom_m$, rather than by the bare molecular diffusivity~$\nu$.

\smallskip

At a high level, the RG state variables are (i) the deterministic scalar effective diffusivity~$\shom_m$ associated with the infrared-truncated coefficient field~$\a_m$, and (ii) a random \emph{defect} variable~$\EthmB(m)$ which in rough terms quantifies how well, on the cube~$\cu_m$, the heterogeneous operator~$\nabla\cdot \a_m\nabla$ is approximated by the constant-coefficient operator~$\nabla\cdot \shom_m\nabla=\shom_m\Delta$. The main proposition describing the \emph{RG flow} gives this defect good tails, uniformly in~$m$, and shows that~$\shom_m$ follows an essentially deterministic growth law.

\smallskip

A key point is that we control this defect \emph{at the same scale as the cutoff}: we compare~$\nabla\cdot \a_m\nabla$ to~$\shom_m\Delta$ on~$\cu_m$, rather than homogenizing~$\a_m$ only on much larger cubes to exploit decorrelation beyond~$3^m$. In this sense the infrared cutoff is primarily a bookkeeping device: once we have a scale-local approximation for~$\a_m$ on~$\cu_m$, the coarse-grained sensitivity estimates allow us to reinsert the longer wavelengths~$\mathbf{j}_n$ with~$n>m$ at a cumulative relative cost of order~$\cgamma^{\nf12}$.

\smallskip

The heuristic is that, from the point of view of~$\cu_m$, the small-scale modes~$\mathbf{j}_n$ with~$m-n\gg 1$ have already been ``integrated out'': their cumulative effect is absorbed into the running diffusivity~$\shom_m$, while the new randomness at scale~$3^m$ comes from the fresh layer~$n\approx m$, which is independent of the previously integrated modes. This layer has typical amplitude~$\simeq 3^{\cgamma m}$, whereas the running diffusivity is already of size~$\shom_m\simeq \cgamma^{-\nf12}3^{\cgamma m}$ (in the disorder-dominated regime), so the new contribution enters at relative size~$\simeq \cgamma^{\nf12}$. This is precisely what makes the iteration perturbative in the Wilsonian sense: the small parameter is the relative size of one newly added shell compared to the running diffusivity, rather than a homogenization-style scale-separation parameter.

\smallskip

Because the environment is self-similar and the new shells are (nearly) independent across scales, there is no terminal scale at which one can ``close'' by appealing to a distinguished microscopic ellipticity parameter: as we zoom out we essentially see the same picture again, with the field resampled at new frequencies.

\smallskip

We implement the induction by repeatedly alternating:~(i) a coarse-graining step that contracts the defect,~(ii) an injection step that turns on the missing layers, and (iii) a stochastic integrability-repair step based on the insensitivity of the coarse-grained matrices to rare, isolated ``bad cubes'' where the environment behaves erratically.

\smallskip

\paragraph{Contraction of the defect.}
Consider the infrared cutoff~$\a_L$ (see~\eqref{e.infrared.cutoff.def}) of the field and coarse-grain it in larger cubes~$\cu_m$ with~$m - L \gg 1$. Because~$\a_L$ has range of dependence~$\lesssim 3^L$, the coarse-grained matrices in~$\cu_m$ should be closer to their homogenized limit by the geometric factor~$3^{-\theta(m-L)}$ for some~$\theta>0$, relative to the assumed size of the defect for~$\a_L$ in~$\cu_L$. The difficulty is to make this quantitative in a way that is compatible with the multifractal regime. This is where the ideas and methods from~\cite{AK.HC} enter: using the coarse-graining theory, we are able to measure the contraction of the defect using the \emph{coarse-grained ellipticity} at the working scale (not the bare microscopic~$\nu$).
As a result, we show that, if the defect is small at the cutoff scale, then after climbing~$O(\left|\log \varepsilon\right|)$ triadic steps it typically shrinks by a factor~$\varepsilon$ (modulo some exceptional bad events of small probability). This argument appears in Section~\ref{ss.homogenization.step}.

\smallskip

The RG step contains many estimates that inflate constants. If we want our induction to close, we need to obtain the same estimate with the same constant at the new scale. Therefore, we must locate one mechanism that contracts \emph{something}, strongly enough to pay for all of that inflation. Here, that mechanism is the contraction of the coarse-graining defect, in other words, the contraction of the  homogenization error.

\paragraph{Injection of the fresh layer.} After the contraction step, we are typically sitting on a cube~$\cu_m$ but still with the coefficient field~$\a_L$ for some~$L < m$. To get back to the real object at scale~$m$, we have to turn on the layers represented by~$\{ \mathbf{j}_{L+1},\ldots,\mathbf{j}_m\}$. This is where the \emph{coarse-grained sensitivity estimates} from Section~\ref{ss.sensitivity} play an important role. These allow us to compare the coarse-grained objects for~$\a_L$ and~$\a_m$, while measuring perturbations in units of the \emph{running diffusivity at the current scale}, rather than in units of~$\nu$. That is the only normalization under which adding a layer of size~$\simeq 3^{\cgamma m}$ produces a relative effect~$\simeq \shom_m^{-1} 3^{\cgamma m} \approx \cgamma^{\nf12}$ instead of something catastrophically large. This is also why we need the lower bounds on~$\shom_m$ to live inside the induction hypothesis: we can only afford to insert the layers back if the running diffusivity is large enough to absorb them. We can also show that updating the homogenized matrix from~$\shom_L$ to~$\shom_m$ produces the same error by a similar mechanism. The injection step appears in Section~\ref{ss.buckle}, when we close the induction.

\smallskip

Taken together, these two steps lead to the heuristic that we should expect~$D_{m+h} \leq \frac12 D_m + C \cgamma^{\nf12}$, where~$D_m$ is the ``defect size'' at scale~$3^m$ and~$h\asymp C|\log\cgamma|$ is the number of triadic scales needed to realize the contraction. This is consistent with~$D_m \leq C \cgamma^{\nf12}$ uniformly in~$m$, and it is indeed essentially what we prove.

\smallskip

But another difficulty needs to be resolved: the RG iteration is not just about keeping the size of the defect small; we also need to keep its stochastic integrability under control as we iterate. Naively, rare bad events at each scale accumulate and destroy the stretched-exponential tails we need to close the induction in a strong Orlicz norm.

\paragraph{Improvement of stochastic integrability.} In Section~\ref{ss.coarse.grained.sensitivity}, in a preemptive move, we improve the stochastic integrability of the coarse-grained ellipticity constants (and thus the tails of the defect) that enter the contraction and injection steps. The price to pay is that we have to enlarge the cube by~$C\left| \log\cgamma\right|$ more triadic scales, and inflate the constants in our estimates. The improvement is accomplished by demonstrating that the coarse-grained matrices are insensitive to outlier, ``bad'' regions of the environment. We show this by adapting the usual subadditivity argument, constructing curl-free and divergence-free fields which are piecewise affine, and constant on a chosen ``bad'' subset of smaller cubes. This effectively replaces the full sum in the subadditivity inequality by a ``trimmed'' sum. The gain in stochastic integrability comes from the simple fact that a trimmed sum of i.i.d.~random variables satisfies much better tail estimates compared to the original sum. Implementing this argument requires that ``chains'' of bad cubes with large diameter do not form, and we rule this out using a percolation-type argument and the approximate locality (in space and in the scale) of our defect random variables.

\paragraph{Closing the induction loop.}
The contraction step improves the deterministic size of the defect but tends to worsen integrability constants; the percolation/trimmed sum step improves integrability but inflates the deterministic size by a large constant. The injection step makes everything worse but gives us back the correct field. As we show in Section~\ref{ss.buckle}, the induction loop closes because these three moves can be alternated (and quantitatively balanced) so that neither the deterministic size, nor the tail parameter drifts off to infinity.

\paragraph{Propagating bounds on the running effective diffusivity.}
Finally, once the defect is under control uniformly in~$m$, we can treat the RG flow of~$\shom_m$ as essentially governed by a one-step perturbative computation which is quite explicit. This yields a one-step recurrence for~$\shom_m$; at a schematic level it takes the form
\begin{equation*}
	\shom_{n+h}
	=
	\shom_{n} + (\log 3)\cstar \shom_{n}^{-1}
	\sum_{k=n+1}^{n+h}3^{2\cgamma k}
	+
	\text{ small error}\,.
\end{equation*}
leading to the sharp asymptotic~$\shom_m \approx (\nu^2 + \cstar \cgamma^{-1} 3^{2\cgamma m})^{\nf12}$ up to relative error~$O(\cgamma^{\nf12}\left| \log\cgamma \right|)$. This step must also be included in the induction loop, for reasons already discussed (we need the lower bound in~$\shom_m$ we get from it), and so it appears in Section~\ref{ss.propagate.diffusivity.bounds}.

\subsubsection{Anomalous regularization (Section~\ref{s.regularity})}

The RG analysis developed in Section~\ref{s.RG.flow} gives us quantitative control of the coarse-grained homogenization error at every scale. In Section~\ref{s.regularity}, we convert this into H\"older regularity estimates for solutions of the equation~$-\nabla \cdot \a \nabla u = 0$, culminating in Theorem~\ref{t.regularity}, followed by the full statement of Theorem~\ref{t.homogenization}. The logical structure is slightly subtle: the excess decay iteration underlying Theorem~\ref{t.regularity} already uses a form of homogenization---specifically, an~$L^2$-type comparison between the heterogeneous solution and a harmonic approximation, valid for homogeneous equations ($\g = 0$). This weaker homogenization estimate follows directly from the coarse-graining machinery of Section~\ref{s.coarse.graining.theory} and the results of the RG analysis in Section~\ref{s.RG.flow}, without any regularity input. The role of the regularity theory is to \emph{upgrade} this to the full~$L^\infty$ estimate of Theorem~\ref{t.homogenization}, which applies to inhomogeneous equations and measures the error in~$L^\infty$ rather than~$L^2$.

\smallskip

The key observation is that, while we cannot expect solutions to be regular on \emph{every} cube (some scales will have anomalously large homogenization error), the \emph{density of bad scales} is small. This is quantified by a ``good event''~$\mathcal{G}(m;s,\varepsilon)$ on which the homogenization error is deterministically bounded by~$\varepsilon$. We show in Section~\ref{ss.scale.local} that these good events satisfy an \emph{approximate scale locality} property:~$\mathcal{G}(m;s,\varepsilon)$ and~$\mathcal{G}(n;s,\varepsilon)$ are nearly independent when~$|m-n|$ is large. This independence comes from the scale decomposition of~$\k$ into independent pieces~$\{\mathbf{j}_k\}_{k\in\Z}$. As a consequence, the proportion of ``bad scales''---those on which the good event fails---can be made as small as~$O(\cgamma^{\nf12})$. The resulting H\"older exponent of~$1 - C\cgamma^{\nf12}$ is expected to be optimal, in view of multifractal heuristics.

\smallskip

The heart of Section~\ref{s.regularity} is an \emph{excess decay iteration} (Section~\ref{ss.excess.decay}). For a solution~$u$, we define the ``excess'' at scale~$3^n$ to be, roughly, the~$L^2$-distance of~$u$ from being affine on~$\cu_n$. We show that, on good scales, the excess contracts by a fixed factor (roughly~$3^{-\nf12}$) when passing from scale~$3^n$ to scale~$3^{n-k}$. On bad scales, we merely have non-expansion. Two effects combine to limit the H\"older exponent to~$\alpha \leq 1 - C\cgamma^{\nf12}$. First, even on good scales, the harmonic approximation error is of order~$\cgamma^{\nf12}$---it does not contract to zero---so the iteration accumulates a scale-independent additive error at each step, which is incompatible with Lipschitz regularity. Second, the density of bad scales is itself of order~$\cgamma^{\nf12}$, further degrading the exponent. The constraint~$\alpha \leq 1 - C\cgamma^{\nf12}$ arises from requiring that the iteration closes despite both of these effects.

\smallskip

The mechanism underlying excess decay is the same one that drives the RG flow: on a good scale, solutions of~$-\nabla \cdot \a \nabla u = 0$ are well-approximated by harmonic functions (solutions of~$-\shom_m \Delta v = 0$), and harmonic functions have excellent regularity. The smallness of the defect~$\EthmB(m)$, established in the RG analysis, translates directly into excess decay via the harmonic approximation lemma (Lemma~\ref{l.harmonic.approximation.good.scales}).

\smallskip

We then define the \emph{minimal scale}~$X_m(\alpha)$ to be the number of scales one must descend below~$m$ before the iteration (centered at every possible point) has accumulated enough good scales to guarantee H\"older regularity with exponent~$\alpha$. By the approximate independence of good events, the minimal scale has the exponentially decaying tails in~\eqref{e.Xm.alpha.integrability}. The energy density estimate~\eqref{e.energy.density.estimate} of Theorem~\ref{t.regularity} is then a direct consequence of the iterated excess decay.

\smallskip

Finally, in Section~\ref{ss.proof.theorem.B}, we prove the full statement of Theorem~\ref{t.homogenization}---the~$L^\infty$ estimate for inhomogeneous equations---by combining the regularity estimates of Theorem~\ref{t.regularity} with the general coarse-graining machinery developed in Section~\ref{s.coarse.graining.theory}. The idea is to use Theorem~\ref{t.regularity} to control the energy density at small scales---this bounds the right side of Proposition~\ref{p.general.coarse.graining} when~$p>2$---and then to convert the resulting~$\Wminusul{-s}{4d}$ gradient bound into an~$L^\infty$ bound on~$u-v$ by the fractional Poincar\'e inequality and the embedding~$W^{1-s,4d}\hookrightarrow L^\infty$, which is valid since~$1-s>d/(4d)=1/4$.

\subsubsection{Quenched superdiffusivity (Section~\ref{s.theprocess})}

In Section~\ref{s.theprocess}, we pass from the generator-level estimates of Theorem~\ref{t.homogenization} to estimates on the stochastic process~$\{X_t\}$ itself, completing the proof of Theorem~\ref{t.superdiffusivity}.

\smallskip

The connection between the two is the probabilistic interpretation of the Dirichlet problem: for a solution~$u$ of~$-\nabla \cdot \a \nabla u = \nabla \cdot \g$ with boundary data~$h$ and smooth~$\g$, we have the probabilistic representation
\begin{equation*}
	u(x) = \mathbf{E}^{\k,x}\biggl[ h(X_{\tau(\cu_m)}) + \int_0^{\tau(\cu_m)} (\nabla \cdot \g)(X_s)\, ds \biggr]\,.
\end{equation*}
Theorem~\ref{t.homogenization} then says that~$u(x)$ is close to the corresponding quantity for a Brownian motion with diffusivity~$\shom_m$, up to a relative error of order~$\EthmB(m) = O(\cgamma^{\nf12}|\log\cgamma|^{\nf72})$.

\smallskip

Moments of the stopped process~$X_{t \wedge \tau(\cu_m)}$ can be computed by choosing appropriate data~$\g$ and~$h$. For instance, taking~$\g = 0$ and~$h(x) = e \cdot x$ gives access to~$\mathbf{E}^{\k,0}[e \cdot X_{\tau(\cu_m)}]$, while taking~$h(x) = (e \cdot x)^2$ and the appropriate right-hand side gives access to~$\mathbf{E}^{\k,0}[(e \cdot X_{\tau(\cu_m)})^2]$. In both cases, It\^o's formula combined with Theorem~\ref{t.homogenization} yields quantitative bounds.

\smallskip

To pass from the stopped process to the unstopped process~$X_t$, we choose the confinement scale~$m$ to depend on~$t$ in such a way that the cube~$\cu_m$ is larger than the intrinsic length scale~$\Rt(t)$ by a factor~$\left|\log\cgamma\right|^{\nf12}$---specifically, we take~$3^m \asymp \left|\log\cgamma\right|^{\nf12} \Rt(t)$. This choice is tuned so that the effective Brownian exit time~$3^{2m}/\shom_m$ is~$t \left| \log \cgamma\right|$ times a large constant. We then show in Section~\ref{ss.exit.time.bounds} that, due to the extra factor of~$\left| \log \cgamma\right|$ and the large constant, the probability of ``early exit''~$\mathbf{P}^{\k,0}[\tau(\cu_m) \leq t]$ is extremely small, of order~$\cgamma^{100}$. The argument proceeds in several steps:
\begin{itemize}
\item First, we use Theorem~\ref{t.homogenization} to show that the expected exit time~$\mathbf{E}^{\k,0}[\tau(\cu_n)]$ is comparable to the Brownian exit time~$\Tr(3^n) \coloneqq 3^{2n}/\nueff(3^n)$.
\item Second, we use a Paley-Zygmund argument to convert this into a lower bound on the survival probability:~$\mathbf{P}^{\k,0}[\tau(\cu_n) \geq c\Tr(3^n)] \geq c$.
\item Third, we iterate across scales using the strong Markov property, together with the percolation estimate of Section~\ref{ss.percolation.on.paths}, which guarantees enough cubes with good exit time bounds along the chain: chaining together~$N \asymp 3^{m-n}$ exit events from cubes of size~$3^n$, we obtain exponential tail bounds for~$\tau(\cu_m)$.
\item Finally, we optimize over the auxiliary scale~$n$ to obtain the displacement tail bound~\eqref{e.early.exit.probability}.
\end{itemize}

\noindent
With the exit probability under control, the error from removing the stopping time is negligible: by the Cauchy-Schwarz inequality and the bound~$|X_{\tau(\cu_m)}|\leq C3^m$,
\begin{align*}
	\lefteqn{ \bigl| \mathbf{E}^{\k,0}[|X_t|^2] - \mathbf{E}^{\k,0}[|X_{t \wedge \tau(\cu_m)}|^2] \bigr| } \qquad & \notag
	\\ & \leq \mathbf{E}^{\k,0}[|X_t|^4]^{\nf12} \, \mathbf{P}^{\k,0}[\tau(\cu_m) \leq t]^{\nf12} + C3^{2m}\mathbf{P}^{\k,0}[\tau(\cu_m) \leq t]
	\leq  C \cgamma^{40} 3^{2m} \,.
\end{align*}
This is much smaller than the homogenization error~$O(\cgamma^{\nf12}\left|\log\cgamma\right|^{\nf92} \Rt(t)^2)$, so the stopped-process estimates transfer directly to the unstopped process, yielding Theorem~\ref{t.superdiffusivity}.

\subsection{Comparison to the borderline case and previous works}
\label{ss.lit.review}

The borderline superdiffusive case~$\xi = 2$ (equivalently~$\cgamma = 0$) corresponds to a drift with a log-correlated stream matrix and leads to the logarithmic-type superdiffusivity predicted on the middle line of~\eqref{e.physics.predictions}, namely a growth of order~$t \sqrt{\log t}$ for the quenched second moment and variance. This prediction was proved rigorously in our previous work~\cite{ABK.SD}, where we obtained the sharper version of the asymptotic (identifying the exact prefactor constant):
\begin{equation}
	\mathbf{E} \bigl[ \bigl| X_t \bigr|^2 \bigr]
	\approx
	\mathbf{E} \bigl[ \bigl| X_t - \mathbf{E}[ X_t] \bigr|^2 \bigr]
	\approx
	2d\cstar^{\nf12} t (\log t)^{\nf 12}
	\,.
	\label{e.ABK.SD}
\end{equation}
We also proved a \emph{central limit theorem} for the process~$\{X_t\}$ with a nonstandard, superdiffusive scaling.
Note that the Gaussian limit is consistent with superdiffusivity only because the rate is very slow---slower than~$t^{1+\ep}$ for every~$\ep>0$.

\smallskip

Prior to that work, the borderline case~$\cgamma=0$ was analyzed in~\cite{TV,CHT,CMOW}. These papers proved growth estimates for the \emph{annealed} second moment consistent with the physics prediction in~\eqref{e.physics.predictions}.
In particular, the works~\cite{CHT,CMOW} established the bound
\begin{equation}
	\E \bigl[ \mathbf{E} \bigl[ \bigl| X_{t} \bigr|^2 \bigr] \bigr]
	=
	O( t \sqrt{\log t} )
	\,,
	\label{e.CMOW}
\end{equation}
up to~$\log\log t$ corrections in the case of~\cite{CHT} which were subsequently removed in~\cite{CMOW}. Very recently, a connection to geometric Brownian motion to analyze intermittent behavior of the averaged Lagrangian coordinate was explored in~\cite{MOW}.

\smallskip

The estimate~\eqref{e.CMOW} should be seen as significantly weaker than~\eqref{e.ABK.SD}. Indeed, as emphasized in~\cite{BG}:~(i) the physics conjecture is \emph{quenched}---the predicted asymptotics should hold for a typical environment, not just after averaging over all environments; and~(ii) even superlinear bounds on the quenched second moment alone do not distinguish between \emph{superdiffusion} (when the quenched variance of~$X_t$ grows superlinearly in~$t$) and \emph{anomalous drift} (when the square of the quenched first moment grows superlinearly but the variance remains diffusive). More concretely, a superlinear lower bound for~$\E[ \mathbf{E}[|X_t|^2]]$ does not tell us whether~$\E[ \mathbf{E}[|X_t - \mathbf{E}[X_t]|^2]]$ is superlinear, or even whether~$\mathbf{E}[|X_t|^2]$ is typically superlinear, since rare environments with anomalously large mean drift could dominate the average.

\smallskip

In this context, we also mention a very recent work~\cite{CMT} which establishes a logarithmic superdiffusive CLT for a critical~2D stochastic Burgers equation, a result which is conceptually very similar to the one of~\cite{ABK.SD}. They formalize an RG flow based on multiscale resolvent estimates for the generator of the SPDE.

\smallskip

In the regime~$\cgamma>0$, there are comparatively few rigorous results for superdiffusion. The only previous work we are aware of is~\cite{KO}, which proves annealed upper and lower bounds, though with exponents that do not match the predicted value.

\smallskip

Theorem~\ref{t.superdiffusivity} asserts a quenched, exponent-matching scaling law in the genuinely superdiffusive regime. It can be viewed as the perturbative continuation of the results of~\cite{ABK.SD} into the regime~$\cgamma>0$. In fact, we can directly compare~\eqref{e.our.main.result} to~\eqref{e.ABK.SD} and observe the~$\sqrt{\log t}$ scaling emerge from the factor of~$\cgamma^{-\nf12}$ in~\eqref{e.our.main.result} by noticing that
\begin{equation*}
	\cgamma^{-1} \asymp \log t \implies
	\bigl( \cgamma^{-\nf12} t \bigr)^{\frac{2}{2-\cgamma}}
	\asymp
	t (\log t)^{\nf12}
	\,.
\end{equation*}
We will not attempt in this paper to re-derive the borderline superdiffusive results of~\cite{ABK.SD} from the present framework---although this would certainly be possible, and the heuristic match is exact at the level of the renormalization group picture. The critically-correlated case lives at the degenerate endpoint~$\cgamma=0$ of the family of strongly correlated drifts treated here, and our power-law superdiffusion result is the genuinely perturbative side of the same RG mechanism. We do not prove a central limit theorem in this regime: unlike the borderline case, where the slow growth of the running diffusivity permits Gaussian fluctuations, power-law superdiffusivity is incompatible with a Gaussian limit. The more natural question in this setting concerns \emph{spontaneous stochasticity}~\cite{BGK,EVE,FGV}: whether, for almost every realization of the drift, the particle trajectories converge to a well-defined limiting process as the molecular diffusivity~$\nu \to 0$ and if so, whether this limit is unique (note that the limiting process will necessarily depend on the realization of the drift). We leave this interesting question for future work.

\smallskip

A related \emph{non-perturbative} RG argument appears in the work~\cite{AV}, which proves anomalous dissipation for a passive scalar equation by iteratively renormalizing the scale-dependent effective diffusivity in a way that is conceptually close to the renormalization of the generator carried out here. The main differences are~(i) the RG argument of~\cite{AV} is non-perturbative in the sense that there is a parameter~$\beta$, which roughly corresponds to our~$\cgamma$, but is allowed to be arbitrarily close to~$\nf43$, which is the Kolmogorov-Onsager threshold; and~(ii) the drift field~$\f(t,x)$ is a carefully engineered ``synthetic turbulence'' built from periodic ingredients, not a stationary random field. We also mention a subsequent paper~\cite{BSW} which modified the construction from~\cite{AV}, using convex integration techniques, to produce an~$\f(t,x)$ which is also a weak solution of the incompressible Euler equation.

\smallskip

The anomalous regularization statement of Theorem~\ref{t.regularity} is related to a distinct circle of ideas originating in the Obukhov-Corrsin theory of passive scalar turbulence~\cite{Obukhov,Corrsin}. That theory predicts a critical relationship between the spatial regularity of a turbulent velocity field and the regularity of passively advected scalars, with the implication that scalars should exhibit H\"older regularity that persists \emph{uniformly as the molecular diffusivity vanishes}. These predictions have been studied most intensively in the Kraichnan model~\cite{Kraichnan68,Kraichnan94,BGK,FGV}, which assumes a velocity field that is Gaussian in space and white-noise in time. The white-in-time structure leads to closed stochastic equations for correlation functions, enabling exact computations of anomalous scaling exponents and a detailed understanding of the Lagrangian mechanisms underlying anomalous dissipation~\cite{GK95,CFKL95,DE}. Very recently, sharp anomalous regularization bounds in this Kraichnan setting were proved in~\cite{GGM24,DGP25,Rowan.OC}, obtaining the predicted Obukhov-Corrsin regularity up to logarithmic corrections. In a different direction,~\cite{CCS,HCR} produced deterministic velocity field constructions exhibiting anomalous dissipation while maintaining uniform-in-diffusivity H\"older regularity of the passive scalar.

\smallskip

The present work establishes anomalous regularization in a qualitatively different setting: the velocity field is \emph{frozen} (time-independent) and genuinely random, rather than white-in-time or deterministically constructed. Our H\"older exponent~$1-C\cgamma^{\nf12}$ depends explicitly on the Hurst parameter~$\cgamma$ and is \emph{quenched}, meaning it holds for almost every realization of the random environment. As in the Obukhov-Corrsin prediction, the bound is uniform in the molecular diffusivity~$\nu$. The mechanism underlying this regularization is the \emph{control} of the homogenization error across scales provided by the renormalization group flow: at most scales that error is of order~$\cgamma^{\nf12}$, which means the heterogeneous operator is well-approximated by a constant-coefficient Laplacian. At such scales, harmonic approximation and a Campanato iteration yield a one-step improvement in regularity, an idea that also underlies the large-scale regularity theory in quantitative stochastic homogenization (see~\cite{AKMBook,GNO}). The scales where this approximation fails are rare, and crucially, the \emph{scale-independence} of the approximation error, obtained from the RG analysis, prevents them from accumulating. This is what allows the iteration to propagate regularity across the full range of scales.

\smallskip

In our previous paper~\cite{ABK.SD}, we proved a weaker version of Theorem~\ref{t.regularity} in the borderline case~$\cgamma=0$. There we obtained a H\"older estimate at a single point, but without scale-independence, we could not obtain sufficient stochastic integrability to deduce an estimate of the global H\"older norm. In effect, this estimate yields oscillation decay only down to a large mesoscopic scale (see~\cite[Theorem C]{ABK.SD} and the discussion following it).

\smallskip

A key methodological feature of the present work is that our arguments do not rely on Gaussian structure or moment closure. Essentially all prior rigorous results on superdiffusion and anomalous transport for stationary random environments---including the Kraichnan literature and the works~\cite{TV,CHT,CMOW,MOW,KO} discussed above---exploit either the white-in-time structure of the velocity field (which yields closed equations for correlation functions) or explicit Gaussian covariance structure. Our approach is different: the coarse-graining theory developed in Section~\ref{s.coarse.graining.theory} is entirely deterministic and applies to arbitrary coefficient fields, while the probabilistic estimates in Section~\ref{s.RG.flow} use only finite-range dependence and moment bounds on the scale decomposition~$\{\mathbf{j}_k\}$. No Gaussian calculus appears anywhere in this paper. As a consequence, our methods are robust in ways that closure-based or exact-formula arguments are not.

\smallskip

More broadly, our proof implements a genuinely multiscale renormalization procedure in the spirit of Wilson~\cite{Wilson71,WilsonKogut74}: one eliminates degrees of freedom scale by scale and tracks the induced flow of a small collection of coarse observables. This viewpoint is classical in rigorous renormalization group approaches to constructive field theory and statistical mechanics; see, for instance, the foundational works~\cite{GK85,GN1,GN2}, the general finite-range multiscale framework developed in~\cite{BBS} and modern implementations in correlated lattice interface models such as the discrete Gaussian model~\cite{BPR1,BPR2}. We also mention the non-perturbative real-space renormalization program originating in Balaban's work (see Dimock's expository account~\cite{DimockBalaban} and references therein).  While these works concern Gibbs measures and correlation functions rather than quenched operators, they share the structural features that are most relevant here: an explicit renormalization map, a quantitative notion of closure, and the control of errors along a long multiscale flow.

\smallskip

From the point of view of infrared renormalization, the dichotomy between the borderline case~$\cgamma=0$ and the supercritical regime~$\cgamma>0$ is reminiscent of the~$\phi^4$ model at and below its upper critical dimension. In~$\phi^4$ theory,~$d=4$ is the marginal dimension: the coupling is (naively) dimensionless and the RG flow produces logarithmic corrections to effective parameters, whereas in dimensions~$d=4-\varepsilon$ the same interaction becomes weakly relevant and leads to power-law scaling governed by a small~$\varepsilon$-dependent fixed point; see for example~\cite{WilsonKogut74,GK85,GN1,GN2,ZinnJustin2002,Cardy96}. In our setting,~$\cgamma=0$ is the marginal, critically-correlated case (with logarithmic-type superdiffusivity), while~$\cgamma>0$ may be viewed as a ``slightly subcritical'' perturbation in which the scale-by-scale injection of infrared layers becomes weakly relevant and drives a power-law growth of the running effective diffusivity. This analogy is meant at the level of RG relevance/marginality rather than a literal identification with a Euclidean field theory.

\smallskip

Renormalization ideas also arise in transport and PDE, and---as in the present work---the natural renormalized quantity is a (running) effective diffusivity rather than a coupling in an action. A prototypical example is the ``exact renormalization'' approach to turbulent transport developed by~\cite{AMajda}, in which one iteratively renormalizes eddy diffusivity across scales (see also~\cite{MajdaKramer99}). The first rigorous proof of power-law superdiffusion using related perpetual homogenization techniques for a multiscale shear flow was obtained in~\cite{BAO}. In the probabilistic literature on quenched random media, related multiscale renormalization schemes have been implemented to prove diffusive scaling limits in perturbative regimes: Bricmont and Kupiainen developed a space--time RG method for (weakly disordered) nonreversible random walks and closely related Markov processes~\cite{BK}; see also the related work of~\cite{SZ}. Related Wilsonian RG viewpoints have also been developed for equations in other settings, for example in~\cite{Kupiainen}. The present paper fits into this landscape but with features specific to scale-invariant random environments: the iteration acts on the generator in a fixed realization, and the coarse-grained diffusivities serve as quenched block variables. The observed flow is \emph{multifractal}: fluctuations do not average out with scale, but remain scale-local and of relative size~$O(\cgamma^{\nf12})$. The outcome is not convergence to a single homogenized operator, but tight control of a running, scale-dependent diffusivity---sufficient to establish both superdiffusivity and anomalous regularity.

\subsection{Formalization of Theorems~\ref{t.superdiffusivity}, \ref{t.homogenization} and~\ref{t.regularity} in the Lean 4 proof assistant}
\label{ss.formalization}

Theorems~\ref{t.superdiffusivity}, \ref{t.homogenization} and~\ref{t.regularity} have been formalized and machine-verified in the Lean~4 proof assistant, in the form in which they are stated above. The only mathematical difference between the formal statements and those presented here are that the former requires a slightly stronger (qualitative) assumption of~$C^2$ regularity for the shells~$\mathbf{j}_k$, in place of~$C^{1,1}$. The formalization is built on the mathlib library, on a formalization of the quantitative homogenization theory of~\cite{AK.HC}, and on a library for Markov processes written for this purpose, over which the diffusion process of Theorem~\ref{t.superdiffusivity} is constructed and characterized through the resolvent of its generator. The development is publicly available at
\begin{center}
\href{https://github.com/scottnarmstrong/Superdiffusion}{\texttt{github.com/scottnarmstrong/Superdiffusion}}\,.
\end{center}
The formal proofs closely follow the proofs presented here, and every intermediate result of Sections~\ref{s.coarse.graining.theory}--\ref{s.theprocess} that they depend on is formalized as well. The three formal theorems reduce to the standard foundational axioms of mathlib; the development contains no unproved placeholders and no custom axioms. A reader who does not wish to read the Lean development can still check what has been proved: each of the three theorems is restated in a self-contained file that uses only mathlib vocabulary, and the fact that this restatement follows from the formal theorem is verified by the Lean kernel. 

\subsection{Notation}
\label{ss.notation}

\subsubsection*{Basic notation}
The Euclidean norm on~$\R^m$ is denoted by~$|\cdot|$. We write~$r\wedge s \coloneqq \min\{ r,s\}$ and~$r\vee s \coloneqq \max\{ r,s\}$. If~$r\in \R$, then~$r_+ \coloneqq r \vee 0$. The H\"older conjugate exponent of~$p\in[1,\infty]$ is denoted by~$p' \coloneqq p(p-1)^{-1}$ for~$p\in(1, \infty)$ and~$p'=1$ (resp.,~$p'=\infty)$ if~$p=\infty$ (resp.,~$p=1$). We let~$\linear_e(x) = e \cdot x$ denote the linear function with slope~$e\in\Rd$. The distance between subsets~$A,B\subseteq\Rd$ is~$\dist(A,B) \coloneqq \inf_{x\in A, y\in B} |x-y|$. Indicator functions---both for events and for subsets of~$\Rd$---are denoted by~$\indc$.

\subsubsection*{Linear algebra}
The set of~$m$-by-$n$ real matrices is~$\R^{m\times n}$, with transpose~$B^t$ for~$B\in \R^{m\times n}$. The~$n$-by-$n$ identity matrix is~$\mathrm{I}_n$, and the~$2d$-by-$2d$ identity matrix is~$\Itwod$. The symmetric and antisymmetric~$n$-by-$n$ matrices are~$\Rsym$ and~$\Rskew$, respectively; the positive definite symmetric matrices are~$\Rsymp$. The Loewner ordering on~$\R^{n\times n}_{\mathrm{sym}}$ is denoted by~$\leq$: for~$A,B\in \R^{n\times n}_{\mathrm{sym}}$,~$A\leq B$ means~$B-A$ has nonnegative eigenvalues. Unless otherwise indicated, the norm~$|A|$ for~$A\in\R^{m\times n}$ is the operator norm (square root of the largest eigenvalue of~$A^tA$). If~$\mathbf{m}$ is a matrix-valued field on~$U\subseteq\Rd$ with entries~$(\mathbf{m}_{ij})_{i,j=1}^d$, the divergence~$\nabla \cdot \mathbf{m}$ is the vector field with components~$(\nabla \cdot \mathbf{m})_j = \sum_{i=1}^d \partial_{x_i} \mathbf{m}_{ij}$.

\subsubsection*{Measure and integration}
The Lebesgue measure of a measurable subset~$U \subseteq\Rd$ is denoted by~$|U|$; for a codimension-$1$ subset~$V$ (such as a boundary~$\partial U$),~$|V|$ denotes the~$(d-1)$-dimensional Hausdorff measure. The cardinality of a finite set~$A$ is also denoted~$|A|$. Volume-normalized integrals and~$L^p$ norms are
\begin{equation*}
	(f)_U \coloneqq \fint_U f(x) \,dx \coloneqq \frac{1}{|U|} \int_U f(x)\,dx
	\qqand
	\| f \|_{\underline{L}^p(U)} \coloneqq \Bigl( \fint_U |f(x)|^p \,dx \Bigr)^{\nf1p}\,.
\end{equation*}
Similarly, for a finite set~$A$ and~$f:A \to \R$, we write
\begin{equation*}
	\avsum_{a \in A} f(a) \coloneqq \frac{1}{|A|} \sum_{ a\in A} f(a)
	\,.
\end{equation*}
For a nonnegative random variable,~$\E[X]$ always denotes its possibly infinite Lebesgue integral. In particular, moment hypotheses involving~$\E[X^p]$ include the corresponding integrability assertion whenever the displayed upper bound is finite.

\subsubsection*{Function spaces}

We use the standard H\"older spaces~$C^{k,\alpha}(U)$ for~$k\in\N_0$ and~$\alpha \in (0,1]$, as well as the Sobolev spaces~$W^{s,p}(U)$ for~$s\in \R$ and~$p\in[1,\infty]$. When~$p = 2$, we write~$H^s \coloneqq W^{s, 2}$. The fractional Sobolev spaces are defined in Appendix~\ref{s.functional.inequalities}. The space~$W^{1,p}_0(U)$ denotes the closure of~$C^\infty_c(U)$ in~$W^{1,p}(U)$. For a domain~$U$,~$X_{\mathrm{loc}}(U)$ denotes the set of functions belonging to~$X(U \cap B_R)$ for every~$R\in [1,\infty)$. We let~$C_0(\Rd)$ denote the space of continuous functions vanishing at infinity, and~$C_c^k(\Rd)$ denotes~$C^k$ functions with compact support. The classical Sobolev norm is
\begin{equation*}
	\| f \|_{{W}^{1,p}(U)}
	\coloneqq
	\Bigl( \| \nabla f \|_{{L}^p(U)}^p+ \| f \|_{{L}^p(U)}^p \Bigr)^{\frac1p} \,,
\end{equation*}
and the volume-normalized variant (for~$|U|<\infty$) is
\begin{equation*}
	\| f \|_{\underline{W}^{1,p}(U)}
	\coloneqq
	\Bigl( \| \nabla f \|_{\underline{L}^p(U)}^p+ |U|^{-\frac pd} \| f \|_{\underline{L}^p(U)}^p \Bigr)^{\frac1p} \,.
\end{equation*}
For~$p=\infty$, the corresponding formulas are
\begin{equation*}
	\|f\|_{W^{1,\infty}(U)}\coloneqq\|\nabla f\|_{L^\infty(U)}+\|f\|_{L^\infty(U)}
	\qqand
	\|f\|_{\underline W^{1,\infty}(U)}\coloneqq\|\nabla f\|_{L^\infty(U)}+|U|^{-1/d}\|f\|_{L^\infty(U)}\,.
\end{equation*}
For~$k\in\N$, the normalized norms of higher order are
\begin{equation*}
	\| f \|_{\underline{W}^{k,\infty}(U)}
	\coloneqq
	\sum_{j=0}^{k} |U|^{-\frac{k-j}{d}} \| \nabla^j f \|_{L^\infty(U)}\,;
\end{equation*}
for~$k=1$ this is the previous display.
The seminorm is~$[ f ]_{\underline{W}^{1,p}(U)} \coloneqq  \| \nabla f \|_{\underline{L}^p(U)}$.
The negative, dual seminorms are defined by
\begin{equation*}
	\bigl[ f \bigr]_{\underline{W}^{-1,p'}(U)} \coloneqq
	\sup
	\biggl\{
	\fint_U f(x) g(x) \,dx \,:\,
	g\in C^\infty_c(U),\ [ g ]_{\underline{W}^{1,p}(U)} \leq 1 \biggr\}
	\,,
\end{equation*}
and, with zero-average test functions,
\begin{equation*}
	\bigl[ f \bigr]_{\Wminusul{-1}{p'}(U)} \coloneqq
	\sup
	\biggl\{
	\fint_U f(x) g(x) \,dx \,:\,
	g\in C^\infty(U),\ [g]_{\underline{W}^{1,p}(U)} \leq 1\,, \ (g)_U=0 \biggr\}\,.
\end{equation*}
When~$p=p'=2$ we write~$\underline{H}^{-1}$ and~$\Hminusul$ in place of~$\underline{W}^{-1,p}$ and~$\Wminusul{-1}{p}(U)$, respectively. For~$s\in(0,1)$ and~$p\in(1,\infty)$, the full and mean-zero negative dual norms are
\begin{align*}
	\|f\|_{\underline{W}^{-s,p}(U)}
	&\coloneqq
	\sup\Bigl\{\bigl|\fint_U f\cdot g\bigr|:
	g\in C^\infty(U),\ \|g\|_{\underline{W}^{s,p'}(U)}\leq1\Bigr\}\,,
	\\
	\|f\|_{\Wminusul{-s}{p}(U)}
	&\coloneqq
	\sup\Bigl\{\bigl|\fint_U f\cdot g\bigr|:
	g\in C^\infty(U),\ \|g\|_{\underline{W}^{s,p'}(U)}\leq1,\ (g)_U=0\Bigr\}\,.
\end{align*}
Here the full volume-normalized fractional norm is
\begin{equation*}
	\|g\|_{\underline{W}^{s,p}(U)}
	\coloneqq
	[g]_{\underline{W}^{s,p}(U)} + |U|^{-\frac sd}\|g\|_{\underline{L}^p(U)}
	\,, \qquad
	\|g\|_{\underline{W}^{s,\infty}(U)}
	\coloneqq
	|U|^{-\frac sd}\|g\|_{L^\infty(U)} + [g]_{C^{0,s}(U)}
	\,;
\end{equation*}
the definitions are componentwise, with~$[F]_{\underline H^1(U)}^2=[F_1]_{\underline H^1(U)}^2+[F_2]_{\underline H^1(U)}^2$ for doubled fields.
The full dual is used for flux pairings; the hatted dual is used when the test field has zero mean.
The H\"older seminorm and norm for~$\alpha\in(0,1]$ are defined by
\begin{equation*}
	[ f ]_{C^{0,\alpha}(U)}
	\coloneqq
	\sup_{x,y\in U,\, x\neq y} \frac{|f(x) - f(y)|}{|x-y|^\alpha}
	\qqand
	\| f \|_{C^{0,\alpha}(U)}
	\coloneqq
	\| f \|_{L^\infty(U)} + [ f ]_{C^{0,\alpha}(U)}\,,
\end{equation*}
with higher-order versions~$C^{k,\alpha}(U)$ defined by requiring~$\partial^\beta f \in C^{0,\alpha}(U)$ for all~$|\beta| \leq k$. The volume-normalized H\"older norm is
\begin{equation*}
	\| f \|_{\underline{C}^{0,\alpha}(U)}
	\coloneqq
	|U|^{-\frac{\alpha}{d}} \| f \|_{L^\infty(U)} + [ f ]_{C^{0,\alpha}(U)}\,.
\end{equation*}
For~$s\in(0,1)$ and~$U$ a cube, the H\"older seminorm~$[\cdot]_{C^{0,s}(U)}$ is equivalent (up to dimensional constants) to the volume-normalized fractional Sobolev norm~$[\cdot]_{\underline{W}^{s,\infty}(U)}$; see~\eqref{e.fractional.Sobolev.seminorm}. We often use these interchangeably.
For~$s<\nf12$,~$C^{0,\nf12}(\cu_m)\subset\underline H^s(\cu_m)$ with norm factor~$C(1-3^{2s-1})^{-\nf12}$; the endpoint is excluded.

\smallskip

The spaces of potential and solenoidal vector fields vanishing on~$\partial U$ are
\begin{equation*}
	\Lpoto(U) \coloneqq \left\{ \nabla u \,:\, u \in H^1_0(U) \right\}\,,
	\quad
	\Lsolo(U) \coloneqq \Bigl\{ \g \in L^2(U;\Rd) \,:\, \forall \phi\in H^1(U) \,, \ \int_U \g\cdot \nabla \phi = 0 \Bigr\}\,.
\end{equation*}
On~$\Omega$,~$\Lpot$ and~$\Lsol$ denote the stationary potential and solenoidal subspaces of~$L^2(\Omega;\Rd)$.

\subsubsection*{Besov-type seminorms}

The (volume-normalized) Besov-type seminorm on~$\cu_m$ for exponents~$s \in (0,1]$ and~$p,q\in [1,\infty)$ is
\begin{equation}
	\label{e.Besov.seminorm}
	[ f ]_{\underline{B}^s_{p,q}(\cu_m)}
	\coloneqq
	\biggl(
	\sum_{n=-\infty}^m
	3^{-nsq}
	\biggl(
	\avsum_{z\in 3^{n-1} \Zd \cap \cu_m, \, z+\cu_n \subseteq \cu_m}
	\bigl\| f - (f)_{z+\cu_n} \bigr\|_{\underline{L}^p(z+\cu_n)}^p
	\biggr)^{\!\nf qp} \biggr)^{\!\nf1q} \,,
\end{equation}
with corresponding norm~$\| f \|_{\underline{B}^s_{p,q}(\cu_m)} \coloneqq [ f ]_{\underline{B}^s_{p,q}(\cu_m)} + 3^{-sm} |(f)_{\cu_m}|$. For~$q=p$, this seminorm is equivalent to the fractional Sobolev seminorm: there exists~$C(d)<\infty$ such that, for every~$s\in (0,1)$ and~$p\in [1,\infty)$,
\begin{equation}
	\label{e.Wsp.vs.Bspp.intro.A}
	C^{-1} [ f ]_{\underline{W}^{s,p}(\cu_m)}
	\leq
	[ f ]_{\underline{B}^s_{p,p}(\cu_m)}
	\leq
	C [ f ]_{\underline{W}^{s,p}(\cu_m)}
	\,.
\end{equation}

\subsubsection*{Triadic cubes and simplices}

The axis-aligned triadic cube with side length~$3^m$ centered at the origin is
\begin{equation*}
	\cu_m \coloneqq \Bigl( -\frac12 3^m, \frac 12 3^m \Bigr)^{\!d}\,.
\end{equation*}
More generally,~$z + \cu_m$ denotes the translate of~$\cu_m$ centered at~$z \in \Rd$.

We let~$\mathcal{P}$ denote the set of permutations of~$\{ 1,\ldots, d\}$. For any permutation~$\pi \in \mathcal{P}$,~$n\in\Z$, and~$z \in\Rd$, the simplex of size~$3^n$ is
\begin{equation}
	\label{e.simplex.def}
	\triangle_n^\pi(z) \coloneqq z + 3^n \Bigl\{ (x_1,\ldots,x_d) \in\Rd \,:\, -\frac12 < x_{\pi(1)} < x_{\pi(2)} < \cdots < x_{\pi(d)} < \frac12 \Bigr\}
	\,.
\end{equation}
We write~$\triangle_n = \triangle_n^\pi$ when~$\pi$ is the identity permutation. A \emph{vertex} of a simplex is an extreme point; each simplex~$\triangle$ has~$d+1$ vertices, and we denote this set by~$V(\triangle)$. For a collection~$S$ of simplices,~$V(S) \coloneqq \bigcup_{\triangle\in S} V(\triangle)$.

Since~$|\mathcal{P}| = d!$, we have~$| \triangle_n^\pi| = 3^{dn}/d!$. Any triadic simplex~$\triangle_m^\pi(z)$ is the disjoint union (up to measure zero) of~$3^{d(m-n)}$ simplices of the form~$\triangle_n^\sigma(z')$ with~$z'\in \triangle_m^\pi(z)$. The cube~$z+\cu_n$ is the disjoint union of~$\{ \triangle_n^\pi(z) \,:\, \pi\in \mathcal{P} \}$. The partition of~$z+\cu_m$ into triadic simplices of size~$3^n$ is denoted
\begin{equation*}
	S_{n}(z{+}\cu_m) \coloneqq
	\bigl\{ z+\triangle_n^\pi(z') \,:\, z'\in 3^n\Zd \,, \, \triangle_n^\pi(z') \subseteq \cu_m \bigr\}\,,
	\quad n\leq m\,.
\end{equation*}

\section{Deterministic coarse-graining}
\label{s.coarse.graining.theory}

In this section, we develop a \emph{deterministic} coarse-graining theory for general elliptic coefficient fields
\begin{equation*}
	\a : \Rd \to \mathbb{R}^{d\times d}
	\,.
\end{equation*}
We split the field into its symmetric part~$\s\coloneqq \frac12(\a+\a^t)$ and antisymmetric part~$\k \coloneqq \frac12(\a-\a^t)$. For convenience, we assume that~$\a$ is \emph{qualitatively} uniformly elliptic in every bounded subset of~$\Rd$. This means that~$\s$ is a positive matrix and
\begin{equation}
	\| \s + \k^t \s^{-1} \k \|_{L^\infty(B_R)}
	+
	\| \s^{-1} \|_{L^\infty(B_R)} < \infty
	\,, \qquad \forall R<\infty \,.
	\label{e.qualitative.ellipticity}
\end{equation}
Throughout the section, the field~$\a$ is a fixed coefficient field; in particular, no randomness is used here. Probabilistic arguments will be important in Sections~\ref{s.RG.flow} and~\ref{s.regularity}, when the deterministic estimates developed here are applied to the random coefficient field~$\a(x)=\nu \Id  + \k(x)$ defined in~\eqref{e.generator.intro}, with~$\k$ being the random stream matrix in Section~\ref{ss.assumptions}.

\smallskip

Our starting point is the coarse-graining framework of~\cite{AK.HC}. Section~\ref{ss.cg.mats} recalls the essential definitions and basic properties of the coarse-grained matrices and ellipticity constants. Section~\ref{ss.cg.ellipticity} states the basic elliptic estimates and functional inequalities from this theory, such as the coarse-grained Caccioppoli and Poincar\'e inequalities.

\smallskip

The proof of our main results will require coarse-graining tools that go beyond the theory developed in~\cite{AK.HC}. For equations with symmetric coefficients, this is relatively straightforward---one can work directly at the level of the scalar energy. In the nonsymmetric case, we need to work with the~$2d\times 2d$ block system and to control both gradients and fluxes simultaneously. This is developed in Section~\ref{ss.cg.RHS}, which is of independent interest.

\smallskip

Section~\ref{ss.sensitivity} develops the \emph{coarse-grained sensitivity} theory: quantitative estimates showing how the coarse-grained matrices and ellipticity constants change under antisymmetric perturbations of the field. We say that these bounds are ``coarse-grained'' because all estimates must measure the perturbations against the \emph{coarse-grained ellipticity at the current scale} rather than the constants of uniform ellipticity.  These sensitivity estimates are the analytic engine that makes a genuine renormalization group iteration possible in Section~\ref{s.RG.flow}: when applied to the setting of our main results, they allow us to understand the effect of individual scales of the stream matrix without ever ``reaching down'' to the microscopic scale to exploit the bare microscopic diffusivity~$\nu>0$.

\smallskip

The renormalization group argument in Section~\ref{s.RG.flow} propagates estimates on the coarse-grained matrices and ellipticity constants across scales. Section~\ref{ss.blackboxes} explains how to turn such information into quantitative homogenization and harmonic-approximation estimates for inhomogeneous equations: we treat these as ``black boxes'' that take as input the multiscale ellipticity bounds and output estimates on Dirichlet problems. Combined with the RG analysis, these lead eventually to the full statement of Theorem~\ref{t.homogenization}.

\subsection{Coarse-grained matrices}
\label{ss.cg.mats}

In this subsection we recall the definitions and basic properties of the coarse-grained matrices associated with a coefficient field~$\a$, following~\cite{AKMBook,AK.HC}.
For each bounded Lipschitz domain~$U\subseteq\Rd$, the coarse-grained matrices~$\s(U)$,~$\s_*(U)$ and~$\k(U)$ together represent the effective linear response of the operator at that scale. These matrices and their associated ellipticity constants play a central role in essentially all subsequent estimates in the paper.
The properties presented here are standard in the theory, can be found in~\cite{AK.HC}, and most have relatively simple proofs. We record them only to fix notation and make the later arguments self-contained.

\smallskip

\subsubsection{Definitions of the coarse-grained matrices}

We denote the linear space of solutions~$u$ of~$-\nabla \cdot \a \nabla u = 0$ in an open set~$U \subseteq \Rd$ by
\begin{equation}
	\mathcal{A}(U; \a)  \coloneqq
	\bigl\{ u \in H^1(U) : -\nabla \cdot \a \nabla u = 0 \quad \mbox{in $U$} \bigr\}
	\,.
	\label{e.local.solutions.general}
\end{equation}
The space of solutions of the adjoint equation is denoted by~$\mathcal{A}^*(U; \a)  \coloneqq \mathcal{A}(U; \a^t)$.
For every~$p,q \in \Rd$ and bounded Lipschitz domain~$U \subseteq \Rd$, we define the quantity
\begin{equation}
	J(U,p,q\,; \a)  \coloneqq  \sup_{v \in \mathcal{A}(U;\a)} \fint_{U} \Bigl( -\frac{1}{2} \nabla v \cdot \s \nabla v - p \cdot \a \nabla v + q \cdot \nabla v \Bigr)
	\,.
	\label{e.J.general}
\end{equation}
The supremum in the variational problem on the right side of~\eqref{e.J.general} is attained by a maximizer in~$H^1(U)$, unique up to additive constants. We denote it by
\begin{equation}
	v(\cdot,U,p,q \,; \a)
	\coloneqq
	\text{maximizer of~\eqref{e.J.general}}
	\,.
	\label{e.argmax.J}
\end{equation}
For the adjoint maximizer, we write
\begin{equation*}
	v^*(\cdot,U,p,q;\a)\coloneqq v(\cdot,U,p,q;\a^t)\,.
\end{equation*}
It is immediate that the mapping~$(p,q) \mapsto J(U,p,q;\a)$ is quadratic. We define the coarse-grained matrices~$\s(U; \a),\s_*(U; \a)$ and~$\k(U;\a)$ via the identity
\begin{equation}
	\label{e.J.coarse.grained}
	J(U,p,q\,; \a) = \frac12 p \cdot \s (U; \a) p + \frac12 (q + \k(U; \a) p) \cdot \s_*^{-1}(U; \a)(q + \k(U; \a) p) - p \cdot q
	\,.
\end{equation}
We also define
\begin{equation}
	\b(U; \a) \coloneqq  ( \s + \k^t\s_*^{-1}\k )(U; \a)
	\,.
	\label{e.b.U.def}
\end{equation}
We further define
\begin{equation}
	\a(U;\a) \coloneqq \s(U;\a) - \k^t(U;\a)
	\qquad \mbox{and} \qquad
	\a_*(U;\a) \coloneqq \s_*(U;\a) - \k^t(U;\a) \,.
	\label{e.a.astar.U.def}
\end{equation}
We will sometimes drop the dependence on the field~$\a$ from the notation if it is clear from context.

\subsubsection{Properties of the coarse-grained matrices}

We list here (without proof) some of the important properties of the coarse-grained matrices and the associated energy functional~$J$. Proofs can be found in~\cite{AK.Book,AK.HC}.
The following assertions are valid for every bounded Lipschitz domain~$U \subseteq \Rd$.

\begin{itemize}

\item Ellipticity bounds: the coarse-grained matrices are bounded by integrals of the field itself:
\begin{equation}
	\label{e.CG.bounds.1}
	\Bigl( \fint_U \s^{-1}(x)\, dx \Bigr)^{\!-1}
	\leq \s_*(U) \leq \s(U)
	\leq
	\b(U)
	\leq
	\fint_U (\s + \k^t\s^{-1} \k )(x)\, dx
\end{equation}

\item {Control of the symmetric part of~$\k(U)$.} We have
\begin{equation}
	\label{e.ksym.by.sss}
	(\k+\k^t)(U) \leq (\s-\s_*)(U)
	\quad \mbox{and} \quad
	-(\k+\k^t)(U) \leq (\s-\s_*)(U)\,.
\end{equation}

\item First variation: for every~$w\in \A(U)$,
\begin{equation}
	q\cdot \fint_U \nabla w - p \cdot \fint_U \a \nabla w
	=
	\fint_U \nabla w \cdot \s \nabla v(\cdot,U,p,q)
	\,.
	\label{e.firstvar}
\end{equation}

\item Second variation and quadratic response: for every~$w\in \A(U)$,
\begin{multline}
	J(U,p,q) - \fint_U \Bigl( -\frac12 \nabla w \cdot \s\nabla w -p\cdot \a\nabla w+ q\cdot \nabla w   \Bigr)
	\\
	=
	\fint_U \frac12 \bigl ( \nabla v(\cdot,U,p,q) - \nabla w \bigr )\cdot \s\bigl ( \nabla v(\cdot,U,p,q) - \nabla w \bigr )\,.
	\label{e.quadresp.nosymm}
\end{multline}

\item Characterization in terms of the energy of the maximizer: for every~$p,q\in\Rd$,
\begin{equation}
	\label{e.Jenergyv.nosymm}
	J(U,p,q) = \fint_U \frac12 \nabla v(\cdot,U,p,q) \cdot \s \nabla v(\cdot,U,p,q)
\end{equation}
and
\begin{equation}
	\label{e.J.by.lin}
	J(U,p,q) = \frac12 \Bigl( q\cdot \fint_U \nabla v(\cdot,U,p,q) - p \cdot \fint_U \a \nabla v(\cdot,U,p,q) \Bigr) \,\,,
\end{equation}

\item Characterization in terms of spatial averages of gradients \& fluxes: for every~$p,q\in\Rd$,
\begin{equation}
	\label{e.v.spatial.averages}
	\left\{
	\begin{aligned}
		& \fint_U
		\nabla v(\cdot,U,p,q)
		=
		- p + \s_{*}^{-1}(U) (q + \k(U) p)
		\\ &
		\fint_U
		\a \nabla v(\cdot,U,p,q)
		=
		(\Id -\k^t  \s_{*}^{-1})(U)  q - \b(U) p
		\,.
	\end{aligned}
	\right.\,,
\end{equation}

\item Subadditivity: for every~$m,n\in\Z$ with~$n\leq m$ and~$p,q\in\Rd$,
\begin{equation}
	\label{e.subaddJ.nosymm}
	J(\cu_m,p,q)
	\leq
	\avsum_{z\in 3^n\Zd \cap \cu_m}
	J(z+\cu_n,p,q)
	\,.
\end{equation}

\item Linear response inequality: for every~$p,q \in \Rd$ and~$w\in \A(U)$
\begin{equation}
	\label{e.fluxmaps}
	\biggl|
	\fint_U \bigl( p \cdot \a \nabla w - q \cdot \nabla w  \bigr)
	\biggr|
	\leq
	\left( \fint_U \nabla w \cdot \s \nabla w \right)^{\nf12}
	\bigl( 2J(U, p , q  )\bigr )^{\nf12}
	\,.
\end{equation}

\item {Coarse-graining inequality:} For every~$p\in\Rd$ and~$w\in \A(U)$,
\begin{equation}
	\label{e.coarse.graining.nosymm}
	\biggl| p\cdot \fint_U ( \a_*(U) - \a ) \nabla w \biggr|
	\leq
	2^{\nf12}
	\bigl( p\cdot( \s(U) - \s_*(U))p \bigr)^{\nf12} \biggl( \fint_U \nabla w \cdot \s \nabla w \biggr)^{\!\nf12}
	\,.
\end{equation}

\item Energy-averaging inequalities: for every~$w\in \A(U; \a)$,
\begin{equation}
	\label{e.energymaps.nonsymm}
	\frac12\left( \fint_U \nabla w \right) \cdot \s_*(U) \left( \fint_U \nabla w \right)
	\leq
	\fint_U \frac12 \nabla w \cdot \s\nabla w
\end{equation}
and
\begin{equation}
	\label{e.energymaps.nonsymm.dual}
	\frac12\left( \fint_U \a\nabla w \right) \cdot  \b^{-1}(U) \left( \fint_U \a\nabla w \right)
	\leq
	\fint_U \frac12 \nabla w \cdot \s\nabla w
	\,.
\end{equation}

\item {Redundancy of the adjoint quantity.} The coarse-grained matrices for the transpose~$\a^t$ of the field are given by
\begin{equation}
	\label{e.adjoint.redundant}
	\s(U;\a^t) = \s(U ; \a)
	\,, \quad
	\s_*(U;\a^t) = \s_*(U ; \a)
	\,, \quad
	\k(U;\a^t) = - \k(U ; \a)
	\,.
\end{equation}

\item Homogeneity in the underlying coefficient field: for every~$\lambda \in (0,\infty)$,
\begin{equation}
	\label{e.homogeneity.in.a}
	\s(U; \lambda \a ) = \lambda\s(U ; \a)
	\,, \quad
	\s_*(U;\lambda\a) = \lambda \s_*(U ; \a)
	\,, \quad
	\k(U;\lambda \a) = \lambda \k(U ; \a)
	\,.
\end{equation}
In terms of the quantity~$J$, this states that, for every~$\lambda\in (0,\infty)$ and~$p,q \in \Rd$,
\begin{equation}
	\label{e.homogeneous.J}
	J(U,p,q\,;\lambda \a) = J(U,\lambda^{\nf12} p, \lambda^{-\nf12}q \,;  \a) \,.
\end{equation}

\item We can write~$J(U,p,q)$ as
\begin{align}
	J(U,p,q)
	&
	=
	\frac 12p \cdot \bigl (\s(U)- \s_* (U) \bigr )p
	+
	\frac12 p \cdot \bigl( \k(U) + \k^t(U) \bigr) p
	\notag \\ &
	\qquad
	+ \frac 12 \bigl (q -(\s_*-\k)(U) p\bigr ) \cdot \s_*^{-1} (U) \bigl (q-(\s_*-\k)(U) p\bigr )
	\,.
	\label{e.J.magic.squares}
\end{align}
\end{itemize}

The formula~\eqref{e.J.magic.squares} is just an algebraic rearrangement of~\eqref{e.J.coarse.grained}, but it is a very useful way to think about~$J$. In particular,~$J$ is small if~$(\s-\s_*)(U)$ is small and~$q$ is close to~$(\s_*-\k)(U)p$. The estimate~\eqref{e.shaking.lambda} is a simple consequence of~\eqref{e.ksym.by.sss},~\eqref{e.homogeneity.in.a},~\eqref{e.homogeneous.J} and~\eqref{e.J.magic.squares}.

\smallskip

Most of the properties listed above have relatively short proofs of only a few lines. The main exceptions are~\eqref{e.ksym.by.sss} and~\eqref{e.adjoint.redundant}. These require the~$2d\times 2d$ block matrix characterization of the coarse-grained matrices and the corresponding variational principles, which we explain in the next subsection.

\smallskip

We also require one additional property of~$J$, which is not proved in~\cite{AK.Book,AK.HC}. It is:
\begin{itemize}
\item For every~$p,q\in\Rd$ and~$\lambda \in (0,\infty)$,
\begin{equation}
	J(U, \lambda^{-\nf12} p,\lambda^{\nf12} q; \a)
	\leq
	\lambda^{-1} \bigl(
	\bigl( J(U, p, q; \a)  \bigr)^{\nf 12}
	+
	2^{-\nf12}
	|\lambda -1|
	\bigl( q \cdot \s_*^{-1} (U;\a) q\bigr)^{\nf 12}
	\bigr)^2
	\,.
	\label{e.shaking.lambda}
\end{equation}
\end{itemize}

\begin{proof}[{Proof of~\eqref{e.shaking.lambda}}]
We combine~\eqref{e.homogeneity.in.a},~\eqref{e.homogeneous.J} and~\eqref{e.J.magic.squares} to obtain
\begin{align*}
	J(U, \lambda^{-\nf12} p,\lambda^{\nf12} q; \a)
	&
	=
	J(U, p,q; \lambda^{-1} \a)
	\notag \\ &
	=
	\frac 12p \cdot \bigl (\s(U;  \lambda^{-1} \a)- \s_* (U;  \lambda^{-1} \a) \bigr )p
	+
	\frac12 p \cdot \bigl( \k(U;  \lambda^{-1} \a) + \k^t(U;  \lambda^{-1} \a) \bigr) p
	\notag \\ &
	\qquad
	+ \frac 12 \bigl (q -(\s_*-\k)(U;  \lambda^{-1} \a) p\bigr ) \cdot \s_*^{-1} (U;  \lambda^{-1} \a) \bigl (q-(\s_*-\k)(U;  \lambda^{-1} \a) p\bigr )
	\notag \\ &
	=
	\lambda^{-1}
	\biggl(
	\frac 12p \cdot \bigl (\s(U;\a)- \s_* (U;\a) \bigr )p
	+
	\frac12 p \cdot \bigl( \k(U;\a) + \k^t(U;\a) \bigr) p
	\notag \\ &
	\qquad\qquad
	+ \frac 12 \bigl (\lambda q -(\s_*-\k)(U;\a) p\bigr ) \cdot \s_*^{-1} (U;\a) \bigl (\lambda q-(\s_*-\k)(U;\a) p\bigr )
	\biggr)
	\,.
\end{align*}
By Young's inequality, we have, for every~$\ep >0$,
\begin{align*}
	\lefteqn{
	\frac12 \bigl(\lambda q -(\s_*-\k)(U;\a) p\bigr) \cdot \s_*^{-1} (U;\a) \bigl(\lambda q-(\s_*-\k)(U;\a) p\bigr)
	} \quad &
	\notag \\ &
	\leq
	\frac{1+\ep}{2}
	\bigl(q -(\s_*-\k)(U;\a) p\bigr) \cdot \s_*^{-1} (U;\a) \bigl(q-(\s_*-\k)(U;\a) p\bigr)
	+
	\frac{1+\ep}{2\ep} (\lambda -1)^2
	q \cdot \s_*^{-1} (U;\a) q
	\,.
\end{align*}
Using~\eqref{e.ksym.by.sss} and combining the above displays, we obtain, for every~$\ep >0$,
\begin{equation*}
	J(U, \lambda^{-\nf12} p,\lambda^{\nf12} q; \a)
	\leq \frac{1+\ep}{\lambda}
	J(U, p, q; \a)
	+
	\frac{1+\ep}{2\ep \lambda}  (\lambda -1)^2
	q \cdot \s_*^{-1} (U;\a) q
	\,.
\end{equation*}
Optimizing in the parameter~$\ep$ yields~\eqref{e.shaking.lambda}.
\end{proof}

\subsubsection{Block matrix characterization}

There is a second, equivalent way to define the coarse-grained matrices using a~$2d\times2d$ block matrix. This formulation exposes the duality and the Dirichlet and Neumann variational structures.

\smallskip

Proofs of the assertions in this subsection appear in~\cite{AK.Book,AK.HC}; see also~\cite[Chapter 10]{AKMBook} for the variational formulation of nonsymmetric divergence-form equations.

\smallskip

Throughout the rest of the paper, we associate, to each coefficient field~$\a$, the~$2d\times 2d$ block matrix-valued field
\begin{equation}
	\label{e.bfA.def}
	\bfA(x) \coloneqq
	\begin{pmatrix}
		( \s + \k^t\s^{-1}\k )(x)
		& -(\k^t\s^{-1})(x)
		\\ - ( \s^{-1}\k )(x)
		& \s^{-1}(x)
	\end{pmatrix}
	\,.
\end{equation}
It follows from our qualitative ellipticity assumption~\eqref{e.qualitative.ellipticity} that~$\bfA(x)$ is a positive matrix. Moreover, it admits the factorization
\begin{equation}
	\label{e.bfA.LDLt}
	\bfA
	=
	\begin{pmatrix}
		\Id & -\k^t \\
		0 & \Id
	\end{pmatrix}
	\begin{pmatrix}
		\s & 0 \\
		0 & \s^{-1}
	\end{pmatrix}
	\begin{pmatrix}
		\Id & 0 \\
		-\k & \Id
	\end{pmatrix}
	=
	\bfG_{-\k}^t \begin{pmatrix}
		\s & 0 \\
		0 & \s^{-1}
	\end{pmatrix}
	\bfG_{-\k}
	\,,
\end{equation}
where we define, for each matrix~$\h \in \R^{d \times d}$
\begin{equation}
	\label{e.Gh.def}
	\bfG_{\h}  \coloneqq
	\begin{pmatrix}
		\Id & 0 \\
		\h & \Id
	\end{pmatrix}
	\,.
\end{equation}
Note that~\eqref{e.bfA.LDLt} implies~$\det \bfA = 1$.
Observe that the~$\mathbf{G}_\h$'s satisfy
\begin{equation}
	\label{e.Gh.additive}
	\bfG_{\h_1 + \h_2} =  \bfG_{\h_1} \bfG_{\h_2} \qquad \forall \h_1,\h_2 \in \R^{d \times d} \,.
\end{equation}
In particular,~$\mathbf{G}_\h^{-1} = \mathbf{G}_{-\h}$ and from this and~\eqref{e.bfA.LDLt} we find that the inverse of~$\bfA$ is given by
\begin{equation}
	\label{e.bfA.inverse}
	\bfA^{-1} =
	\bfG_{\k}
	\begin{pmatrix}
		\s^{-1} & 0 \\
		0 & \s
	\end{pmatrix}
	\bfG_{\k}^t
	=
	\begin{pmatrix}
		\s^{-1} & -\s^{-1} \k \\
		-\k^t \s^{-1} & \s + \k^t \s^{-1} \k
	\end{pmatrix}
	=
	\bfR \bfA \bfR\,,
\end{equation}
where we denote
\begin{equation}
	\bfR \coloneqq \begin{pmatrix}
		0 & \Id \\
		\Id & 0
	\end{pmatrix}
	\,.
	\label{e.R.def}
\end{equation}

For each bounded Lipschitz domain~$U \subseteq\Rd$, we define subspaces of~$L^2(U;\Rd)$ consisting of divergence-free and potential vector fields by:
\begin{equation}
	\left\{
	\begin{aligned}
		& L^2_{\pot}(U) \coloneqq \nabla H^1(U)
		= \bigl\{ \nabla u \,:\, u \in H^1(U) \bigr\}
		\,, \\
		& L^2_{\pot,0}(U) \coloneqq \nabla H^1_0(U)
		= \bigl\{ \nabla u \,:\, u \in H^1_0(U) \bigr\}
		\,, \\
		& L^2_{\sol}(U) \coloneqq
		L^2_{\pot,0}(U)^{\perp}
		= \Bigl\{ \g \in L^2(U;\Rd) \,:\, \forall u \in H^1_0(U) \,, \int_U \g \cdot \nabla u = 0 \Bigr\}
		\,, \\
		& L^2_{\sol,0}(U)
		\coloneqq
		L^2_{\pot}(U)^{\perp}
		= \Bigl\{ \g \in L^2(U;\Rd) \,:\, \forall u \in H^1(U) \,, \int_U \g \cdot \nabla u = 0 \Bigr\}
		\,.
	\end{aligned}
	\right.
	\label{e.sols.and.pots}
\end{equation}
We next define, for each~$P\in\R^{2d}$, the variational quantity
\begin{equation}
	\mu(U,P) \coloneqq
	\inf\left\{ \fint_{U} \frac12 X \cdot \bfA X \,:\, X \in P+L^2_{\pot,0}(U) \times L^2_{\sol,0}(U) \right\}
	\,.
	\label{e.variational.mu.U.P}
\end{equation}
The infimum is once again a minimum, and~$P \mapsto \mu(U,P)$ is evidently quadratic in~$P$. Therefore we can identify a matrix~$\bfA(U)$ such that
\begin{equation}
	\mu(U,P)
	=
	\frac12 P \cdot \bfA(U) P
	\,.
	\label{e.bigA.def}
\end{equation}
The matrix~$\bfA(U)$ defined in~\eqref{e.bigA.def} is linked to the matrices defined in~\eqref{e.J.coarse.grained} and~\eqref{e.b.U.def} by
\begin{equation}
	\bfA(U)
	=
	\begin{pmatrix}
		( \s + \k^t\s_*^{-1}\k )(U)
		& -(\k^t\s_*^{-1})(U)
		\\ - ( \s_*^{-1}\k )(U)
		& \s_*^{-1}(U)
	\end{pmatrix}
	=
	\begin{pmatrix}
		\b(U)
		& -(\k^t\s_*^{-1})(U)
		\\ - ( \s_*^{-1}\k )(U)
		& \s_*^{-1}(U)
	\end{pmatrix}
	\,.
	\label{e.bigA.to.little.A}
\end{equation}
The last three lines give a useful variational formulation of the coarse-grained matrices.

\smallskip

We next define a double-variable version of~$J$: for each~$P,Q\in\R^{2d}$, we set
\begin{equation}
	\label{e.bfJfull.var}
	\bfJ(U, P,Q )
	\coloneqq
	\sup_{X  \in \S(U)}
	\fint_U
	\biggl( -\frac12  X  \cdot \bfA X -  P \cdot \bfA X +Q\cdot X \biggr)
	\,,
\end{equation}
where~$\S(U)$ is the linear subspace of~$\Lpot(U) \times \Lsol(U)$ defined by
\begin{equation}
	\S(U)  \coloneqq
	\biggl\{
	X \in \Lpot(U) \times \Lsol(U) \,:\,
	\int_U Y \cdot \bfA X=0
	\; \;  \forall Y \in \Lpoto(U) \times \Lsolo(U)
	\biggr\}
	\,.
	\label{e.bfS}
\end{equation}
Then this quantity admits the formula
\begin{equation}\label{e.bfcoarsegrainedmatrices}
\bfJ(U,P,Q ) = \frac12 P \cdot \bfA(U) P + \frac12 Q \cdot \bfA_*^{-1}(U) Q - P \cdot Q \,\,,
\end{equation}
where~$\bfA_*^{-1}(U)$ is given by the formula
\begin{equation*}
	\bfA_*^{-1}(U) = \bfR\bfA(U)\bfR
	=
	\begin{pmatrix}
		\s_*^{-1}(U) &
		- ( \s_*^{-1}\k )(U)  \\
		-(\k^t\s_*^{-1})(U) &
		( \s + \k^t\s_*^{-1}\k )(U)
	\end{pmatrix}
	\,.
\end{equation*}
In particular,~$\bfJ(U, P,0 ) = \mu(U,P)$.

\subsubsection{More properties of the coarse-grained matrices}

We list here (without proof) some of the further properties of the coarse-grained matrices, the associated energy functional~$J$ and their double-variable versions. Proofs can be found in~\cite{AK.Book,AK.HC}; we only record the statements needed later.

The following assertions are valid for every bounded Lipschitz domain~$U \subseteq \Rd$.

\begin{itemize}

\item More general version of~\eqref{e.CG.bounds.1}:
\begin{equation}
	\label{e.CG.bounds.2}
	\Bigl( \fint_U\bfA^{-1}(x)\,dx \Bigr)^{\!-1}\leq \bfA_*(U)\leq \bfA(U)\leq\fint_U\bfA(x)\,dx\,.
\end{equation}

\item Formula for~$\bfJ$ in terms of~$J$ and its adjoint: for every~$p,q,p^*,q^* \in \Rd$
\begin{equation} \label{e.bfJ.general}
\bfJ\left(U, \begin{pmatrix} p \\ q \end{pmatrix}, \begin{pmatrix} q^* \\ p^* \end{pmatrix}; \a\right) = \frac{1}{2}J(U, p-p^*, q^*-q\,; \a) + \frac{1}{2}J(U, p^*+p, q^*+q\,; \a^t)\,.
\end{equation}

\item {Formula for the sum of~$J$ and its adjoint.} For every~$p,q,h\in\Rd$,
\begin{align}
	\label{e.JJstar1}
	J(U,p,q-h; \a) + J(U,p,q+h ; \a^t)
	&
	=
	p \cdot (\s -\s_*)(U)p
	\notag \\ & \qquad
	+
	\bigl (q - \s_*(U)p\bigr )\cdot \s_*^{-1}(U) \bigl (q - \s_*(U)p\bigr )
	\notag \\ & \qquad
	+
	\bigl ( h-\k(U)p \bigr )  \cdot \s_*^{-1}(U)\bigl ( h-\k(U) p\bigr )
	\,.
\end{align}

\item {Formula for~$J$ and its adjoint in terms of~$\bfA$.} For every~$p,q\in\Rd$,
\begin{equation}
	\label{e.J.by.means.of.bfA}
	J(U,p,q\,; \a) = \frac12\begin{pmatrix} p \\ - q \end{pmatrix} \cdot \bfA(U; \a) \begin{pmatrix} p \\ -q \end{pmatrix} - p\cdot q
	\,, \quad
	J(U,p,q\,; \a^t) = \frac12\begin{pmatrix} p \\ q \end{pmatrix} \cdot \bfA(U; \a) \begin{pmatrix} p \\ q \end{pmatrix} - p\cdot q\,.
\end{equation}

\item
Characterization of the space~$\mathcal{S}(U)$ defined in~\eqref{e.bfS}:
\begin{equation}
	\label{e.findSfull}
	\S(U)
	=
	\left\{ (\nabla v+\nabla v^*, \a\nabla v - \a^t\nabla v^*) \,:\, v\in \A(U), \ v^*\in \A^*(U)
	\right\}
	\,.
\end{equation}

\item Denote the maximizer of the variational problem in~\eqref{e.bfJfull.var} by~$S(\cdot,U,P,Q)$.
Then, for every~$p,p^*,q,q^*\in\Rd$,
\begin{equation}
	\label{e.maximizers.J.to.bfJ}
	S
	\biggl(\cdot,U, \begin{pmatrix} p  \\ q \end{pmatrix}, \begin{pmatrix} q^* \\ p^* \end{pmatrix} \biggr)
	=
	\frac12
	\begin{pmatrix} \nabla v (\cdot,U,p{-}p^*,q^*{-}q)
		+
		\nabla v^*\bigl(\cdot,U,p^*{+}p,q^*{+}q\bigr) \\
		\a \nabla v (\cdot,U,p{-}p^*,q^*{-}q)
		-
		\a^t \nabla v^*\bigl(\cdot,U,p^*{+}p,q^*{+}q\bigr)
	\end{pmatrix}
	\,.
\end{equation}

\item The minimizer of the variational problem~\eqref{e.variational.mu.U.P} is~$S(\cdot,U,-P,0)$.

\item First variation for~$\bfJ$: for every~$T \in \S(U)$ and~$P,Q \in \R^{2 d}$,
\begin{equation}
	\fint_{U} (Q \cdot T - P \cdot \bfA T ) = \fint_{U} T \cdot \bfA S(\cdot, U, P, Q) \,\,.
	\label{e.firstvar.bfJ}
\end{equation}

\item Subadditivity: for every~$m,n\in\Z$ with~$n<m$,
\begin{equation}
	\label{e.subadda.nosymm}
	\bfA(\cu_m)
	\leq
	\avsum_{z\in 3^n\Zd \cap \cu_m} \!\!\!
	\bfA(z+\cu_n)
	\quad \mbox{and} \quad
	\bfA_*^{-1} (\cu_m)
	\leq
	\avsum_{z\in 3^n\Zd \cap \cu_m} \!\!\!
	\bfA_*^{-1}(z+\cu_n)
	\,.
\end{equation}

\end{itemize}

Note that the formula~\eqref{e.bfcoarsegrainedmatrices} can also be rewritten as
\begin{equation}
	\label{e.bfJ.magic}
	\bfJ(U,P,Q)
	=
	\frac12 P \cdot (\bfA(U) - \bfA_*(U)) P  +   \frac12 (Q  - \bfA_*(U) P) \cdot \bfA_*^{-1}(U) (Q  - \bfA_*(U)  P)
	\,.
\end{equation}
It can also be written as
\begin{equation}
	\label{e.bfJ.magic2}
	\bfJ(U,P,Q)
	=
	\frac12 Q \cdot (\bfA_*^{-1}(U) - \bfA^{-1}(U)) Q  +   \frac12 (P  - \bfA^{-1}(U) Q) \cdot \bfA(U) (P - \bfA^{-1}(U) Q)
	\,.
\end{equation}
We will often need to estimate the normalized matrix~$\mathbf D^{-\nf12}\mathbf E\mathbf D^{-\nf12}$. Conjugating~$\mathbf D$ and~$\mathbf E$ by~$\G_\h$ defined in~\eqref{e.Gh.def} does not change its spectrum, as recorded in the following lemma.

\begin{lemma}
\label{l.conjugating.ratio.Gh}
Let~$\mathbf{D} \in \R^{2d \times 2d}_{\sym,+}$,~$\mathbf{E} \in \R^{2d \times 2d}_{\sym}$,
and~$\h \in \R^{d \times d}$.
Denote~$\mathbf{E}_{\h}  \coloneqq  \G^t_{\h} \mathbf{E} \G_{\h}$ and~$\mathbf{D}_{\h}  \coloneqq  \G^t_{\h} \mathbf{D} \G_{\h}$.
Then
\begin{equation*}
	|\mathbf{D}^{-\nf 12} \mathbf{E} \mathbf{D}^{-\nf 12} | =
	|\mathbf{D}_{\h}^{-\nf 12} \mathbf{E}_{\h} \mathbf{D}_{\h}^{-\nf 12} |
	\qand
	|\mathbf{D}^{-\nf 12} \mathbf{E} \mathbf{D}^{-\nf 12} - \Itwod| =
	|\mathbf{D}_{\h}^{-\nf 12} \mathbf{E}_{\h} \mathbf{D}_{\h}^{-\nf 12} - \Itwod| \,.
\end{equation*}
\end{lemma}
\begin{proof}
The matrices~$\mathbf{D}_{\h}^{-1} \mathbf{E}_{\h}$,~$\mathbf{D}^{-1} \mathbf{E}$ are similar and hence share the same eigenvalues.
Indeed, for any invertible matrix~$\mathbf{P} \in \R^{2d \times 2d}$,
\begin{equation*}
	(\mathbf{P}^t \mathbf{D} \mathbf{P})^{-1} (\mathbf{P}^t \mathbf{E} \mathbf{P}) = \mathbf{P}^{-1} \mathbf{D}^{-1} \mathbf{P}^{-t} \mathbf{P}^t \mathbf{E} \mathbf{P} = \mathbf{P}^{-1} \mathbf{D}^{-1} \mathbf{E} \mathbf{P} \,\,.
\end{equation*}
Similarly, conjugating by~$\mathbf{D}^{\nf 12}$ shows that~$\mathbf{D}^{-1} \mathbf{E}$ and~$\mathbf{D}^{-\nf12} \mathbf{E} \mathbf{D}^{-\nf 12}$ have the same spectrum. The normalized matrices are symmetric, so equality of their spectra gives both asserted operator-norm identities.
\end{proof}

\subsection{Coarse-grained ellipticity constants and functional inequalities}
\label{ss.cg.ellipticity}

Following~\cite[Section~5.4]{AK.HC}, we attach to each cube~$\cu_m$ a pair of multiscale
ellipticity constants~$\lambda_{s,q}(\cu_m;\a)$ and~$\Lambda_{s,q}(\cu_m;\a)$.  These are scale-dependent versions of the usual lower and upper ellipticity constants: they are built from the coarse-grained matrices~$\s_*(U)$ and~$\b(U)$ over all subcubes~$U = z+\cu_k \subseteq \cu_m$, with weights~$3^{-s(m-k)}$ favoring scales close to~$3^m$. The point is that these \emph{coarse-grained ellipticity constants} are scale-local and thus more forgiving than uniform ellipticity, yet they still provide genuine elliptic estimates at scale~$\cu_m$. In particular, we will present \emph{coarse-grained} versions of the Poincar\'e and Caccioppoli inequalities which are expressed directly in terms of~$\lambda_{s,q}$ and~$\Lambda_{s,q}$.

\smallskip

\subsubsection{Definitions and basic properties of the multiscale quantities}

We begin with the definition of~$\lambda_{s,q}(\cu_m;\a)$ and~$\Lambda_{s,q}(\cu_m;\a)$.

\begin{definition}[Coarse-grained ellipticity constants]
Let~$\cs \coloneqq 1-3^{-s}$ so that~$n\mapsto \cs 3^{-ns}$ is a probability mass function on~$\N_0$.
For every~$s\in (0,\infty)$,~$q\in [1,\infty)$,~$m\in\Z$ and coefficient field~$\a: \cu_m \to \R^{d\times d}_+$, we define the multiscale composite quantities
\begin{equation}
	\label{e.coarse.grained.ellipticity}
	\left\{
	\begin{aligned}
		& {\Lambda}_{s,q}(\cu_m\,;\a)
		\coloneqq
		\biggl(
		\css{sq} \sum_{k=-\infty}^{m}
		3^{-sq(m-k)}
		\max_{z\in 3^k\Zd \cap \cu_m}
		\bigl| \b(z+\cu_k; \a) \bigr|^{\nf q2}
		\biggr)^{\!\nf{2}{q}}
		\,, \\  &
		{\lambda}_{s,q}(\cu_m;\a)
		\coloneqq
		\biggl(\css{sq} \sum_{k=-\infty}^{m}
		3^{-sq(m-k)}
		\max_{z\in 3^k\Zd \cap \cu_m}
		\bigl| \s_{*}^{-1}(z+\cu_k; \a) \bigr|^{\nf q2}
		\biggr)^{\!- \nf{2}{q}}
		\,,
	\end{aligned}
	\right.
\end{equation}
and, for~$q=\infty$,
\begin{equation}
	\label{e.coarse.grained.ellipticity.infty}
	\left\{
	\begin{aligned}
		& {\Lambda}_{s,\infty}(\cu_m\,;\a)
		\coloneqq
		\sup_{k\in(-\infty,m]\cap \Z}
		3^{-2s(m-k)}
		\max_{z\in 3^k\Zd \cap \cu_m}
		\bigl| \b(z+\cu_k; \a) \bigr|
		\,, \\  &
		{\lambda}_{s,\infty}(\cu_m;\a)
		\coloneqq
		\biggl(
		\sup_{k\in(-\infty,m]\cap \Z}
		3^{-2s(m-k)}
		\max_{z\in 3^k\Zd \cap \cu_m}
		\bigl| \s_{*}^{-1}(z+\cu_k; \a) \bigr|
		\biggr)^{\!-1}
		\,.
	\end{aligned}
	\right.
\end{equation}
Note that~\eqref{e.coarse.grained.ellipticity.infty} is the limit as~$q \to \infty$ of~\eqref{e.coarse.grained.ellipticity}.
If we do not write the exponent~$q$, then it is implied that~$q=1$, that is,
\begin{equation*}
	\Lambda_{s} (\cu_m;\a)
	=
	\Lambda_{s,1} (\cu_m;\a)
	\qqand
	\lambda_{s} (\cu_m;\a)
	=
	\lambda_{s,1} (\cu_m;\a)
	\,.
\end{equation*}
We extend these definitions to translations~$y+\cu_m$ of the cube~$\cu_m$ in the obvious way.
\end{definition}

\begin{lemma}[Monotonicity of the coarse-grained ellipticity constants]
\label{l.ellipticities.monotone.exponents}
For every~$m\in\Z$,~$0<s\leq t<\infty$ and~$1\leq p\leq q\leq\infty$,
\begin{equation}
	\lambda_{s,p}(\cu_m;\a)
	\leq
	\lambda_{t,q}(\cu_m;\a)
	\qand
	\Lambda_{t,q}(\cu_m;\a)
	\leq
	\Lambda_{s,p}(\cu_m;\a)
	\,.
	\label{e.ellipticities.change.q}
\end{equation}
In particular, for~$q\in[1,\infty]$ and~$0<t<s<\infty$,
\begin{equation}
	\lambda_{t,q}(\cu_m; \a)
	\leq
	\lambda_{s,q} (\cu_m; \a)
	\leq |\s_*^{-1}(\cu_m;\a) |^{-1}
	\leq
	|\b(\cu_m;\a)|
	\leq
	\Lambda_{s,q}(\cu_m; \a)
	\leq
	\Lambda_{t,q}(\cu_m; \a)
	\,.
	\label{e.ellipticities.monotone.ordered}
\end{equation}
\end{lemma}
\begin{proof}
For either of the two quantities in~\eqref{e.coarse.grained.ellipticity}, let~$a_n$ be the square root of the corresponding maximum over the cubes of size~$3^{m-n}$. By subadditivity~\eqref{e.subadda.nosymm}, the sequence~$\{a_n\}_{n\in\N_0}$ is nondecreasing. For~$p<\infty$, let~$N_{s,p}$ have the geometric law
\begin{equation*}
	\P_{\mathrm g}[N_{s,p}=n]
	=
	\css{sp}3^{-spn}
	\,.
\end{equation*}
Then the square roots of~$\Lambda_{s,p}$ and~$\lambda_{s,p}^{-1}$ are the corresponding quantities~$\|a_{N_{s,p}}\|_{L^p}$. Indeed, if~$s\leq t$, then~$\P_{\mathrm g}[N_{t,p}\geq n]=3^{-tpn}\leq3^{-spn}=\P_{\mathrm g}[N_{s,p}\geq n]$ for every~$n\in\N_0$. Since~$a_n$ is nondecreasing, this proves monotonicity in~$s$.
\smallskip
For~$p\leq q<\infty$, we also have, for every~$r\geq0$,
\begin{equation*}
	\P_{\mathrm g}[a_{N_{s,q}}>r]
	=
	\P_{\mathrm g}[a_{N_{s,p}}>r]^{\nf qp}
	\,.
\end{equation*}
Writing~$\alpha=\nf qp$ and~$F(r)=\P_{\mathrm g}[a_{N_{s,p}}^p>r]$, the layer-cake formula and~$rF(r)\leq\int_0^rF$ give
\begin{align*}
	\E_{\mathrm g}[a_{N_{s,q}}^q]
	&=
	\alpha\int_0^\infty r^{\alpha-1}F(r)^\alpha\,dr
	\leq
	\biggl(\int_0^\infty F(r)\,dr\biggr)^{\!\alpha}
	=
	\E_{\mathrm g}[a_{N_{s,p}}^p]^{\nf qp}
	\,.
\end{align*}
\smallskip
If~$p<q=\infty$, the same conclusion follows by bounding each~$3^{-sn}a_n$ by~$\|a_{N_{s,p}}\|_{L^p}$. The case~$p=q=\infty$, including monotonicity in~$s$, is immediate from the supremum definition. This proves~\eqref{e.ellipticities.change.q}. The second and fourth inequalities in~\eqref{e.ellipticities.monotone.ordered} follow from the monotonicity of~$a_n$, and the middle inequality is~\eqref{e.CG.bounds.1}.
\end{proof}

The definition also gives, for every~$k\in (-\infty,m]\cap \Z$,
\begin{equation}
	\label{e.bound.one.cube.by.lambdas}
	\left\{
	\begin{aligned}
		&
		\max_{z\in 3^k\Zd\cap \cu_m}
		\bigl| \b(z+\cu_k; \a) \bigr|
		\leq
		\max_{z\in 3^k\Zd\cap \cu_m}
		\Lambda_{s,q} (z+\cu_k; \a)
		\leq 3^{2s(m-k)} \Lambda_{s,q} (\cu_m;\a)
		\,,
		\\ &
		\max_{z\in 3^k\Zd\cap \cu_m}
		\bigl| \s_{*}^{-1}(z+\cu_k; \a) \bigr|
		\leq
		\max_{z\in 3^k\Zd\cap \cu_m}
		\lambda_{s,q}^{-1} (z+\cu_k; \a)
		\leq 3^{2s(m-k)} \lambda_{s,q}^{-1} (\cu_m;\a)
		\,.
	\end{aligned}
	\right.
\end{equation}
It is easy to check from the Lebesgue differentiation theorem that
\begin{equation*}
	\text{$\a$ is uniformly elliptic}
	\quad \iff \quad
	\lim_{s\to 0} \lambda_s(\cu_m; \a) >0
	\qand
	\lim_{s\to 0} \Lambda_s(\cu_m; \a) <\infty\,,
\end{equation*}
in which case the two limits on the right coincide with the lower and upper constants of uniform ellipticity, respectively, by~\eqref{e.CG.bounds.1} and the Lebesgue differentiation theorem. Note that, if~$\lambda_s(\cu_m; \a)$ is positive for some~$s>0$, then~$\lim_{s\to \infty} \lambda_s(\cu_m; \a) = |\s_*^{-1}(\cu_m;\a) |^{-1}$ and, similarly, if~$\Lambda_s(\cu_m; \a)$ is finite for some~$s>0$, then~$\lim_{s\to \infty} \Lambda_s(\cu_m; \a) = |\b(\cu_m;\a)|$. Thus the pair~$( \lambda_s(\cu_m; \a), \Lambda_s(\cu_m; \a))$ interpolates between the uniform ellipticity constants and~$( |\s_*^{-1}(\cu_m;\a) |^{-1} , |\b(\cu_m;\a)|)$,  the smallest and largest eigenvalues of the coarse-grained matrices.

In parallel with ellipticity, we need a multiscale control of the \emph{error in the linear response} when we compare~$-\nabla\cdot\a\nabla$ to a constant-coefficient operator~$-\nabla\cdot\a_0\nabla$.  This is quantified by the homogenization error~$\mathcal{E}_{s,p,q}(\cu_m;\a,\a_0)$. Its building blocks are the quantities
\begin{equation*}
	\bfJ\bigl(z+\cu_l,\bfA_0^{-\nf12}e,\bfA_0^{\nf12} e\,;\a\bigr)\,,
	\qquad
	e\in\R^{2d}\,, \
	|e|=1\,, \
	l \in \Z \,, \
	z \in 3^l \Zd \cap \cu_m
	\,.
\end{equation*}
Each of these quantities, by the flux-gradient comparison inequality~\eqref{e.fluxmaps}, measures how far the response of~$-\nabla\cdot\a\nabla$ is from that of~$-\nabla\cdot\a_0\nabla$ in the cube~$z+\cu_l$. The quantity~$\mathcal{E}_{s,p,q}(\cu_m;\a,\a_0)$ aggregates these local response errors over all subcubes of~$\cu_m$, with the same geometric discounting weights~$3^{-s(n-l)}$ as in the definition of~$\lambda_{s,q}$ and~$\Lambda_{s,q}$. It is therefore a scale-local version of a homogenization error built from the coarse-grained matrices, measuring not just their size (like~$\lambda_{s,q}$ and~$\Lambda_{s,q}$ do) but their closeness to a particular given matrix~$\a_0$.

\smallskip

In the RG argument,~$\mathcal{E}_{s,p,q}$ is the basic state variable. It allows us to compare the operators at two successive scales and close the iteration purely at the level of coarse-grained coefficients. We turn now to its definition.

\begin{definition}[Homogenization error]
\label{d.mathcal.E}
For every~$s \in (0,\infty)$,~$p,q\in [1,\infty]$,~$m,n \in \Z$ with~$n\leq m$, coefficient field~$\a: \cu_m \to \R^{d\times d}_+$ and matrix~$\a_0 \in \R^{d\times d}_+$, we let~$\bfA_0$ denote the matrix defined in terms of~$\a_0$ by the usual formula~\eqref{e.bfA.def}, that is,
\begin{equation}
	\label{e.form.of.A.naught}
	\bfA_0  \coloneqq
	\begin{pmatrix}
		\s_0 + \k_0^t\s_0^{-1}\k_0
		& -\k_0^t\s_0^{-1}
		\\ - \s_0^{-1}\k_0
		& \s_0^{-1}
	\end{pmatrix}
	\quad \mbox{where} \quad
	\s_0 \coloneqq  \frac12(\a_0+\a_0^t)
	\,,
	\quad
	\k_0 \coloneqq  \frac12(\a_0-\a_0^t)
	\,,
\end{equation}
We then define the multiscale composite quantity
\begin{align*}
	\lefteqn{
	\mathcal{E}_{s,p,q} (\cu_m , n; \a, \a_0 )
	} \ \notag \\ &
	=
	\left\{
	\begin{aligned}
		&
		\biggl(
		\css{sq}
		\sum_{l=-\infty}^n
		3^{-sq(n-l)}
		\biggl(
		\avsum_{z \in 3^l\Zd \cap \cu_m}
		\max_{|e| =1} \bigl( \bfJ\bigl(z {+} \cu_l,\bfA_0^{-\nf 12} e, \bfA_0^{\nf 12} e \,; \a\bigr) \bigr)^{\nf p2}
		\biggr)^{\!\nf qp}
		\biggr)^{\!\nf1q}
		\,,
		\quad  p,q \in [1,\infty)\,,
		\\ &
		\biggl(
		\css{sq}
		\sum_{l=-\infty}^n
		3^{-sq(n-l)} \!\!
		\max_{z \in 3^l\Zd \cap \cu_m}
		\max_{|e| =1} \bigl( \bfJ\bigl(z {+} \cu_l,\bfA_0^{-\nf 12} e, \bfA_0^{\nf 12} e \,; \a\bigr) \bigr)^{\nf q2}
		\biggr)^{\!\nf1q}
		\,,
		\quad p =\infty,\,  q \in [1,\infty)\,,
		\\
		&
		\sup_{l \in (-\infty,n]\cap \Z}
		3^{-s(n-l)}
		\biggl(
		\avsum_{z \in 3^l\Zd \cap \cu_m}
		\max_{|e| =1} \bigl( \bfJ\bigl(z {+} \cu_l,\bfA_0^{-\nf 12} e, \bfA_0^{\nf 12} e \,; \a\bigr) \bigr)^{\nf p2}
		\biggr)^{\!\nf 1p}
		\,,
		\quad p \in [1,\infty),\, q =\infty \,,
		\\ &
		\sup_{l \in (-\infty,n]\cap \Z}
		3^{-s(n-l)}
		\max_{z \in 3^l\Zd \cap \cu_m}
		\max_{|e| =1} \bigl( \bfJ\bigl(z {+} \cu_l,\bfA_0^{-\nf 12} e, \bfA_0^{\nf 12} e \,; \a\bigr) \bigr)^{\nf 12}
		\,,
		\quad p = q = \infty \,.
	\end{aligned}
	\right.
\end{align*}
We extend these definitions to translations~$y+\cu_m$ of~$\cu_m$ in the obvious way.
If we drop the argument~$n$, then it is to be inferred that~$n=m$; that is,
\begin{equation*}
	\mathcal{E}_{s,p,q} (\cu_m ; \a, \a_0 )
	\coloneqq
	\mathcal{E}_{s,p,q} (\cu_m , m; \a, \a_0 )
	\,.
\end{equation*}
\end{definition}

Below we will show, in Lemma~\ref{l.mathcal.E.to.Lambdas}, that~$\mathcal{E}_{s,p,q}$ encodes essentially the same quantitative information as the ellipticity constants~$\Lambda_{s,q}$ and~$\lambda_{s,q}$, while Lemma~\ref{l.coarse.graining.operator} identifies~$\mathcal{E}_{s,p,q}$ as the coarse-level quantity controlling the operator error~$(\a-\a_0)\nabla u$ in negative Besov norms.

\smallskip

We continue with some of the basic properties of~$\mathcal{E}_{s,p,q}$, starting with the following analogues of~\eqref{e.ellipticities.monotone.ordered} and~\eqref{e.bound.one.cube.by.lambdas}.
For every~$m\in\Z$,~$p \in [2,\infty]$,~$q \in [1,\infty]$ and~$t,s\in (0,\infty)$ with~$t<s$, we have
\begin{equation}
	\max_{|e| =1}  \bfJ\bigl(\cu_m,\bfA_0^{-\nf 12} e, \bfA_0^{\nf 12} e \,; \a\bigr)^{\nf12}
	\leq
	\mathcal{E}_{s,p,q} (\cu_m ; \a, \a_0 )
	\leq
	\mathcal{E}_{t,p,q} (\cu_m ; \a, \a_0 )
	\label{e.mathcalE.monotone.ordered}
\end{equation}
and, for every~$k\in (-\infty,m]\cap \Z$,
\begin{equation}
	\label{e.bound.one.cube.by.mathcalE}
	\left\{
	\begin{aligned}
		&
		\max_{z\in 3^k\Zd\cap \cu_m}
		\mathcal{E}_{s,\infty,q} (z+\cu_k;\a, \a_0)
		\leq 3^{s(m-k)} \mathcal{E}_{s,\infty,q}  (\cu_m;\a,\a_0)
		\,,
		\\ &
		\avsum_{z\in 3^k\Zd\cap \cu_m}
		\mathcal{E}_{s,p,q} (z+\cu_k;\a, \a_0)
		\leq 3^{s(m-k)} \mathcal{E}_{s,p,q}  (\cu_m;\a,\a_0)
		\,.
	\end{aligned}
	\right.
\end{equation}
We next observe that by Definition~\ref{d.mathcal.E}, for every~$p,q \in [1,\infty]$ with~$q\leq p$ and~$s \in (\nf{d}{p},1)$,
\begin{equation}
	\label{e.mathcalE.infty.to.q}
	\mathcal{E}_{s, \infty, q} (\cu_m,n ; \a, \a_0 )
	\leq
	3^{\frac{d}{p}(m-n)}
	\Bigl(
	\frac{\css{sq}}{\css{q(s-\nf{d}{p})}}
	\Bigr)^{\nf1q}
	\mathcal{E}_{s - \nf{d}{p}, p, q}(\cu_m,n;\a,\a_0) \,.
\end{equation}
In the case~$q < \infty$,
\begin{align*}
	\mathcal{E}_{s, \infty, q} (\cu_m,n ; \a, \a_0 )^q
	&
	=
	\css{sq}
	\sum_{j=-\infty}^n
	3^{-sq(n-j)} \!\!
	\max_{z \in 3^j\Zd \cap \cu_m}
	\max_{|e| =1} \bigl( \bfJ\bigl(z+\cu_j,\bfA_0^{-\nf 12} e, \bfA_0^{\nf 12} e \,; \a\bigr)\bigr)^{\nf q2}
	\notag \\ &
	\leq
	\css{sq}
	\sum_{j=-\infty}^n
	3^{-sq(n-j)}
	\biggl(
	\sum_{z \in 3^j\Zd \cap \cu_m}
	\!\max_{|e| =1} \bigl( \bfJ\bigl(z+\cu_j,\bfA_0^{-\nf 12} e, \bfA_0^{\nf 12} e \,; \a\bigr)\bigr)^{\frac{p}{2}}
	\biggr)^{\!\nf qp}
	\notag \\ &
	=
	\css{sq}
	\sum_{j=-\infty}^n
	3^{-sq(n-j) + \frac{q}{p} d (m-j)}
	\biggl(
	\avsum_{z \in 3^j\Zd \cap \cu_m}
	\!\max_{|e| =1} \bigl( \bfJ\bigl(z+\cu_j,\bfA_0^{-\nf 12} e, \bfA_0^{\nf 12} e \,; \a\bigr)\bigr)^{\frac{p}{2}}
	\biggr)^{\!\nf qp}
	\notag \\ &
	=
	3^{\frac{q}{p}d(m-n)}
	\frac{\css{sq}}{\css{sq-\nf{dq}{p}}}
	\mathcal{E}_{s - \nf{d}{p}, p, q}(\cu_m,n;\a,\a_0)^q
	\,.
\end{align*}
The case~$q = \infty$ follows by sending~$q \to \infty$ in the above. We also observe that, for every~$q,p \in [1, \infty)$ with~$p \geq q$ and~$s > 0$,
\begin{equation}
	\mathcal{E}_{s, \infty, q}(\cu_m; \a, \a_0)
	\leq
	\mathcal{E}_{\frac{s q}{p}, \infty, p}(\cu_m; \a, \a_0)
	\,.
	\label{e.compareEqs}
\end{equation}
Indeed, by Jensen's inequality,
\begin{align*}
	\mathcal{E}_{s,\infty,q}(\cu_m;\a,\a_0)^q
	&
	=
	\css{sq}
	\sum_{j=-\infty}^m
	3^{-sq(m-j)}
	\max_{z\in3^j\Zd\cap\cu_m}
	\max_{|e|=1}
	\bigl( \bfJ\bigl(z+\cu_j,\bfA_0^{-\nf 12} e, \bfA_0^{\nf 12} e \,; \a\bigr)\bigr)^{\nf q2}
	\notag \\ &
	\leq
	\Bigl(
	\css{sq}
	\sum_{j=-\infty}^m
	3^{-sq(m-j)}
	\max_{z\in3^j\Zd\cap\cu_m}
	\max_{|e|=1}
	\bigl( \bfJ\bigl(z+\cu_j,\bfA_0^{-\nf 12} e, \bfA_0^{\nf 12} e \,; \a\bigr)\bigr)^{\nf p 2}
	\Bigr)^{\nf q p}
	\notag \\ &
	=
	\mathcal{E}_{\frac{s q}{p}, \infty, p}(\cu_m; \a, \a_0)^{q}
	\,.
\end{align*}
Similarly,
\begin{equation}
	\Lambda_{s, q}(\cu_m; \a)
	\leq
	\Lambda_{\frac{s q}{p}, p}(\cu_m; \a)
	\qand
	\lambda_{s, q}^{-1}(\cu_m; \a)
	\leq
	\lambda_{\frac{s q}{p}, p}^{-1}(\cu_m; \a) \,\,.
	\label{e.compareLambdaqs}
\end{equation}
Finally, we note that, for every~$m,n\in\Z$ with~$n\leq m$, every coefficient field~$\a$ on~$\cu_m$, every scalar comparator~$\s_0>0$, every~$t \in (0,1]$ and~$q \in [1,\infty]$,
\begin{equation}
	\label{e.mathcal.E.t.infty.one}
	\max_{z \in 3^n \Zd \cap \cu_m}
	\mathcal{E}_{t, \infty, q} (z{+}\cu_n;\a,\s_0)
	\leq
	\mathcal{E}_{t, \infty, q} (\cu_m,n;\a,\s_0)
	\,.
\end{equation}

We next record a simple stability result for coarse-grained ellipticity constants.

\begin{lemma}
\label{l.lambdas.stability}
There exists a constant~$C(d)<\infty$ such that, for every~$q \in [1,\infty]$,~$s,t \in (0,\nf 12]$ with~$s<t$, and~$x \in \cu_{0}$ such that~$x +\cu_{-1} \subseteq \cu_0$, we have the estimates
\begin{equation}
	\label{e.lambda.stability}
	\lambda_{t,q}^{-1}(x+\cu_{-1}\,;\a)
	\leq
	\frac{C}{1-2s} \Bigl( \frac{t}{t-s}\Bigr)^{\! \nf 2q} \lambda_{s,q}^{-1}(\cu_{0}\,;\a) \,,
\end{equation}
\begin{equation}
	\label{e.big.Lambda.stability}
	\Lambda_{t,q}(x+\cu_{-1}\,;\a)
	\leq
	\frac{C}{1-2s} \Bigl( \frac{t}{t-s}\Bigr)^{\! \nf 2q} \Lambda_{s,q}(\cu_{0}\,;\a)\,.
\end{equation}
For every comparator~$\a_0$ as in Definition~\ref{d.mathcal.E}, we also have
\begin{equation}
	\label{e.mathcalE.stability}
	\mathcal{E}_{t,\infty,2}(x+\cu_{-1}\,;\a,\a_0)
	\leq
	\frac{C}{(1-2s)^{\nf12}} \Bigl( \frac{t}{t-s}\Bigr)^{\!\nf12} \mathcal{E}_{s,\infty,2}(\cu_0\,;\a,\a_0)
	\,.
\end{equation}
\end{lemma}
\begin{proof}
Let~$U$ be any domain in~$\cu_0$. For~$n \in -\N$, we define recursively~$\mathcal{G}_0(U) \coloneqq \emptyset$ and
\begin{equation*}
	\mathcal{G}_n(U) \coloneqq \biggl\{ z \in 3^n \Zd \setminus \bigcup_{k \in (n,0] \cap \Z} \bigcup_{y \in \mathcal{G}_k(U)} (y+\cu_k) \, : \, z+\cu_{n+1} \subset U \biggr\} \,.
\end{equation*}
Then, by subadditivity and~\eqref{e.bound.one.cube.by.lambdas},
\begin{equation*}
	| \s_*^{-1}(U;\a)| \leq \sum_{n=-\infty}^0 \sum_{z \in \mathcal{G}_n(U)} \frac{|\cu_n|}{|U|} | \s_*^{-1}(z+\cu_n\,;\a)|
	\leq \lambda_{s,q}^{-1}(\cu_0\,;\a)  \sum_{n=-\infty}^0
	3^{-2sn}  \sum_{z \in \mathcal{G}_n(U)} \frac{|\cu_n|}{|U|}
	\,.
\end{equation*}
For every~$y \in \cu_{-1}$ and~$k \in -\N$ and~$n<k$, we have~$\sum_{z \in \mathcal{G}_n(y+\cu_{k})} \frac{|\cu_n|}{|\cu_k|} \leq C 3^{n-k}$, and it follows that
\begin{align*}
	\lambda_{t,q}^{-\nf q2}(x+\cu_{-1}\,;\a)
	&
	=
	\css{tq} \sum_{k=-\infty}^{-1} 3^{t q (1+k)}
	\max_{y \in x+(3^k \Zd \cap \cu_{-1})}
	| \s_*^{-1}(y+\cu_k\,;\a)|^{\nf q2}
	\notag \\ &
	\leq
	C^q \css{tq}  \lambda_{s,q}^{-\nf q2}(\cu_0\,;\a)
	\sum_{k=-\infty}^{-1} 3^{t q k}
	\biggl( \sum_{n=-\infty}^k  3^{n-k - 2sn} \biggr)^{\! \nf q2}
	\notag \\ &
	\leq
	\frac{C^q t }{(1-2s)^{\nf q2}(t-s)} \lambda_{s,q}^{-\nf q2}(\cu_0\,;\a)
	\,.
\end{align*}
For~$q<\infty$, taking the power~$2/q$ gives~\eqref{e.lambda.stability}; the case~$q=\infty$ follows from the same decomposition with suprema. The argument for~$\Lambda$ is identical. Applying the same decomposition to~$\bfJ$ and taking the square root gives~\eqref{e.mathcalE.stability}.
\end{proof}

The following lemma is nearly the same as~\cite[Lemma 5.9]{AK.HC}. We present the proof for the reader's convenience.

\begin{lemma}
\label{l.mathcal.E.to.Lambdas}
For every~$m\in\Z$,~$q \in [2,\infty]$,~$s\in (0,\infty)$, coefficient field~$\a: \cu_m \to \R^{d\times d}_+$ and scalar~$\s_0 \in (0,\infty)$,
\begin{align}
	\frac12 \mathcal{E}_{s,\infty,q}
	(\cu_m;\a,\s_0)^2
	&
	\leq
	\max\bigl\{
	\s_0^{-1} \Lambda_{s,q}(\cu_m;\a)
	\,,
	\s_0 \lambda_{s,q}^{-1} (\cu_m;\a)
	\bigr\}
	\notag \\ &
	\leq
	1+
	2 \mathcal{E}_{s,\infty,q} (\cu_m;\a,\s_0)^2
	+
	2^{\nf12} \mathcal{E}_{s,\infty,q}(\cu_m;\a,\s_0)
	\,.
	\label{e.bound.Lambdas.by.Es}
\end{align}
Moreover, in the case~$q = 2$,
\begin{equation}
	\label{e.bound.Lambdas.by.Es.q2}
	1+\frac12 \mathcal{E}_{s,\infty,2}
	(\cu_m;\a,\s_0)^2
	\leq
	\max\bigl\{
	\s_0^{-1} \Lambda_{s,2}(\cu_m;\a)
	\,,
	\s_0 \lambda_{s,2}^{-1} (\cu_m;\a)
	\bigr\} \,\,.
\end{equation}
\end{lemma}
\begin{proof}
The statement in the case~$q = \infty$ follows by sending~$q \to \infty$ in~\eqref{e.bound.Lambdas.by.Es}. Thus, we may assume~$q \in [2,\infty)$. By~\eqref{e.bfcoarsegrainedmatrices}, for every~$e \in \R^{2d}$ with~$|e|=1$,~$\bfA_0 \in \R^{2d \times 2d}_{\sym,+}$ and Lipschitz domain~$U \subseteq \Rd$,
\begin{equation*}
	\bfJ(U,\bfA_0^{-\nf12}e,\bfA_0^{\nf12} e ; \a)
	=
	\frac12 e \cdot \bfA_0^{-\nf12} \bfA(U;\a) \bfA_0^{-\nf12} e
	+
	\frac12 e \cdot \bfA_0^{\nf12} \bfA_*^{-1} (U;\a) \bfA_0^{\nf12} e
	- 1
	\,.
\end{equation*}
In the case that~$\a_0 = \s_0 \in (0,\infty)$ is a scalar matrix and~$\bfA_0$ has the form of~\eqref{e.form.of.A.naught}, we obtain that
\begin{align}
	\label{e.J.by.f}
	\frac12
	\bigl|
	\s_0^{-1} \b(U;\a)
	+
	\s_0 \s_*^{-1} (U;\a)
	-2 \Id
	\bigr|
	&
	\leq
	\max_{|e|=1} \bfJ(U,\bfA_0^{-\nf12}e,\bfA_0^{\nf12} e ; \a)
	\notag \\ &
	\leq
	\bigl|
	\s_0^{-1} \b(U;\a)
	+
	\s_0 \s_*^{-1} (U;\a)
	-2 \Id
	\bigr|
	\,.
\end{align}
Using also that~$\s_0^{-1} \b(U;\a)
+
\s_0 \s_*^{-1} (U;\a)
-2 \Id \geq 0$, which implies
\begin{align*}
	\bigl|
	\s_0^{-1} \b(U;\a)
	+
	\s_0 \s_*^{-1} (U;\a)
	-2 \Id
	\bigr|
	&
	=
	\bigl|
	\s_0^{-1} \b(U;\a)
	+
	\s_0 \s_*^{-1} (U;\a)
	\bigr|
	-2
	\notag \\ &
	\leq
	\s_0^{-1} | \b(U;\a)| + \s_0 |\s_*^{-1} (U;\a)| - 2
	\,,
\end{align*}
positivity of the coarse block matrix and Minkowski's inequality applied to the square roots give, for every scalar matrix~$\s_0\in (0,\infty)$,
\begin{align*}
	\mathcal{E}_{s,\infty,q}
	(\cu_m;\a,\s_0) ^2
	&
	\leq
	\s_0^{-1} \Lambda_{s,q}(\cu_m;\a) +
	\s_0 \lambda_{s,q}^{-1} (\cu_m;\a)
	\leq
	2 \max\bigl\{
	\s_0^{-1} \Lambda_{s,q}(\cu_m;\a)
	\,,
	\s_0 \lambda_{s,q}^{-1} (\cu_m;\a)
	\bigr\}
	\,,
\end{align*}
which implies the first inequality of~\eqref{e.bound.Lambdas.by.Es}.
In the case~$q = 2$, we can replace the first term on the left above by~$\mathcal{E}_{s,\infty,2}
(\cu_m;\a,\s_0) ^2  + 2$. This implies~\eqref{e.bound.Lambdas.by.Es.q2}.

To prove the second inequality, we use that
\begin{align}
	\label{e.xminusonetimesxminusonesquared.sstar}
	\lefteqn{
	\bigl(
	\s_0 \s_*^{-1} (U;\a)
	- \Id
	\bigr)
	\s_0^{-1}  \s_*(U;\a)
	\bigl(
	\s_0 \s_*^{-1} (U;\a)
	- \Id
	\bigr)
	} \qquad &
	\notag \\ &
	\leq
	\s_0^{-1} ( \b(U;\a) - \s_* (U;\a))
	+
	\bigl(
	\s_0 \s_*^{-1} (U;\a)
	- \Id
	\bigr)
	\s_0^{-1}  \s_*(U;\a)
	\bigl(
	\s_0 \s_*^{-1} (U;\a)
	- \Id
	\bigr)
	\notag \\ &
	=
	\s_0^{-1} \b(U;\a)
	+
	\s_0 \s_*^{-1} (U;\a)
	-2 \Id
	\,.
\end{align}
Similarly,
\begin{equation}
	\label{e.xminusonetimesxminusonesquared.b}
	\bigl(
	\s_0^{-1} \b(U;\a)
	- \Id
	\bigr)
	\s_0 \b^{-1}(U;\a)
	\bigl(
	\s_0^{-1} \b(U;\a)
	- \Id
	\bigr)
	\leq \s_0^{-1} \b(U;\a)
	+
	\s_0 \s_*^{-1} (U;\a)
	-2 \Id\,.
\end{equation}
An eigenvector of~$\s_0 \s_*^{-1} (U;\a)$ with eigenvalue~$\mu$ is also an eigenvector of the left side of~\eqref{e.xminusonetimesxminusonesquared.sstar} with eigenvalue~$f(\mu)$, where we set~$f(x)  \coloneqq  x^{-1}(x-1)^2$. It follows that
\begin{equation*}
	f\bigl( \s_0 |\s_*^{-1} (U;\a)| \bigr)
	\leq
	\bigl|
	\s_0^{-1} \b(U;\a)
	+
	\s_0 \s_*^{-1} (U;\a)
	-2 \Id
	\bigr| \,.
\end{equation*}
For every~$\delta,\ep>0$, we have that~$f(x) \leq \delta$ implies~$x-1 \leq \delta + \delta^{\nf12} \leq (1+\ep^{-1} )\delta + \frac14\ep$. In view of the first inequality in~\eqref{e.J.by.f}, we obtain
\begin{equation*}
	\s_0 |\s_*^{-1} (U;\a)| -1
	\leq
	2 (1+\ep^{-1} ) \max_{|e|=1} \bfJ(U,\bfA_0^{-\nf12}e,\bfA_0^{\nf12} e ; \a)
	+
	\frac14 \ep\,.
\end{equation*}
This implies that
\begin{equation*}
	\s_0 \lambda_{s,q}^{-1} (\cu_m;\a) - 1
	\leq
	2 (1+\ep^{-1} ) \mathcal{E}_{s,\infty,q}
	(\cu_m;\a,\s_0)^2
	+
	\frac14 \ep
	\,.
\end{equation*}
For~$q\geq2$, optimizing in~$\ep$ and using the ordinary triangle inequality yields
\begin{equation*}
	\s_0 \lambda_{s,q}^{-1} (\cu_m;\a) - 1
	\leq
	2 \mathcal{E}_{s,\infty,q} (\cu_m;\a,\s_0)^2
	+
	2^{\nf12} \mathcal{E}_{s,\infty,q}(\cu_m;\a,\s_0)
	\,.
\end{equation*}
The upper bound for~$\s_0^{-1} \Lambda_{s,q}(\cu_m;\a)$ is obtained similarly, using~\eqref{e.xminusonetimesxminusonesquared.b} in place of~\eqref{e.xminusonetimesxminusonesquared.sstar}.
\end{proof}

For~$q\in[1,2)$, the ellipticity constants are controlled by the~$q=2$ ellipticities at halved exponent: by Jensen's inequality and the monotonicity~\eqref{e.ellipticities.change.q},~$\Lambda_{s,q}(\cu_m;\a)\leq\Lambda_{\nf s2,2}(\cu_m;\a)$ and~$\lambda_{s,q}^{-1}(\cu_m;\a)\leq\lambda_{\nf s2,2}^{-1}(\cu_m;\a)$; we use this below without further comment.

\subsubsection{Coarse-grained elliptic and functional inequalities}

We next give the statements of two very important coarse-grained estimates which we use many times in the paper. The first of these is the coarse-grained Poincar\'e inequality, which is identical to~\cite[Lemma 2.3]{AK.HC}. We nevertheless present the full proof since it is short and conveys a lot of intuition about how our multiscale arguments work.

\begin{proposition}[Coarse-grained Poincar\'e inequality]
\label{p.coarse.grained.Poincare}
For every~$s\in (0,1]$,~$q\in [1,\infty]$,~$m\in\Z$ and solution~$u\in \mathcal{A}(\cu_m)$,
\begin{equation}
	3^{-sm} [ \nabla u ]_{\Besov{-s}{2}{q}(\cu_m)}
	\leq
	\css{sq}^{-\nf1q}
	\lambda_{s,q}^{-\nf12}  (\cu_m; \a) \| \s^{\nf 12} \nabla u \|_{\underline{L}^2(\cu_m)}
	\label{e.besov.grad.poincare}
\end{equation}
and
\begin{equation}
	3^{-sm} [ \a \nabla u ]_{\Besov{-s}{2}{q}(\cu_m)}
	\leq
	\css{sq}^{-\nf1q}
	\Lambda^{\nf12 }_{s,q} (\cu_m; \a)  \| \s^{\nf 12} \nabla u \|_{\underline{L}^2(\cu_m)}
	\,\,.
	\label{e.besov.flux.poincare}
\end{equation}
\end{proposition}
\begin{proof}
Select~$u\in \A(\cu_m)$ and apply~\eqref{e.energymaps.nonsymm} to get
\begin{align*}
	3^{-sqm} [ \nabla u ]_{\Besov{-s}{2}{q}(\cu_m)}^q
	& =
	\sum_{k=-\infty}^m  3^{sq (k-m)}
	\biggl( \avsum_{z \in 3^{k} \Zd \cap \cu_m}  \bigl|(\nabla u)_{z+\cu_k} \bigr|^2
	\biggr)^{\! \nf q2}
	\notag \\ &
	\leq
	\sum_{k=-\infty}^m  3^{ sq (k-m)}
	\max_{z \in 3^{k} \Zd \cap \cu_m}
	|\s_*^{-1}(z + \cu_k)|^{\nf q2}
	\biggl(
	\avsum_{z \in 3^{k} \Zd \cap \cu_m}  \|\s^{\nf12} \nabla u \|_{\underline{L}^2(z + \cu_k)}^2
	\biggr)^{\! \nf q2}
	\notag \\ &
	=
	\css{sq}^{-1}
	\|\s^{\nf12} \nabla u \|_{\underline{L}^2(\cu_m)}^q
	\css{sq}
	\sum_{k=-\infty}^m  3^{-sq (m-k)} \max_{z \in 3^{k} \Zd \cap \cu_m} |\s_*^{-1}(z + \cu_k)|^{\nf q2}
	\notag \\ &
	=
	\css{sq}^{-1}
	\lambda_{s,q}^{-\nf q2}(\cu_m\,;\a) \| \s^{\nf 12} \nabla u \|_{\underline{L}^2(\cu_m)}^q
	\,.
\end{align*}
This gives us~\eqref{e.besov.grad.poincare}.
The proof of~\eqref{e.besov.flux.poincare} is similar. We apply~\eqref{e.energymaps.nonsymm.dual} to get
\begin{align*}
	3^{-sq m} [ \a \nabla u ]_{\Besov{-s}{2}{q}(\cu_m)}^q
	& =
	\sum_{k=-\infty}^m  3^{sq (k-m)}
	\biggl( \avsum_{z \in 3^{k} \Zd \cap \cu_m}  \bigl|(\a \nabla u)_{z+\cu_k} \bigr|^2
	\biggr)^{\!\nf q2}
	\notag \\ &
	\leq
	\sum_{k=-\infty}^m  3^{sq (k-m)} \max_{z \in 3^{k} \Zd \cap \cu_m} |\b(z + \cu_k)|^{\nf q2}
	\biggl( \avsum_{z \in 3^{k} \Zd \cap \cu_m}  \|\s^{\nf12} \nabla u \|_{\underline{L}^2(z + \cu_k)}^2 \biggr)^{\!\nf q2}
	\notag \\ &
	=
	\css{sq}^{-1} \Lambda_{s,q}^{\nf q2}  (\cu_m\,;\a) \| \s^{\nf 12} \nabla u \|_{\underline{L}^2(\cu_m)}^q
	\,\,.
\end{align*}
This completes the proof.
\end{proof}

\begin{remark}
\label{r.cg.poincare.doubled.variables}
Proposition~\ref{p.coarse.grained.Poincare} can be formulated equivalently in terms of the ``double-variable solutions,'' that is, elements of the space~$\mathcal{S}(\cu_m)$ defined in~\eqref{e.bfS}.
Let~$\bfA_0$ be defined by~\eqref{e.form.of.A.naught} where~$\s_0$ is a scalar matrix and~$\k_0=0$.  The statement is that, for every~$X \in \mathcal{S}(\cu_m)$,
\begin{align}
	\lefteqn{
	3^{-sm} [ \bfA_0^{\nf12} X ]_{\Besov{-s}{2}{q}(\cu_m)}
	+
	3^{-sm} [ \bfA_0^{-\nf12} \bfA X ]_{\Besov{-s}{2}{q}(\cu_m)}
	} \qquad &
	\notag \\ &
	\leq
	4 \css{sq}^{-\nf1q}
	\bigl(
	\s_0^{\nf12} \lambda_{s,q}^{-\nf12}  (\cu_m; \a)
	+
	\s_0^{-\nf12} \Lambda^{\nf12 }_{s,q} (\cu_m; \a)
	\bigr)
	\| \bfA^{\nf 12} X \|_{\underline{L}^2(\cu_m)}
	\,.
	\label{e.CG.Poincare.doubled.vars}
\end{align}
This is an immediate consequence of Proposition~\ref{p.coarse.grained.Poincare}, the identity
\begin{equation}
	\bfA \begin{pmatrix} p + p^* \\ \a p -  \a^t p^* \end{pmatrix}
	=
	\begin{pmatrix} \a p + \a^t p^* \\ p - p^* \end{pmatrix}
	\label{e.bfA.magic.swapping}
\end{equation}
and the characterization of the space~$\mathcal{S}(\cu_m)$ in~\eqref{e.findSfull}.
\end{remark}

We next present the statement of the coarse-grained Caccioppoli inequality, which is a generalization of the version proved in~\cite[Proposition 2.5]{AK.HC}. The argument closely follows that of~\cite{AK.HC}, with a slight adaptation to cubes overlapping the boundary of~$\cu_0$. Testing the equation with~$u\varphi$ and integrating by parts gives
\begin{equation*}
	\int_{\cu_0}
	\varphi\nabla u\cdot \s\nabla u
	=
	- \int_{\cu_0}
	u \nabla \varphi \cdot \a\nabla u
	\,.
\end{equation*}
In the uniformly elliptic setting, the integral on the right side is  estimated by Cauchy-Schwarz and~$\|\a\|_{L^\infty}\leq \Lambda$, yielding the classical Caccioppoli inequality. To obtain a  coarse-grained version, the inequality is bounded differently, by splitting into small subcubes and then using the Besov duality pairing, leading to factors involving the gradients and fluxes in norms of negative regularity.

\begin{proposition}[Coarse-grained Caccioppoli inequality, boundary version]
\label{p.coarse.grained.Caccioppoli.boundary}
There exists a constant~$C(d)<\infty$ such that,
for every~$s,t\in (0,1)$ with~$s+t < 1$, every~$x \in \cu_0$, and every~$u\in H^1(\cu_0)$ satisfying
\begin{equation}
	\left\{
	\begin{aligned}
		& -\nabla \cdot \a\nabla u = 0
		&  \mbox{in} & \ \cu_0 \,,
		\\ &
		u = 0 &  \mbox{on} & \ (\partial \cu_0) \cap (x + \cu_{-1})\,,
	\end{aligned}
	\right.
\end{equation}
we have the estimate
\begin{equation}
	\label{e.coarse.grained.Caccioppoli.boundary}
	\| \s^{\nf12}  \nabla u \|_{\underline{L}^2((x + \cu_{-2}) \cap \cu_0)}^2
	\leq
	\biggl( \frac{C}{1-s-t} \biggr)^{\! 2 + \frac{4s}{1-s-t}}
	\biggl( \frac{\Lambda_{s}(\cu_0; \a)}{\lambda_{t}(\cu_0; \a)} \biggr)^{\!\frac{s}{1-s-t}}
	{\Lambda}_{s}(\cu_0; \a)
	\| u \|_{\underline{L}^2(\cu_{0})}^2
	\,.
\end{equation}
\end{proposition}
\begin{proof}
Fix~$x \in \cu_0$ and denote~$\sigma \coloneqq  1 - s - t$. Note that~$\sigma >0$ implies~$1-s \geq \sigma + t > t$. Fix parameters~$\rho_1, \rho_2 \in [\nf13,1)$ with~$\rho_1 < \rho_2$ and scale separation parameters~$h,k\in\N$ satisfying
\begin{equation*}
	3^{-4} (\rho_2 - \rho_1) \leq 3^{-k} \leq 3^{-3} (\rho_2 - \rho_1) \qand h\geq k+4
	\,.
\end{equation*}
Take a smooth cutoff function~$\varphi \in C^\infty(\Rd)$ satisfying
\begin{equation*}
	\indc_{x + \rho_1 \cu_{-1}}
	\leq
	\varphi \leq
	\indc_{x + \frac12(\rho_1 + \rho_2) \cu_{-1}}
\end{equation*}
and, for~$j \in \{1,2\}$,
\begin{equation*}
	\| \nabla^j\varphi\|_{L^\infty(\Rd)} \leq C 3^{k j} \leq C (\rho_2 - \rho_1)^{-j}
	\,.
\end{equation*}
Since~$u = 0$ on~$(\partial \cu_0) \cap (x + \cu_{-1})$ and~$\supp \varphi \subset x + \cu_{-1}$, we have~$u\varphi \in H^1_0(\cu_0)$. Testing the equation with~$u \varphi$ yields
\begin{equation}
	\label{e.bdy.Caccioppoli.the.beginning}
	\int_{\cu_0}
	\varphi \nabla u \cdot \a \nabla u
	=
	-\int_{\cu_0}
	u \nabla \varphi \cdot \a\nabla u
	\,.
\end{equation}
Denote the set of subcubes of~$\cu_0$ which overlap with the support of~$\nabla \varphi$ by
\begin{equation*}
	S \coloneqq  \bigl\{ z \in 3^{-h} \Zd \cap \cu_0 \,:\, \supp(\nabla \varphi) \cap (z+\cu_{-h}) \neq \emptyset \bigr\}
	\,.
\end{equation*}
Since the cubes~$\{z + \cu_{-h} : z \in 3^{-h}\Zd \cap \cu_0\}$ partition~$\cu_0$, we have
\begin{equation}
	\label{e.bdy.Caccioppoli.sum.over.S}
	\int_{\cu_0}
	u \nabla \varphi \cdot \a\nabla u
	=
	\avsumcube{z}{-h}{0} \indc_{S}(z)
	\fint_{z+\cu_{-h}}
	u \nabla \varphi \cdot \a\nabla u
	\,.
\end{equation}
Using~\eqref{e.duality.for.cubes}, we estimate, for~$z \in S$,
\begin{align*}
	\biggl| \fint_{z+\cu_{-h}}
	u \nabla \varphi \cdot \a\nabla u \biggr|
	&
	\leq
	C \bigl| (u)_{z+\cu_{-h}}\bigr|
	\| \nabla \varphi \|_{\underline{B}_{2,\infty}^{1}(z+\cu_{-h})}
	[ \a \nabla u ]_{\Besov{-1}{2}{1}(z+\cu_{-h})}
	\notag \\ & \qquad
	+
	C
	\bigl\| (u - (u)_{z+\cu_{-h}})  \nabla \varphi \bigr\|_{\underline{B}^{s}_{2,\infty}(z+\cu_{-h})}
	[ \a \nabla u ]_{\Besov{-s}{2}{1}(z+\cu_{-h})}
	\,.
\end{align*}
By Lemma~\ref{l.divcurl} and Proposition~\ref{p.coarse.grained.Poincare}, we have
\begin{align*}
	\bigl\| (u - (u)_{z+\cu_{-h}})  \nabla \varphi \bigr\|_{\underline{B}^{s}_{2,\infty}(z+\cu_{-h})}
	&
	\leq
	C 3^{k} [\nabla u ]_{\Besov{s-1}{2}{1}(z+\cu_{-h})}
	\notag \\ &
	\leq
	\frac{C}{1-s} 3^{k-(1-s)h} \lambda_{1-s,1}^{-\nf12} (z+\cu_{-h};\a) \| \s^{\nf 12} \nabla u \|_{\underline{L}^2(z+\cu_{-h})}
	\,.
\end{align*}
Again by Proposition~\ref{p.coarse.grained.Poincare}, we get, for~$r \in \{s,1\}$,
\begin{equation*}
	[ \a \nabla u ]_{\Besov{-r}{2}{1}(z+\cu_{-h})}
	\leq
	2 r^{-1}
	3^{-r h}\Lambda^{\nf12 }_{r,1} (z + \cu_{-h}; \a)
	\| \s^{\nf 12} \nabla u \|_{\underline{L}^2(z+\cu_{-h})}
	\,.
\end{equation*}
Since~$\| \nabla \varphi \|_{\underline{B}_{2,\infty}^{1}(z+\cu_{-h})} \leq C3^{k+h}$, we may combine the previous three displays to obtain
\begin{align}
	\label{e.bdy.single.cube.est}
	\biggl| \fint_{z+\cu_{-h}}
	u \nabla \varphi \cdot \a\nabla u \biggr|
	&
	\leq
	C 3^{k} \Lambda^{\nf12 }_{1,1} (z + \cu_{-h}; \a)
	\bigl| (u)_{z+\cu_{-h}}\bigr|
	\| \s^{\nf 12} \nabla u \|_{\underline{L}^2(z+\cu_{-h})}
	\notag \\ & \qquad
	+
	\frac{C 3^{k -h} }{s(1-s)} \biggl( \frac{\Lambda_{s}(z+\cu_{-h};\a)}{\lambda_{1-s,1}(z+\cu_{-h};\a)}
	\biggr)^{\! \nf 12} \| \s^{\nf 12} \nabla u \|_{\underline{L}^2(z+\cu_{-h})}^2
	\,.
\end{align}
For~$z \in S$, we estimate the mean value term by
\begin{equation*}
	\bigl|(u)_{z + \cu_{-h}}\bigr| \leq \|u\|_{\underline{L}^2(z + \cu_{-h})}
	\,.
\end{equation*}
We may therefore bound the first term on the right side of~\eqref{e.bdy.single.cube.est}
for~$z \in S$ by
\begin{equation*}
	C 3^k \Lambda_{1,1}^{\nf12}(z + \cu_{-h}; \a) \|u\|_{\underline{L}^2(z + \cu_{-h})} \|\s^{\nf12} \nabla u\|_{\underline{L}^2(z + \cu_{-h})}\,.
\end{equation*}
We now apply~\eqref{e.ellipticities.monotone.ordered} and~\eqref{e.bound.one.cube.by.lambdas} to get
\begin{equation*}
	\max_{z\in 3^{-h} \Zd \cap \cu_0}
	\Lambda_{1,1}^{\nf12} (z+\cu_{-h};\a)
	\leq
	3^{sh} \Lambda_{s}^{\nf12} (\cu_{0};\a)
	\,.
\end{equation*}
For the second term in~\eqref{e.bdy.single.cube.est}, using~$\lambda_{1-s,1}(z+\cu_{-h}) \geq \lambda_t(z+\cu_{-h})$ (since~$1-s -t \geq \sigma > 0$) and then scaling, we obtain
\begin{equation*}
	3^{-h}
	\max_{z\in 3^{-h} \Zd \cap \cu_0}
	\biggl(
	\frac{\Lambda_{s}(z+\cu_{-h};\a)}{\lambda_{1-s,1}(z+\cu_{-h};\a) }
	\biggr)^{\!\nf12}
	\leq
	3^{(s+t-1)h}
	\biggl(
	\frac{\Lambda_{s}(\cu_{0};\a)}{\lambda_{t}(\cu_{0};\a) }
	\biggr)^{\!\nf12}
	\leq
	3^{-\sigma h}
	\biggl(
	\frac{\Lambda_{s}(\cu_{0};\a)}{\lambda_{t}(\cu_{0};\a) }
	\biggr)^{\!\nf12}
	\,.
\end{equation*}
Summing over~$z \in S$ and combining with~\eqref{e.bdy.Caccioppoli.the.beginning} and~\eqref{e.bdy.Caccioppoli.sum.over.S}, we obtain
\begin{align}
	\label{e.bdy.Caccioppoli.almost.done}
	\| \s^{\nf12}  \nabla u \|_{\underline{L}^2((x + \rho_1\cu_{-1}) \cap \cu_0)}^2
	&
	\leq
	\frac{C3^{k}}{s(1-s)}
	3^{-\sigma h}
	\biggl(
	\frac{\Lambda_{s}(\cu_{0};\a)}{\lambda_{t}(\cu_{0};\a) }
	\biggr)^{\!\nf12}
	\| \s^{\nf12}  \nabla u \|_{\underline{L}^2((x + \rho_2\cu_{-1}) \cap \cu_0)}^2
	\notag \\ & \qquad
	+
	C3^{k} 3^{sh}
	\Lambda_{s}^{\nf12}  (\cu_{0};\a)
	\| u \|_{\underline{L}^2(\cu_0)}
	\| \s^{\nf12}  \nabla u \|_{\underline{L}^2((x + \rho_2 \cu_{-1}) \cap \cu_0)}
	\,.
\end{align}
We now select
\begin{equation*}
	h  \coloneqq  \Biggl\lceil \sigma^{-1} \log_3 \biggl( \frac{4C_{\eqref{e.bdy.Caccioppoli.almost.done}}3^{k}}{s(1-s)}
	\biggl(\frac{\Lambda_{s}(\cu_{0};\a)}{\lambda_{t}(\cu_{0};\a) }
	\biggr)^{\!\nf12} \biggr) \Biggr\rceil
	\,.
\end{equation*}
By Young's inequality, we obtain
\begin{align*}
	\| \s^{\nf12}  \nabla u \|_{\underline{L}^2((x + \rho_1\cu_{-1}) \cap \cu_0)}^2
	&
	\leq
	\frac12 \| \s^{\nf12}  \nabla u \|_{\underline{L}^2((x + \rho_2 \cu_{-1}) \cap \cu_0)}^2
	\notag \\ & \quad
	+
	C (\rho_2-\rho_1)^{- \frac{2(s + \sigma)}{\sigma}}
	\biggl(  \frac{C}{s(1-s)} \biggl(\frac{\Lambda_{s}(\cu_{0};\a)}{\lambda_{t}(\cu_{0};\a) }
	\biggr)^{\!\nf12}   \biggr)^{\! \nf{2s}{\sigma}} \Lambda_{s}  (\cu_{0};\a)\| u \|_{\underline{L}^2(\cu_0)}^2
	\,.
\end{align*}
A standard iteration argument (\cite[Lemma C.6]{AKMBook}) then yields
\begin{equation*}
	\| \s^{\nf12}  \nabla u \|_{\underline{L}^2((x + \cu_{-2}) \cap \cu_0)}^2
	\leq
	\bigl(C \sigma^{-1} \bigr)^{\frac{2(s + \sigma)}{\sigma}}
	(s(1-s))^{-\nf {2s}\sigma}
	\biggl(\frac{\Lambda_{s}(\cu_{0};\a)}{\lambda_{t}(\cu_{0};\a) }
	\biggr)^{\! \nf{s}{\sigma}}
	\Lambda_{s}  (\cu_{0};\a)\| u \|_{\underline{L}^2(\cu_0)}^2
	\,.
\end{equation*}
Simplifying by using~$s^{-s} \leq C$ and~$(1-s)^{-1} \leq \sigma^{-1}$ and inserting the result into the previous display, we obtain~\eqref{e.coarse.grained.Caccioppoli.boundary}.
\end{proof}

The next statement is analogous to Proposition~\ref{p.coarse.grained.Poincare}, but is an estimate for the weak norm of the error in the linear response, rather than that of the gradient and flux of a solution. It is very close to~\cite[Lemma 5.5]{AK.HC}.

\begin{lemma}[Coarse-graining an elliptic operator]
\label{l.coarse.graining.operator}
For every~$m\in\Z$,~$s\in (0,1)$, symmetric positive-definite matrix~$\a_0\in\Rsymp$, and~$u\in\A(\cu_m; \a)$, we have
\begin{equation}
	\label{e.CG.elliptic.operator}
	3^{-s m} \| ( \a-\a_0) \nabla u \|_{\Besov{-s}{2}{1}(\cu_m)}
	\leq
	4s^{-1}
	|\a_0|^{\nf12}
	\|\s^{\nf12} \nabla u \|_{\underline{L}^2(\cu_m)}
	\mathcal{E}_ {s,\infty,1} (\cu_m;\a,\a_0)
	\,.
\end{equation}
\end{lemma}
\begin{proof}
According to~\eqref{e.fluxmaps} with~$p = \a_0^{-\nf 12} e$,~$q = \a_0^{\nf12} e$ and~$e \in \Rd$,
\begin{align*}
	\lefteqn{
	3^{-sm} \| ( \a-\a_0) \nabla u \|_{\Besov{-s}{2}{1}(\cu_m)}
	} \qquad &
	\notag \\ &
	=
	\sum_{k=-\infty}^m  3^{s (k-m)}
	\biggl (
	\avsum_{z \in 3^{k} \Zd \cap \cu_m}  \bigl|\bigl(( \a-\a_0) \nabla u \bigr)_{z+\cu_k} \bigr|^2
	\biggr)^{\!\nf12}
	\notag \\ &
	\leq
	2^{\nf12}
	|\a_0|^{\nf12}
	\sum_{k=-\infty}^m
	3^{s (k-m)}
	\biggl(
	\avsum_{z \in 3^{k} \Zd \cap \cu_m}
	\| \s^{\nf12} \nabla u \|_{\underline{L}^2(z+\cu_k)}^2
	\max_{|e| \leq 1}
	J(z+\cu_k, \a_0^{-\nf12} e, \a^{\nf12}_0 e) \biggr)^{\!\nf12}
	\notag \\ &
	\leq
	2^{\nf12}
	\css{s}^{-1}
	|\a_0|^{\nf12}
	\|\s^{\nf12} \nabla u \|_{\underline{L}^2(\cu_m)}
	\css{s}
	\sum_{k=-\infty}^m  3^{s (k-m)}
	\max_{z \in 3^{k} \Zd \cap \cu_m}
	\max_{|e| \leq 1}
	J(z+\cu_k, \a_0^{-\nf12} e, \a^{\nf12}_0 e)^{\nf12}
	\notag \\ &
	\leq
	4s^{-1}  |\a_0|^{\nf12}
	\|\s^{\nf12} \nabla u \|_{\underline{L}^2(\cu_m)}
	\mathcal{E}_{s,\infty,1} (\cu_m;\a,\a_0) \,\,.
\end{align*}
This completes the proof.
\end{proof}

\subsection{Coarse-graining solutions of inhomogeneous equations}
\label{ss.cg.RHS}

In this subsection, we generalize Propositions~\ref{p.coarse.grained.Poincare},~\ref{p.coarse.grained.Caccioppoli.boundary} and Lemma~\ref{l.coarse.graining.operator} to solutions~$u$ of the inhomogeneous equation
\begin{equation*}
	-\nabla \cdot \a \nabla u
	= \nabla \cdot \g
	\,,
\end{equation*}
where the vector field~$\g$ belongs to a space of positive regularity such as~$H^s=B^{s}_{2,2}$ for~$s>0$.
Most of the coarse-graining theory developed in~\cite{AK.HC} is for solutions of a homogeneous equation (zero right-hand side). For equations with symmetric coefficients, the generalization to non-zero right-hand side is rather straightforward. The case of equations with nonsymmetric coefficient fields is more subtle. Most of the results in this subsection are therefore new and of independent interest.

\smallskip

We begin with the generalization of Proposition~\ref{p.coarse.grained.Poincare}.

\begin{lemma}[Coarse-grained Poincar\'e inequality with RHS]
\label{l.coarse.graining.RHS}
There exists~$C(d)<\infty$ such that, for every~$m\in\Z$, every coefficient field~$\a \in L^{\infty}(\cu_m ;\R_{+}^{d \times d})$, every~$s\in (0,1)$,~$\g \in H^s(\cu_m;\Rd)$ and solution~$u \in H^1(\cu_m)$ of the equation
\begin{equation*}
	-\nabla \cdot \a \nabla u = \nabla \cdot \g \quad \mbox{in} \ \cu_m\,,
\end{equation*}
we have the estimates
\begin{equation}
	\label{e.cg.Poincare.with.rhs.grad}
	3^{-sm}\| \nabla u \|_{\Besov{-s}{2}{2}(\cu_m)}
	\leq
	C s^{-\nf32}  \lambda_{\nf s2,2}^{-\nf12}  (\cu_m;\a)
	\| \s^{\nf12} \nabla u\|_{\underline{L}^2(\cu_m)}
	+
	C s^{-3}
	\lambda_{\nf s2,2}^{-1} (\cu_m;\a)
	3^{sm} [ \g ]_{\underline{H}^s(\cu_m) }  \,\,,
\end{equation}
\begin{equation}
	\label{e.cg.Poincare.with.rhs.flux}
	3^{-sm} \| \a \nabla u \|_{\Besov{-s}{2}{2}(\cu_m)}
	\leq
	C s^{-\nf32} \Lambda_{\nf s2,2}^{\nf12}  (\cu_m;\a)
	\| \s^{\nf12} \nabla u\|_{\underline{L}^2(\cu_m)}
	+
	C s^{-\nf 92}
	\frac{\Lambda_{\nf s 2,2}^{\nf12}  (\cu_m;\a)  }{\lambda_{\nf s 2,2}^{\nf12}  (\cu_m;\a)  }
	3^{sm} [ \g ]_{\underline{H}^s(\cu_m) }  \,\,.
\end{equation}
Moreover, for every~$h\in H^{1+s}(\cu_m)$, the solutions~$v,w \in H^1(\cu_m)$  of the Dirichlet and Neumann problems
\begin{equation}
	\label{e.coarse.graining.RHS.eq}
	\left\{
	\begin{aligned}
		& -\nabla \cdot \a\nabla v = \nabla \cdot \g
		& \mbox{in} & \ \cu_m\,, \\
		& v = h & \mbox{on} & \ \partial \cu_m\,,
	\end{aligned}
	\right.
	\qand
	\left\{
	\begin{aligned}
		& -\nabla \cdot \a\nabla w = \nabla \cdot \g
		& \mbox{in} & \ \cu_m\,, \\
		& \mathbf{n} \cdot (\a\nabla w +\g - (\g)_{\cu_m}) = 0 & \mbox{on} & \ \partial \cu_m\,.
	\end{aligned}
	\right.
\end{equation}
satisfy the estimate
\begin{align}
	\label{e.cg.RHS}
	\lefteqn{
	\| \s^{\nf12} \nabla v\|_{\underline{L}^2(\cu_m)}
	+
	\| \s^{\nf12} \nabla w\|_{\underline{L}^2(\cu_m)}
	} \qquad &
	\notag \\ &
	\leq
	Cs^{-3} \lambda_{\nf s 2,2}^{-\nf12} (\cu_m\,;\a)
	3^{sm} [ \g ]_{\underline{H}^s(\cu_m)}
	+ C s^{-\nf32} \Lambda_{\nf s2,2}^{\nf12}  (\cu_m;\a)
	3^{sm}
	\| \nabla h \|_{\underline{H}^s(\cu_m)}
	\,.
\end{align}
\end{lemma}
\begin{proof}
Fix~$m \in \Z$,~$\g \in H^s(\cu_m)$ with~$(\g)_{\cu_m}=0$ and~$u \in H^1(\cu_m)$ which satisfy
\begin{equation*}
	- \nabla \cdot \a\nabla u = \nabla \cdot \g  \quad \mbox{in} \ \cu_m\,\,.
\end{equation*}
Define for each~$n\in (-\infty,m]\cap\Z$ the quantities
\begin{equation*}
	R_{n}  \coloneqq
	\avsum_{z \in 3^{n} \Zd \cap \cu_m}
	3^{-2 sn} \| \nabla u \|_{\Besov{-s}{2}{2}(z+\cu_n)}^2
	\quad
	\qand
	\quad
	D_n
	\coloneqq
	\avsum_{z \in 3^{n} \Zd \cap \cu_m}
	3^{-2 sn} \| \a \nabla u \|_{\Besov{-s}{2}{2}(z+\cu_n)}^2 \,.
\end{equation*}
In Step 1 below, we will show that, for every~$n \in \Z \cap (-\infty,m]$,
\begin{equation}
	R_{n} \leq 3^{-\frac 32 s}  R_{n-1}+ C s^{-2} \lambda_{\nf s2,2}^{-1}(\cu_m)
	3^{s(m-n)}
	\| \s^{\nf12} \nabla u \|_{\underline{L}^2(\cu_m)}^2
	+
	C s^{-5} \lambda_{\nf s2,2}^{-2}(\cu_m) 3^{2sm} \| \g \|_{\underline{B}^s_{2,2}(\cu_m)}^2
	\,\,.
	\label{e.Rn.recurrence}
\end{equation}
An iteration of this inequality, using the crude bound~$\limsup_{n \to -\infty} 3^{ sn } R_n \leq \| \nabla u \|_{L^2(\cu_m)}^2 < \infty$, yields~\eqref{e.cg.Poincare.with.rhs.grad}. Using this, in Step 2, we establish the energy bound~\eqref{e.cg.RHS},
which is then used, in Step 3, to show
\begin{equation}
	D_{n} \leq 3^{-\frac 32 s}  D_{n-1}+ Cs^{-2} \Lambda_{\nf s2,2}(\cu_m)
	3^{s(m-n)}
	\| \s^{\nf12} \nabla u \|_{\underline{L}^2(\cu_m)}^2
	+
	C s^{-8} \frac{\Lambda_{\nf s2,2}(\cu_m)}{\lambda_{\nf s2,2}(\cu_m)} 3^{2sm} \| \g \|_{\underline{B}^s_{2,2}(\cu_m)}^2
	\,\,.
	\label{e.Dn.recurrence}
\end{equation}
Since~$\limsup_{n \to -\infty} 3^{ sn } D_n \leq \| \a \nabla u \|_{L^2(\cu_m)}^2 < \infty$, this yields~\eqref{e.cg.Poincare.with.rhs.flux}.
Note that, since~$(\g)_{\cu_m}=0$, we have that~$\| \g \|_{\underline{H}^s(\cu_m)} \leq C[  \g ]_{\underline{H}^s(\cu_m)}$. Also, we have~$H^s = B^s_{2,2}$ and the norms~$\|\cdot\|_{\underline{H}^s}$ and~$\|\cdot\|_{\underline{B}^s_{2,2}}$ are equivalent, so we use these interchangeably.

\smallskip

\emph{Step 1.}
We prove~\eqref{e.Rn.recurrence} and thus~\eqref{e.cg.Poincare.with.rhs.grad}.
What we show is that, for every~$n \in \Z \cap (-\infty,m]$ and~$z \in 3^{n} \Zd \cap \cu_m$,
\begin{align}
	\label{e.u.besov.norm.RHS.pre}
	\lefteqn{
	3^{-2sn} \| \nabla u \|_{\Besov{-s}{2}{2} (z+\cu_n)}^2
	} \qquad &
	\notag \\ &
	\leq
	3^{- \frac 32 s}
	\avsum_{y \in 3^{n-1} \Zd \cap (z+\cu_n)}
	3^{-2s(n-1)} \| \nabla u \|_{\Besov{-s}{2}{2} (y+\cu_{n-1})}^2
	+ C s^{-2} \lambda_{s,2}^{-1}(z+\cu_n) \| \s^{\nf12} \nabla u \|_{\underline{L}^2(z+\cu_n)}^2
	\notag \\ & \qquad
	+
	C s^{-5} \lambda_{s,2}^{-2}(z+\cu_n) 3^{2sn} [ \g ]_{\underline{B}^s_{2,2}(z+\cu_n)}^2
	\,.
\end{align}
Using~\eqref{e.bound.one.cube.by.lambdas} and~\eqref{e.ellipticities.monotone.ordered}, we have that, for every~$z \in 3^{n} \Zd \cap \cu_m$,
\begin{equation*}
	\lambda_{s,2}^{-1}(z+\cu_n)
	\leq \lambda_{\nf s2, 2}^{-1}(z+\cu_n) \leq 3^{s(m-n)} \lambda^{-1}_{\nf s2,2}(\cu_m)
	\,,
\end{equation*}
and upon substituting this into~\eqref{e.u.besov.norm.RHS.pre} and then summing over~$z\in 3^n\Zd\cap \cu_m$, we obtain~\eqref{e.Rn.recurrence}.

\smallskip

For ease of notation, we consider the case~$z = 0$.
Fix~$\eps \in (0, \nf12]$ to be determined below and
let~$v_n \in H^1_0(\cu_n)$ be the solution of the Dirichlet problem
\begin{equation*}
	\left\{
	\begin{aligned}
		& -\nabla \cdot \a\nabla v_n = \nabla \cdot \g
		& \mbox{in} & \ \cu_n\,, \\
		& v_n = 0 & \mbox{on} & \ \partial \cu_n\,.
	\end{aligned}
	\right.
\end{equation*}
Since~$u-v_n$ is a solution of the equation with zero right side, we may apply~\eqref{e.besov.grad.poincare} to get
\begin{align}
	\| \nabla u - \nabla v_n \|_{\Besov{-s}{2}{2}(\cu_n)}^2
	&
	\leq
	\css{2s}^{-1}
	\lambda_{s,2}^{-1}(\cu_n)  3^{2 sn}
	\| \s^{\nf12} ( \nabla u - \nabla v_n)  \|_{\underline{L}^2(\cu_n)}^2
	\notag \\ &
	\leq
	3 s^{-1} \lambda_{s,2}^{-1}(\cu_n)  3^{2 sn}
	\bigl(
	\| \s^{\nf12} \nabla u \|_{\underline{L}^2(\cu_n)}^2
	+
	\| \s^{\nf12} \nabla v_n \|_{\underline{L}^2(\cu_n)}^2
	\bigr)
	\,.
	\label{e.weaknormbound.inproof}
\end{align}
Testing the equation for~$v_n$ with~$v_n$ and using~\eqref{e.duality.for.cubes}, we get
\begin{align*}
	\| \s^{\nf12} \nabla v_n \|_{\underline{L}^2(\cu_n)}^2
	=
	\fint_{\cu_n} \nabla v_n  \cdot \a \nabla v_n
	&
	=
	- \fint_{\cu_n} ( \g - (\g)_{z+\cu_n} ) \cdot \nabla v_n
	\notag \\ &
	\leq
	3^{d+s}
	[ \g ]_{\underline{B}^s_{2,2}(\cu_n)}
	\| \nabla v_n \|_{\Besov{-s}{2}{2} (\cu_n)}
	\notag \\ &
	\leq
	3^{d+s}
	[ \g ]_{\underline{B}^s_{2,2}(\cu_n)}
	\bigl(
	\| \nabla u \|_{\Besov{-s}{2}{2} (\cu_n)}
	+
	\| \nabla u - \nabla v_n \|_{\Besov{-s}{2}{2} (\cu_n)}
	\bigr)
	\,.
\end{align*}
Combining these and using Young's inequality, we get
\begin{align}
	[ \nabla u - \nabla v_n ]_{\Besov{-s}{2}{2}(\cu_n)}^2
	&
	\leq
	C s^{-1} \lambda_{s,2}^{-1}(\cu_n)
	3^{2 sn}
	\| \s^{\nf12} \nabla u \|_{\underline{L}^2(\cu_n)}^2
	+
	C s^{-1} \lambda_{s,2}^{-1}(\cu_n)
	3^{2 sn}
	[ \g ]_{\underline{B}^s_{2,2}(\cu_n)}
	\| \nabla u \|_{\Besov{-s}{2}{2} (\cu_n)}
	\notag \\ & \qquad
	+
	C s^{-2} \lambda_{s,2}^{-2}(\cu_n)
	3^{4sn}
	[ \g ]_{\underline{B}^s_{2,2}(\cu_n)}^2
	\,.
	\label{e.weaknormbound.inproof.two}
\end{align}
Using that~$(\nabla v_n)_{\cu_n} = 0$, we compute
\begin{align*}
	\| \nabla v_n \|_{\Besov{-s}{2}{2} (\cu_n)}^2
	&
	=
	\sum_{k=-\infty}^n 3^{2s k}
	\avsum_{z \in 3^{k} \Zd \cap \cu_n} \bigl| (\nabla v_n)_{z+\cu_k} \bigr|^{2}
	\notag \\ &
	=
	\sum_{k=-\infty}^{n-1} 3^{2s k}
	\avsum_{z \in 3^{k} \Zd \cap \cu_n} \bigl| (\nabla v_n)_{z+\cu_k} \bigr|^{2}
	=
	\avsum_{y \in 3^{n-1} \Zd \cap \cu_n}
	\| \nabla v_n \|_{\Besov{-s}{2}{2} (y+\cu_{n-1})}^2
	\,.
\end{align*}
By the triangle inequality, for every~$\ep >0$,
\begin{equation}
	\label{e.u.besov.norm.discounted}
	\| \nabla u \|_{\Besov{-s}{2}{2} (\cu_n)}^2
	\leq
	(1+\ep)
	\avsum_{y \in 3^{n-1} \Zd \cap \cu_n}
	\| \nabla u \|_{\Besov{-s}{2}{2} (y+\cu_{n-1})}^2
	+
	(1+\eps^{-1})
	\| \nabla u-\nabla v_n \|_{\Besov{-s}{2}{2} (\cu_n)}^2
	\,.
\end{equation}
Combining~\eqref{e.weaknormbound.inproof.two} and~\eqref{e.u.besov.norm.discounted} and using Young's inequality again, we obtain, for every~$\ep \in (0,\nf12]$,
\begin{align*}
	\| \nabla u \|_{\Besov{-s}{2}{2} (\cu_n)}^2
	&
	\leq
	(1+\ep)
	\avsum_{y \in 3^{n-1} \Zd \cap \cu_n}
	\| \nabla u \|_{\Besov{-s}{2}{2} (y+\cu_{n-1})}^2
	+
	C s^{-1}
	(1+\eps^{-1})
	\lambda_{s,2}^{-1}(\cu_n)
	3^{2 sn}
	\| \s^{\nf12} \nabla u \|_{\underline{L}^2(\cu_n)}^2
	\notag \\ & \qquad
	+
	C s^{-2} \ep^{-1}
	(1+\eps^{-1})^2
	\lambda_{s,2}^{-2}(\cu_n)
	3^{4sn}
	[ \g ]_{\underline{B}^s_{2,2}(\cu_n)}^2
	+
	\ep \| \nabla u \|_{\Besov{-s}{2}{2} (\cu_{n})}^2
	\,.
\end{align*}
We rearrange this inequality to obtain
\begin{align*}
	\| \nabla u \|_{\Besov{-s}{2}{2} (\cu_n)}^2
	&
	\leq
	\frac{1+\ep}{1-\ep}
	\avsum_{y \in 3^{n-1} \Zd \cap \cu_n}
	\| \nabla u \|_{\Besov{-s}{2}{2} (y+\cu_{n-1})}^2
	\notag \\ & \qquad
	+ C s^{-1} \eps^{-1} \lambda_{s,2}^{-1}(\cu_n)  3^{2sn}\| \s^{\nf12} \nabla u \|_{\underline{L}^2(\cu_n)}^2
	+
	C s^{-2} \eps^{-3} \lambda_{s,2}^{-2}(\cu_n)
	3^{4sn} [ \g ]_{\underline{B}^s_{2,2}(\cu_n)}^2
	\,.
\end{align*}
Selecting~$\eps \coloneqq 2^{-2} s$ and then dividing by~$3^{2 sn}$, we get~\eqref{e.u.besov.norm.RHS.pre}.

\smallskip

\emph{Step 2.} We prove~\eqref{e.cg.RHS} in the case~$h=0$. Test~\eqref{e.coarse.graining.RHS.eq} with~$v \in H^1_0(\cu_m)$, and use~\eqref{e.duality.for.cubes}, Young's inequality, and the bound~\eqref{e.cg.Poincare.with.rhs.grad} with~$v$ in place of~$u$ to obtain, for every~$\eps \in (0,1)$,
\begin{align}
	\fint_{\cu_m}
	\nabla v \cdot \s\nabla v
	&
	=
	- \fint_{\cu_m} \g \cdot \nabla v
	\notag \\ &
	\leq
	C
	\| \nabla v \|_{\Besov{-s}{2}{2}(\cu_m)}
	\| \g \|_{\underline{B}^s_{2,2}(\cu_m)}
	\notag \\ &
	\leq
	\eps
	3^{-2 s m}
	\lambda_{\nf s 2,2}(\cu_m)
	[ \nabla v ]^2_{\Besov{-s}{2}{2}(\cu_m)}
	+
	C
	\eps^{-1}
	3^{2 s m}
	\lambda_{\nf s 2,2}^{-1}(\cu_m)
	\| \g \|_{\underline{B}^s_{2,2}(\cu_m)}^2
	\notag \\ &
	\leq
	C
	\eps
	s^{-3}
	\|\s^{\nf12} \nabla v\|_{\underline{L}^2(\cu_m)}^2
	+
	C (\eps^{-1} + \eps s^{-6})
	3^{2 s m}
	\lambda_{\nf s 2,2}^{-1}(\cu_m)
	\| \g \|_{\underline{B}^s_{2,2}(\cu_m)}^2
	\,.
	\label{e.cs.rhs.nabla.v.energy.0}
\end{align}
Taking~$\eps \coloneqq  s^{3} (2 C_{\eqref{e.cs.rhs.nabla.v.energy.0}})^{-1}$ small enough and absorbing the first term on the right side into the left-hand side, we obtain~\eqref{e.cg.RHS} in the case~$h=0$.
We obtain the bound for~$\|\s^{\nf12} \nabla w\|_{\underline{L}^2(\cu_m)}$ by a similar argument.

\smallskip

\emph{Step 3.} We prove~\eqref{e.Dn.recurrence} and thus~\eqref{e.cg.Poincare.with.rhs.flux}. We assert, for every~$n \in \Z \cap (-\infty,m]$ and~$z \in 3^{n} \Zd \cap \cu_m$,
\begin{align}
	\lefteqn{
	3^{-2sn} [ \a \nabla u ]_{\Besov{-s}{2}{2} (z+\cu_n)}^2
	} \qquad &
	\notag \\ &
	\leq
	3^{- \frac 32 s}
	\avsum_{y \in 3^{n-1} \Zd \cap (z+\cu_n)}
	3^{-2s(n-1)} [ \a \nabla u ]_{\Besov{-s}{2}{2} (y+\cu_{n-1})}^2
	+ C s^{-2} \Lambda_{s,2}(z+\cu_n)
	\| \s^{\nf12} \nabla u \|_{\underline{L}^2(z+\cu_n)}^2
	\notag \\ &\qquad
	+
	C s^{-8} \Lambda_{s,2}(z+\cu_n) \lambda_{\nf s2,2}^{-1}(z+\cu_n) 3^{2sn}
	[ \g ]_{\underline{B}^s_{2,2}(z+\cu_n)}^2
	\label{e.u.besov.norm.RHS.pre.flux}
	\,.
\end{align}
Using~\eqref{e.bound.one.cube.by.lambdas} and~\eqref{e.ellipticities.monotone.ordered} we have that, for every~$z \in 3^{n} \Zd \cap \cu_m$,
\begin{equation*}
	\lambda_{\nf s2,2}^{-1}(z+\cu_n) \leq 3^{s (m-n)} \lambda^{-1}_{\nf s2,2}(\cu_m)
	\qand
	\Lambda_{s,2}(z+\cu_n) \leq 3^{s (m-n)} \Lambda_{\nf s2,2}(\cu_m)
	\,.
\end{equation*}
Therefore,~\eqref{e.u.besov.norm.RHS.pre.flux} implies~\eqref{e.Dn.recurrence}.

\smallskip

For ease of notation, set~$z = 0$
and assume that~$(\g)_{\cu_n} = 0$.
Fix~$\eps \in (0, \nf12]$ to be determined below and
let~$w_n \in H^1(\cu_n)$ be the solution of the Neumann problem
\begin{equation*}
	\left\{
	\begin{aligned}
		& -\nabla \cdot \a\nabla w_n = \nabla \cdot \g
		& \mbox{in} & \ \cu_n\,, \\
		&  \mathbf{n} \cdot (\a\nabla w_n +\g) = 0 & \mbox{on} & \ \partial \cu_n\,.
	\end{aligned}
	\right.
\end{equation*}
Since~$u-w_n$ is a solution with zero right side, we may apply~\eqref{e.besov.flux.poincare}
and~\eqref{e.cg.RHS} to obtain
\begin{align*}
	3^{-2 s n} \Lambda^{-1}_{s,2}(\cu_n)
	\| \a \nabla u - \a \nabla w_n \|_{\Besov{-s}{2}{2}(\cu_n)}^2
	&\leq
	Cs^{-1}
	\| \s^{\nf12} ( \nabla u - \nabla w_n)  \|_{\underline{L}^2(\cu_n)}^2
	\notag \\
	&\leq
	Cs^{-1}
	\| \s^{\nf12} \nabla u \|_{\underline{L}^2(\cu_n)}^2
	+
	C s^{-7}
	\lambda_{\nf s2,2}^{-1} (\cu_n)
	3^{2 sn}
	\| \g \|_{\underline{B}^s_{2,2}(\cu_n)}^2
	\,.
\end{align*}
Note that~$(\g)_{\cu_n}=0$ implies~$\| \g \|_{\underline{B}^s_{2,2}(\cu_n)}^2\leq C [ \g ]_{\underline{B}^s_{2,2}(\cu_n)}^2$.
Using~$(\a\nabla w_n)_{\cu_n} = (\a\nabla w_n+\g)_{\cu_n} = 0$, we obtain, by the triangle inequality and Young's inequality, for every~$\eps \in (0,1)$,
\begin{align*}
	\| \a \nabla u \|_{\Besov{-s}{2}{2} (\cu_n)}^2
	\leq
	(1+\ep)
	\avsum_{y \in 3^{n-1} \Zd \cap \cu_n}
	\| \a \nabla u \|_{\Besov{-s}{2}{2} (y+\cu_{n-1})}^2
	+
	(1 + \eps^{-1})
	\| \a \nabla u- \a\nabla w_n\|_{\Besov{-s}{2}{2} (\cu_n)}^2
	\,\,.
\end{align*}
Selecting~$\eps = c s$ for sufficiently small~$c$ in the above display, dividing by~$3^{2sn}$ and combining this with the previous display yields~\eqref{e.u.besov.norm.RHS.pre.flux}. This completes the proof of~\eqref{e.Dn.recurrence} and thus of~\eqref{e.cg.Poincare.with.rhs.flux}.

\smallskip

\emph{Step 4.} We complete the proof of~\eqref{e.cg.RHS}. By linearity, the case of general~$(\g,h)$ is the sum of the already proved case with datum~$(\g,0)$ and the remaining case~$(0,h)$, so it suffices to assume~$\g=0$.
Testing the equation for~$v$ in~\eqref{e.coarse.graining.RHS.eq} with~$v-h$  and then applying~\eqref{e.cg.Poincare.with.rhs.flux} yields the estimate
\begin{align*}
	\fint_{\cu_m}
	\nabla v \cdot \s \nabla v
	=
	\fint_{\cu_m} \nabla h \cdot \a \nabla v
	&
	\leq
	\| \nabla h \|_{\underline{H}^s(\cu_m)}
	\| \a\nabla v \|_{\underline{H}^{-s}(\cu_m)}
	\notag \\ &
	\leq
	C
	3^{sm}
	\| \nabla h \|_{\underline{H}^s(\cu_m)}
	s^{-\nf32} \Lambda_{\nf s2,2}^{\nf12}  (\cu_m;\a)
	\| \s^{\nf12} \nabla v\|_{\underline{L}^2(\cu_m)}
	\notag \\ &
	\leq
	\frac12
	\| \s^{\nf12} \nabla v\|_{\underline{L}^2(\cu_m)}^2
	+
	C
	s^{-3} \Lambda_{\nf s2,2}  (\cu_m;\a)
	3^{2sm}
	\| \nabla h \|_{\underline{H}^s(\cu_m)}^2
	\,.
\end{align*}
This completes the proof of~\eqref{e.cg.RHS}.
\end{proof}

We generalize Proposition~\ref{p.coarse.grained.Caccioppoli.boundary} to solutions of inhomogeneous equations with general right-hand sides and boundary data.

\begin{lemma}[Coarse-grained Caccioppoli with RHS, boundary version]
\label{l.coarse.grained.Caccioppoli.RHS}
There exists a constant~$C(d)<\infty$ such that,
for every~$s\in (0,1)$ and~$t \in (0,\nf 12)$ with~$\sigma \coloneqq  1-s-t >0$, every~$x \in \cu_0$, every
coefficient field~$\a \in L^{\infty}(\cu_0 ;\R_{+}^{d \times d})$, and every~$u \in H^1(\cu_0)$,~$\g \in H^{2t}(\cu_0;\Rd)$, and~$h \in H^{1+2t}(\cu_0)$ with~$(h)_{\cu_0}=0$  satisfying
\begin{equation*}
	\left\{
	\begin{aligned}
		& -\nabla \cdot \a\nabla u = \nabla \cdot \g
		&  \mbox{in} & \ \cu_0 \,,
		\\ &
		u = h &  \mbox{on} & \ (\partial \cu_0) \cap (x + \cu_{-1})\,,
	\end{aligned}
	\right.
\end{equation*}
we have the estimate
\begin{align}
	\lefteqn{
	\| \s^{\nf12} \nabla u \|_{\underline{L}^2((x + \cu_{-2}) \cap \cu_0)}^2
	} \ \ &
	\notag \\ &
	\leq
	\biggl( \frac{C}{\sigma} \biggr)^{\! 2 + \frac{4s}{\sigma}}
	\biggl( \frac{\Lambda_{s}(\cu_0; \a)}{\lambda_{t}(\cu_0; \a)} \biggr)^{\!\frac{s+\sigma}{\sigma}}
	\notag \\ & \qquad \times
	\Bigl(
	{\lambda}_{t}(\cu_0; \a)
	\| u \|_{\underline{L}^2(\cu_{0})}^2
	\notag \\ & \qquad\quad +
	\bigl(1-2t\bigr)^{-1}
	\Bigl(
	t^{-11}
	\lambda_{t}^{-1} (\cu_0\,;\a)
	\| \g \|_{\underline{H}^{2t} (\cu_0)}^2
	+
	t^{-6}
	\Lambda_t(\cu_0; \a)
	\| \nabla h \|_{\underline{H}^{2t}(\cu_0)}^2
	\Bigr)
	\Bigr)
	\,.
	\label{e.cg.Caccioppoli.with.RHS}
\end{align}
The same two-scale estimate holds, by rescaling, for a triadic parent and its child patch.
\end{lemma}
\begin{proof}
Let~$v$ be the solution of the Dirichlet problem
\begin{equation*}
	\left\{
	\begin{aligned}
		& -\nabla \cdot \a\nabla v = \nabla \cdot \g
		& \mbox{in} & \ \cu_0\,, \\
		& v = h & \mbox{on} & \ \partial \cu_0\,.
	\end{aligned}
	\right.
\end{equation*}
We apply Proposition~\ref{p.coarse.grained.Caccioppoli.boundary} to~$u - v$ to obtain
\begin{align}
	\| \s^{\nf12} ( \nabla u - \nabla v) \|_{\underline{L}^2((x + \cu_{-2}) \cap \cu_0)}^2
	&
	\leq
	\biggl( \frac{C}{\sigma} \biggr)^{\! 2 + \frac{4s}{\sigma}}
	\biggl( \frac{\Lambda_{s}(\cu_0; \a)}{\lambda_{t}(\cu_0; \a)} \biggr)^{\!\frac{s}{\sigma}}
	{\Lambda}_{s}(\cu_0; \a)
	\| u - v  \|_{L^2(\cu_{0})}^2
	\,.
	\label{e.cgC.boundary.applied}
\end{align}
By~\eqref{e.cg.RHS} and~\eqref{e.ellipticities.change.q},
\begin{equation}
	\label{e.nabla.v.est.for.cgC.RHS}
	\| \s^{\nf12} \nabla v\|_{L^2(\cu_0)}
	\leq
	Ct^{-3} \lambda_{t}^{-\nf12} (\cu_0\,;\a)
	\| \g \|_{H^{2t} (\cu_0)}
	+
	Ct^{-\nf32} \Lambda_{t}^{\nf12}  (\cu_0;\a)
	\| \nabla h \|_{H^{2t}(\cu_0)}
	\,.
\end{equation}
Since~$v - h \in H_0^1(\cu_0)$ and~$(h)_{\cu_0} = 0$, we deduce by~\eqref{e.fractional.Sobolev.embedding} and~\eqref{e.fractional.Sobolev.embedding.zero.bndr} that
\begin{equation*}
	\| v  \|_{L^2(\cu_{0})}
	\leq
	\| v - h  \|_{L^2(\cu_{0})} + \| h  \|_{L^2(\cu_{0})}
	\leq
	C \| \nabla v  \|_{\Besov{-1}{2}{1}(\cu_0)}
	+
	C \| \nabla h  \|_{\Besov{-1}{2}{1}(\cu_0)}
	\,.
\end{equation*}
For~$t \in (0,\nf 12)$ we have
\begin{align*}
	\| \nabla v  \|_{\Besov{-1}{2}{1}(\cu_0)}
	&
	=
	\sum_{n=-\infty}^0 3^{n} \biggl(\avsum_{z \in 3^n \Zd \cap \cu_0} \bigl| (\nabla v)_{z+\cu_n} \bigr|^2 \biggr)^{\! \nf 12}
	\notag \\ &
	\leq
	\biggl(
	\sum_{n=-\infty}^0 3^{2(1-2t)n}  \biggr)^{\! \nf 12}
	\biggl( \sum_{n=-\infty}^0 3^{4tn} \avsum_{z \in 3^n \Zd \cap \cu_0} \bigl| (\nabla v)_{z+\cu_n} \bigr|^2 \biggr)^{\! \nf 12}
	\notag \\ &
	\leq C\bigl(1-2t\bigr)^{-\nf12} \| \nabla v  \|_{\Besov{-2t}{2}{2}(\cu_0)}
	\,,
\end{align*}
and~$\| \nabla h  \|_{\Besov{-1}{2}{1}(\cu_0)} \leq C \| \nabla h  \|_{\underline{H}^{2t}(\cu_0)}$. Applying~\eqref{e.cg.Poincare.with.rhs.grad} to~$v$ with exponent~$2t$ and using~\eqref{e.cg.RHS} and~\eqref{e.ellipticities.change.q} gives
\begin{equation*}
	\| v  \|_{\underline{L}^2(\cu_{0})}
	\leq
	C\bigl(1-2t\bigr)^{-\nf12}
	\Bigl(
	t^{-\nf{11}2} \lambda_{t}^{-1}  (\cu_0;\a)
	\| \g \|_{\underline{H}^{2t} (\cu_0) }
	+
	t^{-3} \Lambda_{t}^{\nf12}  (\cu_0;\a) \lambda_t^{-\nf12}(\cu_0;\a)
	\| \nabla h \|_{\underline{H}^{2t}(\cu_0)}
	\Bigr)
	\,.
\end{equation*}
By the triangle inequality and the previous display,
\begin{align*}
	\lefteqn{
	\| u - v  \|_{\underline{L}^2(\cu_{0})}^2
	} \quad & \notag
	\\ & \leq
	2\| u \|_{\underline{L}^2(\cu_{0})}^2
	+
	C\bigl(1-2t\bigr)^{-1}
	\Bigl(
	t^{-11} \lambda_{t}^{-2}  (\cu_0;\a)
	\| \g \|_{\underline{H}^{2t} (\cu_0) }^2
	+
	t^{-6} \Lambda_{t}  (\cu_0;\a) \lambda_t^{-1}(\cu_0;\a)
	\| \nabla h \|_{\underline{H}^{2t}(\cu_0)}^2
	\Bigr)
	\,.
\end{align*}
Substituting this into~\eqref{e.cgC.boundary.applied} and using~$\Lambda_s \geq \lambda_t$, we obtain
\begin{align*}
	\lefteqn{
	\| \s^{\nf12} ( \nabla u - \nabla v) \|_{\underline{L}^2((x + \cu_{-2}) \cap \cu_0)}^2
	} \quad &
	\notag \\ &
	\leq
	\biggl( \frac{C}{\sigma} \biggr)^{\! 2 + \frac{4s}{\sigma}}
	\biggl( \frac{\Lambda_{s}}{\lambda_{t}} \biggr)^{\!\frac{s}{\sigma}}
	{\Lambda}_{s}
	\Bigl(
	\| u \|_{\underline{L}^2(\cu_{0})}^2
	+
	\bigl(1-2t\bigr)^{-1}
	\Bigl(
	t^{-11} \lambda_{t}^{-2}
	\| \g \|_{\underline{H}^{2t}(\cu_{0})}^2
	+
	t^{-6} \Lambda_{t} \lambda_t^{-1}
	\| \nabla h \|_{\underline{H}^{2t}(\cu_{0})}^2
	\Bigr)
	\Bigr)
	\notag \\ &
	=
	\biggl( \frac{C}{\sigma} \biggr)^{\! 2 + \frac{4s}{\sigma}}
	\biggl( \frac{\Lambda_{s}}{\lambda_{t}} \biggr)^{\!\frac{s+\sigma}{\sigma}}
	\Bigl(
	\lambda_t \| u \|_{\underline{L}^2(\cu_{0})}^2
	+
	\bigl(1-2t\bigr)^{-1}
	\Bigl(
	t^{-11} \lambda_{t}^{-1}
	\| \g \|_{\underline{H}^{2t}(\cu_{0})}^2
	+
	t^{-6} \Lambda_{t}
	\| \nabla h \|_{\underline{H}^{2t}(\cu_{0})}^2
	\Bigr)
	\Bigr)
	\,.
\end{align*}
In the last line we used~$(\Lambda_s/\lambda_t)^{s/\sigma}\Lambda_s=(\Lambda_s/\lambda_t)^{(s+\sigma)/\sigma}\lambda_t$. Combining~\eqref{e.nabla.v.est.for.cgC.RHS} with the previous display yields~\eqref{e.cg.Caccioppoli.with.RHS} by the triangle inequality.
\end{proof}

We next generalize Lemma~\ref{l.coarse.graining.operator} to solutions of the equation with nonzero right side.

\begin{lemma}[Coarse-graining an elliptic operator with RHS]
\label{l.coarse.graining.operator.RHS}
There exists a constant~$C(d)<\infty$ such that, for every~$m,n\in\Z$ with~$n\leq m$,
coefficient field~$\a \in L^{\infty}(\cu_m ;\R_{+}^{d \times d})$,~$s\in (0,1)$, scalar~$\s_0\in (0,\infty)$,~$\g \in H^s(\cu_m;\Rd)$ and solution~$u \in H^1(\cu_m)$ of the equation
\begin{equation*}
	-\nabla \cdot \a \nabla u = \nabla \cdot \g \quad \mbox{in} \ \cu_m\,,
\end{equation*}
we have the estimate
\begin{align}
	\label{e.fluxmaps.RHS}
	3^{-s m}
	\| ( \a-\s_0) \nabla u \|_{\underline H^{-s} (\cu_m)}
	&
	\leq
	C s^{-1} \s_0^{\nf12}  \| \s^{\nf12}  \nabla u  \|_{\underline{L}^2(\cu_m)}
	\mathcal{E}_{s, \infty, 1} (\cu_m,n;\a,\s_0)
	\notag \\ & \qquad
	+
	C 3^{-s(m-n)} s^{-\nf 9 2} \bigl(1+ \mathcal{E}^2_{\nf s2, \infty, 2} ( \cu_m,n;\a,\s_0) \bigr)
	3^{sm}
	[ \g ]_{\underline{H}^s(\cu_m)}
	\,.
\end{align}
\end{lemma}
\begin{proof}
Fix~$z\in 3^n\Zd\cap \cu_m$.
Let~$v$ be the solution of the Dirichlet problem in~\eqref{e.coarse.graining.RHS.eq} with~$h=0$ in the domain~$z + \cu_n$. According to the triangle inequality, Lemmas~\ref{l.coarse.graining.operator} and~\ref{l.coarse.graining.RHS} and the fact that~$\| \cdot \|_{\underline{H}^{-s}(\cu)} \leq C\| \cdot \|_{\Besov{-s}{2}{2}(\cu)} \leq C\| \cdot \|_{\Besov{-s}{2}{1}(\cu)}$ (see~\eqref{e.Wsp.vs.Bspp.intro}; note that the seminorm~\eqref{e.negative.Besov.def} contains the scale-weighted cube means, so it controls the full dual norm),
\begin{align*}
	\lefteqn{
	3^{- sn } \| ( \a-\s_0) \nabla u \|_{\underline H^{-s} (z+\cu_n)}
	} \qquad &
	\notag \\ &
	\leq
	3^{- sn} \| ( \a-\s_0) (\nabla u - \nabla v) \|_{\underline H^{-s} (z+\cu_n)}
	+
	3^{-sn }\| \a \nabla v \|_{\underline H^{-s} (z+\cu_n)}
	+
	\s_0 3^{-sn}  \|  \nabla v \|_{\underline H^{-s} (z+\cu_n)}
	\notag \\ &
	\leq
	C s^{-1} \s_0^{\nf12}
	\| \s^{\nf12}  \nabla (u-v)  \|_{\underline{L}^2(z+\cu_n)}
	\mathcal{E}_{s, \infty, 1} (z{+}\cu_n;\a,\s_0)
	\notag \\ & \qquad
	+
	C s^{-\nf 9 2}
	\biggl( \s_0 \lambda^{-1}_{\nf s2, 2}(z{+}\cu_n;\a)  + \frac{\Lambda_{\nf s2, 2}^{\nf12}(z{+}\cu_n;\a)}{\lambda_{\nf s2,2}^{\nf12}(z{+}\cu_n;\a) }\biggr)
	3^{sn} [\g ]_{\underline{H}^s(z+\cu_n) }
	\,.
\end{align*}
Using the triangle inequality and Lemma~\ref{l.coarse.graining.RHS} again, we bound the first term on the right side by
\begin{align*}
	\| \s^{\nf12}  \nabla (u-v)  \|_{\underline{L}^2(z+\cu_n)}
	&
	\leq
	\| \s^{\nf12}  \nabla u  \|_{\underline{L}^2(z+\cu_n)}
	+
	\| \s^{\nf12}  \nabla v  \|_{\underline{L}^2(z+\cu_n)}
	\notag \\ &
	\leq
	\| \s^{\nf12}  \nabla u  \|_{\underline{L}^2(z+\cu_n)}
	+
	C s^{-3} \lambda_{\nf s2,2}^{-\nf12} (z{+}\cu_n\,;\a)
	3^{sn} [ \g ]_{\underline{H}^s(z+\cu_n)}   \,\,.
\end{align*}
For the second and third terms on the right, we use~\eqref{e.bound.Lambdas.by.Es} which implies
\begin{equation*}
	\s_0 \lambda_{\nf s2,2}^{-1} (z{+}\cu_n;\a)   +  \frac{\Lambda_{\nf s2, 2}^{\nf12}(z{+}\cu_n;\a)}{\lambda_{\nf s2,2}^{\nf12}(z{+}\cu_n;\a) }
	\leq
	C(1 + \mathcal{E}^2_{\nf s2, \infty, 2} ( z{+}\cu_n;\a,\s_0)) \,.
\end{equation*}
Putting these together, using also~\eqref{e.compareEqs}, we obtain
\begin{align*}
	3^{-s n} \| ( \a-\s_0) \nabla u \|_{\underline H^{-s} (z+\cu_n)}
	&
	\leq
	C  s^{-1} \s_0^{\nf12} \| \s^{\nf12}  \nabla u  \|_{\underline{L}^2(z+\cu_n)}
	\mathcal{E}_{s, \infty, 1} (z{+}\cu_n;\a,\s_0)
	\notag \\ & \qquad
	+
	C s^{-\nf 9 2} \bigl(1+ \mathcal{E}^2_{\nf s2, \infty, 2} ( z{+}\cu_n;\a,\s_0) \bigr)
	3^{sn}
	[ \g ]_{\underline{H}^s(z+\cu_n)}
	\,.
\end{align*}
By squaring this and summing over~$z\in 3^n\Zd\cap \cu_m$, using also~\eqref{e.mathcal.E.t.infty.one} and the superadditivity~\eqref{e.Hs.superadditivity} and subadditivity~\eqref{e.H.minus.s.subadditivity} of the Sobolev seminorms, we get
\begin{align*}
	3^{-s m} \| ( \a-\s_0) \nabla u \|_{\underline H^{-s} (\cu_m)}
	&
	\leq
	C s^{-1} \s_0^{\nf12}  \| \s^{\nf12}  \nabla u  \|_{\underline{L}^2(\cu_m)}
	\mathcal{E}_{s, \infty, 1} (\cu_m,n;\a,\s_0)
	\notag \\ & \qquad
	+
	C 3^{-s(m-n)} s^{-\nf 9 2}  \bigl(1+ \mathcal{E}^2_{\nf s2, \infty, 2} ( \cu_m,n;\a,\s_0)  \bigr)
	3^{sm}
	[ \g ]_{\underline{H}^s(\cu_m)}
	\,.
\end{align*}
This completes the proof.
\end{proof}

We next extend the previous lemma to include a general integrability exponent~$p\in[2,\infty)$.

\begin{lemma}
\label{l.gen.coarse.grain}
Let~$m,n\in\Z$ with~$n<m$,~$s, s_1, s_2 \in (0,1)$ with~$s_1 \in (0,s)$,~$s_2 \in (s,1)$ and~$p \in [2,\infty)$.
There exists a constant~$C(p,d)<\infty$ such that for every coefficient field~$\a \in L^{\infty}( \cu_m ;\R_{+}^{d \times d})$, scalar~$\s_0\in (0,\infty)$,~$\g \in W^{s_2, p}(\cu_m;\Rd)$ and solution~$u \in H^1(\cu_m)$ of the equation
\begin{equation*}
	-\nabla \cdot \a\nabla u
	=
	\nabla \cdot \g
	\quad \mbox{in} \ \cu_m\,,
\end{equation*}
we have the estimate
\begin{align}
	\label{e.gen.coarse.grain}
	\lefteqn{
	3^{-sn}
	\biggl(
	\avsum_{z\in 3^n\Zd \cap \cu_m}
	\| (\a-\s_0) \nabla u \|_{\Besov{-s}{p}{p}(z+\cu_n)}^p
	\biggr)^{\!\nf1p}
	} \quad &
	\notag \\ &
	\leq
	C s^{-1} \s_0^{\nf12}
	\mathcal{E}_{s_1, \infty, 1} (\cu_m, n;\a,\s_0)
	\biggl(
	\sum_{k=-\infty}^n
	3^{-(s-s_1) p (n-k)}
	\avsum_{z\in 3^k\Zd\cap\cu_m}
	\| \s^{\nf12} \nabla u \|_{\underline{L}^2(z+\cu_k)}^{p}
	\biggr)^{\!\nf 1p}
	\notag \\ & \qquad
	+
	C s^{-\nf 9 2}
	(s_2 - s)^{-1}
	\bigl(1 + \mathcal{E}^{2}_{\nf{s_1}{2}, \infty, 2}(\cu_m, n;\a,\s_0) \bigr)
	3^{s_2 n}
	[ \g ]_{\underline{B}^{s_2}_{p,p}(\cu_m)}
	\,.
\end{align}
\end{lemma}
\begin{proof}
By definition,
\begin{align*}
	\avsum_{z\in 3^n\Zd \cap \cu_m}
	\| (\a-\s_0) \nabla u \|_{\Besov{-s}{p}{p}(z+\cu_n)}^p
	&
	=
	\sum_{k=-\infty}^n
	3^{spk}
	\avsum_{z\in 3^{k}\Zd \cap \cu_m}
	\bigl| \bigl( (\a-\s_0) \nabla u \bigr)_{z+\cu_{k}} \bigr|^p
	\,.
\end{align*}
Applying Lemma~\ref{l.coarse.graining.operator.RHS}, we have
\begin{align*}
	\bigl| \bigl( (\a-\s_0) \nabla u \bigr)_{z+\cu_{k}} \bigr|
	&
	\leq
	C s^{-1} \s_0^{\nf12}  \| \s^{\nf12}  \nabla u  \|_{\underline{L}^2(z+\cu_{k})}
	\mathcal{E}_{s, \infty, 1} (z+\cu_{k};\a,\s_0)
	\notag \\ & \qquad
	+
	C  s^{-\nf 9 2}  \bigl(1+ \mathcal{E}^2_{\nf s 2, \infty, 2} ( z+\cu_{k};\a,\s_0) \bigr)
	3^{sk}
	[ \g ]_{\underline{H}^s(z+\cu_{k})}
	\,.
\end{align*}
We compute, using~\eqref{e.mathcalE.monotone.ordered} and~\eqref{e.bound.one.cube.by.mathcalE},
\begin{align*}
	\lefteqn{
	\sum_{k=-\infty}^n
	3^{sp(k-n)}
	\avsum_{z\in 3^{k}\Zd \cap \cu_m}
	\| \s^{\nf12}  \nabla u  \|^p_{\underline{L}^2(z+\cu_{k})}
	\mathcal{E}^p_{s, \infty, 1} (z+\cu_{k};\a,\s_0)
	} \qquad &
	\notag \\ &
	\leq
	\mathcal{E}^p_{s_1, \infty, 1} (\cu_m, n;\a,\s_0)
	\sum_{k=-\infty}^n
	3^{(s -s_1) p (k-n)}
	\avsum_{z\in 3^{k}\Zd \cap \cu_m}
	\| \s^{\nf12}  \nabla u  \|^p_{\underline{L}^2(z+\cu_{k})} \,\,.
\end{align*}
Similarly,
\begin{align*}
	\lefteqn{
	\sum_{k=-\infty}^n
	3^{s p  (k-n)}
	\avsum_{z\in 3^k\Zd\cap\cu_m}
	\bigl(1+ \mathcal{E}^2_{\nf s2, \infty, 2} ( z+\cu_{k};\a,\s_0) \bigr)^p
	3^{s p k}
	[ \g ]_{\underline{H}^s(z+\cu_{k})}^p
	}
	\qquad &
	\notag \\ &
	\leq
	C^p (1 + \mathcal{E}^{2 p}_{\nf{s_1}{2}, \infty, 2}(\cu_m, n;\a,\s_0))
	\sum_{k=-\infty}^n
	3^{(s -s_1) p (k-n)}
	\avsum_{z\in 3^{k}\Zd \cap \cu_m}
	3^{ s p k}
	[ \g ]_{\underline{H}^s(z+\cu_{k})}^p \,\,.
\end{align*}
We estimate the latter term using Jensen's inequality as follows:
\begin{align*}
	\lefteqn{
	3^{-s p n}
	\sum_{k=-\infty}^n
	3^{(s-s_1) p (k-n)}
	\avsum_{z\in 3^{k}\Zd \cap \cu_m}
	3^{ s pk}
	[ \g ]_{\underline{B}^{s}_{2,2}(z+\cu_k)}^p
	} \qquad &
	\notag \\ &
	\leq
	C
	\sum_{k=-\infty}^n
	3^{(2 s -s_1) p (k-n)}
	\avsum_{z\in 3^{k}\Zd \cap \cu_m}
	\Bigl(
	\sum_{j=-\infty}^k
	3^{-2 s j}
	\avsum_{z' \in 3^{j} \Zd \cap z + \cu_k}
	\bigl\| \g - \bigl(\g \bigr)_{z'+\cu_{j}} \bigr\|_{\underline{L}^2(z'+\cu_j)}^2
	\Bigr)^{\nf p 2}
	\notag \\ &
	\leq
	C (s_2 - s)^{-p}
	\sum_{k=-\infty}^n
	3^{(2 s -s_1) p (k-n)}
	3^{ (s_2 - s) p k}
	\avsum_{z\in 3^{k}\Zd \cap \cu_m}
	[ \g ]_{\underline{B}^{s_2}_{p,p}(z + \cu_k)}^p
	\notag \\ &
	\leq
	C (s_2 - s)^{-p}
	3^{(s_2 - s) p n}
	[ \g ]_{\underline{B}^{s_2}_{p,p}(\cu_m)}^p
	\,.
\end{align*}
The estimate~\eqref{e.gen.coarse.grain} follows from the previous displays.
\end{proof}

\subsection{Sensitivity estimates}
\label{ss.sensitivity}

In this subsection, we quantify the dependence of the coarse-grained matrices and the ellipticity lower bound~$\lambda_{s,q}$ on antisymmetric perturbations~$\h$ of~$\a$. The basic message is that these objects are locally Lipschitz (or differentiable) with respect to~$\h$, provided~$\h$ is smooth enough, with Lipschitz constants expressed in terms of coarse-grained ellipticity rather than microscopic uniform ellipticity. This is the mechanism which will allow us to propagate control of the coarse-grained matrices along the RG flow.

\smallskip

We begin with a basic lemma which considers~$L^\infty$-type perturbations for quantitatively uniformly elliptic coefficient fields.

\begin{lemma} [Sensitivity for uniformly elliptic fields]
\label{l.crude.sensitivity}
Let~$U\subseteq\Rd$ be a bounded Lipschitz domain and~$\a_1,\a_2: U \to \R^{d\times d}$ be a pair of uniformly elliptic coefficient fields with~$\a_i = \s_i+\k_i$, where~$\s_i:U\to\R^{d \times d}_{\sym}$ and~$\k_i:U\to\R^{d \times d}_{\skew}$. Denote
\begin{equation*}
	\bfA_i  \coloneqq  \begin{pmatrix}
		\s_i + \k_i^t\s_i^{-1}\k_i
		& -\k_i^t\s_i^{-1}
		\\ -  \s_i^{-1}\k_i
		& \s_i^{-1}
		\end{pmatrix},\qquad i\in\{1,2\}\,\,.
\end{equation*}
Then
\begin{equation}
	\label{e.crude.sensitivity}
	\bfA_1^{-\nf12}(U) \bfA_2(U) \bfA_1^{-\nf12}(U) - \Itwod \leq \| \bfA_1^{-\nf 12} \bfA_2  \bfA_1^{-\nf 12} -  \Itwod\|_{L^{\infty}(U)} \Itwod
	\,\,.
\end{equation}
\end{lemma}
\begin{proof}
Let~$P \in\R^{2d}$, and set~$S\coloneqq S(\cdot,U,-P,0\,;\a_1) \in P+L^2_{\pot,0}(U) \times \Lsolo(U)$. Recall that~$S$ is the minimizer of~$\mu(U,P\,; \a_1) = \frac12P\cdot \bfA_1(U)P$ in~\eqref{e.variational.mu.U.P}. We compute
\begin{align*}
	P \cdot \bfA_2(U) P
	&
	=
	\min_{X \in P+L^2_{\pot,0}(U) \times \Lsolo(U)}
	\fint_{U}
	X \cdot \bfA_2 X
	\notag \\ &
	\leq
	\fint_{U}
	S \cdot \bfA_2 S
	=
	P\cdot \bfA_1(U)P
	+
	\fint_{U}
	S \cdot \bfA_1^{\nf 12} ( \bfA^{-\nf 12}_1 \bfA_2  \bfA^{-\nf 12}_1 -  \Itwod) \bfA_1^{\nf 12} S
	\,.
\end{align*}
Moreover,
\begin{align*}
	\fint_{U}
	S \cdot \bfA_1^{\nf 12} ( \bfA^{-\nf 12}_1 \bfA_2  \bfA^{-\nf 12}_1 -  \Itwod) \bfA_1^{\nf 12} S
	&\leq \| \bfA^{-\nf 12}_1 \bfA_2  \bfA^{-\nf 12}_1 -  \Itwod\|_{L^{\infty}(U)}
	\fint_{U}
	S \cdot \bfA_1 S  \\
	&=  \| \bfA^{-\nf 12}_1 \bfA_2  \bfA^{-\nf 12}_1 -  \Itwod\|_{L^{\infty}(U)}
	P \cdot \bfA_1(U) P
	\,.
\end{align*}
Combining the previous two displays yields~\eqref{e.crude.sensitivity}.
\end{proof}

We next intend to differentiate the coarse-grained matrices with respect to perturbations in the field, with a focus on antisymmetric perturbations in view of the applications to the present paper.
For an arbitrary matrix field~$\h\in L^\infty(U;\R^{d\times d})$, we introduce the symmetric
matrix~$\mathbf{D}_{\h}(U\,;\a) \in \R^{2d\times 2d}$ defined by
\begin{equation}
	\label{e.def.Dh}
	\begin{pmatrix} p \\ q \end{pmatrix} \cdot \mathbf{D}_{\h}(U\,;\a) \begin{pmatrix} p \\ q \end{pmatrix}
	=
	\fint_{U}
	\nabla v(\cdot,U,p,q;\a^t) \cdot  \h \nabla v(\cdot,U,p,-q;\a)
	\,.
\end{equation}
Notice that
\begin{equation*}
	\begin{pmatrix} p \\ - q \end{pmatrix} \cdot \mathbf{D}_{\h^t}(U\,;\a^t) \begin{pmatrix} p \\ - q \end{pmatrix}
	=
	\begin{pmatrix} p \\ q \end{pmatrix} \cdot \mathbf{D}_{\h}(U\,;\a) \begin{pmatrix} p \\ q \end{pmatrix}
	\,.
\end{equation*}

\begin{lemma}[Differentiating~$\bfA(U;\a)$]
\label{l.sensitivity.coarse.grained.general}
Let~$U \subseteq\Rd$ be a bounded Lipschitz domain,~$\a:U\to \R^{d\times d}$ a uniformly elliptic coefficient field with symmetric part~$\s \coloneqq  \frac12(\a+\a^t)$. For every antisymmetric field~$\h\in L^\infty(U;\R^{d\times d}_{\skew})$ and~$P\in\R^{2d}$,
we have
\begin{equation}
	\label{e.sensitivity.coarse.grained.general}
	P\cdot\bigl( \bfA(U; \a + \h) - \bfA(U; \a) - \mathbf{D}_\h(U\,;\a) \bigr)P
	\leq
	\| \s^{-\nf12} \h \s^{-\nf 12} \|_{L^\infty(U)}^2
	P\cdot\bfA(U;\a)P\,.
\end{equation}
\end{lemma}
\begin{proof}
Denote, as usual, the antisymmetric part of~$\a$ by~$\k  \coloneqq  \frac12(\a-\a^t)$.
Let~$P  \coloneqq  ( p, q )$ for~$p,q\in\Rd$. As explained below~\eqref{e.maximizers.J.to.bfJ}, the minimizer of the variational problem in~\eqref{e.variational.mu.U.P} is~$S = S(\cdot,U,-P,0\,;\a)$ given in~\eqref{e.maximizers.J.to.bfJ}. That is,
\begin{equation*}
	S(\cdot,U,-P,0) =
	-\begin{pmatrix}
		\nabla w +\nabla w^* \\
		\a\nabla w - \a^t\nabla w^*
	\end{pmatrix}
	\,,
	\quad \mbox{where} \quad
	w  \coloneqq  \frac12 v(\cdot, U, p, -q\,; \a)
	\,, \quad
	w^*  \coloneqq  \frac12 v(\cdot, U, p, q\,; \a^t) \,\,.
\end{equation*}
Observe that
\begin{align*}
	P \cdot \bfA(U; \a + \h) P
	&
	=
	\inf_{X \in (p+L^2_{\pot,0}(U)) \times (q + \Lsolo(U))}
	\fint_{U} X \cdot \G_{-\h}^t \bfA \G_{-\h} X
	\notag \\ &
	\leq
	\fint_{U} S \cdot \G_{-\h}^t \bfA \G_{-\h} S
	= P \cdot \bfA(U; \a) P
	+
	\fint_{U} S \cdot (\G_{-\h}^t \bfA \G_{-\h} - \bfA) S
	\,.
\end{align*}
We compute using~\eqref{e.bfA.LDLt} and~\eqref{e.Gh.additive},
\begin{align*}
	(\G_{-\h}^t \bfA \G_{-\h} - \bfA)
	&
	=
	\G_{-\k}^t
	\left(
	\G_{-\h}^t
	\begin{pmatrix}
		\s & 0 \\
		0 & \s^{-1}
	\end{pmatrix}
	\G_{-\h}
	-
	\begin{pmatrix}
		\s & 0 \\
		0 & \s^{-1}
	\end{pmatrix}
	\right)
	\G_{-\k}  \\
	&=
	\G_{-\k}^t
	\begin{pmatrix}
		( \h^t\s^{-1}\h )
		& -(\h^t\s^{-1})
		\\ - ( \s^{-1}\h )
		& 0
	\end{pmatrix}
	\G_{-\k}
\end{align*}
and
\begin{equation*}
	\bfG_{-\k}
	\begin{pmatrix} \nabla (w+w^*) \\   \a\nabla w - \a^t \nabla w^*  \end{pmatrix}
	=
	\begin{pmatrix} \Id   & 0 \\ -\k & \Id  \end{pmatrix}
	\begin{pmatrix} \nabla (w + w^*)  \\ (\s+\k) \nabla w  - (\s - \k) \nabla w^* \end{pmatrix}
	=
	\begin{pmatrix} \nabla w + \nabla w^*  \\  \s (\nabla w -\nabla w^*) \end{pmatrix}
	\,.
\end{equation*}
The previous three displays and the assumption that~$\h$ is antisymmetric imply
\begin{align*}
	&P \cdot (\bfA(U; \a + \h) - \bfA(U; \a)) P   \\
	&\qquad \leq \fint_{U} \begin{pmatrix} \nabla w + \nabla w^*  \\  \s (\nabla w -\nabla w^*) \end{pmatrix} \cdot
	\begin{pmatrix}
		( \h^t\s^{-1}\h )
		& -(\h^t\s^{-1})
		\\ - ( \s^{-1}\h )
		& 0
		\end{pmatrix}\begin{pmatrix} \nabla w + \nabla w^*  \\  \s (\nabla w -\nabla w^*) \end{pmatrix}
	\\ &\qquad
	=
	\fint_{U}
	(\nabla w + \nabla w^*) \cdot \left(\h^t \s^{-1} \h\right)(\nabla w + \nabla w^*)
	+ 2 (\nabla w + \nabla w^*)  \cdot  \h (\nabla w - \nabla w^*)
	\,.
\end{align*}
We next use that
\begin{equation*}
	\fint_U 2 (\nabla w + \nabla w^*)  \cdot  \h (\nabla w - \nabla w^*)
	=
	\fint_U 4 \nabla w^* \cdot  \h \nabla w
	=
	\fint_U \nabla v^* \cdot  \h \nabla v
	=P\cdot  \mathbf{D}_{\h}(U\,;\a) P
\end{equation*}
and that, by~\eqref{e.Jenergyv.nosymm} and~\eqref{e.J.by.means.of.bfA},
\begin{equation*}
	\fint_{U}
	|  \s^{\nf 12} \nabla (w + w^*) |^2
	\leq
	2\fint_{U}
	|  \s^{\nf 12} \nabla w |^2
	+
	2\fint_{U}
	|  \s^{\nf 12} \nabla w^* |^2
	=
	P \cdot \bfA(U;\a) P
	\,.
\end{equation*}
Inserting these into the previous display, we obtain
\begin{align*}
	P \cdot (\bfA(U; \a + \h) - \bfA(U; \a) - \mathbf{D}_{\h}(U\,;\a)) P
	&
	\leq
	\| \s^{-\nf12} \h \s^{-\nf 12} \|_{L^\infty(U)}^2 P \cdot \bfA(U;\a) P
	\,.
\end{align*}
This completes the proof of~\eqref{e.sensitivity.coarse.grained.general}.
\end{proof}

Applying the preceding upper estimate at the shifted base field~$\a+t\h$ with increment~$-t\h$, and using continuity of the minimizing fields as~$t\to0$, gives the matching lower derivative. Thus the matrix~$\mathbf{D}_\h(U;\a)$ defined in~\eqref{e.def.Dh} is the derivative of the coarse-grained matrix~$\bfA(U;\a)$ with respect to an antisymmetric perturbation~$\h$ of the field.
By using~\eqref{e.J.by.means.of.bfA}, we can write this in terms of~$J$ rather than~$\bfA$ as
\begin{equation}
	\label{e.derivative.J}
	\frac{\partial}{\partial t}
	J(U,p,q\,;\a+t\h) \Bigr\vert_{t=0}
	=
	\frac12
	\begin{pmatrix} p \\ -q \end{pmatrix} \cdot \mathbf{D}_{\h}(U\,;\a) \begin{pmatrix} p \\ -q \end{pmatrix}
	\,.
\end{equation}
A very crude way to bound the size of the right side is to start from the expression in~\eqref{e.def.Dh} and apply Cauchy-Schwarz and~\eqref{e.Jenergyv.nosymm} to obtain
\begin{align}
	\lefteqn{
	\biggl| \begin{pmatrix} p \\ -q \end{pmatrix} \cdot \mathbf{D}_{\h}(U\,;\a) \begin{pmatrix} p \\ -q \end{pmatrix} \biggr|
	} \qquad &
	\notag \\ &
	\leq
	\| \s^{-\nf12} \h \s^{-\nf 12} \|_{L^\infty(U)}
	\| \s^{\nf12} \nabla v(\cdot,U,p,-q;\a^t) \|_{\underline{L}^2(U)}
	\| \s^{\nf12} \nabla v(\cdot,U,p,q;\a)  \|_{\underline{L}^2(U)}
	\notag \\ &
	\leq
	2\| \s^{-\nf12} \h \s^{-\nf 12} \|_{L^\infty(U)}
	J(U,p,-q;\a^t)^{\nf12}
	J(U,p,q;\a)^{\nf12}
	\,.
	\label{e.bad.derivative.bound}
\end{align}
For our purposes, this estimate is not satisfactory because of the factor of~$\| \s^{-\nf12} \h \s^{-\nf 12} \|_{L^\infty(U)}$ on the right, which measures the ratio of the perturbation and the field~$\s$ in a pointwise/$L^\infty$ way. In most situations, we can do no better than to bound this by~$\| \s^{-1}\|_{L^\infty(U)} \| \h \|_{L^\infty(U)}$, and so we see the inverse of the \emph{microscopic} ellipticity lower bound constant~$\| \s^{-1}\|_{L^\infty(U)}^{-1}$ appearing on the right side. The same microscopic ratio shows up implicitly in the right side of~\eqref{e.crude.sensitivity}. To handle multifractal problems, we must do better.

\smallskip

We next present a \emph{coarse-grained} version of the estimate~\eqref{e.bad.derivative.bound}, which replaces~$\| \s^{-1} \|_{L^\infty(\cu_0)}$ with~$\lambda_s^{-1}(\cu_0\,;\a)$, the inverse of the coarse-grained ellipticity constant. This coarse-grained derivative estimate plays a fundamental role in the analysis in the rest of the paper. The price to pay for the replacement of the microscopic ellipticity by the coarse-grained one is higher derivatives on the perturbation: the right side of~\eqref{e.Dh.bound} has~$W^{2,\infty}$ norms of~$\h$, compared to the~$L^\infty$ norm in~\eqref{e.bad.derivative.bound}.

\begin{lemma}[Coarse-grained derivative bound]
\label{l.Dh.bound}
Let~$\a:\cu_0 \to \R^{d\times d}$ be a uniformly elliptic coefficient field. There exists a constant~$C(d) < \infty$ such that, for every~$\h\in W^{2,\infty}(\cu_0;\R^{d\times d}_{\skew})$  and~$p,q \in \R^{d}$,
\begin{align}
	\label{e.Dh.bound}
	\lefteqn{
	\biggl| \begin{pmatrix} p \\ -q \end{pmatrix} \cdot \mathbf{D}_{\h}(\cu_0\,;\a) \begin{pmatrix} p \\ -q \end{pmatrix} \biggr|
	} \qquad &
	\notag \\ &
	\leq
	C \Bigl(  |p|  \| \h \|_{W^{1,\infty}(\cu_0)}  \lambda_{\nf 38,2}^{-\nf12}(\cu_0\,; \a)  +
	|p\cdot q|^{\nf12} \| \nabla \h \|_{W^{1,\infty}(\cu_0)}
	\lambda_{\nf 38,2}^{-1}(\cu_0\,; \a)
	\Bigr)
	J(\cu_0,p,q\, ; \a)^{\nf 12}
	\notag \\ & \qquad
	+
	C \| \nabla \h \|_{W^{1,\infty}(\cu_0)}   \lambda_{\nf 38,2}^{-1}(\cu_0\,; \a)   J(\cu_0,p,q\, ; \a)
	\,.
\end{align}
\end{lemma}

\begin{proof}
We denote
\begin{equation}
	\label{e.wstarplusw.spatavg}
	P  \coloneqq \begin{pmatrix} p \\ -q \end{pmatrix}\,, \quad
	w  \coloneqq  \frac12 v(\cdot, \cu_0, p, q\,; \a)
	\qand
	w^*  \coloneqq  \frac12 v(\cdot, \cu_0, p, -q\,; \a^t) \,\,.
\end{equation}
By the antisymmetry of~$\h$ (pointwise), and then by integration by parts using~$w+w^*+\ell_p\in H^1_0(\cu_0)$, we find that
\begin{align}
	\label{e.sensitivity.basic.split}
	\fint_{\cu_0}
	\nabla w^* \cdot \h\nabla w
	=
	\fint_{\cu_0}
	\nabla (w^* + w) \cdot \h\nabla w
	&
	=
	-p  \cdot \fint_{\cu_0}  \h \nabla w
	+
	\fint_{\cu_0}
	\nabla (w + w^* + \ell_{p})   \cdot \h \nabla w
	\notag \\ &
	=
	-p  \cdot \fint_{\cu_0}  \h \nabla w
	-
	\fint_{\cu_0}
	(w + w^* + \ell_{p}) (\nabla \cdot \h ) \cdot \nabla w
	\,.
\end{align}
By~\eqref{e.Jenergyv.nosymm}, we have
\begin{equation*}
	\| \s^{\nf 12} \nabla w \|_{\underline{L}^2(\cu_0)}^2  = \frac12 J(\cu_0,p,q\,; \a) \,,
\end{equation*}
thus by the coarse-grained Poincar\'e inequality~\eqref{e.besov.grad.poincare},
\begin{equation}
	\label{e.w.wstar.to.energy}
	\| \nabla w \|_{\Besov{-3/8}{2}{2}(\cu_0)}
	\leq
	C \lambda_{\nf 38,2}^{-\nf 12}(\cu_0; \a) J(\cu_0,p,q\, ; \a)^{\nf 12}    \,\,.
\end{equation}
Similarly, by~\eqref{e.Jenergyv.nosymm} and~\eqref{e.J.by.means.of.bfA}, we have
\begin{equation*}
	\| \s^{\nf 12} \nabla w^* \|_{\underline{L}^2(\cu_0)}^2 + \| \s^{\nf 12} \nabla w \|_{\underline{L}^2(\cu_0)}^2
	=
	\frac12 P \cdot \bfA(\cu_0;\a) P
	\,.
\end{equation*}
The block formula~\eqref{e.bigA.to.little.A} and~$\s_*(\cu_0)\leq\s(\cu_0)$ give
\begin{equation*}
	P\cdot\bfA(\cu_0;\a)P
	=p\cdot\s(\cu_0)p
	+\bigl(q+\k(\cu_0)p\bigr)\cdot\s_*^{-1}(\cu_0)
	\bigl(q+\k(\cu_0)p\bigr)
	\geq p\cdot\s_*(\cu_0)p\,.
\end{equation*}
Consequently, by~\eqref{e.ellipticities.monotone.ordered},
\begin{equation*}
	|p|^2
	\leq
	|\s_*^{\nf12}(\cu_0) p |^2
	|\s_*^{-\nf12}(\cu_0)|^2
	\leq
	C \lambda_{\nf12,2}^{-1}(\cu_0\,; \a) P \cdot \bfA(\cu_0;\a) P  \,\,,
\end{equation*}
which implies
\begin{equation}
	\|  \nabla(w + w^* + \ell_{p})  \|_{ {\Besov{-\nf 12}{2}{2}(\cu_0)} }
	\leq
	C \lambda^{-\nf12}_{\nf 12,2}(\cu_0\,; \a) \bigl( P \cdot \bfA(\cu_0;\a) P \bigr)^{\nf 12}   \,\,.
	\label{e.w.wstar.with.other.stuff.to.energy}
\end{equation}
We estimate the first term on the right side of~\eqref{e.sensitivity.basic.split} using~\eqref{e.w.wstar.to.energy}:
\begin{align}
	\biggl|
	p  \cdot \fint_{\cu_0}  \h \nabla w
	\biggr|
	&
	\leq
	C |p|
	\| \h \|_{\underline{B}_{2,2}^{\nf 38}(\cu_0)}
	\| \nabla w \|_{\Besov{-\nf 38}{2}{2}(\cu_0)}
	\notag \\ &
	\leq
	C |p| \bigl( | (\h)_{\cu_0}| +  \| \nabla \h \|_{L^\infty(\cu_0)} \bigr)   \lambda_{\nf 38,2}^{-\nf12}(\cu_0\,;\a)
	\bigl( J(\cu_0,p,q\, ; \a) \bigr)^{\nf 12}
	\,.
	\label{e.sensitivity.basic.split1}
\end{align}
To estimate the second term, we compute
\begin{equation}
	\biggl| \fint_{\cu_0} (w + w^* + \ell_{p}) (\nabla \cdot \h ) \cdot \nabla w\biggr|
	\leq C
	\bigl\| (w + w^* + \ell_{p}) (\nabla \cdot \h) \bigr\|_{\underline{B}_{2,2}^{\nf 38}(\cu_0)}
	\| \nabla w \|_{\Besov{-\nf 38}{2}{2}(\cu_0)}
	\label{e.split.besov.stuff.begin}
	\,.
\end{equation}
By~\eqref{e.v.in.Bpps.vs.nabla.v.in.B.pp.minus.t}, the strict embedding~$B^{1/2}_{2,\infty}\subseteq B^{3/8}_{2,2}$, and~\eqref{e.w.wstar.with.other.stuff.to.energy}, we obtain
\begin{equation*}
	[ w + w^* + \ell_{p}  ]_{B_{2,2}^{\nf 38}(\cu_0)}
	\leq
	C \| \nabla (w + w^*) + p  \|_{\underline{B}_{2,2}^{-\nf 12 }(\cu_0)}
	\leq
	C\lambda_{\nf12,2}^{-\nf12}(\cu_0\,; \a) \bigl( P \cdot \bfA(\cu_0;\a) P \bigr)^{\nf 12}
	\,.
\end{equation*}
On the other hand, by~\eqref{e.negative.Besov.to.L2.zero.bndr} and~\eqref{e.w.wstar.with.other.stuff.to.energy} we deduce that
\begin{equation*}
	\|  w + w^* + \ell_{p}  \|_{L^{2}(\cu_0)}
	\leq
	C
	\|  \nabla(w + w^*) + p  \|_{ {\Besov{-\nf 12}{2}{2}(\cu_0)} }
	\leq
	C \lambda^{-\nf12}_{\nf 12,2}(\cu_0\,; \a) \bigl( P \cdot \bfA(\cu_0;\a) P \bigr)^{\nf 12}
	\,.
\end{equation*}
Using the above two displays together with~\eqref{e.besov.mult}, we obtain
\begin{align}
	\lefteqn{
	\bigl\| (w + w^* + \ell_{p}) (\nabla \cdot \h) \bigr\|_{\underline{B}_{2,2}^{\nf38}(\cu_0)}
	} \qquad &
	\notag \\ &
	\leq
	C \|  w + w^* + \ell_{p}  \|_{L^{2}(\cu_0)} \| \nabla^2 \h \|_{L^\infty(\cu_0)} + C \|\nabla \cdot \h\|_{L^{\infty}(\cu_0)}
	[ w + w^* + \ell_{p}  ]_{B_{2,2}^{\nf 38}(\cu_0)} \notag \\ &
	\leq
	C \| \nabla \h \|_{W^{1,\infty}(\cu_0)}  \lambda^{-\nf12}_{\nf 12, 2}(\cu_0\,; \a) \bigl( P \cdot \bfA(\cu_0;\a) P \bigr)^{\nf 12}
	\,.
	\label{e.split.besov.stuff}
\end{align}
Substituting the previous display and~\eqref{e.w.wstar.to.energy} into~\eqref{e.split.besov.stuff.begin}, we obtain
\begin{align*}
	\biggl| \fint_{\cu_0} (w + w^* + \ell_{p}) (\nabla \cdot \h ) \cdot \nabla w \biggr|
	\leq
	C  \| \nabla \h \|_{W^{1,\infty}(\cu_0)}   \lambda_{\nf 38,2}^{-1}(\cu_0\,; \a)
	\bigl( P \cdot \bfA(\cu_0;\a) P \bigr)^{\nf 12}   \bigl( J(\cu_0,p,q\, ; \a) \bigr)^{\nf 12}
	\,\,.
\end{align*}
Plugging this and~\eqref{e.sensitivity.basic.split1} into~\eqref{e.sensitivity.basic.split} and applying~\eqref{e.J.by.means.of.bfA} completes the proof.
\end{proof}

We next differentiate the coarse-grained ellipticity lower bound constant~$\lambda_s^{-1}(U;\a)$ with respect to antisymmetric perturbations of~$\a$. Since~$\lambda_s^{-1}(U;\a)$ is a multiscale amalgamation of~$\s_*(\cu)$'s, we obtain our estimate by applying the result of the previous two lemmas, specialized to the bottom right entry of~$\bfA(\cu)$, to all triadic subcubes of~$U$.

Note that the estimate in the next lemma is conditional---it requires the perturbation to be small at the coarse-grained scale, see~\eqref{e.lambda.sensitivity.smallness.condition}. This will be removed later in Lemma~\ref{l.J.sensitivity.no.conditions}.

\begin{lemma}[{Sensitivity estimate for~$\lambda_{s,q}$}]
\label{l.lambda.sensitivity}
Let~$\a:\cu_0 \to \R^{d\times d}$ be a uniformly elliptic coefficient field.  There exists~$C(d) < \infty$ such that, for every~$\h\in W^{2,\infty}(\cu_0;\R^{d\times d}_{\skew})$ satisfying
\begin{equation}
	\| \nabla \h \|_{W^{1,\infty}(\cu_0)}
	\leq
	C^{-1}  \lambda_{\nf38,2}(\cu_0 \, ; \a) \,,
	\label{e.lambda.sensitivity.smallness.condition}
\end{equation}
we have, for every~$s\in (0,\nf 12]$ and~$q\in [1,\infty]$,  the estimate
\begin{equation}
	\label{e.lambda.sensitivity}
	\lambda_{s,q}^{-1} (\cu_0 ; \a + \h)
	\leq
	\bigl(
	1 + C  \| \nabla \h \|_{W^{1,\infty}(\cu_0)}  \lambda_{\nf 38,2}^{-1} (\cu_0 \, ; \a)
	\bigr)
	\lambda_{s,q}^{-1} (\cu_0\, ; \a)
	\,.
\end{equation}
\end{lemma}
\begin{proof}
For each~$n \in - \N_0$ and~$z\in 3^{n}\Zd \cap \cu_0$, we apply~\eqref{e.Dh.bound} and~\eqref{e.sensitivity.coarse.grained.general} with~$p = 0$ and~$q \in \Rd$ with~$|q|=1$, and with~$U=\cu_0$ and then rescale and translate so that~$z+\cu_n$ appears in place of~$\cu_0$. (In other words we apply~\eqref{e.Dh.bound} and~\eqref{e.sensitivity.coarse.grained.general} to the rescaled field~$\a_{z,n}(x)\coloneqq \a(z+3^n x)$ on~$\cu_0$, and then transfer the result back to the cube~$z+\cu_n$.)
We obtain, for every~$\tau \in (-1,1)$,
\begin{align}
	\label{e.sensitivity.coarse.grained.general.bottom.right}
	\lefteqn{
	| \s_*^{-1} (z+\cu_{n}\,; \a + \tau \h) |
	} \quad &
	\notag \\  &
	\leq
	\bigl(1+
	C|\tau| \lambda_{\nf38,2}^{-1}(z {+} \cu_n\,; \a)
	3^{2n} \| \nabla \h \|_{\underline{W}^{1,\infty}(z+ \cu_n)}
	+ \tau^2 \| \s^{-\nf12} \h \s^{-\nf 12} \|_{L^\infty(z + \cu_n)}^2 \bigr) \bigl| \s_*^{-1} (z{+}\cu_{n}; \a) \bigr|
	\notag \\ &
	\leq
	\bigl(1+
	C|\tau| \lambda_{\nf38,2}^{-1}( \cu_0\,; \a)
	\| \nabla \h \|_{\underline{W}^{1,\infty}( \cu_0)}
	+ \tau^2 \| \s^{-\nf12} \h \s^{-\nf 12} \|_{L^\infty(z + \cu_n)}^2 \bigr) \bigl| \s_*^{-1} (z{+}\cu_{n}; \a) \bigr|
	\,,
\end{align}
where in the last line we used~\eqref{e.bound.one.cube.by.lambdas} and that~$\| \nabla \h \|_{\underline{W}^{1,\infty}(z+ \cu_n)} \leq 3^{-n} \| \nabla \h \|_{\underline{W}^{1,\infty}(\cu_0)}$.

\smallskip

From~\eqref{e.sensitivity.coarse.grained.general.bottom.right}
we obtain, for every~$s \in (0,\nf 12]$,~$q\in[1,\infty]$ and~$\tau \in (-1,1)$,
\begin{align*}
	\lambda_{s,q}^{-1}(\cu_0\,;\a+\tau \h)
	\leq
	\bigl(1+
	C|\tau| \lambda_{\nf 38,2}^{-1}(\cu_0\,;\a)   \| \nabla \h \|_{W^{1,\infty}(\cu_{0})} + |\tau|^2 \| \s^{-\nf12} \h \s^{-\nf 12} \|_{L^\infty(\cu_0)}^2 \Bigr)
	\lambda_{s,q}^{-1}(\cu_0\,;\a)
	\,.
\end{align*}
For~$r>0$ and~$q\in[1,\infty]$, set~$\Phi_{r,q}(t)\coloneqq
\lambda_{r,q}^{-1}(\cu_0;\a+t\h)$. Applying the preceding estimate at the base field~$\a+t\h$ gives the one-step ratio
\begin{equation}
	\label{e.lambda.sensitivity.pre}
	\Phi_{s,q}(t+\tau)
	\leq
	\Bigl(1+C|\tau|\|\nabla\h\|_{W^{1,\infty}(\cu_0)}
	\Phi_{\nf38,2}(t)+R\tau^2\Bigr)\Phi_{s,q}(t)\,,
\end{equation}
where~$R<\infty$ is the fixed microscopic quadratic ratio along the segment. We compose this estimate along the uniform partition~$t_k=k/N$. The sum of the quadratic remainders is~$O(N^{-1})$ and hence vanishes as~$N\to\infty$. First applying the argument at~$(r,q)=(\nf38,2)$ and using~\eqref{e.lambda.sensitivity.smallness.condition} gives
\begin{equation*}
	\sup_{t\in[0,1]}\Phi_{\nf38,2}(t)
	\leq 2\Phi_{\nf38,2}(0)\,.
\end{equation*}
The same finite-step argument at general~$(s,q)$ therefore yields
\begin{equation*}
	\lambda_{s,q}^{-1}(\cu_0;\a+\h)
	\leq
	\exp\Bigl(C\|\nabla\h\|_{W^{1,\infty}(\cu_0)}
	\lambda_{\nf38,2}^{-1}(\cu_0;\a)\Bigr)
	\lambda_{s,q}^{-1}(\cu_0;\a)\,.
\end{equation*}
The smallness condition turns the exponential into the factor in~\eqref{e.lambda.sensitivity}.
\end{proof}

\begin{lemma}[{Sensitivity estimate for~$J$}]
\label{l.J.sensitivity}
Let~$\a:\cu_0 \to \R^{d\times d}$ be a uniformly elliptic coefficient field. There exists~$C(d) < \infty$ such that, for every~$\h\in W^{2,\infty}(\cu_0;\R^{d\times d}_{\skew})$ satisfying
\begin{equation}
	\| \nabla \h \|_{W^{1,\infty}(\cu_0)}
	\leq
	C^{-1}  \lambda_{\nf38,2}(\cu_0 \, ; \a) \,,
	\label{e.J.sensitivity.smallness.condition}
\end{equation}
we have, for every~$\delta \in (0,1]$ and~$p,q\in\Rd$,
\begin{align}
	\label{e.J.sensitivity}
	J(\cu_0,p,q\,; \a + \h)
	&
	\leq
	\bigl(1+  \delta + C   \| \nabla \h \|_{W^{1,\infty}(\cu_0)}   \lambda_{\nf 38,2}^{-1}( \cu_0 \,; \a) \bigr)   J(\cu_0,p,q\,; \a)
	\notag \\ & \qquad
	+ C \delta^{-1} \Bigl(  |p|^2  \| \h \|_{W^{1,\infty}(\cu_0)}^2 \lambda_{\nf 38,2}^{-1}(\cu_0\,; \a)
	+
	|p\cdot q| \| \nabla \h \|_{W^{1,\infty}(\cu_0)}^2
	\lambda_{\nf 38,2}^{-2}(\cu_0 \,; \a)  \Bigr)
\end{align}
and, for every~$s \in (0,\nf 38)$,
\begin{align}
	\label{e.big.Lambda.sensitivity}
	\lefteqn{
	\Lambda_{s,2}(\cu_0\,; \a + \h)
	} \quad  &
	\notag \\ &
	\leq
	\bigl(1+\delta+C\|\nabla\h\|_{W^{1,\infty}(\cu_0)}
	\lambda_{\nf38,2}^{-1}(\cu_0;\a)\bigr)\Lambda_{s,2}(\cu_0;\a)
	+C\delta^{-1}(\nf38-s)^{-1}\|\h\|_{W^{1,\infty}(\cu_0)}^2
	\lambda_{s,2}^{-1}(\cu_0;\a)
	\,.
\end{align}

\end{lemma}
\begin{proof}
Set
\begin{equation*}
	E\coloneqq |p|^2\|\h\|_{W^{1,\infty}(\cu_0)}^2
	\lambda_{\nf38,2}^{-1}(\cu_0;\a)
	+|p\cdot q|\|\nabla\h\|_{W^{1,\infty}(\cu_0)}^2
	\lambda_{\nf38,2}^{-2}(\cu_0;\a)\,.
\end{equation*}
Apply Lemmas~\ref{l.Dh.bound} and~\ref{l.lambda.sensitivity} along the uniform chain~$\a+(k/N)\h$,~$k=0,\ldots,N$. The one-step estimate is
\begin{align*}
	J\bigl(\cu_0,p,q;\a+\tfrac{k+1}{N}\h\bigr)
	&\leq
	\Bigl(1+\frac{C}{N}\|\nabla\h\|_{W^{1,\infty}(\cu_0)}
	\lambda_{\nf38,2}^{-1}(\cu_0;\a)+O(N^{-2})\Bigr)
	J\bigl(\cu_0,p,q;\a+\tfrac{k}{N}\h\bigr)
	\\ &\qquad
	+\frac{C}{N}E^{\nf12}
	J\bigl(\cu_0,p,q;\a+\tfrac{k}{N}\h\bigr)^{\nf12}\,.
\end{align*}
The~$O(N^{-2})$ microscopic remainders sum to~$O(N^{-1})$. Composing the chain, letting~$N\to\infty$, and applying Young's inequality to the square-root term with the free parameter~$\delta$ gives~\eqref{e.J.sensitivity}. This finite-step argument avoids any differentiability assumption on the variational supremum defining~$J$.

\smallskip

We turn to the proof of~\eqref{e.big.Lambda.sensitivity}. Since
\begin{equation*}
	\lambda_{\nf38,2}^{-1}(z {+} \cu_n\,; \a)
	3^{2n} \| \nabla \h \|_{\underline{W}^{1,\infty}(z+ \cu_n)}
	\leq
	C\| \nabla\h \|_{W^{1,\infty}(\cu_0)}
	\lambda_{\nf38,2}^{-1}(\cu_0\,; \a)
\end{equation*}
and~$3^{n} \| \h \|_{W^{1,\infty}(z+\cu_n)} \leq C \| \h \|_{W^{1,\infty}(\cu_0)}$,
the estimate~\eqref{e.J.sensitivity} yields
\begin{align*}
	\bigl| \b(z+\cu_n; \a+\h) \bigr|
	&
	\leq
	\bigl(1+\delta+C\|\nabla\h\|_{W^{1,\infty}(\cu_0)}
	\lambda_{\nf38,2}^{-1}(\cu_0;\a)\bigr)
	\bigl| \b(z+\cu_n; \a) \bigr|
	\notag \\ & \qquad
	+ C\delta^{-1}\| \h \|_{W^{1,\infty}(\cu_0)}^2
	\lambda_{\nf38,2}^{-1}(z+\cu_n\,; \a)
	\,.
\end{align*}
By definition~\eqref{e.coarse.grained.ellipticity}, we have
\begin{align*}
	\sum_{n=-\infty}^{0}
	3^{2sn}
	\max_{z\in 3^n\Zd \cap \cu_0}
	\lambda_{\nf 38,2}^{-1}(z+\cu_n\,; \a)
	&
	\leq
	\sum_{n=-\infty}^{0}
	3^{2sn}
	\sum_{k=-\infty}^{n}
	3^{-\frac34(n-k)}
	\max_{z \in 3^k \Zd \cap \cu_0}
	\bigl| \s_{*}^{-1}(z+\cu_k; \a) \bigr|
	\notag \\ &
	=
	\sum_{k=-\infty}^{0}
	3^{2sk}
	\max_{z \in 3^k \Zd \cap \cu_0}
	\bigl| \s_{*}^{-1}(z+\cu_k; \a) \bigr|
	\sum_{n=k}^{0}
	3^{(2s-\frac34)(n-k)}
	\notag \\ &
	\leq
	C (\nf38-s)^{-1} \css{2s}^{-1}  \lambda_{s,2}^{-1}(\cu_0\,; \a)
	\,.
\end{align*}
Combining the previous two displays yields~\eqref{e.big.Lambda.sensitivity}, completing the proof.
\end{proof}

The sensitivity estimates in Lemmas~\ref{l.lambda.sensitivity} and~\ref{l.J.sensitivity} require a smallness assumption of the form~$\|\nabla\h\|_{W^{1,\infty}} \lambda_{\nf38,2}^{-1} \ll 1$. For our purposes, this is sometimes too restrictive: there will be rare ``bad boxes'' where this condition is not satisfied, but for which we still need an estimate. The next lemma removes this smallness assumption. The idea is to choose a mesoscopic scale~$3^{-h}$ depending on~$\|\nabla\h\|_{W^{1,\infty}}\lambda_{s,2}^{-1}$, apply the conditional sensitivity estimates cube-by-cube at scale~$h$, and then patch the resulting bounds via subadditivity. If we are unlucky and the cube is ``bad,'' then we may need to take the scale parameter~$h$ larger, but we will still get an estimate using essentially coarse-scale ellipticity constants. The outcome is an ``unconditional'' sensitivity estimate for~$\lambda_{t,q}$ and~$J$ with a controlled power dependence on the dimensionless quantity~$\|\nabla\h\|_{W^{1,\infty}} \lambda_{s,2}^{-1}$.

\begin{lemma}[Unconditional sensitivity estimates]
\label{l.J.sensitivity.no.conditions}
Let~$\a:\cu_0 \to \R^{d\times d}$ be a uniformly elliptic coefficient field and~$\s_0 \in (0,\infty)$. There exists~$C(d) < \infty$ such that, for every~$\h\in W^{2,\infty}(\cu_0;\R^{d\times d}_{\skew})$ and every~$s,t \in (0,\nf14]$,~$q \in [1,\infty]$ and~$\mu\in (0,\infty)$, we have the estimates
\begin{equation}
	\label{e.lambda.sensitivity.no.conditions}
	\lambda_{t,q}^{-1} (\cu_{0}\, ; \a + \h)
	\leq 6 \Bigl( 1 + C \| \nabla  \h   \|_{\underline{W}^{1,\infty}(\cu_0)} \lambda_{s,2}^{-1}
	(\cu_0; \a ) \Bigr)^{\! \frac{2 t}{1-2s}}\lambda_{t,q}^{-1} (\cu_{0}\, ; \a)
\end{equation}
and, for every~$e \in \Rd$ with~$|e| = 1$,
\begin{align}
	\label{e.J.sensitivity.no.conditions}
	\lefteqn{
	J(\cu_0,(\mu \s_0)^{-\nf12}e,(\mu \s_0)^{\nf 12}e\,;\a+\h)
	} \qquad &
	\notag \\ &
	\leq
	C \mu^{-1} \bigl( 1 +
	\| \nabla  \h   \|_{{W}^{1,\infty}(\cu_0)} \lambda_{s,2}^{-1}
	(\cu_0; \a )
	\bigr)^{\! \frac{2 t}{1-2s} }
	\mathcal{E}_{t,2,2}(\cu_0;\a,\s_0)^2
	\notag \\ & \qquad
	+
	C\mu^{-1}
	\s_0\lambda_{s,2}^{-1} (\cu_0 \,; \a)
	\bigl( 1+
	\| \nabla  \h   \|_{{W}^{1,\infty}(\cu_0)} \lambda_{s,2}^{-1}
	(\cu_0; \a )
	\bigr)^{\frac{2 s}{1-2s}}
	\bigl(
	\s_0^{-2}
	\| \h \|_{\underline{L}^2(\cu_0)}^2 +
	(\mu-1)^2 \bigr)
	\notag \\ & \qquad
	+
	C \mu^{-2} \s_0^{-2}
	\| \nabla  \h   \|_{{W}^{1,\infty}(\cu_0)}^2
	+
	C
	\min\bigl\{1\,,
	\| \nabla  \h   \|_{{W}^{1,\infty}(\cu_0)} \lambda_{s,2}^{-1}
	(\cu_0; \a ) \bigr\}^2
	\,.
\end{align}
The same estimates hold, by rescaling, on every cube~$z+\cu_n$, since~$J$,~$\mathcal{E}_{t,2,2}$ and~$\lambda_{s,2}^{-1}$ are unchanged by the translation and dilation~$x \mapsto z+3^n x$ of the cube together with the corresponding rescaling of~$\a$ and~$\h$.
\end{lemma}
\begin{proof}
We use subadditivity in the form
\begin{equation*}
	J(\cu_{0},(\mu \s_0)^{-\nf 12}e,(\mu \s_0)^{\nf 12}e\,; \a+\h)
	\leq
	\avsumcube{z}{-h}{0}
	J(z+\cu_{-h},(\mu \s_0)^{-\nf 12}e,(\mu \s_0)^{\nf 12}e\,; \a+\h)
	\,,
\end{equation*}
where the scale separation parameter~$h\in\N$ is given explicitly by
\begin{equation*}
	h \coloneqq  \biggl\lceil  \frac1{(1-2s)} \log_3 \Bigl(
	1 +  C_{\eqref{e.lambda.sensitivity.smallness.condition}\&\eqref{e.J.sensitivity.smallness.condition}}  \| \nabla  \h   \|_{\underline{W}^{1,\infty}(\cu_0)} \lambda_{s,2}^{-1}
	(\cu_0; \a )   \Bigr) \biggr\rceil
	\,.
\end{equation*}
The definition of~$h$ ensures that
\begin{equation}
	3^{(1-2s)h-1 }
	\leq
	1 + C_{\eqref{e.lambda.sensitivity.smallness.condition}\&\eqref{e.J.sensitivity.smallness.condition}}   \| \nabla  \h   \|_{\underline{W}^{1,\infty}(\cu_0)} \lambda_{s,2}^{-1}
	(\cu_0; \a )
	\leq
	3^{(1-2s)h}
	\,.
	\label{e.another.way.to.say.h}
\end{equation}
Observe that, by~\eqref{e.ellipticities.monotone.ordered} and~\eqref{e.bound.one.cube.by.lambdas}, for every~$z \in 3^{-h}\Zd \cap \cu_0$,
\begin{equation*}
	\lambda_{\nf 38,2}^{-1}(z+\cu_{-h} \,; \a)
	\leq
	\lambda_{s,2}^{-1}(z+\cu_{-h} \,; \a)
	\leq
	3^{2s h}
	\lambda_{s,2}^{-1}(\cu_0 \,; \a)
\end{equation*}
and thus, by the previous two displays,
\begin{align*}
	3^{-2h} \| \nabla  \h   \|_{\underline{W}^{1,\infty}(z+\cu_{-h})}
	&
	\leq
	3^{-h} \| \nabla  \h   \|_{\underline{W}^{1,\infty}(\cu_{0})}
	\notag \\ &
	\leq
	C_{\eqref{e.lambda.sensitivity.smallness.condition}\&\eqref{e.J.sensitivity.smallness.condition}}^{-1}
	3^{-2sh} \lambda_{s,2} (\cu_0; \a )
	\leq
	C_{\eqref{e.lambda.sensitivity.smallness.condition}\&\eqref{e.J.sensitivity.smallness.condition}} ^{-1} \lambda_{\nf 38,2} (z+\cu_{-h}; \a)
	\,.
\end{align*}
This is the condition we need to apply Lemmas~\ref{l.lambda.sensitivity} and~\ref{l.J.sensitivity} to switch from the field~$\a$ to~$\a+\h$ in each subcube~$z + \cu_{-h}\subseteq \cu_0$. An application of Lemma~\ref{l.lambda.sensitivity} gives us
\begin{equation*}
	\lambda_{t,q}^{-1} (z+\cu_{-h} ; \a + \h)
	\leq
	2 \lambda_{t,q}^{-1} (z+\cu_{-h}\, ; \a) \,.
\end{equation*}
Since, by~\eqref{e.bound.one.cube.by.lambdas},
\begin{equation*}
	\lambda_{t,q}^{-1} (z+\cu_{-h}\, ; \a) \leq 3^{2 th} \lambda_{t,q}^{-1} (\cu_{0}\, ; \a)\,,
\end{equation*}
we then obtain by the subadditivity of~$\lambda^{-1}$ that
\begin{equation*}
	\lambda_{t,q}^{-1} (\cu_{0}\, ; \a + \h)
	\leq 6 \Bigl( 1 + C \| \nabla  \h   \|_{\underline{W}^{1,\infty}(\cu_0)} \lambda_{s,2}^{-1}
	(\cu_0; \a ) \Bigr)^{\! \frac{2 t}{1-2s}}\lambda_{t,q}^{-1} (\cu_{0}\, ; \a)
	\,,
\end{equation*}
which is~\eqref{e.lambda.sensitivity.no.conditions}.

\smallskip

The estimate for~$J$ is similar. Using Lemma~\ref{l.J.sensitivity} with~$\delta = 1$, followed by an application of~\eqref{e.shaking.lambda} and then~\eqref{e.ellipticities.monotone.ordered}, yields, for every~$z \in 3^{-h}\Zd \cap \cu_0$,
\begin{align*}
	J(z+\cu_{-h},(\mu \s_0)^{-\nf 12}e,(\mu \s_0)^{\nf 12}e\,; \a+\h)
	&
	\leq
	6 \mu^{-1} J(z+\cu_{-h},\s_0^{-\nf 12}e,\s_0^{\nf 12}e\,; \a)
	\notag \\ & \qquad
	+ 4\mu^{-1} ( \mu -1 )^2
	\s_0 \lambda_{\nf12,2}^{-1}(z+\cu_{-h}\,; \a)
	\notag \\ & \qquad
	+ C \mu^{-1} \s_0^{-1}
	\lambda_{\nf 38,2}^{-1}(z+\cu_{-h} \,; \a)
	3^{-2h} \| \h \|_{\underline{W}^{1,\infty}(z+\cu_{-h})}^2
	\notag \\ & \qquad
	+
	C
	\lambda_{\nf 38,2}^{-2}(z+\cu_{-h} \,; \a)
	3^{-4h}  \| \nabla  \h  \|_{\underline{W}^{1,\infty}(z+\cu_{-h})}^2
	\,.
\end{align*}
By~\eqref{e.another.way.to.say.h}, for every~$t \in (0,1]$,
\begin{multline*}
	\avsum_{z \in 3^{-h}\Zd \cap \cu_0}
	J(z+\cu_{-h},\s_0^{-\nf 12}e,\s_0^{\nf 12}e\,; \a)
	\\
	\leq
	3^{2 th} \mathcal{E}_{t,2,2}(\cu_0;\a,\s_0)^2
	\leq
	\bigl( 1 + C
	\| \nabla  \h   \|_{\underline{W}^{1,\infty}(\cu_0)} \lambda_{s,2}^{-1}
	(\cu_0; \a )
	\bigr)^{\! \frac{2 t}{1-2s} }
	\mathcal{E}_{t,2,2}(\cu_0;\a,\s_0)^2
	\,.
\end{multline*}
Using~\eqref{e.another.way.to.say.h} again, we get
\begin{align*}
	\mu^{-1}( \mu -1)^2
	\s_0 \lambda_{\nf12,2}^{-1}(z+\cu_{-h}\,; \a)
	&
	\leq
	\mu^{-1}( \mu -1)^2
	\s_0
	3^{2 s h}
	\lambda_{s,2}^{-1}(\cu_0 \,; \a)
	\notag \\ &
	\leq
	\mu^{-1}( \mu -1)^2
	\bigl( 1 +
	C \| \nabla  \h   \|_{\underline{W}^{1,\infty}(\cu_0)} \lambda_{s,2}^{-1}
	(\cu_0; \a )
	\bigr)^{\frac{2 s}{1-2s}}
	\s_0
	\lambda_{s,2}^{-1}(\cu_0 \,; \a)
\end{align*}
and, similarly, using Young's inequality and~\eqref{e.another.way.to.say.h},
\begin{align*}
	\lefteqn{
	\avsum_{z \in 3^{-h}\Zd \cap \cu_0}
	\mu^{-1} \s_0^{-1}
	\lambda_{\nf 38,2}^{-1}(z+\cu_{-h} \,; \a)
	3^{-2h}   \| \h \|_{\underline{W}^{1,\infty}(z+\cu_{-h})}^2
	} \qquad &
	\notag \\ &
	\leq
	\mu^{-1} \s_0^{-1}
	\lambda_{s,2}^{-1}(\cu_0 \,; \a)
	\bigl(
	3^{2 s h}
	\| \h \|_{\underline{L}^2(\cu_0)}^2
	+
	3^{-2(1-s)h}
	\| \nabla  \h  \|_{L^\infty(\cu_0)}^2  \bigr)
	\notag \\ &
	\leq
	\mu^{-1}
	\s_0^{-1}\lambda_{s,2}^{-1} (\cu_0 \,; \a)
	3^{2sh}
	\| \h \|_{\underline{L}^2(\cu_0)}^2
	+
	\frac12\bigl( \mu^{-2} \s_0^{-2} +  \lambda_{s,2}^{-2}(\cu_0;\a) \bigr)
	3^{-2(1-s)h}
	\| \nabla  \h   \|_{L^{\infty}(\cu_0)}^2
	\notag \\ &
	\leq
	C \mu^{-1}
	\s_0^{-1}\lambda_{s,2}^{-1} (\cu_0 \,; \a)
	\bigl( 1 +
	\| \nabla  \h   \|_{\underline{W}^{1,\infty}(\cu_0)} \lambda_{s,2}^{-1}
	(\cu_0; \a )
	\bigr)^{\frac{2 s}{1-2s}}
	\| \h \|_{\underline{L}^2(\cu_0)}^2
	\notag \\ & \qquad
	+
	\mu^{-2} \s_0^{-2}
	\| \nabla  \h   \|_{L^{\infty}(\cu_0)}^2
	+
	C \min\bigl\{ 1 \,,
	\lambda_{s,2}^{-2}
	(\cu_0; \a )
	\| \nabla  \h   \|_{L^{\infty}(\cu_0)}^2
	\bigr\}
\end{align*}
and
\begin{align*}
	\lefteqn{
	\avsum_{z \in 3^{-h}\Zd \cap \cu_0}
	3^{-4 h}  \| \nabla  \h  \|_{\underline{W}^{1,\infty}(z+\cu_{-h})}^2
	\lambda_{\nf 38,2}^{-2}(z+\cu_{-h} \,; \a)
	} \qquad \quad &
	\notag \\ &
	\leq
	3^{-2(1-2s)h}
	\lambda_{s,2}^{-2}(\cu_0 \,; \a)
	\| \nabla  \h   \|_{\underline{W}^{1,\infty}(\cu_0)} ^2
	\leq
	C \min\bigl\{ 1 \,,\,
	\lambda_{s,2}^{-2} (\cu_0; \a )
	\| \nabla  \h   \|_{\underline{W}^{1,\infty}(\cu_0)}^2
	\bigr\}
	\,.
\end{align*}
Combining the above displays yields~\eqref{e.J.sensitivity.no.conditions} and completes the proof.
\end{proof}

\subsection{Homogenization estimates}
\label{ss.blackboxes}

In this subsection we package the coarse-graining machinery into a single homogenization-type estimate. Given a solution~$u$ of the variable-coefficient equation and the corresponding solution~$v$ of the constant-coefficient problem in a large cube~$\cu_m$ with the same right-hand side and boundary data, we show that their gradients are close in~$\Wminusul{-s}{p}$ and their fluxes in~$\underline{W}^{-s,p}$. The error is expressed purely in terms of the coarse-grained ellipticity quantities and the regularity of the (divergence-form) forcing terms. This estimate will be used as a black-box homogenization input in the proof of Theorem~\ref{t.homogenization}.

\begin{proposition}[General coarse-graining estimate]
\label{p.general.coarse.graining}
Let~$m,n\in\Z$ with~$n<m$,~$\a:\cu_m \to \R^{d\times d}_+$ be a uniformly elliptic coefficient field and~$\s_0 \in (0,\infty)$. Fix exponents~$p \in [2,\infty)$ and~$s, s_1, s_2 \in (0,1)$ such that~$s_1 < s < s_2 < 1$ and~$sp' < 1$, where~$p' \coloneqq p/(p-1)$. There exists a constant~$C(p,d)<\infty$ such that, for every vector field~$\g \in W^{s_2,p}(\cu_m;\Rd)$ and pair~$u,v\in H^1(\cu_m)$ satisfying
\begin{equation*}
	\left\{
	\begin{aligned}
		& -\nabla \cdot \a \nabla u  = \nabla \cdot \g  & \mbox{in} & \ \cu_m\,,
		\\
		& -\nabla \cdot \s_0 \nabla v = \nabla \cdot \g  & \mbox{in} & \ \cu_m\,,
		\\
		&
		u-v \in H_0^1(\cu_m)\,,
	\end{aligned}
	\right.
\end{equation*}
we have the estimate
\begin{align}
	\lefteqn{
	3^{-ms} \s_0 \|  \nabla u - \nabla v \|_{\Wminusul{-s}{p} (\cu_m)}
	+
	3^{-ms} \| \a \nabla u - \s_0 \nabla v \|_{\underline{W}^{-s,p} (\cu_m)}
	}
	\qquad &
	\notag \\ &
	\leq
	C s^{-1} \s_0^{\nf12}
	\mathcal{E}_{s_1, \infty, 1} (\cu_m, n;\a,\s_0)
	\biggl(
	\sum_{k=-\infty}^n
	3^{-(s-s_1) p (n-k)}
	\avsum_{z\in 3^k\Zd\cap\cu_m}
	\| \s^{\nf12} \nabla u \|_{\underline{L}^2(z+\cu_k)}^{p}
	\biggr)^{\!\nf 1p}
	\notag \\ & \qquad
	+
	Cs^{-\nf 9 2}
	(s_2 - s)^{-1}
	\bigl(1 + \mathcal{E}^{2}_{\nf{s_1}{2}, \infty, 2}(\cu_m, n;\a,\s_0) \bigr)
	3^{s_2 n}
	[ \g ]_{\underline{W}^{s_2,p}(\cu_m)}
	\,.
	\label{e.general.coarse.graining.estimate}
\end{align}
\end{proposition}

The role of the parameter~$n$ in Proposition~\ref{p.general.coarse.graining} is to provide a mesoscopic scale slightly below~$m$, from which we can take advantage of the regularity of~$\g$. Note that~$\mathcal{E}$ provides no smallness for the last term on the right of~\eqref{e.general.coarse.graining.estimate}. In practice, it is small because the factor
\begin{equation*}
	3^{s_2 n}
	[ \g ]_{\underline{W}^{s_2,p}(\cu_m)}
	=
	3^{-s_2 (m-n)} \cdot 3^{s_2 m} [ \g ]_{\underline{W}^{s_2,p}(\cu_m)}
\end{equation*}
and the smallness comes from the factor of~$3^{-s_2 (m-n)}$. If~$\g=0$, then this term vanishes and the proof of Proposition~\ref{p.general.coarse.graining} applies verbatim with~$n=m$, giving
\begin{align}
	\lefteqn{
	3^{-ms} \s_0 \|  \nabla u - \nabla v \|_{\Wminusul{-s}{p} (\cu_m)}
	+
	3^{-ms} \| \a \nabla u - \s_0 \nabla v \|_{\underline{W}^{-s,p} (\cu_m)}
	}
	\qquad &
	\notag \\ &
	\leq
	C s^{-1} \s_0^{\nf12}
	\mathcal{E}_{s_1, \infty, 1} (\cu_m; \a,\s_0)
	\biggl(
	\sum_{k=-\infty}^m
	3^{-(s-s_1) p (m-k)}
	\avsum_{z\in 3^k\Zd\cap\cu_m}
	\| \s^{\nf12} \nabla u \|_{\underline{L}^2(z+\cu_k)}^{p}
	\biggr)^{\!\nf 1p}
	\,.
	\label{e.general.coarse.graining.estimate.simped}
\end{align}
If we specialize further to the case~$p=2$ and~$\g=0$, then the last factor on the right of~\eqref{e.general.coarse.graining.estimate.simped} is bounded by~$C(s-s_1)^{-\nf12} \| \s^{\nf12} \nabla u \|_{\underline{L}^2(\cu_m)}$. Taking~$s_1=\nf s2$ and weakening the flux norm to its hatted counterpart gives
\begin{multline}
	3^{-ms} \s_0 \|  \nabla u - \nabla v \|_{\Hminusuls{-s}  (\cu_m)}
	+
	3^{-ms} \| \a \nabla u - \s_0 \nabla v \|_{\Hminusuls{-s}(\cu_m)}
	\\
	\leq
	C s^{-\nf32}
	\s_0^{\nf12}
	\mathcal{E}_{\nf s2, \infty, 1} (\cu_m;\a,\s_0)
	\| \s^{\nf12} \nabla u \|_{\underline{L}^2(\cu_m)}
	\,.
	\label{e.general.coarse.graining.estimate.simped.again}
\end{multline}

To prove Proposition~\ref{p.general.coarse.graining}, we use a simple, deterministic duality argument to control the global difference~$u-v$ in terms of the flux defect~$(\a-\s_0)\nabla u$, which is presented in the following lemma. Previously, in Lemma~\ref{l.gen.coarse.grain}, we have estimated the flux defect in terms of the right side of~\eqref{e.general.coarse.graining.estimate}. The combination of these lemmas immediately yields the proposition.

\begin{lemma}
\label{l.cz.flux.estimates}
Let~$m,n\in\Z$ with~$n<m$ and~$\a:\cu_m \to \R^{d\times d}_+$ be a uniformly elliptic coefficient field and~$\s_0 \in (0,\infty)$. Let~$s \in (0,1)$ and~$p \in [2,\infty)$ with~$sp' < 1$, where~$p' \coloneqq p/(p-1)$. There exists a constant~$C(p,d)<\infty$ such that, for every pair~$u,v\in H^1(\cu_m)$ satisfying
\begin{equation*}
	\left\{
	\begin{aligned}
		& \nabla \cdot ( \a \nabla u  - \s_0 \nabla v)  = 0  & \mbox{in} & \ \cu_m\,,
		\\
		&
		u-v \in H_0^1(\cu_m)\,,
	\end{aligned}
	\right.
\end{equation*}
we have
\begin{multline}
	\s_0 \|  \nabla u - \nabla v \|_{\Wminusul{-s}{p} (\cu_m)}
	+
	\| \a \nabla u - \s_0 \nabla v \|_{\underline{W}^{-s,p} (\cu_m)}
	\\
	\leq
	C  3^{s(m-n)}\biggl( \avsum_{z\in 3^n\Zd \cap \cu_m}
	\| (\a-\s_0) \nabla u \|_{\underline{W}^{-s,p} (z+\cu_n)}^p
	\biggr)^{\!\nf1p}
	\label{e.cz.flux.estimates}
	\,.
\end{multline}
\end{lemma}
\begin{proof}
We first prove the bound for gradients:
\begin{equation}
	\s_0 \|  \nabla u - \nabla v \|_{\Wminusul{-s}{p} (\cu_m)}
	\leq
	C  3^{s(m-n)}
	\biggl(
	\avsum_{z \in 3^n \Zd \cap \cu_m}
	\bigl\|  (\a - \s_0)\nabla u \bigr\|_{\underline{W}^{-s,p}(z+\cu_n)}^p
	\biggr)^{\!\nf1p}
	\,.
	\label{e.by.duality.we.show}
\end{equation}
By duality, we have
\begin{equation*}
	\|  \nabla u - \nabla v \|_{\Wminusul{-s}{p} (\cu_m)}
	=
	\sup \biggl\{
	\fint_{\cu_m} (\nabla u  -\nabla v) \cdot \h
	\,:\, \h\in C^\infty(\cu_m ;\Rd)\,, \
	\| \h \|_{\underline{W}^{s,p'} (\cu_m)}^{p'} \leq 1
	\biggr\}
	\,.
\end{equation*}
Fix~$\h \in C^\infty(\cu_m ;\Rd)$ satisfying
\begin{equation}
	\| \h \|_{\underline{W}^{s,p'} (\cu_m)}^{p'}
	=
	3^{-msp'} \| \h \|_{\underline{L}^{p'}(\cu_m)}^{p'}
	+
	[ \h ]_{ \underline{W}^{s,p'} (\cu_m)}^{p'}
	\leq 1
	\label{e.f.normalization}
\end{equation}
and let~$w \in H^1(\cu_m )$ be the solution of the Dirichlet problem
\begin{equation*}
	\left\{
	\begin{aligned}
		& -\nabla \cdot \s_0\nabla w = \nabla \cdot \h
		&  \mbox{in} & \ \cu_m \,, \\
		& w = 0 &  \mbox{on} & \ \partial \cu_m \,.
	\end{aligned}
	\right.
\end{equation*}
By the global Calder\'on-Zygmund~$W^{1,p'}$ estimate, applied in the macroscopic cube~$\cu_m$, we have that
\begin{equation}
	\| \nabla w \|_{\underline{L}^{p'}(\cu_m)}
	\leq
	C\s_0^{-1} \| \h \|_{\underline{L}^{p'}(\cu_m)}
	\leq
	C\s_0^{-1} 3^{sm}
	\,.
	\label{e.CZ.Lp.app.nablaw}
\end{equation}
By the local and global Calder\'on-Zygmund~$W^{1+s,p'}$ estimates, we also obtain
\begin{align}
	[ \nabla w ]_{\underline{W}^{s,p'} (\cu_m)}
	&
	\leq
	C
	\s_0^{-1} [ \h ]_{\underline{W}^{s,p'} (\cu_m)}
	\leq
	C\s_0^{-1}
	\,.
	\label{e.CZ.Wsp.app.nablaw}
\end{align}
Here we use the standard fractional Calder\'on--Zygmund estimate for the Dirichlet problem on a cube; it follows, for example, by a compatible extension and localization argument from the whole-space Fourier multiplier estimates in~\cite[Chapter~V]{Stein}.
Testing the equation~$-\s_0\Delta w=\h$ with~$u-v$ and the difference of the equations for~$u$ and~$v$ with~$w$ yields
\begin{equation}
	\fint_{\cu_m} (\nabla u  -\nabla v) \cdot \h
	=
	- \fint_{\cu_m}
	\nabla w \cdot
	\s_0( \nabla u  - \nabla v)
	=
	\fint_{\cu_m} \nabla w\cdot (\a - \s_0)\nabla u \,.
	\label{e.duality.testing}
\end{equation}
We split the integral on the right side of~\eqref{e.duality.testing} into mesoscopic cubes:
\begin{align}
	\fint_{\cu_m}
	\nabla w
	\cdot
	(\a - \s_0)\nabla u
	&
	=
	\avsum_{z \in 3^n \Zd \cap \cu_m}
	\fint_{z+\cu_n} (\a - \s_0)\nabla u \cdot \bigl( \nabla w - (\nabla w)_{z+\cu_n} \bigr)
	\notag \\ & \qquad \qquad
	+
	\avsum_{z \in 3^n \Zd \cap \cu_m}
	\bigl( (\a - \s_0)\nabla u \bigr)_{z+\cu_n}  \cdot (\nabla w)_{z+\cu_n}
	\,.
	\label{e.Wminusul.subadditivity}
\end{align}
Using H\"older's inequality,~\eqref{e.f.normalization},~\eqref{e.CZ.Lp.app.nablaw} and~\eqref{e.CZ.Wsp.app.nablaw}, together with~\eqref{e.Hs.superadditivity}, we estimate the first sum by
\begin{align*}
	\lefteqn{
	\avsum_{z \in 3^n \Zd \cap \cu_m}
	\fint_{z+\cu_n} (\a - \s_0)\nabla u \cdot \bigl( \nabla w - (\nabla w)_{z+\cu_n} \bigr)
	} \qquad &
	\notag \\ &
	\leq
	\avsum_{z \in 3^n \Zd \cap \cu_m}
	\|  (\a - \s_0)\nabla u \|_{\underline{W}^{-s,p}(z+\cu_n)}
	[ \nabla w ]_{\underline{W}^{s,p'} (z+\cu_n)}
	\notag \\ &
	\leq
	\biggl(
	\avsum_{z \in 3^n \Zd \cap \cu_m}
	\|  (\a - \s_0)\nabla u \|_{\underline{W}^{-s,p}(z+\cu_n)} ^p
	\biggr)^{\!\nf1p}
	\biggl(
	\avsum_{z \in 3^n \Zd \cap \cu_m}
	[ \nabla w ]_{\underline{W}^{s,p'} (z+\cu_n)}^{p'}
	\biggr)^{\!\nf1{p'}}
	\notag \\ &
	\leq
	C
	\biggl(
	\avsum_{z \in 3^n \Zd \cap \cu_m}
	\| (\a - \s_0)\nabla u \|_{\underline{W}^{-s,p}(z+\cu_n)}^p
	\biggr)^{\!\nf1p}
	[ \nabla w ]_{\underline{W}^{s,p'} (\cu_m)}
	\notag \\ &
	\leq
	C \s_0^{-1}
	\biggl(
	\avsum_{z \in 3^n \Zd \cap \cu_m}
	\| (\a - \s_0)\nabla u \|_{\underline{W}^{-s,p}(z+\cu_n)}^p
	\biggr)^{\!\nf1p} \,,
\end{align*}
and the second sum by
\begin{align*}
	\avsum_{z \in 3^n \Zd \cap \cu_m}
	\bigl( (\a - \s_0)\nabla u \bigr)_{z+\cu_n}  \cdot (\nabla w)_{z+\cu_n}
	&
	\leq
	\biggl(
	\avsum_{z \in 3^n \Zd \cap \cu_m}
	\bigl| \bigl( (\a - \s_0)\nabla u \bigr)_{z+\cu_n}  \bigr|^p
	\biggr)^{\!\nf1p}
	\| \nabla w\|_{\underline{L}^{p'}(\cu_m)}
	\notag \\ &
	\leq
	C\s_0^{-1} 3^{sm}
	\biggl(
	\avsum_{z \in 3^n \Zd \cap \cu_m}
	\bigl| \bigl( (\a - \s_0)\nabla u \bigr)_{z+\cu_n}  \bigr|^p
	\biggr)^{\!\nf1p} \,.
\end{align*}
Combining the above yields~\eqref{e.by.duality.we.show}.

\smallskip

We next estimate the fluxes. We observe by the triangle inequality that
\begin{equation*}
	\| \a \nabla u - \s_0\nabla v \|_{\underline{W}^{-s,p} (\cu_m)}
	\leq
	\| (\a-\s_0) \nabla u  \|_{\underline{W}^{-s,p} (\cu_m)}
	+
	\s_0 \| \nabla u - \nabla v \|_{\Wminusul{-s}{p} (\cu_m)}
	\,.
\end{equation*}
Using the subadditivity~\eqref{e.H.minus.s.subadditivity} of the full~$\| \cdot  \|_{\underline W^{-s,p} (\cu_m)}$ norm,
\begin{equation*}
	\| (\a-\s_0) \nabla u  \|_{\underline{W}^{-s,p} (\cu_m)}
	\leq
	C3^{s(m-n)}\biggl( \avsum_{z\in 3^n\Zd \cap \cu_m}
	\| (\a-\s_0) \nabla u \|_{\underline{W}^{-s,p} (z+\cu_n)}^p
	\biggr)^{\!\nf1p} \,.
\end{equation*}
The previous two displays and~\eqref{e.by.duality.we.show} yield
\begin{equation*}
	\| \a \nabla u - \s_0\nabla v \|_{\underline{W}^{-s,p} (\cu_m)}
	\leq
	C  3^{s(m-n)}\biggl( \avsum_{z\in 3^n\Zd \cap \cu_m}
	\| (\a-\s_0) \nabla u \|_{\underline{W}^{-s,p} (z+\cu_n)}^p
	\biggr)^{\!\nf1p}
	\,.
\end{equation*}
This completes the proof of the estimate for the fluxes, and thus of~\eqref{e.cz.flux.estimates}.
\end{proof}

\begin{proof}[{Proof of Proposition~\ref{p.general.coarse.graining}}]
Combine Lemmas~\ref{l.gen.coarse.grain} and~\ref{l.cz.flux.estimates}.
\end{proof}

\section{Renormalization group flow}
\label{s.RG.flow}

In this section we prove the main coarse-graining estimate (stated below in Proposition~\ref{p.induction.bounds}) which asserts that, at each triadic scale~$3^m$, the infrared-truncated operator~$-\nabla\cdot\a_m\nabla$ is quantitatively close to the constant-coefficient Laplacian~$-\shom_m\Delta$.

\smallskip

Throughout the section (and in contrast to the previous one), we always assume~$\P$ is a probability measure satisfying~\ref{a.j.frd},~\ref{a.j.reg},~\ref{a.j.iso} and~\ref{a.j.nondeg}; the parameter~$\nu>0$, the antisymmetric matrix-valued fields~$\k$,~$\mathbf{j}_n$ and the coefficient field~$\a \coloneqq \nu\Id + \k$ are as defined in Section~\ref{ss.assumptions}.

\smallskip

For each integer~$m\in\Z$, we introduce the \emph{infrared cutoffs} of~$\k$ and~$\a$ at scale~$3^m$ by
\begin{equation}
	\k_m(x) \coloneqq  \sum_{n=-\infty}^{m} \mathbf{j}_n(x)\,,
	\qquad
	\a_m(x) \coloneqq  \nu \Id + \k_m(x)
	\,.
	\label{e.infrared.cutoff.def}
\end{equation}
The field~$\k_m$ is locally bounded,~$\P$-almost surely, so~$\a_m$ is locally bounded. Moreover, since~$\k_m$ is antisymmetric, the symmetric part of~$\a_m$ is~$\nu\Id$, and therefore~$\a_m$ is uniformly elliptic (with ellipticity constant~$\nu$) on every bounded domain. Moreover, the moments of~$\k_m$ and hence~$\a_m$ are finite, and these fields are~$\Rd$-stationary with range of dependence of order~$3^m$. It follows from standard qualitative homogenization results and~\ref{a.j.iso} that the operator~$-\nabla \cdot \a_m \nabla$ homogenizes, in the infinite scale limit, to~$-\shom_m\Delta$ for some~$\shom_m\in (0,\infty)$. (We make no use of this result, except as a convenient way to define~$\shom_m$.)

\smallskip

The following proposition asserts that, at scale~$3^m$, the operator~$-\nabla\cdot \a_m \nabla$ is quantitatively close to~$-\shom_m \Delta$. The discrepancy is measured by the random variable~$\mathcal{E}_{s,\infty,2}(\cu_m;\a_m,\shom_m)$ defined in the previous section (see Definition~\ref{d.mathcal.E}), which is shown to have a typical size of~$C\cgamma^{\nf12}$ and stretched exponential tails. The effective diffusivity~$\shom_m$  is shown to satisfy~$\shom_m \approx ( \nu^2 + \cstar \cgamma^{-1} 3^{2\cgamma m} )^{\nf12}$, up to a relative error of order~$C\cstar^{-2} \cgamma^{\nf12}\left| \log \cgamma\right|$.
This entire section is devoted to its proof.

\begin{proposition}[Self-similar coarse-graining estimate]
\label{p.induction.bounds}
There exists a constant~$C(d) < \infty$ such that, if the parameter~$\cgamma$ satisfies~$\cgamma \leq C^{-10} \cstar^{10}$, then, for every~$m \in \Z$ and~$s \in [8\cgamma,1]$,
\begin{equation}
	\label{e.induction.E.bounds}
	\mathcal{E}_{s,\infty,2}(\cu_m; \a_m, \shom_m)
	\leq
	\O_{\Gamma_{2} }
	\bigl( C\cstar^{-1}  s^{-1} \cgamma^{\nf12}  \bigr) +
	\O_{\Gamma_{\nf12}}
	\bigl( \exp( - C^{-1}\cstar^{3}  \cgamma^{-1}) \bigr)
\end{equation}
and
\begin{equation}
	\label{e.formula.for.shom}
	\bigl|
	\shom_{m}
	-\bigl(\nu^2
	+
	\cstar\cgamma^{-1}3^{2\cgamma m} \bigr)^{\nf12}
	\bigr|
	\leq
	C \cstar^{-2} \cgamma^{\nf12}  \left|\log \cgamma\right|
	\,\shom_m
	\,.
\end{equation}
\end{proposition}

We emphasize that the scale of the cube~$\cu_m$ is the same as the scale of the infrared cutoff in the field~$\a_m$. That is, we compare the operator~$\nabla\cdot \a_m \nabla$ to~$\nabla \cdot \shom_m\nabla$ exactly at the scale of the infrared cutoff, without needing to take advantage of the decorrelation of~$\a_m$ on scales larger than~$3^m$ by considering much larger boxes. In combination with the sensitivity estimates we proved in Section~\ref{ss.sensitivity}, this will allow us to ``put back'' the longer wavelengths~$\mathbf{j}_n$ for~$n>m$, whose contribution on~$\cu_m$ has relative size of order~$\cgamma^{\nf12}$. Thus, although the statement of Proposition~\ref{p.induction.bounds} is phrased in terms of the truncated operator~$\nabla \cdot \a_m \nabla$, it implies similar estimates (keeping~$\cu_m$ fixed) for each operator~$\nabla \cdot \a_L\nabla$ with~$L\geq m$, as well as the original, non-cutoff operator~$\nabla \cdot \a\nabla$.

\smallskip

We prove Proposition~\ref{p.induction.bounds} by induction on the scale parameter~$m$, starting from the base scale~$\mstarstar$ (defined below in~\eqref{e.mstarstar}) and propagating to larger~$m$. The induction hypothesis we need to propagate is given in the next definition.

\begin{definition}[Induction hypothesis]
\label{d.mathcalS.def}
For each~$m_0 \in \Z$ and~$E \in [1,\infty)$, we let~$\mathcal{S}(m_0,E)$ refer to the statement that the following two assertions are valid:
\begin{itemize}
\item Bound on the effective diffusivities: for every~$m \in \Z$ with~$m\leq m_0$,
\begin{equation}
	\label{e.shom.h.bounds}
	\frac14\max\bigl\{  \cstar \cgamma^{-1} 3^{2\cgamma m} \,, \, \nu^2 \bigr\}
	\leq \shom_m^2 \leq
	4 \max\bigl\{ \cstar \cgamma^{-1} 3^{2\cgamma m} \,, \, \nu^2 \bigr\}
	\,;
\end{equation}
\item
Smallness of the homogenization errors: for every~$m \in \Z$ with~$m\leq m_0$
and~$s \in [8\cgamma,1]$,
\begin{equation}
	\mathcal{E}_{s,\infty,2}(\cu_m; \a_m, \shom_m)
	\leq
	\O_{\Gamma_{2} }
	\bigl( E s^{-1} \cgamma^{\nf12}  \bigr) +
	\O_{\Gamma_{\nf12}}
	\bigl( s^{-2}\exp( - E^{-3}\cgamma^{-1}) \bigr)
	\,.
	\label{e.new.induction.for.shom}
\end{equation}
\end{itemize}
\end{definition}

The main step in the proof of Proposition~\ref{p.induction.bounds} is to establish the existence of~$C(d)<\infty$ such that, for every~$m\in\Z$,~$\cgamma \in (0,1)$ and~$E \in[1,\infty)$ satisfying
\begin{equation}
	\label{e.xi.delta1.condition}
	E \geq C \cstar^{-1}
	\qand
	\cgamma
	\leq  E^{-10}
	\,,
\end{equation}
we have the implication
\begin{equation}
	\label{e.propagation.of.indyhyp}
	\mathcal{S}(m-1,E)
	\implies
	\mathcal{S}(m,E) \,.
\end{equation}
In the course of propagating the bounds~\eqref{e.shom.h.bounds} on the effective diffusivities, we will manage to prove that the induction hypothesis actually implies the more precise estimate~\eqref{e.formula.for.shom}.

\smallskip

The proof of~\eqref{e.propagation.of.indyhyp} will consist of the following steps:
\begin{itemize}

\item Upper bounds on the effective ellipticity constants (Section~\ref{ss.coarse.grained.sensitivity}). We establish \emph{improved tail bounds} for the coarse-grained ellipticity constants by a trimmed-sum argument: loosely, we show that the coarse-grained matrices are insensitive to rare ``bad'' subcubes, which allows us to upgrade stretched-exponential integrability at the cost of inflating deterministic constants.

\item Propagation of the effective diffusivity bounds (Section~\ref{ss.propagate.diffusivity.bounds}). We derive an approximate recurrence relation for~$\shom_m$ and integrate it to obtain the sharp asymptotics in~\eqref{e.formula.for.shom} which, in particular, implies~\eqref{e.shom.h.bounds}.

\item Propagation of the homogenization estimate (Section~\ref{ss.homogenization.step}). We prove that if the homogenization error at a particular scale is small, then the error contracts by an algebraic factor after moving up sufficiently many scales above the infrared truncation (see Proposition~\ref{p.homogenization.step}). This contraction is the mechanism that ultimately compensates for the various inflationary errors introduced elsewhere.

\item Injection of the fresh shell (Section~\ref{ss.buckle}). We close the induction loop by reinserting the layers~$\mathbf{j}_{m-h+1},\ldots,\mathbf{j}_m$ via the coarse-grained sensitivity estimates from Section~\ref{s.coarse.graining.theory}, showing that their cumulative effect is of relative size~$O(\cgamma^{\nf12} )$.
\end{itemize}

We introduce two distinguished small scales~$\mstar,\mstarstar\in\Z$ with~$\mstarstar < \mstar$ which depend on~$\nu$ and from which we will begin our induction. They are defined as follows:
\begin{itemize}

\item We let~$\mstarstar$ be the largest integer~$m$ such that~$\nu^2 \geq 2 (\log 3) \cstar^{+}  \cgamma^{-3} 3^{2\cgamma m}$, that is,
\begin{equation}
	\label{e.mstarstar}
	\mstarstar  \coloneqq
	\Bigl\lfloor (2\cgamma)^{-1}
	\log_3 \Bigl( \frac{\nu^2 \cgamma^3}{2 (\log3) \cstar^{+}} \Bigr) \Bigr\rfloor
	\,.
\end{equation}

\item We let~$\mstar$ be the largest integer~$m$ such that~$\nu^2 \geq 2 (\log 3) \cstar^{+}  \cgamma^{-1} 3^{2\cgamma m}$, that is,
\begin{equation}
	\label{e.mstar}
	\mstar \coloneqq
	\Bigl\lfloor (2\cgamma)^{-1}
	\log_3 \Bigl( \frac{\nu^2 \cgamma}{2 (\log3) \cstar^{+}} \Bigr) \Bigr\rfloor
	\,.
\end{equation}

\end{itemize}
Notice that the gap between these scales is~$\mstar - \mstarstar \approx \cgamma^{-1} \left|\log_3 \cgamma\right|$.
We should think of~$\mstarstar$ as the smallest scale at which the vector field is active, and~$\mstar$ as the first scale at which the effective diffusivity~$\shom_m$ begins to increase exponentially in~$m$.

\smallskip

To establish the base case for the induction, we essentially argue that the antisymmetric part of the coefficient field does not have much effect on small enough scales.

\begin{itemize}

\item Base case (Section~\ref{ss.basecase}). We show that
$\mathcal{S}(\mstarstar,C_{\mathrm{start}}\cstar^{-\nf12})$ is valid for a
constant~$C_{\mathrm{start}}(d)<\infty$.
\end{itemize}

We begin in the next subsection by collecting some basic probabilistic estimates on the truncated stream matrices~$\k_m$ and their increments~$\k_m-\k_n$ that will be used throughout the section.

\subsection{Stream matrix estimates}
\label{ss.cutoffs}

In this subsection, we collect basic estimates on norms of the field~$\k_m$ and its increments~$\k_m-\k_n$.
We also introduce the infrared cutoff of the~$2d\times 2d$ block matrix~$\bfA$, defined by
\begin{equation}
	\label{e.bfA.ell.def}
	\bfA_m(x)  \coloneqq
	\begin{pmatrix}
		( \nu\Id + \nu^{-1}\k_m^t\k_m )(x)
		& -\nu^{-1}\k_m^t(x)
		\\ - \nu^{-1}\k_m (x)
		& \nu^{-1}\Id
	\end{pmatrix}
	\,.
\end{equation}
For every bounded Lipschitz domain~$U \subset \Rd$, we let~$\s_m(U)$,~$\s_{m,*}(U)$,~$\k_m(U)$ and~$\b_m(U)$ be the associated coarse-grained matrices for the field~$\a_m$. By~\ref{a.j.frd}, the fields~$\a_m$,~$\k_m$ and~$\bfA_m$ are~$\Rd$-stationary and have range of dependence~$\sqrt{d}3^m$.

\smallskip

We will often encounter the sum~$\sum_{k=n+1}^{m} 3^{2\cgamma k}$ throughout the rest of the paper. To aid our computations we introduce the constant
\begin{equation}
	K_{\cgamma}
	\coloneqq \frac{2\log 3}{1-3^{-2\cgamma}}
	\,,
	\label{e.K.cgamma.def}
\end{equation}
so that
\begin{equation}
	2(\log 3) \! \sum_{j=n+1}^m
	3^{2\cgamma j}
	=
	K_{\cgamma}
	\bigl(
	3^{2\cgamma m} - 3^{2\cgamma n}
	\bigr)
	\,.
	\label{e.K.cgamma.sum.identity}
\end{equation}
One may check that
\begin{equation}
	0\leq  K_{\cgamma} - \cgamma^{-1} \leq 4
	\,.
	\label{e.K.cgamma.ineq}
\end{equation}
By~\eqref{e.K.cgamma.sum.identity} and~\eqref{e.K.cgamma.ineq} we see that, for every~$m,n\in\Z$ with~$n<m$,
\begin{equation}
	2(\log 3) \! \sum_{j=n+1}^m
	3^{2\cgamma j}
	\leq
	\min \bigl\{ K_{\cgamma} \,, 2(\log 3) (m-n) \bigr\} 3^{2\cgamma m}
	\leq
	C \min \bigl\{ \cgamma^{-1}  \,, m-n \bigr\} 3^{2\cgamma m}
	\,.
	\label{e.geo.part.sum}
\end{equation}
We also have the more precise bound
\begin{equation}
	\label{e.simple.summing.with.gamma}
	\biggl| 2(\log 3) \! \sum_{j=n+1}^m 3^{2\cgamma j}
	-
	\cgamma^{-1} \bigl( 3^{2\cgamma m} - 3^{2\cgamma n} \bigr)
	\biggr|
	\leq
	4\bigl( 3^{2\cgamma m} - 3^{2\cgamma n} \bigr)
	\leq
	C \min \bigl\{ 1 \,, \cgamma (m-n) \bigr\} 3^{2\cgamma m}
	\,.
\end{equation}
Observe that, by~\eqref{e.simple.summing.with.gamma}, the scaling and
stationarity in~\eqref{e.diff.law.shift}, the independence of the
fields~$\mathbf{j}_n$, and the definition of~$\cstar^+$, we have
\begin{equation}
	\Bigl|
	\E \biggl[\fint_{\cu_l} |\k_m(x)|^2\,dx\biggr]
	-
	\frac12 \cstar^{+} \cgamma^{-1} 3^{2\cgamma m}
	\Bigr|
	\leq
	2 \cstar^{+} 3^{2\cgamma m}
	\qquad \forall l \in \Z \,\,.
	\label{e.km.L2.expectation}
\end{equation}

We next record some estimates for the volume-normalized~$L^p$ norms of the stream matrix increments~$\k_m-\k_n$; note that we can obtain estimates for~$\k_m$ by sending~$n\to -\infty$.

\begin{lemma}[{Estimates for stream matrix increments}]
\label{l.km.kn.Lp.estimates}
There exists a constant~$C(d)<\infty$ such that, for every~$p\in[2,\infty)$ and~$l,m,n \in \Z$ with~$n < m \leq l$,
\begin{equation}
	3^{-2\cgamma m}  \bigl\|\k_m-\k_n\bigr\|_{\underline L^2(\cu_l)}^2
	=
	\cstar^{+} (\log3)
	\sum_{k=n+1}^m
	3^{2\cgamma (k-m)}
	+
	\O_{\Gamma_1}\bigl(
	C
	3^{-\frac d2(l-m)}
	\bigr)
	\,,
	\label{e.km.kn.L2.exact}
\end{equation}
\begin{equation}
	\label{e.km.kn.Lp}
	3^{-2\cgamma m}  \bigl\| \k_m - \k_n \bigr\|_{\underline{L}^{p}(\cu_l)}^2
	\leq
	C p\min\{ \cgamma^{-1} ,  m-n \}
	+ \O_{\Gamma_{1}}\bigl(  C p \min\{ \cgamma^{-1} ,  (m-n) \} 3^{-\frac{d}{p}(l-m)}\bigr)
	\,,
\end{equation}
\begin{equation}
	3^{-2\cgamma m} \bigl\| \k_m - \k_n \bigr\|_{L^\infty(\cu_n)}^2
	\leq \O_{\Gamma_1}\bigl(  C \min\{ \cgamma^{-1} , m-n \}  \bigr)
	\label{e.km.kn.Linfty.smallcube}
\end{equation}
and
\begin{equation}
	3^{-2\cgamma m}  \bigl\| \k_m - \k_n \bigr\|_{{L}^{\infty}(\cu_l)}^2
	\leq
	\O_{\Gamma_1} \bigl(
	C (l-n) \min\{ \cgamma^{-1} , m-n \}
	\bigr)\,.
	\label{e.km.kn.Linfty}
\end{equation}
For every~$m,n\in\Z$ with~$n\leq m$, we also have
\begin{equation}
	3^{-2\cgamma m}\bigl\|\k_m-\k_n\bigr\|_{\underline L^4(\cu_m)}^2
	\leq C(m-n)+\O_{\Gamma_1}(C)
	\,.
	\label{e.km.kn.L4.centered}
\end{equation}
For every~$l,m,n \in \Z$ with~$n\leq m$,
\begin{equation}
	\label{e.W.1.inf.bound}
	\sum_{k=n+1}^m
	\bigl(
	3^{k}  \|  \nabla \mathbf{j}_k \|_{L^{\infty}(\cu_l)}
	+
	3^{2 k}  \|  \nabla^2 \mathbf{j}_k \|_{L^{\infty}(\cu_l)}
	\bigr)
	\leq
	\O_{\Gamma_2}
	\bigl(
	C \max\{1,l-n\}^{\nf 12}
	\min\bigl\{ \cgamma^{-1} , m-n \bigr\} 3^{\cgamma m}
	\bigr)
	\,.
\end{equation}
\end{lemma}
\begin{proof}
By~\eqref{e.diff.law.shift} and~\ref{a.j.reg}, for every~$k\in\Z$,
\begin{equation}
	\| \mathbf{j}_k \|_{L^{\infty}(\cu_k)} = \O_{\Gamma_2}(3^{\cgamma k})
	\label{e.jk.O}
\end{equation}
and
\begin{equation}
	\| \nabla \mathbf{j}_k \|_{L^{\infty}(\cu_k)}
	+
	3^k \| \nabla^2 \mathbf{j}_k \|_{L^{\infty}(\cu_k)} = \O_{\Gamma_2}(C3^{(\cgamma -1) k})
	\,.
	\label{e.nabla.jk.O}
\end{equation}
By the triangle inequality~\eqref{e.Gamma.sigma.triangle} and~\eqref{e.nabla.jk.O}, for every~$m,n\in\Z$ with~$n< m$,
\begin{align}
	\bigl\| \nabla (\k_m - \k_n) \bigr\|_{L^\infty(\cu_n)}
	+
	3^n\bigl\| \nabla^2 (\k_m - \k_n) \bigr\|_{L^\infty(\cu_n)}
	&
	\leq
	\O_{\Gamma_2}
	\bigl(
	C 3^{(\cgamma-1) n}
	\bigr)
	\,.
	\label{e.nabla.km.kn.infty}
\end{align}
By~\eqref{e.jk.O}, the fields~$\mathbf{j}_k(0)$ are independent and have mean zero. Proposition~\ref{p.concentration} therefore gives
\begin{align}
	|(\k_m - \k_n) (0)| &\leq
	\O_{\Gamma_2}\biggl( C \biggl( \sum_{k=n+1}^m 3^{2 \cgamma k} \biggr)^{\!\nf 12}\, \biggr)
	&
	=
	\O_{\Gamma_2}\bigl(
	C\cgamma^{-\nf12 } ( 3^{2\cgamma m} - 3^{2\cgamma n} )^{\nf12}
	\bigr)
	\notag \\ &
	\leq
	\O_{\Gamma_2}
	\bigl( C \min\bigl\{ \cgamma^{-\nf 12} ,  (m-n)^{\nf 12} \bigr\} 3^{\cgamma m} \bigr)
	\,,
	\label{e.k.ell.upscales}
\end{align}
where in the last line we used~\eqref{e.geo.part.sum}. Using also~\eqref{e.nabla.km.kn.infty}, we obtain, for every~$m,n\in\Z$ with~$n< m$,
\begin{align}
	\label{e.k.ell.upscales.infty}
	\bigl\| \k_m - \k_n \bigr\|_{L^\infty(\cu_n)}
	&
	\leq
	\sqrt{d}3^n \,\bigl\| \nabla(\k_m - \k_n) \bigr\|_{L^\infty(\cu_n)}
	+\bigl| (\k_m - \k_n) (0)\bigr|
	\notag \\ &
	\leq \O_{\Gamma_2}\bigl(
	C  3^{\cgamma n}
	+
	C \min\{ \cgamma^{-\nf 12} ,  (m-n)^{\nf 12} \} 3^{\cgamma m} \bigr)
	\notag \\ &
	\leq \O_{\Gamma_2}\bigl(  C \min\bigl\{ \cgamma^{-\nf 12} ,  (m-n)^{\nf 12} \bigr\} 3^{\cgamma m}  \bigr)
	\,.
\end{align}
This yields~\eqref{e.km.kn.Linfty.smallcube}. To obtain~\eqref{e.km.kn.Linfty} we combine~\eqref{e.km.kn.Linfty.smallcube} and~\eqref{e.maxy.bound}.

\smallskip

For~$k'<k$, condition on~$\mathbf j_k$, partition~$\cu_l$ into cubes of size~$3^{k'}$, and use a fixed coloring so that the centered averages of~$\mathbf j_{k'}$ within each color class are independent. Proposition~\ref{p.concentration} and the product estimate~\eqref{e.multGammasig} then give the factor~$3^{-\frac d2(l-k')}$ below; the bounded-gap cases follow from the direct product bound. The resulting two gap sums are geometric uniformly in~$\cgamma$, and we compute, for every~$l,m,n\in\Z$ with~$n<m\leq l$,
\begin{align}
	\lefteqn{
	\fint_{\cu_l} |\k_m(x)-\k_n(x)|^2\,dx
	} \qquad &
	\notag \\ &
	=
	\sum_{k,k'=n+1}^m
	\fint_{\cu_l}
	\mathbf{j}_{k}(x) :
	\mathbf{j}_{k'}(x)
	\,dx
	\notag \\ &
	=
	2
	\sum_{k=n+1}^m
	\sum_{k'=n+1}^{k-1}
	\underbrace{
	\fint_{\cu_l}
	\mathbf{j}_{k}(x) :
	\mathbf{j}_{k'}(x)
	\,dx
	}_{
	\leq\,\O_{\Gamma_1}( C 3^{\cgamma (k+k')} 3^{-\frac d2(l-k')}  )
	}
	+
	\sum_{k=n+1}^m
	\underbrace{
	\fint_{\cu_l}
	\mathbf{j}_{k}(x) :
	\mathbf{j}_{k}(x)
	\,dx
	}_{
	=\, \cstar^{+} (\log3)3^{2\cgamma k} + \O_{\Gamma_1}(C 3^{2\cgamma k} 3^{-\frac d2(l-k)})
	}
	\notag \\ &
	=
	\cstar^{+} (\log3)
	\sum_{k=n+1}^m
	3^{2\cgamma k}
	+
	\O_{\Gamma_1}\biggl(
	C
	\sum_{k=n+1}^m
	3^{2\cgamma k -\frac d2(l-k)}
	\biggr)
	\notag \\ &
	=
	\cstar^{+} (\log3)
	\sum_{k=n+1}^m
	3^{2\cgamma k}
	+
	\O_{\Gamma_1}\bigl(
	C
	3^{2\cgamma m -\frac d2(l-m)}
	\bigr)
	\,.
	\label{e.km.kn.L2.bound}
\end{align}
Multiplying by~$3^{-2\cgamma m}$ yields~\eqref{e.km.kn.L2.exact}.

The proof of~\eqref{e.km.kn.L4.centered} is very similar to the proof of~\eqref{e.km.kn.L2.exact}: we expand the fourth power and apply Proposition~\ref{p.concentration} after grouping the centered spatial averages by their largest scale.

By the assumption of~$\Rd$-stationarity and Lemma~\ref{l.Gamma.sigma.triangle}, we obtain from~\eqref{e.k.ell.upscales} that, for every~$l\in\Z$ and~$m,n\in\Z$ with~$n< m$,
\begin{equation}
	\label{e.kmn.bounds}
	\fint_{\cu_l} \bigl| (\k_m - \k_n) (x)\bigr|^p \,dx
	\leq \O_{\Gamma_{\nf 2 p}}\bigl(  (Cp)^{\nf p2} \min\{ \cgamma^{-\nf p2} ,  (m-n)^{\nf p2} \} 3^{\cgamma p  m} \bigr)  \,.
\end{equation}
Using the finite range of dependence assumption~\ref{a.j.frd} and Proposition~\ref{p.concentration} we can improve the previous bound on scales larger than~$3^{m}$: for every~$l,m,n\in \Z$ with~$n< m\leq l$,
\begin{align}
	\label{e.kl.bounds.large}
	\fint_{\cu_l} \bigl| (\k_m - \k_n) (x)\bigr|^p \,dx
	&
	\leq
	(Cp)^{\nf p2} \min\{ \cgamma^{-\nf p2} ,  (m-n)^{\nf p2} \} 3^{\cgamma p m}
	\notag \\ & \qquad
	+
	\O_{\Gamma_{\nf2p}}\Bigl(  (Cp)^{\nf p 2}\bigl( \min\{ \cgamma^{-\nf p2} ,  (m-n)^{\nf p2} \} 3^{\cgamma p m}    \bigr)  3^{-\frac d2(l-m)}\Bigr)
	\,.
\end{align}
Taking the~$2/p$ power yields~\eqref{e.km.kn.Lp}.

\smallskip

Finally, turning to the proof of~\eqref{e.W.1.inf.bound}, we have, by~\eqref{e.nabla.jk.O} and~\eqref{e.maxy.bound}, for every~$k \in \Z$ with~$k < l$,
\begin{equation}
	3^{k}
	\| \nabla \mathbf{j}_k \|_{L^{\infty}(\cu_l)}
	+
	3^{2k}
	\| \nabla^2 \mathbf{j}_k \|_{L^{\infty}(\cu_l)}
	\leq
	\O_{\Gamma_2}
	\bigl(C(l-k)^{\nf12} 3^{\cgamma k}\bigr)
	\,.
	\label{e.W1inf.jL.bound.smaller}
\end{equation}
Summing the above display for~$k \in \Z \cap [n+1,m\wedge l)$ and then summing~\eqref{e.nabla.jk.O} for~$k\in \Z \cap [ m\wedge l,m]$ yields~\eqref{e.W.1.inf.bound}.
\end{proof}

\paragraph{Annealed coarse-grained matrices.}
Associated to the coarse-grained matrices are the \emph{annealed} matrices denoted by
\begin{equation*}
	\shom_m(U)
	\,, \quad
	\shom_{m,*}(U)
	\,, \quad
	\khom_m(U)
	\,, \quad
	\bhom_m(U)
	\,, \quad
\end{equation*}
and defined by
\begin{equation}
	\label{e.meet.the.homs}
	\left\{
	\begin{aligned}
		& \shom_{m,*}(U)  \coloneqq  \E \bigl[ \s_{m,*}^{-1}(U) \bigr]^{-1} \,,
		\\ &
		\khom_m(U)
		\coloneqq  \shom_{m,*}(U)  \E \bigl[ \s_{m,*}^{-1}(U) \k_m(U) \bigr] \,,
		\\ &
		\bhom_m(U)
		\coloneqq
		\E \bigl[ \b_m(U) \bigr]\,,
		\\ &
		\shom_m(U)
		\coloneqq
		\bhom_m(U)
		-
		\khom_m^t(U) \shom_{m,*}^{-1}(U) \khom_m(U)
		\,.
	\end{aligned}
	\right.
\end{equation}
Equivalently, these are defined in such a way that
\begin{equation}
	\label{e.homs.defs.U.pre}
	\bfAhom_m (U) \coloneqq
	\E \bigl[ \bfA_m(U) \bigr]
	=
	\begin{pmatrix}
		( \shom_m + \khom_m^t\shom_{m,*}^{-1}\khom_m )(U)
		& -(\khom_m^t\shom_{m,*}^{-1})(U)
		\\ - ( \shom_{m,*}^{-1}\khom_m )(U)
		& \shom_{m,*}^{-1}(U)
	\end{pmatrix}
	\,.
\end{equation}
By~\eqref{e.adjoint.redundant} and the assumption~\ref{a.j.iso} that~$\a_m$ has the same law as~$\a_m^t$, we deduce that, for every bounded Lipschitz domain~$U\subseteq\Rd$,
\begin{equation}
	\khom_m(U) = 0\,.
	\label{e.annealed.khom.zero}
\end{equation}
In particular~$\bhom_m (U) = \shom_m(U)$ and~$\bfAhom_m(U)$ is block diagonal,
\begin{equation}
	\label{e.homs.defs.U.diag}
	\bfAhom_m (U)
	=
	\begin{pmatrix}
		\shom_m(U)
		& 0
		\\ 0
		& \shom_{m,*}^{-1}(U)
		\end{pmatrix} \,\,.
\end{equation}
The hypercube symmetry assumption~\ref{a.j.iso} also
implies, in the case that~$U = \cu$ is a cube, that each of the diagonal blocks of~$\bfAhom_m(\cu)$ is a scalar matrix.

By subadditivity, the matrices~$\bfAhom_m(\cu_n)$ are monotone nonincreasing in~$n$ in the Loewner (positive semidefinite) order. We also define the deterministic matrices~$\bfAhom_m$ and~$\shom_m$ as the infinite-volume limits of these:
\begin{equation}
	\label{e.homs.defs}
	\bfAhom_m  \coloneqq
	\begin{pmatrix}
		\shom_m
		& 0
		\\ 0
		& \shom_{m}^{-1}
	\end{pmatrix}
	\coloneqq  \lim_{n\to \infty} \E \bigl[ \bfA_m(\cu_n) \bigr]\,.
\end{equation}
The fact that~$\shom_{m,*}^{-1}(\cu_n)$ converges to~$\shom_{m}^{-1}$ as~$n\to \infty$ is a consequence of qualitative homogenization for the field~$\nu\Id + \k_m$ (see~\cite[Proposition D.2 \& Theorem 3.1]{AK.HC}).
The hypercube symmetry assumption ensures that~$\shom_m$ is a scalar matrix.
Throughout, we will consider~$\shom_m$ to be a positive scalar matrix, or a positive real number, whichever is convenient.

\subsection{The base case of the induction}
\label{ss.basecase}

We check the validity of the induction hypothesis for all sufficiently small scales. Recall that~$\mstar$ and~$\mstarstar$ are defined in~\eqref{e.mstar} and~\eqref{e.mstarstar}, above.

\begin{proposition}[Base case]
\label{p.base.case}
There exists~$C_0(d)<\infty$ such that,
for every~$s\in (0,1)$ and~$m\in \Z$,
\begin{equation}
	\nu \Id
	\leq
	\shom_m
	\leq
	\nu \Bigl( 1 +  \frac1{2\nu^2} \cstar^{+}
	\cgamma^{-1}
	( 1+4\cgamma )
	3^{2 \cgamma m}
	\Bigr) \Id
	\label{e.basecase.shom.m}
\end{equation}
and
\begin{align}
	\mathcal{E}_{s,\infty,2}(\cu_m; \a_m,\shom_m)
	&
	\leq
	\O_{\Gamma_2}\bigl(  C_0 \nu^{-1} \cgamma^{-\nf12 } 3^{\cgamma m}
	\bigl(1+6d(\log3)s^{-1}\bigr)^{\nf12} \bigr)
	\notag \\ & \qquad
	+
	\O_{\Gamma_4}\bigl(  C_0 \nu^{-\nf12} \cgamma^{-\nf14} 3^{\frac12 \cgamma m}
	\bigl(1+6d(\log3)s^{-1}\bigr)^{\nf14} \bigr)
	\,.
	\label{e.basecase.mathcal.E}
\end{align}
Consequently, for every~$s \in (0,1)$,~$m \in \Z \cap (-\infty, \mstarstar] $,
\begin{equation}
	\label{e.basecase.homogenization}
	\nu \Id \leq \shom_m \leq (1+\cgamma^2 ) \nu \Id
	\qqand
	\mathcal{E}_{s,\infty,2}(\cu_m\,;\a_m,\shom_m)
	\leq \O_{\Gamma_2}\bigl(C_0(\cstar^{+})^{-\nf12}\cgamma^{\nf12}
	\bigl(1+6d(\log3)s^{-1}\bigr)^{\nf12}\bigr)
	\,,
\end{equation}
and, for every~$s\in(0,1)$ and~$m \in (-\infty, \mstar] \cap \Z$,
\begin{equation}
	\nu \Id  \leq \shom_m  \leq 2 \nu\Id
	\qqand
	\mathcal{E}_{s,\infty,2}(\cu_m\,;\a_m,\shom_m)
	\leq
	2
	+
	\O_{\Gamma_2}
	\bigl(
	C_0(\cstar^{+})^{-\nf12}\bigl(1+6d(\log3)s^{-1}\bigr)^{\nf12}
	\bigr)
	\,.
	\label{e.basecase.diffusivity}
\end{equation}
In particular, there exists~$C(d)<\infty$ such that~$\mathcal{S}(\mstarstar ,C \cstar^{-\nf12} )$ is valid.
\end{proposition}
\begin{proof}
According to~\eqref{e.CG.bounds.1}, for every~$m,l\in\Z$,
\begin{equation}
	\nu \Id \leq \s_{m,*}(\cu_l) \leq
	\s_m(\cu_l)
	\leq
	\b_m(\cu_l  )
	\leq
	\nu \bigl( 1 + \nu^{-2} \| \k_m \|^2_{\underline{L}^2(\cu_l )} \bigr) \Id
	\,.
	\label{e.basecase.tightbound}
\end{equation}
In particular,
\begin{equation}
	\k_m^t (\cu_l) \s_{m,*}^{-1}(\cu_l) \k_m (\cu_l)
	=
	\b_m(\cu_l) -
	\s_m(\cu_l)
	\leq
	\nu^{-1} \| \k_m \|^2_{\underline{L}^2(\cu_l )} \Id
	\,.
	\label{e.basecase.ksk.term}
\end{equation}
Taking an expectation of~\eqref{e.basecase.tightbound} and using~\eqref{e.annealed.khom.zero} and~\eqref{e.km.L2.expectation} yields that, for every~$m,l\in\Z$,
\begin{equation}
	\nu \Id
	\leq
	\shom_{m,*}(\cu_l)
	\leq
	\shom_m
	\leq
	\shom_{m}(\cu_l)
	\leq
	\nu \Bigl( 1 +  \frac1{2\nu^2} \cstar^{+}
	\cgamma^{-1}
	( 1+4\cgamma )
	3^{2 \cgamma m}
	\Bigr) \Id
	\,.
	\label{e.plateau.region.bound}
\end{equation}
This is~\eqref{e.basecase.shom.m}.
Using the definition~\eqref{e.J.coarse.grained} of~$J$ and the bounds~\eqref{e.basecase.tightbound},~\eqref{e.basecase.ksk.term} and~\eqref{e.plateau.region.bound}, we obtain, for every~$e\in\Rd$ with~$|e|=1$,
\begin{align*}
	\lefteqn{
	J(\cu_l,\shom_m^{-\nf12} e , \shom_m^{\nf12} e \,;\a_m)
	} \ \ &
	\notag \\ &
	=
	\frac12 e\cdot \shom_m^{-\nf12} \b_m (\cu_l) \shom_m^{-\nf12} e
	+
	\frac12 e\cdot \shom_m^{\nf12} \s_{m,*}^{-1}(\cu_l) \shom_m^{\nf12} e
	+
	e \cdot \shom_m^{\nf12}\s_{m,*}^{-1}(\cu_l) \k_m (\cu_l) \shom_m^{-\nf12}e -  1
	\notag \\ &
	\leq
	\frac12 \bigl|  \shom_m^{-\nf12} \b_m (\cu_l) \shom_m^{-\nf12} \bigr|
	+
	\frac12 \bigl| \shom_m^{\nf12} \s_{m,*}^{-1}(\cu_l) \shom_m^{\nf12} \bigr|
	+
	\bigl| \shom_m^{\nf12} \s_{m,*}^{-\nf12}(\cu_l) \bigr|
	\bigl|
	(\k_m^t \s_{m,*}^{-1}\k_m )(\cu_l) \shom_m^{-1}\bigr|^{\nf12}  -  1
	\notag \\ &
	\leq
	\frac12 \bigl( 1 + \nu^{-2} \| \k_m \|^2_{\underline{L}^2(\cu_l )}
	\bigr)
	+
	\frac12
	\Bigl( 1 +  \frac1{2} \nu^{-2} \cstar^{+}
	\cgamma^{-1}
	( 1+4\cgamma )
	3^{2 \cgamma m}
	\Bigr)
	\notag \\ & \qquad
	+
	\Bigl( 1 +  \frac1{2}\nu^{-2} \cstar^{+}
	\cgamma^{-1}
	( 1+4\cgamma )
	3^{2 \cgamma m}
	\Bigr)^{\!\nf12}
	\nu^{-1} \| \k_m \|_{\underline{L}^2(\cu_l )}
	- 1
	\notag \\ &
	=
	\frac12 \nu^{-2} \| \k_m \|^2_{\underline{L}^2(\cu_l )}
	+
	\frac1{4}\nu^{-2} \cstar^{+}
	\cgamma^{-1}
	( 1+4\cgamma )
	3^{2 \cgamma m}
	+
	\Bigl( 1 +  \frac1{2}\nu^{-2} \cstar^{+}
	\cgamma^{-1}
	( 1+4\cgamma )
	3^{2 \cgamma m}
	\Bigr)^{\!\nf12}
	\nu^{-1} \| \k_m \|_{\underline{L}^2(\cu_l )}
	\notag \\ &
	\leq
	\nu^{-2} \| \k_m \|^2_{\underline{L}^2(\cu_l )}
	+
	\nu^{-1} \| \k_m \|_{\underline{L}^2(\cu_l )}
	+
	\frac1{2}\nu^{-2} \cstar^{+}
	\cgamma^{-1}
	( 1+4\cgamma )
	3^{2 \cgamma m}
	\,.
\end{align*}
Applying this estimate cube-wise, using~\eqref{e.km.L2.expectation}, and summing the depth gain~$3^{-d(m-l)/2}$ gives
\begin{align*}
	\mathcal{E}_{s,\infty,2}(\cu_m;\a_m,\shom_m)
	&\leq
	\O_{\Gamma_2}\!\left(
	C\nu^{-1}\cgamma^{-\nf12}3^{\cgamma m}
	\bigl(1+6d(\log3)s^{-1}\bigr)^{\nf12}\right)
	\\ &\qquad+
	\O_{\Gamma_4}\!\left(
	C\nu^{-\nf12}\cgamma^{-\nf14}3^{\frac12\cgamma m}
	\bigl(1+6d(\log3)s^{-1}\bigr)^{\nf14}\right)\,,
\end{align*}
which is~\eqref{e.basecase.mathcal.E}.

\smallskip

The statements~\eqref{e.basecase.homogenization} and~\eqref{e.basecase.diffusivity} follow immediately from~\eqref{e.basecase.shom.m} and~\eqref{e.basecase.mathcal.E} by using the definitions of~$\mstarstar$ and~$\mstar$, see~\eqref{e.mstarstar} and~\eqref{e.mstar}. For~$s\in[8\cgamma,1]$, applying~\eqref{e.basecase.homogenization} with~$\min\{s,\frac12\}$, using monotonicity in the exponent, and observing that~$\min\{s,\frac12\}^{-1}\leq 2s^{-1}$ covers the endpoint~$s=1$. Finally, the last statement about the validity of~$\mathcal{S}(\mstarstar ,C\cstar^{-\nf12} )$ is immediate from~\eqref{e.basecase.homogenization}.
\end{proof}

\subsection{Tail bounds for the coarse-grained ellipticity}
\label{ss.coarse.grained.sensitivity}

In this subsection we prove the following statement, which improves the tail bounds of the estimate~\eqref{e.new.induction.for.shom} assumed in the induction hypothesis.

\begin{proposition}[Improved tail bounds]
\label{p.tail.bounds}
Assume~$m \in \Z$ and~$E \in [1,\infty)$ are such that~$\S(m-1, E)$ is valid. There exists~$C(d)<\infty$ such that, if~$s \in [4 \cgamma ,1]$ and
\begin{equation}
	\label{e.percolation.admissibility.tails}
	\max\bigl\{ C , \cstar^{-1} \bigr\}
	\leq
	E
	\leq \cgamma^{-\nf15}
	\,,
\end{equation}
then
\begin{equation}
	\mathcal{E}_{s,\infty,2}
	(\cu_m;\a_m,\shom_{m-1})
	\leq
	C
	+
	\O_{\Gamma_{2}} \bigl( C\cstar^{-\nf12}s^{-1} \cgamma^{\nf12}  \bigr)
	+
	\O_{\Gamma_{\nf 12}}\bigl( \exp \bigl( - C^{-1} E^{-2} \cgamma^{-1}\bigr)  \bigr)\,.
	\label{e.tail.bounds}
\end{equation}
\end{proposition}

Comparing~\eqref{e.new.induction.for.shom} and~\eqref{e.tail.bounds}, we see that the latter has an additional deterministic constant~$C$ on the right side---so it gives a much worse estimate for the \emph{typical} size of~$\mathcal{E}_{s,\infty,2}
(\cu_m;\a_m,\shom_{m-1})$. On the other hand, the other two random terms are better, because they have each replaced one instance of the parameter~$E$ by~$C$. To close the induction, we will assume~$E$ is very large compared to~$C$, so these terms will be better by at least a large constant compared to~\eqref{e.new.induction.for.shom}. Therefore~\eqref{e.tail.bounds} gives a better bound than~\eqref{e.new.induction.for.shom} for the probability of the rare event that~$\mathcal{E}_{s,\infty,2} (\cu_m;\a_m,\shom_{m-1})$ is larger than a fixed large constant. This estimate will be crucial to closing the induction: see the proof of the estimate for~\eqref{e.local.bad.events.summed} in Section~\ref{ss.buckle}, below. The estimates in this subsection will also be needed in our regularity iteration in Sections~\ref{s.regularity}.

\smallskip

As we have seen in Lemma~\ref{l.mathcal.E.to.Lambdas}, the random variable~$\mathcal{E}_{s,\infty,2} (\cu_m;\a_m,\shom_{m-1})$ is closely related to the coarse-grained ellipticity constants. In fact, we will prove a slightly stronger version of Proposition~\ref{p.tail.bounds} which is stated explicitly in terms of the coarse-grained ellipticity constants, rather than~$\mathcal{E}_{s,\infty,2} (\cu_m;\a_m,\shom_{m-1})$.

\begin{proposition}[Tail bounds for coarse-grained ellipticity constants]
\label{p.cg.ellipticity.bounds}
Assume~$m \in \Z$ and~$E \in [1,\infty)$ are such that~$\S(m-1, E)$ is valid. Let~$\sigma \in (0,\nf12]$.
There exists~$C(d)<\infty$ such that, if
\begin{equation}
	\label{e.percolation.admissibility}
	\max\bigl\{ \exp(C\sigma^{-1} ) , \cstar^{-1} \bigr\}
	\leq
	E
	\leq \cgamma^{-\nf15}
	\,,
\end{equation}
then, for every~$q \in [1,\infty]$ and~$s \in [\cgamma,1]$, we have
\begin{equation}
	\sup_{L \in [m-1,\infty)\cap\Z}  \lambda_{s ,q}^{-1}(\cu_m\, ; \a_L) \shom_{m-1}
	\leq
	C
	+
	\O_{\Gamma_{\frac 12(1 -\sigma)}} \bigl( \exp \bigl( - C^{-1} E^{-2} \cgamma^{-1}\bigr)  \bigr)
	\label{e.cg.ellip.lower}
\end{equation}
and
\begin{equation}
	\Lambda_{s,q}(\cu_m\,;\a_m) \shom_{m-1} ^{-1}
	\leq
	C
	+
	\O_{\Gamma_{1}} \bigl( C\cstar^{-1}s^{-2}\cgamma \bigr)
	+
	\O_{\Gamma_{\frac 13(1 -\sigma)}}\bigl( \exp \bigl( - C^{-1} E^{-2} \cgamma^{-1}\bigr)  \bigr)
	\,.
	\label{e.cg.ellip.upper}
\end{equation}
\end{proposition}

In order to check that Proposition~\ref{p.cg.ellipticity.bounds} implies Proposition~\ref{p.tail.bounds}, we use  Lemma~\ref{l.mathcal.E.to.Lambdas}, which asserts
\begin{equation*}
	\mathcal{E}_{s,\infty,2}
	(\cu_m;\a_m,\shom_{m-1})^2
	\leq
	2  \shom_{m-1} ^{-1} \Lambda_{s,2}(\cu_m;\a_m)
	+
	2\shom_{m-1} \lambda_{s,2}^{-1} (\cu_m;\a_m)
\end{equation*}
and therefore~\eqref{e.cg.ellip.lower} and~\eqref{e.cg.ellip.upper} imply
\begin{equation*}
	\mathcal{E}_{s,\infty,2}
	(\cu_m;\a_m,\shom_{m-1} )^2
	\leq
	C
	+
	\O_{\Gamma_{1}} \bigl( C\cstar^{-1}s^{-2}\cgamma \bigr)
	+
	\O_{\Gamma_{\frac 13(1 -\sigma)}}\bigl( \exp \bigl( - C^{-1} E^{-2} \cgamma^{-1}\bigr)  \bigr)\,.
\end{equation*}
Taking~$\sigma=\frac14$, using~$s\geq4\cgamma$, and taking square roots give~\eqref{e.tail.bounds}, after weakening the stretched-exponential exponents.

\smallskip

The rest of this subsection is devoted to the proof of Proposition~\ref{p.cg.ellipticity.bounds}. The lower bound~\eqref{e.cg.ellip.lower} is immediate for~$m\leq\mstar$ from~\eqref{e.CG.bounds.1} and~\eqref{e.basecase.diffusivity}, uniformly in~$L$. By~\eqref{e.bound.Lambdas.by.Es} and~\eqref{e.basecase.homogenization}, the upper bound, and hence the proposition, is immediate for~$m\leq\mstarstar$. We may therefore assume throughout the section that~$m>\mstarstar$.

\smallskip

The idea of the proof of Proposition~\ref{p.cg.ellipticity.bounds} is to use subadditivity to bound the coarse-grained ellipticity constants by a sum over the same quantities on smaller-scale subcubes and to bound these with the induction hypothesis. The improvement comes from the fact that, by an adaptation of the subadditivity argument, we can avoid a sparse set of cubes for which the coarse-grained ellipticity constants are large. By avoiding the worst cubes, this argument shows that the coarse-graining operation is more insensitive to extreme regions of the coefficient field than a simple average, and hence the sum has better tail bounds compared to the assumption in the induction hypothesis.

\smallskip

We will still select bad triadic cubes at a fixed mesoscopic scale, but we pass to a triangulation~$\SW$ of~$\cu_m$ so that we can build piecewise affine competitors that are exactly constant on the selected bad region. Our notation for simplices and their vertices can be found in the paragraph containing~\eqref{e.simplex.def}, above.

\smallskip

For each slope~$p \in \Rd$, we construct an approximation of the affine function~$\linear_p$ given by~$\hat{\linear}_p$ which satisfies
\begin{equation}
	\label{e.hat.linear.1}
	\hat{\linear}_p \in \linear_p + H^1_0(\cu_m)
\end{equation}
and
\begin{equation}
	\label{e.hat.linear.2}
	\text{
	$\hat{\linear}_p$ is piecewise affine on~$\cu_m$, and~$\hat{\linear}_p \vert_\spx$ is affine for each element~$\spx$ of~$\SW$\,.
	}
\end{equation}
The family of approximate affines will be linear in~$p$, that is,
\begin{equation}
	\label{e.hat.linear.linearity}
	\hat{\linear}_{\alpha p_1+\beta p_2}=\alpha\hat{\linear}_{p_1}+\beta\hat{\linear}_{p_2}
	\qquad \forall \alpha,\beta\in\R,\ p_1,p_2\in\Rd\,.
\end{equation}
We define the antisymmetric matrix~$D_q$ for each~$q\in\Rd$ by
\begin{equation}
	\label{e.dq}
	(D_q)_{ij}  \coloneqq
	\frac1{d-1} \bigl( q_j\linear_{e_i}  - q_i\linear_{e_j}  \bigr)\,.
\end{equation}
It is immediate that\footnote{We use the convention that the divergence of a matrix-valued function~$A=(A_{ij})$ is the vector field with~$j$th coordinate given by~$(\nabla \cdot A)_j = \sum_{i=1}^d \partial_{x_i} A_{ij}$.}
\begin{equation*}
	\nabla \cdot D_q = q\,.
\end{equation*}
Given a family~$\{ \hat{\linear}_p \,:\, p \in\Rd \}$ satisfying~\eqref{e.hat.linear.1}, we define, for each~$q\in\Rd$ an antisymmetric matrix~$\hat{D}_q$ by
\begin{equation}
	\label{e.hatdq}
	(\hat{D}_q)_{ij}  \coloneqq
	\frac1{d-1} \bigl( q_j\hat{\linear}_{e_i}  - q_i\hat{\linear}_{e_j}  \bigr)\,.
\end{equation}
Since~$\hat D_q-D_q$ has antisymmetric entries in~$H^1_0(\cu_m)$, its divergence is divergence-free and has zero normal trace. Therefore,
\begin{equation}
	\nabla \cdot \hat{D}_q \in q + L^2_{\sol,0}(\cu_m)
	\,.
	\label{e.hat.D.bc}
\end{equation}
For each~$\spx\in \SW$, let~$\nabla \hat{\linear}_p (\spx)$ and~$\nabla \cdot \hat{D}_q(\spx)$ denote, respectively, the constant values of~$\nabla \hat{\linear}_p$ and~$\nabla \cdot \hat{D}_q$ in~$\spx$.
\begin{lemma}
\label{l.subadd.betterer}
Let~$\SW$ be a partition of~$\cu_m$ into triadic simplices and let~$\{ \hat{\linear}_p : p\in\Rd \}$ be a family of functions satisfying~\eqref{e.hat.linear.1},~\eqref{e.hat.linear.2} and~\eqref{e.hat.linear.linearity}, and let~$\hat{D}_q$ be defined by~\eqref{e.hatdq}. Then, for every~$p,q\in\Rd$,
\begin{equation}
	\label{e.subadd.betterer}
	\begin{pmatrix}
		p \\
		q
	\end{pmatrix}
	\cdot \bfA(\cu_m; \a)
	\begin{pmatrix}
		p \\
		q
	\end{pmatrix}
	\leq
	\sum_{\spx \in \SW}
	\frac{|\spx|}{|\cu_m|}
	\begin{pmatrix}
		\nabla \hat{\linear}_p (\spx) \\
		(\nabla \cdot \hat{D}_q)(\spx)
	\end{pmatrix}
	\cdot
	\bfA(\spx; \a)
	\begin{pmatrix}
		\nabla \hat{\linear}_p (\spx) \\
		(\nabla \cdot \hat{D}_q)(\spx)
	\end{pmatrix}
	\,.
\end{equation}
\end{lemma}
\begin{proof}
 Using that~$\nabla \hat{\linear}_p \in p+L^2_{\mathrm{pot},0}(\cu_m)$ and~$\nabla \cdot \hat{D}_q \in q + L^2_{\sol,0}(\cu_m)$, we have by~\eqref{e.variational.mu.U.P} and~\eqref{e.bigA.def},
\begin{align*}
	\begin{pmatrix}
		p \\
		q
		\end{pmatrix}  \cdot \bfA(\cu_m; \a) \begin{pmatrix}
		p \\
		q
	\end{pmatrix}
	&
	=
	\min_{X \in (p+L^2_{\pot,0}(\cu_m)) \times (q + \Lsolo(\cu_m))}
	\fint_{\cu_m} X \cdot \bfA  X
	\notag \\ &
	=
	\min_{ u \in \hat{\linear}_p + H^1_0(\cu_m) \, ,  \,  \mathbf{g} \in \nabla \cdot \hat{D}_q + L^2_{\sol,0} (\cu_m)}
	\fint_{\cu_m} \begin{pmatrix} \nabla u \\ \mathbf{g} \end{pmatrix} \cdot \bfA   \begin{pmatrix} \nabla u \\ \mathbf{g} \end{pmatrix}
	\notag \\ &
	\leq
	\sum_{\spx \in \SW}
	\frac{|\spx|}{|\cu_m|}
	\min_{ u \in \hat{\linear}_p + H^1_0(\spx) \, ,  \,  \mathbf{g} \in \nabla \cdot \hat{D}_q + L^2_{\sol,0} (\spx)}
	\fint_{\spx} \begin{pmatrix} \nabla u \\ \mathbf{g} \end{pmatrix} \cdot \bfA   \begin{pmatrix} \nabla u \\ \mathbf{g} \end{pmatrix}
	\notag \\ &
	=
	\sum_{\spx \in \SW}
	\frac{|\spx|}{|\cu_m|}
	\begin{pmatrix}
		\nabla \hat{\linear}_p(\spx) \\
		(\nabla \cdot \hat{D}_q)(\spx)
		\end{pmatrix}  \cdot \bfA(\spx; \a) \begin{pmatrix}
		\nabla \hat{\linear}_p(\spx) \\
		(\nabla \cdot \hat{D}_q)(\spx)
	\end{pmatrix}
	\,.
\end{align*}
The countable sum is justified by applying the gluing argument to the first~$N$ simplex layers and the uncovered remainder, and then sending~$N\to\infty$; the remainder has vanishing measure and the corresponding competitor fields converge in~$L^2$.
This completes the proof.
\end{proof}

Note that the previous lemma reduces to the usual subadditivity inequality if we specialize to~$\hat{\linear}_p=\linear_p$
and~$\hat{D}_q = D_q$.
The flexibility of choosing more general functions~$\hat{\linear}_p$ and~$\hat{D}_q$ is that it allows us to avoid the ``bad'' elements~$\spx\in \SW$ from the sum on the right of~\eqref{e.subadd.betterer}, since these terms vanish if~$\nabla\hat{\linear}_p(\spx)=0$ for all~$p$, as in this case~$(\nabla \cdot \hat{D}_q) (\spx) = 0$ as well.

\subsubsection{Bad cubes and the minimal scale separation}
\label{sss.bad.cubes}

In this subsubsection we identify a sparse random set of ``bad'' triadic cubes and establish diameter and density bounds for its connected components. The main output is a random variable~$\hsep \in \N_0$, which we call the \emph{minimal scale separation}, with stretched-exponential tails. This parameter will control the gap between adjacent layers in the Whitney-type partition constructed in the next subsubsection: for scale separations~$h \geq \hsep$, bad clusters have small diameter and low density. The construction here does not yet involve any partition of~$\cu_m$; the partition will be built using the bounds established below.

\smallskip

We work under the assumption that~$\mathcal{S}(m-1,E)$ is valid for some~$m \in \Z$ and~$E \in [1,\infty)$.

\smallskip

For each~$j \in \Z$ and~$z \in 3^j \Zd$, we define two types of bad events for the cube~$z + \cu_j$. The first captures the event that the oscillation of the stream matrix at wavelengths above scale~$3^j$ is too large:
\begin{equation}
	\label{e.Bosc.def}
	\mathcal{B}_{\mathrm{osc}}(z + \cu_j)
	\coloneqq
	\biggl\{
	\sup_{L \geq j} 3^{2 j} \|\nabla (\k_{L} - \k_{j}) \|_{\underline{W}^{1,\infty} (z + \cu_j)}
	>  \delta_{\mathrm{osc}}
	\cstar^{\nf12} \cgamma^{-\nf12} 3^{\cgamma j}
	\biggr\}\,,
\end{equation}
where~$\delta_{\mathrm{osc}} > 0$ is a small constant depending only on~$d$, chosen so that the coarse-grained sensitivity estimates (Lemmas~\ref{l.lambda.sensitivity} and~\ref{l.J.sensitivity}) can be applied when~$\mathcal{B}_{\mathrm{osc}}(z + \cu_j)$ does not occur. Specifically, we require
\begin{equation}
	\delta_{\mathrm{osc}} \leq 2^{-10} C_{\eqref{e.J.sensitivity.smallness.condition}}^{-1}
	\,.
	\label{e.delta.osc.choice}
\end{equation}
The event~$\mathcal{B}_{\mathrm{osc}}(z + \cu_j)$ admits a decomposition by scale; for each~$L > j$, if we define
\begin{equation}
	\label{e.BoscL.def}
	\mathcal{B}_{\mathrm{osc},L}(z+\cu_j) \coloneqq \Bigl\{ 3^{2j} \| \nabla \mathbf{j}_L \|_{\underline{W}^{1,\infty}(z + \cu_j)} > \tfrac15 \delta_{\mathrm{osc}} \cstar^{\nf12} \cgamma^{-\nf12} 3^{\cgamma j} 3^{-\frac15(L-j)} \Bigr\}\,\,,
\end{equation}
then we have
\begin{equation}
	\label{e.Bosc.decomp}
	\mathcal{B}_{\mathrm{osc}}(z + \cu_j) \subseteq \bigcup_{L > j} \mathcal{B}_{\mathrm{osc},L}(z+\cu_j) \,\,.
\end{equation}
The second bad event captures the failure of the local coarse-grained ellipticity bounds:
\begin{equation}
	\label{e.Bloc.def}
	\mathcal{B}_{\mathrm{loc}}(z + \cu_j)
	\coloneqq
	\Bigl\{
	\lambda^{-1}_{\nf14,2}(z + \cu_j; \a_{j})  > 10 \shom^{-1}_j
	\Bigr\} \cup
	\Bigl\{
	\Lambda_{\nf14,2}(z + \cu_j; \a_{j})  > 10 \shom_j
	\Bigr\}
	\,.
\end{equation}
We say that~$z + \cu_j$ is \emph{bad} if the event
\begin{equation}
	\label{e.B.def}
	\mathcal{B}(z+\cu_j)  \coloneqq  \mathcal{B}_{\mathrm{osc}}(z + \cu_j) \cup \mathcal{B}_{\mathrm{loc}}(z + \cu_j)
\end{equation}
occurs, and \emph{good} otherwise.

\smallskip

To control the geometry of bad clusters in a partition that involves cubes of varying sizes (as will be constructed in the next subsubsection), we need to ensure that a cube and all of its descendants up to a fixed number of generations are good. For this reason, we introduce the \emph{extended} bad events. For each~$j \in \Z$ and~$z \in 3^j \Zd$, we define
\begin{equation}
	\label{e.Bext.def}
	\mathcal{B}^{*}(z + \cu_j)
	\coloneqq
	\bigcup_{i \in \{0,1,\ldots,9\}} \; \bigcup_{z' \in 3^{j-i}\Zd \cap (z + \cu_j) }
	\mathcal{B}(z' +\cu_{j-i})
	\,.
\end{equation}
In words,~$\mathcal{B}^{*}(z + \cu_j)$ is the event that~$z + \cu_j$ or any of its triadic descendants up to nine generations is bad. The choice of nine generations ensures that, in the Whitney partition constructed below, bad clusters cannot span more than two adjacent layers.

\smallskip

We next establish probability bounds for the bad events defined above.

\begin{lemma}[Probability of bad events]
	\label{l.bad.event.probabilities}
	Assume~$m \in \Z$ and~$E \in [1,\infty)$ are such that~$\mathcal{S}(m-1,E)$ is valid. There exists~$c(d),C(d) < \infty$ such that, if
	\begin{equation}
		\label{e.bad.event.admissibility.prelim}
		E \geq C \cstar^{-1}
		\qand
		\cgamma \leq E^{-5}\,,
	\end{equation}
	then for every~$j \in \Z$ with~$j \leq m-1$ and every~$z \in 3^j \Zd$, we have
	\begin{equation}
		\label{e.Bloc.prob}
		\P\bigl[ \mathcal{B}_{\mathrm{loc}}(z + \cu_j) \bigr]
		\leq
		\exp\bigl( - c E^{-2} \cgamma^{-1} \bigr)\,.
	\end{equation}
	Moreover, for every~$L > j$,
	\begin{equation}
		\label{e.BoscL.prob}
		\P\bigl[ \mathcal{B}_{\mathrm{osc},L}(z + \cu_j) \bigr]
		\leq
		\exp\bigl( - c E^{-2} \cgamma^{-1} 3^{\frac32(L-j)} \bigr)\,.
	\end{equation}
\end{lemma}
\begin{proof}
\emph{Step 1.} We prove~\eqref{e.Bloc.prob}. By~\eqref{e.bound.Lambdas.by.Es}, we have
\begin{equation*}
	\shom_j^{-1} \Lambda_{\nf14,2}(z + \cu_j; \a_j) \leq 4 + 4 \bigl( \mathcal{E}_{\nf14,\infty,2}(z+\cu_j; \a_j, \shom_j) \bigr)^2
\end{equation*}
and similarly for~$\shom_j \lambda_{\nf14,2}^{-1}(z + \cu_j; \a_j)$. Thus, by Markov's inequality and the induction hypothesis~\eqref{e.new.induction.for.shom},
\begin{equation*}
	\P\bigl[ \mathcal{B}_{\mathrm{loc}}(z + \cu_j) \bigr]
	\leq
	\P\Bigl[ \bigl( \mathcal{E}_{\nf14,\infty,2}(z+\cu_j; \a_j, \shom_j) \bigr)^2 > C \Bigr]
	\leq
	\exp\bigl( - c E^{-2} \cgamma^{-1} \bigr)\,.
\end{equation*}

\emph{Step 2.} To prove~\eqref{e.BoscL.prob},  we observe that by~\eqref{e.nabla.jk.O} and~\eqref{e.maxy.bound}, for every~$L \geq j$,
\begin{equation*}
	3^{2j} \| \nabla \mathbf{j}_L \|_{\underline{W}^{1,\infty}(z + \cu_j)}
	\leq
	\O_{\Gamma_2}\bigl( C 3^{\cgamma L - (L-j)} \bigr)\,.
\end{equation*}
Therefore, comparing with the threshold, the probability of~$\mathcal{B}_{\mathrm{osc},L}(z+\cu_j)$ is bounded by
\begin{equation*}
	\P\bigl[ \mathcal{B}_{\mathrm{osc},L}(z+\cu_j) \bigr]
	\leq
	\exp\Bigl( - c \cstar \cgamma^{-1} 3^{2(\frac45 - \cgamma)(L-j)} \Bigr)
	\leq
	\exp\Bigl( - c E^{-2} \cgamma^{-1} 3^{\frac32(L-j)} \Bigr)\,,
\end{equation*}
where we used~$\cgamma \leq \frac{1}{20}$ in the last step.
\end{proof}

\smallskip

We next apply Lemma~\ref{l.percolation.bound.general} to control the geometry of connected components of bad cubes. For each~$h \in \N$, we consider the lattice~$3^{m-h} \Zd \cap \cu_m$ and the random subset
\begin{equation}
	\label{e.bad.cluster.def}
	\mathcal{B}^{*}_{m-h}  \coloneqq  \bigl\{ z + \cu_{m-h} \,:\, z \in 3^{m-h} \Zd \cap \cu_m\,,\; \mathcal{B}^{*}(z + \cu_{m-h}) \text{ occurs} \bigr\}\,.
\end{equation}
 A \emph{2-connected component} (or \emph{2-cluster}) of~$\mathcal{B}^{*}_{m-h}$ is a maximal subset~$A \subseteq \mathcal{B}^{*}_{m-h}$ such that any two cubes in~$A$ can be connected by a chain
 of successive cubes~$z+\cu_{m-h},z'+\cu_{m-h}\in A$ with
$\|3^{-(m-h)}(z-z')\|_{\ell^\infty}\leq2$. Thus adjacency is measured in
lattice units and is invariant under a change of scale. For~$z \in 3^{m-h} \Zd \cap \cu_m$, we denote by~$\mathcal{B}^{*}_{m-h}(z;2)$ the 2-connected component containing~$z + \cu_{m-h}$ if~$\mathcal{B}^{*}(z + \cu_{m-h})$ occurs, and the empty set otherwise.

\begin{lemma}[Diameter and density bounds for bad clusters]
\label{l.percolation.bad.clusters}
Assume~$m \in \Z$ and~$E \in [1,\infty)$ are such that~$\mathcal{S}(m-1,E)$ is valid. Let~$\sigma \in (0,\frac12]$ and~$b \in (0,\frac19]$. There exist~$c(d)>0$ and~$C(d) < \infty$ such that, if
\begin{equation}
	\label{e.percolation.admissibility.bad.clusters}
	E \geq \exp(C \sigma^{-1}) \vee C\cstar^{-1}\vee Cb^{-1}
	\qand
	\cgamma \leq E^{-5}\,,
\end{equation}
then for every~$h \in \N$, we have the following bounds.
\begin{itemize}
\item \emph{Diameter bound:}
\begin{equation}
	\label{e.diameter.bound.bad.clusters}
	\P\Bigl[ \max_{z \in 3^{m-h} \Zd \cap \cu_m} \diam\bigl( \mathcal{B}^{*}_{m-h}(z;2) \bigr) \geq 3^{bh} \cdot 3^{m-h} \Bigr]
	\leq
	\exp\bigl( - 3^{(1-\sigma) b h} \bigr)\,.
\end{equation}

\item \emph{Density bound:}
\begin{equation}
	\label{e.density.bound.bad.clusters}
	\P\Biggl[ \avsum_{z \in 3^{m-h} \Zd \cap \cu_m} \indc_{\mathcal{B}^{*}(z + \cu_{m-h})} > \exp\bigl( - c E^{-2} \cgamma^{-1} \bigr) + 3^{-\frac{h}{2}} \Biggr]
	\leq
	\exp\bigl( - 3^{\frac h 2} \bigr)\,.
\end{equation}
\end{itemize}
\end{lemma}

\begin{proof}
	For each~$z \in 3^{m-h} \Zd \cap \cu_m$ and~$L \in \N_0$, define the scale~$L$ bad event
	\begin{equation*}
		\mathrm{B}_L(z + \cu_{m-h}) \coloneqq
		\begin{cases}
			\displaystyle\bigcup_{i=0}^{9} \bigcup_{z' \in 3^{m-h-i}\Zd \cap (z + \cu_{m-h})} \mathcal{B}_{\mathrm{loc}}(z' + \cu_{m-h-i}) & \text{if } L = 0\,, \\[12pt]
			\displaystyle\bigcup_{i=0}^{9} \bigcup_{z' \in 3^{m-h-i}\Zd \cap (z + \cu_{m-h})} \mathcal{B}_{\mathrm{osc},m-h-i+L}(z' + \cu_{m-h-i}) & \text{if } L \geq 1
		\end{cases}
	\end{equation*}
	and note that
	\begin{equation*}
		\mathcal{B}^{*}(z + \cu_{m-h}) \subseteq \bigcup_{L \in \N_0} \mathrm{B}_L(z + \cu_{m-h})\,.
	\end{equation*}
	By Lemma~\ref{l.bad.event.probabilities}, with~$T \coloneqq (c_{\eqref{e.BoscL.prob}} E^{-2} \cgamma^{-1})^{\nf12}$ these events satisfy the probability bounds
		\begin{equation*}
			\P\bigl[ \mathrm{B}_L(z + \cu_{m-h}) \bigr] \leq \exp\bigl( - T^2 3^{\frac32 L} \bigr) \qquad \forall L \in \N_0 \,\,,
		\end{equation*}
		By~\ref{a.j.frd}, for each~$L$ the cumulative families~$\{\mathrm B_{L'}(z):L'\leq L\}$ at sets of sites separated by more than~$C(d)3^L$ are independent. After a dimension-only shift of~$L$, this is the independence hypothesis of Lemma~\ref{l.percolation.bound.general}. The admissibility condition~\eqref{e.percolation.admissibility.bad.clusters} ensures that~$T \geq C \exp(C\sigma^{-1})$. Applying that lemma with tail exponent~$a=\frac32$ gives~\eqref{e.density.bound.bad.clusters}; its maximal-diameter estimate with diameter parameter~$(1-\sigma)b$ gives~\eqref{e.diameter.bound.bad.clusters} after rescaling.
\end{proof}

\smallskip

We next introduce a random variable~$\hsep$ as the smallest scale separation above which both the diameter and density bounds hold.

\begin{proposition}[Minimal scale separation for bad clusters]
\label{p.minimal.scale.separation}
Assume~$m \in \Z$ and~$E \in [1,\infty)$ are such that~$\mathcal{S}(m-1,E)$ is valid. Let~$\sigma \in (0,\frac12]$ and~$b \in (0,\frac19]$. There exist~$c(d)>0$ and~$C(d) < \infty$ such that, if~\eqref{e.percolation.admissibility.bad.clusters} holds, then there exists~$\hsep \in \N_0$ satisfying
\begin{equation}
	\label{e.hsep.tails}
	3^{b \hsep} \leq \O_{\Gamma_{1-\sigma}}\bigl(C(\sigma,b)\bigr)
\end{equation}
such that, for every~$h \geq \hsep$, we have
\begin{equation}
	\label{e.diameter.deterministic}
	\max_{z \in 3^{m-h} \Zd \cap \cu_m} \diam\bigl( \mathcal{B}^{*}_{m-h}(z;2) \bigr) < 3^{bh} \cdot 3^{m-h}
\end{equation}
and
\begin{equation}
	\label{e.density.deterministic}
	\avsum_{z \in 3^{m-h} \Zd \cap \cu_m} \indc_{\mathcal{B}^{*}(z + \cu_{m-h})} \leq \exp\bigl( - c E^{-2} \cgamma^{-1} \bigr) + 3^{-\frac{h}{2}}\,.
\end{equation}
\end{proposition}

\begin{proof}
Let~$\mathsf G(h)$ denote the conjunction of
\eqref{e.diameter.deterministic} and~\eqref{e.density.deterministic} at
scale~$h$, and define the random variable
\begin{equation*}
	\widehat h_{\mathrm{sep}}
	\coloneqq
	\inf\bigl\{j\in\N_0:\ \mathsf G(h)\text{ holds for every }h\geq j\bigr\}
	\in\N_0\cup\{\infty\}\,.
\end{equation*}
By a union bound over~$h \geq j$ and Lemma~\ref{l.percolation.bad.clusters}, we have
\begin{equation*}
	\P\bigl[ \widehat h_{\mathrm{sep}} > j \bigr]
	\leq
	\sum_{h=j}^{\infty}\left(
	\exp\bigl( - 3^{(1-\sigma) b h} \bigr)
	+\exp\bigl(-3^{h/2}\bigr)\right)
	\leq
	C(\sigma,b) \exp\bigl( - c(\sigma,b)3^{(1-\sigma) b j} \bigr)\,.
\end{equation*}
The preceding tail estimate shows that~$\widehat h_{\mathrm{sep}}$ is finite. Setting~$\hsep=\widehat h_{\mathrm{sep}}$ gives~\eqref{e.hsep.tails}, and the two bounds hold for every~$h\geq\hsep$ by construction.
\end{proof}

\subsubsection{The Whitney partition}
\label{sss.whitney.partition}

We now construct a Whitney-type partition of~$\cu_m$ into triadic cubes of varying sizes. The partition is designed so that cubes become smaller as they approach the boundary~$\partial \cu_m$, ensuring that no element of the partition touches the boundary. This is necessary to enforce the boundary condition~\eqref{e.hat.linear.1} for the approximate affine functions constructed in the next subsubsection.

\smallskip

The partition depends on the minimal scale separation~$\hsep$ from Proposition~\ref{p.minimal.scale.separation}, as well as two parameters~$b \in (0,\frac19]$ and~$k_0 \in \N$ to be selected below. We define the \emph{scale separation sequence}~$\{h_n\}_{n \in \N_0}$ by
\begin{equation}
	\label{e.hn.def}
	h_n  \coloneqq  \bigl\lceil b (1-b)^{-1} n \bigr\rceil + \hsep + k_0\,.
\end{equation}
Observe that for~$b \leq \frac18$, consecutive terms satisfy
\begin{equation}
	\label{e.hn.gap}
	h_{n+1} - h_n \leq \bigl\lceil b(1-b)^{-1} \bigr\rceil \leq 1\,,
\end{equation}
so the scale separation grows slowly. This bound, together with the definition of~$\mathcal{B}^*$ which looks at nine generations of descendants, will ensure that bad clusters cannot span more than two adjacent layers of the partition.

\smallskip

\paragraph{Construction of the partition.}
The Whitney partition~$\W(\cu_m)$ is constructed by the following recursive procedure.
\begin{itemize}
\item \emph{Layer~$1$.} Divide~$\cu_m$ into its~$3^d$ triadic subcubes of size~$3^{m-1}$. The central cube~$\cu_{m-1}$ does not touch the boundary~$\partial \cu_m$; subdivide it into its~$3^{d h_1}$ triadic subcubes of size~$3^{m-1-h_1}$ and denote these finalized interior cubes by~$\W(\cu_m,1)$. The remaining~$3^d - 1$ cubes touch the boundary and are left undivided for now.

\item \emph{Layer~$n \geq 2$.} Take each cube from layer~$n-1$ that touches~$\partial \cu_m$ and subdivide it into its~$3^d$ triadic subcubes of size~$3^{m-n}$. For each resulting cube of size~$3^{m-n}$ that does not touch~$\partial \cu_m$, further subdivide it into its~$3^{d h_n}$ triadic subcubes of size~$3^{m-n-h_n}$; these finalized interior cubes form~$\W(\cu_m,n)$.

\item \emph{Infinite-layer limit.} Run the recursion for every~$n\in\N$ and
let~$\W(\cu_m)=\bigcup_{n\geq1}\W(\cu_m,n)$.  No member of~$\W(\cu_m)$
touches~$\partial\cu_m$, and the cubes cover~$\cu_m$ up to the null limiting
boundary layer. After~$N$ layers the uncovered remainder has relative
volume at most~$2d\,3^{-N}$.
\end{itemize}
The cubes in~$\W(\cu_m)$ have sizes in the set~$\{ 3^{m-n-h_n} : n \in \N \}$. We denote the cubes of a given size by
\begin{equation}
	\label{e.W.n.def}
	\W(\cu_m, n)  \coloneqq  \bigl\{ \cu \in \W(\cu_m) \,:\, \size(\cu) = 3^{m-n-h_n} \bigr\}\,.
\end{equation}
For every~$Q\in\W(\cu_m,n)$,~$3^{m-n}\leq\dist(Q,\partial\cu_m)<3^{m-n+1}$.

\smallskip

\paragraph{Bad cubes in the partition.}
A cube~$\cu \in \W(\cu_m)$ is \emph{bad} if the bad event~$\mathcal{B}(\cu)$ defined in~\eqref{e.B.def} occurs. We collect the bad cubes into
\begin{equation}
	\label{e.bad.cubes.in.partition}
	\mathcal{I}  \coloneqq  \bigl\{ \cu \in \W(\cu_m) \,:\, \mathcal{B}(\cu) \text{ occurs} \bigr\}\,,
\end{equation}
and, for each~$n \in \N$, we write~$\mathcal{I}_n  \coloneqq  \mathcal{I} \cap \W(\cu_m, n)$ for the bad cubes at layer~$n$. The \emph{neighborhood} of a subset~$\mathcal{J} \subseteq \W(\cu_m)$ is defined by
\begin{equation}
	\label{e.neighborhood.def}
	\mathcal{N}(\mathcal{J})  \coloneqq  \bigl\{ \cu \in \W(\cu_m) \,:\, \dist(\cu, \mathcal{J}) = 0 \bigr\}\,,
\end{equation}
where~$\dist(\cu, \mathcal{J})  \coloneqq  \inf\{ \| x - y \|_{\ell^\infty} : x \in \cu,\, y \in \bigcup_{\cu' \in \mathcal{J}} \cu' \}$.

\smallskip

Two Whitney cubes are called adjacent if their closures touch, equivalently if their~$\ell^\infty$ distance is zero. A \emph{connected component} of~$\mathcal{I}$ is a maximal subset~$\mathcal{C} \subseteq \mathcal{I}$ such that any two cubes in~$\mathcal{C}$ can be connected by a chain of pairwise adjacent cubes in~$\mathcal{I}$. We say that a connected component~$\mathcal{C}$ \emph{intersects layer~$n$} if~$\mathcal{C} \cap \mathcal{I}_n \neq \emptyset$.

\smallskip

The following lemma shows that bad clusters have small diameter and cannot span more than two adjacent layers. The key observation is that any connected component~$\mathcal{C}$ of bad cubes in the partition---even one that zigzags between two adjacent layers---can be ``lifted'' to a connected component of extended bad cubes~$\mathcal{B}^*$ at the coarsest scale in the region, to which the percolation bounds of Proposition~\ref{p.minimal.scale.separation} apply.

\begin{lemma}[Geometry of bad clusters in the partition]
\label{l.bad.clusters.geometry}
Under the hypotheses of Proposition~\ref{p.minimal.scale.separation}, assume~$k_0\geq3$, let~$\mathcal C$ be a connected component of~$\mathcal I$ and set
\begin{equation*}
	j_{\mathcal C}\coloneqq
	\min\{j:\mathcal N(\mathcal C)\cap\W(\cu_m,j)\ne\varnothing\}\,.
\end{equation*}
\begin{enumerate}
\item The touching neighborhood occupies two consecutive layers:
\begin{equation*}
	\mathcal N(\mathcal C)\subseteq
	\W(\cu_m,j_{\mathcal C})\cup\W(\cu_m,j_{\mathcal C}+1)\,.
\end{equation*}

\item In the~$\ell^\infty$ metric,
\begin{equation}
	\label{e.cluster.diameter.bound}
	\diam_{\ell^\infty}(\mathcal N(\mathcal C))
	<4\,3^{b(j_{\mathcal C}+h_{j_{\mathcal C}})}
	3^{m-j_{\mathcal C}-h_{j_{\mathcal C}}}\,.
\end{equation}

\end{enumerate}
\end{lemma}

\begin{proof}
Put~$r\coloneqq j_{\mathcal C}$ and lift each Whitney cube to its ancestor at scale~$3^{m-r-h_r}$. The enlargement by nine generations in~\eqref{e.Bext.def} and~\eqref{e.hn.gap} ensure that every adjacent bad chain lifts to a 2-connected~$\mathcal B^*$-bad family. If~$\mathcal N(\mathcal C)$ met layer~$r+2$, the first subchain crossing two layer interfaces would have diameter at least~$c3^{m-r}$, whereas~\eqref{e.diameter.deterministic} bounds its lift by~$3^{m-r}3^{-(1-b)k_0}<c3^{m-r}$. Thus only layers~$r,r+1$ occur. Applying the same lift to all of~$\mathcal C$ and enlarging by the touching cubes gives
\begin{align*}
	\diam_{\ell^\infty}\bigl(\mathcal N(\mathcal C)\bigr)
	& <
	4\,3^{b(r+h_r)} \cdot 3^{m-r-h_r}
	\leq 4\,3^{m-r} \cdot 3^{-(1-b)k_0}\,,
\end{align*}
where the second inequality uses~$h_r\geq b(1-b)^{-1}r+k_0$.
This is~\eqref{e.cluster.diameter.bound}.
\end{proof}

\smallskip

The next lemma gives the density bound for bad cubes.

\begin{lemma}[Density of bad cubes]
\label{l.bad.cubes.density}
Under the hypotheses of Proposition~\ref{p.minimal.scale.separation}, for every~$n \in \N$, we have
\begin{equation}
	\label{e.bad.cubes.density}
	\sum_{\cu \in \W(\cu_m,n)\cap\mathcal I}
	\frac{|\cu|}{|\cu_m|}
	\leq
	\exp\bigl( - c E^{-2} \cgamma^{-1} \bigr) + 3^{-\frac{n + h_n}{2}}\,.
\end{equation}
\end{lemma}

\begin{proof}
Since~$\mathcal{B}(\cu) \subseteq \mathcal{B}^*(\cu)$, it suffices to bound the density of~$\mathcal{B}^*$-bad cubes. The cubes in the partition~$\W(\cu_m, n)$ have size~$3^{m-n-h_n}$ and are a subset of~$\{ z + \cu_{m-n-h_n} : z \in 3^{m-n-h_n} \Zd \cap \cu_m \}$. Applying~\eqref{e.density.deterministic} with~$h = n + h_n \geq \hsep$ (which holds since~$h_n \geq \hsep + k_0$), we obtain~\eqref{e.bad.cubes.density}.
\end{proof}

\subsubsection{Construction of the affine approximation}
\label{sss.affine.approximation}

We now construct piecewise affine approximations of affine functions that are constant on the bad cubes identified in the previous subsubsection. The construction proceeds in two steps: first, we refine the Whitney partition~$\W(\cu_m)$ into a simplex partition~$\SW(\cu_m)$; second, we define a family of piecewise affine functions subordinate to this partition.

\smallskip

\paragraph{The simplex partition.}
We say that a function is \emph{subordinate} to a simplex partition~$\SW$ if its restriction to each element of the partition is affine. The simplex partition~$\SW(\cu_m)$ is obtained by subdividing each cube in~$\W(\cu_m)$ into simplices as follows. For each cube~$\cu \in \W(\cu_m, n)$ of size~$3^{m-n-h_n}$, we subdivide~$\cu$ into~$3^{d(1 + h_{n+1} - h_n)} d!$ simplices of size~$3^{m-(n+1)-h_{n+1}}$. We denote the collection of all simplices obtained from the cubes in~$\W(\cu_m,n)$ by~$\SW(\cu_m,n)$. That is,
\begin{equation}
	\label{e.SW.def}
	\SW(\cu_m)
	\coloneqq
	\bigcup_{n \in \N} \; \bigcup_{\cu \in \W(\cu_m, n)}
	S_{m-(n+1)-h_{n+1}}(\cu)\,,
\end{equation}
where we recall from~\eqref{e.simplex.def} that~$S_j(\cu)$ denotes the standard simplicial decomposition of the cube~$\cu$ into simplices of size~$3^j$.

Across unequal adjacent cells, affine traces are restricted from the common coarse mesh to the finer mesh.

\smallskip

The next lemma allows overlapping neighborhoods by summing the components active at each Whitney cube.

\begin{lemma}[Approximate affines avoiding a bad set]
\label{l.piecewise.affine.approx}
Let~$\mathcal I\subseteq\W(\cu_m)$, and denote by~$\mathfrak C$ the family of its connected components. Suppose that every touching neighborhood occupies two consecutive layers and satisfies~\eqref{e.cluster.diameter.bound}. For each component define
\begin{equation*}
	j_{\mathcal C}\coloneqq
	\min\{j:\mathcal N(\mathcal C)\cap\W(\cu_m,j)\ne\varnothing\}\,,
	\qquad
	\rho_{\mathcal C}\coloneqq
	3^{m-(j_{\mathcal C}+1)-h_{j_{\mathcal C}+1}}\,,
\end{equation*}
and, for~$R\in\W(\cu_m)$, define the active set
\begin{equation*}
	\mathcal A_R\coloneqq
	\{\mathcal C\in\mathfrak C:R\in\mathcal N(\mathcal C)\}\,.
\end{equation*}

Then there exists a family of piecewise affine functions
$\{\hat{\linear}_p\}_{p\in\Rd}$ satisfying~\eqref{e.hat.linear.linearity}
such that~$\hat{\linear}_p$ is subordinate to~$\SW(\cu_m)$,
$\hat{\linear}_p-\linear_p\in H^1_0(\cu_m)$, and
\begin{align}
	\label{e.hat.linear.properties}
	\nabla\hat{\linear}_p&=0
	&&\text{in every }Q\in\mathcal I,\notag\\
	\nabla\hat{\linear}_p&=p
	&&\text{on every simplex outside }\mathcal N(\mathcal I),\notag\\
	|\nabla\hat{\linear}_p-p|
	&\leq C(d)|p|\sum_{\mathcal C\in\mathcal A_R}
	\frac{\diam_{\ell^\infty}(\mathcal C)}{\rho_{\mathcal C}}
	&&\text{on every simplex contained in }R\in\W(\cu_m)\,.
\end{align}
Moreover, the Whitney geometry gives the dimension-only overlap bound
$|\mathcal A_R|\leq M(d)$.
\end{lemma}

\begin{proof}
For each~$\mathcal C$, interpolate on the common mesh of scale~$\rho_{\mathcal C}$ a correction~$\psi_{\mathcal C,p}$ supported in~$\mathcal N(\mathcal C)$, equal to a constant minus~$\linear_p$ on~$\mathcal C$, and zero on the boundary of this neighborhood; then restrict it to the final fine simplices. The inverse affine estimate gives~$|\nabla\psi_{\mathcal C,p}|\leq C(d)|p|\diam_{\ell^\infty}(\mathcal C)/\rho_{\mathcal C}$. The finite sums over components with~$j_{\mathcal C}\leq N$ belong to~$H^1_0(\cu_m)$, and the overlap bound and summability of~$3^{-j}3^{2b(j+h_j)}$ make them Cauchy in~$H^1$. Their limit defines~$\hat\linear_p=\linear_p+\sum_{\mathcal C}\psi_{\mathcal C,p}$. Since~$\mathcal A_R$ is finite for every~$R$, the affine slopes stabilize locally and the asserted properties pass to the limit. Whitney comparability gives~$|\mathcal A_R|\leq M(d)$.
\end{proof}

\begin{remark}
\label{r.gradient.bound.simplified}
If~$R$ lies in layer~$n$,~\eqref{e.cluster.diameter.bound} and~\eqref{e.hn.gap} give~$|\nabla\hat{\linear}_p-p|\leq C(d)|p|3^{b(n+h_n)}$ on every simplex in~$R$.
\end{remark}

\subsubsection{Coarse-grained multiscale estimate}
\label{sss.cg.multiscale}

We prove the following estimate, from which Proposition~\ref{p.cg.ellipticity.bounds} easily follows.

\begin{proposition}
\label{p.bfA.multiscalebound}
Assume that~$m \in\Z$ with~$m-1 \geq \mstarstar$ and~$E \in[1,\infty)$ are such that~$\S(m-1,E)$ is valid. Let~$\sigma \in (0,\nf12]$. There exist~$c(d)>0$ and~$C(d)<\infty$ such that, if
\begin{equation}
	\cgamma \leq E^{-5}
	\qand
	E \geq  \exp(C \sigma^{-1}) \vee C\cstar^{-1}\,,
	\label{e.Esigmares.in.proof}
\end{equation}
then for every~$n \in \Z$ with~$n \leq m-1$,
\begin{equation}
	\sup_{L \in [m-1, \infty) \cap \Z} 3^{-\cgamma (m-n)}  \shom_{m-1} |\s_{L,*}^{-1}(\cu_n)| \leq  C+  \O_{\Gamma_{\frac{1}{2}(1-\sigma)}}(\exp(-c E^{-2} \cgamma^{-1}))
	\label{e.slstar.multiscale}
\end{equation}
and, for every~$n\leq m-1$ and every~$L \in \Z$ with~$L \geq m-1$,
\begin{align}
	\shom^{-1}_{m-1} |\b_L(\cu_n)| &\leq C + \O_{\Gamma_1}(C \cstar^{-1}\cgamma (L-n+1) 3^{2 \cgamma(L-n) - \cgamma(m-n)})
	\notag \\ & \qquad
	+ \O_{\Gamma_{\frac{1}{3}(1 - \sigma)}}( 3^{2 \cgamma(L-n) - \cgamma(m-n)} \exp(-c E^{-2} \cgamma^{-1}))
	\,.
	\label{e.bL.multiscale}
\end{align}
\end{proposition}

\begin{proof}[Proof of Proposition~\ref{p.cg.ellipticity.bounds} assuming Proposition~\ref{p.bfA.multiscalebound}]
We apply Proposition~\ref{p.bfA.multiscalebound} with~$\sigma \mapsto \sigma/4$.
By~\eqref{e.ellipticities.change.q}, it is enough to prove both estimates for~$q=1$.
The estimate~\eqref{e.slstar.multiscale} yields, by~\eqref{e.maxy.bound}, for every~$n \in \Z$ with~$n \leq m-1$,
\begin{align*}
	&\sup_{L \in [m-1, \infty) \cap \Z}
	\max_{z \in 3^n \Zd \cap \cu_m}
	\bigl|\s_{L,*}^{-1}(z+\cu_n) \bigr|   \shom_{m-1}
	\\
	&\qquad
	\leq
	3^{\cgamma (m-n)} \Bigl(
	C
	+
	\O_{\Gamma_{\frac{1}{2}(1-  \sigma/4 )}} \Bigl(C   (m-n+1)^{2 + \sigma}
	\exp( - c E^{-2} \cgamma^{-1})  \Bigr)
	\Bigr)
	\,\,,
\end{align*}
By the one-step subadditivity~\eqref{e.subadda.nosymm}, the top-scale blocks with~$n=m$ are bounded by the maximum of their scale~$m-1$ children. Thus the lower top block is covered by the preceding~$n=m-1$ estimate, and the upper one by the analogous estimate below.
Summing over~$n$ before taking the square in the definition of~$\lambda_{s,1}^{-1}$, the deterministic part contributes~$C(s/(2s-\cgamma))^2\leq C$, while the exponentially small term is multiplied by a fixed power of~$(2s-\cgamma)^{-1}$.
Similarly~\eqref{e.bL.multiscale} yields for every~$n \in \Z$ with~$n \leq m-1$,
\begin{align*}
	\max_{z \in 3^n \Zd \cap \cu_m}
	\bigl|\b_{m}(z+\cu_n) \bigr| \shom_{m-1}^{-1}
	&
	\leq
	C
	+
	\O_{\Gamma_1}\bigl(C \cstar^{-1}\cgamma (m-n+1)^2   3^{\cgamma (m-n)} \bigr)
	\notag \\ &  \qquad
	+
	\O_{\Gamma_{\frac{1}{3}(1-  \sigma/4 )}} \Bigl(C  (m-n+1)^{3 + \sigma}
	\exp( -c E^{-2} \cgamma^{-1})  3^{\cgamma(m-n)} \Bigr)
	\,,
\end{align*}
The same summation for~$\Lambda_{s,1}$ gives a deterministic contribution bounded by~$C$, a~$\Gamma_1$-term of scale~$C\cstar^{-1}s\cgamma(2s-\cgamma)^{-3}\leq C\cstar^{-1}s^{-2}\cgamma$, and an exponentially small term multiplied by a fixed power of~$(2s-\cgamma)^{-1}$. Finally, after increasing~$C$ in~\eqref{e.percolation.admissibility}, that condition gives~$E^{-2}\cgamma^{-1}\geq\cgamma^{-\frac35}$ and~$\exp(-C^{-1}E^{-2}\cgamma^{-1})\leq\frac12\cgamma$. The remaining powers of~$(2s-\cgamma)^{-1}$ are therefore absorbed into the exponential. This proves both estimates for~$q=1$ and hence, by~\eqref{e.ellipticities.change.q}, the proposition.
\end{proof}

\begin{proof}[Proof of Proposition~\ref{p.bfA.multiscalebound}]
Fix~$P = (p,q)$ with~$p, q \in \Rd$ and~$|p|,|q| \leq 1$, and fix~$L,n \in \Z$ with~$L \geq m-1 \geq n$.
Fix~$b\coloneqq2^{-20}$ and a dimensional parameter~$\eta\in(0,1)$ to be fixed below, and set
\begin{equation*}
	k_0\coloneqq3\vee\left\lceil\eta E^{-2}\cgamma^{-1}\right\rceil\in\N\,.
\end{equation*}
The hypothesis~$\mathcal S(m-1,E)$ contains~$\mathcal S(n-1,E)$. After increasing the constant in~\eqref{e.Esigmares.in.proof} to include the fixed condition~$E\geq Cb^{-1}$, apply Proposition~\ref{p.minimal.scale.separation} with root~$\cu_n$ and denote its output by~$\hsep$. With this root-$n$ variable and~$k_0$, define~$\{h_k\}_{k \in \N_0}$ by~\eqref{e.hn.def} and construct~$\W(\cu_n)$ and~$\SW(\cu_n)$.

\smallskip

Let~$\mathcal{I}  \coloneqq  \{ \cu \in \W(\cu_n) : \mathcal{B}(\cu) \text{ occurs} \}$ be the set of bad cubes defined in~\eqref{e.bad.cubes.in.partition}. For~$\spx\in\SW(\cu_n)$, let~$Q(\spx)$ be its unique parent Whitney cube, write~$\mathcal B(\spx)\coloneqq\mathcal B(Q(\spx))$, and interpret~$\spx\in\mathcal N(\mathcal I)$ as~$Q(\spx)\in\mathcal N(\mathcal I)$.
We partition~$\mathcal I$ into its connected components.  Lemma~\ref{l.bad.clusters.geometry}
supplies the two-layer neighborhood and diameter hypotheses of
Lemma~\ref{l.piecewise.affine.approx}.  We thus obtain the overlap-safe family
$\{\hat{\linear}_p\}_{p\in\Rd}$ satisfying~\eqref{e.hat.linear.properties}.
Applying Lemma~\ref{l.subadd.betterer} yields
\begin{align}
	\label{e.sum-of-a-decomp}
	&\bigl|
	\bfAhom_{n-1}^{-\nf12} P \cdot
	\bfA_L(\cu_n)
	\bfAhom_{n-1}^{-\nf12}
	P
	\bigr|
	\notag  \\
	&\quad \leq
	\sum_{k\in\N}\sum_{\spx\in\SW(\cu_n,k)}
	\frac{|\spx|}{|\cu_n|}
	\bigl|
	\bfA^{\nf12}_L(\spx)
	\bfAhom_{n-1}^{-\nf12}
	P
	\bigr|^2
	\indc_{\{\neg \mathcal{B}(\spx)\}}
	\notag \\ & \qquad
	+
	\sum_{k\in\N}\sum_{\spx\in\SW(\cu_n,k)}
	\frac{|\spx|}{|\cu_n|}
	\Biggl|
	\bfA^{\nf12}_L(\spx)
	\bfAhom_{n-1}^{-\nf12}
	\begin{pmatrix}
		\nabla \hat{\linear}_p (\spx) \\
		(\nabla \cdot \hat{D}_q)(\spx)
	\end{pmatrix}
	\Biggr|^2
	\indc_{\{\neg \mathcal{B}(\spx)\}}
	\indc_{\{\spx \in \mathcal{N}(\mathcal{I})\}}
	\,\,,
\end{align}
where~$\hat{D}_q$ is defined in~\eqref{e.hatdq}. In the remainder of the proof, which we split into four steps, we bound the terms appearing on the right above.

\smallskip

\emph{Step 1. Good simplices.}
In this step we show that, for every~$k\in\N$,~$\spx\in\SW(\cu_n,k)$ and~$\hat P=(\hat p,\hat q)\in\R^{2d}$,
\begin{align}
	\label{e.good.simplex.consequence}
	\! \!
	\bigl|
	\bfA^{\nf12}_L(\spx)
	\bfAhom_{n-1}^{-\nf12}
	\hat{P}
	\bigr|^2
	\indc_{\neg \mathcal{B}(\spx)}
	&\leq C|\hat q|^2 3^{\cgamma(k+h_k)}
	+C|\hat p|^2\left(1+
	\frac{|(\k_L-\k_{n-k-h_k})_{\spx}|^2}
	{\shom_{n-k-h_k}\shom_{n-1}}\right)
	\notag\\
	&\leq  C  |\hat{q}|^2 3^{\cgamma(k + h_k)}
	+ C  |\hat{p}|^2  \bigl(1 + \cstar^{-1}\cgamma 3^{-2 \cgamma(n -k-h_k)}  | (\k_L - \k_{n-k-h_k})_{\spx}|^2  \bigr)
	\,\,.
\end{align}
Let~$z+\cu_j\in\W(\cu_n,k)$ be the cube containing~$\spx$, where~$j\coloneqq n-k-h_k$. By the induction hypothesis~\eqref{e.shom.h.bounds}
and the definition of~$\mathcal{B}_{\mathrm{osc}}$ in~\eqref{e.Bosc.def}, we have that the hypothesis~\eqref{e.J.sensitivity.smallness.condition} is satisfied for~$\spx \subset z + \cu_j$ with~$\neg \mathcal{B}(\spx)$, taking~$\h  \coloneqq  (\k_L - \k_j)$  and~$\a  \coloneqq  \a_j$. Indeed, by the definition of~$\mathcal{B}_{\mathrm{osc}}$ and~\eqref{e.delta.osc.choice}, when~$\neg \mathcal{B}_{\mathrm{osc}}(z + \cu_j)$ holds we have
\begin{align}
	3^{2 j} \|\nabla (\k_{L} - \k_{j}) \|_{\underline{W}^{1,\infty} (z + \cu_j)}  \indc_{\{\neg \mathcal{B}(\spx)\}}
	&
	\leq  \delta_{\mathrm{osc}} \cstar^{\nf12}  \cgamma^{-\nf12} 3^{\cgamma j}
	\notag \\ &
	\leq  2^{-5}   C^{-1}_{\eqref{e.J.sensitivity.smallness.condition} }  \shom_j
	\leq C_{\eqref{e.J.sensitivity.smallness.condition}}^{-1} \lambda_{\nf14,2}(z + \cu_j; \a_j) \,\,,
	\label{e.we.can.apply.cg}
\end{align}
where in the last inequality we used~$\neg \mathcal{B}_{\mathrm{loc}}(z + \cu_j)$. A triadic Whitney decomposition of~$\spx$, followed by subadditivity and~\eqref{e.bound.one.cube.by.lambdas}, bounds~$J(\spx,p,0;\a)$ and~$J(\spx,0,q;\a)$ by~$C|p|^2\Lambda_{\nf14,2}(z+\cu_j;\a)$ and~$C|q|^2\lambda_{\nf14,2}^{-1}(z+\cu_j;\a)$, respectively. For every cube~$U$ in this decomposition and~$\h=\k_L-\k_j$, the estimate~$|(\h)_U|\leq|(\h)_\spx|+C3^j\|\nabla\h\|_{L^\infty(z+\cu_j)}$ converts the cube-based sensitivity error into the simplex-mean error below. The boundary-layer sum converges since~$\nf14<\nf12$, and the scale loss is absorbed into~$C(d)$ since~$1+h_{k+1}-h_k\leq2$. The triangle inequality and~\eqref{e.J.by.means.of.bfA} give
\begin{align*}
	\bigl|
	\bfA^{\nf12}_L(\spx)
	\bfAhom_{n-1}^{-\nf12}
	\hat{P}
	\bigr|^2
	\indc_{\neg \mathcal{B}(\spx)}
	& \leq
	2 \biggl|
	\bfA_L^{\nf12}(\spx)
	\begin{pmatrix} \shom_{n-1}^{-\nf12} \hat{p} \\ 0
		\end{pmatrix} \biggr|^2
	\indc_{\neg \mathcal{B}(\spx)}
	+
	2 \biggl|
	\bfA_L^{\nf12}(\spx)
	\begin{pmatrix} 0 \\ \shom_{n-1}^{\nf12} \hat{q}
		\end{pmatrix} \biggr|^2
	\indc_{\neg \mathcal{B}(\spx)}
	\\
	&= 4 \bigl( J(\spx,\shom_{n-1}^{-\nf12} \hat{p}  , 0 \,; \a_L) +  J(\spx,0 , \shom_{n-1}^{\nf12} \hat{q} \,; \a_L) \bigr)\indc_{\neg \mathcal{B}(\spx)}
	\,\,.
\end{align*}
Using~\eqref{e.ellipticities.monotone.ordered} and the event~$\neg \mathcal{B}(z + \cu_j)$ (which by~\eqref{e.Bloc.def} ensures~$\Lambda_{\nf14,2}(z+\cu_j; \a_j) \leq 10 \shom_j$ and~$\lambda_{\nf14,2}^{-1}(z+\cu_j; \a_j) \leq 10 \shom_j^{-1}$), we have by subadditivity and Lemma~\ref{l.J.sensitivity}, together with~\eqref{e.Esigmares.in.proof}
and~\eqref{e.we.can.apply.cg}, that
\begin{align*}
	\lefteqn{
	J(\spx,\shom_{n-1}^{-\nf12} \hat{p}  , 0 \,; \a_L) \indc_{\neg \mathcal{B}(\spx)}
	}   \qquad &
	\notag \\ &
	\leq
	C \shom_{n-1}^{-1} |\hat{p}|^2 \Lambda_{\nf 14,2}(z+\cu_{j}\, ; \a_{L}) \indc_{\neg \mathcal{B}(\spx)}
	\notag \\ &
	\leq
	C  \shom_{n-1}^{-1} |\hat{p}|^2 \Lambda_{\nf 14,2}(z+\cu_{j}\, ; \a_{j})\indc_{\neg \mathcal{B}(\spx)}
	\notag \\ & \qquad
	+
	C \shom_{n-1}^{-1} |\hat{p}|^2
	\bigl( | (\k_L - \k_j)_{\spx}| + 3^{j}\| \nabla (\k_L - \k_j) \|_{L^\infty(z + \cu_{j})} \bigr)^2 \lambda_{\nf14,2}^{-1}(z + \cu_{j}\,; \a_{j})
	\indc_{\neg \mathcal{B}(\spx)}
	\notag \\ &
	\leq
	C  |\hat{p}|^2 \shom_j \shom_{n-1}^{-1} \bigl( 1
	+
	\shom_j^{-2} | (\k_L - \k_j)_{\spx}|^2 \bigr)
	\,.
\end{align*}
We also have
\begin{align*}
	J(\spx,0 , \shom_{n-1}^{\nf12} \hat{q} \,; \a_L) \indc_{\neg \mathcal{B}(\spx)}
	&
	\leq
	C |\hat{q}|^2 \lambda_{\nf 14,2}^{-1}(z+\cu_{j}\, ; \a_{L}) \indc_{\neg \mathcal{B}(\spx)}  \shom_{n-1}
	\notag \\ &
	\leq C |\hat q|^2\lambda_{\nf14,2}^{-1}(z+\cu_j;\a_j)
	\indc_{\neg\mathcal B(\spx)}\shom_{n-1}
	\notag \\ &
	\leq
	C |\hat{q}|^2 \shom_{n-1} \shom^{-1}_j
	\,\,.
\end{align*}
By the induction hypothesis~\eqref{e.shom.h.bounds} if~$n -1 \geq \mstar$ and~\eqref{e.basecase.diffusivity} otherwise,
\begin{equation*}
	\shom_{n-1} \shom^{-1}_j \leq C 3^{\cgamma(n-1-j)} \leq C 3^{\cgamma (k + h_k)}
	\qand
	\shom_j \shom_{n-1}^{-1}\leq C\,.
\end{equation*}
The previous four displays yield~\eqref{e.good.simplex.consequence}.

\smallskip

\emph{Step 2. Layerwise H\"older estimate.}
For every~$k\in\N$ and~$\spx\in\SW(\cu_n,k)$ with~$\spx\in\mathcal N(\mathcal I)$, Remark~\ref{r.gradient.bound.simplified} and~\eqref{e.hat.linear.properties} give
\begin{equation}
	| (\nabla \cdot \hat{D}_q)(\spx)|^2\leq C|q|^2 3^{2b(k+h_k)}
	\qand
	|\nabla\hat{\linear}_p(\spx)|^2\leq C|p|^2 3^{2b(k+h_k)}\,.
	\label{e.bounds.on.slopes.when.bad}
\end{equation}
For each~$k\in\N$, set
\begin{equation*}
	R_k\coloneqq\Biggl(\sum_{\spx\in\SW(\cu_n,k)}\frac{|\spx|}{|\cu_n|}
	\bigl(1+3^{8b(\hsep+k_0)}\indc_{\{\spx\in\mathcal N(\mathcal I)\}}\bigr)\Biggr)^{\! \nf14}\,.
\end{equation*}
By the construction of~$\W(\cu_n)$, the simplices in~$\SW(\cu_n,k)$ triangulate a boundary layer of relative volume at most~$C3^{-k}$. Hence, by H\"older's inequality,
\begin{equation*}
	\sum_{\spx\in\SW(\cu_n,k)}\frac{|\spx|}{|\cu_n|}
	\bigl(1+3^{2b(\hsep+k_0)}\indc_{\{\spx\in\mathcal N(\mathcal I)\}}\bigr)
	\leq C3^{-\frac 34 k}R_k
	\,.
\end{equation*}
By increasing the constant in~\eqref{e.Esigmares.in.proof}, we may assume~$\frac{4\cgamma+4b}{1-b}\leq\frac1{32}$. Applying H\"older's inequality and Jensen's inequality to~\eqref{e.sum-of-a-decomp}, and then using~\eqref{e.good.simplex.consequence},~\eqref{e.bounds.on.slopes.when.bad} and the preceding estimate, yields
\begin{align}
	\lefteqn{
	\bigl|\bfAhom_{n-1}^{-\nf12}P\cdot\bfA_L(\cu_n)\bfAhom_{n-1}^{-\nf12}P\bigr|
	}\quad&
	\notag\\ &
	\leq C|p|^2+C\bigl(3^{\cgamma(\hsep+k_0)}|q|^2+|p|^2\bigr)
	\biggl( \sum_{k\in\N}3^{-\frac1{16} k}R_k \biggr)
	\notag\\ &\quad
	+C\cstar^{-1}\cgamma |p|^2
	\biggl( \sum_{k\in\N}
	3^{-\frac1{16} k}R_k^4 \biggr)^{\! \nf 14}
	\biggl( \sum_{k\in\N} 3^{-\frac1{16} k}
	\bigl( 3^{-\cgamma(n-k-h_k)}
	\|\k_L-\k_{n-k-h_k}\|_{\underline L^4(\cu_n)} \bigr)^4 \biggr)^{\! \nf 12}
	\,.
	\label{e.sum.of.a.decomp.holdered}
\end{align}
In the last line we used H\"older's inequality with exponents~$4$,~$2$ and~$4$; the fourth power of its third factor is bounded by the summable series~$C\sum_{k\in\N}3^{-(\frac{13}{16}-\frac{8b}{1-b})k}$. We estimate the three sums in~\eqref{e.sum.of.a.decomp.holdered} in Steps~3 and~4.

\smallskip

\emph{Step 3. Wave oscillations.}
In this step, we show that
\begin{multline*}
	\Biggl(\sum_{k\in\N}3^{-\frac1{16} k}
	\bigl(3^{-\cgamma(n-k-h_k)}\|\k_L-\k_{n-k-h_k}\|_{\underline L^4(\cu_n)}\bigr)^4\Biggr)^{\!\nf12}
	\\
	\leq 3^{2\cgamma\hsep}\Bigl(
	C(k_0+1)3^{2\cgamma k_0}
	+\O_{\Gamma_1}\bigl(C(1+\min\{\cgamma^{-1},L-n\})3^{2\cgamma(L-n+k_0)}\bigr)
	+\O_{\Gamma_{(1+\sigma)^{-1}}}\bigl(C\sigma^{-2}\cgamma^{-1}3^{-\frac d4k_0}\bigr)
	\Bigr)\,.
\end{multline*}
We use the centered deep-band estimate~\eqref{e.km.kn.L4.centered}.

For each~$k\in\N$, set the deterministic index
\begin{equation*}
	\ell_k\coloneqq n-k-h_k+\hsep
	=n-k-\lceil b(1-b)^{-1}k\rceil-k_0\,.
\end{equation*}
We split
\begin{equation*}
	\k_L-\k_{\ell_k-\hsep}
	=(\k_L-\k_n)+(\k_n-\k_{\ell_k})+(\k_{\ell_k}-\k_{\ell_k-\hsep})\,.
\end{equation*}
The case~$L=n$ is trivial. Otherwise, by~\eqref{e.kmn.bounds} and~\eqref{e.powerofGammasigma},
\begin{equation*}
	3^{-2\cgamma n}\|\k_L-\k_n\|_{\underline L^4(\cu_n)}^2
	\leq\O_{\Gamma_1}\bigl(C\min\{\cgamma^{-1},L-n\}3^{2\cgamma(L-n)}\bigr)\,.
\end{equation*}
Finally, decompose according to~$\{\hsep=h\}$ and use~\eqref{e.km.kn.Lp},~\eqref{e.hsep.tails},~\eqref{e.indc.O.sigma},~\eqref{e.multGammasig} and~\eqref{e.powerofGammasigma}. Since~$n-\ell_k\geq k_0$, this yields
\begin{equation*}
	3^{-2\cgamma(\ell_k-\hsep)}\|\k_{\ell_k}-\k_{\ell_k-\hsep}\|_{\underline L^4(\cu_n)}^2
	\leq 3^{2\cgamma\hsep}\Bigl(
	C\hsep+\O_{\Gamma_{(1+\sigma)^{-1}}}\bigl(C\sigma^{-2}\cgamma^{-1}3^{-\frac d4k_0}\bigr)\Bigr)\,.
\end{equation*}
Since~$n-\ell_k=k+\lceil b(1-b)^{-1}k\rceil+k_0$, the smallness imposed after~\eqref{e.bounds.on.slopes.when.bad} gives
\begin{equation*}
	\Biggl(\sum_{k\in\N}3^{-\frac1{16}k}(n-\ell_k)^2 3^{4\cgamma(n-\ell_k)}\Biggr)^{\!\nf12}
	\leq C(k_0+1)3^{2\cgamma k_0}\,.
\end{equation*}
The same summation applied to the centered parts combines them into~$\O_{\Gamma_1}(C(1+\min\{\cgamma^{-1},L-n\})3^{2\cgamma(L-n+k_0)})$. Finally,~\eqref{e.hsep.tails} gives~$\hsep=\O_{\Gamma_1}(C)$, so the term~$C\hsep$ is absorbed into this random term. The generalized triangle inequality now gives the claim.

\smallskip

\emph{Step 4. Layer estimates and conclusion.}
The relative-volume bound~$C3^{-k}$ for the~$k$-th boundary layer and Lemma~\ref{l.bad.cubes.density}, together with the bounded overlap of touching Whitney cubes, give
\begin{equation*}
	R_k^4
	\leq C3^{-k}\Biggl(1+3^{8b(\hsep+k_0)}
	\min\Bigl\{1,3^k\bigl(\exp(-cE^{-2}\cgamma^{-1})+3^{-\frac12(k+h_k)}\bigr)\Bigr\}\Biggr)
	\,.
\end{equation*}
Taking fourth roots and summing, and also summing the preceding estimate before taking the fourth root, using~$h_k=\lceil b(1-b)^{-1}k\rceil+\hsep+k_0$, yield, after reducing~$c$,
\begin{align*}
	\sum_{k\in\N}3^{-\frac1{16}k}R_k
	+\Biggl(\sum_{k\in\N}3^{-\frac1{16}k}R_k^4\Biggr)^{\!\nf14}
	&\leq C\Bigl(1+3^{2b(\hsep+k_0)}
	\bigl(\exp(-cE^{-2}\cgamma^{-1})+3^{-\frac18 k_0}\bigr)\Bigr)
	\,.
\end{align*}
We now fix~$\eta=\eta(d)>0$ sufficiently small. Since~$b=2^{-20}$ and~$\cgamma k_0\leq\eta E^{-2}+3\cgamma\leq CE^{-2}$, after reducing~$c$ we have
\begin{equation*}
	3^{2bk_0}\left(\exp(-cE^{-2}\cgamma^{-1})+3^{-\frac18 k_0}\right)
	\leq C\exp(-cE^{-2}\cgamma^{-1})\,.
\end{equation*}
By~\eqref{e.hsep.tails} and~\eqref{e.powerofGammasigma}, after reducing~$c$ and increasing the constant in~\eqref{e.Esigmares.in.proof}, this also gives
\begin{equation*}
	3^{2b(\hsep+k_0)}\left(\exp(-cE^{-2}\cgamma^{-1})+3^{-\frac18 k_0}\right)
	\leq\O_{\Gamma_{\frac12(1-\sigma)}}\bigl(C\exp(-cE^{-2}\cgamma^{-1})\bigr)\,.
\end{equation*}
By~\eqref{e.shom.h.bounds} and~\eqref{e.basecase.diffusivity}, we have~$3^{-\cgamma(m-n)}\shom_{m-1}\shom_{n-1}^{-1}\leq C$.
Setting~$p=0$ in~\eqref{e.sum.of.a.decomp.holdered}, using the first term on the left side of the layer estimate above and taking the supremum over~$|q|\leq1$ therefore yield
\begin{align*}
	3^{-\cgamma(m-n)}\shom_{m-1}|\s_{L,*}^{-1}(\cu_n)|
	&\leq C3^{\cgamma(\hsep+k_0)}
	\Bigl(1+3^{2b(\hsep+k_0)}\bigl(\exp(-cE^{-2}\cgamma^{-1})+3^{-\frac18k_0}\bigr)\Bigr)
	\\
	&\leq C3^{\cgamma(\hsep+k_0)}
	\Bigl(1+3^{2b\hsep}\exp(-cE^{-2}\cgamma^{-1})\Bigr)\,.
\end{align*}
For a sufficiently small dimensional constant~$\kappa>0$, a split at~$\cgamma\hsep=\kappa\cstar$ and~\eqref{e.hsep.tails} give
\begin{equation*}
	3^{\cgamma\hsep}\leq C+3^{\cgamma\hsep}\indc_{\{\cgamma\hsep>\kappa\cstar\}}
	\leq C+\O_{\Gamma_2}\bigl(\exp(-cE^{-2}\cgamma^{-1})\bigr)\,.
\end{equation*}
Here we used~\eqref{e.indc.O.sigma} and~\eqref{e.multGammasig}. Applying~\eqref{e.hsep.tails} once more to~$3^{2b\hsep}$ and relabelling~$\sigma$, the two previous estimates yield~\eqref{e.slstar.multiscale}.

\smallskip

For the upper-left block, retain the first inequality in~\eqref{e.good.simplex.consequence}. If~$j=n-k-h_k$, then~\eqref{e.shom.h.bounds} and~\eqref{e.basecase.diffusivity} give
\begin{equation*}
	\frac1{\shom_j\shom_{m-1}}
	\leq C\cstar^{-1}\cgamma3^{-\cgamma(j+m-1)}
	\leq C\cstar^{-1}\cgamma
	3^{-\cgamma(m-n)}3^{-2\cgamma j}\,,
\end{equation*}
where the last inequality uses~$j\leq n-1$. Thus the wave part of the layerwise H\"older estimate involves the same~$R_k^4$-sum and quartic wave sum as~\eqref{e.sum.of.a.decomp.holdered}, with the additional factor~$3^{-\cgamma(m-n)}$.

Set~$q=0$, take the supremum over~$|p|\leq1$, and use the wave estimate from Step~3 and the layer estimate above. On the event~$\{\cgamma\hsep\leq\kappa\cstar\}$, the factor~$3^{2\cgamma\hsep}$ is bounded, and the deterministic deep-band term and the terminal random term are controlled by
\begin{equation*}
	C\cstar^{-1}\cgamma(k_0+1)3^{2\cgamma k_0}\leq C
	\qand
	C\cstar^{-1}\sigma^{-2}3^{-\frac d4k_0}\leq C\exp(-cE^{-2}\cgamma^{-1})\,.
\end{equation*}
The part of the layer estimate without the exponentially small multiplier combines with the~$\Gamma_1$-term in Step~3 to give the first random term below. For the remaining part, the preceding~$\Gamma_{\frac12(1-\sigma)}$ estimate and~\eqref{e.multGammasig} give an exponent~$\frac{1-\sigma}{3-\sigma}\geq\frac13(1-\sigma)$. Its scale is uniform in~$L$ since~$\cgamma(1+\min\{\cgamma^{-1},L-n\})\leq C$, and~$\cstar^{-1}$ is absorbed into the exponential using~$\cstar^{-1}\leq CE$ and~$E^{-2}\cgamma^{-1}\geq E^3$. On the complementary event,~\eqref{e.hsep.tails} absorbs all three wave terms into the same random term.

The nonwave terms are bounded since~$\shom_{n-1}\shom_{m-1}^{-1}\leq C$. Applying~\eqref{e.Gamma.sigma.triangle},~\eqref{e.multGammasig} and~\eqref{e.powerofGammasigma}, and relabelling~$\sigma$, we obtain
\begin{align*}
	\shom^{-1}_{m-1}|\b_L(\cu_n)|
	&\leq C
	+\O_{\Gamma_1}\bigl(C\cstar^{-1}\cgamma(L-n+1)3^{2\cgamma(L-n)-\cgamma(m-n)}\bigr)
	\notag\\ &\qquad
	+\O_{\Gamma_{\frac13(1-\sigma)}}\bigl(3^{2\cgamma(L-n)-\cgamma(m-n)}\exp(-cE^{-2}\cgamma^{-1})\bigr)
	\,.
\end{align*}
This is~\eqref{e.bL.multiscale}.
\end{proof}

\subsubsection{The event that enables coarse-grained sensitivity}
\label{sss.cg.sensitivity.event}

Later in this section, we will need to use the coarse-grained sensitivity estimates, Lemmas~\ref{l.lambda.sensitivity} and~\ref{l.J.sensitivity}.
To that end, we define, for each~$m,n \in \Z$ and~$z \in \Rd$, the event
\begin{multline}
	\mathcal{Q}(m,n,z) \coloneqq
	\Bigl\{
	C_{\eqref{e.J.sensitivity.smallness.condition}\&\eqref{e.lambda.sensitivity.smallness.condition}}
	\sup_{L \geq n} 3^{2m} \| \nabla (\k_L-\k_{n}) \|_{\underline{W}^{1,\infty}(z+\cu_m)}
	\leq
	\frac{1}{2} C_{\eqref{e.cg.ellip.lower}}^{-1}
	3^{-\frac12(n-m)_+}
	\shom_{n-1}
	\\
	\leq
	3^{-\frac14(n-m)_+}
	\lambda_{\nf18,2}(z{+}\cu_m;\a_{n})
	\Bigr\}
	\,.
	\label{e.good.local.events}
\end{multline}
The event is a modification of the one defined in~\eqref{e.Bosc.def} and allows us to switch from the field~$\a_L$ to~$\a_{n}$ in the cube~$z+\cu_m$ for any~$L \geq n$, using the coarse-grained sensitivity estimates.

\begin{lemma}[Estimate of the bad event]
\label{l.bad.event.lemma}
Assume~$m_0 \in\Z$ and~$E \in [1,\infty)$ are such that~$\S(m_0-1, E)$ is valid.
There exists~$c(d)>0$ such that, if
\begin{equation}
	\label{e.bad.event.admissibility}
	c^{-1}
	\cstar^{-1}
	\leq
	E
	\leq \cgamma^{-\nf15}
	\,,
\end{equation}
then for each~$m,n \in \Z$ with~$n \leq m_0 - 1$ and~$z \in \Rd$ such
that, whenever~$n<m$,
\begin{equation*}
	3^{5(m-n)}\leq c\cstar\cgamma^{-1}
	\quad\hbox{and}\quad
	d(m-n)\log3
	+\exp\Bigl(\frac18 C_{\eqref{e.cg.ellip.lower}}^{-1}E^{-2}\cgamma^{-1}\Bigr)
	\leq\exp\Bigl(\frac14 C_{\eqref{e.cg.ellip.lower}}^{-1}E^{-2}\cgamma^{-1}\Bigr)\,,
\end{equation*}
we have
\begin{equation}
	\label{e.bad.event.Q.estimate}
	\P \bigl[ \neg \mathcal{Q}(m,n,z) \bigr]
	\leq
	\exp \Bigl( - c \cstar \cgamma^{-1} 3^{-5 (m-n)_+} 3^{(n-m)_+} \Bigr)
	+
	\exp\bigl(-\exp\bigl(c E^{-2} \cgamma^{-1}\bigr)  \bigr)  \,\,.
\end{equation}
\end{lemma}
\begin{proof}
By stationarity, it suffices to prove the case~$z = 0$. Using the lower bound for~$\shom_{n-1}$ in~\eqref{e.shom.h.bounds}, we decompose the event
\begin{equation*}
	\neg \mathcal{Q}(m,n,0)
	\subseteq \mathcal{B}_{\mathrm{osc}}(m,n)
	\cup
	\Bigl\{
	\lambda^{-1}_{\nf18,2}(\cu_m;\a_{n})
	\shom_{n-1}
	>2 \cdot 3^{\frac14 (n-m)_+} C_{\eqref{e.cg.ellip.lower}}
	\Bigr\} \,\,,
\end{equation*}
where we define
\begin{equation*}
	\mathcal{B}_{\mathrm{osc}}(m,n)
	\coloneqq
	\biggl\{
	\sup_{L \geq n} 3^{2 m} \|\nabla (\k_{L} - \k_{n}) \|_{\underline{W}^{1,\infty} (\cu_m)}
	>  \delta_{\mathrm{osc}}
	\cstar^{\nf12} \cgamma^{-\nf12} 3^{\cgamma n}
	3^{-\frac12 (n- m)_+}
	\biggr\} \,\,.
\end{equation*}
We now bound the probability of the second event. After reducing~$c$ in~\eqref{e.bad.event.admissibility}, that condition and~$\cstar\leq\frac32$ ensure the admissibility of Proposition~\ref{p.cg.ellipticity.bounds} with~$\sigma=\nf12$ and also give~$\cgamma\leq\nf1{16}$. We consider two cases.

\smallskip

\emph{Case~$n \geq m$:} We have~$\cu_{m \vee n} = \cu_n$. By Proposition~\ref{p.cg.ellipticity.bounds} applied at scale~$n$ with~$\sigma=\nf12$ and~$s=\nf18$ (which is valid since~$n \leq m_0 - 1$ and~$\S(m_0-1, E)$ implies~$\S(n-1, E)$), we obtain
\begin{equation*}
	\P\bigl[ \lambda^{-1}_{\nf18,2}(\cu_n;\a_{n})
	\shom_{n-1}
	> 2 C_{\eqref{e.cg.ellip.lower}} \bigr]
	\leq
	\exp\bigl(-\exp\bigl(c E^{-2} \cgamma^{-1}\bigr)  \bigr) \,\,.
\end{equation*}
Since~$(n-m)_+ = n - m \geq 0$, the factor~$3^{\frac14(n-m)_+}$ provides additional room, and the desired bound follows.

\smallskip

\emph{Case~$n < m$:} We have~$\cu_{m \vee n} = \cu_m$ and~$(n-m)_+ = 0$. By the sharp finite-grid subadditivity estimate,
\begin{equation*}
	\lambda^{-1}_{\nf18,2}(\cu_m;\a_{n})
	\leq
	2\max_{z \in 3^n \Zd \cap \cu_m} \lambda^{-1}_{\nf1{16},2}(z + \cu_n;\a_{n})\,.
\end{equation*}
For each~$z \in 3^n \Zd \cap \cu_m$, Proposition~\ref{p.cg.ellipticity.bounds} at scale~$n$ with~$\sigma=\nf12$ and~$s=\nf1{16}$ (valid since~$\S(n-1, E)$ holds) gives
\begin{equation*}
	\shom_{n-1} \lambda^{-1}_{\nf1{16},2}(z + \cu_n;\a_{n})
	\leq
	C_{\eqref{e.cg.ellip.lower}}
	+
	\O_{\Gamma_{\nf14}} \bigl( \exp \bigl( - c E^{-2} \cgamma^{-1}\bigr)  \bigr)\,.
\end{equation*}
Taking the maximum over~$|3^n \Zd \cap \cu_m| = 3^{d(m-n)}$ cubes and using~\eqref{e.maxy.bound}, we obtain
\begin{equation*}
	\P\bigl[ \lambda^{-1}_{\nf18,2}(\cu_m;\a_{n})
	\shom_{n-1}
	> 2 C_{\eqref{e.cg.ellip.lower}} \bigr]
	\leq
	\exp\bigl(-\exp\bigl(c E^{-2} \cgamma^{-1}\bigr)  \bigr) \,\,,
\end{equation*}
where the factor~$3^{d(m-n)}$ from the union bound is absorbed using the
hypothesis on~$m-n$ in the statement.

\smallskip

Combining both cases, we conclude that
\begin{equation*}
	\P\bigl[ \lambda^{-1}_{\nf18,2}(\cu_{m \vee n};\a_{n})
	\shom_{n-1}
	> 2 \cdot 3^{\frac14 (n-m)_+} C_{\eqref{e.cg.ellip.lower}} \bigr]
	\leq
	\exp\bigl(-\exp\bigl(c E^{-2} \cgamma^{-1}\bigr)  \bigr) \,\,.
\end{equation*}
We claim that for every~$m,n \in \Z$,
\begin{equation}
	\P[\mathcal{B}_{\mathrm{osc}}(m,n)]  \leq
	\exp \Bigl( - c \cstar \cgamma^{-1} 3^{-5 (m-n)_+} 3^{(n-m)_+} \Bigr)
	\,.
	\label{e.oscillation.bound}
\end{equation}
It remains to prove~\eqref{e.oscillation.bound}. Set~$r\coloneqq m\vee n$ and decompose
into the genuine lower shells and shells strictly above~$r$:
\begin{equation*}
	\mathcal{B}_{\mathrm{osc}}(m,n)
	\subseteq
	\biggl\{
	3^{2m}\sum_{J=n+1}^{m}
	\|\nabla\mathbf j_J\|_{\underline W^{1,\infty}(\cu_m)}
	>  \frac12\delta_{\mathrm{osc}}
	\cstar^{\nf12} \cgamma^{-\nf12} 3^{\cgamma n}
	\biggr\}
	\cup
	\bigcup_{L \in (r,\infty) \cap \Z}
	\mathcal{B}_{\mathrm{osc},L}(m,n) \,\,,
\end{equation*}
where the lower-shell event is empty when~$n\geq m$, and
where
\begin{equation*}
	\mathcal{B}_{\mathrm{osc},L}(m,n)
	\coloneqq
	\bigl\{
	3^{2m} \|\nabla \mathbf{j}_L \|_{\underline{W}^{1,\infty} (\cu_m )}
	>   \frac1{10}\delta_{\mathrm{osc}}
	\cstar^{\nf12}
	\cgamma^{-\nf12} 3^{\cgamma n} 3^{-\frac15(L-r)}
	3^{-\frac12 (n-m)_+}
	\bigr\} \,\,.
\end{equation*}
When~$n \geq m$ the first event on the right is empty. When~$n < m$, by~\eqref{e.W.1.inf.bound} and the first additional gate in the statement we have
\begin{align*}
	\P\bigl[
	3^{2m}\sum_{J=n+1}^{m}
	\|\nabla\mathbf j_J\|_{\underline W^{1,\infty}(\cu_m)}
	>  \frac12\delta_{\mathrm{osc}}
	\cstar^{\nf12} \cgamma^{-\nf12} 3^{\cgamma n}
	\bigr]
	&\leq
	\exp \Bigl( - c \cstar \cgamma^{-1}
	3^{-(4 + 2 \cgamma) (m-n)}
	(m-n)^{-2}
	\Bigr) \\
	&\leq
	\exp \Bigl( - c \cstar \cgamma^{-1}
	3^{-5 (m-n)}
	\Bigr) \,\,.
\end{align*}
On the other hand, for~$L>r$, by~\eqref{e.nabla.jk.O} and~\eqref{e.maxy.bound},
\begin{equation*}
	\P\bigl[
	\mathcal{B}_{\mathrm{osc},L}(m,n)
	\bigr]
	\leq
	\exp \Bigl( - c \cstar \cgamma^{-1}
	3^{2(1 - \cgamma - \frac{1}{5})(L-r) }
	3^{(n-m)_+}
	\Bigr)
	\leq
	\exp \Bigl( - c \cstar \cgamma^{-1}
	3^{(L-r)}3^{(n-m)_+}
	\Bigr) \,\,.
\end{equation*}
Summing over~$L>r$ proves~\eqref{e.oscillation.bound} and completes the proof.
\end{proof}

\subsection{Propagation of the effective diffusivity bounds}
\label{ss.propagate.diffusivity.bounds}

In this subsection we propagate the upper and lower bound on~$\shom_m$ by proving the implication
\begin{equation}
	\mathcal{S}(m_0-1,E)
	\quad\implies \quad
	\text{\eqref{e.shom.h.bounds} is valid for~$m = m_0$}\,.
	\label{e.propagation.of.indyhyp.shombounds}
\end{equation}
In fact, we obtain a much more precise estimate for~$\shom_m$ than the crude bound~\eqref{e.shom.h.bounds} needed to propagate the induction hypothesis.
In Lemma~\ref{l.approximate.recurrence.formula}, we derive an approximate recurrence formula for the sequence~$\shom_m$ which roughly states that
\begin{equation*}
	\shom_{n+h}
	=
	\shom_{n} + (\log 3) \cstar  \shom_{n}^{-1} \sum_{k=n+1}^{n+h} 3^{2\cgamma k}
	+
	\text{ small error}\,.
\end{equation*}
This recurrence is then carefully analyzed in Lemma~\ref{l.integrate.approx.recurrence}, where we integrate it to obtain the asymptotic
\begin{equation*}
	\shom_{m}
	=
	\bigl( \nu^2
	+
	\cstar\cgamma^{-1}3^{2\cgamma m}
	\bigr)^{\nf12}
	+
	\text{ small error}  \,.
\end{equation*}
The precise version of the latter estimate is stated in the following proposition.

\begin{proposition}[{Asymptotics for~$\shom_m$}]
\label{p.propagate.diffusivity.lower.bound}
Assume~$m_0 \in (\mstarstar,\infty) \cap \Z$ and~$E\in [1,\infty)$ are such that~$\mathcal{S}(m_0-1,E)$ is valid.
There exists a constant~$C(d)<\infty$ such that, if
\begin{equation}
	\label{e.lower.bound.cgamma.cond}
	E \geq C \cstar^{-1}
	\qand
	\cgamma \leq E^{-10}\,,
\end{equation}
then
\begin{equation}
	\bigl| \shom_m - (\nu^2 + \cstar \cgamma^{-1} 3^{2\cgamma m})^{\nf12} \bigr|
	\leq
	C\cstar^{-1} E\cgamma^{\nf12}\left|\log \cgamma\right|
	\shom_m
	\,, \quad \forall m\in\Z \cap (-\infty,m_0]\,.
	\label{e.shom.m.flow}
\end{equation}
Moreover,~\eqref{e.shom.h.bounds} is valid for~$m = m_0$.
\end{proposition}

The approximate recurrence formula is presented below in Lemma~\ref{l.approximate.recurrence.formula}. Considered together, the two one-sided bounds in~\eqref{e.what.do.we.have} can be compared with~\cite[Proposition 7.2]{ABK.SD}, where we proved an analogous recurrence formula in the case~$\cgamma=0$.

\smallskip

The lemma is proved by a quantitative version of the recurrence heuristic for the running diffusivity: we compare the operator at scale~$m$ to the one at scale~$m-h$, where the only new randomness is the ``fresh shell''~$\h\coloneqq\k_m-\k_{m-h}$. The step size~$h$ can be selected, but must be no larger than~$\cstar\cgamma^{-1}$. On scales~$\gg 3^{m-h}$, the field~$\a_{m-h}$ has already been absorbed into the constant operator~$\shom_{m-h}\Delta$, so inserting the shell amounts (at leading order) to perturbing~$\shom_{m-h}\Id$ by the antisymmetric matrix field~$\h$. If we look for the coarse solution with affine boundary data~$\linear_e$ in the form~$\linear_e+w$, then linearizing
\begin{equation*}
	-\nabla\!\cdot\!\bigl((\shom_{m-h}\Id+\h)\nabla(\linear_e+w)\bigr)=0
	\,.
\end{equation*}
around~$\shom_{m-h}\Delta$ gives the Poisson equation~$-\Delta w=\nabla\!\cdot(\shom_{m-h}^{-1}\h e)$. This is exactly why the proof uses both the Dirichlet and Neumann perturbative correctors, denoted respectively by~$w_{D,e}^{(K)}$ and~$w_{N,e'}^{(K)}$ in~\eqref{e.def.w}: they are the explicit first-order responses to the newly added shell, with the large auxiliary scale~$K\gg m$ ensuring that boundary effects are negligible. The quantity which drives the flow is second order in~$\h$ and is encoded by the energies of these two correctors. For~$|e|=|e'|=1$, the computations~\eqref{e.perturb.assumption} and~\eqref{e.perturb.assumption.two} show that both expected energies converge, as~$K\to\infty$, to
\begin{equation*}
	(\log 3)\cstar\,\shom_{m-h}^{-2}\sum_{k=m-h+1}^m 3^{2\cgamma k}\,.
\end{equation*}
This is precisely the leading increment appearing in~\eqref{e.what.do.we.have}.

\smallskip

The proof of the lemma turns this linear-response computation into two one-sided inequalities for~$\shom_m$ itself. One probes the block energy with the vector~$P=\bfAhom_{m-h}^{-\nf12}\binom{e'}{e}$ (see~\eqref{e.recurrence.P.def}) and then localizes the corresponding minimization problem by tiling~$\cu_K$ into mesoscopic cubes~$z+\cu_n$, where~$n=m-h-O(|\log\cgamma|)$ creates a buffer scale between~$3^{m-h}$ and the localization scale. On each cube, the test field is split as~$X_z=S_z+\tilde S_z$, where~$S_z$ is the local minimizer at the almost affine vector~$P_z$ and~$\tilde S_z$ carries only the mean-zero fluctuation~$\mathbf F_z$. Step~2 shows that this fluctuation cost is negligible thanks to the scale separation (cf.~\eqref{e.lower.bound.oscillations}). The principal contribution is then reduced to~$\E[P_z\cdot \bfA_m(z+\cu_n)P_z]$, and here the sensitivity input enters in its sharpest form: on the good event~$\mathcal Q_z=\mathcal Q(n,m-h,z)$ one may replace~$\a_m$ by~$\a_{m-h}$ \emph{after gauging out the local average of the inserted shell}. Concretely, this is where the conjugation by~$\G_{-(\h)_{z+\cu_n}}$ appears in Step~3 and~\eqref{e.lower.bound.principal.one}: it implements, at the block-matrix level, the fact that subtracting a constant antisymmetric matrix corresponds to a simple linear transformation of the~$(p,q)$ variables. After averaging and switching~$\bfA_{m-h}(z+\cu_n)$ to its annealed limit~$\bfAhom_{m-h}$, with a controlled error, the proof boils down to estimating deterministic quadratic forms of
$\bfAhom_{m-h}^{\nf12}\G_{-(\h)_{z+\cu_n}}\bfAhom_{m-h}^{-\nf12}$ on the two pure inputs~$\binom{e'}0$ and~$\binom0e$, perturbed respectively by the averaged flux of the Neumann corrector and the averaged gradient of the Dirichlet corrector. Step~4 treats the first input; the gauge removes the explicit shell term from the Neumann flux, and~\eqref{e.perturb.assumption.two} yields the upper inequality for~$\shom_m\shom_{m-h}^{-1}$. Step~5 treats the second input, where~\eqref{e.what.nablaw.really.is} produces the compensating negative contribution~$-\|\nabla w_{D,e}^{(K)}\|_{\underline L^2}^2$ and hence the lower inequality for~$\shom_{m-h}\shom_m^{-1}$.


\begin{lemma}[Approximation of stationary fields by local test fields]
\label{l.approximation.stationary.by.local}
Assume that the translation action on~$\Omega$ is ergodic. Let~$\mathbf{p} \in L^2(\Omega;\Rd) \cap \Lpot$ be a stationary random potential field and~$\mathbf{j} \in L^2(\Omega;\Rd) \cap \Lsol$ a stationary random solenoidal field such that
\begin{equation*}
	\E[\mathbf p]=\E[\mathbf j]=0\,.
\end{equation*}
There exists~$C(d)<\infty$ such that the following hold.

\begin{enumerate}
\item[\rm{(i)}] Potential approximation. For every~$L\in\N$, there exists a family~$\{ v_{K,L} \}_{K \in \N}$ with~$ v_{K,L} \in H^1_0(\cu_K)$ such that
\begin{equation}
	\label{e.approximation.stationary.potential.by.dirichlet}
	\limsup_{K\to\infty}
	\E\bigl[
	\|\nabla v_{K,L}-\mathbf{p}\|^2_{\underline L^2(\cu_K)}
	\bigr]
	\leq
	C 3^{-L}\E[|\mathbf{p}|^2] \,.
\end{equation}

\item[\rm{(ii)}] Solenoidal approximation.
For every~$L\in\N$, there exists a family~$\{ \mathbf{j}_{K,L}\}_{K \in \N}$ with~$\mathbf{j}_{K,L} \in \Lsolo(\cu_K)$ such that
\begin{equation}
	\label{e.approximation.stationary.solenoidal.by.neumann}
	\limsup_{K\to\infty}
	\E\bigl[
	\|\mathbf{j}_{K,L}-\mathbf{j}\|^2_{\underline L^2(\cu_K)}
	\bigr]
	\leq
	C3^{-L}\E[|\mathbf{j}|^2]
	\,.
\end{equation}
\end{enumerate}
\end{lemma}

\begin{proof}
By the stationary Weyl decomposition, see~\cite[Chapter~7]{JKO94} and~\cite[Section~1.3]{AK.Book}, there exist random fields~$\phi$ and~$S$, whose realizations belong to~$H^1_{\loc}(\Rd)$ and~$H^1_{\loc}(\Rd;\R^{d\times d})$, respectively, such that~$S$ is antisymmetric,~$\nabla S$ is stationary and square integrable, and, almost surely,
\begin{equation*}
	\nabla\phi=\mathbf p\,,
	\qquad
	\mathbf j_i=\sum_{m=1}^d\partial_mS_{im}\,,
	\quad i\in\{1,\ldots,d\}\,.
\end{equation*}
The representation of~$\mathbf j$, rather than the sublinearity, is where the ergodic input enters.
The sublinearity of mean-zero stationary cocycles gives
\begin{equation*}
	\lim_{K\to\infty}3^{-2K}\E\biggl[\fint_{\cu_K}\Bigl(\bigl|\phi-(\phi)_{\cu_K}\bigr|^2+\bigl|S-(S)_{\cu_K}\bigr|^2\Bigr)\biggr]=0\,.
\end{equation*}
Indeed, the cited sublinearity is initially almost sure. Poincar\'e's inequality bounds its normalized left side by~$C\fint_{\cu_K}(|\mathbf p|^2+|\nabla S|^2)$, and these spatial averages are uniformly integrable by stationarity, Jensen's inequality, and the de la Vall\'ee--Poussin criterion. Vitali's theorem therefore gives the displayed convergence in expectation.
Fix~$K,L\in\N$ and choose~$\eta_{K,L}\in C_c^\infty(\cu_K)$ such that~$0\leq\eta_{K,L}\leq1$,~$\eta_{K,L}(x)=1$ whenever~$\dist(x,\partial\cu_K)\geq C3^{K-L}$, and~$|\nabla\eta_{K,L}|\leq C3^{L-K}$. Set
\begin{equation*}
	v_{K,L}\coloneqq\eta_{K,L}\bigl(\phi-(\phi)_{\cu_K}\bigr)
	\qand
	\mathbf j_{K,L,i}\coloneqq\sum_{m=1}^d\partial_m\Bigl(\eta_{K,L}\bigl(S_{im}-(S_{im})_{\cu_K}\bigr)\Bigr)\,.
\end{equation*}
Then~$v_{K,L}\in H^1_0(\cu_K)$ and~$\mathbf j_{K,L}\in\Lsolo(\cu_K)$. The boundary-strip estimate and stationarity give
\begin{align*}
	\E\bigl[\|\nabla v_{K,L}-\mathbf p\|_{\underline L^2(\cu_K)}^2\bigr]
	&\leq C3^{-L}\E\bigl[|\mathbf p|^2\bigr]+C3^{2(L-K)}\E\biggl[\fint_{\cu_K}\bigl|\phi-(\phi)_{\cu_K}\bigr|^2\biggr]\,,
	\\
	\E\bigl[\|\mathbf j_{K,L}-\mathbf j\|_{\underline L^2(\cu_K)}^2\bigr]
	&\leq C3^{-L}\E\bigl[|\mathbf j|^2\bigr]+C3^{2(L-K)}\E\biggl[\fint_{\cu_K}\bigl|S-(S)_{\cu_K}\bigr|^2\biggr] \,.
\end{align*}
Taking~$K\to\infty$ proves both assertions.
\end{proof}


\begin{lemma}
\label{l.corrector.limit}
Assume that the translation action on~$\Omega$ is ergodic. Let~$\f \in L^2(\Omega; \Rd)$ be a stationary random vector field satisfying~$\E[\f]=0$, and let~$w_{D}^{(K)} \in H^1_0(\cu_K)$ and~$w_{N}^{(K)} \in H^1(\cu_K)$ be, respectively, the weak Dirichlet solution and a weak Neumann solution, the latter unique modulo constants, of
\begin{equation}
	\label{e.corrector.limit.pde}
	\left\{
	\begin{aligned}
		& -\Delta w_{D}^{(K)} = \nabla \cdot \f & \text{in} & \ \cu_K\,,\\
		& w_{D}^{(K)} = 0 & \text{on} & \ \partial\cu_K\,,
	\end{aligned}
	\right.
	\qquad
	\left\{
	\begin{aligned}
		& -\Delta w_{N}^{(K)} = \nabla \cdot \f & \text{in} & \ \cu_K\,,\\
		& \mathbf{n}_{\cu_K}\cdot(\nabla w_{N}^{(K)} + \f) = 0
		& \text{on} & \ \partial\cu_K\,,
	\end{aligned}
	\right.\,,
\end{equation}
Here~$\mathbf n_{\cu_K}$ denotes the outward unit normal, and the Neumann boundary condition is understood through the weak formulation.
Let~$w$ be the unique (up to constants) solution of
$-\Delta w = \nabla \cdot \f$ in~$\Rd$
with~$\nabla w \in \Lpot$ stationary and~$\E[\nabla w]=0$.
Then the limits below exist and
\begin{equation}
	\label{e.corrector.limit}
	\lim_{K \to \infty}
	\E\bigl[\|\nabla w_{D}^{(K)}\|^2_{\underline{L}^2(\cu_K)}\bigr]
	=
	\lim_{K \to \infty}
	\E\bigl[\|\nabla w_{N}^{(K)}\|^2_{\underline{L}^2(\cu_K)}\bigr]
	=
	\E\bigl[|\nabla w|^2\bigr]
	\,.
\end{equation}
\end{lemma}

\begin{proof}
Set~$\mathbf p\coloneqq\nabla w$ and~$\mathbf j\coloneqq\f+\mathbf p$. Then~$\mathbf p\in\Lpot$,~$\mathbf j\in\Lsol$,~$\E[\mathbf p]=\E[\mathbf j]=0$, and
\begin{equation}
	\label{e.pot.sol.decomposition.corrector.limit}
	\E[\f\cdot\mathbf p]=-\E[|\mathbf p|^2]\,.
\end{equation}
For~$v\in H^1(\cu_K)$, set
\begin{equation*}
	I_K(v)\coloneqq\fint_{\cu_K}\biggl(\frac12|\nabla v|^2+\f\cdot\nabla v\biggr)\,.
\end{equation*}
The functions~$w_D^{(K)}$ and~$w_N^{(K)}$ minimize~$I_K$ over~$H^1_0(\cu_K)$ and~$H^1(\cu_K)/\R$, respectively, and
\begin{equation}
	\label{e.energy.identities.corrector.limit}
	I_K(w_D^{(K)})=-\frac12\|\nabla w_D^{(K)}\|^2_{\underline L^2(\cu_K)}
	\qqand
	I_K(w_N^{(K)})=-\frac12\|\nabla w_N^{(K)}\|^2_{\underline L^2(\cu_K)}
	\,.
\end{equation}
In particular,
\begin{equation}
	\label{e.ordering.corrector.limit}
	\|\nabla w_D^{(K)}\|^2_{\underline L^2(\cu_K)}
	\leq
	\|\nabla w_N^{(K)}\|^2_{\underline L^2(\cu_K)}
	\,.
\end{equation}
Fix~$L\in\N$, and let~$v_{K,L}$ be given by Lemma~\ref{l.approximation.stationary.by.local}\rm{(i)} with~$\mathbf p=\nabla w$. By minimality,
\begin{equation*}
	-\frac12\|\nabla w_D^{(K)}\|^2_{\underline L^2(\cu_K)}=I_K(w_D^{(K)})\leq I_K(v_{K,L})\,.
\end{equation*}
By stationarity, H\"older's inequality,~\eqref{e.pot.sol.decomposition.corrector.limit}, and~\eqref{e.approximation.stationary.potential.by.dirichlet},
\begin{equation*}
	\lim_{L\to\infty}\limsup_{K\to\infty}\biggl|\E\bigl[I_K(v_{K,L})\bigr]+\frac12\E\bigl[|\mathbf p|^2\bigr]\biggr|=0\,.
\end{equation*}
It follows that
\begin{equation}
	\label{e.lower.bound.corrector.limit}
	\liminf_{K\to\infty}\E\bigl[\|\nabla w_D^{(K)}\|^2_{\underline L^2(\cu_K)}\bigr]\geq\E\bigl[|\mathbf p|^2\bigr]
	\,.
\end{equation}
Let~$\mathbf j_{K,L}$ be given by Lemma~\ref{l.approximation.stationary.by.local}\rm{(ii)}. Since~$\mathbf j_{K,L}$ and~$\nabla w_N^{(K)}+\f$ belong to~$\Lsolo(\cu_K)$, orthogonality to gradients gives
\begin{equation*}
	\|\nabla w_N^{(K)}\|^2_{\underline L^2(\cu_K)}=-\fint_{\cu_K}\f\cdot\mathbf p-\fint_{\cu_K}(\mathbf j-\mathbf j_{K,L})\cdot\nabla w_N^{(K)}\,.
\end{equation*}
Taking expectations and using H\"older's and Young's inequalities yields, for every~$\delta\in(0,1)$,
\begin{equation*}
	\E\bigl[\|\nabla w_N^{(K)}\|^2_{\underline L^2(\cu_K)}\bigr]\leq(1+\delta)\E\bigl[|\mathbf p|^2\bigr]+C\delta^{-1}\E\bigl[\|\mathbf j-\mathbf j_{K,L}\|^2_{\underline L^2(\cu_K)}\bigr]
	\,.
\end{equation*}
By~\eqref{e.approximation.stationary.solenoidal.by.neumann}, we may take first~$K\to\infty$, then~$L\to\infty$, and finally~$\delta\to0$ to obtain
\begin{equation}
	\label{e.upper.bound.neumann.corrector.limit}
	\limsup_{K\to\infty}\E\bigl[\|\nabla w_N^{(K)}\|^2_{\underline L^2(\cu_K)}\bigr]\leq\E\bigl[|\mathbf p|^2\bigr]
	\,.
\end{equation}
Combining~\eqref{e.ordering.corrector.limit},~\eqref{e.lower.bound.corrector.limit}, and~\eqref{e.upper.bound.neumann.corrector.limit} proves~\eqref{e.corrector.limit}.
\end{proof}

\begin{lemma}[{Approximate recurrence formula for~$\shom_m^2$}]
\label{l.approximate.recurrence.formula}
Assume~$m_0 \in (\mstarstar,\infty) \cap \Z$ and~$E\in [1,\infty)$ are such that~$\mathcal{S}(m_0-1,E)$ is valid.
There exists a constant~$C(d)<\infty$ such that, if
\begin{equation}
	\label{e.lower.bound.cgamma.cond.again.0}
	E \geq C \cstar^{-1}
	\qand
	\cgamma \leq E^{-5}\,,
\end{equation}
then, for every~$n \in (-\infty,m_0] \cap \Z$ and~$h \in \N$ with~$h \leq 6 \cstar \cgamma^{-1}$ and~$\mstarstar \leq n \leq m_0-h$,
we have
\begin{equation}
	\left\{
	\begin{aligned}
		&
		\shom_{n+h} \shom_{n}^{-1}
		\leq 1 + (\log 3) \cstar  \shom_{n}^{-2} \sum_{k=n+1}^{n+h} 3^{2\cgamma k}
		+ C E^2  \left|\log \cgamma\right|^{2}  \cgamma \,,
		\\
		& \shom_{n} \shom_{n+h}^{-1}
		\leq
		1 - (\log 3) \cstar \shom_{n}^{-2} \sum_{k=n+1}^{n+h} 3^{2\cgamma k}
		+
		C h^2 \shom_{n}^{-4}  3^{4\cgamma (n+h)}
		+ C E^2  \left|\log \cgamma\right|^{2}  \cgamma \,.
	\end{aligned}
	\right.
	\label{e.what.do.we.have}
\end{equation}
\end{lemma}
\begin{proof}
\emph{Step 1.}
We set up the proof using~$m$ for the upper scale, so that the lower scale denoted by~$n$ in the statement is~$m-h$. Let~$m \in \Z$ with~$m \leq m_0$ and~$h \in \N$ with~$h \leq 6 \cstar \cgamma^{-1}$ be given.
Assume that
\begin{equation}
	\cgamma \leq E^{-5} \qand E \geq M \cstar^{-1}\,,
	\label{e.cgamma.constraints}
\end{equation}
for a constant~$M = M(d) < \infty$ to be determined below. Fix parameters~$K, n \in \Z$ which satisfy
\begin{equation}
	K \geq m +10^{10}  \cgamma^{-1}
	\qand
	n \coloneqq m-h-32\lceil|\log_3\cgamma|\rceil\,.
	\label{e.recurrence.params}
\end{equation}
Note that selecting~$M$ large enough in~\eqref{e.cgamma.constraints} ensures that~$m-n\leq C\cgamma^{-1}$.

\smallskip

Denote~$\h \coloneqq  \k_{m} - \k_{m-h}$, fix~$e,e' \in \Rd$ with~$|e|,|e'| \leq 1$ , and let
\begin{equation}
	P  \coloneqq  \bfAhom_{m-h}^{-\nf 12} \begin{pmatrix} e'  \\ e \end{pmatrix} \in \R^{2d}\,.
	\label{e.recurrence.P.def}
\end{equation}
Let~$w_{D,e}^{(K)} \in H^1_0(\cu_K)$ and~$w_{N,e'}^{(K)} \in H^1(\cu_K)$ be the solutions of
\begin{equation}
	\label{e.def.w}
	\left\{
	\begin{aligned}
		& - \Delta w_{D,e}^{(K)}  =
		\nabla \cdot ( \shom_{m-h}^{-1} \h e)
		& \mbox{in} & \ \cu_{K} \,, \\
		&
		w_{D,e}^{(K)} = 0  & \mbox{on} & \ \partial \cu_{K} \,,
	\end{aligned}
	\right.
	\quad
	\left\{
	\begin{aligned}
		& - \Delta w_{N,e'}^{(K)} =
		\nabla \cdot ( \shom_{m-h}^{-1} \h e')
		& \mbox{in} & \ \cu_{K} \,, \\
		&
		\mathbf{n}_{\cu_K} \cdot (\nabla w_{N,e'}^{(K)}+ \shom_{m-h}^{-1} \h e')  = 0  & \mbox{on} & \ \partial \cu_{K} \,,
	\end{aligned}
	\right.\,,
\end{equation}
respectively.
Testing the equations for~$w_{D,e}^{(K)}$ and~$w_{N,e'}^{(K)}$ with themselves yields the identities
\begin{equation}
	\label{e.what.nablaw.really.is}
	\|\nabla w_{D,e}^{(K)} \|_{\underline{L}^2(\cu_{K})}^2
	=
	-\shom_{m-h}^{-1} \fint_{\cu_K} \h e \cdot \nabla w_{D,e}^{(K)}
	=
	\fint_{\cu_K} \shom_{m-h}^{-1} e  \cdot \h\nabla w_{D,e}^{(K)}
	\,,
\end{equation}
and
\begin{equation}
	\label{e.what.nablaw.really.is.two}
	\|\nabla w_{N,e'}^{(K)} \|_{\underline{L}^2(\cu_{K})}^2
	=
	-\shom_{m-h}^{-1} \fint_{\cu_K} \h e' \cdot \nabla w_{N,e'}^{(K)}
	=
	\fint_{\cu_K} \shom_{m-h}^{-1} e'  \cdot \h\nabla w_{N,e'}^{(K)}
	\,.
\end{equation}
By Calder\'on-Zygmund estimates,~\eqref{e.km.kn.Lp}, and~\eqref{e.shom.h.bounds}, increasing~$M$ in~\eqref{e.cgamma.constraints}, if necessary, we have
\begin{align}
	\label{e.nablaw.in.L.eight}
	\|\nabla w_{D,e}^{(K)} \|_{\underline{L}^8(\cu_{K})}^2
	+
	\|\nabla w_{N,e'}^{(K)} \|_{\underline{L}^8(\cu_{K})}^2
	&\leq C \shom^{-2}_{m-h}  \|\k_m {-} \k_{m-h}\|_{\underline{L}^8(\cu_{K})}^2
	\notag \\ &
	\leq
	C \shom_{m-h}^{-2} h 3^{2\cgamma m}
	{+}
	\O_{\Gamma_{1}}\Bigl(  C \shom_{m-h}^{-2}  3^{2\cgamma m}  h 3^{-\frac{(K-m)d}{8}}\Bigr)
	\notag \\ &
	\leq
	C + \O_{\Gamma_1}(\cgamma^{100})
	\,.
\end{align}
The finite-range assumptions imply ergodicity of the translation action. Moreover,~\eqref{e.diff.law.shift} and~\ref{a.j.frd} show that both forcings are stationary and mean-zero. Using Lemma~\ref{l.corrector.limit}, together with independence and assumption~\ref{a.j.nondeg}, we obtain
\begin{align}
	\label{e.perturb.assumption}
	\lim_{K \to \infty}
	\E \Bigl[
	\|\nabla w_{D,e}^{(K)} \|^2_{\underline{L}^2(\cu_{K})}
	\Bigr]
	&=
	\shom_{m-h}^{-2} \sum_{k = m-h+1}^m \E\Bigl[
	\Bigl| \bigl(\nabla \Delta^{-1} \nabla \cdot (\mathbf{j}_k(\cdot) e)\bigr)(0) \Bigr|^2
	\Bigr]
	\notag \\ & =
	\shom_{m-h}^{-2} \cstar (\log 3) |e|^2
	\sum_{k = m-h+1}^m 3^{2\cgamma k}
	\,.
\end{align}
In the last equality above we used~\ref{a.j.nondeg} and~\eqref{e.diff.law.shift}. Similarly,
\begin{equation}
	\label{e.perturb.assumption.two}
	\lim_{K \to \infty}
	\E \Bigl[
	\|\nabla w_{N,e'}^{(K)} \|^2_{\underline{L}^2(\cu_{K})}
	\Bigr]
	=
	\shom_{m-h}^{-2} \cstar (\log 3) |e'|^2
	\sum_{k = m-h+1}^m 3^{2\cgamma k}
	\,.
\end{equation}

\smallskip

The Calder\'on--Zygmund estimates for the two problems in~\eqref{e.def.w} give
\begin{equation}
	\label{e.nablaw.hessian.L.eight}
	\|\nabla^2w_{D,e}^{(K)}\|_{\underline L^8(\cu_K)}
	+
	\|\nabla^2w_{N,e'}^{(K)}\|_{\underline L^8(\cu_K)}
	\leq
	C\shom_{m-h}^{-1}\|\nabla\h\|_{\underline L^8(\cu_K)}
	\,.
\end{equation}
By stationarity,~\ref{a.j.reg},~\eqref{e.nablaw.hessian.L.eight} and~\eqref{e.W.1.inf.bound},
\begin{multline}
	\label{e.lower.bound.oscillations}
	3^n\E\Bigl[\|\nabla^2w_{D,e}^{(K)}\|_{\underline L^8(\cu_K)}^8+\|\nabla^2w_{N,e'}^{(K)}\|_{\underline L^8(\cu_K)}^8\Bigr]^{\!\nf18}
	+\shom_{m-h}^{-1}3^n\E\Bigl[\|\nabla\h\|_{L^\infty(\cu_{m-h})}^4\Bigr]^{\!\nf14}
	\\
	\leq C\shom_{m-h}^{-1}\cgamma^{-1}3^{-(m-h-n)}3^{\cgamma m}
	\leq\cgamma^{15}
	\,.
\end{multline}
The Poincar\'e inequality on each cube~$z+\cu_n$, followed by summation, gives, for~$w\in\{w_{D,e}^{(K)},w_{N,e'}^{(K)}\}$,
\begin{equation}
	\label{e.lower.bound.gradient.oscillations}
	\biggl(\avsum_{z\in3^n\Zd\cap\cu_K}\E\Bigl[\|\nabla w-(\nabla w)_{z+\cu_n}\|_{\underline L^8(z+\cu_n)}^8\Bigr]\biggr)^{\!\nf18}
	\leq C\shom_{m-h}^{-1}\cgamma^{-1}3^{-(m-h-n)}3^{\cgamma m}
	\leq\cgamma^{15}
	\,.
\end{equation}

Next, for each~$z \in 3^n \Zd \cap \cu_K$, we define
\begin{equation}
	P_z = \begin{pmatrix} p_z \\ q_z \end{pmatrix}
	\coloneqq
	\bfAhom_{m-h}^{-\nf 12} \begin{pmatrix} e' +(  \nabla w_{D,e}^{(K)})_{z+\cu_n}  \\  e + ( \nabla w_{N,e'}^{(K)}+ \shom_{m-h}^{-1} \h e')_{z+\cu_n}
	\end{pmatrix}
	\,,
	\label{e.Pz.def}
\end{equation}
and
\begin{equation}
	\mathbf{F}_z  \coloneqq
	\bfAhom_{m-h}^{-\nf 12}\begin{pmatrix}  \nabla w_{D,e}^{(K)} - (  \nabla w_{D,e}^{(K)})_{z+\cu_n}
		\\
		\nabla w_{N,e'}^{(K)}+ \shom_{m-h}^{-1} \h e' - ( \nabla w_{N,e'}^{(K)}+ \shom_{m-h}^{-1} \h e' )_{z+\cu_n}
	\end{pmatrix}
	\,.
	\label{e.Fz.def}
\end{equation}
Consider the minimizers of the variational problem in~\eqref{e.variational.mu.U.P} in the domain~$U = z+\cu_n$,
\begin{equation*}
	S_z \coloneqq S\left( \cdot,z+\cu_n, -P_z, 0 \, ; \a_{m} \right)
	\,,
\end{equation*}
and
\begin{multline*}
	X_z \coloneqq \argmin\Biggl\{\fint_{z+\cu_n}\frac12 X\cdot\bfA_mX:\ X\in
	\bfAhom_{m-h}^{-\nf12}
	\begin{pmatrix}e'+\nabla w_{D,e}^{(K)}\\ e+\nabla w_{N,e'}^{(K)}+\shom_{m-h}^{-1}\h e'\end{pmatrix}
	\\
	+(L^2_{\pot,0}\times\Lsolo)(z+\cu_n)\Biggr\}\,,
\end{multline*}
and the perturbation
\begin{equation*}
	\tilde{S}_z
	\coloneqq
	\argmin \left\{
	\fint_{z+\cu_n} \tfrac12\,X\!\cdot\!\bfA_m X
	:\; X \in \mathbf{F}_z + (L^2_{\pot,0} \times \Lsolo)(z+\cu_n)
	\right\}\,.
\end{equation*}
Observe that~$X_z = S_{z} + \tilde{S}_z$. Since~$\nabla w_{D,e}^{(K)} \in L^2_{\pot,0}(\cu_K)$ and~$ \nabla w_{N,e'}^{(K)}+ \shom_{m-h}^{-1} \h e'   \in \Lsolo(\cu_K)$, we have that~$ \sum_{z \in 3^n \Zd \cap \cu_K} X_z \indc_{z+\cu_n} \in P + L^2_{\pot,0}(\cu_K) \times \Lsolo(\cu_K)$, and we may insert
it into the minimization problem in~\eqref{e.variational.mu.U.P} for~$\bfA_m(\cu_K)$. We obtain
\begin{align}
	\label{e.lower.bound.basic.split}
	\! \!
	P \cdot \bfA_m(\cu_{K}) P
	&
	\leq
	\avsum_{z \in 3^n \Zd \cap \cu_K}  \fint_{z+\cu_n} X_z \cdot \bfA_m X_z
	\notag \\ &
	=
	\avsum_{z \in 3^n \Zd \cap \cu_K}
	\fint_{z+\cu_n} \Bigl(
	S_z \cdot  \bfA_{m}  S_z
	+
	2 S_{z} \cdot  \bfA_{m}  \tilde S_z
	+
	\tilde S_z  \cdot  \bfA_{m}  \tilde S_z  \Bigr)
	\notag \\ &
	=
	\avsum_{z \in 3^n \Zd \cap \cu_K}  \! \!
	\biggl( P_z \cdot  \bfA_{m}(z+\cu_n)  P_z
	+
	2 P_{z} \cdot \fint_{z + \cu_n}
	\bfA_{m}  \tilde S_z
	+
	\fint_{z + \cu_n}
	\tilde S_z  \cdot  \bfA_{m}  \tilde S_z  \biggr)
	\,.
\end{align}
In the next three steps we will estimate the terms on the right in the above display.
The first term is the principal one; the latter two terms are small because of the scale separation between~$m-h$ and~$n$.

\smallskip

\emph{Step 2.} We show that for a sufficiently large choice of~$M$ in~\eqref{e.cgamma.constraints},
\begin{equation}
	\label{e.lower.bound.localization.terms}
	\E\Biggl[ \avsum_{z \in 3^n \Zd \cap \cu_K}
	\biggl(
	2 P_{z} \cdot \fint_{z + \cu_n}
	\bfA_{m}  \tilde S_z
	+
	\fint_{z + \cu_n}
	\tilde S_z  \cdot  \bfA_{m}  \tilde S_z  \biggr)\Biggr]
	\leq
	\cgamma^{10}
	\,.
\end{equation}
For each~$z \in 3^n \Zd \cap \cu_K$, since~$\mathbf{F}_z$ has mean zero in~$z + \cu_n$,
\begin{equation*}\fint_{z+\cu_n}  \tilde{S}_z \cdot \bfA_{m} \tilde{S}_z
=
\fint_{z+\cu_n}  \mathbf{F}_z \cdot \bfA_{m} \tilde{S}_z
\leq
[ \bfAhom_{m-h}^{\nf12}  \mathbf{F}_z  ]_{\underline{H}^1(z+\cu_n)}
[ \bfAhom_{m-h}^{-\nf12} \bfA_{m} \tilde{S}_z]_{\Hminusul(z+\cu_n)}
\,,
\end{equation*}
and the last term can be estimated,  using the coarse-grained Poincar\'e inequality~\eqref{e.CG.Poincare.doubled.vars} for elements of~$\mathcal{S}(z+\cu_n)$ from Remark~\ref{r.cg.poincare.doubled.variables}, by
\begin{multline*}
	3^{-n} [ \bfAhom_{m-h}^{-\nf12} \bfA_{m} \tilde{S}_z]_{\Hminusul(z+\cu_n)}
	\\
	\leq
	C  \bigl( \shom_{m-h}^{-1} \Lambda_{1,2}(z+\cu_n\, ;\a_m) + \shom_{m-h} \lambda_{1,2}^{-1}(z+\cu_n\, ; \a_m)
	\bigr)^{\nf12}
	\| \bfA_{m}^{\nf12} \tilde{S}_z\|_{L^2(z + \cu_n)}
	\,.
\end{multline*}
For each~$z \in 3^n \Zd \cap \cu_K$, let~$\cu_m^{(z)}$ be its
$z$-dependent scale-$m$ triadic ancestor.  Using
\eqref{e.ellipticities.monotone.ordered},~\eqref{e.bound.one.cube.by.lambdas},
and~$m-n\leq C\cgamma^{-1}$,
\begin{equation*}
	\bigl( \shom_{m-h}^{-1} \Lambda_{1,2}(z+\cu_n\, ;\a_m) + \shom_{m-h} \lambda_{1,2}^{-1}(z+\cu_n\, ; \a_m)
	\bigr)
	\leq
	C  \bigl( \shom_{m-h}^{-1} \Lambda_{\cgamma,2}(\cu_m^{(z)}\, ;\a_m) + \shom_{m-h} \lambda_{\cgamma,2}^{-1}(\cu_m^{(z)}\, ; \a_m)
	\bigr)\,.
\end{equation*}
The previous two displays yield that for each~$z \in 3^n \Zd \cap \cu_K$,
\begin{equation*}\| \bfA_{m}^{\nf12} \tilde{S}_z\|_{\underline{L}^2(z+\cu_n)}^2
\leq
C   \bigl( \shom_{m-h}^{-1} \Lambda_{\cgamma,2}(\cu_m^{(z)}\, ; \a_m) + \shom_{m-h} \lambda_{\cgamma,2}^{-1}(\cu_m^{(z)}\, ; \a_m)
 \bigr)
3^{2n} [ \bfAhom_{m-h}^{\nf12}  \mathbf{F}_z  ]_{\underline{H}^1(z+\cu_n)}^2
\,.
\end{equation*}
Similarly, using also the above display, for each~$z \in 3^n \Zd \cap \cu_K$,
\begin{equation*}
	P_{z} \cdot \fint_{z+\cu_n}
	\bfA_{m}  \tilde S_z
	\leq
	C  | \bfAhom_{m-h}^{\nf 12}  P_{z} |
	\bigl( \shom_{m-h}^{-1} \Lambda_{\cgamma,2}(\cu_m^{(z)}\, ; \a_m) + \shom_{m-h} \lambda_{\cgamma,2}^{-1}(\cu_m^{(z)}\, ; \a_m)
	\bigr)
	3^n [ \bfAhom_{m-h}^{\nf12}  \mathbf{F}_z  ]_{\underline{H}^1(z+\cu_n)}
	\,.
\end{equation*}
By stationarity, the fourth moment of the ellipticity factor above does not depend on~$z$ and equals that of the corresponding expression on~$\cu_m$; the finite-domain correctors are handled by the averaged estimates~\eqref{e.nablaw.in.L.eight} and~\eqref{e.lower.bound.oscillations}. Proposition~\ref{p.cg.ellipticity.bounds} with~$\sigma=\nf12$ and~$s=\cgamma$,~\eqref{e.shom.h.bounds}, and
$h\leq6\cstar\cgamma^{-1}\leq9\cgamma^{-1}$ therefore give
\begin{equation}
	\label{e.lower.bound.ellipticity}
	\E\Bigl[ \bigl( \shom_{m-h}^{-1} \Lambda_{\cgamma,2}(\cu_m\, ; \a_m) + \shom_{m-h} \lambda_{\cgamma,2}^{-1}(\cu_m\, ; \a_m)
	\bigr)^4 \Bigr]^{\nf 14}
	\leq
	C \cstar^{-1}\cgamma^{-1} \shom_{m-1} \shom^{-1}_{m-h}
	\leq
	C \cstar^{-1}\cgamma^{-1}\,,
\end{equation}
and we get, by~\eqref{e.nablaw.in.L.eight},
\begin{equation*}
	\biggl(\avsum_{z \in 3^n \Zd \cap \cu_K}
	\E\Bigl[  | \bfAhom_{m-h}^{\nf 12}  P_{z} |^4 \Bigr] \biggr)^{\! \nf 14}
	\leq C\,,
\end{equation*}
and, by~\eqref{e.lower.bound.oscillations},
\begin{equation*}
	3^n \biggl(\avsum_{z \in 3^n \Zd \cap \cu_K}    \E\Bigl[  [ \bfAhom_{m-h}^{\nf12}  \mathbf{F}_z  ]_{\underline{H}^1(z+\cu_n)}^4 \Bigr] \biggr)^{\! \nf 14}
	\leq
	C\cgamma^{15}
	\,.
\end{equation*}
Combining the above estimates with H\"older's inequality yields~\eqref{e.lower.bound.localization.terms}; here~$E\geq C\cstar^{-1}$ and~$\cgamma\leq E^{-5}$ absorb the extra factor~$\cstar^{-1}$ into the stated remainder.

\smallskip

\emph{Step 3.}  We estimate the first term on the right side of~\eqref{e.lower.bound.basic.split} by showing that for a sufficiently large choice of~$M$ in~\eqref{e.cgamma.constraints}
\begin{align}
	\label{e.lower.bound.principal.one}
	\lefteqn{
	\avsum_{z \in 3^n \Zd \cap \cu_K}  \E\bigl[ P_z \cdot \bfA_m(z+\cu_n) P_{z} \bigr]
	} \quad &
	\notag \\ &
	\leq
	\bigl(1+C E^2  \left|\log \cgamma\right|^{2}  \cgamma \bigr)
	\avsum_{z \in 3^n \Zd \cap \cu_K} \E\bigl[ \G_{-(\h)_{z+\cu_n}} P_z \cdot \bfAhom_{m-h} \G_{-(\h)_{z+\cu_n}} P_{z}  \bigr]
	+ \cgamma^6
	\,.
\end{align}
To derive this estimate, we use the event~$\mathcal{Q}_z \coloneqq \mathcal{Q}(n, m-h, z)$ in~\eqref{e.good.local.events} which allows us to switch from the field~$\a_m$ to~$\a_{m-h}$ in the cube~$z+\cu_{n}$. We obtain by Lemma~\ref{l.J.sensitivity}, applied with~$\delta  \coloneqq  \frac12 3^{-\frac14 (m-h-n)}$,  and~\eqref{e.J.by.means.of.bfA} that
\begin{align*}
	P_z \cdot \bfA_m(z+\cu_n) P_{z}  \indc_{ \mathcal{Q}_z }
	&
	=
	\G_{-(\h)_{z+\cu_n}} P_z \cdot \bfA\bigl(z+\cu_n\,; \a_m - (\h)_{z+\cu_n} \bigr) \G_{-(\h)_{z+\cu_n}} P_{z} \indc_{ \mathcal{Q}_z }
	\notag \\ &
	=
	\Bigl( 2 J\bigl(z+\cu_n, -p_z,q_z - (\h)_{z+\cu_n} p_z \,; \a_m - (\h)_{z+\cu_n} \bigr) - 2 p_z \cdot q_z \Bigr)\indc_{ \mathcal{Q}_z }
	\notag \\ &
	\leq
	2 \bigl(1+ \tfrac32 3^{-\frac14 (m-h-n)} \bigr)
	J\bigl(z+\cu_n, -p_z,q_z - (\h)_{z+\cu_n} p_z \,; \a_{m-h} \bigr)\indc_{ \mathcal{Q}_z }
	\notag \\ & \qquad
	+ C 3^{\frac14 (m-h-n)} \Bigl(\bigl(\shom_{m-h} |p_z|^2 + | p_z \cdot q_z| \bigr)  3^{-\frac12 (m-h-n)}  \Bigr)\indc_{ \mathcal{Q}_z }
	- 2 p_z \cdot q_z\indc_{ \mathcal{Q}_z }
	\notag \\ &
	\leq
	\bigl(1+\tfrac32 3^{-\frac14 (m-h-n)}  \bigr) \G_{-(\h)_{z+\cu_n}} P_z \cdot \bfA_{m-h}(z+\cu_n) \G_{-(\h)_{z+\cu_n}} P_{z}
	\notag \\ &
	\qquad
	+ C 3^{-\frac14 (m-h-n)} \bigl(
	\shom_{m-h} |p_z|^2 +  \shom_{m-h}^{-1 }| q_z|^2
	\bigr)
	\,.
\end{align*}
By~\eqref{e.nablaw.in.L.eight} we have
\begin{equation*}\avsum_{z \in 3^n \Zd \cap \cu_K}
\E \bigl[
 \shom_{m-h} |p_z|^2 +  \shom_{m-h}^{-1 }| q_z|^2
 \bigr] \leq C
 \,.
\end{equation*}
Therefore, by independence of~$\bfA_{m-h}$ and~$\G_{-(\h)_{z+\cu_n}} P_z$ (the latter is a function of~$\k_m - \k_{m-h}$), increasing~$M$ in~\eqref{e.cgamma.constraints} if necessary,
\begin{align}
	\label{e.lower.bound.principal.one.pre}
	\lefteqn{
	\avsum_{z \in 3^n \Zd \cap \cu_K}
	\E
	\bigl[ P_z \cdot \bfA_m(z+\cu_n) P_{z}  \indc_{ \mathcal{Q}_z } \bigr]
	} \qquad &
	\notag \\ &
	\leq
	( 1+ \cgamma^6)
	\avsum_{z \in 3^n \Zd \cap \cu_K} \E\Bigl[ \G_{-(\h)_{z+\cu_n}} P_z \cdot \bfAhom_{m-h}(z+\cu_n) \G_{-(\h)_{z+\cu_n}} P_{z} \Bigr]
	+
	\frac12 \cgamma^{6}  \,.
\end{align}
To estimate the quantity on the bad event, we first note that by Proposition~\ref{p.cg.ellipticity.bounds} with~$s=\cgamma$ and~\eqref{e.ellipticities.monotone.ordered},
\begin{equation*}
	\E \Bigl[ \bigl| \bfAhom_{m-1}^{-\nf 12} \bfA_m(\cu_n) \bfAhom_{m-1}^{-\nf 12} \bigr|^8 \Bigr]^{\nf 18}
	\leq
	C
	\E \Bigl[  \bigl( \shom_{m-1}^{-1} \Lambda_{\cgamma ,1}(\cu_m;\a_m) + \shom_{m-1} \lambda_{\cgamma ,1}^{-1}(\cu_m;\a_m)  \bigr)^8 \Bigr]^{\nf 18}
	\leq C \cgamma^{-1}\,,
\end{equation*}
and
\begin{equation*}\biggl( \avsum_{z \in 3^n \Zd \cap \cu_K} \E\Bigl[ \bigl| \bfAhom_{m-1}^{\nf 12} P_z \bigr|^4 \Bigr]
\biggr)^{\! \nf 14 }
\leq
C
\,.
\end{equation*}
Thus, by H\"older's inequality, using Lemma~\ref{l.bad.event.lemma} to estimate the bad event, recalling that~$\mathcal{Q}_z \coloneqq \mathcal{Q}(n, m-h, z)$ is defined above in~\eqref{e.good.local.events},
\begin{equation*}\avsum_{z \in 3^n \Zd \cap \cu_K} \E \bigl[ P_z \cdot \bfA_m(z+\cu_n) P_{z}  \indc_{\neg  \mathcal{Q}_z } \bigr]
\leq
C \cgamma^{-2} \exp(-\cgamma^{-1})
\leq
\frac12 \cgamma^{6}
 \,.
\end{equation*}
It remains to switch the term~$\bfAhom_{m-h}(\cu_n)$ to~$\bfAhom_{m-h}$; for this, we claim that
\begin{equation}
	\label{e.use.also.for.the.upper.bound}
	\E\Bigl[ \bigl| \bfAhom_{m-h}^{\nf12}(z+\cu_n) \G_{-(\h)_{z+\cu_n}}  P_z \bigr|^2\Bigr]
	\leq
	\bigl(1+C E^2  \left|\log \cgamma\right|^{2}  \cgamma \bigr)
	\E\Bigl[\bigl|\bfAhom_{m-h}^{\nf12} \G_{-(\h)_{z+\cu_n}}P_z \bigr|^2
	\Bigr]
	\,.
\end{equation}
Combining the above two displays with~\eqref{e.lower.bound.principal.one.pre} leads to~\eqref{e.lower.bound.principal.one}.

\smallskip

To prove~\eqref{e.use.also.for.the.upper.bound}, cube symmetry shows that the two diagonal blocks of~$\bfAhom_{m-h}(\cu_n)$ are scalar. Denote the corresponding eigenvalues after normalization by~$\bfAhom_{m-h}$ by~$x$ and~$y$. The monotonicity in~\eqref{e.homs.defs} gives~$x,y\geq1$, and~\eqref{e.bfJ.magic} gives
\begin{equation*}
	(x-1)+(y-1)
	=2\max_{|e|=1}\E\bigl[\bfJ(\cu_n,\bfAhom_{m-h}^{-\nf12}e\,,
	\bfAhom_{m-h}^{\nf12}e\,;\a_{m-h})\bigr]\,.
\end{equation*}
Therefore, by~\eqref{e.bound.one.cube.by.mathcalE},~\eqref{e.mathcalE.monotone.ordered}, and~\eqref{e.new.induction.for.shom}, for~$t\coloneqq\left|\log\cgamma\right|^{-1}$, using that~$m-h-n\leq C\left|\log\cgamma\right|$ from the definition of~$n$ in~\eqref{e.recurrence.params},
\begin{align*}
	\bigl| \bfAhom_{m-h}^{-\nf 12}  \bfAhom_{m-h}(\cu_n) \bfAhom_{m-h}^{-\nf 12} - \Itwod \bigr|
	&=\max\{x-1,y-1\}
	\leq
	2
	\max_{|e|=1}
	\E\bigl[ \bfJ(\cu_n, \bfAhom_{m-h}^{-\nf 12} e ,  \bfAhom_{m-h}^{\nf 12} e\,;
	\a_{m-h})\bigr]  \\
	&\leq
	C 3^{2t(m-h-n)} \E\bigl[ \mathcal{E}_{t, \infty, 2} (\cu_{m-h}; \a_{m-h}, \shom_{m-h})^2 \bigr]
	\leq
	C E^2  \left|\log \cgamma\right|^{2}  \cgamma
	\,.
\end{align*}
This implies~\eqref{e.use.also.for.the.upper.bound}, completing the proof of~\eqref{e.lower.bound.principal.one}.

\smallskip

\emph{Step 4.}
We show that, for every~$|e'| = 1$,
\begin{multline}
	\limsup_{K\to \infty}
	\avsum_{z \in 3^n \Zd \cap \cu_K} \E\Biggl[ \, \Biggl| \bfAhom_{m-h}^{\nf12} \G_{-(\h)_{z+\cu_n}}  \bfAhom_{m-h}^{-\nf 12} \begin{pmatrix} e'  \\  ( \nabla w_{N,e'}^{(K)}+ \shom_{m-h}^{-1} \h e')_{z+\cu_n}
		\end{pmatrix} \Biggr|^2\Biggr]
	\\
	\leq
	1+  (\log 3) \cstar \shom_{m-h}^{-2} \sum_{k = m-h+1}^m 3^{2\cgamma k}
	\,.
	\label{e.lower.bound.pre1}
\end{multline}
Averaging the cube-wise Jensen inequalities over~$z \in 3^{n} \Zd \cap \cu_K$ gives
\begin{multline}
	\avsum_{z \in 3^n \Zd \cap \cu_K}
	\Biggl| \bfAhom_{m-h}^{\nf12} \G_{-(\h)_{z+\cu_n}}  \bfAhom_{m-h}^{-\nf 12} \begin{pmatrix} e'  \\  ( \nabla w_{N,e'}^{(K)}+ \shom_{m-h}^{-1} \h e')_{z+\cu_n}
		\end{pmatrix} \Biggr|^2
	\\
	=
	|e'|^2 +
	\avsum_{z \in 3^n \Zd \cap \cu_K} \bigl| (\nabla w_{N,e'}^{(K)})_{z+\cu_n} \bigr|^2
	\leq
	|e'|^2 +
	\bigl\| \nabla w_{N,e'}^{(K)} \bigr\|_{\underline{L}^2(\cu_K)}^2
	\,.
\end{multline}
The above display, together with~\eqref{e.perturb.assumption.two},  implies~\eqref{e.lower.bound.pre1}.

\smallskip

\emph{Step 5.}
We show that for every~$|e| = 1$,
\begin{multline}
	\label{e.lower.bound.pre2}
	\limsup_{K \to \infty}
	\avsum_{z \in 3^n \Zd \cap \cu_K} \E\Biggl[ \biggl| \bfAhom_{m-h}^{\nf12} \G_{-(\h)_{z+\cu_n}}  \bfAhom_{m-h}^{-\nf 12} \begin{pmatrix} (  \nabla w_{D,e}^{(K)})_{z+\cu_n}  \\  e
		\end{pmatrix} \biggr|^2\Biggr]
	\\
	\leq
	1 - (\log 3) \cstar \shom_{m-h}^{-2}|e|^2 \sum_{k = m-h+1}^m 3^{2\cgamma k}
	+
	C h^2 \shom_{m-h}^{-4}  3^{4\cgamma m}
	+
	C \cgamma^{15}
	\,.
\end{multline}
For each~$z \in 3^{n} \Zd \cap \cu_K$, we have
\begin{align*}
	\lefteqn{
	\biggl| \bfAhom_{m-h}^{\nf12} \G_{-(\h)_{z+\cu_n}}  \bfAhom_{m-h}^{-\nf 12} \begin{pmatrix} (  \nabla w_{D,e}^{(K)})_{z+\cu_n}  \\  e
		\end{pmatrix} \biggr|^2
	} \quad &
	\notag \\ &
	=
	|e|^2 + |(\nabla w_{D,e}^{(K)})_{(z+\cu_n)}|^2
	+ \shom_{m-h}^{-2} |(\h)_{z+\cu_n}(\nabla w_{D,e}^{(K)})_{(z+\cu_n)}|^2
	-2 \shom_{m-h}^{-1} e \cdot (\h)_{z+\cu_n}(\nabla w_{D,e}^{(K)})_{(z+\cu_n)}\,.
\end{align*}
By~\eqref{e.what.nablaw.really.is},~\eqref{e.km.kn.Lp},~\eqref{e.nablaw.in.L.eight},~\eqref{e.lower.bound.oscillations},~\eqref{e.lower.bound.gradient.oscillations}, and H\"older's inequality,
\begin{multline*}
	\avsum_{z \in 3^n \Zd \cap \cu_K}
	\!\!\!\!
	\E\Bigl[ |(\nabla w_{D,e}^{(K)})_{(z+\cu_n)}|^2
	+ \shom_{m-h}^{-2} |(\h)_{z+\cu_n}(\nabla w_{D,e}^{(K)})_{(z+\cu_n)}|^2 -2 \shom_{m-h}^{-1} e \cdot (\h)_{z+\cu_n}(\nabla w_{D,e}^{(K)})_{(z+\cu_n)}
	\Bigr] \\
	\leq
	\shom_{m-h}^{-2} \E[\| \h \nabla w_{D,e}^{(K)} \|^2_{\underline{L}^2(\cu_{K})}] - \E[\| \nabla w_{D,e}^{(K)} \|^2_{\underline{L}^2(\cu_{K})}] + \cgamma^{15}\,.
\end{multline*}
By~\eqref{e.perturb.assumption},~\eqref{e.km.kn.Lp}, and H\"older's inequality
\begin{multline*}
	\limsup_{K \to \infty}
	\Bigl(\shom_{m-h}^{-2} \E\Bigl[ \| \h \nabla w_{D,e}^{(K)} \|^2_{\underline{L}^2(\cu_{K})} \Bigr]
	-
	\E\Bigl[\| \nabla w_{D,e}^{(K)} \|^2_{\underline{L}^2(\cu_{K})} \Bigr] \Bigr)
	\\
	\leq
	-(\log 3) \cstar \shom_{m-h}^{-2} \sum_{k = m-h+1}^m 3^{2\cgamma k}
	+ C \shom_{m-h}^{-4} h^2 3^{4\cgamma m}\,.
\end{multline*}
Combining the above three displays gives~\eqref{e.lower.bound.pre2}.

\smallskip

\emph{Step 6.} The conclusion. By combining~\eqref{e.lower.bound.basic.split},~\eqref{e.lower.bound.localization.terms},~\eqref{e.lower.bound.principal.one} and the definition of~$P_z$ in~\eqref{e.Pz.def} and passing~$K \to \infty$, we obtain
\begin{align*}
	\lefteqn{
	\begin{pmatrix} e'  \\ e \end{pmatrix}
	\cdot
	\bfAhom_{m-h}^{-\nf 12} \bfAhom_m \bfAhom_{m-h}^{-\nf 12}
	\begin{pmatrix} e'  \\ e \end{pmatrix}
	} \qquad &
	\notag \\ &
	=
	P \cdot \bfAhom_m P
	\notag \\ &
	= \lim_{K\to \infty}
	\E \bigl[ P \cdot \bfA_m(\cu_{K}) P \bigr]
	\notag \\ &
	\leq
	C \cgamma^6
	+
	\bigl(1+C E^2  \left|\log \cgamma\right|^{2}  \cgamma \bigr)
	\limsup_{K\to \infty}
	\avsum_{z \in 3^n \Zd \cap \cu_K} \E\bigl[ \G_{-(\h)_{z+\cu_n}} P_z \cdot \bfAhom_{m-h} \G_{-(\h)_{z+\cu_n}} P_{z}  \bigr]
	\notag \\ &
	=
	C \cgamma^6 +
	\bigl(1+C E^2  \left|\log \cgamma\right|^{2}  \cgamma \bigr)
	\notag \\ & \qquad
	\times
	\limsup_{K \to \infty}
	\avsum_{z \in 3^n \Zd \cap \cu_K} \E\Biggl[ \biggl| \bfAhom_{m-h}^{\nf12} \G_{-(\h)_{z+\cu_n}}  \bfAhom_{m-h}^{-\nf 12} \begin{pmatrix} e'+(  \nabla w_{D,e}^{(K)})_{z+\cu_n}  \\  e
		+ ( \nabla w_{N,e'}^{(K)}+ \shom_{m-h}^{-1} \h e')_{z+\cu_n}
		\end{pmatrix} \biggr|^2\Biggr]
	\,.
\end{align*}
We specialize first to~$e=0$ and then to~$e'=0$.
If~$e=0$, then we use~\eqref{e.lower.bound.pre1} to bound the expression on the right side, and this leads to the first estimate of~\eqref{e.what.do.we.have}. Similarly, if~$e'=0$, then we use~\eqref{e.lower.bound.pre2} to obtain the second estimate of~\eqref{e.what.do.we.have}. This completes the proof.
\end{proof}

We integrate the recurrence relation~\eqref{e.what.do.we.have} in the following lemma.

\begin{lemma}[Integration of approximate recurrence]
\label{l.integrate.approx.recurrence}
Let~$m_0\in \Z$ and~$\{ s_n \}_{n\in\Z} \subseteq (0,\infty)$ be a sequence satisfying, for every~$m\in\Z$,
\begin{equation}
	\nu \leq s_m \leq \bigl(1 + \nu^{-2} \cstar^{+} \cgamma^{-1} 3^{2\cgamma m} \bigr)
	\nu
	\label{e.small.scale.bound}
	\,,
\end{equation}
and, for every~$m,n\in\Z$ with
$\mstarstar\leq n\leq m\leq m_0$ and
$m\leq n+\cstar\cgamma^{-1}$,
\begin{align}
	& s_m \leq (1 + E\cgamma) s_n + A_{n,m} s_n^{-1} \,, \quad\text{and}
	\label{e.assump.upper}
	\\ &
	s_m \geq (1-E\cgamma)s_n + A_{n,m}s_n^{-2} s_m
	- F(m-n)^2 s_n^{-4} s_m 3^{4\cgamma m} \,,
	\label{e.assump.lower.modified}
\end{align}
where~$E,F \in [1,\infty)$ are given constants and
\begin{equation}
	\label{e.A.def}
	A_{n,m}\coloneqq \cstar(\log 3)\sum_{k=n+1}^m 3^{2\cgamma k}
	\,.
\end{equation}
Then there exists a universal~$C<\infty$ such that, if~$\cgamma \leq C^{-1} E^{-1} (1+\cstar^{-2} F+\cstar^{-2}\cstar^{+})^{-1}$, then
\begin{equation}
	\bigl| s_m - (\nu^2 + \cstar \cgamma^{-1} 3^{2\cgamma m})^{\nf12} \bigr|
	\leq
	C\bigl( E (1+\cstar^{-2} F+\cstar^{-2}\cstar^{+}) \cgamma \bigr)^{\nf 12}  s_m
	\,, \quad \forall m\in\Z \cap [\mstarstar,m_0]\,.
	\label{e.concl.root.gamma.inv}
\end{equation}
\end{lemma}
\begin{proof}
We use~$\cstar\leq\frac32$,~$\cstar\leq\cstar^{+}$,~$E\cgamma\leq\frac14$, and~$\cgamma\leq\frac12\cstar$, all consequences of the standing assumptions.
Put
\begin{equation*}
	\bar c\coloneqq\min\{\cstar,1\}\,,\qquad
	q\coloneqq1+\cstar^{-2}F+(\cstar\bar c)^{-1}\cstar^{+}\,,
	\qquad
	\delta\coloneqq\Bigl(\frac{E\cgamma}{q}\Bigr)^{\!\nf12}\,.
\end{equation*}
Since~$\cstar\leq\frac32$,
\begin{equation*}
	q\leq4\bigl(1+\cstar^{-2}F+\cstar^{-2}\cstar^{+}\bigr)
	\,.
\end{equation*}
Let~$K_{\cgamma}$ be defined as in~\eqref{e.K.cgamma.def}.
Denote~$T_m\coloneqq \nu^2+\cstar K_{\cgamma}\,3^{2\cgamma m}$. Since~$\cgamma^{-1}  \leq K_{\cgamma} \leq \cgamma^{-1} + C$ by~\eqref{e.K.cgamma.ineq}, we have
\begin{equation*}
	\bigl| T_m  - (\nu^2 + \cstar \cgamma^{-1} 3^{2\cgamma m}) \bigr|
	\leq
	C \cstar 3^{2\cgamma m}
	\leq
	C \cgamma T_m
	\,.
\end{equation*}
Therefore to prove~\eqref{e.concl.root.gamma.inv}, it suffices to demonstrate, for every~$m\in[\mstarstar,m_0]\cap\Z$,
\begin{equation}
	\bigl| s_m^2 - T_m \bigr|
	\leq
	C\delta q T_m\,.
	\label{e.sm.Tm.wts}
\end{equation}
The advantage of working with~$T_m$ is that it satisfies the exact increment identity
\begin{equation}\label{e.T.increment.exact}
T_m-T_n
=2A_{n,m}
\,, \quad \forall m,n\in\Z\,,\ n\leq m\,.
\end{equation}
Therefore the deviation~$d_m \coloneqq s_m^2 - T_m$ satisfies, for every~$m,n\in\Z$ with~$n\leq m$,
\begin{equation}\label{e.d.increment.identity}
d_m-d_n = (s_m^2-s_n^2) - (T_m-T_n) = s_m^2 - s_n^2-2A_{n,m}
\,.
\end{equation}
Let~$m_{\circ}$ be the largest integer such that
\begin{equation}
	3^{2\cgamma m_\circ}
	\leq\delta\nu^2\cgamma(\cstar\bar c)^{-1}
	\,.
	\label{e.mcirc.def}
\end{equation}
After increasing the universal constant in the smallness assumption, we have
\begin{equation*}
	m_\circ-\mstarstar\geq\bar c\cgamma^{-1}\,,\qquad
	\delta\leq\bar c\,,\qquad
	\delta\cgamma^{-1}\geq2\,,\qquad
	\delta\log3\leq1\,,\qquad
	C\delta q\leq\frac12
	\,.
\end{equation*}
Moreover,~$q\geq(\cstar\bar c)^{-1}\cstar^+\geq\bar c^{-1}$, and therefore~\eqref{e.small.scale.bound} and~\eqref{e.mcirc.def} imply, for every~$m\leq m_\circ$,
\begin{equation*}
	\nu^{-2}\cstar^+\cgamma^{-1}3^{2\cgamma m}
	\leq\delta(\cstar\bar c)^{-1}\cstar^+
	\leq\delta q
	\qand
	\cstar K_\cgamma3^{2\cgamma m}
	\leq C\delta\bar c^{-1}\nu^2
	\leq C\delta q\nu^2
	\,.
\end{equation*}
It follows that
\begin{equation}
	|s_m^2-T_m|\leq C\delta qT_m\,,\qquad m\leq m_\circ\,.
	\label{e.small.m.is.done}
\end{equation}
This is the base case for~\eqref{e.sm.Tm.wts}.

\smallskip

\emph{Step 1.} We rearrange~\eqref{e.assump.upper} and~\eqref{e.assump.lower.modified} to get the following  bound for the right side of~\eqref{e.d.increment.identity}:
\begin{equation}
	\bigl| s_m^2 - s_n^2-2A_{n,m} \bigr|
	\leq
	\Bigl( 3E\cgamma + \bigl(s_n^{-2}A_{n,m}\bigr)^2 +2F(m-n)^2 s_n^{-4}3^{4\cgamma m} \Bigr) \bigl(s_n^2+2A_{n,m}\bigr)
	\,.
	\label{e.local.square.increments}
\end{equation}
For simplicity, denote~$x\coloneqq s_m$,~$y\coloneqq s_n$,~$\eps\coloneqq E\cgamma$,~$D \coloneqq F(m-n)^2 s_n^{-4} 3^{4\cgamma m}$ and~$A\coloneqq A_{n,m}$ so that the desired bound~\eqref{e.local.square.increments} can be written as
\begin{equation}\label{e.double.local.square.increment}
|x^2-y^2-2A|
\leq \bigl(3\eps + (A y^{-2})^2  +2D \bigr)(y^2+2A)\,.
\end{equation}
The lower bound~\eqref{e.assump.lower.modified} can be written as
\begin{equation*}
	\bigl( 1- A y^{-2}  +D) x \geq (1-\eps)y\,.
\end{equation*}
Since~$x,y>0$, we deduce that~$1 - Ay^{-2} +D>0$. Squaring both sides of the previous display and applying the elementary inequality
\begin{equation*}
	(1-u)^{-2} \geq 1+2u\,, \quad \forall u \in (-\infty,1)\,,
\end{equation*}
with~$u= Ay^{-2} -D$, we obtain
\begin{align*}
	x^2
	&
	\geq
	(1-\ep)^2 (1 + 2 Ay^{-2} -2D) y^2
	=
	(1-\eps)^2y^2 + 2(1-\eps)^2A -2(1-\eps)^2Dy^2
	\,.
\end{align*}
Rearranging this expression and using that~$1 - (1-\ep)^2 \leq 2\ep$, we obtain
\begin{equation}
	x^2 - y^2 - 2A \geq -2(\eps+D)(y^2+2A)
	\,.
	\label{e.lower.bound.square.incr}
\end{equation}
This is one side of~\eqref{e.local.square.increments}. For the other side, we square the upper bound~\eqref{e.assump.upper} to get
\begin{equation*}
	x^2
	\leq
	(1+\eps)^2y^2
	+ 2(1+\eps)A + A^2 y^{-2}
	\,.
\end{equation*}
Rearranging this, using that~$(1+\eps)^2-1=2\eps+\eps^2\le 3\eps$, we get
\begin{align}
	x^2-y^2-2A
	&
	\leq
	\bigl((1+\eps)^2-1\bigr)y^2 + 2\eps A + A^2y^{-2}
	\leq
	\bigl( 3\eps + (A y^{-2})^2\bigr) (y^2+A) \,.
	\label{e.upper.bound.square.incr}
\end{align}
Combining~\eqref{e.lower.bound.square.incr} and~\eqref{e.upper.bound.square.incr} yields~\eqref{e.double.local.square.increment}.

\smallskip

\emph{Step 2.} We prove the following statement:
for every~$m,n\in\Z$ with~$\mstarstar\leq n\leq m\leq m_0$ satisfying
\begin{equation}
	\frac12 T_n \leq s_n^2 \leq 2T_n
	\,,
	\label{e.boot.sass}
\end{equation}
and
\begin{equation}
	n \leq m \leq n + \cstar \cgamma^{-1}
	\,,
	\label{e.step.sass}
\end{equation}
we have the estimate
\begin{equation}
	\bigl| s_m^2 - s_n^2-2A_{n,m} \bigr|
	\leq
	C \cgamma \Bigl( E + (1+ \cstar^{-2} F) \min \bigl\{ \nu^{-2}  \cstar \cgamma^{-1} 3^{2\cgamma n} \,, 1 \bigr\}^2  \cgamma (m-n)^2
	\Bigr)T_m
	\,.
	\label{e.chunk.increment.under.sass}
\end{equation}
The case~$m=n$ is immediate, so we assume~$n<m$ below.
Using the assumption~\eqref{e.boot.sass} in combination with~\eqref{e.local.square.increments}, we obtain
\begin{align*}
	\bigl| s_m^2 - s_n^2-2A_{n,m} \bigr|
	&
	\leq
	\Bigl(3E\cgamma + \bigl(2T_n^{-1}A_{n,m}\bigr)^2 +8F(m-n)^2 T_n^{-2}3^{4\cgamma m} \Bigr) \bigl(2T_n+2A_{n,m}\bigr)
	\notag \\ &
	\leq
	2 \bigl(3E\cgamma + \bigl(2T_n^{-1}A_{n,m}\bigr)^2 +8F(m-n)^2 T_n^{-2}3^{4\cgamma m} \bigr) T_m
	\,.
\end{align*}
For every~$m,n\in\Z$ satisfying~\eqref{e.step.sass}, we have
\begin{equation*}\frac{2 A_{n,m}}{T_n(m-n)}
\leq
\frac{2  \cstar (\log 3) 3^{2 \cgamma m}}{\nu^2 + \cstar \cgamma^{-1} 3^{2\cgamma n}}
=
\frac{2 (\log 3) 3^{2\cgamma(m-n)} }{\nu^{2}\cstar^{-1}  3^{-2\cgamma n} +\cgamma^{-1} }
\leq
C\min\bigl\{ \nu^{-2}  \cstar 3^{2\cgamma n} \,, \cgamma \bigr\}
\,,
\end{equation*}
and
\begin{equation*}T_n^{-1}3^{2\cgamma m}
\leq 3^{2 \cgamma(m-n)} \frac{3^{2\cgamma n}}{\nu^{2} + \cstar \cgamma^{-1} 3^{2\cgamma n}}
\leq \frac{C \cstar^{-1}}{\nu^{2}\cstar^{-1}  3^{-2\cgamma n} +\cgamma^{-1} }
\leq
C\cstar^{-1} \min\bigl\{ \nu^{-2}  \cstar 3^{2\cgamma n} \,, \cgamma \bigr\}
 \,.
\end{equation*}
Combining these yields~\eqref{e.chunk.increment.under.sass}, completing the proof of the claim.

\smallskip

\emph{Step 3.} We complete the proof by strong induction. For~$k>m_\circ$, set
\begin{equation*}
	a_k\coloneqq
	\max\bigl\{\delta\nu^2\cstar^{-1}3^{-2\cgamma k},\delta\cgamma^{-1}\bigr\}
	\qand
	h_k\coloneqq\lfloor a_k\rfloor
	\,.
\end{equation*}
The definition of~$m_\circ$ and the preceding smallness bounds give
\begin{equation}
	2\leq a_k\,,\qquad
	1\leq h_k\leq a_k\leq\bar c\cgamma^{-1}
	\leq\min\{\cstar\cgamma^{-1},\cgamma^{-1}\}
	\,.
	\label{e.adaptive.step.cap}
\end{equation}
In particular,~$k-h_k\geq\mstarstar$, so the recurrence can be applied between~$k-h_k$ and~$k$. If
\begin{equation*}
	\mu_j\coloneqq\min\bigl\{\nu^{-2}\cstar\cgamma^{-1}3^{2\cgamma j},1\bigr\}
	\,,
\end{equation*}
then monotonicity and~$h_k\leq a_k$ give
\begin{equation*}
	\mu_{k-h_k}h_k\leq\mu_k a_k=\delta\cgamma^{-1}
	\,.
\end{equation*}
Assuming the induction hypothesis at all lower indices, the bound~$C\delta q\leq\frac12$ gives~\eqref{e.boot.sass} at~$k-h_k$. Hence Step~2 and~$q\geq1+\cstar^{-2}F$ yield
\begin{equation}
	|d_k-d_{k-h_k}|
	\leq C\cgamma\Bigl(E+(1+\cstar^{-2}F)\mu_{k-h_k}^2\cgamma h_k^2\Bigr)T_k
	\leq CE\cgamma T_k
	\,.
	\label{e.capped.step.increment}
\end{equation}
Since~$a_k\geq2$, we also have~$h_k\geq\frac12a_k$. Using the two factorizations of~$T_k-T_{k-h_k}$, the bounds~$3^r-1\geq r\log3$ for~$r\geq0$ and~$1-3^{-\delta}\geq\frac12\delta\log3$, and~\eqref{e.adaptive.step.cap}, we obtain
\begin{equation}
	T_k-T_{k-h_k}
	\geq\frac{\log3}{18}\delta T_k
	\,.
	\label{e.capped.step.drop}
\end{equation}
Indeed, the first term in the maximum defining~$a_k$ controls the~$\nu^2$ part of~$T_k$, while the second controls~$\cstar K_\cgamma3^{2\cgamma k}$. Since~$E\cgamma=\delta^2q$,~\eqref{e.capped.step.increment} and~\eqref{e.capped.step.drop} imply
\begin{equation*}
	|d_k-d_{k-h_k}|
	\leq C\delta q\bigl(T_k-T_{k-h_k}\bigr)
	\,.
\end{equation*}
Combining this with the induction hypothesis at~$k-h_k$ gives
\begin{equation*}
	|d_k|
	\leq C\delta qT_{k-h_k}
	+C\delta q\bigl(T_k-T_{k-h_k}\bigr)
	=C\delta qT_k
	\,.
\end{equation*}
Together with the base estimate~\eqref{e.small.m.is.done}, this proves~\eqref{e.sm.Tm.wts} and completes the proof.
\end{proof}

We next finish the proof of Proposition~\ref{p.propagate.diffusivity.lower.bound}.

\begin{proof}[{Proof of Proposition~\ref{p.propagate.diffusivity.lower.bound}}]

Apply Lemma~\ref{l.integrate.approx.recurrence} from~$\mstarstar$, using the plateau below it and Lemma~\ref{l.approximate.recurrence.formula} above it, with~$E$ replaced by~$CE^2|\log\cgamma|^2$ and~$F$ by~$C$, using~$\cstar \leq \cstar^{+} \leq C(d)$ to absorb the term~$\cstar^{-2}\cstar^{+}$ into the constant. The forward bound controls~$\shom_m/\shom_n$, so the inverse-ratio bound yields~\eqref{e.assump.lower.modified} after multiplication by~$\shom_m$. Thus~\eqref{e.shom.m.flow} holds at every scale.
In view of~\eqref{e.lower.bound.cgamma.cond} and~\eqref{e.shom.m.flow}, upon taking~$E$ sufficiently large, we have that~\eqref{e.shom.h.bounds} is valid for every~$m\leq m_0$ and in particular for~$m=m_0$.
\end{proof}

We conclude this subsection with a continuity estimate for the map~$m\mapsto \shom_m$, which is needed when we switch from~$\shom_m$ to~$\shom_{m+h}$ in our induction argument.

\begin{lemma}[{Continuity of~$m\mapsto \shom_m$}]
\label{l.shom.continuity}
There exists~$C(d) < \infty$ such that the following holds.
Assume~$m_0 \in (\mstarstar,\infty) \cap \Z$ and~$E\in [C\cstar^{-1},\infty)$ are such that~$\mathcal{S}(m_0-1,E)$ is valid and~$\cgamma \leq E^{-10}$.
Then, for every~$m,n \in (-\infty,m_0] \cap \Z$ with~$m > n$,
\begin{equation}
	\label{e.shom.m.vs.shom.n}
	\shom_m
	\geq
	\frac14\shom_n
	\qand
	\bigl|
	\shom_{m} \shom_{n}^{-1} - 1
	\bigr|
	+
	\bigl|
	\shom_{n} \shom_{m}^{-1} - 1
	\bigr|
	\leq
	C \min\bigl\{ 1 \,, \cgamma (m-n) + E^2  \cgamma \left|\log \cgamma\right|^{2} \bigr\}  3^{\cgamma (m-n)}
	\,.
\end{equation}
\end{lemma}
\begin{proof}
We first claim that, for~$h \in \N$ with~$h \leq \frac12 \min\{1,\cstar^{2}\} \cgamma^{-1}$ and~$n + h \leq m_0$, we have
\begin{equation}
	\label{e.shom.continuity.pre}
	\bigl| \shom_{n+h} \shom_{n}^{-1} - 1 \bigr|
	+
	\bigl| \shom_{n} \shom_{n+h}^{-1} - 1 \bigr|
	\leq
	C \cgamma \bigl( h + E^2 \left|\log \cgamma\right|^{2} \bigr)\,.
\end{equation}
If~$n+h\leq\mstarstar$, the plateau gives~\eqref{e.shom.continuity.pre}. If~$n\geq\mstarstar$,~\eqref{e.what.do.we.have} bounds both
$|\shom_{n+h}\shom_n^{-1}-1|$ and
$|\shom_n\shom_{n+h}^{-1}-1|$ by
\begin{equation*}
	(\log 3) \cstar \shom_{n}^{-2} \sum_{k=n+1}^{n+h} 3^{2\cgamma k}
	+ C h^2 \shom_{n}^{-4} 3^{4\cgamma (n+h)}
	+ C E^2 \cgamma \left|\log \cgamma\right|^{2}\,.
\end{equation*}
For the first term, using~\eqref{e.shom.h.bounds}, we get
\begin{equation*}
	\cstar \shom_{n}^{-2} \sum_{k=n+1}^{n+h} 3^{2\cgamma k}
	\leq
	\cstar  \bigl( 4 \cstar^{-1} \cgamma 3^{-2\cgamma n}\bigr)  h  3^{2\cgamma (n+h)}
	\leq 4 \cgamma h 3^{2 \cgamma h} \leq 12 \cgamma h\,,
\end{equation*}
since~$3^{2\cgamma h} \leq 3$ for~$h \leq \frac12 \cgamma^{-1}$.
For the second term, using again~\eqref{e.shom.h.bounds}, we have
\begin{equation*}
	h^2 \shom_{n}^{-4} 3^{4\cgamma (n+h)}
	= h^2 \bigl( \shom_n^{-2} 3^{2\cgamma n} \bigr)^2 3^{4\cgamma h}
	\leq C \cstar^{-2} \cgamma^2 h^2  3^{4\cgamma h}
	\leq
	C \cstar^{-2} \cgamma^2 h^2
	\leq
	C\cgamma h
	\,.
\end{equation*}
This proves~\eqref{e.shom.continuity.pre} above the landmark; a mixed interval is split at~$\mstarstar$. If~$h \coloneqq m-n \leq \frac12 \min\{1,\cstar^{2}\}  \cgamma^{-1}$ and~$n<m \leq m_0$, then~\eqref{e.shom.continuity.pre} gives
\begin{equation*}
	\bigl| \shom_{m} \shom_{n}^{-1} - 1 \bigr| + \bigl| \shom_{n} \shom_{m}^{-1} - 1 \bigr|
	\leq C \cgamma ( h + E^2 |\log\cgamma|^2 )\,,
\end{equation*}
and since~$3^{\cgamma h} \leq 3$, this implies~\eqref{e.shom.m.vs.shom.n} for~$0< m-n \leq \frac12 \min\{1,\cstar^{2}\}  \cgamma^{-1}$.

\smallskip

For~$m-n \geq  \frac12 \min\{1,\cstar^{2}\}  \cgamma^{-1}$ and~$m \leq m_0$, we use~\eqref{e.lower.bound.cgamma.cond} and obtain, for~$k \in \{m,n\}$,
\begin{equation*}
	\bigl| \shom_k^2 - \bigl( \nu^2 + \cstar \cgamma^{-1} 3^{2\cgamma k} \bigr) \bigr|
	\leq
	C\cstar^{-1} E\cgamma^{\nf12}\left|\log \cgamma\right|
	\shom_k^2
	\,.
\end{equation*}
Thus, by the triangle inequality,
\begin{align*}
	\bigl| \shom_m^2 - \shom_n^2
	\bigr|
	&
	\leq
	\cstar \cgamma^{-1} 3^{2\cgamma n} \bigl| 3^{2\cgamma (m-n)} - 1 \bigr|
	+
	C\cstar^{-1} E\cgamma^{\nf12}\left|\log \cgamma\right|
	\bigl( \shom_m^2 + \shom_n^2 \bigr)
	\,.
\end{align*}
Since~$\cgamma\leq E^{-10}$ implies~$E^4\cgamma^{1/2}|\log\cgamma|\leq16$, put~$\eta=C\cstar^{-1}E\cgamma^{1/2}|\log\cgamma|\leq\frac14\min\{1,\cstar^2\}$. If~$\cgamma(m-n)\leq\eta$, apply~\eqref{e.shom.continuity.pre}; otherwise absorb the defect into~$C\cgamma(m-n)$. Thus, for
$\frac12\min\{1,\cstar^2\}\cgamma^{-1}\leq m-n\leq\cgamma^{-1}$,
\begin{equation*}
	\bigl| \shom_m \shom_n^{-1} - 1
	\bigr|
	+
	\bigl| \shom_n \shom_m^{-1} - 1
	\bigr|
	\leq C \cgamma (m-n)\,.
\end{equation*}
For~$m-n \geq \cgamma^{-1}$,
\begin{equation*}
	\bigl| \shom_m \shom_n^{-1} - 1
	\bigr|
	+
	\bigl| \shom_n \shom_m^{-1} - 1
	\bigr|
	\leq
	C 3^{\cgamma (m-n)}
	\,.
\end{equation*}
Combining the above two cases yields~\eqref{e.shom.m.vs.shom.n} for~$ m - n \geq \frac12 \min\{1,\cstar^{2}\}  \cgamma^{-1}$ and~$m \leq m_0$.

\smallskip

Finally, from~\eqref{e.shom.h.bounds}, for any~$m, n \leq m_0$ with~$m > n$, we have
\begin{equation*}
	\shom_m^2
	\geq
	\tfrac14 \max\{ \cstar \cgamma^{-1} 3^{2\cgamma m}, \nu^2 \}
	\geq
	\tfrac14 \max\{ \cstar \cgamma^{-1} 3^{2\cgamma n}, \nu^2 \}
	\geq
	\tfrac{1}{16} \shom_n^2\,,
\end{equation*}
where the second inequality uses~$m > n$ and the third uses the upper bound in~\eqref{e.shom.h.bounds}. Taking square roots gives~$\shom_m \geq \frac14 \shom_n$.
\end{proof}

\subsection{Contraction of the coarse-graining defect}
\label{ss.homogenization.step}

In this subsection we establish the following bound on the homogenization error, at scales above the infrared cutoff.

\begin{proposition}[Contraction of the coarse-graining defect]
\label{p.homogenization.step}
Assume that~$m\in\Z$ and~$E\geq1$ are such that
\eqref{e.new.induction.for.shom} is valid for every~$k\leq m-1$ and every
$s\in[8\cgamma,1]$. There exists~$C(d)\in[1,\infty)$ such that for every
$\ep\in(0,\nf12]$,
under the assumptions
\begin{equation}
	\label{e.gamma.condition.homog}
	15\cstar^{-1}\leq E\,,
	\qquad \cgamma\leq E^{-10}\,,
	\qquad \cgamma \leq C^{-1} E^{-2} \ep \,\,,
\end{equation}
we have
\begin{equation}
	\label{e.homogenization.step.bound}
	\sup_{L \in \Z \cap (-\infty, m- C \left| \log \ep \right|]}
	\,
	\max_{|e|=1}
	J\bigl(\cu_m,\shom_{L}^{-\nf12} e, \shom_{L}^{\nf12} e \,; \a_L\bigr)
	\leq
	\O_{\Gamma_{1}}\bigl(\ep E^2 \cgamma\bigr)
	+ \O_{\Gamma_{\nf{1}{4} }}\bigl(\ep \exp\bigl(-2E^{-3} \cgamma^{-1}\bigr) \bigr)
	\,.
\end{equation}
\end{proposition}

We follow the renormalization argument given in~\cite[Chapter 5]{AKMBook}, with several ``coarse-graining'' modifications and improvements introduced in~\cite{AK.HC,ABK.SD}.

\smallskip

We begin by giving an upper bound for~$J$ in terms of the weak norm of its maximizer (in a slightly larger cube) and the subadditivity defect. This is a deterministic estimate, that is, it is valid for an arbitrary elliptic coefficient field~$\a(\cdot)$.

\begin{lemma}[{Estimate of~$J$ by weak norms}]
\label{l.J.bound.by.Besov}
There exists~$C(d)<\infty$ such that,
for every~$m\in\Z$, every uniformly elliptic coefficient field~$\a:\cu_m \to \R^{d\times d}_+$,~$p,q\in\Rd$, and~$s \in (0,\nf14]$,
\begin{align}
	\label{e.J.L.bound}
	J(\cu_{m-1},p,q\,;\a)
	&
	\leq
	C {\Lambda}_{s,1}(\cu_m;\a)
	\biggl(\frac{{\Lambda}_{s,1}(\cu_m\,;\a)}{{\lambda}_{s,1}(\cu_m\,;\a)} \biggr)^{\frac{s}{1-2s}}
	3^{-2m} \| \nabla v(\cdot,\cu_m,p,q\,;\a) \|_{\underline H^{-1}(\cu_m)}^2
	\notag \\ &  \qquad \quad
	+
	\sum_{z\in 3^{m-1}\Zd\cap \cu_m}
	2( J(z+\cu_{m-1},p,q\,;\a) - J (\cu_m,p,q\,;\a) )
	\,.
\end{align}

\end{lemma}
\begin{proof}
For convenience, in this proof we simplify the notation by writing~$v \coloneqq  v(\cdot,\cu_m,p,q\,;\a)$ and~$v_{n,z} \coloneqq  v(\cdot,z+\cu_n,p,q\,;\a)$ for each~$n\in\Z$ and~$z\in 3^n\Zd$.
We compute, using the triangle inequality,~\eqref{e.Jenergyv.nosymm} and~\eqref{e.quadresp.nosymm},
\begin{align}
	\label{e.subadd.dig.room}
	J(\cu_{m-1},p,q\,;\a)
	&
	=
	\fint_{\cu_{m-1}}
	\frac12
	\nabla v_{m-1,0}
	\cdot \s\nabla v_{m-1,0}
	\notag \\ &
	\leq
	\fint_{\cu_{m-1}}
	\nabla v
	\cdot \s\nabla v
	+
	\fint_{\cu_{m-1}}
	\left( \nabla v
	-
	\nabla v_{m-1,0} \right)
	\cdot
	\s\left( \nabla v
	-
	\nabla v_{m-1,0}
	\right)
	\notag \\ &
	\leq
	\fint_{\cu_{m-1}}
	\nabla v
	\cdot \s\nabla v
	+
	\sum_{z\in 3^{m-1}\Zd\cap \cu_m}
	\fint_{z+\cu_{m-1}}
	\left( \nabla v
	-
	\nabla v_{m-1,z} \right)
	\cdot
	\s\left( \nabla v
	-
	\nabla v_{m-1,z}
	\right)
	\notag \\ &
	=
	\fint_{\cu_{m-1}}
	\nabla v
	\cdot \s\nabla v
	+
	\sum_{z\in 3^{m-1}\Zd\cap \cu_m}
	\!\!
	2( J(z+\cu_{m-1},p,q\,;\a) - J (\cu_m,p,q\,;\a) )
	\,.
\end{align}
Applying the localized coarse-grained Caccioppoli inequality (Proposition~\ref{p.coarse.grained.Caccioppoli.boundary}) on a finite cover of~$\cu_{m-1}$ by cubes~$z+\cu_{m-2}$ and summing the resulting bounds yields
\begin{equation*}
	\fint_{\cu_{m-1}}
	\nabla v
	\cdot \s \nabla v
	\leq
	C {\Lambda}_{s,1}(\cu_m;\a)
	\biggl(\frac{{\Lambda}_{s,1}(\cu_m\,;\a)}{{\lambda}_{s,1}(\cu_m\,;\a)} \biggr)^{\frac{s}{1-2s}}
	3^{-2m}  \fint_{\cu_m}  | v-(v)_{\cu_m}|^2
	\,.
\end{equation*}
Using~$\| v-(v)_{\cu_m} \|_{\underline{L}^2(\cu_m)}
\leq
C \| \nabla v \|_{\underline H^{-1}(\cu_m)}$ and combining the above displays, we obtain~\eqref{e.J.L.bound}.
\end{proof}

We continue with an estimate of the weak norm of the gradient of the maximizer of~$J$. Like the previous lemma, it is a deterministic estimate and so we write it for general coefficient fields.

\begin{lemma}[Negative norms estimates]
\label{l.Besov.norms}
There exists~$C(d)<\infty$ such that, for every~$m \in \Z$, every uniformly elliptic coefficient field~$\a:\cu_m \to \R^{d\times d}_+$,~$s\in (0,\nf12]$,~$L_0 \in \Z$ with~$L_0 \leq m$,~$p,q\in\Rd$, and~$r \in [1,\infty]$,
\begin{align}
	\label{e.Besov.norms.gradient}
	\lefteqn{
	3^{-m} \bigl\| \nabla v(\cdot,\cu_m,p,q\,;\a) -
	\bigl( \s_*^{-1} (\cu_m)
	( q + \k(\cu_m)p  ) - p \bigr)
	\bigr\|_{\Hminusul(\cu_m)}
	}
	\quad &
	\notag \\ &
	\leq
	C \lambda_{s,r}^{-\nf12}(\cu_m;\a)
	\sum_{n=L_0}^m
	3^{-(1-s)(m-n) }
	\biggl( \avsum_{z\in3^n\Zd \cap \cu_m}
	\!\!
	J(z+\cu_n,p,q\,;\a) - J(\cu_m,p,q\,;\a) \biggr)^{\!\nf12}
	\notag \\ & \quad
	+
	C\sum_{n=L_0}^m
	3^{-(m-n)}
	\biggl(\avsum_{z\in 3^n\Zd \cap \cu_m}
	\bigl|
	\s_*^{-1} (\cu_m) ( q + \k(\cu_m)p  )
	- \s_*^{-1} (z+\cu_n) ( q + \k(z+\cu_n)p  )
	\bigr|^2
	\biggr)^{\!\nf12}
	\notag \\ & \quad
	+C\lambda_{s,r}^{-\nf12}(\cu_m;\a)  3^{-(1-s) (m-L_0)}
	\| \s^{\nf12} \nabla v \|_{\underline{L}^2(\cu_m)}
	\notag \\ & \quad
	+
	C 3^{-(m-L_0)}
	\bigl| \s_*^{-1} (\cu_m)( q + \k(\cu_m)p  ) - p \bigr|
	\,.
\end{align}
\end{lemma}
\begin{proof}
Denote~$v \coloneqq  v(\cdot,\cu_m,p,q\,;\a)$ and~$v_{n,z} \coloneqq  v(\cdot,z+\cu_n,p,q\,;\a)$ for short.
Using~\eqref{e.v.spatial.averages}, we may apply the multiscale Poincar\'e inequality (see~\cite[Proposition~1.10]{AK.Book}) to get
\begin{align}
	\label{e.msp.plus.split.begin}
	\lefteqn{
	3^{-m} \bigl\| \nabla v
	-
	\bigl( \s_*^{-1} (\cu_m)( q + \k(\cu_m)p  ) - p \bigr)
	\bigr\|_{\Hminusul(\cu_m)}
	} \quad &
	\notag \\ &
	\leq
	C\sum_{n=-\infty}^m
	3^{-(m-n)}
	\biggl(\avsum_{z\in 3^n\Zd \cap \cu_m}
	\bigl| (\nabla v)_{z+\cu_n} - \bigl( \s_*^{-1} (\cu_m)( q + \k(\cu_m)p  ) - p \bigr)
	\bigr|^2
	\biggr)^{\!\nf12} \,\,.
\end{align}
We split the sum on the right side in two. By~\eqref{e.energymaps.nonsymm} and~\eqref{e.bound.one.cube.by.lambdas}, for every~$n < L_0$ and~$r \in [1,\infty]$,
\begin{align*}
	\avsum_{z\in 3^n\Zd \cap \cu_m}
	\bigl| (\nabla v)_{z+\cu_n} \bigr|^2
	&
	\leq
	\max_{z\in 3^n\Zd \cap \cu_m}
	\bigl| \s_{*}^{-1}(z+\cu_n) \bigr|
	\| \s^{\nf12} \nabla v \|_{\underline{L}^2(\cu_m)}^2
	\notag \\ &
	\leq
	3^{2s(m-n)} {\lambda}_{s,r}^{-1}(\cu_m; \a)
	\| \s^{\nf12} \nabla v \|_{\underline{L}^2(\cu_m)}^2
	\,.
\end{align*}
Summing this over~$n\in\Z$ with~$n<L_0$ and using the triangle inequality, we get
\begin{align}
	\label{e.msp.bound.for.scalesep}
	\lefteqn{
	\sum_{n=-\infty}^{L_0-1}
	3^{-(m-n)}
	\biggl(\avsum_{z\in 3^n\Zd \cap \cu_m}
	\bigl| (\nabla v)_{z+\cu_n} - \bigl( \s_*^{-1} (\cu_m)( q + \k(\cu_m)p  ) - p \bigr)
	\bigr|^2
	\biggr)^{\!\nf12}
	} &
	\notag \\ &
	\leq
	C3^{-(1 - s) (m-L_0)}  {\lambda}_{s,r}^{-\nf12}(\cu_m;\a)
	\| \s^{\nf12} \nabla v \|_{\underline{L}^2(\cu_m)} +
	C 3^{-(m-L_0)}
	\bigl| \s_*^{-1} (\cu_m)( q + \k(\cu_m)p  ) - p \bigr|
	\,.
\end{align}
We split the sum over the larger scales~$n\geq L_0$ as
\begin{align}
	\label{e.msp.plus.split}
	\lefteqn{
	\sum_{n=L_0}^m
	3^{-(m-n)}
	\biggl(\avsum_{z\in 3^n\Zd \cap \cu_m}
	\bigl| (\nabla v)_{z+\cu_n} - \bigl( \s_*^{-1} (\cu_m)( q + \k(\cu_m)p  ) - p \bigr)
	\bigr|^2
	\biggr)^{\!\nf12}
	} \quad &
	\notag \\ &
	\leq
	C\sum_{n=L_0}^m
	3^{-(m-n)}
	\biggl(\avsum_{z\in 3^n\Zd \cap \cu_m}
	\bigl| (\nabla v)_{z+\cu_n} - \bigl( \s_*^{-1} (z+\cu_n)( q + \k(z+\cu_n)p  ) - p \bigr)
	\bigr|^2
	\biggr)^{\!\nf12}
	\notag \\ & \quad
	+
	C\sum_{n=L_0}^m
	3^{-(m-n)}
	\biggl(\avsum_{z\in 3^n\Zd \cap \cu_m}
	\bigl|
	\s_*^{-1} (\cu_m)( q + \k(\cu_m)p )
	-
	\s_*^{-1} (z+\cu_n)( q + \k(z+\cu_n)p  )
	\bigr|^2
	\biggr)^{\!\nf12}
	\,.
\end{align}
Using~\eqref{e.v.spatial.averages} and~\eqref{e.energymaps.nonsymm}, we have that, for each~$n\leq m$ and~$z\in3^n\Zd\cap\cu_m$,
\begin{align*}
	\lefteqn{
	\Bigl|
	\s_*^{\nf12} (z+\cu_n)
	\bigl(
	(\nabla v)_{z+\cu_n} - \bigl( \s_*^{-1} (z+\cu_n)( q + \k(z+\cu_n)p  ) - p \bigr)
	\bigr)
	\Bigr|^2
	} \qquad &
	\notag \\ &
	=
	\Bigl|
	\s_*^{\nf12} (z+\cu_n)
	\bigl( \nabla v - \nabla v_{n,z} \bigr)_{z+\cu_n}
	\Bigr|^{2}
	\leq
	\fint_{z+\cu_n}
	( \nabla v - \nabla v_{n,z} )\cdot \s ( \nabla v - \nabla v_{n,z} )\,.
\end{align*}
Summing this over~$z\in 3^n\Zd \cap \cu_m$ and applying~\eqref{e.quadresp.nosymm} yields, for each~$n\leq m$,
\begin{align*}
	\lefteqn{
	\avsum_{z\in 3^n\Zd \cap \cu_m}
	\Bigl|
	\s_*^{\nf12} (z+\cu_n)
	\bigl(
	(\nabla v)_{z+\cu_n} - \bigl( \s_*^{-1} (z+\cu_n)( q + \k(z+\cu_n)p  ) - p \bigr)
	\bigr)
	\Bigr|^2
	}   &
	\notag \\ &
	\leq
	\avsum_{z\in 3^n\Zd \cap \cu_m}
	\fint_{z+\cu_n}
	( \nabla v - \nabla v_{n,z} )\cdot \s ( \nabla v - \nabla v_{n,z} )
	=
	2 \!\!
	\avsum_{z\in3^n\Zd \cap \cu_m}
	\!\!
	\bigl( J(z{+}\cu_n,p,q\,;\a) - J(\cu_m,p,q\,;\a) \bigr)
	\,.
\end{align*}
By~\eqref{e.ellipticities.monotone.ordered} and~\eqref{e.bound.one.cube.by.lambdas}, for every~$P  \in \Rd$,~$n\leq m$ and~$z\in 3^n\Zd \cap \cu_m$,
\begin{align*}
	\bigl| P \bigr|^2
	&
	\leq
	\bigl| \s_*^{-1} (z+\cu_n)
	\bigr|
	\bigl|
	\s_*^{\nf12} (z+\cu_n)
	P \bigr|^{2}
	\leq
	3^{2s(m-n)} {\lambda}_{s,r}^{-1}(\cu_m;\a)
	\bigl|
	\s_*^{\nf12} (z+\cu_n)
	P \bigr|^{2}\,.
\end{align*}
Combining the previous three displays, we obtain
\begin{align*}
	\lefteqn{
	\sum_{n=L_0}^m
	3^{-(m-n)}
	\biggl(\avsum_{z\in 3^n\Zd \cap \cu_m}
	\bigl| (\nabla v)_{z+\cu_n} - \bigl( \s_*^{-1} (z+\cu_n)( q + \k(z+\cu_n)p  ) - p \bigr)
	\bigr|^2
	\biggr)^{\!\nf12}
	} \qquad &
	\notag \\ &
	\leq
	C \lambda_{s,r}^{-\nf12}(\cu_m;\a)
	\sum_{n=L_0}^m
	3^{-(1-s)(m-n) }
	\biggl( \avsum_{z\in3^n\Zd \cap \cu_m}
	\!\!
	J(z+\cu_n,p,q\,;\a) - J(\cu_m,p,q\,;\a) \biggr)^{\!\nf12}
	\,.
\end{align*}
Combining the above display,~\eqref{e.msp.plus.split},~\eqref{e.msp.bound.for.scalesep}, and~\eqref{e.msp.plus.split.begin},  yields~\eqref{e.Besov.norms.gradient}.
\end{proof}

\begin{lemma}[Variance estimate]
\label{l.var.bounds}
For every~$m,n,k\in\Z$ with~$k\leq n$ and matrix~$\bfB \in \R^{2 d\times 2 d}_{\sym,+}$,
\begin{align}
	\label{e.var.a.star}
	\var\big[ \bfB^{\nf12}\bfA_{m,*}^{-1}(\cu_n) \bfB^{\nf12} \bigr]
	&
	\leq
	2\var\Biggl[
	\avsum_{z\in 3^k\Zd \cap \cu_n} \!\!\!\!
	\bfB^{\nf12}\bfA_{m,*}^{-1} (z+\cu_k) \bfB^{\nf12}\Biggr]
	\notag \\ & \qquad
	+
	8 \E \Biggl[
	\biggl( \avsum_{z\in 3^k\Zd\cap \cu_n}
	\sum_{i=1}^{2d}
	\bfJ(z+\cu_k, \bfB^{-\nf12}  e_i, \bfB^{\nf12} e_i\,;\a_m)
	\biggr)^{\!2}
	\Biggr]
	\,.
\end{align}
\end{lemma}
\begin{proof}
By the triangle inequality, for every~$k\leq n$,
\begin{align}
	\label{e.var.splitting}
	\var\big[ \bfB^{\nf12}\bfA_{m,*}^{-1}(\cu_n) \bfB^{\nf12} \bigr]
	&
	\leq
	2\var\Biggl[  \avsum_{z\in 3^k\Zd \cap \cu_n}
	\!\!\!
	\bfB^{\nf12} \bfA_{m,*}^{-1}(z + \cu_k ) \bfB^{\nf12} \Biggr]
	\notag \\ & \qquad
	+
	2\E \Biggl[
	\biggl|
	\avsum_{z\in 3^k\Zd \cap \cu_n}
	\!\!\!
	\bfB^{\nf12}\bigl(
	\bfA_{m,*}^{-1} (z+\cu_k)
	-
	\bfA_{m,*}^{-1} (\cu_n)
	\bigr)\bfB^{\nf12}
	\biggr|^2
	\Biggr]\,.
\end{align}
To estimate the second term on the right-side we use the fact that for every sequence~$\b_1,\ldots,\b_N$ of symmetric positive definite matrices, if we denote their mean by~$\mathbf{m}  \coloneqq  \avsum_{i=1}^N \b_i$ and their harmonic mean by~$\mathbf{h} \coloneqq  \bigl(\avsum_{i=1}^N \b_i^{-1} \bigr)^{\!\!-1}$, then we have, for any other symmetric matrix~$\g$ of the same size,
\begin{align}
	\label{e.no.harm.just.mean}
	\mathbf{m}
	&
	=
	\mathbf{h}+\avsum_{i=1}^N  (\b_i - \g ) \b_i^{-1} (\b_i - \g)
	-
	\underbrace{(\mathbf{h} - \g) \mathbf{h}^{-1} (\mathbf{h} - \g)}_{\geq 0}
	\leq
	\mathbf{h}+\avsum_{i=1}^N  (\b_i - \g ) \b_i^{-1} (\b_i - \g)
	\,.
\end{align}
Applying the above display with~$\{\b_i\}$ given by~$\{\bfA_m^{-1}(z+\cu_k)\,:\, z\in 3^k\Zd\cap\cu_n \}$, we deduce that for any~$\bfB \in \R^{2 d \times 2d}_{\sym,+}$,
\begin{align*}
	\avsum_{z\in 3^k\Zd\cap\cu_n} \!\!\!
	\bfA_m^{-1}(z+\cu_k)
	&
	\leq
	\biggl( \avsum_{z\in 3^k\Zd\cap\cu_n} \bfA_m(z+\cu_k)  \biggr)^{\!-1}
	\notag \\  & \qquad
	+
	\avsum_{z\in 3^k\Zd\cap\cu_n} \!\!\! (\bfB^{-1} - \bfA_m^{-1}(z+\cu_k)  ) \bfA_m(z+\cu_k) (\bfB^{-1} - \bfA_m^{-1}(z+\cu_k) )\,.
\end{align*}
Using~\eqref{e.bfJ.magic2} we have that for any~$\bfB \in \R^{2 d\times 2 d}_{\sym,+}$,~$P \in \R^{2 d}$ and Lipschitz domain~$U \subseteq \Rd$,
\begin{equation*}
	\bfJ(U, \bfB^{-1} P, P;\a_m)
	=
	\frac12 P \cdot \bigl( \bfA_{m,*}^{-1}(U) - \bfA_m^{-1}(U) \bigr) P
	+
	\frac12 P \cdot \bigl(\bfB^{-1} - \bfA_m^{-1}(U) \bigr)  \bfA_m(U) \bigl(\bfB^{-1} - \bfA_m^{-1}(U)  \bigr) P
	\,.
\end{equation*}
Applying~\eqref{e.no.harm.just.mean} with~$\g=\bfB^{-1}$ and using~\eqref{e.bfJ.magic2} to control both~$\bfA_{m,*}^{-1}-\bfA_m^{-1}$ and the quadratic remainder, we obtain from~\eqref{e.CG.bounds.2} that, for any~$\bfB \in \R^{2 d \times 2d}_{\sym,+}$,
\begin{align*}
	\lefteqn{
	\bfB^{\nf12}\bfA_m^{-1}(\cu_n)\bfB^{\nf12}
	} \qquad &
	\notag \\ &
	\leq
	\bfB^{\nf12}\bfA_{m,*}^{-1}(\cu_n)\bfB^{\nf12}
	\notag \\ &
	\leq
	\avsum_{z \in 3^k \Zd \cap \cu_n} \bfB^{\nf12}\bfA_{m,*}^{-1}(z + \cu_k) \bfB^{\nf12}
	\notag \\ &
	\leq
	\bfB^{\nf12}
	\biggl( \avsum_{z\in 3^k\Zd\cap\cu_n} \bfA_m(z+\cu_k)  \biggr)^{\!-1} \bfB^{\nf12}
	+
	2 \sum_{i=1}^{2 d}
	\avsum_{z \in 3^k \Zd \cap \cu_n}
	\!\!\!
	\bfJ(z + \cu_k, \bfB^{-\nf12} e_i, \bfB^{\nf12}e_i;\a_m) \Itwod
	\notag \\ &
	\leq
	\bfB^{\nf12} \bfA_m^{-1}(\cu_n)\bfB^{\nf12}
	+
	2 \sum_{i=1}^{2 d} \avsum_{z \in 3^k \Zd \cap \cu_n}  \bfJ(z + \cu_k, \bfB^{-\nf12} e_i, \bfB^{\nf12}e_i;\a_m) \Itwod \,\,.
\end{align*}
We deduce from the above display that
\begin{align*}
	0 &
	\leq
	\bfB^{\nf12}\Biggl( \avsum_{z\in 3^k\Zd \cap \cu_n}
	\bfA_{m,*}^{-1}(z+\cu_k)
	-
	\bfA_{m,*}^{-1} (\cu_n) \Biggr) \bfB^{\nf12}
	\notag \\ &
	\leq
	2\biggl (
	\avsum_{z\in 3^k\Zd\cap \cu_n}
	\sum_{i=1}^{2 d}
	\bfJ(z+\cu_k, \bfB^{-\nf12}  e_i, \bfB^{\nf12}e_i ;\a_m)\biggr )\Itwod\,,
\end{align*}
and therefore
\begin{align*}
	\lefteqn{
	\E \Biggl[
	\biggl|
	\bfB^{\nf12} \Biggl( \avsum_{z\in 3^k\Zd \cap \cu_n}
	\bfA_{m,*}^{-1} (z+\cu_k)
	-
	\bfA_{m,*}^{-1} (\cu_n)
	\Biggr) \bfB^{\nf12}
	\biggr|^2 \Biggr]
	} \qquad &
	\notag \\ &
	\leq
	4 \E \Biggl[
	\biggl( \avsum_{z\in 3^k\Zd\cap \cu_n}
	\sum_{i=1}^{2d}
	\bfJ(z+\cu_k, \bfB^{-\nf12}  e_i, \bfB^{\nf12} e_i\,;\a_m)
	\biggr)^{\!2}
	\Biggr]
	\,.
\end{align*}
Combining this estimate with~\eqref{e.var.splitting} yields~\eqref{e.var.a.star}.
\end{proof}

We next combine the previous four lemmas to obtain an upper bound on the expectation of~$J(\cu_m,p,q\,;\a_L)$, assuming a moment bound  on~$\mathcal{E}(\cu_L; \a_L, \shom_L)$.

\begin{proposition}
\label{p.combine.under.S}
Let~$L,m\in\Z$ with~$L<m$, and let~$E\geq1$ be such that the full preceding
hypothesis~$\mathcal S(m-1,E)$ is valid. Assume
\begin{equation*}
	15\cstar^{-1}\leq E,\qquad \cgamma\leq E^{-10}\,,
	\qquad
	\delta_1\coloneqq10^9E^2\cgamma\leq\tfrac12\,.
\end{equation*}
There exists~$C(d)<\infty$ such that, if
$C3^{-\frac34(m-L)}\leq\nf14$, then for every~$p,q\in\Rd$ satisfying
\begin{equation}
	\label{e.pq.normed}
	q = \shom_{L} p
	\qqand
	|q|^2 = \shom_{L} \,,
\end{equation}
we have the estimate
\begin{align}
	\label{e.J.bound.by.indyhyp}
	\lefteqn{
	\E \bigl[ J(\cu_m,p,q \,;\a_L)
	\bigr]
	} \qquad &
	\notag \\ &
	\leq
	C
	\delta_1 \bigl( \delta_1 + 3^{-(m-L)} \bigr)
	+
	C \!\!
	\sum_{n=-\infty}^m
	3^{-(m-n)}
	\E \bigl[  J(\cu_n,p,q\,;\a_L) - J(\cu_m,p,q\,;\a_L) \bigr]
	\,.
\end{align}
\end{proposition}
\begin{proof}
Set
$s\coloneqq\nf14$,~$t\coloneqq\nf18$, and~$u\coloneqq\nf1{16}$.
The full hypothesis~$\mathcal S(m-1,E)$, evaluated at scale~$L$ and at
discount~$s$, gives
\begin{equation}
	\E\bigl[
	\mathcal{E}_{s, \infty,2}(\cu_L ; \a_L, \shom_L)^4 \bigr]^{\nf12} \leq \delta_1 \,\,.
	\label{e.moment.bound.onmathcalE}
\end{equation}
Define a (``good'') event~$\mathsf{G}$ by
\begin{equation}
	\label{e.goodevent.combine}
	\mathsf{G}  \coloneqq
	\Bigl\{
	\frac 12 \shom_{L}
	\leq
	{\lambda}_{s,1}(\cu_m;\a_L)
	\leq
	{\Lambda}_{s,1}(\cu_m;\a_L)
	\leq
	2 \shom_{L}
	\Bigr\}\,.
\end{equation}
Let~$\mathsf{B}$ denote the complement of~$\mathsf{G}$.
We note that~\eqref{e.moment.bound.onmathcalE} implies~$\shom_L(\cu_L) \leq 2 \shom_L \leq 4 \shom_{L,*}(\cu_L)$.

\smallskip

\emph{Step 1.} We combine the previous lemmas to obtain an upper bound for~$\E \bigl[ J(\cu_{m},p,q\,;\a_L) \bigr]$. We show that there exists~$C(d,s)<\infty$ such that, if~$C 3^{-(1-s)(m-L)}\leq \nf14$, then
\begin{align}
	\label{e.initial.JL.bound}
	\E \bigl[ J(\cu_{m},p,q\,;\a_L) \bigr]
	&
	\leq
	C
	\sum_{n=L}^m 3^{-(m-n)}
	\E \bigl[  J (\cu_n,p,q\,;\a_L) - J (\cu_m,p,q\,;\a_L) \bigr]
	\notag \\ & \qquad
	+
	C \shom^2_{L}
	\sum_{n=L}^m 3^{-(m-n)}
	\Bigl(
	\var
	\bigl[  \s_{L,*}^{-1}(\cu_n) \bigr]
	+ \bigl| \shom_{L,*}^{-1} (\cu_m)  - \shom_{L,*}^{-1} (\cu_n) \bigr|^2 \Bigr)
	\notag \\ & \qquad
	+
	C\sum_{n=L}^m 3^{-(m-n)}
	\var\bigl[ \s_{L,*}^{-1}(\cu_n ) \k_L(\cu_n) \bigr]
	\notag \\ & \qquad
	+
	C\delta_1^2
	+
	C \E \bigl[ J (\cu_{m-1},p,q\,;\a_L) \indc_{\mathsf{B}}  \bigr]
	\,.
\end{align}
Using the subadditivity of~$J$ and then taking the expectation of~\eqref{e.J.L.bound} and using stationarity, we obtain
\begin{align*}
	\E \bigl[ J(\cu_{m},p,q\,;\a_L) \bigr]
	&
	\leq
	\E \bigl[ J (\cu_{m-1},p,q\,;\a_L) \bigr]
	\notag \\ &
	=
	\E \bigl[ J (\cu_{m-1},p,q\,;\a_L) \indc_{\mathsf{G}}  \bigr]
	+
	\E \bigl[ J (\cu_{m-1},p,q\,;\a_L) \indc_{\mathsf{B}}  \bigr]
	\notag \\ &
	\leq
	C \E \biggl[  {\Lambda}_{s,1}(\cu_m;\a_L)
	\biggl( \frac{{\Lambda}_{s,1}(\cu_m;\a_L)}{{\lambda}_{s,1}(\cu_m;\a_L)} \biggr)^{\frac{s}{1-2s}}
	3^{-2m}\| \nabla v(\cdot,\cu_m,p,q\,;\a_L) \|_{\underline H^{-1}(\cu_m)}^2 \indc_{\mathsf{G}}  \biggr]
	\notag \\ &  \qquad
	+
	C\E \bigl[  J (\cu_{m-1},p,q\,;\a_L) - J (\cu_m,p,q\,;\a_L) \bigr]
	+
	C \E \bigl[ J (\cu_{m-1},p,q\,;\a_L) \indc_{\mathsf{B}}  \bigr]
	\,.
\end{align*}
Applying Lemma~\ref{l.Besov.norms} with~$r=1$, using the triangle inequality and the fixed-scale negative-norm bound for the constant centering vector, squaring both sides of~\eqref{e.Besov.norms.gradient}, taking expectations and using the fact that
\begin{equation*}
	\biggl( \frac{{\Lambda}_{s,1}(\cu_m;\a_L)}{{\lambda}_{s,1}(\cu_m;\a_L)} \biggr)^{1+\frac{s}{1-2s}}
	\indc_{\mathsf{G}}
	\leq
	C
	\qquad \mbox{and} \qquad
	{\Lambda}_{s,1}(\cu_m;\a_L) \indc_{\mathsf{G}}  \leq 2 \shom_{L}\,,
\end{equation*}
we obtain
\begin{align*}
	\lefteqn{
	\E \biggl[  {\Lambda}_{s,1}(\cu_m;\a_L)
	\biggl( \frac{{\Lambda}_{s,1}(\cu_m;\a_L)}{{\lambda}_{s,1}(\cu_m;\a_L)} \biggr)^{\frac{s}{1-2s}}
	3^{-2m}\| \nabla v(\cdot,\cu_m,p,q\,;\a_L) \|_{\underline H^{-1}(\cu_m)}^2 \indc_{\mathsf{G}}  \biggr]
	}
	\quad &
	\notag \\ &
	\leq
	C
	\sum_{n=L}^m
	3^{-(m-n)}
	\E \bigl[  J (\cu_n,p,q\,;\a_L) - J (\cu_m,p,q\,;\a_L) \bigr]
	\notag \\ & \quad
	+
	C
	\shom_L
	\sum_{n=L}^m
	3^{-(m-n)}
	\E \biggl[
	\avsum_{z\in 3^n\Zd \cap \cu_m}
	\bigl|
	\s_{L,*}^{-1} (\cu_m) ( q + \k_L(\cu_m)p  )
	- \s_{L,*}^{-1} (z+\cu_n) ( q + \k_L(z+\cu_n)p  )
	\bigr|^2 \biggr]
	\notag \\ & \quad
	+
	C
	\shom_L \E \Bigl[ \bigl| \s_{L,*}^{-1} (\cu_m)
	( q + \k_L(\cu_m)p  ) - p \bigr|^2 \Bigr]
	+ C 3^{-(1-s)(m-L)}
	\E \biggl[ \| \s_L^{\nf12} \nabla v \|_{\underline{L}^2(\cu_m)}^2 \biggr] \,\,.
\end{align*}
By taking~$C$ large enough in the constraint~$C 3^{-(1-s)(m-L)}\leq \frac14$ and using~\eqref{e.Jenergyv.nosymm}, we have that the last term on the right side is at most~$\frac12 \E \bigl[ J(\cu_{m},p,q\,;\a_L) \bigr]  $;  this term can be reabsorbed. By stationarity, the triangle inequality, and~\eqref{e.annealed.khom.zero}, we have
\begin{align*}
	\lefteqn{
	\E \biggl[
	\avsum_{z\in 3^n\Zd \cap \cu_m}
	\bigl|
	\s_{L,*}^{-1} (\cu_m) ( q + \k_L(\cu_m)p  )
	- \s_{L,*}^{-1} (z+\cu_n) ( q + \k_L(z+\cu_n)p  )
	\bigr|^2 \biggr]
	} \qquad &
	\notag \\ &
	\leq
	2 |q|^2 \E \biggl[
	\avsum_{z\in 3^n\Zd \cap \cu_m}
	\bigl|
	\s_{L,*}^{-1} (\cu_m) -
	\s_{L,*}^{-1} (z+\cu_n)
	\bigr|^2 \biggr]
	\\ &\qquad +
	2 |p|^2 \E \biggl[
	\avsum_{z\in 3^n\Zd \cap \cu_m}
	\bigl|
	(\s_{L,*}^{-1}\k_L)(\cu_m)
	-
	(\s_{L,*}^{-1} \k_L)(z+\cu_n)
	\bigr|^2 \biggr]
	\notag \\ &
	\leq
	6|q|^2 \left(
	\var\bigl[ \s_{L,*}^{-1}(\cu_m ) \bigr]
	+ \var\bigl[ \s_{L,*}^{-1}(\cu_n ) \bigr]
	+
	\bigl| \shom_{L,*}^{-1} (\cu_m) - \shom_{L,*}^{-1} (\cu_n) \bigr|^2
	\right)
	\\ & \qquad
	+
	6 |p|^2
	\left(
	\var\bigl[ (\s_{L,*}^{-1}\k_L)(\cu_m) \bigr]
	+
	\var\bigl[( \s_{L,*}^{-1} \k_L)(\cu_n) \bigr]
	\right)
	\,.
\end{align*}
Combining the previous two displays, using the normalization~$|q|^2 = \shom_L$,~$|p|^2 = \shom^{-1}_L$ assumed in~\eqref{e.pq.normed}, and the fact that~$\shom_{L,*}(\cu_m) \leq \shom_L \leq 2 \shom_{L,*}(\cu_m)$, gives all the displayed terms in~\eqref{e.initial.JL.bound}. The remaining~$C\delta_1^2$ is the residual obtained by taking the limit in Step~2 below; it does not follow from the preceding two-sided bound.

\smallskip

\emph{Step 2.} Using~\eqref{e.moment.bound.onmathcalE}, we prove that there exists a universal constant~$C < \infty$ so that
\begin{equation}
	\label{e.means.downscale.by.defect}
	\bigl| \shom_{L,*}^{-1} (\cu_m)  - \shom_{L,*}^{-1} (\cu_n) \bigr|^2
	\leq
	C \shom_{L}^{-2} \delta_1^2
	\,.
\end{equation}
Taking the expectation of~\eqref{e.J.coarse.grained} and using~\eqref{e.subaddJ.nosymm} and~\eqref{e.annealed.khom.zero} gives, for every~$n\in\Z$ with~$L\leq n \leq m$ and~$e,e'\in \Rd$ with~$|e|=1$,
\begin{align*}
	0\leq
	\shom_{L,*}^{-1} (\cu_n) - \shom_{L,*}^{-1} (\cu_m)
	&
	=
	2 \E \bigl[ J(\cu_n,0,e\,;\a_L)  - J(\cu_m,0,e\,;\a_L) \bigr]
	\notag \\ &
	\leq
	2 \E \bigl[ J(\cu_n,e',e\,;\a_L) - J(\cu_m,e',e\,;\a_L) \bigr]
	\,.
\end{align*}
In view of the normalization in~\eqref{e.pq.normed}, we deduce that
\begin{equation*}
	0\leq
	\shom_{L,*}^{-1} (\cu_n)  - \shom_{L,*}^{-1} (\cu_m)
	\leq
	2\shom_{L,*}^{-1}(\cu_m) \E \bigl[ J (\cu_n,p,q\,;\a_L) - J (\cu_m,p,q\,;\a_L) \bigr]
	\,.
\end{equation*}
Squaring the previous display and using the assumption of~\eqref{e.moment.bound.onmathcalE}, subadditivity,~\eqref{e.homogeneous.J} and~$n\geq L$, we get
\begin{equation*}
	\bigl| \shom_{L,*}^{-1} (\cu_m)  - \shom_{L,*}^{-1} (\cu_n) \bigr|^2
	\leq
	4\shom_{L,*}^{-2}(\cu_m)   \E \bigl[ J (\cu_n,p,q\,;\a_L) - J (\cu_m,p,q\,;\a_L) \bigr]^2
	\leq
	C \delta_1^2 \shom_{L,*}^{-2}(\cu_m)
	\,.
\end{equation*}
Since~$\delta_1\leq 1$, using the assumption of~\eqref{e.moment.bound.onmathcalE} again, we have~$\shom^{-1}_{L,*}(\cu_m) \leq 2 \shom^{-1}_L$ and so this implies~\eqref{e.means.downscale.by.defect}.

\smallskip

\emph{Step 3.}
Using the variance estimate in Lemma~\ref{l.var.bounds}, the independence assumption and the assumption of~\eqref{e.moment.bound.onmathcalE}, we show that, for every~$n\geq L$,
\begin{equation}
	\label{e.var.bound.astar.withS}
	\shom_{L}^{2} \var\big[ \s_{L,*}^{-1}(\cu_n)  \bigr]
	+
	\var\big[ \s_{L,*}^{-1}(\cu_n) \k_L(\cu_n)  \bigr]
	\leq
	C  \delta_1 \bigl( \delta_1 + 3^{-d(n-L)} \bigr)
	\,.
\end{equation}
We bound the left side of~\eqref{e.var.bound.astar.withS} using
\begin{equation*}
	\begin{pmatrix}
		\lambda^{\nf12}\Id & 0 \\
		0 & \lambda^{-\nf12}  \Id
	\end{pmatrix}
	\bfA_{L,*}^{-1} (\cu_n)
	\begin{pmatrix}
		\lambda^{\nf12}\Id & 0 \\
		0 & \lambda^{-\nf12}\Id
	\end{pmatrix}
	=
	\begin{pmatrix}
		\lambda \s_{L,*}^{-1}(\cu_n) & -\s_{L,*}^{-1}(\cu_n) \k_L(\cu_n) \\
		\ast  & \ast
		\end{pmatrix} \qquad \forall \lambda > 0 \,\,,
\end{equation*}
and applying this with~$\lambda = \shom_L$ yields, in view of~\eqref{e.homs.defs},
\begin{equation*}
	\shom_L^2 \var\bigl[\s_{L,*}^{-1} (\cu_n)\bigr] + \var\bigl[ \s_{L,*}^{-1}(\cu_n) \k_L(\cu_n) \bigr]
	\leq
	\var\Bigl[
	\bfAhom_L^{\nf12}
	\bfA_{L,*}^{-1} (\cu_n)
	\bfAhom_L^{\nf12} \Bigr]
	\,.
\end{equation*}
By Lemma~\ref{l.var.bounds}, applied with~$\bfB = \bfAhom_L$ as in~\eqref{e.homs.defs}, we have, for every~$n \geq L$,
\begin{align*}
	\var\Bigl[
	\bfAhom_L^{\nf12}
	\bfA_{L,*}^{-1} (\cu_n)
	\bfAhom_L^{\nf12} \Bigr]
	&
	\leq
	2\var\Biggl[
	\avsum_{z\in 3^L\Zd \cap \cu_n}
	\bfAhom_L^{\nf12} \bfA_{L,*}^{-1} (z+\cu_L) \bfAhom_L^{\nf12} \Biggr]
	\notag \\ & \qquad
	+
	8 \E \Biggl[
	\biggl( \avsum_{z\in 3^L\Zd\cap \cu_n}
	\sum_{i=1}^{2 d}
	\bfJ(z+\cu_L,\bfAhom_L^{-\nf12} e_i, \bfAhom_L^{\nf12}  e_i\,;\a_L)
	\biggr)^{\!2}
	\Biggr]
	\,.
\end{align*}
By assumption~\ref{a.j.frd} and the definition of~$\a_L$, for every~$y,z\in 3^L\Zd$ with~$|y-z|> d 3^{L+1}$, the random variables~$ \bfA_{L,*}^{-1} (z+\cu_L)$ and~$\bfA_{L,*}^{-1} (y+\cu_L)$ are independent.
Consequently, by the assumption of~\eqref{e.moment.bound.onmathcalE}, we have that
\begin{align*}\var\Biggl[
\avsum_{z\in 3^{L}\Zd \cap \cu_n} \!\!\!
\bfAhom_L^{\nf12} \bfA_{L,*}^{-1} (z+\cu_{L})\bfAhom_L^{\nf12}  \Biggr]
&
\leq
C 3^{-d(n-L)}
\var\bigl[\bfAhom_L^{\nf12} \bfA_{L,*}^{-1} (\cu_{L}) \bfAhom_L^{\nf12} \bigr]
\notag \\ &
\leq
C 3^{-d(n-L)}
\E \Bigl[ \bigl|\bfAhom_L^{\nf12} \bfA_{L,*}^{-1}(\cu_{L})\bfAhom_L^{\nf12}  - \Itwod  \bigr|^2 \Bigr]
\notag \\ &
\leq
C  3^{-d(n-L)} \delta_1
\,.
\end{align*}
Similarly,
\begin{equation*}\var \Biggl[
\avsum_{z\in 3^{L}\Zd\cap \cu_n}
\sum_{i=1}^{2d}
\bfJ(z+\cu_{L},\bfAhom_{L}^{-\nf12} e_i, \bfAhom_L^{\nf12} e_i\,;\a_L)
\Biggr]
\leq
C3^{-d(n-L)} \delta_1^2
\end{equation*}
and
\begin{equation*}
	\sup_{|e|=1}
	\E \Bigl[ \bfJ(\cu_{L},\bfAhom_{L}^{-\nf12} e, \bfAhom_L^{\nf12} e\,;\a_L)\Bigr]^2
	\leq C\delta_1^2
	\,.
\end{equation*}
Combining the previous three displays yields~\eqref{e.var.bound.astar.withS}.

\smallskip

\emph{Step 4.}
We prove that
\begin{equation}
	\label{e.badeventbound.sect42}
	\E \bigl[ J (\cu_{m-1},p,q\,;\a_L) \indc_{\mathsf{B}}  \bigr]
	\leq
	C\delta_1^2 \,.
\end{equation}
Set~$s=1/4$,~$t=1/8$, and~$u=1/16$. The change from~$q=1$ to~$q=2$ and~\eqref{e.bound.Lambdas.by.Es} give~$\mathsf B\subseteq\{\mathcal E_{t,\infty,2}(\cu_m;\a_L,\shom_L)\geq1/5\}$. Splitting the depth sum at scale~$L$, applying subadditivity above it and the strict gap~$u<t$ below it, the hypothesis~$\mathcal S(m-1,E)$ yields
\begin{equation*}
	\mathcal E_{t,\infty,2}(\cu_m;\a_L,\shom_L)^2
	\leq\O_{\Gamma_{1/4}}(C(d)\delta_1)\,.
\end{equation*}
Markov's inequality now yields
\begin{equation}
	\label{e.bad.B.bound}
	\P[\mathsf B]
	\leq
	\P\bigl[\mathcal E_{t,\infty,2}(\cu_m;\a_L,\shom_L)\geq\nf15\bigr]
	\leq C\delta_1^2\,.
\end{equation}
By~\eqref{e.bfJ.general},~\eqref{e.mathcalE.monotone.ordered},~\eqref{e.bound.one.cube.by.mathcalE},~\eqref{e.pq.normed} and the preceding Orlicz estimate, we have
\begin{equation*}
	\E \bigl[ J (\cu_{m-1},p,q\,;\a_L)^2 \bigr]
	\leq C\delta_1^2
	\,.
\end{equation*}
Combining this with~\eqref{e.bad.B.bound}, we obtain
\begin{equation*}
	\E \bigl[ J (\cu_{m-1},p,q\,;\a_L) \indc_{\mathsf{B}}  \bigr]
	\leq
	\E \bigl[ J (\cu_{m-1},p,q\,;\a_L)^2 \bigr]^{\nf12}
	\P [ \mathsf{B} ]^{\nf12}
	\leq
	C \delta_1^2
	\,.
\end{equation*}
This completes the proof of~\eqref{e.badeventbound.sect42}.

\smallskip

\emph{Step 5.} The conclusion. Combining~\eqref{e.initial.JL.bound},~\eqref{e.means.downscale.by.defect},~\eqref{e.var.bound.astar.withS} and~\eqref{e.badeventbound.sect42} gives the desired bound~\eqref{e.J.bound.by.indyhyp} and completes the proof of the proposition.
\end{proof}

We next use~\eqref{e.J.bound.by.indyhyp} together with a concentration inequality to show that~$J$ contracts above the infrared cutoff scale.

\begin{proof}[Proof of Proposition~\ref{p.homogenization.step}]
Since the response is a nonnegative quadratic form in~$e$, there is a dimensional constant~$C_d$ such that every unit direction is bounded by~$C_d$ times the maximum over the coordinate directions. We apply the fixed-direction argument below simultaneously to those directions, with~$\ep$ replaced by~$C_d^{-1}\ep$; finite maxima preserve the two Orlicz classes at a dimensional cost. To avoid extra notation, we continue to write~$\ep$ for the reduced tolerance.
The assumed estimate~\eqref{e.new.induction.for.shom}, applied with~$s=\nf14$, gives
\begin{equation*}
\E\bigl[\mathcal{E}_{\nf14,\infty,2}(\cu_L;\a_L,\shom_L)^4\bigr]^{\nf12} \leq C E^2 \cgamma\,, 
\qquad \forall L \in (-\infty, m-1] \cap \Z  \,.
\end{equation*}
The mean-contraction argument of Steps~1--4 above has a standalone form at this hypothesis: its error-moment input is the preceding display, while the diffusivity comparisons are supplied by the already established propagation estimates under~$15\cstar^{-1}\leq E$ and~$\cgamma\leq E^{-10}$. Applying that form with~$\delta_1\simeq E^2\cgamma$, and then using subadditivity, gives for every~$|e|=1$,
\begin{equation}
	\E[J(\cu_n, \shom_{L}^{-\nf12} e, \shom_{L}^{\nf12} e \,; \a_L)] \leq C E^2 \cgamma\,,  \qquad \forall L,n \in \Z \quad \mbox{with~$L\leq m-1$ and $L\leq n\leq m$} \,\,.
	\label{e.iter.init}
\end{equation}
Fix~$L\leq m-1$ and~$|e|=1$, put~$T=m-L$ and
\begin{equation*}
	F(j)\coloneqq\E\bigl[J(\cu_{L+j},\shom_L^{-\nf12}e,\shom_L^{\nf12}e\,;\a_L)\bigr]\,,
	\qquad \delta_1\coloneqq C_0E^2\cgamma\,,
\end{equation*}
where~$C_0(d)$ is fixed large enough. By~\eqref{e.iter.init}, positivity and subadditivity, we have~$0\leq F(j)\leq\delta_1$ for~$0\leq j\leq T$. The finite version of Steps~1--4 in the proof of~\eqref{e.J.bound.by.indyhyp} gives constants~$A,C,j_0$ depending only on~$d$ such that
\begin{equation}
	\label{e.iter.corridor}
	F(j)\leq A\delta_1\bigl(\delta_1+3^{-j}\bigr)+C\sum_{r=1}^j3^{-r}\bigl(F(j-r)-F(j)\bigr)
\end{equation}
whenever~$j_0<j\leq T$ and~$C|\log\ep|\leq j$. Set~$t_\ep\coloneqq\max\{j_0,\lceil C|\log\ep|\rceil\}$. After shifting by~$t_\ep$, the terms in~\eqref{e.iter.corridor} which cross the entrance scale satisfy
\begin{equation*}
	\sum_{r=j+1}^{t_\ep+j}3^{-r}\bigl(F(t_\ep+j-r)-F(t_\ep+j)\bigr)
	\leq\delta_1\sum_{r=j+1}^\infty3^{-r}
	=\frac12\delta_1 3^{-j}\,.
\end{equation*}
Choose~$\rho\in(0,1)$, depending only on~$d$, sufficiently close to one that
\begin{equation*}
	C\sum_{r=1}^j3^{-r}\bigl(\rho^{j-r}-\rho^j\bigr)\leq\frac12\rho^j
	\qquad\forall j\in\N\,.
\end{equation*}
A strong induction in~\eqref{e.iter.corridor} then yields~$F(t_\ep+j)\leq K\delta_1(\delta_1+\rho^j)$ whenever~$t_\ep+j\leq T$, for a constant~$K(d)<\infty$. The last assumption in~\eqref{e.gamma.condition.homog} and a choice~$j_\ep\leq C|\log\ep|$ therefore give~$K\delta_1(\delta_1+\rho^{j_\ep})\leq\frac12\ep E^2\cgamma$. Setting~$k_1\coloneqq t_\ep+j_\ep$, we obtain
\begin{equation}
	\label{e.iter.post.new}
	\E \bigl[ J (\cu_n,\shom_{L}^{-\nf12} e, \shom_{L}^{\nf12} e \,;\a_L) \bigr]
	\leq
	\frac12 \ep E^2 \cgamma \,,
	\quad \forall L,n \in \Z\,, \ \ L\leq m-1\,,\quad L+k_1\leq n\leq m \,\,.
\end{equation}
By Proposition~\ref{p.concentration} and~\eqref{e.new.induction.for.shom}, we have, for every~$|e|=1$ and~$L,n\in \Z$ with~$L\leq n\leq m-1$,
\begin{multline*}
	\biggl|
	\avsum_{z\in 3^n\Zd \cap \cu_m}
	\bigl( J (z+\cu_n,\shom_{L}^{-\nf12} e, \shom_{L}^{\nf12} e \,;\a_L) - \E \bigl[ J (\cu_n,\shom_{L}^{-\nf12} e, \shom_{L}^{\nf12} e \,;\a_L) \bigr] \bigr)  \biggr|
	\\
	\leq
	\O_{\Gamma_{1}}(E^2 \cgamma 3^{-\frac d2(m-n)})
	+ \O_{\Gamma_{\nf{1}{4} }}(\exp(-2E^{-3} \cgamma^{-1}) 3^{-\frac d2(m-n)})
	\,.
\end{multline*}
We put~$k_2 \coloneqq  \lceil \log_3 (10\ep^{-1}) \rceil$,~$k \coloneqq k_1+k_2$, choose~$n  \coloneqq  L +k_1$ and then sum over~$L \in \Z$ with~$L\leq m-k$ to obtain
\begin{multline}
	\label{e.union.bound.fluct}
	\sup_{L \in \Z \cap (-\infty,m-k]}
	\biggl|
	\avsum_{z\in 3^{L +k_1} \Zd \cap \cu_m}
	\bigl( J (z+\cu_{L +k_1},\shom_{L}^{-\nf12} e, \shom_{L}^{\nf12} e \,;\a_L) - \E \bigl[ J (\cu_{L +k_1},\shom_{L}^{-\nf12} e, \shom_{L}^{\nf12} e \,;\a_L) \bigr] \bigr)  \biggr|
	\\
	\leq
	\O_{\Gamma_{1}}\Bigl(\frac12 \ep E^2 \cgamma \Bigr)
	+ \O_{\Gamma_{\nf{1}{4} }}
	\bigl(\ep \exp(-2E^{-3} \cgamma^{-1})  \bigr)
	\,.
\end{multline}
The finite-direction reduction at the start of the proof now gives the supremum over~$|e|\leq1$. Subadditivity and~\eqref{e.iter.post.new}--\eqref{e.union.bound.fluct} give
\begin{align*}
	\sup_{L \in \Z \cap (-\infty,m-k]}
	J (\cu_m,\shom_{L}^{-\nf12} e, \shom_{L}^{\nf12} e \,;\a_L)
	&
	\leq
	\sup_{L \in \Z \cap (-\infty,m-k]}
	\avsum_{z\in 3^{L +k_1} \Zd \cap \cu_m}
	J (z+\cu_{L +k_1},\shom_{L}^{-\nf12} e, \shom_{L}^{\nf12} e \,;\a_L)
	\notag \\ &
	\leq
	\frac12 \ep E^2 \cgamma
	+
	\O_{\Gamma_{1}}\Bigl(\frac12 \ep E^2 \cgamma \Bigr)
	+ \O_{\Gamma_{\nf{1}{4} }}
	\Bigl(\ep \exp(-2E^{-3} \cgamma^{-1})  \Bigr)
	\notag \\ &
	\leq
	\O_{\Gamma_{1}}\bigl(\ep E^2 \cgamma \bigr)
	+ \O_{\Gamma_{\nf{1}{4} }}
	\bigl(\ep \exp(-2E^{-3} \cgamma^{-1})  \bigr)
	\,.
\end{align*}
As~$k=k_1+k_2\leq C|\log\ep|$, restoring the original tolerance completes the proof.
\end{proof}

\subsection{Propagation of the multiscale estimate}
\label{ss.buckle}

In this section we assemble the main ingredients of the previous four subsections to complete the proof of Proposition~\ref{p.induction.bounds}. The main remaining task is to propagate the first part of the induction hypothesis, namely the multiscale bound on the defect variable~$\mathcal{E}_{s,\infty,2}(\cu_m;\a_m,\shom_m)$, from scale~$m-1$ to scale~$m$. Proposition~\ref{p.multiscale.estimate} is essentially this propagation step: assuming~$\mathcal{S}(m-1,E)$, it produces the same type of Orlicz-tail bound for~$\mathcal{E}_{s,\infty,2}(\cu_m;\a_m,\shom_m)$, up to a controllable loss factor~$\ep$ which we later fix as an absolute constant. In other words, this proposition is the ``buckle'' that closes the induction for the homogenization error.

\begin{proposition}
\label{p.multiscale.estimate}
Assume that~$m \in\Z$ and~$E \in [1,\infty)$ are such that~$\S(m-1, E)$ is valid.
There exists~$C(d)<\infty$ such that, for every~$\ep \in (0,\nf12]$, if
\begin{equation}
	\label{e.param.conditions.in.main}
	E \geq \cstar^{-1} \ep^{-C}
	\qand
	\cgamma \leq
	E^{-10}
	\,,
\end{equation}
then, for every~$s \in [8\cgamma,1]$,
\begin{equation}
	\label{e.complete.wrapping}
	\mathcal{E}_{s, \infty, 2}(\cu_m; \a_m, \shom_m)
	\leq
	\O_{\Gamma_{2} }
	\bigl( C E s^{-1} (\ep\cgamma)^{\nf12}  \bigr) +
	\O_{\Gamma_{\nf12}}
	\bigl( C \ep s^{-2}
	\exp\bigl( - E^{-3}\cgamma^{-1} \bigr) \bigr)
	\,.
\end{equation}
For each~$s$, the two terms may be realized by the same pair of random
witnesses. If, in addition,~$s\leq C\ep$, that same~$\Gamma_2$ witness satisfies
the sharper bound
\begin{equation}
	\mathcal E_{s,\infty,2}(\cu_m;\a_m,\shom_m)
	\leq
	\O_{\Gamma_2}\bigl(C\ep E s^{-1}\cgamma^{\nf12}\bigr)
	+\O_{\Gamma_{\nf12}}\bigl(C\ep s^{-2}
	\exp(-E^{-3}\cgamma^{-1})\bigr)\,.
	\label{e.complete.wrapping.sharp}
\end{equation}
\end{proposition}

The proof of Proposition~\ref{p.multiscale.estimate} is an assembly argument. The contraction input from the previous subsection controls the functional~$J$ for a deeper cutoff~$\a_{l-h}$ on subcubes~$z+\cu_l$ after climbing~$h\asymp |\log\ep|$ scales. The remaining task is then to transfer this information back to the genuine coefficients~$\a_m$ and the correct normalization~$\shom_m$ at scale~$\cu_m$.

\smallskip

This is accomplished in the following lemma, which is the point in our RG argument in which we \emph{undo} the infrared cutoff and return to the true coefficients at the working scale. Concretely, we localize the global defect~$\mathcal{E}_{s,\infty,2}(\cu_m;\a_m,\shom_m)$ into contributions coming from subcubes~$z+\cu_l$ and from a \emph{lower} cutoff field~$\a_{l-h}$, at the matching normalization~$\shom_{l-h}$ (which we can control via the results from the previous subsection) and errors caused by injecting the  missing layers~$\k_m-\k_{l-h}$.

\smallskip

This lemma is also where the coarse-grained sensitivity theory from Section~2 is used in an essential way. The comparison between~$\a_m$ and~$\a_{l-h}$ on each cube~$z+\cu_l$ is obtained by applying Lemma~\ref{l.J.sensitivity} on the local good event~$\mathcal{Q}(l,l-h,z)$ defined in~\eqref{e.good.local.events}. We separate the contribution due to the bad events and later control it using the improved tail bounds from Section~\ref{ss.coarse.grained.sensitivity}.

\begin{lemma}[Localization lemma]
\label{l.localization.mathcalE}
Assume that~$m \in\Z$ and~$E \in [1,\infty)$ are such that~$\S(m-1, E)$ is valid. There exists~$C(d,\cstar)<\infty$ such that, if
\begin{equation}
	\label{e.localization.lemma.conditions}
	\cgamma\leq E^{-10} \,,
	\qand
	E \geq C\cstar^{-1} \,,
\end{equation}
then we have, for every~$h\in\N \cap [1, \cgamma^{-1}]$ and~$s\in [8\cgamma,1]$,
\begin{align}
	\mathcal{E}_{s,\infty,2 } (\cu_m ; \a_m, \shom_m )^2
	&
	\leq
	C \css{2s}
	\sum_{l=-\infty}^m
	3^{- s (m-l)}
	\biggl( \avsum_{z \in 3^l\Zd \cap \cu_m}
	\max_{|e| =1} \bigl( J(z +\cu_l,\shom_{l-h}^{-\nf 12} e, \shom_{l-h}^{\nf 12} e \,; \a_{l-h} ) \bigr)^{\frac{d}{s}}
	\biggr)^{\!\frac{s}{d}}
	\notag \\ & \quad
	+
	C \css{2s}
	\sum_{l=-\infty}^m
	3^{-s (m-l)}
	\biggl( \avsum_{z \in 3^l\Zd \cap \cu_m}
	\max_{|e| =1} \bigl( J(z +\cu_l,\shom_{l-h}^{-\nf 12} e, \shom_{l-h}^{\nf 12} e \,; \a_{l-h}^t ) \bigr)^{\frac{d}{s}}
	\biggr)^{\!\frac{s}{d}}
	\notag \\ & \quad
	+
	C\cstar^{-1} \cgamma
	s
	\sum_{l=-\infty}^m
	3^{-s (m-l)}
	\biggl(
	\avsum_{z \in 3^l\Zd \cap \cu_m}
	\bigl(
	3^{2l-\cgamma l}\| \k_m-\k_{l-h} \|_{\underline{W}^{2,\infty}(z+\cu_l)}
	\bigr)^{\frac{2d}{s}}
	\biggr)^{\!\frac{s}{d}}
	\notag \\ & \quad
	+
	8\,
	\mathcal{E}_{\nf s4,\infty,2 }(\cu_m;\a_m,\shom_m )^2
	\max_{l\leq m}
	3^{-\frac s2 (m-l)}
	\biggl( \avsum_{z \in 3^l\Zd \cap \cu_m}
	\indc_{\neg  \mathcal{Q}(l,l-h,z) }
	\biggr)^{\!\frac{s}{d}}
	\notag \\ & \quad
	+
	C\cgamma^2 \bigl( h^2 + s^{-2} +  E^4 \left|\log \cgamma\right|^4 \bigr)
	\,.
	\label{e.localization.mathcalE.estimate}
\end{align}
\end{lemma}
\begin{proof}
We apply~\eqref{e.mathcalE.infty.to.q} with~$p = 2ds^{-1}$ and~$q =2$ to get
\begin{align}
	\mathcal{E}_{s,\infty,2 } (\cu_m ; \a_m, \shom_m )^2
	&
	\leq
	\css{2s}
	\sum_{l=-\infty}^m
	3^{-s (m-l)}
	\biggl( \avsum_{z \in 3^l\Zd \cap \cu_m}
	\max_{|e| =1} \bigl( J(z {+} \cu_l,\shom_m^{-\nf 12} e, \shom_m^{\nf 12} e \,; \a_m \bigr) \bigr)^{\frac{d}{s}}
	\biggr)^{\!\frac{s}{d}}
	\notag \\ & \qquad
	+
	\css{2s}
	\sum_{l=-\infty}^m
	3^{-s (m-l)}
	\biggl( \avsum_{z \in 3^l\Zd \cap \cu_m}
	\max_{|e| =1} \bigl( J(z {+} \cu_l,\shom_m^{-\nf 12} e, \shom_m^{\nf 12} e \,; \a_m^t \bigr) \bigr)^{\frac{d}{s}}
	\biggr)^{\!\frac{s}{d}}
	\,.
	\label{e.mathcal.E.breakdown}
\end{align}
Set~$\Delta\coloneqq m-l+h$. We first switch from~$\shom_m$'s to~$\shom_{l-h}$'s using~\eqref{e.shaking.lambda}, applied with~$\lambda = \shom_m\shom_{l-h}^{-1}$.
If~$m\leq\mstarstar$, then~\eqref{e.basecase.homogenization} at~$m$ and~$l-h$ gives ratios at most~$2$ and an error at most~$\cgamma^2$, which is absorbed by the displayed remainder since~$m-l+h\geq1$; hence we may assume~$m>\mstarstar$.
For this we recall that~\eqref{e.shom.m.vs.shom.n} gives
\begin{equation}
	|\shom_{l'}\shom_{l}^{-1}-1|
	+
	|\shom_{l}\shom_{l'}^{-1}-1|
	\leq
	C \min\bigl\{  \cgamma |l-l'| + E^2 \cgamma \left|\log \cgamma\right|^2, 1    \bigr\} 3^{\cgamma |l-l'|}
	\qquad \forall l, l' \in (-\infty,m]\cap\Z
	\,\,,
	\label{e.switchtheshoms.betterer.again}
\end{equation}
and that~$\shom_m^{-1}\shom_{l-h} \leq 4$. We obtain
\begin{align*}
	\lefteqn{
	J(z+\cu_l,\shom_m^{-\nf12}e,\shom_m^{\nf12}e\,; \a_m )
	\indc_{ \mathcal{Q}(l,l-h,z) }
	} \quad &
	\notag \\ &
	\leq
	C\shom_m^{-1}\shom_{l-h}
	J(z+\cu_l,\shom_{l-h}^{-\nf12}e,\shom_{l-h}^{\nf12}e\,; \a_m )
	+
	C\shom_m^{-1}\shom_{l-h} \bigl (\shom_m\shom_{l-h}^{-1}-1\bigr)^2
	\shom_{l-h} |\s_{m,*}^{-1} (z+\cu_l)|\indc_{ \mathcal{Q}(l,l-h,z) }
	\notag \\ &
	\leq
	C
	J(z+\cu_l,\shom_{l-h}^{-\nf12}e,\shom_{l-h}^{\nf12}e\,; \a_m )
	+
	C\min\{\cgamma^2\Delta^2,1\}3^{2\cgamma\Delta}
	+CE^4\cgamma^2|\log\cgamma|^4\,3^{2\cgamma\Delta}
	\,.
\end{align*}
The last inequality on the right side was obtained thanks to Lemma~\ref{l.lambda.sensitivity}, which yields, in view of~\eqref{e.ellipticities.monotone.ordered}, the definition of~$ \mathcal{Q}(l,l-h,z) $ in~\eqref{e.good.local.events} and~\eqref{e.switchtheshoms.betterer.again},
\begin{equation*}
	|\s_{m,*}^{-1} (z+\cu_l)|\indc_{ \mathcal{Q}(l,l-h,z) }
	\leq
	\lambda_{\nf18,2}^{-1}  (z+\cu_l;\a_m)
	\indc_{ \mathcal{Q}(l,l-h,z) }
	\leq
	2 \lambda_{\nf18,2}^{-1}  (z+\cu_l;\a_{l-h} )
	\indc_{ \mathcal{Q}(l,l-h,z) }
	\leq
	C \shom_{l-h}^{-1} \,.
\end{equation*}
We next apply the sensitivity estimate of
Lemma~\ref{l.J.sensitivity} with~$\delta=1$ to obtain, for every~$|e|=1$,~$l\in (-\infty,m]\cap \Z$ and~$z\in 3^l\Zd\cap \cu_m$,
\begin{align}
	\lefteqn{
	J(z+\cu_l,\shom_{l-h}^{-\nf12}e,\shom_{l-h}^{\nf12}e\,; \a_m )
	\indc_{ \mathcal{Q}(l,l-h,z) }
	} \qquad  &
	\notag \\ &
	\leq
	3 J(z+\cu_l,\shom_{l-h}^{-\nf12}e,\shom_{l-h}^{\nf12}e\,; \a_{l-h})
	+
	C
	\shom_{l-h-1}^{-2}
	3^{4l}\| \k_m-\k_{l-h} \|_{\underline{W}^{2,\infty}(z+\cu_l)}^2
	\notag \\ &
	\leq
	3 J(z+\cu_l,\shom_{l-h}^{-\nf12}e,\shom_{l-h}^{\nf12}e\,; \a_{l-h})
	+
	C\cstar^{-1} \cgamma
	3^{4l-2\cgamma l}\| \k_m-\k_{l-h} \|_{\underline{W}^{2,\infty}(z+\cu_l)}^2
	\,,
	\label{e.J.sensitivity.apppp}
\end{align}
where in the last line we used~\eqref{e.shom.h.bounds} from the induction hypothesis and~$h\cgamma\leq 1$.
Combining these inequalities, summing over~$z\in 3^l\Zd \cap \cu_m$ and applying the triangle inequality, we obtain
\begin{align*}
	\lefteqn{
	\css{2s}
	\sum_{l=-\infty}^m
	3^{-s (m-l)}
	\biggl( \avsum_{z \in 3^l\Zd \cap \cu_m}
	\max_{|e| =1} \bigl( J(z {+} \cu_l,\shom_m^{-\nf 12} e, \shom_m^{\nf 12} e \,; \a_m \bigr) \bigr)^{\frac{d}{s}}
	\indc_{ \mathcal{Q}(l,l-h,z) }
	\biggr)^{\!\frac{s}{d}}
	} \quad &
	\notag \\ &
	\leq
	\css{2s}
	\sum_{l=-\infty}^m
	3^{-s (m-l)}
	\biggl( \avsum_{z \in 3^l\Zd \cap \cu_m}
	\max_{|e| =1} \bigl( J(z +\cu_l,\shom_{l-h}^{-\nf 12} e, \shom_{l-h}^{\nf 12} e \,; \a_{l-h} ) \bigr)^{\frac{d}{s}}
	\biggr)^{\!\frac{s}{d}}
	\notag \\ & \qquad
	+
	C\cstar^{-1} \cgamma
	\css{2s}
	\sum_{l=-\infty}^m
	3^{-s (m-l)}
	\biggl( \avsum_{z \in 3^l\Zd \cap \cu_m}
	\bigl(
	3^{4l-2\cgamma l}\| \k_m-\k_{l-h} \|_{\underline{W}^{2,\infty}(z+\cu_l)}^2
	\bigr)
	^{\frac{d}{s}}
	\biggr)^{\!\frac{s}{d}}
	\notag \\ & \qquad
	+
	C \css{2s}
	\sum_{l=-\infty}^m
	3^{-s(m-l)}
	\Bigl(
	\min\{\cgamma^2(m-l+h)^2,1\}3^{2\cgamma(m-l+h)}
	+E^4\cgamma^2|\log\cgamma|^4 3^{2\cgamma(m-l+h)}
	\Bigr)
	\,.
\end{align*}
Since~$s \geq 8\cgamma$ and~$3^{-s\cgamma^{-1}} \leq C s^{-2} \cgamma^2$, we have
\begin{equation*}
	\css{2s}
	\sum_{l=-\infty}^m
	3^{-s (m-l)}
	\min\{ \cgamma^2 (m-l+h)^2 ,1 \} 3^{2 \cgamma (m-l+h)}
	\leq
	C\bigl( h^2 + s^{-2} \bigr)\cgamma^2
	\,.
\end{equation*}
For the bad events, we estimate
\begin{align*}
	\lefteqn{
	\css{2s} \sum_{l=-\infty}^m
	3^{-s (m-l)}
	\biggl( \avsum_{z \in 3^l\Zd \cap \cu_m}
	\max_{|e| =1} \bigl( J(z {+} \cu_l,\shom_m^{-\nf 12} e, \shom_m^{\nf 12} e \,; \a_m \bigr) \bigr)^{\frac{d}{s}}
	\indc_{\neg  \mathcal{Q}(l,l-h,z) }
	\biggr)^{\!\frac{s}{d}}
	} \qquad &
	\notag \\ &
	\leq
	\css{2s} \sum_{l=-\infty}^m
	3^{-s (m-l)}
	\max_{z\in 3^l\Zd \cap \cu_m}
	\max_{|e| =1}
	\bigl( J(z {+} \cu_l,\shom_m^{-\nf 12} e, \shom_m^{\nf 12} e \,; \a_m) \bigr)
	\biggl( \avsum_{z \in 3^l\Zd \cap \cu_m}
	\indc_{\neg  \mathcal{Q}(l,l-h,z) }
	\biggr)^{\!\frac{s}{d}}
	\notag \\ &
	\leq
	\frac{\css{2s}}{\css{\nf s2}}\,
	\mathcal{E}_{\nf s4,\infty,2}(\cu_m;\a_m,\shom_m )^2
	\max_{l\leq m}
	3^{-\frac 12 s (m-l)}
	\biggl( \avsum_{z \in 3^l\Zd \cap \cu_m}
	\indc_{\neg  \mathcal{Q}(l,l-h,z) }
	\biggr)^{\!\frac{s}{d}}
	\,.
\end{align*}
Here the weighted maximum is extracted before the~$\mathcal E_{s/4,\infty,2}^2$ scale series is summed. Since~$\css{2s}/\css{\nf s2}=(1+3^{-s/2})(1+3^{-s})\leq4$, the preceding display controls the first frame leg on the right side of~\eqref{e.mathcal.E.breakdown} with coefficient at most~$4$. The transposed leg follows from the adjoint redundancy (equivalently, samplewise~$\mathbf j\mapsto-\mathbf j$), so the two legs have total coefficient at most~$8$. Combining this with~\eqref{e.J.sensitivity.apppp} completes the proof.
\end{proof}

We are now ready to close the RG induction loop.

\begin{proof}[Proof of Proposition~\ref{p.multiscale.estimate}]
By~\eqref{e.localization.mathcalE.estimate}, it suffices to prove the following three bounds:
\begin{multline}
	\label{e.what.homogenization.gives}
	s\sum_{l=-\infty}^m
	3^{- s (m-l)}
	\biggl( \avsum_{z \in 3^l\Zd \cap \cu_m}
	\max_{|e| =1} \bigl( J(z +\cu_l,\shom_{l-h}^{-\nf 12} e, \shom_{l-h}^{\nf 12} e \,; \a_{l-h} ) \bigr)^{\frac{d}{s}}
	\biggr)^{\!\frac{s}{d}}
	\\
	\leq
	\O_{\Gamma_{1}}(C s^{-1} \ep E^2 \cgamma)
	+ \O_{\Gamma_{\nf{1}{4} }}(C s^{-4} \ep^2 \exp(-2E^{-3} \cgamma^{-1}))\,,
\end{multline}
\begin{equation}
	\label{e.wave.sizes}
	s\sum_{l=-\infty}^m
	3^{-s (m-l)}
	\biggl(
	\avsum_{z \in 3^l\Zd \cap \cu_m}
	\bigl(
	3^{(2-\cgamma )l}\| \k_m-\k_{l-h} \|_{\underline{W}^{2,\infty}(z+\cu_l)}
	\bigr)^{\frac{2d}{s}}
	\biggr)^{\!\frac{s}{d}}
	\leq
	\O_{\Gamma_1}
	( s^{-2} \ep^{-C} )
\end{equation}
and
\begin{multline}
	\label{e.local.bad.events.summed}
	\mathcal{E}_{\nf s4,\infty,2}(\cu_m;\a_m,\shom_m )^2
	\max_{l\leq m}
	3^{-\frac 12 s (m-l)}
	\biggl( \avsum_{z \in 3^l\Zd \cap \cu_m}
	\indc_{\neg  \mathcal{Q}(l,l-h,z) }
	\biggr)^{\!\frac{s}{d}}
	\\
	\leq
	\O_{\Gamma_{1}}
	\bigl( C\cstar^{-1} s^{-2} \ep^{-C}  \cgamma\bigr)
	+
	\O_{\Gamma_{\nf 27}}\bigl(   \exp \bigl( - C^{-1} E^{-2} \cgamma^{-1}\bigr)  \bigr)
	\,.
\end{multline}
The first estimate~\eqref{e.what.homogenization.gives} follows from Proposition~\ref{p.homogenization.step}, applied with~$\ep^2$, the triangle inequality~\eqref{e.Gamma.sigma.triangle},~$s\geq 8\cgamma$ and~$\ep^2\leq\ep$, provided that~$h\coloneqq\max\{1,\lceil C|\log\ep|\rceil\}$ for~$C$ large enough. Retaining the factor~$\ep^2$ gives the sharper estimate with the same witnesses; in its rare lane, the cube maximum and the weighted scale sum contribute the factor~$s^{-4}$.

For~\eqref{e.wave.sizes}, decompose~$\k_m-\k_{l-h}$ layerwise and apply~\eqref{e.W.1.inf.bound} with~$n=k-1$,~$m=k$. Transport to scale~$l$ gives the summable weights~$\max\{3^{l-k},3^{2(l-k)}\}$ for~$k>l$;~\eqref{e.Gamma.sigma.triangle} then gives~\eqref{e.wave.sizes}.

\smallskip

We turn to the proof of the third estimate~\eqref{e.local.bad.events.summed}.
According to~\eqref{e.bound.Lambdas.by.Es} and Proposition~\ref{p.cg.ellipticity.bounds} with~$\sigma=\frac17$ and parameter~$s/4\geq2\cgamma$, which is valid by the assumption of~\eqref{e.param.conditions.in.main},
we have
\begin{align}
	\mathcal{E}^2_{\nf s4,\infty,2}(\cu_m; \a_m, \shom_m)
	&
	\leq
	2 \lambda_{\nf s4,2}^{-1} (\cu_m;\a_m) \shom_m
	+
	2
	\Lambda_{\nf s4,2} (\cu_m;\a_m) \shom_m^{-1}
	\notag \\ &
	\leq
	C
	+
	\O_{\Gamma_{1}} \bigl( C\cstar^{-1} s^{-2} \cgamma  \bigr)
	+
	\O_{\Gamma_{\nf 27}}\bigl( \exp \bigl( - C^{-1} E^{-2} \cgamma^{-1}\bigr)  \bigr)
	\,.
	\label{e.mathcalE.crude.bound}
\end{align}
Here Lemma~\ref{l.shom.continuity} permits the change from~$\shom_{m-1}$ to~$\shom_m$.
Using~\eqref{e.indc.O.sigma} and~\eqref{e.bad.event.Q.estimate} and our constraints in~\eqref{e.param.conditions.in.main}, we have that, for every~$\rho \in (0,\infty)$,
\begin{equation}
	\indc_{\neg  \mathcal{Q}(l,l-h,z) }
	\leq
	\O_{\Gamma_{\rho}} \bigl( \bigl| \log \P[\neg  \mathcal{Q}(l,l-h,z) ] \bigr|^{-\nf1\rho} \bigr)
	\leq
	\O_{\Gamma_{\rho}} \bigl(
	C^{\nf1\rho}
	\cstar^{-\nf1\rho}
	3^{\nf {5h}\rho}
	\cgamma^{\nf1\rho}
	\bigr)
	\,.
	\label{e.badevent.indicator.propagate}
\end{equation}
Take~$\rho=s/d$.  For
$Y_l\coloneqq\avsum_z\indc_{\neg\mathcal Q(l,l-h,z)}$ we have
$0\leq Y_l\leq1$.  Averaging~\eqref{e.badevent.indicator.propagate} first and
then applying Markov's inequality to~$Y_l$ shows that raising~$Y_l$ to the
power~$s/d$ converts the tail to~$\Gamma_1$ without an Orlicz triangle loss.
Thus
\begin{equation*}
	\biggl( \avsum_{z \in 3^l\Zd \cap \cu_m}
	\indc_{\neg  \mathcal{Q}(l,l-h,z) } \biggr)^{\!\frac{s}{d}}
	\leq
	\O_{\Gamma_{1}} \bigl(
	C s^{-1} \cstar^{-1}
	3^{5h}
	\cgamma
	\bigr)\,.
\end{equation*}
Combining the above, using~$3^{5h} \leq C\ep^{-C}$, and taking the
weighted scale sum supplies the second factor~$s^{-1}$.  Extracting
the maximum before this sum yields~\eqref{e.local.bad.events.summed}.
Taking square roots retains the factor~$s^{-2}$ in the~$\Gamma_{\nf12}$ terms of~\eqref{e.complete.wrapping} and~\eqref{e.complete.wrapping.sharp}.
\end{proof}

We conclude this section with the proof of Proposition~\ref{p.induction.bounds}.

\begin{proof}[{Proof of Proposition~\ref{p.induction.bounds}}]
Proposition~\ref{p.propagate.diffusivity.lower.bound} asserts the existence of~$C(d)<\infty$ such that, if
\begin{equation}
	\label{e.lower.bound.cgamma.cond.again}
	E \geq C \cstar^{-1}
	\qand
	\cgamma \leq E^{-10}
	\,,
\end{equation}
then
\begin{equation}
	\mathcal{S}(m_0-1,E)
	\quad\implies \quad
	\text{\eqref{e.shom.h.bounds} is valid for~$m = m_0$}\,.
	\label{e.propagation.of.indyhyp.shombounds.again}
\end{equation}
Proposition~\ref{p.multiscale.estimate} asserts, upon choosing~$\ep \coloneqq  (2C_{\eqref{e.complete.wrapping}})^{-2}$, under the same condition~\eqref{e.lower.bound.cgamma.cond.again} with a possibly larger~$C$, that
\begin{equation}
	\mathcal{S}(m_0-1,E)
	\quad\implies \quad
	\text{\eqref{e.new.induction.for.shom} is valid for~$m = m_0$}\,.
	\label{e.propagation.of.indyhyp.mathcalEbounds}
\end{equation}
Put together, the previous two displays imply that, under the condition~\eqref{e.lower.bound.cgamma.cond.again},
\begin{equation*}
	\mathcal{S}(m_0-1,E)
	\implies
	\mathcal{S}(m_0,E) \,.
\end{equation*}
According to Proposition~\ref{p.base.case},
$\mathcal{S}(\mstarstar ,C_{\mathrm{start}}\cstar^{-\nf12})$ is valid.
After enlarging the dimensional constant, we may fix
$E=C\cstar^{-1}\geq C_{\mathrm{start}}\cstar^{-\nf12}$; the preceding
implication and induction then give~$\mathcal S(m,E)$ for every~$m\in\Z$.
Since~$s\geq8\cgamma$ and~$\cgamma\leq C^{-10}\cstar^{10}$, after enlarging~$C$ once more,
\begin{equation*}
	s^{-2}\exp\bigl(-C^{-3}\cstar^3\cgamma^{-1}\bigr)
	\leq\exp\bigl(-(2C^3)^{-1}\cstar^3\cgamma^{-1}\bigr)
	\,.
\end{equation*}
For any prescribed~$m$, reapply
Proposition~\ref{p.propagate.diffusivity.lower.bound} with
$m_0=\max\{m,\mstarstar+1\}$. Its estimate~\eqref{e.shom.m.flow}, with this
choice of~$E$, gives~\eqref{e.formula.for.shom}, while the multiscale part of
$\mathcal S(m,E)$ and the previous display give~\eqref{e.induction.E.bounds}. This completes the
proof.
\end{proof}

\section{Anomalous regularization}
\label{s.regularity}

In this section we prove Theorem~\ref{t.regularity}, which establishes H\"older estimates for solutions of the equation~$-\nabla \cdot \a \nabla u = \nabla \cdot \f$ with regularity exponent~$\alpha$ for any~$\alpha < 1 - C\cgamma^{\nf12}$. The argument is a multiscale iteration that combines probabilistic and deterministic ingredients. The probabilistic component, developed in Section~\ref{ss.scale.local}, identifies ``good scales'' on which the homogenization error is deterministically small; due to an approximate scale locality property, the proportion of bad scales is shown to be at most~$O(\cgamma^{\nf12})$. In Section~\ref{ss.minimal.scale.sep}, we convert this proportion bound into a minimal scale separation result, identifying a random scale~$Z_m$ below which both the density of bad scales and the cumulative homogenization error are controlled. The deterministic component, presented in Section~\ref{ss.excess.decay}, is an excess decay iteration: on good scales, solutions are well-approximated by harmonic functions, yielding a contraction of the excess, while on bad scales we accept a bounded loss. Section~\ref{ss.proof.regularity} combines these ingredients to complete the proof of Theorem~\ref{t.regularity}.

\smallskip

Finally, in Section~\ref{ss.proof.theorem.B}, we prove Theorem~\ref{t.homogenization} by combining the regularity estimates with the coarse-graining machinery developed in  Sections~\ref{s.coarse.graining.theory} and~\ref{s.RG.flow}.

\smallskip

For later reference, we record the induction-state bookkeeping from Section~\ref{s.RG.flow}. We fix~$E=C\cstar^{-1}$ as in the proof of Proposition~\ref{p.induction.bounds} and assume that~$\cgamma\leq E^{-10}$. The conclusion of that proof gives~$\mathcal S(m,E)$ for every~$m\in\Z$. In particular,~\eqref{e.shom.m.vs.shom.n} is valid for every~$m,n\in\Z$ with~$m>n$.

\subsection{Scale locality of the homogenization error}
\label{ss.scale.local}

In this subsection we introduce an explicit event, denoted by~$\mathcal{G}(m;s,\ep)$, on which the homogenization error can be shown to be deterministically small. We say that~$m$ represents a ``good scale'' if this event is valid and we show that the \emph{proportion} of good scales is large: depending on our choice of parameters, this proportion can be as large as~$1-C \cgamma^{\nf12}$. This relies on the \emph{approximate scale locality} of the events~$\mathcal{G}(m;s,\ep)$; in other words,~$\mathcal{G}(m;s,\ep)$ and~$\mathcal{G}(n;s,\ep)$ are nearly independent if~$|m-n|$ is large. We prove this using the scale decomposition of our fields, the relatively explicit form of the homogenization errors~$\mathcal{E}_{s,q}(\cu_m)$ and their subadditivity property, and the coarse-grained sensitivity estimates in Section~\ref{ss.sensitivity}.

\begin{definition}[Good events]
\label{d.good.event.for.lambda}
For each~$m\in\Z$,~$s\in (0,1]$ and~$\ep,T>0$, we denote
\begin{equation*}
	\left\{
	\begin{aligned}
		&
		\mathcal{G}_0(m)
		\coloneqq
		\biggl\{
		\sup_{k \in \Z \cap (-\infty,m-1]}  3^{-\frac14 \cgamma (m-k)} \max_{z\in 3^{k-2}\Zd \cap (\cu_{m} \setminus \cu_{k})} \bigl( \shom_{k-3} \lambda_{\cgamma,2}^{-1}(z+\cu_{k-2};\a_{k-2}) - C_{\eqref{e.cg.ellip.lower}} \bigr)_+
		\leq
		1
		\biggr\}
		\,,
		\\
		& \mathcal{G}_1(m \, ; s,T)
		\coloneqq
		\biggl\{
		\sum_{k=m}^\infty
		3^{(2-\cgamma)m}
		\|  \nabla \mathbf{j}_k  \|_{\underline{W}^{1,\infty}(\cu_m)}
		\leq
		T
		\biggr\}
		\notag \\  & \qquad \qquad \qquad \qquad
		\cap
		\biggl\{
		\sum_{n=-\infty}^{m}
		3^{-\frac14 s(m-n)}
		\sum_{k=n-1}^m
		\bigl(
		3^{(2-\cgamma)k}
		\max_{z \in 3^k \Zd \cap \cu_m}
		\|  \mathbf{j}_k  \|_{\underline{W}^{2,\infty}(z+\cu_k)}
		\bigr)^2
		\leq
		T^2
		\biggr\}
		\,,
		\\
		& \mathcal{G}_2(m \, ; s,\ep)
		\coloneqq
		\biggl\{
		s
		\sum_{j=-\infty}^{m} \sum_{n=-\infty}^{j-1}
		3^{-\frac14 s(m-n)}
		\max_{z \in 3^n \Zd \cap (\cu_j \setminus \cu_{j-1})}
		\mathcal{E}_{s,2,2}(z+\cu_n;\a_{n-2}
		,\shom_{n-2})^2
		\leq
		\ep^2
		\biggr\}
		\,.
	\end{aligned}
	\right.
\end{equation*}
We also set
\begin{equation}
	\label{e.lambda.good.events}
	\mathcal{G}(m\,;s,\ep)
	\coloneqq
	\mathcal{G}_0(m)
	\cap
	\mathcal{G}_1\bigl(m \, ; s, s \ep \cstar^{\nf 12}\cgamma^{-\nf12} \bigr)
	\cap \mathcal{G}_2(m \, ; s,\ep)
	\,.
\end{equation}
For~$y \in \Rd$, we let~$\mathcal{G}(m,y\,;s,\ep)$ denote the translation of the above quantity by~$y$, that is, the cube~$\cu_m$ is replaced by~$y+\cu_m$.
\end{definition}

\begin{proposition}[Homogenization error estimate on good event]
\label{p.mathcalE.annular.decomp}
There exists a constant~$C(d)<\infty$ such that if~$\cgamma \leq C^{-1} \cstar^{10}$, then, for every~$s \in [8\cgamma,\nf14]$ and~$m\in\Z$, we have
\begin{align}
	\lefteqn{
	\sup_{L \geq m}
	\mathcal{E}_{s,\infty,2}
	(\cu_m;\a_L - (\k_L - \k_m)_{\cu_m},\shom_m)^2
	\indc\bigl\{
	\mathcal{G}_0(m)
	\cap
	\mathcal{G}_1(m \, ;  s,\cstar^{\nf12} \cgamma^{-\nf12})\bigr\}
	} \qquad &
	\notag \\ &
	\leq
	Cs\sum_{j=-\infty}^{m} \sum_{n=-\infty}^{j-1}
	3^{-s(m-n)}
	\max_{z \in 3^n \Zd \cap (\cu_j \setminus \cu_{j-1})}
	\mathcal{E}_{s,2,2}(z+\cu_n;\a_{n-2}
	,\shom_{n-2})^2
	\notag \\ & \qquad
	+
	C s^{-3} \cstar^{-4} \cgamma^2 \left|\log \cgamma\right|^{4}
	+
	C s^{-2} \cstar^{-1} \cgamma
	\Biggl(
	\sum_{k=m}^\infty
	3^{(2-\cgamma)m}
	\|  \nabla \mathbf{j}_k  \|_{\underline{W}^{1,\infty}(\cu_m)}
	\Biggr)^2
	\notag \\ & \qquad
	+
	C s^{-2} \cstar^{-1} \cgamma
	\sum_{n=-\infty}^{m}
	3^{-\frac12 s(m-n)}
	\sum_{k=n-1}^m
	\bigl(
	3^{(2-\cgamma)k}
	\max_{z \in 3^k \Zd \cap \cu_m}
	\|  \mathbf{j}_k  \|_{\underline{W}^{2,\infty}(z+\cu_k)}
	\bigr)^2
	\,.
	\label{e.mathcalE.annular.decomp}
\end{align}
In particular, for every~$L\geq m$, if~$\ep\in (0,\nf12]$ and~$ \cgamma \left|\log \cgamma\right|^{2}  \leq s^{\frac32} \cstar^2 \ep $, then
\begin{equation}
	\label{e.mathcalE.annular.decomp.good.event}
	\mathcal{E}_{s,\infty,2}
	(\cu_m;\a_L- (\k_L - \k_m)_{\cu_m},\shom_m) \indc\bigl\{\mathcal{G}(m \, ; s,\ep)\bigr\}
	\leq
	C \ep
	\,.
\end{equation}
\end{proposition}

\begin{proposition}[Proportion of good scales]
\label{p.independence.between.scales}
There exists~$C(d)<\infty$ such that, for every~$s,\ep,\theta\in (0,\nf12]$ satisfying
\begin{equation}
	\label{e.theta.condition}
	\cgamma \leq C^{-1} \min\{ \cstar^{10}\,, \cstar^{2} s^5 \ep^2\}
	\,,\quad
	s \geq 8\cgamma
	\,,\quad
	\theta \geq C \cstar^{-2}  s^{-6} \ep^{-2} \cgamma
\end{equation}
and every~$m_0\in \Z$ and~$M \in \N_0$, we have the estimate
\begin{equation}
	\label{e.proportion.of.good.scales}
	\P \Biggl[
	\avsum_{m=m_0}^{m_0+M}
	\indc{\bigl\{ \mathcal{G} (m\,; s ,\ep ) \bigr\}}
	\leq
	1- \theta \Biggr]
	\leq
	6\exp\biggl( - \frac{\cstar^2 s^7 \ep^2 \theta}{C \cgamma} (M+1) \biggr)
	\,.
\end{equation}
\end{proposition}

We begin by proving a deterministic, uniform lower bound on the coarse-grained ellipticity constants under the events~$\mathcal{G}_0(m)\cap \mathcal{G}_1(m \, ; 1,\cstar^{\nf12}  \cgamma^{-\nf12})$.

\begin{lemma}
\label{l.localize.lambdas.for.regularity}
There exists a constant~$C(d)<\infty$ such that, if~$\cgamma\leq\nf18$, then,
for every~$\widetilde s\in(\cgamma,\nf12]$,
\begin{multline}
	\sup_{n\in\Z\cap(-\infty,m-2]}
	3^{-\cgamma(m-n)}
	\max_{z\in3^n\Zd\cap\cu_m}
	\shom_n\lambda_{\widetilde s,2}^{-1}(z+\cu_n;\a_n)
	\indc\bigl\{\mathcal G_0(m)\cap
	\mathcal G_1(m;1,\cstar^{\nf12}\cgamma^{-\nf12})\bigr\}
	\\
	\leq
	C\frac{1-3^{-2\widetilde s}}
	{1-3^{-2(\widetilde s-\cgamma)}}\,.
	\label{e.localize.lambdas.for.regularity}
\end{multline}
In particular, at~$\widetilde s=2\cgamma$ the quotient on the right is at
most~$2$.
\end{lemma}
\begin{proof}
We work on the event in the statement and first record the estimate for cubes away from the centered column. By monotonicity in the exponent, the two-scale cube-raising consequence of~\eqref{e.subadda.nosymm}, the definition of~$\mathcal G_0(m)$, Lemma~\ref{l.J.sensitivity.no.conditions}, the event~$\mathcal G_1(m;1,\cstar^{\nf12}\cgamma^{-\nf12})$ and~\eqref{e.shom.m.vs.shom.n}, we have, for every~$k\leq N\leq m-2$ and~$z\in3^k\Zd\cap(\cu_m\setminus\cu_k)$,
\begin{equation}
	\shom_N\lambda_{\widetilde s,2}^{-1}(z+\cu_k;\a_N)
	\leq C3^{\cgamma(m-N)}3^{2\cgamma(N-k)}
	\,.
	\label{e.localize.lambdas.off.center}
\end{equation}
Indeed, the first three inputs give
\begin{equation*}
	\shom_{k-2}\lambda_{\cgamma,2}^{-1}(z+\cu_k;\a_{k-2})\leq C3^{\frac14\cgamma(m-k)}\,,
\end{equation*}
while the sensitivity factor is bounded on~$\mathcal G_1(m;1,\cstar^{\nf12}\cgamma^{-\nf12})$ by~$C3^{\frac{1+\cgamma}{4}(m-k)}$, before it is raised to the power~$2\cgamma/(1-2\cgamma)$. Since
\begin{equation*}
	\frac14+\frac{1+\cgamma}{2(1-2\cgamma)}\leq1
\end{equation*}
for~$\cgamma\leq\nf18$, and~\eqref{e.shom.m.vs.shom.n} gives~$\shom_N\shom_{k-2}^{-1}\leq C3^{\cgamma(N-k)}$, this proves~\eqref{e.localize.lambdas.off.center}.

\smallskip

It remains to control the centered column. Fix~$N\leq m-2$ and set
\begin{equation*}
	X_N\coloneqq\css{2\widetilde s}\sum_{k=-\infty}^N3^{-2\widetilde s(N-k)}\bigl|\s_*^{-1}(\cu_k;\a_N)\bigr|\,.
\end{equation*}
This sum is finite by the uniform ellipticity of~$\a_N$. The matrix subadditivity~\eqref{e.subadda.nosymm} gives, for every~$k\leq N$,
\begin{equation*}
	\bigl|\s_*^{-1}(\cu_k;\a_N)\bigr|\leq3^{-d}\bigl|\s_*^{-1}(\cu_{k-1};\a_N)\bigr|+3^{-d}\sum_{z\in3^{k-1}\Zd\cap(\cu_k \setminus \cu_{k-1})}\bigl|\s_*^{-1}(z+\cu_{k-1};\a_N)\bigr|\,.
\end{equation*}
Multiplying by~$\css{2\widetilde s}3^{-2\widetilde s(N-k)}$ and summing over~$k\leq N$, we obtain
\begin{align*}
	X_N
	&\leq3^{2\widetilde s-d}X_N+\css{2\widetilde s}3^{-d}\sum_{k=-\infty}^N3^{-2\widetilde s(N-k)}\sum_{z\in3^{k-1}\Zd\cap(\cu_k \setminus \cu_{k-1})}\bigl|\s_*^{-1}(z+\cu_{k-1};\a_N)\bigr| \\
	&\leq3^{2\widetilde s-d}X_N+\css{2\widetilde s}\sum_{k=-\infty}^N3^{-2\widetilde s(N-k)}\max_{z\in3^{k-1}\Zd\cap(\cu_k \setminus \cu_{k-1})}\lambda_{\widetilde s,2}^{-1}(z+\cu_{k-1};\a_N)\,.
\end{align*}
Here we used that the inner sum has~$3^d-1$ terms and the one-cube bound~$|\s_*^{-1}(z+\cu_{k-1};\a_N)|\leq\lambda_{\widetilde s,2}^{-1}(z+\cu_{k-1};\a_N)$. Since~$\widetilde s\leq\nf12$ and~$d\geq2$, we have~$3^{2\widetilde s-d}\leq3^{1-d}\leq\nf13$. We may therefore reabsorb the first term on the right and reindex the remaining sum to get
\begin{equation*}
	X_N\leq C\css{2\widetilde s}\sum_{k=-\infty}^{N-1}3^{-2\widetilde s(N-k)}\max_{z\in3^k\Zd\cap(\cu_N\setminus\cu_k)}\lambda_{\widetilde s,2}^{-1}(z+\cu_k;\a_N)\,.
\end{equation*}
By the definition of~$\lambda_{\widetilde s,2}^{-1}$, separating at every depth the centered cube from the other cubes and using the same one-cube bound, we conclude that
\begin{equation*}
	\lambda_{\widetilde s,2}^{-1}(\cu_N;\a_N)\leq C\css{2\widetilde s}\sum_{k=-\infty}^{N-1}3^{-2\widetilde s(N-k)}\max_{z\in3^k\Zd\cap(\cu_N\setminus\cu_k)}\lambda_{\widetilde s,2}^{-1}(z+\cu_k;\a_N)\,.
\end{equation*}
Applying~\eqref{e.localize.lambdas.off.center} to the right side and using~$\widetilde s>\cgamma$ yields
\begin{align*}
	\shom_N\lambda_{\widetilde s,2}^{-1}(\cu_N;\a_N)
	&\leq C\css{2\widetilde s}3^{\cgamma(m-N)}\sum_{k=-\infty}^{N-1}3^{-2(\widetilde s-\cgamma)(N-k)} \leq C\frac{1-3^{-2\widetilde s}}{1-3^{-2(\widetilde s-\cgamma)}}3^{\cgamma(m-N)}\,.
\end{align*}
For~$z\in3^N\Zd\cap\cu_m$ with~$z\neq0$, estimate~\eqref{e.localize.lambdas.off.center} with~$k=N$ gives the same bound, since the quotient on the right is at least~$1$. Taking the maximum over~$z$, and then the supremum over~$N\leq m-2$, proves~\eqref{e.localize.lambdas.for.regularity}. For~$\widetilde s=2\cgamma$, the quotient equals~$1+3^{-2\cgamma}\leq2$.
\end{proof}
We next present the proof of Proposition~\ref{p.mathcalE.annular.decomp}.

\begin{proof}[Proof of Proposition~\ref{p.mathcalE.annular.decomp}]
Denote~$\tilde{\a}_{L,m} \coloneqq \a_L - (\k_L - \k_m)_{\cu_m}$ and set
\begin{equation*}
	\tilde{\mathcal{G}}_m(s)\coloneqq \mathcal{G}_0(m)
	\cap
	\mathcal{G}_1(m \, ;  s,\cstar^{\nf 12} \cgamma^{-\nf12})
	\,.
\end{equation*}
By the definition of~$\mathcal{E}_{s,\infty,2}$ and~\eqref{e.bfJ.general}, we have
\begin{align}
	\mathcal{E}_{s,\infty,2}
	(\cu_m;\tilde{\a}_{L,m},\shom_m)^2
	&
	\leq
	C
	\css{2s}
	\!\!
	\sum_{n=-\infty}^m
	\!\!3^{-2s(m-n)} \max_{z \in 3^n \Zd \cap \cu_m} \max_{|e| = 1} J(z{+}\cu_n,\shom_m^{-\nf12}e,\shom_m^{\nf 12}e\,;\tilde{\a}_{L,m})
	\notag \\ & \qquad
	+
	C \css{2s}
	\!\!
	\sum_{n=-\infty}^m
	\!\!3^{-2s(m-n)} \max_{z \in 3^n \Zd \cap \cu_m} \max_{|e| = 1} J(z{+}\cu_n,\shom_m^{-\nf12}e,\shom_m^{\nf 12}e\,;\tilde{\a}_{L,m}^t )
	\,.
	\label{e.mathcal.E.eight.ess.bound}
\end{align}
The transposed line follows from~$\mathbf j\mapsto-\mathbf j$, and the
$k>m$ tail is absorbed into the fifth term below.
We first show a preliminary version of~\eqref{e.mathcalE.annular.decomp}.
The claim is that
\begin{align}
	\lefteqn{
	\css{2s}
	\!\!
	\sum_{n=-\infty}^m
	\!\!3^{-2s(m-n)} \max_{z \in 3^n \Zd \cap \cu_m} \max_{|e| = 1} J(z{+}\cu_n,\shom_m^{-\nf12}e,\shom_m^{\nf 12}e\,;\tilde{\a}_{L,m})
	\indc_{\tilde{\mathcal{G}}_m(s)}
	} \qquad &
	\notag \\ &
	\leq
	Cs\sum_{j=-\infty}^{m} \sum_{n=-\infty}^{j-1}
	3^{-s(m-n)}
	\max_{z \in 3^n \Zd \cap (\cu_j \setminus \cu_{j-1})}
	\mathcal{E}_{s,2,2}(z+\cu_n;\a_{n-2}
	,\shom_{n-2})^2
	\notag \\ & \qquad
	+
	Cs \sum_{n=-\infty}^{m} 3^{-s(m-n)}
	|\shom_m\shom_{n-2}^{-1}-1|^2
	\notag \\ & \qquad
	+
	C\cstar^{-1} \cgamma s
	\sum_{j=-\infty}^{m} \sum_{n=-\infty}^{j-1}
	3^{-s(m-n)}
	3^{2(2-\cgamma)n}
	\max_{z \in 3^n \Zd \cap (\cu_j \setminus \cu_{j-1})}
	\| \nabla (\k_L{-}\k_{n-2})   \|_{\underline{W}^{1,\infty}(z+\cu_n)}^2
	\notag \\ & \qquad
	+
	C\cstar^{-1} \cgamma s \sum_{j=-\infty}^{m}  \sum_{n=-\infty}^{j-1}
	3^{-s(m-n)}
	3^{-2\cgamma n} \max_{z \in 3^n \Zd \cap (\cu_j \setminus \cu_{j-1})}  \| \k_m{-}\k_{n-2}  \|_{\underline{L}^2(z+\cu_n )}^2
	\notag \\ & \qquad
	+
	C\cstar^{-1} s^{-2} \cgamma
	\bigl(
	3^{(2-\cgamma)m}  \| \nabla (\k_L{-}\k_{m})   \|_{\underline{W}^{1,\infty}(\cu_m)}    \bigr)^2
	\,,
	\label{e.mathcalE.annular.decomp.pre.zero}
\end{align}
with the same bound holding for the second term on the right side of~\eqref{e.mathcal.E.eight.ess.bound} since the right side of~\eqref{e.mathcalE.annular.decomp.pre.zero} does not change if we replace the field with its transpose. We prove this in Steps~1 and~2 and then use this estimate to complete the proof of~\eqref{e.mathcalE.annular.decomp} of the proposition in Step~3.

\smallskip

\emph{Step 1.} Annular decomposition.
We show that there exists a constant~$C(d)<\infty$ such that, for every~$m \in \Z$ and~$s \in (0,\nf14]$, we have
\begin{align}
	\lefteqn{
	\sum_{n=-\infty}^m 3^{-2s(m-n)} \max_{z \in 3^n \Zd \cap \cu_m} \max_{|e| = 1}  J(z+\cu_n,\shom_m^{-\nf12}e,\shom_m^{\nf 12}e\,;\tilde{\a}_{L,m})
	} \quad &
	\notag \\ &
	\leq
	C \sum_{j=-\infty}^{m} \sum_{n=-\infty}^{j-1} 3^{-2s(m-n)}
	\max_{z \in 3^n \Zd \cap (\cu_j \setminus \cu_{j-1})}
	\max_{|e| = 1} J(z+\cu_n,\shom_m^{-\nf12}e,\shom_m^{\nf 12}e\,;\tilde{\a}_{L,m})
	\,.
	\label{e.mathcalE.annular.decomp.pre}
\end{align}
For each~$n \leq m$, we decompose~$3^n\Zd \cap \cu_m$ into annular regions:
\begin{equation*}
	3^n\Zd \cap \cu_m = \{0\} \cup \bigcup_{j=n+1}^{m} \bigl(3^n\Zd \cap (\cu_j \setminus \cu_{j-1})\bigr)\,.
\end{equation*}
Therefore,
\begin{align*}
	\lefteqn{
	\sum_{n=-\infty}^m 3^{-2s(m-n)} \max_{z \in 3^n \Zd \cap \cu_m} J(z+\cu_n)
	} \qquad &
	\\ &
	\leq \sum_{n=-\infty}^m 3^{-2s(m-n)} \sum_{j=n+1}^{m} \max_{z \in 3^n \Zd \cap (\cu_j \setminus \cu_{j-1})} J(z+\cu_n)
	+ \sum_{n=-\infty}^m 3^{-2s(m-n)} J(\cu_n)
	\,.
\end{align*}
For the first term, we interchange the order of summation:
\begin{equation*}
	\sum_{n=-\infty}^m \sum_{j=n+1}^{m}
	= \sum_{j=-\infty}^{m} \sum_{n=-\infty}^{j-1}
	\,.
\end{equation*}
For the second term, we use subadditivity. For any~$n \leq m$, we have
\begin{equation*}
	J(\cu_n)
	\leq
	\avsum_{z' \in 3^{n-1}\Zd \cap \cu_n} J(z' + \cu_{n-1})
	\leq
	\max_{z' \in 3^{n-1}\Zd \cap (\cu_n \setminus \cu_{n-1})} J(z' + \cu_{n-1}) + 3^{-d} J(\cu_{n-1})
	\,.
\end{equation*}
Iterating this, we obtain
\begin{equation*}
	J(\cu_n)
	\leq
	\sum_{k=-\infty}^{n-1} 3^{-d(n-1-k)} \max_{z' \in 3^{k}\Zd \cap (\cu_{k+1} \setminus \cu_{k})} J(z' + \cu_{k})
	\,.
\end{equation*}
Substituting back and interchanging the order of summation yields
\begin{align*}
	\sum_{n=-\infty}^m 3^{-2s(m-n)} J(\cu_n)
	&
	\leq
	\sum_{n=-\infty}^m 3^{-2s(m-n)} \sum_{k=-\infty}^{n-1} 3^{-d(n-1-k)} \max_{z' \in 3^{k}\Zd \cap (\cu_{k+1} \setminus \cu_{k})} J(z' + \cu_{k})
	\\ &
	=
	\sum_{k=-\infty}^{m-1} \max_{z' \in 3^{k}\Zd \cap (\cu_{k+1} \setminus \cu_{k})} J(z' + \cu_{k}) \sum_{n=k+1}^{m} 3^{-2s(m-n) - d(n-1-k)}
	\,.
\end{align*}
Since~$s \leq \nf12$, we have
\begin{equation*}
	\sum_{n=k+1}^{m} 3^{-2s(m-n) - d(n-1-k)}
	=
	3^{-2s(m-k)+d} \sum_{n=k+1}^{m} 3^{(2s-d)(n-k)}
	\leq
	C 3^{-2s(m-k)}
	\,.
\end{equation*}
Relabeling~$k \to n$ and~$k+1 \to j$, the second term becomes
\begin{equation*}
	C \sum_{j=-\infty}^{m} 3^{-2s(m-j+1)} \max_{z' \in 3^{j-1}\Zd \cap (\cu_{j} \setminus \cu_{j-1})} J(z' + \cu_{j-1})
	\,,
\end{equation*}
which is of the form appearing on the right-hand side of~\eqref{e.mathcalE.annular.decomp.pre} with~$n = j-1$. Combining both terms yields~\eqref{e.mathcalE.annular.decomp.pre}.

\smallskip

\emph{Step 2.} Sensitivity estimates. In view of~\eqref{e.mathcalE.annular.decomp.pre}, we localize~$J$ by switching the field~$\tilde{\a}_{L,m}$ to~$\a_{n-2}$ and~$\shom_m$ to~$\shom_{n-2}$. To this end, we apply Lemma~\ref{l.J.sensitivity.no.conditions} with~$\mu \coloneqq \shom_m \shom^{-1}_{n-2}$,~$\s_0 \coloneqq \shom_{n-2}$ and~$\h \coloneqq \k_L - (\k_L - \k_m)_{\cu_m} - \k_{n-2}$ to each cube~$z+\cu_n$. We obtain
\begin{align}
	\lefteqn{
	J(z{+}\cu_n,\shom_m^{-\nf12}e,\shom_m^{\nf 12}e\,;\tilde{\a}_{L,m})
	} \ \ &
	\notag \\ &
	\leq
	C\shom_m^{-1}\shom_{n-2} \bigl( 1 +
	3^{2n} \| \nabla  (\k_L {-} \k_{n-2}) \|_{\underline{W}^{1,\infty}(z+\cu_n)} \lambda_{2\cgamma,2}^{-1}
	(z{+}\cu_n; \a_{n-2} )
	\bigr)^{\! \frac{2s}{1-4\cgamma} }
	\mathcal{E}_{s,2,2}(z{+}\cu_n;\a_{n-2},\shom_{n-2})^2
	\notag \\ & \quad
	+
	C \shom_m^{-1} \shom_{n-2}^2
	\lambda_{2\cgamma,2}^{-1} (z{+}\cu_n \,; \a_{n-2})
	\bigl( 1+
	3^{2n}  \| \nabla (\k_L - \k_{n-2}) \|_{\underline{W}^{1,\infty}(z+\cu_n)} \lambda_{2\cgamma,2}^{-1}
	(z{+}\cu_n; \a_{n-2} )
	\bigr)^{\frac{4\cgamma}{1-4\cgamma}}
	\notag \\ &  \qquad \qquad
	\times
	\bigl(
	\shom_{n-2}^{-2} \| \k_L - (\k_L - \k_m)_{\cu_m} - \k_{n-2} \|_{\underline{L}^2(z+\cu_n )}^2 +
	|\shom_m \shom_{n-2}^{-1}-1|^2 \bigr)
	\notag \\ &  \quad
	+
	C
	\Bigl( \bigl( \shom_m^{-1}   +   \lambda_{2\cgamma,2}^{-1}
	(z{+}\cu_n; \a_{n-2} ) \bigr) 3^{2n}  \| \nabla (\k_L - \k_{n-2}) \|_{\underline{W}^{1,\infty}(z+\cu_n)}    \Bigr)^2
	\,.
	\label{e.ugly.estimate.for.J.pre}
\end{align}
By Lemma~\ref{l.localize.lambdas.for.regularity},~\eqref{e.shom.m.vs.shom.n} and subadditivity, we have
\begin{equation*}
	\shom_{n-2} \lambda_{2\cgamma,2}^{-1}(z+\cu_{n};\a_{n-2}) \indc_{\tilde{\mathcal{G}}_m(s)}
	\leq
	C 3^{\cgamma(m-n)}
	\,.
\end{equation*}
By the definition of the good event~$\mathcal{G}_1(m \, ; s,T)$  in Definition~\ref{d.good.event.for.lambda},
\begin{align*}
	\lefteqn{
	\shom_{n-2}^{-1}
	3^{2n} \sup_{L \geq m} \| \nabla (\k_L{-}\k_{n-2})   \|_{\underline{W}^{1,\infty}(z+\cu_n)} \indc_{\tilde{\mathcal{G}}_m(s)}
	} \quad
	\notag \\ &
	\leq
	\shom_{n-2}^{-1}
	3^{2n} \sum_{k = n-1}^{m-1} \| \nabla \mathbf{j}_k   \|_{\underline{W}^{1,\infty}(z+\cu_n)}
	\indc_{\tilde{\mathcal{G}}_m(s)}
	+
	\shom_{n-2}^{-1}
	3^{2n} \sum_{k = m}^\infty \| \nabla \mathbf{j}_k   \|_{\underline{W}^{1,\infty}(z+\cu_n)}
	\indc_{\tilde{\mathcal{G}}_m(s)}
	\leq
	C 3^{\frac14 s(m-n)}
	\,.
\end{align*}
For the sum over~$k \geq m$, we use the restriction bound~$\| \cdot \|_{\underline{W}^{1,\infty}(z+\cu_n)} \leq 3^{m-n} \| \cdot \|_{\underline{W}^{1,\infty}(\cu_m)}$, the bound~$\shom_{n-2}^{-1} \leq C \cstar^{-\nf 12} \cgamma^{\nf 12} 3^{-\cgamma n}$ and the first condition in~$\mathcal{G}_1$:
\begin{align*}
	\lefteqn{
	\shom_{n-2}^{-1} 3^{2n} \sum_{k=m}^{\infty}
	\| \nabla \mathbf{j}_k \|_{\underline{W}^{1,\infty}(z+\cu_n)} \indc_{\mathcal{G}_1(m\,; s,\cstar^{\nf 12} \cgamma^{-\nf 12})}
	} \qquad &
	\notag \\ &
	\leq
	C \cstar^{-\nf 12} \cgamma^{\nf 12} 3^{-(1-\cgamma)(m-n)}
	\sum_{k=m}^{\infty} 3^{(2-\cgamma)m}\| \nabla \mathbf{j}_k \|_{\underline{W}^{1,\infty}(\cu_m)}
	\indc_{\mathcal{G}_1(m\,; s,\cstar^{\nf 12} \cgamma^{-\nf 12})}
	\leq
	C 3^{-(1-\cgamma)(m-n)}
	\,.
\end{align*}
For the sum over~$k \in [n-1, m-1]$, we use the restriction bound~$\| \cdot \|_{\underline{W}^{1,\infty}(z+\cu_n)} \leq 3^{k-n} \| \cdot \|_{\underline{W}^{1,\infty}(z+\cu_k)}$ and the second condition in~$\mathcal{G}_1$:
\begin{align*}
	\lefteqn{
	\shom_{n-2}^{-1} 3^{2n} \sum_{k=n-1}^{m-1} \| \nabla \mathbf{j}_k \|_{\underline{W}^{1,\infty}(z+\cu_n)}
	\indc_{\mathcal{G}_1(m\,; s,\cstar^{\nf 12} \cgamma^{-\nf 12})}
	} \quad &
	\notag \\ &
	\leq
	C \cstar^{-\nf 12} \cgamma^{\nf 12}
	\sum_{k=n-1}^{m-1} 3^{(1-\cgamma)(n-k)}
	\max_{z \in 3^k \Zd \cap \cu_m}
	\bigl( 3^{(2-\cgamma)k}
	\| \mathbf{j}_k \|_{\underline{W}^{2,\infty}(z+\cu_k)}
	\bigr)
	\indc_{\mathcal{G}_1(m\,; s,\cstar^{\nf 12} \cgamma^{-\nf 12})}
	\leq
	C 3^{\frac14 s(m-n)}
	\,.
\end{align*}
Since~$(1-\cgamma)(m-n) \geq \frac14 s(m-n)$ for~$s \leq 1$ and~$\cgamma \leq \nf 12$, the claim follows.

\smallskip

Furthermore, since~$\shom_{n-2} \geq \frac14 \cstar^{\nf12}\cgamma^{-\nf12} 3^{\cgamma n}$ and~$\shom_m^{-1} \shom_{n-2}\leq C$, we deduce that
\begin{align}
	\lefteqn{
	J(z+\cu_n,\shom_m^{-\nf12}e,\shom_m^{\nf 12}e\,;\tilde{\a}_{L,m})
	\indc_{\tilde{\mathcal{G}}_m(s)}
	} \quad &
	\notag \\ &
	\leq
	C 3^{s(m-n)}
	\mathcal{E}_{s,2,2}(z+\cu_n;\a_{n-2},\shom_{n-2}) ^2
	+
	C 3^{s(m-n)} |\shom_m\shom_{n-2}^{-1}-1|^2
	\notag \\ &  \qquad
	+
	C \cstar^{-1} \cgamma 3^{s(m-n)}
	\bigl(
	3^{-\cgamma n} \| \k_m{-}\k_{n-2}  \|_{\underline{L}^2(z+\cu_n )} \bigr)^2
	\notag \\ &  \qquad
	+
	C \cstar^{-1} \cgamma 3^{s(m-n)}
	\bigl(
	3^{(2-\cgamma)n}  \| \nabla (\k_L{-}\k_{n-2})   \|_{\underline{W}^{1,\infty}(z+\cu_n)}    \bigr)^2
	\notag \\ &  \qquad
	+
	C \cstar^{-1} \cgamma 3^{(s+\cgamma)(m-n)}
	\bigl(
	3^{(2-\cgamma)m}  \| \nabla (\k_L{-}\k_{m})   \|_{\underline{W}^{1,\infty}(\cu_m)}    \bigr)^2
	\,.
	\label{e.ugly.estimate.for.J}
\end{align}
Indeed, by the Poincar\'e inequality,
\begin{align*}
	\shom_m^{-1} \shom_{n-2}^{-1}
	\| \k_L - \k_m - (\k_L-\k_m)_{\cu_m}\|_{\underline{L}^2(z+\cu_n)}^2
	&
	\leq
	C \cstar^{-1} \cgamma 3^{-\cgamma (m+n)} 3^{2m} \| \nabla (\k_L - \k_m )\|_{L^\infty(\cu_m)}^2
	\notag \\ &
	\leq
	C \cstar^{-1} \cgamma  3^{\cgamma(m-n)} \bigl( 3^{(2-\cgamma)m} \| \nabla (\k_L - \k_m )\|_{\underline{W}^{1,\infty}(\cu_m)}\bigr)^2
	\,,
\end{align*}
which gives the last term in~\eqref{e.ugly.estimate.for.J}. Plugging~\eqref{e.ugly.estimate.for.J} into~\eqref{e.mathcalE.annular.decomp.pre} and using~$s \geq 8 \cgamma$ leads to~\eqref{e.mathcalE.annular.decomp.pre.zero}.

\smallskip

\emph{Step 3.} Conclusion. To conclude, we need to estimate the last three terms on the right in~\eqref{e.mathcalE.annular.decomp.pre.zero}.  By the triangle inequality and interchanging the order of summation, we obtain
\begin{align*}\lefteqn{
\sum_{j=-\infty}^{m} \sum_{n=-\infty}^{j-1}
3^{-s(m-n)}
 3^{2(2-\cgamma)n}
 \max_{z \in 3^n \Zd \cap (\cu_j \setminus \cu_{j-1})}
 \| \nabla (\k_m{-}\k_{n-2})   \|_{\underline{W}^{1,\infty}(z+\cu_n)}^2
} \quad &
\notag \\ &
\leq
\sum_{k=-\infty}^m
\sum_{n=-\infty}^{m}
\sum_{j=-\infty}^{m}
\indc_{\{ n \leq j-1, k \geq n-1\}}
(m-n+1)
3^{-s(m-n)}
 3^{2(2-\cgamma)n}
 \max_{z \in 3^n \Zd \cap \cu_m}
 \| \nabla \mathbf{j}_k  \|_{\underline{W}^{1,\infty}(z+\cu_n)}^2
\notag \\ &
\leq
\sum_{k=-\infty}^m
\sum_{n=-\infty}^{k+1}
(m-n+1)^2
3^{-s(m-n)}
3^{2(2-\cgamma)n}
 \max_{z \in 3^n \Zd \cap \cu_m}
 \| \nabla \mathbf{j}_k  \|_{\underline{W}^{1,\infty}(z+\cu_n)}^2
 \notag \\ &
\leq
\sum_{k=-\infty}^m
\bigl(
3^{(2-\cgamma)k}
\max_{z \in 3^k \Zd \cap \cu_m}
\| \nabla \mathbf{j}_k  \|_{\underline{W}^{1,\infty}(z+\cu_k)}
\bigr)^2
\sum_{n=-\infty}^{k+1}
(m-n+1)^2
3^{-s(m-n)-(1-\cgamma)(k-n)}
\notag \\ &
\leq
C s^{-3}
\sum_{k=-\infty}^m
3^{-\frac12 s (m-k)}
\bigl(
3^{(2-\cgamma)k}
\max_{z \in 3^k \Zd \cap \cu_m}
\| \nabla \mathbf{j}_k  \|_{\underline{W}^{1,\infty}(z+\cu_k)}
\bigr)^2
\,.
\end{align*}
Here we used that the cubes~$z+\cu_k$ with~$z \in 3^k\Zd \cap \cu_m$ tile~$\cu_m$, since~$3^{m-k}$ is odd, so that the maximum over~$z$ of the norms on~$z+\cu_k$ is the norm on~$\cu_m$, up to the factor~$2$ from the two terms of the norm, together with
\begin{align*}3^{(2-\cgamma)n}
\max_{z \in 3^n \Zd \cap \cu_m}
\| \nabla \mathbf{j}_k  \|_{\underline{W}^{1,\infty}(z+\cu_n)}
&
\leq
3^{(1-\cgamma)(n-k)}
\bigl( 3^{(1-\cgamma)k}
\| \nabla \mathbf{j}_k  \|_{L^{\infty}(\cu_m)}
+
3^{(2-\cgamma)k}
\| \nabla^2 \mathbf{j}_k  \|_{L^{\infty}(\cu_m)} \bigr)
\notag \\ &
\leq
2 \cdot
3^{(1-\cgamma)(n-k)}
\bigl( 3^{(2-\cgamma)k}
\max_{z \in 3^k \Zd \cap \cu_m}
\| \nabla \mathbf{j}_k  \|_{\underline{W}^{1,\infty}(z+\cu_k)}
\bigr)
  \,.
\end{align*}
Similarly,
\begin{align*}\lefteqn{
\sum_{j=-\infty}^{m}  \sum_{n=-\infty}^{j-1}
3^{-s(m-n)}
3^{-2\cgamma n} \max_{z \in 3^n \Zd \cap (\cu_j \setminus \cu_{j-1})}  \| \k_m{-}\k_{n-2}  \|_{\underline{L}^2(z+\cu_n )}^2
} \qquad &
\notag \\ &
\leq
\sum_{j=-\infty}^{m}
\sum_{n=-\infty}^{j-1}
\sum_{k=n-1}^m
(m-n+1) 3^{-s(m-n)}
3^{-2\cgamma n} \max_{z \in 3^n \Zd \cap (\cu_j \setminus \cu_{j-1})}  \| \mathbf{j}_k  \|_{\underline{L}^2(z+\cu_n)}^2
\notag \\ &
\leq
\sum_{k=-\infty}^m
\bigl( 3^{-\cgamma k}  \| \mathbf{j}_k  \|_{L^\infty(\cu_m)}\bigr)^2
\sum_{n=-\infty}^{m}
\sum_{j=-\infty}^{m}
\indc_{\{ n \leq j-1, k \geq n-1\}}
(m-n+1) 3^{-s(m-n)+2\cgamma (k-n)}
\notag \\ &
\leq
C
\sum_{k=-\infty}^m 3^{-\frac12 s(m-k)}
\bigl( 3^{-\cgamma k}  \| \mathbf{j}_k  \|_{L^\infty(\cu_m)}\bigr)^2
\sum_{n=-\infty}^{k+1}
\sum_{j=n}^{m}
(m-n+1)
3^{-\frac12s(m-n) + 2\cgamma(k-n)}
\notag \\ &
\leq
C (s-4\cgamma)^{-3}
\sum_{k=-\infty}^m 3^{-\frac12 s(m-k)}
\bigl( 3^{-\cgamma k}  \| \mathbf{j}_k  \|_{L^\infty(\cu_m)}\bigr)^2
\,.
\end{align*}
By~\eqref{e.shom.m.vs.shom.n}, we obtain
\begin{align*}\sum_{n=-\infty}^{m} 3^{-s(m-n)}
|\shom_m\shom_{n-2}^{-1}-1|^2
&
\leq
C
\sum_{n=-\infty}^{m} 3^{-s(m-n)}  \min\bigl\{ \cgamma (m-n) + \cstar^{-2} \cgamma \left|\log \cgamma\right|^{2}  , 1 \bigr\}^2  3^{2\cgamma (m-n)}
\notag \\ &
\leq
\frac{C \cgamma^2}{(s-2\cgamma)^3}  + \frac{C \cgamma^2 \left|\log \cgamma\right|^{4}}{\cstar^4(s-2\cgamma)}
+ C s^{-1} 3^{-s \cgamma^{-1}}
\leq
C s^{-3} \cstar^{-4} \cgamma^2 \left|\log \cgamma\right|^{4}
\,.
\end{align*}
Using the above four displays together with~\eqref{e.mathcalE.annular.decomp.pre.zero} yields~\eqref{e.mathcalE.annular.decomp} and finishes the proof.
\end{proof}

\smallskip

The next three lemmas, together with Lemma~\ref{l.localize.lambdas.for.regularity}, give us the tools needed to prove Proposition~\ref{p.independence.between.scales}. Each of the following lemmas is based on an application of a concentration inequality for indicator functions of ``rare'' events, which is proved in Appendix~\ref{s.concentration}: see Proposition~\ref{p.concentration.for.scales}.

\begin{lemma}
\label{l.good.scales.ratio.lambda}
There exists a constant~$C(d)<\infty$ such that if~$\theta \in (0,1]$ is such that
\begin{equation}
	\label{e.good.scales.ratio.lambda.cond}
	\cgamma \leq C^{-1} \cstar^{10}
	\qand
	\theta
	\geq
	\exp \bigl( - C^{-1} \cstar^2 \cgamma^{-1}\bigr)
	\,,
\end{equation}
then, for every~$m_0 \in \Z$ and~$M \in \N_0$, we have the estimate
\begin{equation}
	\P\biggl[
	\avsum_{m=m_0}^{m_0+M} \indc\bigl\{ \mathcal{G}_0(m) \bigr\} \leq 1 - \theta
	\biggr]
	\leq
	\exp \Bigl( - \theta
	\exp \bigl( - C^{-1} \cstar^2 \cgamma^{-1}\bigr) (M+1)
	\Bigr)
	\,.
	\label{e.good.scales.ratio.lambda}
\end{equation}
\end{lemma}

\begin{proof}
We assume~$\cgamma \leq K^{-1} \cstar^{5}$ and~$\theta
\geq
\exp ( - K^{-1} \cstar^2 \cgamma^{-1})$ for suitably large~$K(d)$ to be determined.
Proposition~\ref{p.cg.ellipticity.bounds}, applied with~$s=\cgamma$, yields that, for~$k \in \Z$,
\begin{equation*}\bigl( \shom_{k-3} \lambda_{\cgamma,2}^{-1}(z+\cu_{k-2};\a_{k-2}) -C_{\eqref{e.cg.ellip.lower}} \bigr)_+
\leq
\O_{\Gamma_{\nf 13}} \bigl( \exp \bigl( - C^{-1} \cstar^2 \cgamma^{-1}\bigr)  \bigr)
 \,,
\end{equation*}
which also gives by~\eqref{e.maxy.bound} that, for~$k,n \in \Z$ with~$k \leq n$,
\begin{equation*}\max_{z\in 3^{k-2}\Zd \cap (\cu_{n} \setminus \cu_{n-1})}
\bigl( \shom_{k-3} \lambda_{\cgamma,2}^{-1}(z+\cu_{k-2};\a_{k-2}) -C_{\eqref{e.cg.ellip.lower}} \bigr)_+
\leq
\O_{\Gamma_{\nf 13}} \bigl( C(n-k+1)^{3} \exp \bigl( - C^{-1} \cstar^2 \cgamma^{-1}\bigr)  \bigr) \,.
\end{equation*}
Thus, by~$\cgamma \leq K^{-1} \cstar^{5}$ for a large enough constant~$K(d)$, we get
\begin{multline*}\mathcal{Y}_n \coloneqq
\cgamma^{-4}
\sum_{k=-\infty}^n
3^{-\frac1{4} \cgamma(n-k)} \max_{z\in 3^{k-2}\Zd \cap (\cu_{n} \setminus \cu_{n-1})} \bigl( \shom_{k-3} \lambda_{\cgamma,2}^{-1}(z+\cu_{k-2};\a_{k-2}) -C_{\eqref{e.cg.ellip.lower}} \bigr)_+
\\
\leq
\O_{\Gamma_{\nf 13}} \bigl( \exp \bigl( - C^{-1} \cstar^2 \cgamma^{-1}\bigr)
\bigr) \,.
\end{multline*}
By union bounds and interchanging the order of summation, we then obtain
\begin{align*}
	\lefteqn{
	\sup_{k \in \Z \cap (-\infty,m]}  3^{-\frac14 \cgamma (m-k)} \max_{z\in 3^{k-2}\Zd \cap (\cu_{m} \setminus \cu_{k})} \bigl( \shom_{k-3} \lambda_{\cgamma,2}^{-1}(z+\cu_{k-2};\a_{k-2}) - C_{\eqref{e.cg.ellip.lower}} \bigr)_+
	} \qquad &
	\\ \notag &
	\leq
	\sum_{k=-\infty}^m  3^{-\frac14 \cgamma (m-k)} \sum_{n= k+1}^m \max_{z\in 3^{k-2}\Zd \cap (\cu_{n} \setminus \cu_{n-1})} \bigl( \shom_{k-3} \lambda_{\cgamma,2}^{-1}(z+\cu_{k-2};\a_{k-2}) - C_{\eqref{e.cg.ellip.lower}} \bigr)_+
	\\ \notag &
	\leq
	\sum_{n=-\infty}^m 3^{-\frac14 \cgamma (m-n)}  \sum_{k = -\infty}^{n} 3^{-\frac14 \cgamma (n-k)}
	\max_{z\in 3^{k-2}\Zd \cap (\cu_{n} \setminus \cu_{n-1})} \bigl( \shom_{k-3} \lambda_{\cgamma,2}^{-1}(z+\cu_{k-2};\a_{k-2}) - C_{\eqref{e.cg.ellip.lower}} \bigr)_+
	\\ \notag &
	=
	\cgamma^4
	\sum_{n=-\infty}^m
	3^{-\frac14 \cgamma (m-n)}
	\mathcal{Y}_n
	\,.
\end{align*}
Set~$r(d)\coloneqq\max\left\{2\,,\left\lceil\log_3\left(3+\frac23\sqrt d\right)\right\rceil\right\}$. 
Finite range dependence and annular separation imply that the sequence~$\{\mathcal Y_n\}$
is~$r(d)$-dependent. The preceding estimate also gives almost sure convergence
of the series. We have
\begin{equation*}\E\bigl[\mathcal{Y}_n^p\bigr]^{\nf 1p}
\leq
C p^3 \exp \bigl( - C^{-1} \cstar^2 \cgamma^{-1}\bigr)
\,.
\end{equation*}
Therefore, since~$\cgamma \leq C^{-1} \cstar^{5}$, by taking~$p \coloneqq  \exp \bigl(K^{-1} \cstar^2 \cgamma^{-1}\bigr)$ with~$K(d)$ large enough, we obtain that~$\E\bigl[\mathcal{Y}_n^p\bigr] \leq 1$. Thus, by assuming~$p\geq 1$ and~$\cgamma \geq 4p^{-1}$, Proposition~\ref{p.concentration.for.scales} yields, for~$\theta \geq p^{-1}$,
\begin{equation}
	\label{e.no.bad.scales.applied.for.lambdas}
	\P \Biggl[ \avsum_{k=m_0}^{m_0+M} \indc{\biggl\{ \sum_{n=-\infty}^k
	3^{-\frac14 \cgamma (k-n)}
	\mathcal{Y}_n > C\cgamma^{-1}  \biggr\}} > \tfrac12 \theta \Biggr]
	\leq
	\exp \biggl(- \frac{  \cgamma p \theta}{C} (M+1)  \biggr) \,.
\end{equation}
The event in~\eqref{e.good.scales.ratio.lambda} is that the bad density is at least~$\theta$, which is contained in the event that it exceeds~$\tfrac12\theta$. By taking~$K(d)$ larger, if necessary, we arrive at~\eqref{e.good.scales.ratio.lambda} by assuming~\eqref{e.good.scales.ratio.lambda.cond}.
\end{proof}

\begin{lemma}
\label{l.ratio.of.good.scales.for.k}
There exists~$C(d)<\infty$ such that, for every~$s \in (0,1]$,~$T \in (0,\infty)$ and~$\theta \in (0,1)$ satisfying
\begin{equation}
	T^2 \geq C s^{-4}\theta^{-1}
	\,,
	\label{e.theta.cond.two}
\end{equation}
and every~$m_0 \in \Z$ and~$M \in \N_0$, we have
\begin{equation}
	\P\biggl[  \avsum_{m=m_0}^{m_0+M} \indc\bigl\{ \mathcal{G}_1(m \, ; s,T)\bigr\} \leq 1-\theta
	\biggr]
	\leq
	\exp\bigl( - C^{-1} s^5 T^2 \theta (M+1) \bigr)
	\,.
	\label{e.ratio.of.good.scales.for.k}
\end{equation}
\end{lemma}

\begin{proof}
The event~$\mathcal{G}_1(m \, ; s,T)$ defined in Definition~\ref{d.good.event.for.lambda} is the intersection of two events, which we estimate separately.

\smallskip

\emph{Step 1.} We prove that there exists~$C(d)<\infty$ such that, for every~$T\in [1,\infty)$ and~$\theta \in (0,1)$ satisfying~$T \geq C$ and~$\theta \geq C T^{-2}$, we have
\begin{equation}
	\P \Biggl[ \avsum_{m=m_0}^{m_0+M} \indc{\biggl\{ \sum_{k=m}^\infty
	3^{(2-\cgamma)m}
	\|  \nabla \mathbf{j}_k  \|_{\underline{W}^{1,\infty}(\cu_m)}
	> T \biggr\}} > \frac12 \theta \Biggr]
	\leq
	\exp \biggl(- \frac{ T^2 \theta}{C} (M+1)  \biggr)
	\,.
	\label{e.large.waves.scale.counting}
\end{equation}
We apply Proposition~\ref{p.concentration.for.scales} with
\begin{equation*}
	X_{k} \coloneqq
	K T^{-1}
	3^{(2-\cgamma)k}
	\|  \nabla \mathbf{j}_{k} \|_{\underline{W}^{1,\infty}(\cu_{k})}
	\,,
\end{equation*}
where~$K \in[1,\infty)$ is a large constant selected below. Using that~$\| \cdot \|_{\underline{W}^{1,\infty}(\cu_{m})} \leq C 3^{k-m} \| \cdot \|_{\underline{W}^{1,\infty}(\cu_{k})}$ for~$m\leq k$, we may write
\begin{align}
	\sum_{k=m}^\infty
	3^{(2-\cgamma)m}
	\|  \nabla \mathbf{j}_k  \|_{\underline{W}^{1,\infty}(\cu_m)}
	&
	\leq
	C \sum_{k=m}^\infty
	3^{(2-\cgamma)m}
	3^{k-m}
	\|  \nabla \mathbf{j}_k  \|_{\underline{W}^{1,\infty}(\cu_k)}
	\leq
	C T K^{-1} \! \sum_{k=m}^{\infty} 3^{-(1-\cgamma)(k-m)}
	X_{k}
	\,.
	\label{e.counting.scales.G1.Xk.link}
\end{align}
By~\ref{a.j.reg} and~\eqref{e.moments.OGamma2} we have, for every~$p \in [1,\infty)$ and~$k\in\Z$,
\begin{equation}
	\E\bigl[ X_{k}^p \bigr]^{\nf 1p}
	\leq
	C K T^{-1} p^{\nf12}
	\,.
	\label{e.G1.scale.count.Xk.bound}
\end{equation}
Taking~$p \coloneqq C_{\eqref{e.G1.scale.count.Xk.bound}}^{-2} K^{-2} T^{2}$ and assuming that~$T \geq 2C_{\eqref{e.G1.scale.count.Xk.bound}}K$, we have~$\E[ X_{k}^p ] \leq 1$ and~$p \geq 4$. By further assuming that~$\theta \geq 2p^{-1} = CK^2T^{-2}$ and taking~$K(d)$ large enough, we may apply  Proposition~\ref{p.concentration.for.scales} to obtain
\begin{equation}
	\P \Biggl[ \avsum_{m=m_0}^{m_0+M} \indc{\biggl\{  \sum_{k=m}^\infty 3^{-(1-\cgamma)(k-m)} X_{k}
	> 16 C^{\nf 1p} \theta^{-\nf 1p} \biggr\}} > \frac12 \theta \Biggr]
	\leq
	\exp \biggl(- \frac{p \theta}{32} (M+1)  \biggr) \,.
	\label{e.counting.scales.G1.step1.pre}
\end{equation}
In view of the condition on~$\theta$ and the fact that~$p\geq 4$, by enlarging~$K$, if necessary, we can ensure that~$K \geq  16C_{\eqref{e.G1.scale.count.Xk.bound}}C_{\eqref{e.counting.scales.G1.step1.pre}}^{\nf1p}\theta^{-\nf1p}$. We therefore obtain~\eqref{e.large.waves.scale.counting} by combining~\eqref{e.counting.scales.G1.Xk.link} and~\eqref{e.counting.scales.G1.step1.pre}.

\smallskip

\emph{Step 2.}
We show that there exists a constant~$C(d)<\infty$ such that, for every~$T\in [1,\infty)$ satisfying~$T^2 \geq C s^{-4}$ and~$\theta\in (0,1)$ satisfying~$\theta \geq C s^{-3} T^{-2}$, we have
\begin{multline}
	\P \Biggl[ \avsum_{m=m_0}^{m_0+M} \indc{
	\biggl\{
	\sum_{n=-\infty}^{m}
	3^{-\frac14 s(m-n)}
	\sum_{k=n-1}^m
	\bigl(
	3^{(2-\cgamma)k}
	\max_{z \in 3^k \Zd \cap \cu_m}
	\|  \mathbf{j}_k  \|_{\underline{W}^{2,\infty}(z+\cu_k)}
	\bigr)^2
	>
	T^2
	\biggr\}
	} > \frac12 \theta \Biggr]
	\\
	\leq
	\exp \biggl(- \frac{  s^5 T^2 \theta}{C} (M+1)  \biggr)
	\,.
	\label{e.small.waves.scale.counting}
\end{multline}
Next we define, for every~$m,k\in\Z$,
\begin{equation*}X_{m,k} \coloneqq K s^{-2} T^{-2} 3^{-\frac18 s(m-k)}
\bigl(
3^{(2-\cgamma)k}
\max_{z \in 3^k \Zd \cap \cu_m}
\|  \mathbf{j}_k  \|_{\underline{W}^{2,\infty}(z+\cu_k)}
\bigr)^2
\indc_{\{ k \leq m \}}
\,,
\end{equation*}
where~$K \in [1,\infty)$ is a large constant selected below.
Observe that
\begin{align*}
	\lefteqn{
	\sum_{n=-\infty}^{m}
	3^{-\frac14 s(m-n)}
	\sum_{k=n-1}^m
	\bigl(
	3^{(2-\cgamma)k}
	\max_{z \in 3^k \Zd \cap \cu_m}
	\|  \mathbf{j}_k  \|_{\underline{W}^{2,\infty}(z+\cu_k)}
	\bigr)^2
	}
	\qquad &
	\notag \\ &
	=
	\sum_{k=-\infty}^m
	\sum_{n=-\infty}^{(k+1) \wedge m}
	3^{-\frac14 s(m-n)}
	\bigl(
	3^{(2-\cgamma)k}
	\max_{z \in 3^k \Zd \cap \cu_m}
	\|  \mathbf{j}_k  \|_{\underline{W}^{2,\infty}(z+\cu_k)}
	\bigr)^2
	\notag \\ &
	\leq
	4 s^{-1}
	\sum_{k=-\infty}^m
	3^{-\frac14 s(m-k)}
	\bigl(
	3^{(2-\cgamma)k}
	\max_{z \in 3^k \Zd \cap \cu_m}
	\|  \mathbf{j}_k  \|_{\underline{W}^{2,\infty}(z+\cu_k)}
	\bigr)^2
	\notag \\ &
	=
	4s K^{-1} T^2   \sum_{k=-\infty}^m 3^{-\frac18 s(m-k)} X_{m,k}
\end{align*}
It therefore suffices to show that
\begin{equation}
	\P \Biggl[ \avsum_{m=m_0}^{m_0+M} \indc{
	\biggl\{
	\sum_{k=-\infty}^m 3^{-\frac18 s(m-k)} X_{m,k}
	>
	\frac 14Ks^{-1}
	\biggr\}
	} > \frac12 \theta \Biggr]
	\leq
	\exp \biggl(- \frac{  s^5 T^2 \theta}{C} (M+1)  \biggr)
	\,.
	\label{e.small.waves.scale.counting.Xmk}
\end{equation}
By~\ref{a.j.reg},~\eqref{e.moments.OGamma2} and~\eqref{e.maxy.bound} we have, for every~$p \in [1,\infty)$,
\begin{equation}
	\E\bigl[ X_{m,k}^p \bigr]^{\nf 1p}
	\leq
	Cs^{-3} K T^{-2} p
	\,.
	\label{e.Xmk.moment.bound}
\end{equation}
The extra factor~$s^{-1}$ in~\eqref{e.Xmk.moment.bound} is the
lattice-maximum penalty absorbed by the discount. Take
$p\coloneqq C^{-1}s^3K^{-1}T^2$. If~$T^2\geq Cs^{-4}K$, then
$\E[X_{m,k}^p]\leq1$ and~$p\geq8s^{-1}$. Requiring also
$\theta\geq p^{-1}=Cs^{-3}KT^{-2}$ permits the application of
Proposition~\ref{p.concentration.for.scales}. Its exponent is at least
$c s^4T^2\theta(M+1)$, which we harmlessly weaken to
$c s^5T^2\theta(M+1)$ since~$s\leq1$; this proves
\eqref{e.small.waves.scale.counting.Xmk}.

\smallskip

\emph{Step 3.} The conclusion. The hypothesis~\eqref{e.theta.cond.two} ensures that the conditions needed to apply the results of Steps~1 and~2 are valid, provided that~$C_{\eqref{e.theta.cond.two}}$ is chosen sufficiently large.
The lemma statement is obtained by combining~\eqref{e.large.waves.scale.counting} and~\eqref{e.small.waves.scale.counting}.
\end{proof}

\begin{lemma}
\label{l.ratio.of.good.scales.for.mathcal.E}
There exists a constant~$C(d)< \infty$ such that, for every~$s \in (0,1]$,~$\ep \in (0,1]$ and~$\theta \in (0,1]$ satisfying
\begin{equation}
	\cgamma \leq C^{-1} \min\{ \cstar^{5}\,, \cstar^{2} s^5 \ep^2\}
	\,,\quad
	s \geq 8\cgamma
	\,,\quad
	\theta \geq C \cstar^{-2}  s^{-4} \ep^{-2} \cgamma
	\,,
	\label{e.theta.cond.one}
\end{equation}
we have, for every~$m_0 \in \Z$ and~$M \in \N_0$, the estimate
\begin{equation}
	\P\Biggl[  \avsum_{m=m_0}^{m_0+M} \indc\bigl\{ \mathcal{G}_2(m \, ; s,\ep)\bigr\}  \leq 1- \theta
	\Biggr]
	\leq
	\exp\biggl( - \frac{\cstar^2 s^5 \ep^2 \theta}{C \cgamma} (M+1) \biggr)
	\,.
	\label{e.ratio.of.good.scales.for.mathcal.E}
\end{equation}
\end{lemma}

\begin{proof}
Let~$K_0,K_1,K_2\in [1,\infty)$ be large constants to be selected below. Assume that
\begin{equation}
	\label{e.K.and.p.choices}
	p \coloneqq K_1^{-1} \cstar^2 s^{4} \ep^2 \cgamma^{-1}\,,
	\quad
	\cgamma \leq K_0^{-1} \min\{ \cstar^{5}\,, \cstar^{2} s^5 \ep^2\}
	\qand
	\theta \geq K_1 \cstar^{-2}  s^{-4} \ep^{-2} \cgamma = p^{-1}
	\,.
\end{equation}
We have~$p \geq K_1^{-1} K_0 s^{-1}$ and taking~$K_0 \geq 4 K_1$ ensures~$p \geq 4 s^{-1}$.

\smallskip

\emph{Step 1.} By Proposition~\ref{p.induction.bounds}, subadditivity of~$\mathcal{E}$,~\eqref{e.powerofGammasigma} and~\eqref{e.maxy.bound}, for every~$n, j\in \Z$ with~$n<j$,
\begin{align}
	\label{e.mathcalE.squared.max.bound}
	\lefteqn{
	\max_{z \in 3^n \Zd \cap (\cu_j \setminus \cu_{j-1})}
	\mathcal{E}_{s,2,2}(z+\cu_n;\a_{n-2},\shom_{n-2})^2
	} \qquad &
	\notag \\ &
	\leq
	\O_{\Gamma_{1} }
	\bigl( C(j-n+1) \cstar^{-2}  s^{-2} \cgamma  \bigr) +
	\O_{\Gamma_{\nf14}}
	\bigl(  C (j-n+1)^{4} \exp( - C^{-1}\cstar^{3}  \cgamma^{-1}) \bigr)
	\,.
\end{align}

\smallskip

\emph{Step 2.} For each~$j \in \Z$, define
\begin{equation*}
	X_{j} \coloneqq
	K_2 \ep^{-2}
	\sum_{n=-\infty}^{j-1} 3^{-\frac14 s(j-n)}
	\max_{z \in 3^n \Zd \cap (\cu_j \setminus \cu_{j-1})}
	\mathcal{E}_{s,2,2}(z+\cu_n;\a_{n-2}
	,\shom_{n-2})^2
	\,.
\end{equation*}
The same finite-range argument shows that this sequence is
$r(d)$-dependent (the sharper geometry here gives two-dependence for
$d\leq728$). By~\eqref{e.mathcalE.squared.max.bound} and~\eqref{e.Gamma.sigma.triangle},
\begin{equation}
	\label{e.Xj.Orlicz.bound}
	X_j \leq
	\O_{\Gamma_{1} }
	\bigl( CK_2 \cstar^{-2}  s^{-4} \ep^{-2} \cgamma  \bigr)
	+
	\O_{\Gamma_{\nf14}}
	\bigl(  C K_2 s^{-5} \ep^{-2} \exp( - C^{-1}\cstar^{3}  \cgamma^{-1}) \bigr)
	\,.
\end{equation}

\smallskip

\emph{Step 3.} By~\eqref{e.Xj.Orlicz.bound}, we have
\begin{equation}
	\label{e.Xj.moment.bound}
	\E\bigl[
	X_j^p
	\bigr]^{\nf 1p}
	\leq
	C K_2 \cstar^{-2}  s^{-4} \ep^{-2} \cgamma p
	+
	C K_2 s^{-5} \ep^{-2} \exp( - C^{-1}\cstar^{3}  \cgamma^{-1}) p^{4}
	\,.
\end{equation}
With~\eqref{e.K.and.p.choices}, the first term equals
\begin{equation*}
	C K_2 \cstar^{-2}  s^{-4}\ep^{-2}  \cgamma p
	=
	C K_2 \bigl( \cstar^{-2}  s^{-4} \ep^{-2} \cgamma \bigr)  \bigl( K_1^{-1} \cstar^2 s^{4} \ep^2 \cgamma^{-1} \bigr)
	=
	C K_2 K_1^{-1}
	\leq \nf12
\end{equation*}
for~$K_2 K_1^{-1}$ small. For the second term, using~$p \leq s^4 \cgamma^{-1}$ and~$\cgamma \leq K_0^{-1} \cstar^5$, we get
\begin{equation*}
	C K_2 \ep^{-2} s^{-5} \exp( - C^{-1}\cstar^{3}  \cgamma^{-1}) p^4
	\leq
	C  K_2 K_1^{-4} \cgamma^{-4} \exp( - C^{-1}K_0^{\nf 35} \cgamma^{-\nf 25})
	\leq \nf12
\end{equation*}
for~$K_0$ large enough.  Thus~$\E[ X_j^p] \leq 1$.

\smallskip

\emph{Step 4.} Define~$Y_m \coloneqq \sum_{j=-\infty}^{m} 3^{-\frac14 s (m-j)} X_{j}$. By interchanging summation,
\begin{equation*}
	s \sum_{j=-\infty}^{m} \sum_{n=-\infty}^{j-1}
	3^{-\frac14 s(m-n)}
	\max_{z \in 3^n \Zd \cap (\cu_j \setminus \cu_{j-1})}
	\mathcal{E}_{s,2,2}(z+\cu_n;\a_{n-2},\shom_{n-2})^2
	=
	K_2^{-1} s \ep^2 Y_m
	\,,
\end{equation*}
so~$\mathcal{G}_2(m\,;s,\ep)^c = \{ Y_m > K_2 s^{-1} \}$.

\smallskip

The event in~\eqref{e.ratio.of.good.scales.for.mathcal.E} is that the bad density is at least~$\theta$, which is contained in the event that it exceeds~$\tfrac12\theta$. We apply Proposition~\ref{p.concentration.for.scales} with~$X_{m,j} \coloneqq X_j$,~$\frac14 s$ in place of~$s$, and~$r=r(d)$. Since~$p \geq 4 s^{-1}$ and~$\E[X_j^p] \leq 1$, by~\eqref{e.no.bad.scales.twosided}, we obtain
\begin{equation}
	\label{e.concentration.applied}
	\P \biggl[ \avsum_{m=m_0}^{m_0+M} \indc{\bigl\{ Y_m > 24 s^{-1} C^{\nf 1p} (\tfrac12\theta)^{-\nf 1p} \bigr\}} > \tfrac12 \theta \biggr]
	\leq
	\exp \biggl(- \frac{sp\theta}{128r(d)}(M+1)\biggr)
	\,.
\end{equation}
Here~$p\geq1/\theta$, so
$(\tfrac12\theta)^{-1/p}=\exp(p^{-1}|\log(\theta/2)|)\leq6$.
Consequently~$24 s^{-1} C^{\nf 1p} (\tfrac12\theta)^{-\nf 1p} \leq K_2 s^{-1}$ after taking~$K_2(d)$ large enough. This fixes~$K_2$. After~$K_2$ has been fixed, we may fix~$K_1$ as in Step 3, and then~$K_0$ large enough as in the beginning of the proof and Step 3.

\smallskip

Inserting~\eqref{e.K.and.p.choices} into~\eqref{e.concentration.applied} gives us
\begin{equation*}
	\exp \biggl(- \frac{sp\theta}{128r(d)}(M+1)\biggr)
	=
	\exp \biggl(- \frac{K_1^{-1}\cstar^2s^5\ep^2\theta}{128r(d)\cgamma}(M+1)\biggr)
	\leq
	\exp \biggl(- \frac{  \cstar^{2}  s^{5} \ep^2 \theta}{C \cgamma} (M+1)  \biggr)
	\,,
\end{equation*}
completing the proof.
\end{proof}

\begin{proof}[Proof of Proposition~\ref{p.independence.between.scales}]
By Definition~\ref{d.good.event.for.lambda},
\begin{equation*}
	\indc\{ \mathcal{G}(m\,;s,\ep)^c \}
	\leq
	\indc\{ \mathcal{G}_0(m)^c \}
	+ \indc\{ \mathcal{G}_1(m\,;s, s \ep \cstar^{\nf 12}\cgamma^{-\nf12})^c \}
	+ \indc\{ \mathcal{G}_2(m\,;s,\ep)^c \}\,.
\end{equation*}
Averaging over~$m \in \{m_0, \ldots, m_0+M\}$ and applying Lemmas~\ref{l.good.scales.ratio.lambda},~\ref{l.ratio.of.good.scales.for.k}, and~\ref{l.ratio.of.good.scales.for.mathcal.E} with~$\theta$ replaced by~$\theta/3$, we obtain~\eqref{e.proportion.of.good.scales}, provided the hypotheses of all three lemmas are satisfied.

\smallskip

Observe that the conditions~\eqref{e.theta.condition} with large enough~$C(d)$ imply the hypotheses of Lemma~\ref{l.good.scales.ratio.lambda} (namely~\eqref{e.good.scales.ratio.lambda.cond}), Lemma~\ref{l.ratio.of.good.scales.for.k} (namely~\eqref{e.theta.cond.two} with~$T = s \ep \cstar^{\nf 12}\cgamma^{-\nf12}$), and Lemma~\ref{l.ratio.of.good.scales.for.mathcal.E} (namely~\eqref{e.theta.cond.one}).

\smallskip

For the~$\mathcal G_0$ lane we use
\eqref{e.no.bad.scales.applied.for.lambdas} directly, with rate
$\cgamma p\theta/C$ at
$p=\exp(K^{-1}\cstar^2/\cgamma)$, which is sharper
than~\eqref{e.good.scales.ratio.lambda}. Taking a union bound over
these three estimates gives~\eqref{e.proportion.of.good.scales}.
\end{proof}

\subsection{Estimates on the minimal scale separation for regularity}
\label{ss.minimal.scale.sep}

In this subsection, we convert the proportion-of-good-scales estimate from Proposition~\ref{p.independence.between.scales} into a minimal scale separation result suitable for the excess decay iteration. The main output is the following statement, which identifies a random scale~$Z_m(s,\delta)$ below which two key properties hold simultaneously: the Ces\`aro average of the homogenization error (restricted to good scales) is at most~$\delta$, and the proportion of bad scales is at most~$\delta$. Crucially, these bounds hold uniformly over all cube centers~$z$, which is necessary for the pointwise regularity estimates we seek.

\begin{proposition}[Minimal scale separation]
\label{p.minimal.scale.separation.sec4}
There exists a constant~$C(d)<\infty$ such that, if~$s\in(0,\nf14]$,~$\delta\in(0,\nf12]$, and the parameters satisfy
\begin{equation}
	\label{e.scale.sep.cond}
	\cgamma \leq C^{-10} \min\{ \cstar^{10} \,,  \delta^2 \cstar^2 s^9\}
	\qqand
	s \geq 8 \cgamma
	\,,
\end{equation}
then, for every~$m\in\Z$, there exists a~$\N$-valued random variable~$Z_m(s,\delta)$ satisfying, for every~$N\in\N$,
\begin{equation}
	\P\bigl[N\leq Z_m(s,\delta)\bigr]
	\leq
	C\exp\biggl(-\frac{(N-1)s^9\cstar^2\delta^2}{C\cgamma}\biggr)
	\label{e.scale.sep.bounds}
\end{equation}
and such that, for every~$n\in\Z$ with~$n\leq m$ and~$Z_m(s,\delta)\leq m-n$, we have
\begin{equation}
	\max_{z\in 3^{n-1} \Zd\cap \cu_m}
	\avsum_{k=n}^m
	\sup_{L \geq k}
	\mathcal{E}_{s,\infty,2}
	(z+\cu_k;\a_L {-} (\k_L {-} \k_k)_{z+\cu_k},\shom_k) \,
	\indc_{ \mathcal{G}(k,z ; s,s \delta^{\nf 12})}
	\leq
	\delta
	\label{e.scale.sep.for.mathcal.E}
\end{equation}
and
\begin{equation}
	\max_{z\in 3^{n-1} \Zd\cap \cu_m}
	\avsum_{k=n}^m
	\bigl(1 - \indc_{ \mathcal{G}(k,z ; s,s \delta^{\nf 12})} \bigr)
	\leq
	\delta
	\,.
	\label{e.scale.sep.for.good.events}
\end{equation}
\end{proposition}

The following lemma is the main step in the proof of Proposition~\ref{p.minimal.scale.separation.sec4}.

\begin{lemma}
\label{l.minimal.scale.sep}
There exists a constant~$C(d) < \infty$ such that, if the parameter~$\cgamma$ satisfies~$\cgamma \leq C^{-10} \cstar^{10}$, then, for every~$m,n\in\Z$ with~$n<m$ and~$s \in [8\cgamma,\nf14]$,
\begin{align}
	\label{e.good.scale.kicking}
	\lefteqn{
	\avsum_{k=n}^m
	\,
	\sup_{L \geq k}
	\,
	\mathcal{E}_{s,\infty,2}
	(\cu_k;\a_L - (\k_L - \k_k)_{\cu_k},\shom_k) \,
	\indc_{ \mathcal{G}(k \, ; s,1)}
	} \qquad &
	\notag \\ &
	\leq
	C  \cstar^{-1} s^{-\nf 72} \cgamma^{\nf12}
	+
	\O_{\Gamma_{2} }
	\bigl( C \cstar^{-1} s^{-\nf 92}  \cgamma^{\nf12}
	(m-n)^{-\nf12}
	\bigr)
	\,.
\end{align}
\end{lemma}
\begin{proof}
By Proposition~\ref{p.mathcalE.annular.decomp}, we have that, for every~$k\in\Z$,
\begin{equation*}
	\sup_{L \geq k}
	\mathcal{E}_{s,\infty,2}
	(\cu_k;\a_L {-} (\k_L {-} \k_k)_{\cu_k},\shom_k) \,
	\indc_{ \mathcal{G}(k \, ; s,1)}
	\leq
	C\bigl(
	D_1(k) {+} D_2(k) {+} D_3(k)
	+ \cstar^{-2} s^{-\nf32}  \cgamma \left|\log \cgamma\right|^{2}\bigr)
	\,,
\end{equation*}
where
\begin{equation*}
	\left\{
	\begin{aligned}
		&
		D_1(k)
		\coloneqq
		s^{\nf12} \sum_{j=-\infty}^{k} \sum_{n=-\infty}^{j-1}
		3^{-\frac12 s(k-n)}
		\max_{z \in 3^n \Zd \cap (\cu_j \setminus \cu_{j-1})}
		\mathcal{E}_{s,2,2}(z+\cu_n;\a_{n-2}
		,\shom_{n-2}) \indc_{\mathcal{G}(k \, ; s ,1)}
		\,,
		\notag \\ &
		D_2(k)
		\coloneqq
		\cstar^{-\nf12} s^{-1} \cgamma^{\nf12}
		\sum_{l=k}^\infty
		3^{(2-\cgamma)k}
		\|  \nabla \mathbf{j}_l  \|_{\underline{W}^{1,\infty}(\cu_k)}
		\,,
		\notag \\ &
		D_3(k)
		\coloneqq
		\cstar^{-\nf12} s^{-1} \cgamma^{\nf12}
		\sum_{l=-\infty}^{k}
		3^{-\frac14 s(k-l)}
		\sum_{i=l-1}^k
		3^{(2-\cgamma)i}
		\max_{z \in 3^i \Zd \cap \cu_k}
		\|  \mathbf{j}_i  \|_{\underline{W}^{2,\infty}(z+\cu_i)}
		\,.
	\end{aligned}
	\right.
\end{equation*}
In view of~$\cstar^{-2} s^{-\nf32}  \cgamma \left|\log \cgamma\right|^{2} \leq \cstar^{-1} s^{-\nf72} \cgamma^{\nf12}$, guaranteed by the assumption~$\cgamma \leq C^{-10} \cstar^{10}$ for~$C$ chosen sufficiently large, it suffices to show that
\begin{align}
	\avsum_{k=n}^m
	\bigl( D_1(k) + D_2(k) + D_3(k)\bigr)
	&
	\leq
	C  \cstar^{-1} s^{-\nf72} \cgamma^{\nf12}
	+
	\O_{\Gamma_{2} }
	\bigl( C \cstar^{-1}  s^{-\nf92} \cgamma^{\nf12}
	(m-n)^{-\nf12}
	\bigr)
	\,.
	\label{e.minimal.scale.D.sum}
\end{align}

\smallskip

\emph{Step 1.}
We prove the estimate~\eqref{e.minimal.scale.D.sum} for~$D_1$. From the definition of the good event (see Definition~\ref{d.good.event.for.lambda}), we have
\begin{equation*}
	\sum_{i =-\infty}^{j-1}
	3^{-\frac12 s(j-i)}
	\max_{z \in 3^i \Zd \cap (\cu_j \setminus \cu_{j-1})}
	\mathcal{E}_{s,2,2}(z+\cu_i;\a_{i-2}
	,\shom_{i-2}) \indc_{\mathcal{G}(k \, ; s ,1)}
	\leq
	C s^{-\nf32} 3^{\frac 1{8} s (k-j)}
	\,.
\end{equation*}
According to~\eqref{e.induction.E.bounds},~\eqref{e.Gamma.sigma.triangle},~\eqref{e.maxy.bound} and the above bound, we have
\begin{equation*}
	\sum_{i =-\infty}^{j-1}
	3^{-\frac12 s(j-i)}
	\max_{z \in 3^i \Zd \cap (\cu_j \setminus \cu_{j-1})}
	\mathcal{E}_{s,2,2}(z+\cu_i;\a_{i-2}
	,\shom_{i-2}) \indc_{\mathcal{G}(k \, ; s ,1)}
	\leq
	\X_{j}^{(1)} + \min\bigl\{ \X_{j}^{(2)}\,, C s^{-\nf32} 3^{\frac 1{8} s (k-j)} \bigr\} \,,
\end{equation*}
where, for each~$i,i' \in \{1,2\}$, the sequences are~$r(d)$-dependent
(with the dependence loss absorbed into~$C(d)$) and, for~$\cgamma \leq C^{-1} \cstar^{5}$ and~$s \geq 8\cgamma$,
\begin{equation*}
	\X_{j}^{(1)}
	\leq
	\O_{\Gamma_2}\bigl(C \cstar^{-1} s^{-\nf 52} \cgamma^{\nf 12}\bigr)
	\qqand
	\X_{j}^{(2)}
	\leq
	\O_{\Gamma_{\nf 12}}\bigl( C \exp(-C^{-1} \cstar^3 \cgamma^{-1}) \bigr)
	\,.
\end{equation*}
Thus,
\begin{align*}
	D_1(k)
	&
	\leq
	s^{\nf12}
	\sum_{j=-\infty}^k
	3^{-\frac12 s(k-j)}
	\X_{j}^{(1)}
	+
	s^{\nf12}
	\sum_{j=-\infty}^k
	3^{-\frac12 s(k-j)}
	\bigl(\X_{j}^{(2)} \bigr)^{\nf 14} \bigl( C s^{-\nf32}  3^{\frac 1{8} s (k-j)}  \bigr)^{\nf 34}
	\notag \\ &
	\leq
	s^{\nf12}
	\sum_{j=-\infty}^k
	3^{-\frac12 s(k-j)}
	\X_{j}^{(1)}
	+
	s^{\nf 12}
	\sum_{j=-\infty}^k
	3^{-\frac14 s(k-j)}
	\bigl( C s^{-\nf 92} \X_{j}^{(2)}\bigr)^{\nf 14}
	\,.
\end{align*}
We have shown that
\begin{equation*}
	D_1(k)
	\leq
	s^{\nf12}
	\sum_{j=-\infty}^k
	3^{-\frac14 s(k-j)}
	\X_{j}
	\,,
\end{equation*}
where we define~$\X_j \coloneqq \X_{j}^{(1)}  + \bigl( C   s^{-\nf 92} \X_{j}^{(2)}\bigr)^{\nf 14}$. By~\eqref{e.powerofGammasigma},~$s \geq 8\cgamma$ and~$\cgamma \leq C^{-1} \cstar^5$, we have
\begin{equation*}
	\X_j \leq \O_{\Gamma_2}(C \cstar^{-1}  s^{-\nf 52}\cgamma^{\nf 12})
\end{equation*}
and, in particular,
\begin{equation}
	\E [ \X_j ]
	\leq
	C \cstar^{-1} s^{-\nf 52}  \cgamma^{\nf12}
	\,.
	\label{e.minimal.scale.Xj.expectation}
\end{equation}
Summing over~$k \in [n,m]\cap \Z$ and rearranging the sums, we obtain
\begin{align}
	s^{-\nf12} \sum_{k=n}^m
	D_1(k)
	\leq
	\sum_{k=n}^m
	\sum_{j=-\infty}^k
	3^{-\frac14 s(k-j)}
	\X_{j}
	&
	=
	\sum_{j=-\infty}^m
	\sum_{k= j \vee n}^m
	3^{-\frac14s(k-j)}
	\X_{j}
	\notag \\ &
	=
	\sum_{j=-\infty}^n
	\sum_{k= n}^m
	3^{-\frac14s(k-j)}
	\X_{j}
	+
	\sum_{j=n+1}^m
	\sum_{k= j }^m
	3^{-\frac14s(k-j)}
	\X_{j}
	\notag \\ &
	\leq
	C s^{-1}
	\sum_{j=-\infty}^m
	3^{-\frac14s(n-j)_+}
	\X_j
	\,.
	\label{e.minimal.scale.D1.bound}
\end{align}
The pair~$\bigl(\X_j^{(1)},\X_j^{(2)}\bigr)$ is a measurable function of the shell increments carried by the annulus at scale~$j$, so the sequence of pairs is~$r(d)$-dependent, and so is~$\{\X_j\}$. Splitting the indices into the~$r(d)$ residue classes modulo~$r(d)$, applying Proposition~\ref{p.concentration} within each class, and summing the classes, we obtain
\begin{equation*}
	(m-n)^{-\nf12}
	\sum_{j=n}^m
	\bigl( \X_j - \E [ \X_j ] \bigr)
	\leq
	\O_{\Gamma_{2} }
	\bigl( C \cstar^{-1}  s^{-\nf52}  \cgamma^{\nf12}
	\bigr)
	\,.
\end{equation*}
Using the triangle inequality in~\eqref{e.Gamma.sigma.triangle} again, we also have
\begin{equation*}
	\sum_{j=-\infty}^n
	3^{-\frac14s(n-j)}
	\X_j
	\leq
	\O_{\Gamma_{2} }
	\bigl( C \cstar^{-1} s^{-\nf72} \cgamma^{\nf12}
	\bigr)
	\,.
\end{equation*}
Adding the previous two displays and~\eqref{e.minimal.scale.Xj.expectation} and inserting the result into~\eqref{e.minimal.scale.D1.bound}, we obtain the bound for~$D_1$ in~\eqref{e.minimal.scale.D.sum}.

\smallskip
\emph{Step 2.}
We show that
\begin{equation}
	\avsum_{k=n}^m
	D_2(k)
	\leq
	C \cstar^{-\nf12} s^{-1}  \cgamma^{\nf12}
	+
	\O_{\Gamma_2} \bigl( C \cstar^{-\nf12} s^{-1} \cgamma^{\nf12} (m-n)^{-\nf12} \bigr)
	\,.
	\label{e.minimal.scale.D2.bound}
\end{equation}
By the definition of~$D_2(k)$,
the normalized restriction estimate
\begin{equation*}
	\|\nabla\mathbf j_l\|_{\underline W^{1,\infty}(\cu_k)}
	\leq3^{l-k}
	\|\nabla\mathbf j_l\|_{\underline W^{1,\infty}(\cu_l)}\,,
	\qquad k\leq l\,,
\end{equation*}
gives
\begin{align*}
	\cstar^{\nf12} s \cgamma^{-\nf12}
	\sum_{k=n}^m
	D_2(k)
	&
	=
	\sum_{k=n}^m
	\sum_{l=k}^\infty
	3^{(2-\cgamma)k}
	\|  \nabla \mathbf{j}_l  \|_{\underline{W}^{1,\infty}(\cu_k)}
	\notag \\ &
	=
	\sum_{l=n}^\infty
	\sum_{k=n}^{m\wedge l}
	3^{(2-\cgamma)k}
	\|  \nabla \mathbf{j}_l  \|_{\underline{W}^{1,\infty}(\cu_k)}
	\notag \\ &
	\leq
	C
	\sum_{l=n}^\infty
	3^{(1-\cgamma)((m\wedge l)-l)}
	3^{(2-\cgamma)l}
	\|  \nabla \mathbf{j}_l  \|_{\underline{W}^{1,\infty}(\cu_l)}
	\notag \\ &
	=
	C
	\sum_{l=n}^{m-1}
	3^{(2-\cgamma) l }
	\|  \nabla \mathbf{j}_l  \|_{\underline{W}^{1,\infty}(\cu_l)}
	+
	C
	\sum_{l=m}^\infty
	3^{(1-\cgamma)(m-l) }
	3^{(2-\cgamma)l }
	\|  \nabla \mathbf{j}_l  \|_{\underline{W}^{1,\infty}(\cu_l)}
	\,.
\end{align*}
By~\eqref{e.diff.law.shift} and~\ref{a.j.reg}, we have that
\begin{equation*}
	3^{(2-\cgamma)n}
	\| \nabla \mathbf{j}_n \|_{\underline{W}^{1,\infty}(\cu_n)}
	\leq \O_{\Gamma_2}(C)\,.
\end{equation*}
Thus,
by the independence of the sequence~$\{ \mathbf{j}_n \}_{n\in\Z}$ and Proposition~\ref{p.concentration},
\begin{equation*}
	\sum_{l=n}^{m-1}
	3^{(2-\cgamma) l }
	\|  \nabla \mathbf{j}_l  \|_{\underline{W}^{1,\infty}(\cu_l)}
	\leq
	C(m-n) +
	\O_{\Gamma_2}\bigl(C(m-n)^{\nf12} \bigr)
	\,.
\end{equation*}
By the triangle inequality in~\eqref{e.Gamma.sigma.triangle},
\begin{equation*}
	\sum_{l=m}^\infty
	3^{(1-\cgamma)(m-l) }
	3^{(2-\cgamma)l }
	\|  \nabla \mathbf{j}_l  \|_{\underline{W}^{1,\infty}(\cu_l)}
	\leq
	\O_{\Gamma_2} (C) \,.
\end{equation*}
Combining these yields~\eqref{e.minimal.scale.D2.bound}.

\smallskip

\emph{Step 3.}
We show that
\begin{equation}
	\avsum_{k=n}^m
	D_3(k)
	\leq
	C \cstar^{-\nf12} s^{-\nf72}  \cgamma^{\nf12}
	+
	\O_{\Gamma_2}
	\bigl(
	C \cstar^{-\nf12} s^{-\nf92} \cgamma^{\nf12}(m-n)^{-\nf12}
	\bigr)
	\,.
	\label{e.minimal.scale.D3.bound}
\end{equation}
We start from the definition of~$D_3(k)$ and rearrange the sums to find that
\begin{align*}
	\cstar^{\nf12} s \cgamma^{-\nf12}
	\sum_{k=n}^m
	D_3(k)
	&
	=
	\sum_{k=n}^m
	\sum_{l=-\infty}^{k}
	3^{-\frac14 s(k-l)}
	\sum_{i=l-1}^k
	3^{(2-\cgamma)i}
	\max_{z \in 3^i \Zd \cap \cu_k}
	\|  \mathbf{j}_i  \|_{\underline{W}^{2,\infty}(z+\cu_i)}
	\notag \\ &
	\leq
	\sum_{k=n}^m
	\sum_{i=-\infty}^k
	3^{-\frac14 s(k-i)}
	\sum_{l=-\infty}^{i+1}
	3^{-\frac14 s(i-l)}
	3^{(2-\cgamma)i}
	\max_{z \in 3^i \Zd \cap \cu_k}
	\|  \mathbf{j}_i  \|_{\underline{W}^{2,\infty}(z+\cu_i)}
	\notag \\ &
	\leq
	Cs^{-1}
	\sum_{k=n}^m
	\sum_{i=-\infty}^k
	3^{-\frac14 s(k-i)}
	3^{(2-\cgamma)i}
	\max_{z \in 3^i \Zd \cap \cu_k}
	\|  \mathbf{j}_i  \|_{\underline{W}^{2,\infty}(z+\cu_i)}
	\notag \\ &
	\leq
	Cs^{-1}
	\sum_{i=n}^m
	\sum_{k=i}^m
	3^{-\frac14 s(k-i)}
	3^{(2-\cgamma)i}
	\max_{z \in 3^i \Zd \cap \cu_k}
	\|  \mathbf{j}_i  \|_{\underline{W}^{2,\infty}(z+\cu_i)}
	\notag \\ & \qquad
	+ Cs^{-1}
	\sum_{i=-\infty}^{n-1} \sum_{k=n}^m
	3^{-\frac14 s(k-i)}
	3^{(2-\cgamma)i}
	\max_{z \in 3^i \Zd \cap \cu_k}
	\|  \mathbf{j}_i  \|_{\underline{W}^{2,\infty}(z+\cu_i)}
	\,.
\end{align*}
By~\eqref{e.diff.law.shift} and~\ref{a.j.reg}, we have that
\begin{equation*}
	3^{(2-\cgamma)i}
	\| \mathbf{j}_i \|_{\underline{W}^{2,\infty}(\cu_i)}
	\leq \O_{\Gamma_2}(C)
\end{equation*}
and thus, by~\eqref{e.maxy.bound},
\begin{equation*}
	\max_{z\in 3^i\Zd\cap \cu_k}
	3^{(2-\cgamma)i}
	\| \mathbf{j}_i \|_{\underline{W}^{2,\infty}(z+\cu_i)}
	\leq \O_{\Gamma_2}(C(k+1-i)^{\nf12} )
	\,.
\end{equation*}
Applying~\eqref{e.Gamma.sigma.triangle}, we obtain
\begin{equation*}
	\sum_{k=i}^m
	3^{-\frac14 s(k-i)}
	3^{(2-\cgamma)i}
	\max_{z \in 3^i \Zd \cap \cu_k}
	\|  \mathbf{j}_i  \|_{\underline{W}^{2,\infty}(z+\cu_i)}
	\leq \O_{\Gamma_2} \biggl( C \sum_{k=i}^m (k+1-i)^{\nf12}
	3^{-\frac14 s(k-i)} \biggr)
	=\O_{\Gamma_2}(C s^{-\nf 32})
\end{equation*}
and, similarly,
\begin{equation*}
	\sum_{i=-\infty}^{n-1} \sum_{k=n}^m
	3^{-\frac14 s(k-i)}
	3^{(2-\cgamma)i}
	\max_{z \in 3^i \Zd \cap \cu_k}
	\|  \mathbf{j}_i  \|_{\underline{W}^{2,\infty}(z+\cu_i)}
	\leq
	\O_{\Gamma_2}(C s^{-\nf 52})
	\,.
\end{equation*}
Therefore, by independence, we get~\eqref{e.minimal.scale.D3.bound}.
\end{proof}

\begin{proof}[{Proof of Proposition~\ref{p.minimal.scale.separation.sec4}}]
For each~$m\in\Z$ and~$\delta \in (0,\nf12]$, we define a~$\N$-valued random variable
\begin{multline*}
	Z^{(1)}_m(s,\delta)
	\coloneqq
	1+
	\sup \biggl\{
	l \in \N \,:\,
	\\
	\max_{z\in 3^{m-l-1} \Zd\cap \cu_m}
	\avsum_{k=m-l}^m
	\sup_{L \geq k}
	\mathcal{E}_{s,\infty,2}
	(z+\cu_k;\a_L {-} (\k_L {-} \k_k)_{z+\cu_k},\shom_k) \,
	\indc_{ \mathcal{G}(k,z ; s,s \delta^{\nf 12})}
	>
	\delta
	\biggr\}
	\,.
\end{multline*}
Put~$\mu_{\rm det}=C\cstar^{-1}s^{-\nf72}\cgamma^{\nf12}$ and
$\mu_{\rm fluc}=C\cstar^{-1}s^{-\nf92}\cgamma^{\nf12}$.
Lemma~\ref{l.minimal.scale.sep}, exponential Markov, and the union over
centres give, for~$\delta\geq2\mu_{\rm det}$ and every~$N\in\N$,
\begin{align*}
	\P \bigl[ Z^{(1)}_m(s,\delta) \geq N \bigr]
	&
	\leq
	\sum_{n=N-1}^\infty
	3^{d(n+1)}
	\P
	\Biggl[
	\avsum_{k=m-n}^m
	\sup_{L \geq k}
	\mathcal{E}_{s,\infty,2}
	(\cu_k;\a_L {-} (\k_L {-} \k_k)_{\cu_k},\shom_k) \,
	\indc_{ \mathcal{G}(k \, ; s,s \delta^{\nf 12})}
	>
	\delta
	\Biggr]
	\notag \\ &
	\leq
	\sum_{n=N-1}^\infty
	\exp\biggl( d(\log 3)(n+1) - \frac{\cstar^2 s^9 \delta^2 (n+1) }{C \cgamma }  \biggr)
	\leq
	C\exp\biggl( - \frac{\cstar^{2} s^9\delta^2 }{C \cgamma } (N-1) \biggr)
	\,.
\end{align*}
This proves the claimed shifted exponential tail for~$Z^{(1)}_m$.

\smallskip

We next define the measurable~$(\N\cup\{\infty\})$-valued variable
\begin{equation*}
	Z^{(2)}_m(s,\delta)  \coloneqq
	1+
	\sup \biggl\{
	n \in \N \,:\,
	\max_{z\in 3^{m-n-1} \Zd\cap \cu_m}
	\avsum_{k=m-n}^m
	\bigl(1 - \indc_{ \mathcal{G}(k,z ; s,s \delta^{\nf 12})}\bigr)
	>
	\delta
	\biggr\}
	\,.
\end{equation*}
We apply Proposition~\ref{p.independence.between.scales} with~$\theta = \delta$,~$\ep = s \delta^{\nf 12}$,~$m_0 = m-n$, and~$M = n$. The conditions in~\eqref{e.theta.condition} then read as
\begin{equation*}
	\cgamma \leq C^{-1} \min\{ \cstar^{10}\,, \cstar^{2} s^5 \ep^2\} = C^{-1} \min\{ \cstar^{10}\,, \cstar^{2} s^7 \delta \}
	\,,\quad
	s \geq 8\cgamma
	\,,\quad
	\delta \geq C \cstar^{-2}  s^{-6} \ep^{-2} \cgamma =  C \cstar^{-2}  s^{-8} \delta^{-1} \cgamma\,,
\end{equation*}
and these are satisfied under our hypotheses~\eqref{e.scale.sep.cond} by taking~$C_{\eqref{e.scale.sep.cond}}$ large enough. We thus obtain
\begin{equation*}
	\P \Biggl[
	\avsum_{k=m-n}^m
	\bigl(1 - \indc_{ \mathcal{G}(k ; s,s\delta^{\nf 12})}\bigr)
	>
	\delta
	\Biggr]
	\leq
	\exp\biggl( - \frac{\cstar^2 s^9 \delta^2}{C \cgamma} (n+1) \biggr)
	\,.
\end{equation*}
Taking a union bound over~$z \in 3^{m-n-1}\Zd \cap \cu_m$ and using~$\delta \geq C s^{-\nf92} \cstar^{-1} \cgamma^{\nf12}$,
\begin{equation*}
	\P \Biggl[
	\max_{z\in 3^{m-n-1} \Zd\cap \cu_m}
	\avsum_{k=m-n}^m
	\bigl(1 - \indc_{ \mathcal{G}(k,z ; s,s\delta^{\nf 12})}\bigr)
	>
	\delta
	\Biggr]
	\leq
	\exp\biggl( - \frac{\cstar^2 s^9 \delta^2}{C \cgamma} (n+1) \biggr)
	\,.
\end{equation*}
By a union bound over~$n \geq N-1$ for~$N \in \N$, we get
\begin{equation*}
	\P \bigl[ Z^{(2)}_m(s,\delta) \geq N \bigr]
	\leq
	\sum_{n=N-1}^\infty
	\exp\biggl( - \frac{\cstar^2 s^9 \delta^2}{C \cgamma} (n+1) \biggr)
	\leq
	C\exp\biggl( - \frac{\cstar^2 s^9 \delta^2}{C \cgamma}  (N-1)  \biggr)
	\,.
\end{equation*}
This proves the same shifted exponential tail for~$Z^{(2)}_m$.

\smallskip

To conclude, we define
\begin{equation*}
	Z_m(s,\delta)
	\coloneqq
	\max\{  Z^{(1)}_m(s,\delta) \,,\,  Z^{(2)}_m(s,\delta) \}\,.
\end{equation*}
Taking the union of the two tails establishes~\eqref{e.scale.sep.bounds}.
By construction, we obtain both~\eqref{e.scale.sep.for.mathcal.E}
and~\eqref{e.scale.sep.for.good.events}.
\end{proof}

\subsection{The excess decay iteration}
\label{ss.excess.decay}

In this subsection, we establish the core analytic engine for the regularity theory: an excess decay iteration that converts the scale-local homogenization estimates of Sections~\ref{ss.scale.local} and~\ref{ss.minimal.scale.sep} into H\"older regularity. The excess of a function~$u$ on a cube~$\cu_m$ measures, roughly, how far~$u$ is from being affine at scale~$3^m$. On good scales---those where the homogenization error~$\mathcal{E}$ is small---solutions of~$-\nabla\cdot\a\nabla u = 0$ are well-approximated by harmonic functions, which enjoy classical regularity. This approximation yields a contraction of the excess by a fixed factor when passing from scale~$3^m$ to scale~$3^{m-k}$. On bad scales, we have no such contraction, but the excess can only increase by a volume factor. Since the density of bad scales is at most~$O(\cgamma^{\nf12})$ by Proposition~\ref{p.minimal.scale.separation.sec4}, the iteration closes and yields H\"older regularity with exponent~$1-C\cgamma^{\nf12}$.

\smallskip

Let~$\mathbb{L}$ denote the set of affine functions in~$\Rd$. We define the~\emph{excess} of a function~$u\in L^2(W)$ on a bounded domain~$W\subseteq\Rd$ by
\begin{equation}
	\label{e.excess.def}
	E(u,W)  \coloneqq
	|W|^{-\nf1d}
	\min_{\ell \in \mathbb{L}}
	\| u - \ell \|_{\underline{L}^2(W)}
	\qand
	\ell (u,W)   \coloneqq  \argmin_{\ell \in \mathbb{L}}
	\| u - \ell \|_{\underline{L}^2(W)}
	\,.
\end{equation}
For~$W = \cu_m$, this is
\begin{equation}
	\label{e.excess.def.cubes}
	E(u,\cu_m)  \coloneqq
	3^{-m}
	\min_{\ell \in \mathbb{L}}
	\| u - \ell \|_{\underline{L}^2(\cu_m)}
	\qand
	\ell (u,\cu_m)   \coloneqq  \argmin_{\ell \in \mathbb{L}}
	\| u - \ell \|_{\underline{L}^2(\cu_m)}
	\,.
\end{equation}

In the next lemma, we combine Proposition~\ref{p.general.coarse.graining} with~$p=2$ and the results of the previous subsection (Propositions~\ref{p.mathcalE.annular.decomp} and~\ref{p.independence.between.scales}) to obtain harmonic approximation at every good scale. The following lemma handles both interior cubes and cubes that intersect the boundary of the domain.

\begin{lemma}[Harmonic approximation on good scales]
\label{l.harmonic.approximation.good.scales}
There exists~$C(d)<\infty$ such that the following holds. Let~$s\in[64\cgamma,\nf12]$ and suppose
\begin{equation*}
	\cgamma\leq C^{-1}\cstar^{10}\qquad\mbox{and}\qquad\cgamma|\log\cgamma|^2\leq C^{-1}\cstar^2s^{\nf{11}{2}}\,.
\end{equation*}
Let~$L,m,n\in\Z$ satisfy~$L\geq m$ and~$n+3\leq m$, and let~$x,z\in\cu_m$ satisfy~$x+\cu_n\subset(z+\cu_{n+1})\cap\cu_m$. Suppose that~$u\in H^1(\cu_m)$,~$\g\in H^s(\cu_m;\Rd)$, and~$h\in H^{1+s}(\cu_m)$ satisfy
\begin{equation*}
	\left\{
	\begin{aligned}
		-\nabla\cdot\a_L\nabla u&=\nabla\cdot\g&&\mbox{in }\cu_m\,,\\
		u&=h&&\mbox{on }\partial\cu_m\,.
	\end{aligned}
	\right.
\end{equation*}
Set
\begin{equation*}
	\mathcal{E}\coloneqq\mathcal{E}_{\nf{s}{8},\infty,2}(z+\cu_{n+3};\a_L-(\k_L-\k_{n+3})_{z+\cu_{n+3}} ,\shom_{n+3})\,.
\end{equation*}
There exists a harmonic function~$v\in H^1(x+\cu_n)$ such that~$v=u=h$ on~$\partial(x+\cu_n)\cap\partial\cu_m$ and
\begin{align}
	\lefteqn{
	\|u-v\|_{\underline L^2(x+\cu_n)}
	\indc_{\mathcal G(n+3,z;\nf{s}{8},C^{-1}s^4)}
	} \qquad &
	\notag\\ &
	\leq Cs^{-4}\mathcal{E}
	\|u-(u)_{(z+\cu_{n+3})\cap\cu_m}\|_{\underline L^2((z+\cu_{n+3})\cap\cu_m)}
	+Cs^{-7}\mathcal{E}3^n|(\nabla h)_{(z+\cu_{n+2})\cap\cu_m}|
	\notag\\ &\qquad
	+Cs^{-7}\shom_n^{-1}3^{(1+s)n}[\g]_{\underline H^s((z+\cu_{n+3})\cap\cu_m)}
	+Cs^{-6}3^{(1+s)n}
	[\nabla h]_{\underline H^s((z+\cu_{n+3})\cap\cu_m)}\,.
	\label{e.good.scale.homog}
\end{align}
If, in addition,~$(z+\cu_{n+2})\cap\partial\cu_m=\emptyset$, then we have
\begin{align}
	\lefteqn{
	\|u-v\|_{\underline L^2(x+\cu_n)}
	\indc_{\mathcal G(n+3,z;\nf{s}{8},C^{-1}s^4)}
	} \quad &
	\notag\\ &
	\leq Cs^{-4}\mathcal{E}
	\|u-(u)_{(z+\cu_{n+2})\cap\cu_m}\|_{\underline L^2((z+\cu_{n+2})\cap\cu_m)}
	+Cs^{-7}\shom_n^{-1}3^{(1+s)n}
	[\g]_{\underline H^s((z+\cu_{n+2})\cap\cu_m)}\,.
	\label{e.good.scale.homog.interior}
\end{align}
\end{lemma}
\begin{proof}
Set~$\tilde{\a}\coloneqq\a_L-(\k_L-\k_{n+3})_{z+\cu_{n+3}}$ and, for~$r\in\{2,3\}$, write
\begin{equation*}
	W_r\coloneqq(z+\cu_{n+r})\cap\cu_m\,.
\end{equation*}
If~$x+\cu_{n+1}\subset\cu_m$, set~$Q\coloneqq x+\cu_n$. Otherwise, set, for~$i\in\{1,\ldots,d\}$,
\begin{equation*}
	y_i\coloneqq\max\bigl\{-\tfrac12(3^m-3^{n+1}),\min\{\tfrac12(3^m-3^{n+1}),x_i\}\bigr\}\,,
\end{equation*}
and set~$Q\coloneqq y+\cu_{n+1}$. The coordinatewise construction gives
\begin{equation*}
	x+\cu_n\subset Q\subset W_3\,.
\end{equation*}
Let~$v\in u+H^1_0(Q)$ be harmonic in~$Q$. Its restriction to~$x+\cu_n$ is harmonic and satisfies~$v=u=h$ on~$\partial(x+\cu_n)\cap\partial\cu_m$.

\smallskip

The subtracted matrix is constant and antisymmetric, and therefore~$-\nabla\cdot\tilde{\a}\nabla u=\nabla\cdot\g$. Let~$v_{\g}$ have the same boundary values as~$u$ on~$Q$ and solve~$-\shom_{n+3}\Delta v_{\g}=\nabla\cdot\g$ in~$Q$. Apply the translated~$p=2$ argument from Proposition~\ref{p.general.coarse.graining} in~$Q$ with~$s_1=s/3$, middle exponent~$s/2$, and~$s_2=s$. Cauchy--Schwarz across the scales trades the~$q=1$ quantities at index~$s/3$ for~$q=2$ quantities at index~$s/6$, and Lemma~\ref{l.lambdas.stability}, in particular~\eqref{e.mathcalE.stability}, then bounds every homogenization error used by the argument by~$C\mathcal{E}$. The zero-boundary estimate~\eqref{e.negative.Besov.to.L2.zero.bndr}, with~$\theta=s/2$, converts the resulting~$\underline B^{-\theta}_{2,2}$ seminorm to~$L^2$, and the same estimate applied to~$v_{\g}-v$ gives
\begin{equation}
	\|u-v\|_{\underline L^2(Q)}
	\indc_{\mathcal G(n+3,z;\nf s8,C^{-1}s^4)}
	\leq Cs^{-4}\mathcal{E}\shom_{n+3}^{-\nf12}3^n
	\nu^{\nf12}\|\nabla u\|_{\underline L^2(Q)}
	+Cs^{-7}\shom_{n+3}^{-1}3^{(1+s)n}
	[\g]_{\underline H^s(W_3)}\,.
	\label{e.harmonic.approx.local.general}
\end{equation}

\smallskip

Suppose first that~$Q=y+\cu_{n+1}$. Let~$Q_1,\ldots,Q_{9^d}$ be the triadic grandchildren of~$Q$, each of side length~$3^{n-1}$, and write~$3Q_j$ for the concentric cube of side length~$3^n$. For every~$j$, choose a clamped cube~$\widetilde Q_j$ of side length~$3^{n+1}$ such that
\begin{equation*}
	Q=\bigcup_{j=1}^{9^d}Q_j\,,\qquad Q_j\subset\widetilde Q_j\subset W_3\,,\qquad \partial\widetilde Q_j\cap3Q_j\subset\partial\cu_m\,,
\end{equation*}
where the union is understood up to the boundaries of the cubes. The existence of~$\widetilde Q_j$ follows from the same coordinatewise clamping argument used to define~$Q$. We keep the good event in~\eqref{e.harmonic.approx.local.general} understood throughout this boundary argument. On this event, the Cauchy--Schwarz index trade, Lemma~\ref{l.lambdas.stability}, and Lemma~\ref{l.mathcal.E.to.Lambdas}, including~\eqref{e.mathcalE.stability}, give
\begin{equation*}
	\Lambda_{\nf12}(\widetilde Q_j;\tilde{\a})+\Lambda_{\nf s3}(\widetilde Q_j;\tilde{\a})\leq C\shom_{n+3}\qquad\mbox{and}\qquad\lambda_{\nf s3}^{-1}(\widetilde Q_j;\tilde{\a})\leq C\shom_{n+3}^{-1}\,,
\end{equation*}
and bound every required homogenization error by~$C\mathcal{E}$. Since the divergence is unchanged by constants, in the Caccioppoli estimate on~$\widetilde Q_j$ we replace~$\g$ by~$\g-(\g)_{\widetilde Q_j}$; fractional Poincar\'e and interpolation then bound its full norm by the displayed seminorm of~$\g$. We may therefore apply Lemma~\ref{l.coarse.grained.Caccioppoli.RHS} on each~$\widetilde Q_j$, with~$Q_j$ as the inner patch, parameters~$(\nf12,\nf s3)$, and~$u-(h)_{\widetilde Q_j}$ and~$h-(h)_{\widetilde Q_j}$ in place of~$u$ and~$h$, respectively. Summing the~$9^d$ estimates and using the fixed-scale mean comparisons gives
\begin{align}
	\shom_{n+3}^{-\nf12}3^n\nu^{\nf12}\|\nabla u\|_{\underline L^2(Q)}
	&\leq C
	\|u-(u)_{W_3}\|_{\underline L^2(W_3)}
	+C\bigl|(u-h)_Q\bigr|
	\notag\\
	&\qquad+Cs^{-3}3^{(1+s)n}
	[\nabla h]_{\underline H^s(W_3)}
	+Cs^{-3}3^n
	|(\nabla h)_{W_2}|
	\notag\\
	&\qquad+Cs^{-\nf{11}{2}}\shom_{n+3}^{-1}3^{(1+s)n}
	[\g]_{\underline H^s(W_3)}\,.
	\label{e.harmonic.approx.energy.general}
\end{align}
The clamping places a full face~$F$ of~$Q$ in~$\partial\cu_m$. After rescaling~$Q$ to~$\cu_0$, choose a smooth function~$\eta$ supported in the relative interior of~$F$ with~$\int_F\eta=1$, and a smooth function~$\Phi$ such that~$\Phi=0$ on~$\partial Q$,~$\partial_{\mathbf n}\Phi=\eta$ on~$F$, and~$\partial_{\mathbf n}\Phi=0$ on the other faces. For~$w=v-(v)_Q$, harmonicity and Green's identity give~$\int_F\eta w=\int_Qw\Delta\Phi$. Since~$v=h$ on~$F$, this estimate and the trace--Poincar\'e inequality for~$h$ give
\begin{equation*}
	\bigl|(v-h)_Q\bigr|
	\leq C\|v-(v)_Q\|_{\underline L^2(Q)}
	+C3^n\|\nabla h\|_{\underline L^2(Q)}\,.
\end{equation*}
The triangle inequality and the bounded ratio of the volumes of~$Q$ and~$W_3$ then give
\begin{equation}
	\bigl|(u-h)_Q\bigr|\leq C\|u-v\|_{\underline L^2(Q)}+C\|u-(u)_{W_3}\|_{\underline L^2(W_3)}+C3^n\|\nabla h\|_{\underline L^2(W_3)}\,.
	\label{e.harmonic.approx.scalar.general}
\end{equation}
Insert~\eqref{e.harmonic.approx.scalar.general} into~\eqref{e.harmonic.approx.energy.general}, and then insert the result into~\eqref{e.harmonic.approx.local.general}. After increasing the constant in the definition of the good event, Proposition~\ref{p.mathcalE.annular.decomp} gives~$Cs^{-4}\mathcal{E}\leq\nf12$, so the resulting term~$Cs^{-4}\mathcal{E}\|u-v\|_{\underline L^2(Q)}$ is absorbed into the left side. The fractional Poincar\'e split,~$s^{-4}\mathcal{E}\leq1$, and~$\shom_{n+3}^{-1}\leq4\shom_n^{-1}$ by~\eqref{e.shom.m.vs.shom.n} give~\eqref{e.good.scale.homog}. Finally,~$x+\cu_n\subset Q$ and the ratio of the volumes is bounded by~$3^d$, so the normalized-volume comparison gives the estimate in the form stated.

\smallskip

It remains to consider~$Q=x+\cu_n$. Set~$\widehat Q\coloneqq x+\cu_{n+1}$, and let~$P_1,\ldots,P_{3^d}$ be the triadic children of~$Q$, each of side length~$3^{n-1}$. Then
\begin{equation*}
	Q=\bigcup_{j=1}^{3^d}P_j\,,\qquad P_j\subset3P_j\subset\widehat Q\subset W_2\,,\qquad \partial\widehat Q\cap3P_j=\emptyset\,.
\end{equation*}
On the same good event, the stability estimates give
\begin{equation*}
	\Lambda_{\nf12}(\widehat Q;\tilde{\a})+\Lambda_{\nf s3}(\widehat Q;\tilde{\a})\leq C\shom_{n+3}\qquad\mbox{and}\qquad\lambda_{\nf s3}^{-1}(\widehat Q;\tilde{\a})\leq C\shom_{n+3}^{-1}\,.
\end{equation*}
For each~$j$, apply Lemma~\ref{l.coarse.grained.Caccioppoli.RHS} on~$\widehat Q$, with~$P_j$ as the inner patch and parameters~$(\nf12,\nf s3)$, to~$u-(u)_{\widehat Q}$, using~$\g-(\g)_{\widehat Q}$ as the forcing field. Since~$\partial\widehat Q\cap3P_j=\emptyset$, the boundary condition in the lemma is vacuous and we may take~$h=0$. Fractional Poincar\'e and interpolation bound the full forcing norm by~$[\g]_{\underline H^s(W_2)}$. Summing the estimates and inserting the result into~\eqref{e.harmonic.approx.local.general}, using the good-event cap, and adding the nonnegative boundary-data terms proves~\eqref{e.good.scale.homog} in this case.

\smallskip

If~$(z+\cu_{n+2})\cap\partial\cu_m=\emptyset$, then~$x+\cu_{n+1}\subset z+\cu_{n+2}\subset\cu_m$, so~$v$ is the harmonic replacement in~$x+\cu_n$. On~$\mathcal G(n+3,z;\nf s8,C^{-1}s^4)$, the stability estimates permit us to apply Lemma~\ref{l.coarse.grained.Caccioppoli.RHS} directly on~$z+\cu_{n+2}$, with~$x+\cu_n$ as the inner patch and parameters~$(\nf12,\nf s3)$, to~$u-(u)_{z+\cu_{n+2}}$, using~$\g-(\g)_{z+\cu_{n+2}}$ as the forcing field. Since~$\partial(z+\cu_{n+2})\cap(x+\cu_{n+1})=\emptyset$, the boundary condition is vacuous and we may take~$h=0$. Fractional Poincar\'e and interpolation again reduce the full forcing norm to the displayed seminorm. This gives
\begin{align*}
	\shom_{n+3}^{-\nf12}3^n\nu^{\nf12}
	\|\nabla u\|_{\underline L^2(x+\cu_n)}
	\indc_{\mathcal G(n+3,z;\nf s8,C^{-1}s^4)}
	&\leq C\|u-(u)_{z+\cu_{n+2}}\|_{\underline L^2(z+\cu_{n+2})}
	\notag\\
	&\qquad+Cs^{-\nf{11}{2}}\shom_{n+3}^{-1}3^{(1+s)n}
	[\g]_{\underline H^s(z+\cu_{n+2})}\,.
\end{align*}
Repeating the local comparison at the same slot with threshold~$C^{-1}s^4$, combining the two estimates, using~$s^{-4}\mathcal{E}\leq C$, and applying~$\shom_{n+3}^{-1}\leq4\shom_n^{-1}$ by~\eqref{e.shom.m.vs.shom.n} gives~\eqref{e.good.scale.homog.interior}. Since~$s\leq\nf12$, both forcing contributions are bounded by the displayed~$s^{-7}$ term.
\end{proof}
The next lemma gives the excess decay for solutions at good scales.

\begin{lemma}[Excess decay on good scales]
\label{l.excess.decay.good.scales}
There exist constants~$C(d)<\infty$ and~$k_0(d)\geq4$ such that the following holds. Let
$s\in[64\cgamma,\nf14]$,~$\delta\in(0,\nf12]$,~$k\in\N$ with~$k\geq k_0$,~$L,m,n\in\Z$ with~$L\geq m$ and~$n\leq m-1$, and~$z\in\cu_m$. Assume that
\begin{equation*}
	\cgamma\leq C^{-1}\cstar^{10}
	\qquad\mbox{and}\qquad
	\cgamma|\log\cgamma|^2\leq C^{-1}\cstar^2s^{\nf{5}{2}}\delta^{\nf12}\,,
\end{equation*}
and
\begin{equation*}
	\delta\leq C^{-1}s^6\,.
\end{equation*}
Suppose that~$u\in H^1(\cu_m)$,~$\g\in H^s(\cu_m;\Rd)$, and~$h\in C^{1,\nf12}(\cu_m)$ satisfy
\begin{equation*}
	\left\{
	\begin{aligned}
		& -\nabla \cdot \a_L \nabla u = \nabla \cdot \g  & \mbox{in} & \ \cu_m\,, \\
		& u = h & \mbox{on} & \ \partial \cu_m\,.
	\end{aligned}
	\right.
\end{equation*}
For~$j\in\Z$, set~$U_j\coloneqq(z+\cu_j)\cap\cu_m$, and set
\begin{equation*}
	\mathcal{E}_n\coloneqq
	\mathcal{E}_{\nf s8,\infty,2}
	(z+\cu_n;\a_L-(\k_L-\k_n)_{z+\cu_n},\shom_n)\,.
\end{equation*}
Then
\begin{align}
	E(u,U_{n-k})
	\indc_{\mathcal{G}(n,z \, ; \frac 18 s, \frac 18 s \delta^{\nf 12})}
	&
	\leq
	C \bigl( 3^{-\frac12k} +
	3^{(1+\frac d2) k}
	s^{-3}
	\delta^{\nf 12} \bigr)
	E(u,U_n)
	\notag \\ & \qquad
	+C3^{(1+\frac d2)k}s^{-7}
	\Bigl(
	\mathcal{E}_n\bigl|\nabla\ell(u,U_n)\bigr|
	+\shom_n^{-1}3^{sn}[\g]_{\underline H^s(U_n)}
	\notag \\ & \qquad \qquad
	+\bigl(
	\mathcal{E}_n\bigl|(\nabla h)_{U_{n-1}}\bigr|
	+
	3^{\frac 12 n}
	[\nabla h]_{W^{\nf12,\infty}(U_n)}
	\bigr)
	\indc_{\{z\notin\cu_{m-1}\}}
	\Bigr)
	\,.
	\label{e.excess.decay.one.step}
\end{align}
\end{lemma}
\begin{proof}
Set~$\ep=(s/8)\delta^{1/2}$. After enlarging~$C(d)$, the last two hypotheses give
\begin{equation*}
	\ep\leq C^{-1}s^4
	\qquad\mbox{and}\qquad
	\cgamma|\log\cgamma|^2\leq C^{-1}\cstar^2s^{\nf{11}{2}}\,.
\end{equation*}
Thus Lemma~\ref{l.harmonic.approximation.good.scales} applies, and the annular estimate on the good event in the statement gives
\begin{equation*}
	\mathcal{E}_n\leq Cs\delta^{1/2}\,.
\end{equation*}

Since~$s\leq\nf14$, the uniform embedding~$C^{0,\nf12}\subset H^s$ gives~$h\in H^{1+s}(\cu_m)$. If~$z\in\cu_{m-1}$, then~$z+\cu_n\subset\cu_m$. Apply the interior estimate in Lemma~\ref{l.harmonic.approximation.good.scales} at scale~$n-3$ with~$x=z$; its event and homogenization error are both at scale~$n$. Harmonic decay, followed by the triangle inequality and normalized-volume comparison, gives the claimed estimate without the terms carrying the boundary indicator.

Suppose that~$z\notin\cu_{m-1}$. Choose the coordinatewise-clamped center~$y$ so that
\begin{equation*}
	U_{n-3}\subset y+\cu_{n-3}\subset U_{n-2}\,,
\end{equation*}
and apply the general estimate in Lemma~\ref{l.harmonic.approximation.good.scales} at scale~$n-3$, with~$x=y$, to obtain~$v$; its event and homogenization error are both at scale~$n$, while~$v$ is harmonic in~$y+\cu_{n-3}$ and satisfies~$v=h$ on~$\partial(y+\cu_{n-3})\cap\partial\cu_m$. If~$\overline{y+\cu_{n-3}}\cap\partial\cu_m=\emptyset$, apply the ordinary interior harmonic-decay estimate to~$v$ and add the nonnegative boundary terms. Otherwise, let~$\ell_h$ be the affine lift of~$h$ based at~$z$ with slope~$(\nabla h)_{U_{n-3}}$, and let~$v_1$ be harmonic on~$U_{n-3}$ with boundary datum~$h-\ell_h$. Then~$V=v-v_1-\ell_h$ has zero trace on every met face. The met coordinate reflections commute, and the iterated odd extension of~$V$ is harmonic on a reflected box containing a fixed neighborhood of~$U_{n-k}$. Boundary rigidity controls the odd-affine defect by~$E(V,U_{n-3})$, so the reflected excess is at most~$CE(V,U_{n-3})$, and the interior harmonic estimate gives the required decay, with the two displayed boundary-data terms. The estimate for~$v_1$ follows from the maximum principle and~$[\nabla h]_{C^{0,1/2}}$.

The usual Campanato triangle estimate costs~$C3^{(1+d/2)k}$. Insert the two harmonic-approximation bounds, use~$\mathcal{E}_n\leq Cs\delta^{1/2}$ in the excess term, compare the averages of~$\nabla h$ on~$U_{n-3}$ and~$U_{n-1}$, and invoke~\eqref{e.shom.m.vs.shom.n} to replace the finite shift~$\shom_{n-3}$ by~$\shom_n$. The same uniform embedding prices the fractional datum term. Since~$s\leq\nf14$, the remaining powers~$s^{-4}$,~$s^{-6}$, and~$s^{-7}$ are all bounded by the common factor~$s^{-7}$ displayed in the statement. Absorbing the one-scale offsets proves~\eqref{e.excess.decay.one.step}.
\end{proof}

We next record a discrete Gronwall-type inequality which is needed below in the proof of Lemma~\ref{l.iteration.lemma}.

\begin{lemma}[Discrete Gronwall]
\label{l.discrete.gronwall}
Suppose that~$H \geq 0$,~$\{\tilde{\ep}_k\}_{k=1}^{i} \subset \R_+$ and~$\{a_k\}_{k=1}^{i}$ satisfy
\begin{equation*}
	a_l \leq H + \sum_{k=1}^{l-1} \tilde{\ep}_k a_k\,, \qquad \forall l \in \{1,\ldots,i\}\,.
\end{equation*}
Then
\begin{equation*}
	a_i \leq H \exp\Bigl(\sum_{k=1}^{i-1} \tilde{\ep}_k\biggr)
	\,.
\end{equation*}
\end{lemma}
\begin{proof}
Define~$S_l \coloneqq \sum_{k=1}^{l} \tilde{\ep}_k a_k$. The hypothesis gives~$\tilde{\ep}_l a_l \leq \tilde{\ep}_l(H + S_{l-1})$, and thus
\begin{equation*}
	S_l + H
	\leq
	S_{l-1}(1 + \tilde{\ep}_l) + H\tilde{\ep}_l + H
	=
	(1 + \tilde{\ep}_l)(S_{l-1}+H)\,.
\end{equation*}
Iterating this yields~$S_{i-1} + H \leq H \prod_{k=1}^{i-1}(1+\tilde{\ep}_k)$, and thus
\begin{equation*}
	a_i \leq H + S_{i-1} \leq H\prod_{k=1}^{i-1}(1+\tilde{\ep}_k) \leq H\exp\biggl(\sum_{k=1}^{i-1}\tilde{\ep}_k\biggr)
	\,. \qedhere
\end{equation*}
\end{proof}

We next formulate the main result of this section, which is an excess decay estimate formalized with Lemma~\ref{l.excess.decay.good.scales} in mind (compare the right sides of~\eqref{e.excess.decay.one.step} and~\eqref{e.Ej.decay.assumption}). We will combine these lemmas in the next subsection, in the proof of Theorem~\ref{t.regularity}.

\begin{lemma}
\label{l.iteration.lemma}
Let~$h \in \N$ and~$\theta \in (0,1)$ with~$\theta^h \in (0, \nf 35)$.
Let~$n,m \in \Z$ with~$n < m$. Let~$\{ U_j \}_{j \in \Z}$ be bounded measurable sets such that, for every~$j \leq m$,~$U_{j-1} \subset U_j$ and there exist~$x_j,y_j \in \Rd$ such that~$x_j + \cu_{j-2} \subset U_j \subset y_j + \cu_j$. Let~$u \in L^2(U_m)$.
Assume that~$\mathcal{B} \subseteq \Z \cap [n,m]$ and~$\{\ep_j\}_{j \in \Z \cap [n,m]}$,~$\{\delta_j\}_{j \in \Z \cap [n,m]}$ are sequences of non-negative numbers  such that
\begin{equation}
	\label{e.Ej.decay.assumption}
	E(u,U_{j-h})
	\leq \theta^h E(u,U_{j}) + \ep_j |\nabla \ell(u,U_j)| + \delta_j\,, \qquad \forall j \in \Z \cap [n,m] \setminus \mathcal{B}\,.
\end{equation}
Then there exists a constant~$C(d)<\infty$ such that, with
\begin{equation*}
	A \coloneqq C(h+1)(|\mathcal B|+1)+C\sum_{j=n}^m\ep_j\,,
\end{equation*}
we have
\begin{equation}
	\label{e.iteration.slope.lemma}
	3^{-n} \| u - (u)_{U_n} \|_{\underline{L}^2(U_n)}
	\leq
	e^A \biggl(3^{-m} \| u - (u)_{U_m} \|_{\underline{L}^2(U_{m})} + \sum_{j=n}^m \delta_j\biggr)\,.
\end{equation}
Moreover, for every~$M\geq0$ such that~$\ep_j\leq M$ for every~$j\in[n,m]\cap\Z$,
\begin{equation}
	\label{e.excess.decay.lemma}
	E(u,U_n)
	\leq
	\theta^{-C(h+1)(|\mathcal{B}|+1)}e^A
	\biggl(
	\theta^{m-n} E(u,U_m) 
	+M3^{-m} \| u - (u)_{U_m} \|_{\underline{L}^2(U_{m})}
	+(1+M)\sum_{j=n}^m \delta_j
	\biggr)\,.
\end{equation}
\end{lemma}

\begin{proof}
We abbreviate~$E_j \coloneqq E(u,U_j)$ and~$\ell_j \coloneqq \ell(u,U_j)$. Define~$\mathcal{G} \coloneqq [n,m] \cap \Z \setminus \mathcal{B}$.

\smallskip

\emph{Step 1: Gradient stability.}
The sandwich assumption and the inverse affine estimate, followed by the
triangle inequality, give, for every~$k\leq m$,
\begin{equation}
	\label{e.grad.stability}
	|\nabla \ell_k - \nabla \ell_{k-1}| \leq C 3^{-k} \|\ell_k - \ell_{k-1}\|_{\underline{L}^2(U_{k-1})} \leq C(E_k + E_{k-1})\,.
\end{equation}

\smallskip

\emph{Step 2: Decay on good runs.}
Set~$\kappa_d\coloneqq3^{3+3d/2}$. Partition each maximal interval of
good indices into its~$h$ residue classes. Summing
\eqref{e.Ej.decay.assumption} along a residue class and using
$\theta^h<3/5$ gives, for every good run~$[a,b]$ and~$i\in[a,b]$,
\begin{align}
	\label{e.sum.excess.bound}
	& \sum_{k=i}^{b}E_k
	\leq C\kappa_d^h E_b
	+C\sum_{k=i}^{b}\ep_k|\nabla\ell_k|
	+C\sum_{k=i}^{b}\delta_k\,,\\ & 
	\label{e.grad.bound.good}
	E_i+|\nabla\ell_i|
	\leq
	\exp\biggl(C(h+1)+C\sum_{k=i}^{b}\ep_k\biggr)
	\biggl(E_b+|\nabla\ell_b|+\sum_{k=i}^{b}\delta_k\biggr)\,.
\end{align}
Indeed, the first estimate uses the scale sandwich only to compare the at most
$h$ terminal excesses, with total cost~$\kappa_d^h$.  Combining it with
\eqref{e.grad.stability} gives a triangular inequality for
$|\nabla\ell_i-\nabla\ell_b|$, to which
Lemma~\ref{l.discrete.gronwall} applies and yields the second estimate.

\smallskip

\emph{Step 3: Bad scales and iteration.}
For any~$k$, the estimate~\eqref{e.grad.stability} and the definition of excess imply~$E_{k-1} + |\nabla \ell_{k-1}| \leq C(E_k + |\nabla \ell_k|)$. There are at most~$|\mathcal B|+1$ good runs, and each run can meet at most~$h$ incomplete residue classes. Iterating~\eqref{e.grad.bound.good} through the good runs and across the bad indices therefore gives
\begin{equation}
	\label{e.combined.bound}
	E_n + |\nabla \ell_n| \leq e^A \biggl(E_m + |\nabla \ell_m| + \sum_{j=n}^{m} \delta_j\biggr)\,.
\end{equation}
Since~$E_m + |\nabla \ell_m| \leq C  3^{-m} \|u - (u)_{U_m}\|_{\underline{L}^2(U_m)}$ and~$3^{-n}\|u - (u)_{U_n}\|_{\underline{L}^2(U_n)} \leq C(E_n + |\nabla \ell_n|)$, this proves~\eqref{e.iteration.slope.lemma}.

\smallskip

\emph{Step 4: Excess decay.}
Retaining the contraction instead of summing it away gives, on a good run~$[a,b]$,
\begin{equation*}
	E_i \leq
	\theta^{-C(h+1)}\theta^{b-i}E_b
	+C M\max_{k\in[i,b]}|\nabla\ell_k|
	+C\sum_{k=i}^{b}\delta_k\,.
\end{equation*}
The total number of bad crossings and incomplete residue classes is bounded by
$C(h+1)(|\mathcal B|+1)$.  Iterating the preceding display and using
\eqref{e.combined.bound} to control the gradients gives precisely
\eqref{e.excess.decay.lemma}.
\end{proof}

\subsection{H\"older regularity: the proof of Theorem~\ref{t.regularity}}
\label{ss.proof.regularity}

We now combine the results of the previous subsections to give the proof of Theorem~\ref{t.regularity}. The argument proceeds by a multiscale iteration. Starting from a solution~$u$ in~$\cu_m$, we zoom in toward an arbitrary point~$x \in \cu_{m-1}$ by considering a sequence of cubes~$z_j + \cu_j$ containing~$x$ at each scale~$j \leq m-1$. On good scales (where the event~$\mathcal{G}(j)$ holds), Lemma~\ref{l.excess.decay.good.scales} provides excess decay with an error proportional to the homogenization error~$\mathcal{E}$. On bad scales, we simply use stability estimates and accept a bounded loss. The iteration lemma (Lemma~\ref{l.iteration.lemma}) combines these into an oscillation bound, with the exponential prefactor controlled by the number of bad scales and the cumulative homogenization error. Proposition~\ref{p.minimal.scale.separation.sec4} ensures that, below the random minimal scale~$Z_m$, both the proportion of bad scales and the average homogenization error are bounded by a small parameter~$\delta \sim (1-\alpha)$. This yields an oscillation decay of~$3^{(1-\alpha)(m-n)}$. Finally, to convert oscillation to gradient control, we apply the coarse-grained Caccioppoli inequality at the nearest good scale, paying a one-time volume factor that can be absorbed by a slight adjustment of the H\"older exponent.

\begin{proof}[{Proof of Theorem~\ref{t.regularity}}]
We split the argument into several steps.

\smallskip

\emph{Step 1.} Parameter selection. Let~$\alpha \in (0,1)$ be given. We select the parameters~$k$,~$s$ and~$\delta$ in terms of~$\alpha$,~$\cgamma$, and~$\cstar$, as follows:
\begin{equation}
	k \coloneqq \max\left\{k_0,\bigl\lceil4\log_3
	\bigl(2C_{\eqref{e.excess.decay.one.step}}\bigr)\bigr\rceil\right\}\,,
	\qquad
	s \coloneqq \nf14
	\qqand
	\delta \coloneqq C_1^{-1}(1-\alpha) \,,
	\label{e.parameter.choices.regularity}
\end{equation}
where~$k_0$ is given by Lemma~\ref{l.excess.decay.good.scales} and~$C_1(d,\cstar)$ is a constant to be chosen below, sufficiently large depending only on~$d$ and~$\cstar$. Note that~$k\in\N$ is a constant depending only on~$d$.

\smallskip

We verify the conditions required by Proposition~\ref{p.minimal.scale.separation.sec4} and to guarantee the applicability of Lemma~\ref{l.excess.decay.good.scales} and~\eqref{e.mathcalE.annular.decomp.good.event}:
\begin{itemize}

\item With the minimal-scale parameter equal to~$s/8$, its hypotheses are
\begin{equation*}
	\cgamma\leq C^{-10}\min\{\cstar^{10}\,,
	\delta^2\cstar^2(s/8)^9\}
	\qquad\hbox{and}\qquad s/8\geq8\cgamma\,.
\end{equation*}
They follow from~$\cgamma\leq\cgamma_0(d,\cstar)$ after decreasing
$\cgamma_0$.  The harmonic and annular inputs consume the additional
condition
\begin{equation*}
	\cgamma|\log\cgamma|^2
	\leq(s/8)^{3/2}\cstar^2(s/8)\delta^{1/2}\,,
\end{equation*}
which follows from the same choice of~$\cgamma_0$.

\item The condition~$\delta \geq C s^{-\nf92} \cstar^{-1} \cgamma^{\nf12}$ in Proposition~\ref{p.minimal.scale.separation.sec4} requires~$C_1^{-1}(1-\alpha) \geq C \cstar^{-1} \cgamma^{\nf12}$ after the fixed choice~$s=1/4$. This is equivalent, up to a dimensional constant, to~$1-\alpha \geq C C_1 \cstar^{-1} \cgamma^{\nf12}$, precisely the hypothesis~$\alpha \leq 1 - C(d,\cstar) \cgamma^{\nf12}$ stated in the theorem.

\item We require that~$\delta\leq\delta_0$, where~$\delta_0(d)>0$ is chosen so that
\begin{equation*}
	\delta_0\leq C_{\eqref{e.excess.decay.one.step}}^{-1}s^6
	\qquad\hbox{and}\qquad
	C_{\eqref{e.excess.decay.one.step}}
	3^{(1+d/2)k}s^{-3}\delta_0^{\nf12}
	\leq\frac12 3^{-k/4}\,.
\end{equation*}
The integer~$k=k(d)$ has already been fixed, so this is guaranteed by
enlarging~$C_1$ by a factor depending only on~$d$.

\end{itemize}
For every~$L\geq m$, we denote
\begin{equation*}
	\ep_j(z)
	\coloneqq
	\mathcal{E}_{\nf s8, \infty, 2}(z + \cu_j; \a_L - (\k_L-\k_j)_{z +\cu_j}, \shom_j)\indc_{\mathcal{G}(j, z\,; \frac18 s, \frac18 s \delta^{\nf 12})}
\end{equation*}
and notice that, by~\eqref{e.mathcalE.annular.decomp.good.event},
\begin{equation*}
	\ep_j(z)\leq
	C_{\eqref{e.mathcalE.annular.decomp.good.event}}s\delta^{\nf12}\,.
\end{equation*}

\smallskip

\emph{Step 2.} Definition of the minimal scale.
For each~$m \in \Z$, we define the random variable
\begin{equation}
	\label{e.def.X.m.alpha}
	X_m(\alpha) \coloneqq Z_m\bigl( \tfrac 18 s,  \delta \bigr) + k + 3\,,
\end{equation}
where~$Z_m(s,\delta)$ is given by Proposition~\ref{p.minimal.scale.separation.sec4}. The additive constant~$k+3$ accounts for the scale loss in the excess decay step and the interior regularity argument. By~\eqref{e.scale.sep.bounds} and our parameter choices, we have
\begin{equation}
	\label{e.X.m.stochastic.integrability}
	X_m(\alpha)-(k+3)
	=
	\O_{\Gamma_1} \bigl( C\cstar^{-2} (1-\alpha)^{-2} \cgamma  \bigr)\,.
\end{equation}
Equivalently, after changing~$C(d,\cstar)$ to absorb~$k+3$,
\begin{equation*}
	\P[X_m(\alpha)\geq N]
	\leq C\exp\left(-\frac{(1-\alpha)^2(N-C)}{C\cgamma}\right)\,,
\end{equation*}
which is~\eqref{e.Xm.alpha.integrability}.

\smallskip

\emph{Step 3.} Preparation for the oscillation estimate.
Fix~$m \in \Z$,~$L \geq m$, and let~$u \in H^1(\cu_m)$ be a solution of~$-\nabla \cdot \a_L \nabla u = \nabla \cdot \g$ in~$\cu_m$ with~$u = h$ on~$\partial \cu_m$, where~$\g \in C^{0,\nf12}(\cu_m;\Rd)$ and~$h \in C^{1,\nf12}(\cu_m)$.
Fix~$n \in \Z$ with~$n \leq m - X_m(\alpha)$.
Our goal is to establish an oscillation bound for~$u$ in~$(z + \cu_{n}) \cap \cu_m$ with~$z \in 3^{n}\Zd \cap \cu_m$. Denote~$U_j \coloneqq (z + \cu_j) \cap \cu_m$ for the (possibly truncated) cube at scale~$j \in \Z$. For each~$j \in [n, m-1] \cap \Z$, we define the excess
\begin{equation*}
	E_j \coloneqq E(u, U_j)
	=|U_j|^{-1/d}\min_{\ell \in \mathbb{L}}
	\|u-\ell\|_{\underline L^2(U_j)}\,.
\end{equation*}
We also denote~$\ell_j \coloneqq \ell(u, U_j)$ for the optimal affine function. Define the set of bad scales
\begin{equation*}
	\mathcal{B}_z \coloneqq \bigl\{ j \in [n, m-1] \cap \Z \,:\, \mathcal{G}(j, z \,; \tfrac 18 s, \tfrac 18 s  \delta^{\nf 12}) \ \text{is not valid} \bigr\}\,.
\end{equation*}
By~\eqref{e.scale.sep.for.good.events} and the definition of~$X_m(\alpha)$, since~$n \leq m - X_m(\alpha) \leq m - Z_m - k - 2$, we have
\begin{equation}
	\label{e.bad.scale.proportion.bound}
	|\mathcal{B}_z|
	\leq
	\max_{z'\in 3^{n} \Zd\cap \cu_m}
	\sum_{j=n}^m
	\Bigl(1 - \indc_{ \mathcal{G}(j,z' ; \frac18 s, \frac18 s \delta^{\nf 12})} \Bigr)
	\leq
	\delta (m - n+1)
	=
	C_1^{-1}(1-\alpha)(m-n+1)\,.
\end{equation}

\smallskip

\emph{Step 4.} Excess decay. An application of Lemma~\ref{l.excess.decay.good.scales} yields, for every~$j \in [n+k, m-1] \cap \Z \setminus \mathcal{B}_z$,
\begin{align}
	E(u, U_{j-k})
	&\leq
	3^{-\frac14 k} E(u, U_j)
	+
	C \ep_j(z) \bigl( \bigl| \nabla \ell_j \bigr| + \| \nabla h\|_{L^\infty(\cu_m)} \indc_{\{ z \notin \cu_{m-1} \}} \bigr)
	\notag \\ & \qquad
	+
	C \shom_j^{-1} 3^{\frac12 j} [ \g ]_{\underline{W}^{\nf12,\infty}(\cu_m)}
	+
	C 3^{\frac12  j} [ \nabla h ]_{\underline{W}^{\nf12,\infty}(\cu_m)} \indc_{\{ z \notin \cu_{m-1} \}}
	\,.
	\label{e.excess.decay.with.f}
\end{align}
Indeed, the definition of~$k$ and the choice of~$\delta_0$ make the two excess coefficients in~\eqref{e.excess.decay.one.step} at most~$\frac12 3^{-k/4}$ each. Since~$s=1/4$ and~$k=k(d)$, the remaining powers of~$s$ and~$3^k$ are absorbed into~$C(d)$. We also used
\begin{equation*}
	[\g ]_{\underline{H}^{\nf14}(U_j)}
	\leq
	C 3^{\frac{1}{4}j}
	[\g ]_{\underline{W}^{\nf12,\infty}(\cu_m) }
	\qand
	[\nabla h ]_{\underline{W}^{\nf12,\infty}(U_j)}
	\leq
	[ \nabla h ]_{\underline{W}^{\nf12,\infty}(\cu_m)}
	\,.
\end{equation*}

\smallskip

\emph{Step 5.} Application of the iteration lemma.
Set
\begin{equation*}
	\widehat{\mathcal B}_z
	\coloneqq
	\bigl(\mathcal B_z\cup([n,n+k)\cap\Z)\cup\{m\}\bigr)\cap[n,m]\,.
\end{equation*}
Then~$|\widehat{\mathcal B}_z|\leq|\mathcal B_z|+k+1$. We apply Lemma~\ref{l.iteration.lemma} with:
\begin{itemize}
\item~$h \coloneqq k$ (a constant),~$\theta \coloneqq 3^{-\nf14}$,
\item~$\ep_j \coloneqq C \ep_j(z)$ for~$j\leq m-1$ and~$\ep_m\coloneqq0$,
\item~$\delta_j \coloneqq C 3^{\frac 12 j} \shom_j^{-1} [\g]_{\underline{W}^{\nf12,\infty}(\cu_m)} + C \bigl( 3^{\frac 12 j} [ \nabla h ]_{\underline{W}^{\nf12,\infty}(\cu_m)}
+ \ep_j \| \nabla h\|_{L^\infty(\cu_m)} \bigr) \indc_{\{ z \notin \cu_{m-1} \}}
$ for~$j\leq m-1$, and~$\delta_m\coloneqq0$,
\item~$\mathcal{B} \coloneqq \widehat{\mathcal B}_z$.
\end{itemize}
By~\eqref{e.scale.sep.for.mathcal.E} and the definition of~$X_m(\alpha)$, we have
\begin{equation}
	\label{e.sum.eps.j.bound}
	\sum_{j=n}^{m-1} \ep_j
	\leq
	C \delta (m-n)
	=
	C C_1^{-1} (1-\alpha)(m-n)\,.
\end{equation}
Summing the series for the~$\delta_j$ terms and using~$3^{\frac 12 j}  \shom_j^{-1} \leq C 3^{-(\frac12 -\cgamma)(m-j)} 3^{\frac12 m} \shom_m^{-1}$ yields
\begin{align}
	\label{e.sum.delta.j.bound}
	\sum_{j=n}^{m-1} \delta_j
	&
	\leq
	C 3^{\frac12 m} \shom_m^{-1} [\g ]_{\underline{W}^{\nf12,\infty}(\cu_m)}
	\notag \\ & \qquad
	+ C  \bigl( 3^{\frac12 m}  [ \nabla h ]_{\underline{W}^{\nf12,\infty}(\cu_m)} + C_1^{-1}(1-\alpha) (m-n)\| \nabla h\|_{L^\infty(\cu_m)} \bigr) \indc_{\{ z \notin \cu_{m-1} \}}
	\,.
\end{align}
For each sub-window~$[n',m']\subset[n,m]$ with~$n'<m'$, we reapply Lemma~\ref{l.iteration.lemma} with
\begin{equation*}
	\left(
	\bigl(\mathcal B_z\cap[n',m']\bigr)
	\cup\bigl([n',n'+k)\cap\Z\bigr)
	\cup\{m'\}
	\right)\cap[n',m']
\end{equation*}
as its bad set. Its cardinality is at most~$|\mathcal B_z|+k+1$, while its sums of~$\ep_j$ and~$\delta_j$ are bounded by the corresponding full-window sums in~\eqref{e.sum.eps.j.bound} and~\eqref{e.sum.delta.j.bound}. Since~$k=k(d)$, the padding contributes only a dimensional prefactor, which we absorb into~$C$. The case~$n'=m'$ is immediate.
Consequently, for every~$n',m' \in \Z$ with~$n \leq n' \leq m' \leq m$,
\begin{align}
	\lefteqn{
	3^{-n'} \| u - (u)_{U_{n'}} \|_{\underline{L}^2(U_{n'})}
	} \quad &
	\notag \\ &
	\leq
	C\exp\biggl( C(k+1)(|\mathcal{B}_z|+1) + C \sum_{j=n}^{m-1} \ep_j \biggr)
	\biggl(
	3^{-m'} \| u - (u)_{U_{m'}} \|_{\underline{L}^2(U_{m'})}
	+
	\sum_{j=n}^{m-1} \delta_j
	\biggr)
	\notag \\ &
	\leq
	C \exp\bigl( C_1^{-1} C(1-\alpha) (m-n) \bigr)
	\Bigl(
	3^{-m'} \| u - (u)_{U_{m'}} \|_{\underline{L}^2(U_{m'})}
	+
	C 3^{\frac12 m} \shom_m^{-1} [ \g ]_{\underline{W}^{\nf12,\infty}(\cu_m)} \Bigr)
	\notag \\ & \qquad
	+
	C \exp\bigl( C_1^{-1} C(1-\alpha) (m-n) \bigr) 3^{\frac12 m}
	\| \nabla h \|_{\underline{W}^{\nf12,\infty}(\cu_m)} \indc_{\{ z \notin \cu_{m-1} \}}  \,,
	\label{e.oscillation.iteration.result}
\end{align}
where we used~\eqref{e.sum.eps.j.bound} and that~$k\leq C$,~$|\mathcal{B}_z| \leq \delta(m-n)$, and~$U_{m-1} \subseteq \cu_m$, by absorbing the prefactor~$C_1^{-1} C(1-\alpha) (m-n) $ into the exponential.

\smallskip

\emph{Step 6.} Verifying the H\"older exponent.
We now verify that the exponential prefactor in~\eqref{e.oscillation.iteration.result} is bounded by~$C 3^{\frac12(1-\alpha)(m-n)}$.
Choosing~$C_1 \geq 4C(k+1)/\log 3$, we have
\begin{equation*}
	\exp\bigl( C_1^{-1} C(1-\alpha) (m-n) \bigr)
	\leq
	3^{\frac14 (1-\alpha)(m-n)}\,.
\end{equation*}
We conclude that, for every~$n',m' \in \Z$ with~$n \leq n' \leq m' \leq m$,
\begin{align}
	\lefteqn{
	3^{-n'} \| u - (u)_{(z+\cu_{n'}) \cap \cu_m} \|_{\underline{L}^2((z+\cu_{n'}) \cap \cu_m)}
	} \qquad &
	\notag \\ &
	\leq
	C 3^{\frac12(1-\alpha)(m-n)}   3^{-{m'}} \| u - (u)_{(z+\cu_{m'}) \cap \cu_m} \|_{\underline{L}^2((z+\cu_{m'}) \cap \cu_m)}
	\notag \\ &  \qquad
	+
	C  3^{\frac12(1-\alpha)(m-n)}   3^{\frac12 m} \bigl( \shom_m^{-1} [ \g ]_{\underline{W}^{\nf12,\infty}(\cu_m)} + \| \nabla h \|_{\underline{W}^{\nf12,\infty}(\cu_m)} \indc_{\{ z \notin \cu_{m-1} \}} \bigr) \,.
	\label{e.oscillation.Holder.bound}
\end{align}

\smallskip

\emph{Step 7.} Conclusion.
It remains to convert the oscillation bound into a gradient bound.
For this, we use the coarse-grained Caccioppoli inequality (Lemma~\ref{l.coarse.grained.Caccioppoli.RHS}) at the nearest good scale. Fix~$x \in \cu_{m}$ and~$n \leq m - X_m(\alpha)$ as above, and choose~$z\in3^n\Zd\cap\cu_m$ by rounding each coordinate of~$x$ toward the origin. Then
\begin{equation*}
	(x+\cu_n)\cap\cu_m
	\subseteq
	(z+\cu_{n+1})\cap\cu_m
	\qand
	\indc_{\{z\notin\cu_{m-1}\}}
	\leq
	\indc_{\{x\notin\cu_{m-1}\}}\,.
\end{equation*}
Define~$n'$ and~$m'$ to be the smallest and the largest good scale, respectively, on~$[n,m-1] \cap \Z$:
\begin{equation*}
	n' \coloneqq \min \bigl\{ j \in [n+3, m-3] \cap \Z \,:\, j \notin \mathcal{B}_z \bigr\}
	\qand
	m' \coloneqq \max \bigl\{ j \in [n+3, m-3] \cap \Z \,:\, j \notin \mathcal{B}_z \bigr\}
	\,.
\end{equation*}
Observe that
$|\mathcal B_z|+7\leq m-n$ and
$\max\{ n' - n,m-m'\}  \leq |\mathcal{B}_z|+6$; in particular,
$n'<m'$ and~$m'-1\in[n',m']$.  Choose~$z'$ and~$z''$ to be the coordinatewise-clamped centers at the scales~$m'$ and~$n'$, respectively, as in the
proof of Lemma~\ref{l.excess.decay.good.scales}.  Then no boundary-indicator
band enlargement is needed and, for each
$(y',k') \in \{(z',m') , (z'',n') \}$,
$(z + \cu_{k'-1}) \cap \cu_m \subset y' + \cu_{k'-1}
\subset (z + \cu_{k'})\cap \cu_m$.
Set~$s_*\coloneqq1/16$ and define the flux-corrected family
\begin{equation*}
	\a_{L,k'}^{(z)}\coloneqq
	\a_L-(\k_L-\k_{k'})_{z+\cu_{k'}}\,.
\end{equation*}
On~$\mathcal{G}(k',z;1/32,\delta^{\nf12}/32)$, the good event first gives the
$q=2$ ellipticity bounds on the outer cube at~$s_*/2=1/32$.  By~\eqref{e.compareLambdaqs} and monotonicity, this gives the
$q=1$ and~$q=2$ bounds there at~$s_*=1/16$.  Lemma~\ref{l.lambdas.stability},
applied from this outer slot to the child slots, and~$\lambda\leq\Lambda$ give
\begin{equation}
	\left\{
	\begin{aligned}
		& \lambda_{1/8,2}^{-1}(y'+\cu_{k'-1};\a_{L,k'}^{(z)})
		\leq
		C \shom_{k'}^{-1}\,,
		\\
		& \lambda_{1/8,1}^{-1}(y'+\cu_{k'-1};\a_{L,k'}^{(z)})
		\leq
		C\shom_{k'}^{-1}\,,
		\qquad
		\lambda_{1/8,1}(y'+\cu_{k'-1};\a_{L,k'}^{(z)})
		\leq
		C\shom_{k'}\,,
		\\
		& \Lambda_{r,1}(y'+\cu_{k'-1};\a_{L,k'}^{(z)})
		\leq C\shom_{k'}
		\qquad\text{for }r\in\{1/8,1/4\}\,.
	\end{aligned}
	\right.
	\label{e.lambda.stability.applied}
\end{equation}
Since~$n'$ is a good scale, we may apply the coarse-grained Caccioppoli inequality with~$s=1/4$ and~$t=1/8$, so that~$\sigma=5/8$, and the oscillation bound~\eqref{e.oscillation.Holder.bound} at scale~$n'$ to obtain a gradient bound. Since~$(z + \cu_{n+1}) \cap \cu_m \subseteq U_{n'-2}$, passing to the smaller set costs a volume factor~$3^{\frac{d}{2}(n'-n)} \leq 3^{\frac d2 (|\mathcal{B}_z|+6)}$. By~\eqref{e.bad.scale.proportion.bound}, this is at most~$C 3^{\frac14(1-\alpha)(m-n)}$.  Before applying the estimate we use
$\shom_{n'}\leq4\shom_m$; taking square roots costs at most~$2$.
We obtain by~\eqref{e.oscillation.Holder.bound},~\eqref{e.lambda.stability.applied} and Lemma~\ref{l.coarse.grained.Caccioppoli.RHS} that
\begin{align}
	\lefteqn{
	\nu^{\nf12}\| \nabla u \|_{\underline{L}^2((z + \cu_{n+1}) \cap \cu_m)}
	} \quad &
	\notag \\ &
	\leq
	C \shom_{n'}^{\nf12} \cdot 3^{\frac34 (1-\alpha)(m-n)}
	3^{-(m'-1)} \| u - (u)_{(z+\cu_{m'-1}) \cap \cu_m} \|_{\underline{L}^2((z+\cu_{m'-1}) \cap \cu_m)}
	\notag \\ & \qquad
	+
	C 3^{\frac34 (1-\alpha)(m-n)}
	3^{\frac12 m} \bigl( \shom_m^{-\nf12} [ \g ]_{\underline{W}^{\nf12,\infty}(\cu_m)} + \shom_m^{\nf12}  \| \nabla h \|_{\underline{W}^{\nf12,\infty}(\cu_m)} \indc_{\{ x \notin \cu_{m-1} \}} \bigr)
	\,.
	\label{e.gradient.with.shom}
\end{align}
On the other hand, apply Lemma~\ref{l.coarse.graining.RHS} at~$t=1/4$. The exact depth conversion is
\begin{equation*}
	3^{-(m'-1)}\|\nabla u\|_{\Besov{-1}{2}{1}(z'+\cu_{m'-1})}
	\leq(1-3^{-3/2})^{-1/2}3^{-\frac14(m'-1)}
	\|\nabla u\|_{\Besov{-\nf14}{2}{2}(z'+\cu_{m'-1})}
	\,.
\end{equation*}
Its constant is less than~$1.12$; the harmless one-scale factor converts~$3^{-(m'-1)}$ to the normalization below. Together with~\eqref{e.lambda.stability.applied}, this yields
\begin{align}
	\label{e.cg.Poincare.with.rhs.grad.applied}
	\lefteqn{
	3^{-m'}\| \nabla u \|_{\Besov{-1}{2}{1}(z'+\cu_{m'-1})}
	\indc_{\mathcal{G}(m', z \,; \frac 18 s, \frac 18 s  \delta^{\nf 12})}
	} \qquad &
	\notag \\ &
	\leq
	C  \shom_{m'}^{-\nf12} \nu^{\nf 12}
	\| \nabla u\|_{\underline{L}^2(z'+\cu_{m'-1})}
	+
	C
	\shom_{m'}^{-1}
	3^{sm'} [ \g ]_{\underline{H}^s(z'+\cu_{m'-1}) }
	\notag \\ &
	\leq
	C 3^{\frac12 (d + \cgamma) (m-m')}
	\bigl(
	\shom_{m}^{-\nf12} \nu^{\nf 12}
	\| \nabla u\|_{\underline{L}^2(\cu_{m})}
	+
	\shom_{m}^{-1}
	3^{sm} [ \g ]_{\underline{H}^s(\cu_{m}) }
	\bigr)
	\,\,.
\end{align}
Thus, by~\eqref{e.fractional.Sobolev.embedding}, we deduce that
\begin{align*}
	\lefteqn{
	3^{-(m'-1)} \| u - (u)_{(z+\cu_{m'-1}) \cap \cu_m} \|_{\underline{L}^2((z+\cu_{m'-1}) \cap \cu_m)}  \indc_{\mathcal{G}(m', z \,; \frac 18 s, \frac 18 s  \delta^{\nf 12})}
	} \qquad &
	\notag \\ &
	\leq
	C 3^{-(m'-1)} \| u - (u)_{(z'+\cu_{m'-1})} \|_{\underline{L}^2(z'+\cu_{m'-1})}  \indc_{\mathcal{G}(m', z \,; \frac 18 s, \frac 18 s  \delta^{\nf 12})}
	\notag \\ &
	\leq
	C \shom_{m}^{-\nf12}  3^{\frac12 (d + \cgamma) (m-m')}
	\bigl(
	\nu^{\nf 12} \| \nabla u\|_{\underline{L}^2(\cu_{m})}
	+
	\shom_{m}^{-\nf 12} 3^{sm} [ \g ]_{\underline{H}^s(\cu_{m}) }
	\bigr)
	\,.
\end{align*}
Since~$m-m' \leq |\mathcal{B}_z| +3 \leq C_1^{-1} (1-\alpha)(m-n) + 3$, the previous display and~\eqref{e.gradient.with.shom}, together with~$\shom_{n'} \leq 4 \shom_m$, give the desired estimate on~$(z+\cu_{n+1})\cap\cu_m$. The first inclusion above and normalized-volume comparison transfer it to~$(x+\cu_n)\cap\cu_m$ at a cost of at most~$3^d$, proving~\eqref{e.energy.density.estimate} for~$\a_L$ uniformly in~$L\geq m$.

It remains to remove the infrared cutoff. Subtract the constant antisymmetric matrix~$(\k_L)_{\cu_m}$ from~$\a_L$; this does not change the equation. By the construction of~$\k$ modulo constants, the resulting coefficients converge uniformly on~$\cu_m$ to a representative of~$\a$. The symmetric part remains~$\nu\Id$, so the corresponding Dirichlet solutions converge strongly in~$H^1(\cu_m)$ by testing the difference equation. We may therefore let~$L\to\infty$ in the preceding uniform estimate, which proves the theorem for~$\a$.
\end{proof}

\subsection{Harmonic approximation in \texorpdfstring{$L^\infty$}{L-infty}: the proof of Theorem~\ref{t.homogenization}}
\label{ss.proof.theorem.B}

In this subsection we complete the proof of Theorem~\ref{t.homogenization}. The strategy is to combine the general coarse-graining estimate of Proposition~\ref{p.general.coarse.graining} (applied at the fixed exponent~$4d$) with Theorem~\ref{t.regularity} (which bounds the energy density at small scales) and the bounds on the homogenization error from Proposition~\ref{p.induction.bounds}.

\smallskip

To prepare for the proof of Theorem~\ref{t.homogenization}, we use the results of Section~\ref{s.RG.flow} and the sensitivity estimate in Lemma~\ref{l.J.sensitivity.no.conditions} to estimate the random variable appearing on the right side of~\eqref{e.general.coarse.graining.estimate}. In this subsection, we allow all constants~$C$ and~$c$ to depend on a lower bound for~$\cstar$ in addition to~$d$.

\begin{lemma}
\label{l.bounds.mathcal.E.aL}
There exists a constant~$C(d, \cstar) < \infty$ such that, if the parameter~$\cgamma$ satisfies~$\cgamma \leq C^{-10} $, then, for every~$m,n\in\Z$ with~$n \leq m \leq n + \cgamma^{-1}$,~$s \in [C^2\cgamma^{\nf12},\nf14]$, and~$p\in [2 ds^{-1},C^{-1} \cgamma^{-1}s ]$,
\begin{equation}
	\E \biggl[
	\biggl( \sup_{L\geq m}
	\mathcal{E}_{s,\infty,2}
	(\cu_m , n;\a_L - (\k_L - \k_m)_{\cu_m} , \shom_m)
	\biggr)^{\!p}
	\biggr]^{\nf1p}
	\leq
	C 3^{\frac12s(m-n)} p^{\nf12}
	\bigl( s^{-1}  + (m-n)^{\nf12} \bigr)
	\cgamma^{\nf12}
	\,.
	\label{e.bounds.mathcal.E.aL}
\end{equation}
\end{lemma}
\begin{proof}
Denote~$\tilde{\a}_{L,m} \coloneqq \a_L - (\k_L - \k_m)_{\cu_m}$ and put~$D\coloneqq m-n$. For~$j\leq n$, set
\begin{align*}
	G_j
	&\coloneqq
	\biggl(
	\avsum_{z\in3^j\Zd\cap\cu_m}
	\sup_{L\geq m}\max_{|e|=1}
	\bigl(J(z+\cu_j,\shom_m^{-\nf12}e,\shom_m^{\nf12}e\,;\tilde{\a}_{L,m})\bigr)^{\nf p2}
	\biggr)^{\!\nf2p}\,,
	\\
	G_j^t
	&\coloneqq
	\biggl(
	\avsum_{z\in3^j\Zd\cap\cu_m}
	\sup_{L\geq m}\max_{|e|=1}
	\bigl(J(z+\cu_j,\shom_m^{-\nf12}e,\shom_m^{\nf12}e\,;\tilde{\a}_{L,m}^t)\bigr)^{\nf p2}
	\biggr)^{\!\nf2p}
	\,.
\end{align*}
Applying~\eqref{e.mathcalE.infty.to.q} with exponents~$p$ and~$2$, using~$p\geq2ds^{-1}$, and returning to the unraised~$q=2$ scale sum gives
\begin{equation*}
	\biggl(\sup_{L\geq m}\mathcal E_{s,\infty,2}(\cu_m,n;\tilde{\a}_{L,m},\shom_m)\biggr)^2
	\leq C3^{sD}\css{2s}\sum_{j=-\infty}^n3^{-s(n-j)}\bigl(G_j+G_j^t\bigr)
	\,.
\end{equation*}
As in~\eqref{e.ugly.estimate.for.J.pre} in Step~2 of the proof of Proposition~\ref{p.mathcalE.annular.decomp}, applying the sensitivity estimate~\eqref{e.J.sensitivity.no.conditions} of Lemma~\ref{l.J.sensitivity.no.conditions} in the cube~$\cu_j$ with~$\mu \coloneqq \shom_m \shom^{-1}_{j-2}$,~$\s_0 \coloneqq \shom_{j-2}$ and~$\h \coloneqq \k_L  - \k_{j-2} - (\k_L - \k_m)_{\cu_m}$ yields
\begin{align}
	\lefteqn{
	J(\cu_j,\shom_m^{-\nf12}e,\shom_m^{\nf 12}e\,;\tilde{\a}_{L,m})
	} \quad &
	\notag \\ &
	\leq
	C\shom_m^{-1}\shom_{j-2} \bigl( 1 +
	3^{2j} \| \nabla  (\k_L - \k_{j-2}) \|_{\underline{W}^{1,\infty}(\cu_j)} \lambda_{\cgamma,2}^{-1}
	(\cu_j; \a_{j-2} )
	\bigr)^{\! \frac{2s}{1-2\cgamma} }
	\mathcal{E}_{s,2,2}(\cu_j;\a_{j-2},\shom_{j-2})^2
	\notag \\ & \qquad
	+
	C \shom_m^{-1} \shom_{j-2}^2
	\lambda_{\cgamma,2}^{-1} (\cu_j \,; \a_{j-2})
	\bigl( 1+
	\lambda_{\cgamma,2}^{-1}
	(\cu_j; \a_{j-2} )
	3^{2j}  \| \nabla (\k_L - \k_{j-2}) \|_{\underline{W}^{1,\infty}(\cu_j)}
	\bigr)^{\frac{2\cgamma}{1-2\cgamma}}
	\notag \\ &  \qquad \qquad
	\times
	\bigl(
	\shom_{j-2}^{-2} \| \k_L - (\k_L - \k_m)_{\cu_m} - \k_{j-2} \|_{\underline{L}^2(\cu_j )}^2 +
	|\shom_m \shom_{j-2}^{-1}-1|^2 \bigr)
	\notag \\ &  \qquad
	+
	C\bigl( \shom_m^{-2} + \lambda_{\cgamma,2}^{-2}
	(\cu_j; \a_{j-2} ) \bigr)
	3^{4j}  \| \nabla (\k_L - \k_{j-2}) \|_{\underline{W}^{1,\infty}(\cu_j)}^2
	\,.
	\label{e.apply.sensitivity.J.aL}
\end{align}
We will bound the~$\nf{p}2$th moments of each of the factors appearing on the right side of~\eqref{e.apply.sensitivity.J.aL}. For many of these, we use Lemma~\ref{l.moments.gamma.psi} to convert~$\O_{\Gamma_\sigma}$-type estimates into moment bounds.
The terms are bounded using the following:
\begin{itemize}

\item By~\eqref{e.formula.for.shom}, we have~$\shom_m^{-1} \shom_j \leq C$ for every~$j\leq m$.

\item
By Lemma~\ref{l.shom.continuity},
we have, for every~$j\leq m$,
\begin{equation*}
	|\shom_m \shom_{j-2}^{-1}-1|^2 \leq
	C\min\bigl\{ 1 \,, \cgamma^2 (m-j)^2 + C  \cgamma^2 \left|\log \cgamma\right|^{4} \bigr\} 3^{2 \cgamma (m-j)}
	\,.
\end{equation*}

\item We have~$3^{\cgamma j}\shom_{j-2}^{-1} \leq C  \cgamma^{\nf12}$ by~\eqref{e.formula.for.shom};

\item By Proposition~\ref{p.induction.bounds}, we have
\begin{equation*}
	\mathcal{E}_{s,2,2}(\cu_j;\a_{j-2},\shom_{j-2})
	\leq
	\O_{\Gamma_{2} }
	\bigl( C  s^{-1} \cgamma^{\nf12}  \bigr) +
	\O_{\Gamma_{\nf12}}
	\bigl( \exp( - C^{-1}  \cgamma^{-1}) \bigr)
\end{equation*}
and hence, for every~$q \in [1,\infty)$,
\begin{equation*}
	\E \bigl[  \bigl( \mathcal{E}_{s,2,2}(\cu_j;\a_{j-2},\shom_{j-2}) \bigr)^q \bigr]^{\nf1q}
	\leq
	C s^{-1} q^{\nf12}  \cgamma^{\nf12}
	+
	C q^2 \exp( - C^{-1} \cgamma^{-1} )
	\,.
\end{equation*}

\item Partition~$\cu_j$ into its~$9^d$ descendants at scale~$j-2$. Apply Proposition~\ref{p.cg.ellipticity.bounds} with~$s=\cgamma$ to each descendant at the matched cutoff~$j-2$, raise the bound to~$\cu_j$ by the two-step subadditivity estimate, and then use the finite-maximum estimate and the fixed-scale comparison of~$\shom_{j-3}$ and~$\shom_{j-1}$. This gives~$\lambda_{\cgamma,2}^{-1} (\cu_j \,; \a_{j-2}) \leq
C \shom_{j-1}^{-1}
+
\O_{\Gamma_{\nf13}} ( C \shom_{j-1}^{-1} \exp(-C^{-1}  \cgamma^{-1} )  ) $ and hence, for every~$q\in[1,\infty)$,
\begin{equation*}
	\E \bigl[ \lambda_{\cgamma,2}^{-q} (\cu_j \,; \a_{j-2}) \bigr]^{\nf1q}
	\leq
	C \shom_{j-1}^{-1}
	\bigl( 1
	+
	C q^3 \exp(-C^{-1}  \cgamma^{-1} )
	\bigr)
	\,.
\end{equation*}

\item By~\eqref{e.nabla.km.kn.infty}, we have~$3^{2j}  \| \nabla (\k_L - \k_{j-2}) \|_{\underline{W}^{1,\infty}(\cu_j)} \leq
\O_{\Gamma_2} \bigl(  C 3^{\cgamma j} \bigr)$ and hence,
for every~$q\in[1,\infty)$,
\begin{equation*}
	\E \bigl[ \bigl( 3^{2j}  \| \nabla (\k_L - \k_{j-2}) \|_{\underline{W}^{1,\infty}(\cu_j)}  \bigr)^q \bigr]^{\nf1q}
	\leq
	C q^{\nf12} 3^{\cgamma j} \,.
\end{equation*}

\item By~\eqref{e.nabla.km.kn.infty} and~\eqref{e.k.ell.upscales.infty}, we have
\begin{equation*}
	\| \k_L - \k_{j-2}  - (\k_L - \k_{m} )_{\cu_m} \|_{\underline{L}^2(\cu_j )} \leq \O_{\Gamma_2}\bigl(  C \bigl(1+\min\bigl\{ \cgamma^{-\nf 12} ,  (m-j)^{\nf 12} \bigr\}\bigr) 3^{\cgamma m}  \bigr)
	\,.
\end{equation*}
and hence, for every~$q\in[1,\infty)$,
\begin{equation*}
	\E \bigl[ \| \k_L - \k_{j-2}  - (\k_L - \k_{m} )_{\cu_m} \|_{\underline{L}^2(\cu_j )}^q \bigr]^{\nf1q}
	\leq
	C q^{\nf12} \bigl(1+\min\bigl\{ \cgamma^{-\nf 12} ,  (m-j)^{\nf 12} \bigr\}\bigr) 3^{\cgamma m}
	\,.
\end{equation*}
\end{itemize}

The right side of~\eqref{e.apply.sensitivity.J.aL} is independent of~$e$. We apply it in each translated cube~$z+\cu_j$, and apply the same estimate to the transposed field. Expanding the increments into shells, dominating the finite upper tails uniformly in~$L\geq m$ by the corresponding summable infinite tails, and using the preceding moment bounds and H\"older's inequality gives
\begin{equation*}
	\E\bigl[G_j^{\nf p2}\bigr]^{\nf2p}
	+\E\bigl[(G_j^t)^{\nf p2}\bigr]^{\nf2p}
	\leq Cp\cgamma\Bigl(s^{-2}+(m-j)3^{2\cgamma(m-j)}+3^{2\cgamma(m-j)}\Bigr)
	\,.
\end{equation*}
We now apply Minkowski's inequality in~$L^{p/2}$ to the unraised scale sum. Since~$4\cgamma\leq s$,~$\cgamma D\leq1$ and~$s\leq1$, the geometric scale sums yield
\begin{align*}
	\E\biggl[\biggl(\sup_{L\geq m}\mathcal E_{s,\infty,2}(\cu_m,n;\tilde{\a}_{L,m},\shom_m)\biggr)^p\biggr]^{\nf2p}
	&\leq C3^{sD}\css{2s}\sum_{j=-\infty}^n3^{-s(n-j)}
	\Bigl(\E\bigl[G_j^{\nf p2}\bigr]^{\nf2p}+\E\bigl[(G_j^t)^{\nf p2}\bigr]^{\nf2p}\Bigr)
	\\
	&\leq Cp\cgamma3^{sD}\bigl(s^{-2}+D\bigr)
	\,.
\end{align*}
Taking square roots proves~\eqref{e.bounds.mathcal.E.aL}.
\end{proof}

\begin{proof}[{Proof of Theorem~\ref{t.homogenization}}]
The estimate~\eqref{e.shom.m.bounds.thmB} for~$\shom_m$ is a restatement of~\eqref{e.formula.for.shom} from Proposition~\ref{p.induction.bounds}. It remains to prove~\eqref{e.mathcal.E.estimate} and~\eqref{e.intro.homogenization.estimate}. We also recall that~$[ \cdot ]_{W^{s,\infty}(\cu_m)} = [ \cdot ]_{\underline{W}^{s,\infty}(\cu_m)} $ is equivalent, up to dimensional constants, to~$[\cdot ]_{C^{0,s}(\cu_m)}$ and so we use these interchangeably.

\smallskip

We first prove the result for~$\a_L$, uniformly in~$L\geq m$, and then pass
to~$L\to\infty$. As in the theorem statement,
we fix data~$\g \in C^{0,\nf12}(\cu_m;\Rd)$ and~$h\in C^{1,\nf12}( \cu_m)$ and consider a pair~$u,v\in H^1(\cu_m)$ satisfying
\begin{equation}
	\left\{
	\begin{aligned}
		& -\nabla\cdot \a_L  \nabla u = \nabla \cdot \g
		& \mbox{in} & \ \cu_m\,, \\
		& u = h & \mbox{on} & \ \partial \cu_m\,,
	\end{aligned}
	\right.
	\qquad
	\mbox{and}
	\qquad
	\left\{
	\begin{aligned}
		& -\shom_m \Delta v = \nabla \cdot \g
		& \mbox{in} & \ \cu_m\,, \\
		& v = h & \mbox{on} & \ \partial \cu_m\,.
	\end{aligned}
	\right.
	\label{e.dirichlet.pairs.proof}
\end{equation}
We will apply the homogenization estimate of Proposition~\ref{p.general.coarse.graining} and use Theorem~\ref{t.regularity} and Lemma~\ref{l.bounds.mathcal.E.aL} to bound the terms appearing on the right side of~\eqref{e.general.coarse.graining.estimate}. Note that the first equation in~\eqref{e.dirichlet.pairs.proof} is unchanged if we replace~$\a_L$ by~$\tilde{\a}_{L,m} \coloneqq \a_L - (\k_L-\k_m)_{\cu_m}$.  We also assume, without loss of generality, that
\begin{equation*}
	[ \g ]_{\underline{W}^{\nf12,\infty}(\cu_m)} \leq 1 \qand \| \nabla h \|_{\underline{W}^{\nf12,\infty}(\cu_m)} \leq 1 \,\,.
\end{equation*}
We select the following parameters:
\begin{equation}
	s\coloneqq \left| \log \cgamma \right|^{-1} \,, \qquad
	s_1 \coloneqq \frac18s \,,
	\qquad
	s_2 \coloneqq \frac{49}{100}\,,
	\qquad
	\alpha \coloneqq 1 - \frac12 s
	\,,
	\qquad
	k \coloneqq \lceil 10 \left| \log \cgamma \right| \rceil
	\,.
	\label{e.homog.parameter.choices}
\end{equation}
Observe that, if~$\cgamma$ is sufficiently small, then~$s\leq1/4$ and
$(s_2-s)^{-1}\leq5$, and~$\alpha\in (0,1-C_0\cgamma^{\nf12}]$, where~$C_0(d,\cstar)$ is the constant from Theorem~\ref{t.regularity}. These choices satisfy the constraints~$s_1 < s < s_2 < 1$ required by Proposition~\ref{p.general.coarse.graining}.
The mesoscopic scale parameter~$n$ in Proposition~\ref{p.general.coarse.graining} will be chosen to be
\begin{equation}
	\label{e.mesoscale.def}
	n \coloneqq m - k\,.
\end{equation}
We let~$X_m(\alpha)$ be the minimal scale random variable from Theorem~\ref{t.regularity}. Fix a sufficiently large constant~$C_B(d,\cstar)$, which remains fixed below, and define the random variable
\begin{align*}
	\EthmB(m)
	&\coloneqq
	s^{-2} 3^{(1-\alpha) X_m(\alpha)}
	\sup_{L \geq m}
	\mathcal{E}_{\nf{s_1}{2},\infty,2}
	\bigl (\cu_m,n ;\tilde{\a}_{L,m},\shom_m \bigr)
	\\[-2pt]
	&\qquad {}\times
	\bigl(1 + \sup_{L \geq m}
	\mathcal E_{\nf18,\infty,2}(\cu_m;\tilde\a_{L,m},\shom_m)\bigr)
	\\[-2pt]
	&\qquad {}\times
	\bigl(1 + \sup_{L \geq m}
	\mathcal E_{\nf18,\infty,2}(\cu_m;\tilde\a_{L,m},\shom_m)\bigr)
	\\
	&\qquad \qquad {}
	+ C_B s^{-1}\cgamma^{5}
	\Bigl(1 + \sup_{L \geq m}
	\mathcal{E}^{2}_{\nf{s_1}{4}, \infty, 2}
	(\cu_m, n;\tilde{\a}_{L,m},\shom_m) \Bigr)  \,\,.
\end{align*}
\emph{Step 1.} We establish the moment bound~\eqref{e.mathcal.E.estimate}. The tail estimate in Theorem~\ref{t.regularity} gives, for every~$q\leq C^{-1}(1-\alpha)\cgamma^{-1}$,
\begin{equation*}
	\bigl\|3^{(1-\alpha)X_m(\alpha)}\bigr\|_{L^q}
	\leq C\,.
\end{equation*}
For the exponent~$p$ in~\eqref{e.mathcal.E.estimate}, put~$q\coloneqq\max\{4p,64ds^{-1}\}$. After decreasing the constant in the upper bound for~$p$ and~$\cgamma_0$, we have~$q\leq C^{-1}\cgamma^{-1}s$. Applying Lemma~\ref{l.bounds.mathcal.E.aL} at exponent~$q$, then using monotonicity of~$L^p$ norms and H\"older's inequality, yields~\eqref{e.mathcal.E.estimate}; here~$q^{1/2}\leq C(p^{1/2}+s^{-1/2})$ and~$s=\left|\log\cgamma\right|^{-1}$.

\smallskip

\emph{Step 2.} The energy density estimate.
For brevity set
\begin{align*}
	\mathfrak D_{L,m}
	&\coloneqq
	C_B\bigl(1+\mathcal E_{\nf18,\infty,2}
	(\cu_m;\tilde\a_{L,m},\shom_m)\bigr)
	\shom_m^{-\nf12}3^{m/2}
	[\g]_{\underline W^{\nf12,\infty}(\cu_m)}
	\\&\quad+
	C_B\bigl(1+\mathcal E_{\nf18,\infty,2}
	(\cu_m;\tilde\a_{L,m},\shom_m)\bigr)
	\shom_m^{\nf12}3^{m/2}
	\|\nabla h\|_{\underline W^{\nf12,\infty}(\cu_m)}\,.
\end{align*}
We show that
\begin{align}
	\sup_{j \leq n} \, \sup_{z \in 3^j \Zd \cap \cu_m}
	3^{-(s - s_1)(n-j)}
	\nu^{\nf12} \| \nabla u \|_{\underline{L}^2(z + \cu_j)}
	\leq
	C
	3^{(1-\alpha) X_m(\alpha)}
	\mathfrak D_{L,m}
	\,.
	\label{e.energy.bound.ell.infty}
\end{align}
By the energy bound~\eqref{e.cg.RHS} and~\eqref{e.bound.Lambdas.by.Es}, we have
\begin{align}
	\label{e.energy.from.data}
	\nu^{\nf12} \| \nabla u \|_{\underline{L}^2(\cu_m)}
	&\leq
	C\bigl(1+\mathcal E_{\nf18,\infty,2}
	(\cu_m;\tilde\a_{L,m},\shom_m)\bigr)
	\shom_m^{-\nf12}3^{m/2}
	[\g]_{\underline W^{\nf12,\infty}(\cu_m)}
	\notag\\
	&\qquad+C\bigl(1+\mathcal E_{\nf18,\infty,2}
	(\cu_m;\tilde\a_{L,m},\shom_m)\bigr)
	\shom_m^{\nf12}3^{m/2}
	\|\nabla h\|_{\underline W^{\nf12,\infty}(\cu_m)}
	\,.
\end{align}
Here~\eqref{e.cg.RHS} is applied at~$t=1/4$, so both error indices are~$1/8$ and the powers of~$t$ are absorbed into~$C$.
After bounding the data terms in Theorem~\ref{t.regularity} by this right-hand
side (using~$[\g]_{\underline W^{1/2,\infty}}\leq1$ for its Besov term),
we obtain, for every~$j \leq m - X_m(\alpha)$ and every~$z \in 3^j \Zd \cap \cu_m$,
\begin{equation*}
	\nu^{\nf12} \| \nabla u \|_{\underline{L}^2(z + \cu_j)}
	\leq
	C 3^{(1-\alpha)(m-j)}
	\mathfrak D_{L,m}
	\,.
\end{equation*}
If~$m-X_m(\alpha)<j\leq m$, partition~$z+\cu_j$ into its descendants at scale~$m-X_m(\alpha)$ and average their normalized energies. Applying the preceding estimate on every descendant and using~$3^{(1-\alpha)(m-j)}\geq1$ gives the same bound.
Hence, for every~$j\leq m$ and~$z \in 3^j \Zd \cap \cu_m$,
\begin{equation*}
	\nu^{\nf12} \| \nabla u \|_{\underline{L}^2(z + \cu_j)}
	\leq
	C3^{(1-\alpha) X_m(\alpha)} 3^{(1-\alpha)(m-j)}
	\mathfrak D_{L,m}
	\,.
\end{equation*}
By the definitions of~$s_1$ and~$\alpha$, we have~$s-s_1=7s/8\geq s/2=1-\alpha$.
Since~$3^{(1-\alpha)(m-n)} \leq 3^{5+s/2}$ (and hence is bounded by
$3^{41/8}$ when~$s\leq1/4$), the previous display yields
\eqref{e.energy.bound.ell.infty}.

\smallskip

\emph{Step 3.} The proof of~\eqref{e.intro.homogenization.estimate}. We apply Proposition~\ref{p.general.coarse.graining} with~$p=4d$,~$\a=\tilde{\a}_{L,m}$,~$\s_0 = \shom_m$, and the parameters~$s$,~$s_1$,~$s_2$ and~$n$ defined in~\eqref{e.homog.parameter.choices} and~\eqref{e.mesoscale.def}. The weighted~$\ell^{4d}$-to-supremum conversion contributes the factor~$Cs^{-1/(4d)}$, and~$s_2<1/2$ permits the H\"older--Gagliardo conversion.
Substituting~\eqref{e.energy.bound.ell.infty} into~\eqref{e.general.coarse.graining.estimate}, we obtain
\begin{align}
	\lefteqn{
	3^{-ms} \shom_m \| \nabla u - \nabla v \|_{\Wminusul{-s}{4d}(\cu_m)}
	} \qquad &
	\notag \\ &
	\leq
	C s^{-1-\nf1{4d}} \shom_m^{\nf12}
	\mathcal{E}_{s_1, \infty, 1} (\cu_m, n;\tilde{\a}_{L,m},\shom_m)
	\sup_{k \leq n}
	3^{-(s-s_1) (n-k)}
	\max_{z\in 3^k\Zd\cap\cu_m}
	\nu^{\nf12 } \| \nabla u \|_{\underline{L}^2(z+\cu_k)}
	\notag \\ & \qquad
	+
	Cs^{-\nf 9 2}
	(s_2 - s)^{-1}
	\bigl(1 + \mathcal{E}^{2}_{\nf{s_1}{2}, \infty, 2}(\cu_m, n;\tilde{\a}_{L,m},\shom_m) \bigr)
	3^{s_2 n}
	\notag \\ &
	\leq
	C s^{-1-\nf1{4d}}\shom_m^{\nf12}
	3^{(1-\alpha) X_m(\alpha)}
	\mathcal{E}_{s_1, \infty, 1} (\cu_m, n;\tilde{\a}_{L,m},\shom_m)
	\mathfrak D_{L,m}
	\notag \\ & \qquad
	+
	Cs^{-\nf92}3^{s_2(n-m)}
	\bigl(1 + \mathcal{E}^{2}_{\nf{s_1}{2}, \infty, 2}(\cu_m, n;\tilde{\a}_{L,m},\shom_m) \bigr)
	3^{\frac12 m}
	\,.
	\label{e.gradient.homog.bound.pre}
\end{align}
Using the parameter choices in~\eqref{e.homog.parameter.choices} and that~$m-n \geq k\geq 10\left| \log\cgamma\right|$, we find that
\begin{equation*}
	s^{-\nf92}3^{s_2(n-m)}
	\leq
	C \cgamma^{10s_2 \log 3} \left| \log \cgamma\right|^{\nf92}
	\leq C \cgamma^5\,.
\end{equation*}
By~\eqref{e.compareEqs},
\begin{equation*}
	\mathcal E_{s_1,\infty,1}(\cu_m,n;\tilde\a_{L,m},\shom_m)
	\leq C\mathcal E_{\nf{s_1}{2},\infty,2}(\cu_m,n;\tilde\a_{L,m},\shom_m)
	\,.
\end{equation*}
Monotonicity in the first index also gives
\begin{equation*}
	\mathcal E_{\nf{s_1}{2},\infty,2}(\cu_m,n;\tilde\a_{L,m},\shom_m)
	\leq C\mathcal E_{\nf{s_1}{4},\infty,2}(\cu_m,n;\tilde\a_{L,m},\shom_m)
	\,.
\end{equation*}
Substituting this into the right side of~\eqref{e.gradient.homog.bound.pre} and comparing to the definition of~$\EthmB(m)$, while retaining one of its two factors of~$s^{-1}$, we find that
\begin{equation}
	3^{-ms} \| \nabla u - \nabla v \|_{\Wminusul{-s}{4d}(\cu_m)}
	\leq
	C s^{1-\nf1{4d}}
	\EthmB(m)
	\bigl(
	\shom_m^{-1}
	3^{\frac12 m}
	+
	3^{\frac12 m}
	\bigr)
	\,.
	\label{e.weak.gradient.bound}
\end{equation}
Since~$u-v\in H^1_0(\cu_m)$, its average gradient vanishes. Moreover,~$1-s>d/(4d)=1/4$, so the fractional Ne\v{c}as inequality~\eqref{e.fractional.Necas} and the fractional Morrey inequality give, for the continuous representative,
\begin{equation*}
	3^{-m}\|u-v\|_{L^\infty(\cu_m)}
	\leq Cs^{-\nf1{4d}}3^{-ms}\|\nabla u-\nabla v\|_{\Wminusul{-s}{4d}(\cu_m)}
	\,.
\end{equation*}
Combining this estimate with~\eqref{e.weak.gradient.bound}, using~$s^{1-\nf1{2d}}\leq1$, and applying the normalization of the data gives~\eqref{e.intro.homogenization.estimate}.

\smallskip

\emph{Step 4.} The proof of~\eqref{e.intro.energies}. By the same application of Proposition~\ref{p.general.coarse.graining} as in Step~3, with the coefficient~$\tilde{\a}_{L,m}$, we obtain the following bound for the difference of the fluxes:
\begin{equation}
	3^{-ms} \shom_m^{-1}  \| \tilde{\a}_{L,m}\nabla u -  \shom_m  \nabla v \|_{\underline{W}^{-s,4d}(\cu_m)}
	\leq
	C
	\EthmB(m)
	\bigl(
	\shom_m^{-1}
	3^{\frac12 m}
	+
	3^{\frac12 m}
	\bigr)
	\,.
	\label{e.weak.flux.bound.again}
\end{equation}
The constant matrix~$\a_L-\tilde{\a}_{L,m}$ is antisymmetric, so it does not change either the equation or the quadratic energy. We next compute
\begin{align*}
	\biggl| \fint_{\cu_m}
	\bigl(  \nabla u \cdot \a_L\nabla u - \nabla v \cdot \shom_m \nabla v \bigr)
	\biggr|
	&
	=
	\biggl| \fint_{\cu_m}
	\bigl(
	\nabla v \cdot ( \tilde{\a}_{L,m} \nabla u -  \shom_m  \nabla v )
	+
	(\nabla u - \nabla v) \cdot \tilde{\a}_{L,m} \nabla u
	\bigr)
	\biggr|
	\,.
\end{align*}
Testing the equations~\eqref{e.dirichlet.pairs.proof} with~$u-v \in H^1_0(\cu_m)$ and subtracting yields
\begin{equation*}
	\fint_{\cu_m} (\nabla u - \nabla v) \cdot \tilde{\a}_{L,m} \nabla u
	=
	\fint_{\cu_m} \nabla v \cdot \shom_m (\nabla u - \nabla v)
	\,.
\end{equation*}
We therefore obtain
\begin{align*}
	\biggl| \fint_{\cu_m}
	\bigl(  \nabla u \cdot \a_L\nabla u - \nabla v \cdot \shom_m \nabla v \bigr)
	\biggr|
	&
	\leq
	\biggl|
	\fint_{\cu_m}
	\nabla v \cdot ( \tilde{\a}_{L,m} \nabla u -  \shom_m  \nabla v )
	\biggr|
	+
	\biggl|
	\fint_{\cu_m}
	\nabla v \cdot \shom_m ( \nabla u -  \nabla v )
	\biggr|
	\notag \\ &
	\leq
	\| \nabla v \|_{\underline{W}^{s,(4d)'}(\cu_m)}
	\| \tilde{\a}_{L,m} \nabla u -  \shom_m  \nabla v \|_{\underline{W}^{-s,4d}(\cu_m)}
	\notag \\ & \qquad
	+
	\shom_m \| \nabla v \|_{\underline{W}^{s,(4d)'}(\cu_m)}
	\| \nabla u -  \nabla v \|_{\Wminusul{-s}{4d}(\cu_m)}
	\,.
\end{align*}
The first pairing uses the full dual norm, while the second may use the hatted norm since~$(\nabla u-\nabla v)_{\cu_m}=0$.
By Schauder estimates for the Poisson equation in a cube (by odd reflection), we have
\begin{align*}
	\| \nabla v \|_{\underline{W}^{s,(4d)'}(\cu_m)}
	&
	\leq
	C 3^{-ms}
	\bigl(
	\shom_m^{-1}
	3^{\frac12 m}
	+
	3^{\frac12 m}
	\bigr)\,.
\end{align*}
Combining the previous two displays with~\eqref{e.weak.flux.bound.again} and~\eqref{e.weak.gradient.bound} yields~\eqref{e.intro.energies} for~$\a_L$, uniformly in~$L\geq m$.

\smallskip

It remains to remove the infrared cutoff. Subtracting the constant antisymmetric matrix~$(\k_L)_{\cu_m}$ does not change the equation or the energy, and the resulting coefficients converge uniformly on~$\cu_m$ to a representative of~$\a$. Since their common symmetric part is~$\nu\Id$, the corresponding Dirichlet solutions converge strongly in~$H^1(\cu_m)$ by testing the difference equation. The energy estimate therefore passes to the limit; after taking an almost-everywhere convergent subsequence, the uniform~$L^\infty$ estimate passes as well. This proves the theorem.
\end{proof}

\section{Superdiffusivity of the process}
\label{s.theprocess}

In this section, we complete the proof of Theorem~\ref{t.superdiffusivity} by passing from estimates on the infinitesimal generator of the process~\eqref{e.SDE} to estimates on the process itself. The main quantitative estimates we proved for the infinitesimal generator in the previous section are stated in terms of elliptic (time-independent) Dirichlet problems in large cubes. These estimates can be translated relatively easily into information about the first and second moments of the \emph{stopped} process~$X_{t\wedge\tau(\cu_m)}$.
Here~$\tau(U)$ denotes, for each open subset~$U \subseteq\Rd$, the time of the first exit from~$U$,
\begin{equation*}
	\tau(U)
	\coloneqq
	\inf \bigl\{ s \geq 0 \,:\, X_s \not \in U \bigr\}\,.
\end{equation*}
If we let~$u_1,u_2\in H^1(\cu_m)$ denote, respectively, the solutions of the Dirichlet problems
\begin{equation*}
	\left\{
	\begin{aligned}
		& -\nabla \cdot \a \nabla u_1 = 0 & \mbox{in} & \ \cu_m\,, \\
		& u_1 = e\cdot x & \mbox{on} & \ \partial \cu_m\,,
	\end{aligned}
	\right.
	\qquad
	\text{and}
	\qquad
	\left\{
	\begin{aligned}
		& -\nabla \cdot \a \nabla u_2 = -2\shom_m & \mbox{in} & \ \cu_m\,, \\
		& u_2 = (e \cdot x)^2 & \mbox{on} & \ \partial \cu_m\,,
	\end{aligned}
	\right.
\end{equation*}
then applying It\^o's formula gives us
\begin{equation*}
	\mathbf{E}^{\k,0}[u_1(X_{t \wedge \tau(\cu_m) })] = u_1(0)
	\,,
\end{equation*}
and
\begin{align}
	\mathbf{E}^{\k,0}[u_2(X_{t \wedge \tau(\cu_m)})]
	&
	=
	u_2(0) + 2\shom_m \mathbf{E}^{\k,0}[t \wedge \tau(\cu_m)]
	\notag \\ &
	=
	u_2(0) + 2\shom_m t
	-
	2\shom_m
	\mathbf{E}^{\k,0}\bigl[ ( t - \tau(\cu_m)) \indc_{\{\tau(\cu_m) < t \} } \bigr]
	\,.
	\label{e.u2.ito}
\end{align}
Theorem~\ref{t.homogenization} says that if we replace~$u_1$ by~$e\cdot x$ and~$u_2$ by~$(e\cdot x)^2$ in these expressions, then we make only a small relative error.

\smallskip

To get information on the original process~$\{ X_t\}$, we will choose the scale parameter~$m$ slightly larger than the length scale associated to~$t$, so that, with high probability, the process does not exit from~$\cu_m$ by time~$t$; that is,~$\mathbf{P}^{\k,0}[ \tau(\cu_m) \leq t ]$ is very small. This allows us to replace~$t \wedge \tau(\cu_m)$ by~$t$ in these expressions, after making another small relative error, which leads to the statement of Theorem~\ref{t.superdiffusivity}.
The remaining work therefore is to derive an upper bound estimate for the probability of an early exit from a cube. This appears in Proposition~\ref{p.early.exit.bound} in Section~\ref{ss.exit.time.bounds}. The proof of this estimate involves a \emph{chaining argument}, in which we demonstrate that an early exit from a big cube~$\cu_m$ would necessarily involve many early exits from smaller-scale cubes of the form~$z+\cu_n$. The probability of this unlikely event is controlled by a stopping time argument, using bounds for the probability of exiting these smaller cubes obtained from Theorem~\ref{t.homogenization}. We collect some preliminary estimates needed for this argument in Sections~\ref{ss.localization.in.L.infty} and~\ref{ss.percolation.on.paths}. The details of the proof of Theorem~\ref{t.superdiffusivity} are presented in Section~\ref{ss.proof.theorem.A}.

\smallskip

We introduce the time scale associated with a length scale~$r$ by
\begin{equation}
	\Tr(r) \coloneqq r^2 \bigl( \nu^2 + \cstar \cgamma^{-1} r^{2\cgamma} \bigr)^{-\nf12}
	=
	r^2 \nueff(r)^{-1}
	\,.
	\label{e.Tr.r.def}
\end{equation}
This is essentially (up to a small error) the inverse of the length scale associated to a given time scale~$t$, which was defined in~\eqref{e.Rt.def} as
\begin{equation}
	\Rt(t)
	\coloneqq
	\bigl(
	(\nu t)^{2-\cgamma}
	+
	\cstar \cgamma^{-1} t^2
	\bigr)^{\frac1{2(2-\cgamma)}}
	\,.
	\label{e.Rt.def.again}
\end{equation}
Throughout this section, unless explicitly stated to the contrary, we allow constants~$C$ to depend on a lower bound for the parameter~$\cstar$, in addition to~$d$. Explicit dependence on~$\cstar$ can be extracted, if desired, from the arguments. We also assume throughout this section that~$\cgamma \leq \cgamma_0(d,\cstar)$, the standing regime of Theorems~\ref{t.homogenization} and~\ref{t.regularity}, whose estimates are used at every scale.

\subsection{Localization of homogenization errors in~\texorpdfstring{$L^\infty$}{L-infty}}
\label{ss.localization.in.L.infty}

In this subsection we introduce local events controlling homogenization errors in~$L^\infty$ and a notion of ``good cube'' which ensures that homogenization and regularity estimates are valid, uniformly in the infrared cutoff, for arbitrary solutions in the cube.

\smallskip

For each~$m,n \in\Z$,~$y \in \Rd$,~$\g \in C^{0,\nf12}(y+\cu_n;\Rd)$, we let~$u_m(\cdot,y+\cu_n,\g)$ and~$v(\cdot,y+\cu_n,\g)$ denote, respectively, the solutions of the Dirichlet problems
\begin{equation}
	\left\{
	\begin{aligned}
		& -\nabla\cdot \a_m \nabla u_m(\cdot,y+\cu_n,\g) = \nabla \cdot \g
		& \mbox{in} & \ y+\cu_n\,, \\
		& u_m(\cdot,y+\cu_n,\g) = 0 & \mbox{on} & \ \partial (y+\cu_n)\,,
	\end{aligned}
	\right.
	\label{e.dirichlet.pairs.Sfive.u}
\end{equation}
and
\begin{equation}
	\left\{
	\begin{aligned}
		& -\shom_n \Delta v(\cdot,y+\cu_n,\g) = \nabla \cdot \g
		& \mbox{in} & \ y+\cu_n \,, \\
		& v(\cdot,y+\cu_n,\g) = 0 & \mbox{on} & \ \partial (y+\cu_n)\,,
	\end{aligned}
	\right.
	\label{e.dirichlet.pairs.Sfive.v}
\end{equation}
We define the random variables
\begin{multline}
	E(y+\cu_n) \coloneqq
	\sup
	\Bigl \{
	\shom_n
	3^{-n} \| u_n(\cdot,y{+}\cu_n,\g) - v(\cdot,y{+}\cu_n,\g) \|_{L^\infty(y+\cu_n)}
	\\
	\,:\,
	\g \in C^{0,\nf12}(y{+}\cu_n;\Rd) \,, \
	3^{\frac12n} [ \g ]_{C^{0,\nf12}(y+\cu_n)}
	\leq 1
	\Bigr\}
	\label{e.error.rv.Sfive}
\end{multline}
and
\begin{multline}
	X(y+\cu_n)
	\\
	\coloneqq
	\sup
	\Bigl \{
	\shom_n 3^{-\frac12 n}
	[ u_n(\cdot,y{+}\cu_n,\g) ]_{C^{0,\nf12}(y+\cu_n)}
	\,:\,
	\g \in C^{0,\nf12}(y{+}\cu_n;\Rd) \,, \
	3^{\frac12n} [ \g ]_{C^{0,\nf12}(y+\cu_n)}
	\leq 1
	\Bigr\}
	\,.
	\label{e.regularity.rv.Sfive}
\end{multline}
In the proof of Theorem~\ref{t.homogenization} in Section~\ref{ss.proof.theorem.B}, we proved that, under the assumptions of that theorem, for every~$p \in[1, C^{-1} \cgamma^{-1} \left| \log\cgamma \right|^{-1} ]$,
\begin{equation}
	\E \bigl[ E(y+\cu_n) ^p \bigr]^{\nf1p}
	\leq C \bigl(p^{\nf12}+\left|\log\cgamma\right|^{\nf12}\bigr)
	\cgamma^{\nf12} \left| \log\cgamma \right|^3
	\,.
	\label{e.error.rv.Sfive.moment.bound}
\end{equation}
In the proof of Theorem~\ref{t.regularity} in Section~\ref{ss.proof.regularity}, we proved, under the assumptions of that theorem,
\begin{equation}
	\E \bigl[ X(y+\cu_n)^{p} \bigr]^{\nf1p} \leq C \,,\quad \mbox{where} \ p = c\cgamma^{-1} \left| \log\cgamma\right|^{-6}
	\,.
	\label{e.reg.rv.Sfive.moment.bound}
\end{equation}
We also introduce the good event defined for each~$\ep>0$ by
\begin{equation}
	\mathcal{G}(y+\cu_n,\ep)
	\coloneqq
	\bigl\{
	E(y+\cu_n) \leq \ep
	\bigr\}
	\cap
	\bigl\{ X(y+\cu_n) \leq 2C_{\eqref{e.reg.rv.Sfive.moment.bound}} \bigr\}
	\,.
	\label{e.event.G}
\end{equation}
Consequently, by Markov's inequality, for every~$p \in[1, C^{-1} \cgamma^{-1} \left| \log\cgamma \right|^{-1} ]$,
\begin{equation*}
	\P \bigl[ \mathcal{G}(y+\cu_n,\ep) \bigr]
	\geq
	1 -
	\bigl( C \ep^{-1} (p^{\nf12}+\left|\log\cgamma\right|^{\nf12})
	\cgamma^{\nf12} \left| \log \cgamma \right|^3 \bigr)^p
	-
	2^{-c\cgamma^{-1} \left| \log\cgamma\right|^{-6}}
	\,.
\end{equation*}
Restricting to~$\ep\in [ C\cgamma^{\nf12} \left| \log \cgamma \right|^{\nf72}, \nf14]$ and then optimizing the previous estimate in~$p$, we arrive at the bound
\begin{equation*}
	\P \bigl[ \mathcal{G}(y+\cu_n,\ep) \bigr]
	\geq
	1 -
	\exp \bigl( - c \ep^2 \cgamma^{-1} \left| \log \cgamma \right|^{-6} \bigr)
	\,.
\end{equation*}
The point of introducing~$E(y+\cu_n)$ and~$X(y+\cu_n)$, rather than just working with the random variables introduced in the arguments in Sections~\ref{ss.proof.regularity} and~\ref{ss.proof.theorem.B}, is that~$E(y+\cu_n)$ and~$X(y+\cu_n)$ are evidently local: they depend only on~$\a_n$ restricted to~$y+\cu_n$. As a result, we have that, for every pair of subsets~$Y,Z \subseteq \Rd$ and~$\ep,\delta>0$,
\begin{align}
	\lefteqn{
	\min_{y \in Y, z\in Z} |y-z|_\infty \geq (2+\sqrt d)\,3^{n}
	} \quad &
	\notag \\ &
	\implies
	\text{$\{ (E,X) (y+\cu_n) \,:\, y\in Y\}$ \ and \ $\{ (E,X)(z+\cu_n) \,:\, z\in Z\}$ are $\P$-independent}
	\label{e.FRD.Sfive} \\ &
	\implies
	\text{$\{ \mathcal{G}(y+\cu_n,\ep) \,:\, y\in Y\}$ \ and \ $\{ \mathcal{G}(z+\cu_n,\delta) \,:\, z\in Z\}$ are $\P$-independent.}
	\label{e.FRD.events.Sfive}
\end{align}
Here~$| x |_\infty \coloneqq \max_{i\in\{ 1,\ldots,d\}} |e_i\cdot x|$. We also need to define a ``large-scale'' good event (cf.~Definition~\ref{d.good.event.for.lambda})
\begin{equation*}
	\mathcal{J}(y+\cu_n,\theta )
	\coloneqq
	\biggl\{
	\sum_{k=n}^\infty
	3^{(2-\cgamma)n}
	\|  \nabla \mathbf{j}_k  \|_{\underline{W}^{1,\infty}(y+\cu_n)}
	\leq
	\theta
	\biggr\}
	\,.
\end{equation*}
This event is not local, since it depends on~$\mathbf{j}_k$ for large~$k$; however, as we will show (and as we have already seen earlier in the paper), it is approximately local since any contribution of~$\mathbf{j}_k$ for large~$k$ is very unlikely.

\smallskip

We defined the random variables in~\eqref{e.error.rv.Sfive} and~\eqref{e.regularity.rv.Sfive} using the field~$\a_n$ in order that these random variables, as well as the corresponding good events in~\eqref{e.event.G}, are local. This gives us enough independence that we can obtain the percolation estimate presented in the next subsection. However, we need that the estimates are still valid with~$\a_m$ for arbitrary~$m\geq n$. In the next proposition we show, by treating~$\a_m - \a_n$ as a perturbation and using an iteration argument, that we can indeed replace~$\a_n$ by~$\a_m$ for any~$m >n$ and have essentially the same estimates.

\begin{proposition}
\label{p.injection.in.L.infty}
There exists a constant~$C(d,\cstar)<\infty$ such that, for every~$\ep \in (0,\nf14]$,~$y\in\Rd$,~$n\in\Z$ and~$\g \in C^{0,\nf12}(y+\cu_n;\Rd)$,
\begin{multline}
	\sup_{m \geq n}
	\shom_n
	3^{-n} \| u_m(\cdot,y+\cu_n,\g) - v(\cdot,y+\cu_n,\g) \|_{L^\infty(y+\cu_n)}
	\indc \bigl\{ \mathcal{G}(y+\cu_n,\ep) \cap \mathcal{J}(y+\cu_n,C^{-1}\ep\cgamma^{-\nf12} ) \bigr\}
	\\
	\leq
	C \ep 3^{\frac12n} [ \g ]_{C^{0,\nf12}(y+\cu_n)}
	\,.
	\label{e.injection.L.infty}
\end{multline}
and
\begin{multline}
	\sup_{m \geq n}
	\shom_n
	3^{-\frac12 n}
	[ u_m(\cdot,y+\cu_n,\g) ]_{C^{0,\nf12}(y+\cu_n)}
	\indc \bigl\{ \mathcal{G}(y+\cu_n,\ep) \cap \mathcal{J}(y+\cu_n,C^{-1} \ep\cgamma^{-\nf12} ) \bigr\}
	\\
	\leq
	C 3^{\frac12n} [ \g ]_{C^{0,\nf12}(y+\cu_n)}\,.
	\label{e.injection.regularity}
\end{multline}
\end{proposition}
\begin{proof}
We fix~$m\geq n$ and assume~$y=0$ and~$\indc_{\mathcal{G}(\cu_n,\ep)} \indc_{\mathcal{J}(\cu_n,\theta)} = 1$ where~$\theta$ will be selected below.
We compute~$u_m(\cdot,\cu_n,\g)$ by an iteration scheme.

\smallskip

We recursively define~$w_0 \coloneqq u_n(\cdot,\cu_n,\g)$ and, for each~$j\in\N$, we define~$h_j \in H^1_0(\cu_n)$ as the solution of
\begin{equation}
	\left\{
	\begin{aligned}
		& -\Delta h_j =
		-\nabla \cdot \bigl( w_j \nabla \cdot (\k_m-\k_n) \bigr)
		& \mbox{in} & \ \cu_n\,, \\
		& h_j  = 0 & \mbox{on} & \ \partial \cu_n\,,
	\end{aligned}
	\right.
	\label{e.iteration.scheme}
\end{equation}
and then set
\begin{equation*}
	w_{j+1} \coloneqq u_n \bigl(\cdot,\cu_n, \nabla h_j \bigr)
	\,.
\end{equation*}
Using that~$\k_m-\k_n$ is anti-symmetric, we find that
\begin{equation*}
	\nabla \cdot \bigl( w_j \nabla \cdot (\k_m-\k_n) \bigr)
	=
	\nabla w_j \cdot ( \nabla \cdot (\k_m-\k_n))
	=
	\nabla \cdot \bigl( (\k_m-\k_n) \nabla w_j \bigr)\,.
\end{equation*}
Hence~$\nabla \cdot \nabla h_j = \nabla \cdot \bigl( (\k_m-\k_n) \nabla w_j \bigr)$ and therefore
\begin{equation*}
	w_{j+1} = u_n \bigl(\cdot,\cu_n,(\k_m-\k_n) \nabla w_j \bigr)
	\,.
\end{equation*}
We first show that~$\sum_j w_j$ converges, and then identify its sum.
By the definition of the good event~$\mathcal{G}(\cu_n,\ep)$ and~$w_{j+1}$, we have that
\begin{equation*}
	\shom_n 3^{-\frac12 n}
	[ w_{j+1} ]_{C^{0,\nf12}(\cu_n)}
	\leq
	C 3^{\frac12 n} [ \nabla h_j ]_{C^{0,\nf12}(\cu_n)}
	\,.
\end{equation*}
By the Schauder estimates for the Laplacian in a cube (by odd reflection) and the definition of the good event~$\mathcal{J}(\cu_n,\theta)$, we have
\begin{align*}
	[ \nabla h_j ]_{C^{0,\nf12}(\cu_n)}
	&
	\leq
	C [ w_j \nabla \cdot (\k_m-\k_n)  ]_{C^{0,\nf12}(\cu_n)}
	\notag \\ &
	\leq
	C [ w_j  ]_{C^{0,\nf12}(\cu_n)}
	\| \nabla (\k_m-\k_n) \|_{L^\infty(\cu_n)}
	+
	C [ \nabla (\k_m-\k_n) ]_{C^{0,\nf12}(\cu_n)}
	\| w_j  \|_{L^\infty(\cu_n)}
	\notag \\ &
	\leq
	C\theta
	3^{(\cgamma-1)n}
	[ w_j  ]_{C^{0,\nf12}(\cu_n)}
	+
	C\theta
	3^{(\cgamma-\frac32)n}
	\| w_j  \|_{L^\infty(\cu_n)}
	\notag \\ &
	\leq
	C\theta
	3^{(\cgamma-1)n}
	[ w_j  ]_{C^{0,\nf12}(\cu_n)}
	\,.
\end{align*}
Combining these, we deduce that
\begin{equation*}
	[ w_{j+1} ]_{C^{0,\nf12}(\cu_n)}
	\leq
	C \shom_n^{-1}  3^{n} [ \nabla h_j ]_{C^{0,\nf12}(\cu_n)}
	\leq
	C\theta
	\shom_n^{-1}
	3^{\cgamma n}
	[ w_j  ]_{C^{0,\nf12}(\cu_n)}
	\leq
	C \theta \cgamma^{\nf12}
	[ w_j  ]_{C^{0,\nf12}(\cu_n)}
	\,.
\end{equation*}
Here we used~\eqref{e.shom.m.bounds.thmB}, which gives~$\shom_n^{-1} 3^{\cgamma n} \leq C \cgamma^{\nf12}$.
Imposing the requirement that~$\theta < c \ep \cgamma^{-\nf12}$, we deduce that
\begin{equation*}
	[ w_{j+1} ]_{C^{0,\nf12}(\cu_n)}
	\leq
	\ep [ w_j  ]_{C^{0,\nf12}(\cu_n)}
	\,.
\end{equation*}
Iterating this inequality, we obtain, for every~$j \in\N$,
\begin{equation*}
	[ w_{j} ]_{C^{0,\nf12}(\cu_n)}
	\leq
	\ep^j
	\cdot
	[ u_n(\cdot,\cu_n,\g) ]_{C^{0,\nf12}(\cu_n)}
	\leq
	C
	\ep^j
	\shom_n^{-1}
	3^n [ \g ]_{C^{0,\nf12}(\cu_n)}
	\,.
\end{equation*}
Since every~$w_j$ vanishes on~$\partial\cu_n$, the preceding estimate also controls~$\|w_j\|_{L^\infty(\cu_n)}$ geometrically. Testing the equations for~$h_j$ and~$w_{j+1}$ gives
\begin{equation*}
	\nu\|\nabla w_{j+1}\|_{L^2(\cu_n)}
	\leq\|\nabla h_j\|_{L^2(\cu_n)}
	\leq\|\nabla\cdot(\k_m-\k_n)\|_{L^\infty(\cu_n)}
	\|w_j\|_{L^2(\cu_n)}
	\leq C\ep^j
	\,.
\end{equation*}
Here the last constant may depend on the fixed~$m,n,\g$ and realization. Thus~$\sum_jw_j$ converges in~$H^1_0(\cu_n)$ and uniformly. Summing the weak equations and passing to the limit shows that its sum solves the equation for~$u_m$; uniqueness gives
\begin{equation}
	u_m(\cdot,\cu_n,\g)
	=
	\sum_{j=0}^\infty w_j\,.
	\label{e.u.m.iteration.scheme}
\end{equation}
Summing over~$j$ yields, in view of~\eqref{e.u.m.iteration.scheme},
\begin{equation*}
	[ u_m(\cdot,\cu_n,\g) - u_n(\cdot,\cu_n,\g)  ]_{C^{0,\nf12}(\cu_n)}
	\leq
	C \ep
	\shom_n^{-1}
	3^n [ \g ]_{C^{0,\nf12}(\cu_n)}
	\,.
\end{equation*}
In particular,
\begin{equation*}
	\|  u_m(\cdot,\cu_n,\g) - u_n(\cdot,\cu_n,\g) \|_{C^{0,\nf12}(\cu_n)}
	\leq
	C \ep
	\shom_n^{-1}
	3^{\frac32 n} [ \g ]_{C^{0,\nf12}(\cu_n)}
	\,.
\end{equation*}
Combining the last two displays with the definition of the good event~$\mathcal{G}(\cu_n,\ep)$ and the triangle inequality, we complete the proof.
\end{proof}

\subsection{Percolation estimates for bad cubes}
\label{ss.percolation.on.paths}

To prepare for the stopping time/chaining argument, we present a percolation-type estimate which says that any chain of cubes which are translations of~$\cu_n$ which connects~$\cu_{m-1}$ to~$\partial \cu_m$ must contain at least~$c3^{m-n}$ good cubes, where the notion of ``good'' event here is defined by
\begin{equation}
	\mathcal{Q}(z+\cu_n,\ep)
	\coloneqq
	\mathcal{G}(z+\cu_n,\ep) \cap \mathcal{J}(z+\cu_n,C_{\mathrm{Prop\, \ref{p.injection.in.L.infty}}} ^{-1}\ep\cgamma^{-\nf12} )
	\,, \qquad n\in\Z\,, \ z \in\Rd\,, \ep\in (0,\nf14]\,.
	\label{e.event.Q}
\end{equation}

\begin{proposition}[Chains of good cubes]
\label{p.percolation.GF}
There exist~$c(d,\cstar)>0$ and~$C(d,\cstar) <\infty$ such that, for every~$\ep\in [C\cgamma^{\nf12}|\log\cgamma|^{\nf72}, \nf14]$ and~$m\in\Z$, there exists a random variable~$Y_m(\ep)$ satisfying, for every~$N\in\N$,
\begin{equation}
	\P \bigl[ Y_m(\ep) > N \bigr]
	\leq \exp \bigl( - c \ep^2 \cgamma^{-1} |\log\cgamma|^{-6} \cdot 3^N \bigr)
	\label{e.Ym.integrability}
\end{equation}
such that, for every~$n\in\Z$ with~$n \leq m - Y_m(\ep)$ and every path~$\Gamma\subseteq 3^{n-1} \Zd$ from~$3^{n-1} \Zd \cap \cu_{m-1}$ to~$3^{n-1} \Zd \setminus \cu_{m}$, we have
\begin{equation*}
	\sum_{z \in \Gamma}
	\indc_{\mathcal{Q}(z+\cu_n,\ep)}
	\geq
	\frac34 \cdot 3^{m-n}\,.
\end{equation*}
\end{proposition}
\begin{proof}
We apply Lemma~\ref{l.percolation.bound.general} with parameters~$a \coloneqq \frac74 - 2\cgamma$ and~$T^2 \coloneqq c\ep^2 \cgamma^{-1} |\log\cgamma|^{-6}$, after rescaling so that~$3^{n-1}\Zd$ becomes~$\Zd$. For~$z \in \Zd$, we define
\begin{equation*}
	\begin{aligned}
		\mathrm{B}_0(z) &\coloneqq \mathcal{G}(3^{n-1}z+\cu_n,\ep)^c \cup \Bigl\{ 3^{(2-\cgamma)n} \| \nabla \mathbf{j}_{n} \|_{\underline{W}^{1,\infty}(3^{n-1}z+\cu_n)} > C_0 \ep \cgamma^{-\nf12} \Bigr\}\,, \\[4pt]
		\mathrm{B}_L(z) &\coloneqq \Bigl\{ 3^{(2-\cgamma)n} \| \nabla \mathbf{j}_{n+L} \|_{\underline{W}^{1,\infty}(3^{n-1}z+\cu_n)} > C_0 \ep \cgamma^{-\nf12} \cdot 3^{-\frac18 L} \Bigr\}\,, \quad L \geq 1\,.
	\end{aligned}
\end{equation*}
Since~$\sum_{L=0}^\infty 3^{-\frac18 L}<\infty$, we may choose~$C_0(d,\cstar)$ small enough that~$\mathcal{Q}(3^{n-1}z+\cu_n,\ep)^c \subseteq \mathrm{B}(z) \coloneqq \bigcup_{L \in \N_0} \mathrm{B}_L(z)$ by~\eqref{e.event.Q}. By scaling, the normalized gradient and Hessian parts have sizes~$3^{-(1-\cgamma)L}$ and~$3^{-(2-\cgamma)L}$, respectively. Thus~\eqref{e.error.rv.Sfive.moment.bound},~\eqref{e.reg.rv.Sfive.moment.bound} and~\ref{a.j.reg} give
\begin{equation*}
	\P[\mathrm{B}_0(z)] \leq \exp\bigl(-c \ep^2 \cgamma^{-1} |\log\cgamma|^{-6}\bigr)
	\qquad \text{and} \qquad
	\P[\mathrm{B}_L(z)] \leq \exp\bigl(-c \ep^2 \cgamma^{-1} \cdot 3^{(\frac74-2\cgamma)L}\bigr)\,.
\end{equation*}
A fixed dimension-dependent shift of the level index gives the finite-range hypothesis by~\ref{a.j.frd} and~\eqref{e.FRD.events.Sfive}, changing~$T$ only by a dimensional factor. Since~$\frac54\leq a<\frac74<d$ and~$a - 1 \geq \frac14$, the condition~$T \geq C\lambda^{-1}(a-1)^{-1}$ with~$\lambda = \nf14$ is satisfied by the assumption~$\ep \geq C \cgamma^{\nf12} |\log\cgamma|^{\nf72}$. The path bound~\eqref{e.path.bound} at scale~$k = m-n$ then yields
\begin{equation}
	\label{e.path.bound.applied}
	\P\Bigl[\exists\ \text{path } \Gamma \text{ from } \Zd \cap \cu_k \text{ to } \Zd \setminus \cu_{k+1} \text{ with } \sum_{z \in \Gamma} \indc_{\mathrm{B}^c(z)} < \tfrac34 \cdot 3^k \Bigr]
	\leq \exp\bigl(-c \ep^2 \cgamma^{-1} |\log\cgamma|^{-6} \cdot 3^{k}\bigr)\,.
\end{equation}
For each~$k$, we read the event in~\eqref{e.path.bound.applied} with~$n=m-k$, and define~$Y_m(\ep)$ as the smallest~$N \in \N$ such that this event fails to occur for every~$k \geq N$. A union bound gives
\begin{equation*}
	\P\bigl[Y_m(\ep) > N\bigr]
	\leq \sum_{k=N}^\infty \exp\bigl(-c \ep^2 \cgamma^{-1} |\log\cgamma|^{-6} \cdot 3^k\bigr) \leq \exp(-c \ep^2 \cgamma^{-1} |\log\cgamma|^{-6} 3^N)\,,
\end{equation*}
which is~\eqref{e.Ym.integrability}.
\end{proof}

\subsection{Tail bounds on the displacement}
\label{ss.exit.time.bounds}

Theorems~\ref{t.homogenization} and~\ref{t.regularity} provide quantitative control of elliptic boundary value problems for the generator~$L$ on cubes at all relevant scales, in the form of~$L^\infty$-accurate comparisons with a scale-dependent Laplacian. The purpose of the next proposition is to convert this finite-volume elliptic information into genuinely parabolic control in the form of an off-diagonal estimate for the parabolic Green function. We state this, equivalently, as an estimate for the probability of an early exit of the process from a cube.

\begin{proposition}[Displacement tail bound]
\label{p.early.exit.bound}
There exists~$C(d,\cstar)<\infty$ such that, for each~$m\in\Z$, there is a random variable~$S_m$ satisfying
\begin{equation}
	\E \bigl[ S_m^{p} \bigr]^{\nf1p }
	\leq
	C\,, \quad \mbox{for} \  p = C^{-1} \cgamma^{-1} \left| \log \cgamma\right|^{-6} \,,
	\label{e.estimate.St}
\end{equation}
and such that, for every~$t \in (0,\infty)$ satisfying
\begin{equation}
	\min \biggl\{
	\frac{3^{2m}}{\nu t}
	\,,\,
	\bigl( \cstar^{-1} \cgamma  \bigr)^{\nf{1}{2}}
	\Bigl(
	\frac{(3^m)^{2-\cgamma}}{t}
	\Bigr)^{\frac{1}{1-\cgamma} }
	\biggr\}
	\geq
	S_m\,,
	\label{e.tail.displacement.minscale.condition}
\end{equation}
we have the estimate
\begin{equation}
	\sup_{x\in \cu_{m-1}}
	\mathbf{P}^{\k,x}
	\bigl[
	\tau(\cu_m) \leq
	t
	\bigr]
	\leq
	\exp
	\biggl(
	-\frac1{C} \min \biggl\{
	\frac{3^{2m}}{\nu t}
	\,,\,
	\bigl( \cstar^{-1} \cgamma  \bigr)^{\nf{1}{2}}
	\Bigl(
	\frac{(3^m)^{2-\cgamma}}{t}
	\Bigr)^{\frac{1}{1-\cgamma} }
	\biggr\}
	\biggr)
	\,.
	\label{e.early.exit.probability}
\end{equation}
\end{proposition}

Observe that, since
\begin{equation*}
	\bigl\{ |X_t| \geq \sqrt{d}3^m \bigr\}
	\subseteq \bigl\{ \tau(\cu_m) \leq t \bigr\} \,,
\end{equation*}
the estimate~\eqref{e.early.exit.probability} implies the following tail estimate for the displacement for every~$t\in (0,\infty)$ and~$m\in\Z$ satisfying~\eqref{e.tail.displacement.minscale.condition}:
\begin{equation}
	\sup_{x\in \cu_{m-1}}
	\mathbf{P}^{\k,x} \bigl[ |X_t| > \sqrt d\,3^m \bigr]
	\leq
	\exp\biggl(
	-\,\frac1C \min \biggl\{
	\frac{3^{2m}}{\nu t}
	\,,\,
	\bigl( \cstar^{-1} \cgamma  \bigr)^{\nf{1}{2}}
	\Bigl(
	\frac{(3^m)^{2-\cgamma}}{t}
	\Bigr)^{\frac{1}{1-\cgamma} }
	\biggr\}
	\biggr)
	\,.
	\label{e.displacement.tail.estimate}
\end{equation}
This is an integrated form of an off-diagonal bound on the Green function for the parabolic operator~$\partial_t - \nabla \cdot \a\nabla$. We expect that this upper bound estimate is sharp, up to the constant~$C$, since it matches the functional form of transition densities for random walks with anomalous walk dimension (see for example~\cite{Barlow1998,GrigoryanTelcs2012}).

\begin{remark}
\label{r.higher.moments.of.displacement}
By integrating the displacement tail estimate~\eqref{e.displacement.tail.estimate}, using~\eqref{e.estimate.St}, we obtain an upper bound for higher moments of the displacement, as announced in~\eqref{e.moment.bounds}: for a constant~$C(d,\cstar) <\infty$,
\begin{equation}
	\label{e.moment.bounds.proof}
	\E\bigl[\mathbf{E}^{\k,0}[|X_t|^p]\bigr]^{\nf1p} \leq C p^{\nf12}\Rt(t)
	\qquad \text{for } p \in [2, C^{-1} \cgamma^{-1} |\log\cgamma|^{-6}]\,.
\end{equation}
The matching lower bound is provided by the statement of Theorem~\ref{t.superdiffusivity}.
\end{remark}

\smallskip

In the classical uniformly parabolic setting, the passage from near-diagonal information to global off-diagonal bounds is achieved by chaining local estimates along a path, an idea already present in Aronson's work on Gaussian bounds for fundamental solutions~\cite{Aronson1967,Aronson1968}. To prove Proposition~\ref{p.early.exit.bound} we use a probabilistic implementation of the same principle due to Barlow~\cite{Barlow1998} (see also~\cite{GrigoryanTelcs2012}): we first bound the Laplace transform of the exit time from a single (small) cube, and then iterate this contraction along a chain of translated cubes using stopping times and the strong Markov property. We crucially use Proposition~\ref{p.percolation.GF}, which ensures that we can find enough cubes with good exit time bounds along any such chain.

\smallskip

For every bounded convex domain~$U \subseteq \Rd$, the expected exit time from~$U$ for the diffusion process starting at~$x\in U$, denoted by
\begin{equation}
	w(x) \coloneqq \mathbf{E}^{\k,x}[\tau(U)] \,, \quad x\in U
	\,,
	\label{e.exit.time.characterization}
\end{equation}
is the solution of the Dirichlet problem
\begin{equation}
	\label{e.exit.time.pde}
	\left\{
	\begin{aligned}
		& -\nabla \cdot \a \nabla w = 1 & \mbox{in} & \ U\,, \\
		& w = 0 & \mbox{on} & \ \partial U\,.
	\end{aligned}
	\right.
\end{equation}
This provides the link between local elliptic problems and stopping times.

\begin{proof}[{Proof of Proposition~\ref{p.early.exit.bound}}]
Fix~$m\in\Z$ and~$t \in (0,\infty)$.
Let~$\ep,\delta \in (0,\nf14]$ be parameters to be selected below (as~$\ep(d)$ and~$\delta(d,\cstar)$). Define
\begin{equation}
	S_m \coloneqq C(d) \delta^{-\frac{1}{1-\cgamma}} 3^{Y_m(\ep)}\,,
	\label{e.Sm.displ.def}
\end{equation}
where~$Y_m(\ep)$ is as in Proposition~\ref{p.percolation.GF} and~$C(d)$ is taken larger than~$3^{d+4}$. The estimate~\eqref{e.estimate.St} is immediate from~\eqref{e.Ym.integrability}, provided that we eventually take~$\delta=c(d,\cstar)>0$.

We turn our attention to the proof of~\eqref{e.early.exit.probability}. Take~$n\in\Z$ to be the scale parameter satisfying
\begin{equation*}
	3^n \leq
	\max \biggl\{ \frac{\nu t}{ \delta 3^m}
	\,,\,
	(\cstar \cgamma^{-1} )^{\nf12}
	\Bigl(
	\frac{t}{\delta 3^m}
	\Bigr)^{\frac1{1-\cgamma}}
	\biggr\}
	< 3^{n+1}
	\,.
\end{equation*}
With~$\mathcal{Q}$ defined as in Proposition~\ref{p.percolation.GF}, define
\begin{equation*}
	Q_{m,n} \coloneqq \bigl\{ z \in 3^{n-1} \Zd \,:\,
	z+\cu_{n} \subseteq \cu_m\setminus \cu_{m-1}\,, \ \indc_{\mathcal{Q}(z+\cu_{n},\ep)} = 1 \bigr\}\,.
\end{equation*}
Observe that
\begin{equation}
	\delta^{-1} t \Tr(3^n)^{-1}
	\asymp
	3^{m-n}
	\asymp
	\delta
	\min \biggl\{
	\frac{3^{2m}}{\nu t}
	\,,\,
	\delta^{\frac{\cgamma}{1-\cgamma}}
	\bigl( \cstar^{-1} \cgamma  \bigr)^{\nf{1}{2}}
	\Bigl(
	\frac{(3^m)^{2-\cgamma}}{t}
	\Bigr)^{\frac{1}{1-\cgamma} }
	\biggr\}
	\,.
	\label{e.scale.sep.gives}
\end{equation}
In particular, in view of the definition~\eqref{e.Sm.displ.def}, the condition~\eqref{e.tail.displacement.minscale.condition} implies~$n \leq m - Y_m(\ep)$. The passage from the lattice statement of Proposition~\ref{p.percolation.GF} to the geometric one below costs a dimensional factor: reading a continuous path from~$\cu_{m-1}$ to~$\partial \cu_m$ as a path in~$3^{n-1}\Zd$ and correcting its two endpoints costs at most~$2d$ sites; discarding the sites within~$3^n$ of~$\partial \cu_{m-1}$ and of~$\partial \cu_m$ removes the two boundary layers; and keeping every third site of~$3^{n-1}\Zd$ in each coordinate, which makes the cubes~$y_i+\cu_n$ pairwise disjoint, costs a further~$3^d$. These losses are absorbed into~$c$. This allows us to apply Proposition~\ref{p.percolation.GF} to obtain that
\begin{equation}
	\left\{
	\begin{aligned}
		& \text{there exist~$N \geq c 3^{m-n}$ and, for every continuous path~$\Gamma$ from~$\cu_{m-1}$ to~$\partial \cu_{m}$,}
		\\ &
		\text{a finite sequence~$\{ y_1,\ldots,y_N \} \subseteq Q_{m,n}$ satisfying: (i) $\Gamma \cap (y_i + \cu_{n-1}) \neq \emptyset$ for }
		\\ &
		\text{every~$i\in\{1,\ldots,N\}$; and (ii) the cubes~$y_i + \cu_n$ are pairwise disjoint. }
	\end{aligned}
	\right.
	\label{e.what.percolation.gives}
\end{equation}
We will use~\eqref{e.what.percolation.gives} in our stopping time argument in Step~4.

\smallskip

\emph{Step 1.}
We show that, if~$\ep$ is sufficiently small, then there exists~$C(d,\cstar)<\infty$ such that, for every~$z\in\Rd$,
\begin{equation}
	C^{-1} \Tr(3^n)
	\indc_{\mathcal{Q}(z+\cu_n,\ep) }
	\leq
	\inf_{x\in z+ \overline{\cu}_{n-1}}
	\mathbf{E}^{\k,x}
	\bigl[ \tau(z+\cu_n) \bigr]
	\label{e.exit.time.lower.bound}
\end{equation}
and
\begin{equation}
	\sup_{x \in z+\cu_n}
	\mathbf{E}^{\k,x}
	\biggl[ \exp\biggl(
	\frac{ \tau(z+\cu_n) }{{C\Tr(3^n)}}
	\biggr) \biggr]
	\indc_{\mathcal{Q}(z+\cu_n,\ep) }
	\leq 2
	\,.
	\label{e.exit.time.exp.moment}
\end{equation}
Assume~$\indc_{\mathcal{Q}(z+\cu_n,\ep) } = 1$. Denote
\begin{equation*}
	w(x) \coloneqq \mathbf{E}^{\k,x}
	\bigl[ \tau(z+\cu_n) \bigr]
	\,,
\end{equation*}
which satisfies~\eqref{e.exit.time.pde} with~$U = z+\cu_n$. To compare with~\eqref{e.exit.time.pde}, we let~$\tilde{w}$ denote the solution of the corresponding  Dirichlet problem for the deterministic equation,
\begin{equation*}
	\left\{
	\begin{aligned}
		& -\shom_n \Delta \tilde{w} = 1 & \mbox{in} & \ z+\cu_n\,, \\
		& \tilde{w} = 0 & \mbox{on} & \ \partial (z+ \cu_n)\,.
	\end{aligned}
	\right.
\end{equation*}
Writing~$\tilde{w}(z+3^n x) = \shom_n^{-1} 3^{2n} W(x)$, where~$W$ solves~$-\Delta W = 1$ in~$\cu_0$ with zero boundary values, we observe, since~$W$ is positive in~$\cu_0$ and therefore bounded below on~$\overline{\cu}_{-1}$, that~$\tilde{w}$ satisfies the estimate
\begin{equation*}
	C^{-1}
	\| \tilde{w} \|_{L^\infty(z+\cu_n)}
	\leq
	\shom_n^{-1} 3^{2n}
	\leq
	C \inf_{x \in z+ \overline{\cu}_{n-1}}
	\tilde{w} (x)
	\,.
\end{equation*}
In view of~\eqref{e.shom.m.bounds.thmB} and~\eqref{e.Tr.r.def},
this estimate can be written equivalently as
\begin{equation}
	\label{e.homog.exit.time.bounds}
	C^{-1}
	\| \tilde{w} \|_{L^\infty(z+\cu_n)}
	\leq
	\Tr(3^n)
	\leq
	C \inf_{x \in z+ \overline{\cu}_{n-1}}
	\tilde{w} (x)
	\,.
\end{equation}
Letting~$m\to\infty$ in Proposition~\ref{p.injection.in.L.infty}, by the cutoff-removal argument at the end of the proof of Theorem~\ref{t.homogenization}, gives the same estimates for the full coefficient~$\a$.
Write~$1 = \nabla \cdot \g$ with~$\g(x) \coloneqq x_1 e_1$, so that the definitions of~$\mathcal{Q}(z+\cu_n,\ep)$ in~\eqref{e.event.Q} and Proposition~\ref{p.injection.in.L.infty} yield
\begin{equation}
	\label{e.exit.time.homog.error}
	3^{-n} \| w - \tilde{w} \|_{L^\infty(z+\cu_n)}
	\leq
	C \ep \shom_n^{-1} 3^{\frac12 n}  [\g]_{\underline{W}^{\nf12,\infty}(z+\cu_n)}\,.
\end{equation}
For the linear vector field~$\g(x) = x_1 e_1$, we have~$[\g]_{\underline{W}^{\nf12,\infty}(z+\cu_n)} \leq C 3^{\frac12 n} $. Rearranging the previous inequality then yields
\begin{equation}
	\label{e.exit.time.Linfty.error}
	\| w - \tilde{w} \|_{L^\infty(z+\cu_n)} \leq C \ep \Tr(3^n)\,.
\end{equation}
Restricting to~$\ep\leq c$ for a sufficiently small constant~$c(d)>0$, we deduce from~\eqref{e.homog.exit.time.bounds} and~\eqref{e.exit.time.Linfty.error} that
\begin{equation}
	\label{e.exit.time.comparison}
	\frac12 \tilde{w}(x) \leq w(x) \leq 2 \tilde{w}(x) \,, \qquad \forall x \in z+\overline{\cu}_{n-1}\,.
\end{equation}
Combining~\eqref{e.homog.exit.time.bounds} and~\eqref{e.exit.time.comparison} yields~\eqref{e.exit.time.lower.bound}. From~\eqref{e.homog.exit.time.bounds} and~\eqref{e.exit.time.Linfty.error} we deduce the upper bound
\begin{equation}
	\sup_{x\in z+ \cu_{n}}
	\mathbf{E}^{\k,x}
	\bigl[ \tau(z+\cu_n) \bigr]
	\leq
	C \Tr(3^n)
	\,.
	\label{e.exit.time.upper.moment}
\end{equation}
The exponential bound~\eqref{e.exit.time.exp.moment} follows from~\eqref{e.exit.time.upper.moment} and the strong Markov property. Indeed,~\eqref{e.exit.time.upper.moment} implies a uniform bound~$\mathbf{P}^{\k,x}[\tau(z+\cu_n) > C\Tr(3^n)] \leq \frac12$ for some large constant~$C$; the strong Markov property then yields geometric decay~$\mathbf{P}^{\k,x}[\tau(z+\cu_n) > C k\Tr(3^n)] \leq 2^{-k}$, which implies~\eqref{e.exit.time.exp.moment}.

\smallskip

\emph{Step 2.} We show that there exists~$c(d,\cstar)>0$ such that, for every~$z\in\Rd$ with~$\indc_{\mathcal{Q}(z+\cu_n,\ep)}=1$,
\begin{equation}
	\inf_{x\in z+\overline{\cu}_{n-1}}
	\mathbf{P}^{\k,x}
	\bigl[ \tau(z+\cu_n) \geq c \Tr(3^n) \bigr]
	\geq c
	\,.
	\label{e.survival.probability}
\end{equation}
Select~$x\in z+\overline{\cu}_{n-1}$ and apply~\eqref{e.exit.time.lower.bound},~\eqref{e.exit.time.exp.moment} and the Paley-Zygmund inequality to obtain
for every~$\rho\in(0,1)$,
\begin{align*}
	\mathbf{P}^{\k,x}\Bigl[ \tau(z+\cu_n) \geq \rho C^{-1} \Tr(3^n) \Bigr]
	&
	\geq
	\mathbf{P}^{\k,x}\Bigl[ \tau(z+\cu_n) \geq \rho \mathbf{E}^{\k,x}[\tau(z+\cu_n)]\Bigr]
	\notag \\ &
	\geq
	(1-\rho)^2\frac{\mathbf{E}^{\k,x}[\tau(z+\cu_n)]^2}{\mathbf{E}^{\k,x}[\tau(z+\cu_n)^2]}
	\geq
	c (1-\rho)^2
	\,.
\end{align*}
Taking~$\rho=\nf12$, we obtain~\eqref{e.survival.probability}.

\smallskip

\emph{Step 3.}  We show that there exists~$c(d,\cstar)>0$ such that, for every~$z\in\Rd$ with~$\indc_{\mathcal{Q}(z+\cu_n,\ep)}=1$,
\begin{equation}
	\sup_{x\in z+\overline{\cu}_{n-1}}
	\mathbf{E}^{\k,x}
	\bigl[ \exp\bigl( - \theta \tau(z+\cu_n) \bigr) \bigr]
	\leq
	1 - c\bigl( 1 - \exp\bigl( -c \theta \Tr(3^n)  \bigr) \bigr)
	\,.
	\label{e.exit.time.exp.moment.pre}
\end{equation}
Assume~$z=0$ for simplicity and fix~$x\in \overline{\cu}_{n-1}$. Set~$a\coloneqq c_{\eqref{e.survival.probability}} \Tr(3^n)$ and split according to whether or not the event~$\{\tau(\cu_n) \geq a \}$ is valid.
\begin{align*}
	\mathbf{E}^{\k,x}\bigl[ \exp\bigl( -\theta\tau(\cu_n) \bigr) \bigr]
	&
	=
	\mathbf{E}^{\k,x}\bigl[\exp\bigl( -\theta\tau(\cu_n) \bigr) \indc_{\{\tau(\cu_n)<a\}}\bigr]
	+
	\mathbf{E}^{\k,x}\bigl[\exp\bigl( -\theta\tau(\cu_n) \bigr) \indc_{\{\tau(\cu_n)\geq a\}}\bigr]
	\\
	&\leq
	\mathbf{P}^{\k,x}[\tau(\cu_n) < a]
	+
	\exp(-\theta a)
	\mathbf{P}^{\k,x}[\tau(\cu_n) \geq a]
	\\
	&=
	1 - \mathbf{P}^{\k,x}[\tau(\cu_n) \geq a]
	\bigl(1-\exp(-\theta a)\bigr)
	\,.
	\\
	& \leq
	1-c \bigl(1-\exp( -\theta c \Tr(3^n))\bigr)\,,
\end{align*}
where in the last line we used~\eqref{e.survival.probability}.
Taking the supremum over~$x\in\overline{\cu}_{n-1}$ yields~\eqref{e.exit.time.exp.moment.pre}.

\smallskip

\emph{Step 4.}
We show that there exists~$c(d,\cstar)>0$ such that, for every~$\theta \geq \Tr(3^n)^{-1}$,
\begin{equation}
	\sup_{x\in \cu_{m-1}}
	\mathbf{E}^{\k,x}
	\bigl[ \exp\bigl( - \theta \tau(\cu_m) \bigr) \bigr]
	\leq
	\exp\bigl( - c3^{m-n}\bigr)
	\,.
	\label{e.exit.time.tail.bounds}
\end{equation}
The argument uses~\eqref{e.exit.time.exp.moment.pre}, the strong Markov property and a chaining argument. Denote
\begin{equation*}
	D_{m,n} \coloneqq
	\bigcup_{z\in Q_{m,n}}
	(z + \overline{\cu}_{n-1})
	\subseteq \overline{\cu}_{m} \setminus \cu_{m-1} \,.
\end{equation*}
We define a sequence of stopping times~$\{ T_j \}_{j=1}^N$ and points~$\{ z_j \}_{j=1}^N\subseteq Q_{m,n}$ as follows. Define
\begin{equation*}
	T_1 \coloneqq \inf \{s \geq 0 \,:\, X_s \in D_{m,n} \}
\end{equation*}
and let~$z_1 \in Q_{m,n}$ be such that~$X_{T_1} \in z_1 + \overline{\cu}_{n-1}$. Continuing recursively, let
\begin{equation*}
	\tilde{T}_{j} \coloneqq
	\inf \bigl\{s \geq T_j \,:\, X_s \not\in (z_j+ \cu_n) \bigr\}
\end{equation*}
be the time of the first exit from~$z_j + \cu_n$ after~$T_j$, and
\begin{equation*}
	T_{j+1}
	\coloneqq
	\inf \biggl\{ s \geq \tilde{T}_{j} \,:\, X_s \in D_{m,n} \setminus \bigcup_{i=1}^j (z_i + \overline{\cu}_{n-1}) \biggr\}
\end{equation*}
be the first time after~$\tilde{T}_j$ that the diffusion visits a good subcube it has not visited before. Let~$z_{j+1} \in Q_{m,n}$ be such that~$X_{T_{j+1}} \in z_{j+1}  + \overline{\cu}_{n-1}$.

\smallskip

By~\eqref{e.what.percolation.gives} and the continuity of the trajectories of the process~$\{ X_t \}$, every cube~$y_i + \cu_{n-1}$ is visited before the exit from~$\cu_m$. Each cube~$y_i+\cu_n$ which meets~$z_j+\cu_n$ is contained in~$z_j+\cu_{n+1}$; since the cubes~$y_i+\cu_n$ are pairwise disjoint, at most~$3^d$ of them meet~$z_j+\cu_n$, and therefore at most~$3^d$ of the~$y_i$ are consumed by any single step of the chain. Hence we obtain, for every~$x\in \cu_{m-1}$,
\begin{equation}
	\tau(\cu_m) \geq T_{N'}
	\,,\qquad
	N' \coloneqq \bigl\lceil 3^{-d} N \bigr\rceil
	\,,\qquad \text{$\mathbf{P}^{\k,x}$-a.s.}
	\label{e.step4.TN.le.tau}
\end{equation}
Thus, for every~$x\in \cu_{m-1}$,
\begin{equation}
	\mathbf{E}^{\k,x}\bigl[ \exp(-\theta\tau(\cu_m))
	\bigr]
	\leq
	\mathbf{E}^{\k,x}\bigl[\exp(-\theta T_{N'})
	\bigr]
	\,.
	\label{e.step4.reduce.to.TN}
\end{equation}
We estimate~$\mathbf{E}^{\k,x}[\exp( -\theta T_{N'})]$ by iterating Step~3.
By the strong Markov property at time~$T_j$ and the definition of the stopping times,~$T_{j+1}-T_j\geq \tilde{T}_j-T_j$ and
the increment~$\tilde{T}_j-T_j$ has the same law as the exit time~$\tau(z_{j}+\cu_n)$ under~$\mathbf{P}^{\k,X_{T_j}}$.
Since~$X_{T_j}\in z_j + \overline{\cu}_{n-1}$ and~$\theta \geq \Tr(3^n)^{-1}$, the estimate~\eqref{e.exit.time.exp.moment.pre} yields
\begin{align*}
	\mathbf{E}^{\k,X_{T_j}}\bigl[\exp( -\theta (T_{j+1}-T_j))\bigr]
	&
	\leq
	\mathbf{E}^{\k,X_{T_j}}\bigl[\exp( -\theta \tau(z_j+\cu_n)) \bigr]
	\leq
	1-c\bigl(1-\exp(-c\theta\Tr(3^n))\bigr)
	\leq
	1 - c
	\,.
\end{align*}
Therefore, for every~$x\in \cu_{m-1}$,
\begin{align*}
	\mathbf{E}^{\k,x}\bigl[\exp(-\theta T_{j+1}) \bigr]
	=
	\mathbf{E}^{\k,x}\Bigl[\exp( -\theta T_j)
	\mathbf{E}^{\k,X_{T_j}}\bigl[\exp(-\theta (T_{j+1}-T_j))\bigr]\Bigr]
	\leq
	(1-c)
	\mathbf{E}^{\k,x}\bigl[\exp(-\theta T_j)\bigr]
	\,.
\end{align*}
Iterating this inequality yields
\begin{equation}
	\mathbf{E}^{\k,x}
	\bigl[\exp(- \theta T_{N'})\bigr]
	\leq
	(1-c)^{N'}
	\leq
	\exp( - c N')
	\leq
	\exp
	\bigl( - c3^{m-n} \bigr)
	\,.
	\label{e.step4.q.power.N}
\end{equation}
This completes the proof of~\eqref{e.exit.time.tail.bounds}.

\smallskip

\emph{Step 5.}
The conclusion.
For every~$\theta \geq \Tr(3^n)^{-1}$, Markov's inequality and~\eqref{e.exit.time.tail.bounds} give
\begin{align}
	\sup_{x\in \cu_{m-1}}
	\mathbf{P}^{\k,x}\bigl[\tau(\cu_m) \leq t\bigr]
	& =
	\sup_{x\in \cu_{m-1}}
	\mathbf{P}^{\k,x}
	\Bigl[\exp(-\theta \tau(\cu_m))\geq \exp(-\theta t)\Bigr]
	\notag \\ &
	\leq
	\exp(\theta t)
	\sup_{x\in \cu_{m-1}}
	\mathbf{E}^{\k,x}\bigl[\exp(-\theta \tau(\cu_m))\bigr]
	\notag \\ &
	\leq
	\exp\bigl(\theta t - c 3^{m-n } \bigr)
	\,.
	\label{e.step5.markov}
\end{align}
Inserting the choice~$\theta \coloneqq \Tr(3^n)^{-1}$ into the previous display, using~\eqref{e.scale.sep.gives} and taking~$\delta$ to be a sufficiently small constant, we obtain
\begin{align*}
	\sup_{x\in \cu_{m-1}}
	\mathbf{P}^{\k,x}\bigl[\tau(\cu_m) \leq t\bigr]
	&
	\leq
	\exp\bigl( t \Tr(3^n)^{-1} - c 3^{m-n } \bigr)
	\notag \\ &
	\leq
	\exp\bigl( - c 3^{m-n } \bigr)
	\notag \\ &
	\leq
	\exp
	\biggl(
	-c
	\min \biggl\{
	\frac{3^{2m}}{\nu t}
	\,,\,
	\bigl( \cstar^{-1} \cgamma  \bigr)^{\nf{1}{2}}
	\Bigl(
	\frac{(3^m)^{2-\cgamma}}{t}
	\Bigr)^{\frac{1}{1-\cgamma} }
	\biggr\}
	\biggr)
	\,.
\end{align*}
This completes the proof of~\eqref{e.early.exit.probability}.
\end{proof}

\subsection{Quenched first and second moment bounds for the displacement}
\label{ss.proof.theorem.A}

In this subsection we complete the proof of Theorem~\ref{t.superdiffusivity} by proving the quantitative, quenched second-moment and variance estimates for the displacement of the process. The key idea is to compare~$X_t$ to the stopped process~$X_{t\wedge\tau(\cu_m)}$ on a cube~$\cu_m$ whose sidelength is chosen slightly larger than the intrinsic scale~$\Rt(t)$.  On the one hand, elliptic (Dirichlet) homogenization estimates give accurate control of~$\mathbf{E}^{\k,0}[(e\cdot X_{t\wedge\tau(\cu_m)})^2]$ and~$\mathbf{E}^{\k,0}[X_{t\wedge\tau(\cu_m)}]$. On the other hand, the exit-time tail bounds
proved in the previous subsection imply that~$\mathbf{P}^{\k,0}[\tau(\cu_m)\leq t]$ is extremely small when the size~$3^m$ of the cube is only
a factor of~$C\left|\log\cgamma \right|^{\nf12}$ larger than~$\Rt(t)$. Removing the stopping time therefore produces an error much smaller than the homogenization error.

\begin{proposition}
\label{p.displacement.superdiffusive.bounds}
Let~$\EthmB(m)$ denote the random variable in the statement of Theorem~\ref{t.homogenization}. There exists a constant~$C(d,\cstar) < \infty$ and, for every~$t\in (0,\infty)$, a random scale parameter~$m_t \in\Z$ with~$3^{m_t } \geq \Rt(t)$ and
\begin{equation}
	\E \bigl[ (3^{m_t })^p \bigr]^{\nf1p}
	\leq
	C \left| \log \cgamma\right|^{\nf12}
	\Rt(t)
	\,, \qquad \text{for} \  p=C^{-1} \cgamma^{-1} \left| \log \cgamma\right|^{-6}
	\,,
	\label{e.estimate.confinement.scale}
\end{equation}
such that, for every~$e\in \Rd$ with~$|e|=1$,
\begin{equation}
	\bigl|
	\mathbf{E}^{\k,0}\bigl[ ( e\cdot X_{t} )^2 \bigr]
	-
	2\Rt(t)^2
	\bigr|
	\leq
	C\bigl( \EthmB(m_t  )
	+
	\cgamma^{\nf12} \left| \log \cgamma\right|
	\bigr) 3^{2m_t  }
	\label{e.second.moment.ThmA}
\end{equation}
and
\begin{equation}
	\bigl|
	\mathbf{E}^{\k,0}
	\bigl[ X_{t} \bigr]
	\bigr|^2
	\leq
	C\bigl(  \EthmB(m_t  )^2
	+
	\cgamma^{80}
	\bigr) 3^{2m_t  }
	\,.
	\label{e.displacement.squared.ThmA}
\end{equation}
\end{proposition}
\begin{proof}
Fix~$t>0$. Fix a parameter~$K \in [1,\infty)$ which will be selected below (and depend only on~$d$ and~$\cstar$). The exponent~$40$ below is not optimal: the early exit probability~\eqref{e.exit.prob.small} is of order~$\cgamma^{100}$, and the square roots taken below produce~$\cgamma^{40}$; the slack in~\eqref{e.displacement.squared.ThmA} is then~$(\cgamma^{40})^2 = \cgamma^{80}$.
Let~$\{S_k\}_{k \in \Z}$ be the random variables from the statement of Proposition~\ref{p.early.exit.bound} and define~$\tilde{S}_n \coloneqq \sup_{k \geq n} S_k 3^{-2(k-n)}$ and the (random) confinement scale~$m_t$ by
\begin{equation}
	m_t  \coloneqq
	\inf
	\Bigl\{ n \in \Z \,:\,
	n \geq
	\bigl\lceil \log_3 \bigl( K (1+\tilde{S}_n)\left| \log \cgamma \right|^{\nf12} \Rt(t) \bigr) \bigr\rceil \Bigr\} \,.
	\label{e.m.confinement.scale}
\end{equation}
By definition, we have~$3^{m_t } \geq K(1+\tilde{S}_{m_t  }) \left| \log \cgamma \right|^{\nf12} \Rt(t)$. Therefore, by Proposition~\ref{p.early.exit.bound}, we have
\begin{equation}
	\mathbf{P}^{\k,0}\bigl[\tau(\cu_{m_t} )\leq t\bigr]
	\leq
	\exp\bigl(-c K^2 \left| \log \cgamma \right| \bigr)
	\,.
	\label{e.exit.prob.small.pre}
\end{equation}
Now choose~$K$ large enough that the right side is at most~$\cgamma^{100}$, so that we obtain
\begin{equation}
	\mathbf{P}^{\k,0}\bigl[\tau(\cu_{m_t} )\leq t\bigr]
	\leq
	\cgamma^{100}
	\,.
	\label{e.exit.prob.small}
\end{equation}
By Minkowski's inequality,~$\sup_n\|\tilde S_n\|_{L^p}\leq C$ in the range of~\eqref{e.estimate.St}, and the same tail-summation argument as in the proof of Proposition~\ref{p.minimal.scale.separation} gives the finiteness of~$m_t$ and~\eqref{e.estimate.confinement.scale}, after a fixed reduction of the admissible moment.

\smallskip

We henceforth write~$m=m_t$.
By~\eqref{e.formula.for.shom}, the definition~\eqref{e.Rt.def.again} of~$\Rt(t)$ and~$3^m \geq \Rt(t)$, we have
\begin{equation}
	| \Rt(t)^2 - \shom_m t | \leq
	C \cgamma^{\nf12} \left| \log \cgamma\right|
	3^{2m}
	\,.
	\label{e.Rtt.to.shom.t}
\end{equation}
By the definition of~$\tilde{S}_n$ and of~$m_t$ in~\eqref{e.m.confinement.scale}, the condition~\eqref{e.tail.displacement.minscale.condition} holds at every scale~$k \geq m$, so~\eqref{e.displacement.tail.estimate} may be integrated over all radii and gives
\begin{equation}
	\mathbf{E}^{\k,0}\bigl[|X_t|^p\bigr]
	\leq (C p)^{\nf p2} 3^{ p m} \qquad \forall p \in [1, \infty)\,.
	\label{e.pth.moment.process}
\end{equation}

\emph{Step 1.} Estimate for the mean displacement. We show that
\begin{equation}
	\bigl|
	\mathbf{E}^{\k,0}
	\bigl[ X_{t\wedge \tau(\cu_m) } \bigr]
	\bigr|^2
	\leq
	C3^{2m} \EthmB(m)^2
	\,.
	\label{e.displacement.squared.stopped}
\end{equation}
Let~$\linear_e(x) \coloneqq e \cdot x$. Define~$u$ to be the solution of the heterogeneous problem
\begin{equation*}
	\left\{
	\begin{aligned}
		& -\nabla \cdot \a \nabla u = 0 & \mbox{in} & \ \cu_m\,, \\
		& u = \linear_e & \mbox{on} & \ \partial \cu_m\,,
	\end{aligned}
	\right.
\end{equation*}
and let~$v \coloneqq \linear_e$ be the solution of the corresponding homogenized problem (since~$\linear_e$ is harmonic, it solves~$-\shom_m \Delta v = 0$ with the same boundary data). By Theorem~\ref{t.homogenization} with~$\g = 0$ and~$h = \linear_e$, noting~$3^{\frac12 m}\| \nabla \linear_e \|_{\underline{W}^{\nf12,\infty}(\cu_m)} \leq C$,
\begin{equation}
	\label{e.u.minus.v.mean}
	\| u - \linear_e \|_{L^\infty(\cu_m)}
	\leq C3^m \EthmB(m)
	\,.
\end{equation}
Since~$u$ is~$\a$-harmonic, the process~$u(X_t)$ is a local martingale. By optional stopping,
\begin{equation*}
	\mathbf{E}^{\k,0}[u(X_{t \wedge \tau(\cu_m) })] = u(0)\,.
\end{equation*}
Therefore, by~\eqref{e.u.minus.v.mean},
\begin{align*}
	\bigl| \mathbf{E}^{\k,0}[e \cdot X_{t \wedge \tau(\cu_m) }] \bigr|
	&
	=
	\bigl| \mathbf{E}^{\k,0}[\linear_e(X_{t \wedge \tau(\cu_m) })] \bigr|
	\\ &
	\leq
	\bigl| \mathbf{E}^{\k,0}[u(X_{t \wedge \tau(\cu_m) })] \bigr|
	+
	\bigl| \mathbf{E}^{\k,0}[(\linear_e - u)(X_{t \wedge \tau(\cu_m) })] \bigr|
	\\ &
	\leq
	|u(0)| + \| u - \linear_e \|_{L^\infty(\cu_m)}
	\\ &
	=
	|u(0) - \linear_e(0)| + \| u - \linear_e \|_{L^\infty(\cu_m)}
	\leq
	C 3^m \EthmB(m) \,.
\end{align*}
This yields~\eqref{e.displacement.squared.stopped}.

\smallskip

\emph{Step 2.}
We remove the stopping time from the mean displacement estimate~\eqref{e.displacement.squared.stopped} to complete the proof of~\eqref{e.displacement.squared.ThmA}.
We use
\begin{equation*}
	\mathbf{E}^{\k,0}[X_t]
	=
	\mathbf{E}^{\k,0}[X_{t\wedge\tau(\cu_m)}]
	+
	\mathbf{E}^{\k,0}
	\bigl[(X_t-X_{\tau(\cu_m)}) \indc_{\{\tau(\cu_m)\leq t\}} \bigr]
\end{equation*}
together with~\eqref{e.pth.moment.process} with~$p=2$ and the exit time bound~\eqref{e.exit.prob.small} to obtain
\begin{align}
	\bigl|\mathbf{E}^{\k,0}[X_t]-\mathbf{E}^{\k,0}[X_{t\wedge\tau(\cu_m)}]\bigr|
	&
	\leq
	\bigl(
	\mathbf{E}^{\k,0}\bigl[|X_t|^2\bigr]^{\nf12}
	+
	\mathbf{E}^{\k,0}\bigl[|X_{\tau(\cu_m)}|^2\bigr]^{\nf12}
	\bigr)
	\mathbf{P}^{\k ,0}\bigl[\tau(\cu_m)\leq t\bigr]^{\nf12}
	\notag\\ &
	\leq
	C \cgamma^{40} 3^m\,.
	\label{e.remove.stop.mean}
\end{align}
The combination of the previous display and~\eqref{e.displacement.squared.stopped} yields~\eqref{e.displacement.squared.ThmA}.

\smallskip

\emph{Step 3.} Estimate for the second moment of the stopped process. We show that
\begin{equation}
	\bigl|
	\mathbf{E}^{\k,0}[(e \cdot X_{t \wedge \tau(\cu_m)})^2] - 2 \shom_m t
	\bigr|
	\leq
	C 3^{2m} \EthmB(m)
	+
	2\shom_m t \mathbf{P}^{\k,0}[\tau(\cu_m)< t]
	\,.
	\label{e.second.moment.stopped}
\end{equation}
Let~$v(x) \coloneqq (e \cdot x)^2$. Then~$-\shom_m \Delta v = -2\shom_m$, so~$v$ solves the homogenized equation with constant right-hand side~$-2\shom_m$ and boundary data~$h(x) = (e\cdot x)^2$. Define~$u$ to be the solution of the heterogeneous problem with the same right-hand side and boundary data:
\begin{equation*}
	\left\{
	\begin{aligned}
		& -\nabla \cdot \a \nabla u = -2\shom_m & \mbox{in} & \ \cu_m\,, \\
		& u = (e \cdot x)^2 & \mbox{on} & \ \partial \cu_m\,.
	\end{aligned}
	\right.
\end{equation*}
Writing~$-2\shom_m = \nabla \cdot \g$ with~$\g(x) = -2\shom_m (e\cdot x)e$, and applying Theorem~\ref{t.homogenization}, we obtain
\begin{equation}
	\label{e.u.minus.v.second}
	\| u - v \|_{L^\infty(\cu_m)}
	\leq
	C 3^{2m}\EthmB(m)
	\,.
\end{equation}
Here we used~$[\g]_{\underline{W}^{\nf12,\infty}(\cu_m)} \leq C \shom_m 3^{\frac12 m}$ and~$\| \nabla h \|_{\underline{W}^{\nf12,\infty}(\cu_m)} \leq C 3^{\frac12 m}$.
By It\^o's formula, we have
\begin{equation}
	\mathbf{E}^{\k,0}[u(X_{t \wedge \tau(\cu_m)})] = u(0) + 2\shom_m \mathbf{E}^{\k,0}[t \wedge \tau(\cu_m)]
	\,,
	\label{e.u.ito}
\end{equation}
and therefore
\begin{align*}
	\mathbf{E}^{\k,0}[(e \cdot X_{t \wedge \tau(\cu_m)})^2]
	&
	=
	\mathbf{E}^{\k,0}[v(X_{t \wedge \tau(\cu_m)})]
	\\ &
	=
	\mathbf{E}^{\k,0}[u(X_{t \wedge \tau(\cu_m)})] + \mathbf{E}^{\k,0}[(v - u)(X_{t \wedge \tau(\cu_m)})]
	\\&=
	u(0) + 2\shom_m \mathbf{E}^{\k,0}[t \wedge \tau(\cu_m)]
	+ \mathbf{E}^{\k,0}[(v - u)(X_{t \wedge \tau(\cu_m)})]
	\,.
\end{align*}
Using~$v(0) = 0$,~\eqref{e.u.minus.v.second} gives~$|u(0)| = |u(0) - v(0)| \leq C  3^{2m}\EthmB(m)$, we deduce from the previous display that
\begin{align*}
	\bigl|
	\mathbf{E}^{\k,0}[(e \cdot X_{t \wedge \tau(\cu_m)})^2] - 2 \shom_m t
	\bigr|
	&
	\leq
	C 3^{2m} \EthmB(m) +
	2\shom_m \mathbf{E}^{\k,0}[t - t \wedge \tau(\cu_m)]
	\notag \\ &
	\leq
	C 3^{2m} \EthmB(m)
	+
	2\shom_m t \mathbf{P}^{\k,0}[\tau(\cu_m)< t]
	\,.
\end{align*}
This completes the proof of~\eqref{e.second.moment.stopped}.

\smallskip

\emph{Step 4.}
We remove the stopping time from the second moment estimate~\eqref{e.second.moment.stopped} to complete the proof of~\eqref{e.second.moment.ThmA}. We use~$(e\cdot X_{t\wedge\tau(\cu_m)})=(e\cdot X_t)$ on~$\{t<\tau(\cu_m)\}$ and~$(e\cdot X_{t\wedge\tau(\cu_m)})=(e\cdot X_{\tau(\cu_m)})$ on~$\{\tau(\cu_m)\leq t\}$, to obtain
\begin{equation*}
	\mathbf{E}^{\k,0}[(e\cdot X_t)^2]
	=
	\mathbf{E}^{\k,0}[(e\cdot X_{t\wedge\tau(\cu_m)})^2]
	+
	\mathbf{E}^{\k,0}\bigl[\bigl((e\cdot X_t)^2-(e\cdot X_{\tau(\cu_m)})^2\bigr)\indc_{\{\tau(\cu_m)\leq t\}}\bigr]
	\,.
\end{equation*}
Using the previous display,~\eqref{e.pth.moment.process} with~$p = 4$ and~\eqref{e.exit.prob.small} and the choice of~$m$, we obtain
\begin{align*}
	\lefteqn{
	\bigl|\mathbf{E}^{\k,0}[(e\cdot X_t)^2]-\mathbf{E}^{\k,0}[(e\cdot X_{t\wedge\tau(\cu_m)})^2] \bigr|
	} \qquad &
	\notag \\ &
	\leq
	\mathbf{E}^{\k,0}\bigl[(e\cdot X_t)^2\indc_{\{\tau(\cu_m)\leq t\}}\bigr]
	+
	\mathbf{E}^{\k,0}\bigl[(e\cdot X_{\tau(\cu_m)})^2\indc_{\{\tau(\cu_m)\le t\}}\bigr]
	\notag \\ &
	\leq
	\mathbf{E}^{\k,0}\bigl[|X_t|^4\bigr]^{\nf12}
	\mathbf{P}^{\k ,0}\bigl[\tau(\cu_m)\leq t\bigr]^{\nf12}
	+
	C3^{2m}\mathbf{P}^{\k ,0}\bigl[\tau(\cu_m)\le t\bigr]
	\notag \\ &
	\leq
	C\cgamma^{40}
	3^{2m}
	\,.
\end{align*}
This inequality, combined with~\eqref{e.Rtt.to.shom.t} and~\eqref{e.second.moment.stopped} and the triangle inequality,
yields~\eqref{e.second.moment.ThmA}.
\end{proof}

We conclude this subsection with the proof of the first two estimates in Theorem~\ref{t.superdiffusivity}.

\begin{proof}[{Proof of Theorem~\ref{t.superdiffusivity}}]
Let~$m_t$ denote the random scale from~\eqref{e.m.confinement.scale} of Proposition~\ref{p.displacement.superdiffusive.bounds} and let
\begin{equation*}
	\tilde{m} \coloneqq \lceil \log_3(K |\log\cgamma|^{\nf12} \Rt(t)) \rceil
\end{equation*}
denote its deterministic lower bound.
Let~$C$ denote the larger of the constants in~\eqref{e.mathcal.E.estimate} and~\eqref{e.estimate.confinement.scale}. For~$p \in [1, (2C)^{-1} \cgamma^{-1} |\log\cgamma|^{-6}]$, by the Cauchy-Schwarz inequality, the uniform-in-$m$ bound~\eqref{e.mathcal.E.estimate}, and~\eqref{e.estimate.confinement.scale},
\begin{align}
	\E \bigl[ \EthmB(m_t)^p \bigr]
	&= \E \biggl[ 3^{p(m_t - \tilde{m})} \sum_{k=0}^{\infty} \EthmB(\tilde{m}+k)^p \, 3^{-pk} \, \indc_{\{m_t = \tilde{m}+k\}} \biggr] \notag \\
	&\leq \E \bigl[ 3^{2p(m_t - \tilde{m})} \bigr]^{\nf12} \biggl( \sum_{k=0}^{\infty} 3^{-2pk} \, \E \bigl[ \EthmB(\tilde{m}+k)^{2p} \bigr] \biggr)^{\!\nf12} \notag \\
	&\leq \bigl( 3^{-2p\tilde{m}} \E \bigl[ 3^{2pm_t} \bigr] \bigr)^{\nf12} \biggl( \sum_{k=0}^{\infty} 3^{-2pk} \biggr)^{\!\nf12}
	\bigl( C (p^{\nf12}+|\log\cgamma|^{\nf12}) \cgamma^{\nf12} |\log\cgamma|^3 \bigr)^{p} \notag \\
	&\leq \bigl( C (p^{\nf12}+|\log\cgamma|^{\nf12}) \cgamma^{\nf12} |\log\cgamma|^3 \bigr)^{p}\,.
	\label{e.EthmB.mt.bound}
\end{align}
By the above display and~\eqref{e.second.moment.ThmA}, for~$p \in[1, C^{-1} \cgamma^{-1} \left| \log\cgamma \right|^{-6} ]$, we have
\begin{align*}
	\E \Bigl[
	\bigl|
	\mathbf{E}^{\k,0}\bigl[ ( e\cdot X_{t} )^2 \bigr]
	-
	2\Rt(t)^2
	\bigr| ^p
	\Bigr]^{\nf1p}
	&
	\leq
	C
	\E \Bigl[
	\bigl( \EthmB(m_t)
	+
	\cgamma^{\nf12} \left| \log \cgamma \right|
	\bigr)^p 3^{2m_t p}
	\Bigr]^{\nf1p}
	\notag \\ &
	\leq
	C \bigl(p^{\nf12} + \left| \log \cgamma \right|^{\nf12}\bigr) \cgamma^{\nf12}
	\left| \log \cgamma \right|^4
	\Rt(t)^2
	\,.
\end{align*}
Since~$|X_t|^2 = \sum_{i=1}^d ( e_i \cdot X_t )^2$ and~$2d\Rt(t)^2 = \sum_{i=1}^d 2\Rt(t)^2$, applying the previous display at~$e=e_1,\ldots,e_d$ and summing gives the first estimate in~\eqref{e.thm.superdiffusivity}.
Similarly, for the square of the expected displacement we have, by~\eqref{e.displacement.squared.ThmA}, for the same~$p$,
\begin{align*}
	\E \Bigl[
	\bigl|
	\mathbf{E}^{\k,0}
	\bigl[ X_{t} \bigr]
	\bigr|^{2p} \Bigr]^{\nf1p}
	&
	\leq
	C
	\E \Bigl[
	C\bigl(  \EthmB(m_t  )
	+
	\cgamma^{40}
	\bigr)^{2p}
	3^{2m_t p}
	\Bigr]^{\nf1p}
	\leq
	C \bigl(p + \left| \log \cgamma \right|\bigr) \cgamma \left| \log \cgamma \right|^7
	\Rt(t)^2
	\,.
\end{align*}
This completes the proof of~\eqref{e.thm.superdiffusivity}.
\end{proof}

\appendix
\section{Orlicz notation for tail bounds}
\label{s.Orlicz}

Throughout the paper we track tail bounds of random variables as follows: for~$A > 0$ and an increasing function~$\Psi: \R_+ \to [1, \infty)$ satisfying
\begin{equation}
	\label{e.Psi.int}
	\int_1^\infty
	\frac{t}{\Psi(t)}\,dt  < \infty
	\,,
\end{equation}
and a random variable~$X$, we write
\begin{equation} \label{e.orlicz.tail}
X \leq \O_{\Psi}(A)
\end{equation}
to mean
\begin{equation} \label{e.orlicztailbound.bound}
\P[X > t A] \leq \frac{1}{\Psi(t)} \, , \quad \forall t \in [1,\infty) \,\,.
\end{equation}
We write
\begin{equation} \label{e.orlicz.tail.symmetric}
X = \O_{\Psi}(A)
\end{equation}
to mean that both~$X \leq \O_{\Psi}(A)$ and~$-X \leq \O_{\Psi}(A)$, i.e., we have a two-sided tail bound.

The constraint~$t \geq 1$ in~\eqref{e.orlicztailbound.bound} can be removed: since~$\Psi$ is increasing and~$\Psi(1) \geq 1$, we have for~$t < 1$ that~$\P[X > tA] \leq \P[X > A] \leq 1/\Psi(1) \leq \Psi(1)/\Psi(t)$. Thus
\begin{equation*}
	X \leq \O_{\Psi}(A)
	\implies
	\P[X > t A] \leq \frac{\Psi(1)}{\Psi(t)}\,,  \quad \forall t > 0 \,\,.
\end{equation*}

As we will see below, this notation induces a useful algebra which allows us to essentially add and multiply tail bounds. We write
\begin{equation*}
	X \leq \O_{\Psi_1}(A_1) + \O_{\Psi_2}(A_2)
\end{equation*}
to mean that there exist random variables~$Y$ and~$Z$ such that~$X \leq Y + Z$,~$Y \leq \O_{\Psi_1}(A_1)$, and~$Z \leq \O_{\Psi_2}(A_2)$. More generally, we write
\begin{equation*}
	X \leq \O_{\Psi_1}(A_1) + \ldots +  \O_{\Psi_n}(A_n)
\end{equation*}
to mean that
\begin{equation*}
	X \leq X_1 + \cdots + X_n \,, \qquad \mbox{where} \ X_i \leq \O_{\Psi_i}(A_i)\,, \quad \forall i\in\{1,\ldots,n\}\,.
\end{equation*}

\smallskip

Although we introduce this notation for a general increasing function~$\Psi$ satisfying~\eqref{e.Psi.int}, for most of this paper we will use~$\Psi = \Gamma_{\sigma}$, where for~$\sigma \in (0, \infty)$,
\begin{equation}
	\label{e.Gamma.sigma}
	\Gamma_{\sigma}(t)  \coloneqq  \exp(t^{\sigma}) \,\,.
\end{equation}
This describes random variables with stretched exponential tails. The important case~$\sigma = 2$ specifies Gaussian tails, and~\eqref{e.kn.reg.ass} can be rewritten as
\begin{equation}
	\label{e.kn.reg.ass.with.gamma}
	\| \mathbf{j}_n \|_{L^{\infty}(\cu_n)}
	+
	\sqrt{d} 3^n \| \nabla \mathbf{j}_n \|_{L^{\infty}(\cu_n)}
	+
	d 3^{2n} \| \nabla^2 \mathbf{j}_n  \|_{L^{\infty}(\cu_n)}
	\leq \O_{\Gamma_2}(1) \,\,.
\end{equation}

We recall some basic properties and refer to~\cite[Appendix C]{AK.HC} and~\cite[Appendix A]{AKMBook} for an in-depth discussion.

\begin{lemma}[Generalized triangle inequality]
\label{l.Gamma.sigma.triangle}
There exists a universal constant~$C<\infty$ such that, for every~$\sigma \in (0,\infty)$ and sequence~$\{X_k\}_{k \in \N}$ of random variables with~$\sum_{k \in \N} a_k < \infty$,
\begin{equation}
	\label{e.Gamma.sigma.triangle}
	X_k \leq \O_{\Gamma_\sigma}(a_k)
	\quad \implies \quad
	\sum_{k \in \N} X_k \leq \O_{\Gamma_\sigma}\biggl( \bigl(1 + C\sigma^{-\nf1\sigma}\indc_{\{\sigma < 1\}} \bigr) \sum_{k \in \N} a_k \biggr) \,.
\end{equation}
\end{lemma}
\begin{proof}
The inequality for~$\sigma \geq 1$ is proved in~\cite[Lemma A.4]{AKMBook}. The proof of that lemma also gives the statement for~$0<\sigma < 1$.
\end{proof}

\begin{lemma}[Multiplication property]
\label{l.o.gamma2.mult}
For every~$\sigma_1, \sigma_2 \in (0, \infty)$, if~$X_1$ and~$X_2$ are positive random variables,
then
\begin{equation}
	\label{e.multGammasig}
	X_1 \leq \O_{\Gamma_{\sigma_1}}(A_1) \qand X_2 \leq \O_{\Gamma_{\sigma_2}}(A_2) \implies X_1 X_2 \leq \O_{\Gamma_{\frac{\sigma_1 \sigma_2}{\sigma_1 + \sigma_2}}}\bigl( (1+\log 2)^{\frac{\sigma_1+\sigma_2}{\sigma_1\sigma_2}} A_1 A_2 \bigr) \,\,.
\end{equation}
In particular, for every~$\sigma, p, K \in (0, \infty)$ and positive random variable~$X$,
\begin{equation}
	\label{e.powerofGammasigma}
	X \leq \O_{\Gamma_{\sigma}}(K)
	\iff
	X^{p} \leq \O_{\Gamma_{\sigma/p}} (K^p)\,.
\end{equation}
\end{lemma}
\begin{proof}
This is~\cite[Lemma A.3]{AKMBook}.
\end{proof}

\begin{lemma}[Maximum of~$\O_{\Gamma_\sigma}$ random variables]
\label{l.maximums.Gamma.s}
Suppose that~$\sigma,A>0$ and~$X_1,\ldots,X_N$ is a sequence of random variables satisfying~$X_i \leq \O_{\Gamma_\sigma}(A)$
for~$N \geq 2$. Then
\begin{equation}
	\label{e.maxy.bound}
	\max_{1\leq i \leq N} X_i \leq \O_{\Gamma_\sigma} \bigl( (3 \log N)^{\nf1\sigma} A\bigr)\,.
\end{equation}
\end{lemma}
\begin{proof}
For every~$t \geq 1$ we use a union bound to estimate
\begin{align*}
	\P \Bigl[ \max_{1\leq i \leq N} X_i > A (3 \log N)^{\nf{1}{\sigma}} t  \Bigr]
	&\leq \sum_{i=1}^N \P \bigl[ X_i > A (3 \log N)^{\nf{1}{\sigma}} t\bigr] \\
	&\leq N \exp\bigl( -3 t^{\sigma} \log N \bigr)
	= \exp\bigl(-t^{\sigma} (3 \log N - \log N)  \bigr)
	\leq  \exp(-t^\sigma) \,\,,
\end{align*}
where in the last inequality we used~$2 \log N \geq 1$ for~$N \geq 2$.
\end{proof}

The indicator function~$\indc_E$ of an event~$E$ with~$0<\P[E] <1$ satisfies, for every~$\sigma\in (0,\infty)$,
\begin{equation}
	\label{e.indc.O.sigma}
	\indc_{E} \leq \O_{\Gamma_\sigma} \bigl( \bigl| \log \P[E] \bigr|^{-\nf 1\sigma} \bigr)
	\,.
\end{equation}
This is immediate from definitions~\eqref{e.orlicz.tail}-\eqref{e.orlicztailbound.bound}.

\begin{lemma}[Moment bounds]
\label{l.moments.gamma.psi}
There exists a universal~$C<\infty$ such that, for every~$K,\sigma \in (0,\infty)$ and random variable~$X$ satisfying~$X = \O_{\Gamma_{\sigma} }(K)$, we have, for every~$p\in [1,\infty)$,
\begin{equation}
	\label{e.moments.OGamma2}
	\E\bigl[|X|^p \bigr]^{\nf1p}
	\leq
	C p^{\nf1\sigma} K
	\,.
\end{equation}
\end{lemma}
\begin{proof}
By replacing~$X$ with~$X^{\sigma}$ and applying~\eqref{e.powerofGammasigma}, we can reduce to the case~$\sigma =1$. By homogeneity, we may reduce to the case~$K=\nf12$. For~$X = \O_{\Gamma_1}(\nf12 )$ and~$m\in\N$, we have
\begin{equation*}
	\frac1{m!} \E\bigl[ |X|^m \bigr]
	\leq
	\E\bigl[ \exp( |X| )  \bigr]
	\leq
	C\,.
\end{equation*}
Taking the~$m$th roots and using~$m! \leq m^m$ yields the result for integer~$m$. The result for~$p \in (m,m+1)$ follows from H\"older's inequality and the result for~$m$ and~$m+1$.
\end{proof}

We will use the following concentration inequality for sums of centered, independent random variables with stretched exponential tails.
For a proof, see~\cite[Lemma~C.2 \& Corollary C.4]{AK.HC}.

\begin{proposition}[Concentration for $\O_{\Gamma_\sigma}$]
\label{p.concentration}
There exists a universal constant~$C <\infty$ such that, for every~$\sigma \in (0,2]$,~$m\in\N$, and finite sequence~$X_1,\ldots,X_m$ of independent random variables satisfying
\begin{equation}
	|X_k| \leq \O_{\Gamma_\sigma} (1)
	\quad \mbox{and} \quad
	\E[ X_k ] = 0 \, , \quad \forall k\in \{1,\ldots,m\}\,,
\end{equation}
we have the estimate
\begin{equation}
	\label{e.concentration}
	\sum_{k=1}^m X_k \leq
	\O_{\Gamma_\sigma}
	\Bigl( \bigl( ( C\sigma^{-1} )^{\nf 1\sigma} + C|\log(\sigma-1)|^{\nf1\sigma} \indc_{\{ \sigma >1\}} \bigr)
	m^{\nf12}
	\Bigr)
	\,.
\end{equation}
\end{proposition}

When the random variables~$\{ X_z \}_{z\in 3^n\Zd \cap \cu_m}$, for some~$m,n\in \N$ with~$n<m$, have the property that~$X_z$ and~$X_{z'}$ are independent provided that the corresponding subcubes do not touch, then Proposition~\ref{p.concentration} implies
\begin{equation*}
	\sum_{z\in 3^n\Zd \cap \cu_m} X_z
	\leq
	\O_{\Gamma_\sigma} \bigl(C(\sigma) 3^{\frac d2(m-n)} \bigr)
	\,,
\end{equation*}
where~$C(\sigma)$ is the constant from~\eqref{e.concentration}.

\section{Percolation estimates}
\label{s.percolation.estimate}

In this section, distances between subsets of~$\Zd$ are measured in the~$\ell^\infty$ metric,~$\dist(A,B) \coloneqq \inf_{x\in A, y\in B} |x-y|_\infty$. We define a~\emph{path} in~$\Zd$ to be a finite sequence~$\{ x_0, \ldots, x_n \}$ with~$\dist(x_i, x_{i+1}) \leq 1$ for~$i = 0, 1,\ldots,n-1$. For~$A,B\subseteq \Zd$, we say that~$\{ x_0,\ldots,x_n\}$ is a~\emph{path from~$A$ to~$B$} if~$\{ x_0,\ldots,x_n\}$ is a path and~$x_0 \in A$,~$x_n\in B$.

\begin{lemma}
\label{l.percolation.bound.general}
Let~$m \in \N$,~$a\in (0,d]$ and~$T \geq 1$.
Let~$\{ \mathrm{B}_L(z) : z \in \Zd \, , L \in \N_0 \}$ be a collection of events satisfying
\begin{itemize}

\item For each~$L\in \N_0$ and disjoint subsets~$\mathcal{Z},\mathcal{Z}' \subseteq \Zd$ with~$\dist(\mathcal{Z},\mathcal{Z}^\prime)>3^L$, the families~$\{\mathrm{B}_{L'}(z) \,:\, z\in \mathcal{Z}, \ L' \in [0,L] \cap \N_0  \}$ and~$\{ \mathrm{B}_{L'}(z) \,:\, z\in\mathcal{Z}^\prime, \ L' \in [0,L] \cap \N_0  \}$ are independent;

\item For every~$z\in\Zd$ and~$L\in\N_0$,
\begin{equation}
	\P\bigl[\mathrm{B}_L(z)\bigr]
	\leq
	\exp\bigl( - T^2 3^{a L} \bigr) \, .\label{e.badevent.upperbound.appendix}
\end{equation}

\end{itemize}
For each~$z \in \Zd$, define
\begin{equation}
	\mathrm{B}(z)  \coloneqq   \bigcup_{L \in \N_0} \mathrm{B}_L(z)  \,\,.
	\label{e.badevent.decomp}
\end{equation}
Then the following statements are valid:
\begin{enumerate}
\item \underline{Density bound.}
There exist constants~$c(d) > 0$ and~$C(d)<\infty$ such that for every~$t \geq C\bigl(a^{-1}+(d-a)^{-1}\bigr) T^{-1} 3^{-\frac a2 m}$, if~$a < d$,
\begin{equation}
	\P\Bigl[ \avsum_{z \in \cu_m \cap \Zd} (\indc_{\mathrm{B}(z)} - \P[\mathrm{B}(z)])  > t \Bigr] \leq
	\exp\bigl(-ct^2 a^2 (d-a)^2 T^2 3^{a m }\bigr)
	\label{e.density.bound} \,\,.
\end{equation}

\item \underline{Path bound.}
There exists a constant~$C(d) < \infty$ so that for every~$\lambda \in (0,1)$, if
\[
a > 1 \qand T \geq C \lambda^{-1} (a-1)^{-1}
\]
then for every~$k \in \N$, we have
\begin{multline}
	\label{e.path.bound}
	\P\Biggl[\exists\ \textup{a path}\ \Gamma \subseteq \Zd\ \textup{from}\ \Zd \cap \cu_k\ \textup{to}\ \Zd \setminus \cu_{k+1}\ \textup{with}\ \sum_{z \in \Gamma} \indc_{\mathrm{B}^c(z)} < (1-\lambda) 3^k \Biggr]
	\\
	\leq \exp\bigl(-c(a-1)^2 T^2 \lambda^2 \cdot 3^k\bigr) \,\,.
\end{multline}

\item \underline{Maximal diameter bound.}
There exists a constant~$C(d)<\infty$ such that, for every~$b \in (0,1) \cap (0,a]$ satisfying
\begin{equation*}
	16 \, d \exp\bigl(C(1-b)^{-1}\bigr)
	\leq
	b \, T^2
	\,,
\end{equation*}
and every~$k\in\N$, we have
\begin{multline}
	\label{e.diameter.bound}
	\P\Bigl[\exists\ \textup{a path }\Gamma\subset \Zd\ \textup{from }\Zd\cap\cu_k \ \textup{to }\Zd \setminus \cu_{k+1}
	\ \textup{such that }\mathrm{B}(z) \ \textup{occurs for every }z\in\Gamma\Bigr]
	\\
	\leq
	\exp\Bigl(- \frac12 \exp\bigl(-C(1-b)^{-1}\bigr)T^2 3^{bk} \Bigr)
	\,.
\end{multline}

\end{enumerate}
\end{lemma}
\begin{proof}
We prove the three statements in order.

\smallskip

\emph{Step 1: Density bound~\eqref{e.density.bound}.}
This follows from a simple union bound and independence.
By~\eqref{e.indc.O.sigma} and~\eqref{e.badevent.upperbound.appendix} for every~$L \in \N_0$
and~$z \in \Zd$,
\begin{equation*}
	\indc_{\mathrm{B}_L(z)} \leq \O_{\Gamma_{2}}(|\log \P[\mathrm{B}_L(z)]|^{-\nf 1 2}) \leq
	\O_{\Gamma_{2}}(T^{-1} 3^{-\frac{aL}{2}}) \,.
\end{equation*}
Using~\eqref{e.concentration}, we obtain
\begin{equation*}
	\avsum_{z \in \Zd \cap \cu_m} \bigl(\indc_{\mathrm{B}_L(z)} - \P[\mathrm{B}_L(z)] \bigr)
	=
	\O_{\Gamma_{2}}\Bigl(
	C
	3^{-\frac d2 (m-L)_+}
	T^{-1} 3^{-\frac{aL}{2}}
	\Bigr)  \qquad \forall L \in \N_0 \,.
\end{equation*}
Using~$\indc_{\mathrm{B}(z)}\le \sum_{L\ge0}\indc_{\mathrm{B}_L(z)}$ and~$\P[\mathrm{B}(z)]\le\sum_{L\ge0}\P[\mathrm{B}_L(z)]$, we may sum over~$L$ and apply~\eqref{e.Gamma.sigma.triangle} to obtain
\begin{equation*}
	\avsum_{z \in \cu_m \cap \Zd} (\indc_{\mathrm{B}(z)} - \P[\mathrm{B}(z)])
	\leq
	\O_{\Gamma_2}\bigl(C (a^{-1} + (d-a)^{-1}) T^{-1} 3^{-\frac{a}{2} m}\bigr)
	\,,
\end{equation*}
which implies~\eqref{e.density.bound}.

\smallskip

\emph{Step 2: Path bound~\eqref{e.path.bound}.}
Define the inflated events,
\[
\widetilde{\mathrm{B}}(z) \coloneqq \bigcup_{L=0}^{\infty} \bigcup_{u \in 3^{L} \Zd : z \in (u + \cu_L)} \widetilde{\mathrm{B}}_L(u)
\qquad \mbox{where} \qquad
\widetilde{\mathrm{B}}_L(u) \coloneqq
 \bigcup_{y \in (u + \cu_L)} \mathrm{B}_L(y) \, .
\]
Observe that for any path~$\gamma$ from~$\Zd \cap \cu_k$ to~$\Zd \setminus \cu_{k+1}$ of length~$|\gamma| = 3^k + m$,
\begin{equation}
	\label{e.path.event.containment}
	\biggl\{ \sum_{z \in \gamma} \indc_{\mathrm{B}^c(z)} < (1-\lambda) 3^k \biggr\} \subseteq \biggl\{ \sum_{z \in \gamma} \indc_{\widetilde{\mathrm{B}}(z)} > \lambda 3^k + m \biggr\} \,\,.
\end{equation}
Also, a union bound and~\eqref{e.badevent.upperbound.appendix} give~$\P[\widetilde{\mathrm{B}}_L(u)] \leq \exp(-\tfrac{1}{2}T^2 3^{aL})$ for~$T \geq C(d)$, hence, by~\eqref{e.indc.O.sigma},
\begin{equation}
	\label{e.orlicz.for.B.tilde}
	\indc_{\widetilde{\mathrm{B}}_L(u)} - \P[\widetilde{\mathrm{B}}_L(u)] \leq \O_{\Gamma_2}(T^{-1} 3^{-aL/2}) \,\,.
\end{equation}
Among all paths from~$\Zd \cap \cu_k$ to~$\Zd \setminus \cu_{k+1}$, let~$\Gamma$ be one which minimizes~$ \sum_{z \in \Gamma} \indc_{\widetilde{\mathrm{B}}(z)^c}$ and among such minimizers, has minimal length. Write~$|\Gamma| = 3^k + m$ for some~$m \geq 0$. If~$\Gamma$ were to exit a cube~$u + \cu_L$ where~$\widetilde{\mathrm{B}}_L(u)$ occurs and later re-enter, we could replace the excursion by a geodesic inside the cube, not increasing the sum. Hence~$\Gamma$ visits each such cube in a single segment of length at most~$C 3^L$. Defining, for any path~$\gamma$,
\begin{equation*}
	N_L(\gamma) \coloneqq \#\bigl\{u \in 3^L\Zd : \widetilde{\mathrm{B}}_L(u) \text{ occurs and } (u + \cu_L) \cap \gamma \neq \emptyset \bigr\} \,\,,
\end{equation*}
the single-visit property of the minimizer~$\Gamma$ implies that the right-hand side of~\eqref{e.path.event.containment} satisfies
\begin{equation}
	\label{e.S.threshold}
	\biggl\{ \sum_{z \in \Gamma} \indc_{\widetilde{\mathrm{B}}(z)} > \lambda 3^k + m \biggr\}
	\subseteq
	\Biggl\{ \sum_{L=0}^{\infty} 3^L N_L(\Gamma) > c (\lambda 3^k + m) \Biggr\} \,\,.
\end{equation}
We now estimate the probability that a deterministic path~$\gamma$ satisfies the event on the right above.

For any fixed path~$\gamma$ of length~$3^k + m$, we first handle large scales. For~$3^L > 3^k + m$, the path intersects at most~$C$ cubes at scale~$L$, so
\begin{equation*}
	\P\biggl[\sum_{L : 3^L > 3^k + m} 3^L N_L(\gamma) \geq 1\biggr] \leq \sum_{L : 3^L > 3^k + m} C \, e^{-cT^2 3^{aL}} \leq e^{-cT^2 (3^k+m)^a} \,\,.
\end{equation*}
For~$3^L \leq 3^k + m$, the number of cubes of the form~$u + \cu_L$ intersecting~$\gamma$ is at most~$C(3^k + m)3^{-L}$. The event~$\widetilde{\mathrm{B}}_L(u)$ depends only on the sites of~$u+\cu_L$, and these cubes touch; we therefore colour~$3^L \Zd$ by the residue of~$3^{-L}u$ modulo~$3$, so that two distinct sites~$u,u'$ of one colour class satisfy~$|u-u'|_\infty \geq 3^{L+1}$ and hence~$\dist(u+\cu_L,u'+\cu_L) \geq 2\cdot 3^L > 3^L$, and the events of that class are independent. Applying Proposition~\ref{p.concentration} and~\eqref{e.orlicz.for.B.tilde} within each of the~$3^d$ colour classes and summing the classes by Lemma~\ref{l.Gamma.sigma.triangle}, we obtain
\begin{equation*}
	3^L \sum_{u : (u+\cu_L) \cap \gamma \neq \emptyset} \bigl(\indc_{\widetilde{\mathrm{B}}_L(u)} - \P[\widetilde{\mathrm{B}}_L(u)]\bigr) \leq \O_{\Gamma_2}\bigl(CT^{-1} (3^k + m)^{1/2} \cdot 3^{(1-a)L/2}\bigr) \,\,.
\end{equation*}
Here~$C = C(d)$ also carries the factor~$3^d$ coming from the colouring.
Since~$1 < a \leq d$, the prefactors~$3^{(1-a)L/2}$ are summable with~$\sum_{L=0}^{\infty} 3^{(1-a)L/2} = \bigl( 1 - 3^{-(a-1)/2} \bigr)^{-1} \leq C(d)(a-1)^{-1}$. Lemma~\ref{l.Gamma.sigma.triangle} therefore gives
\begin{equation}
	\label{e.sum.minus.expectation}
	\sum_{L : 3^L \leq 3^k + m} 3^L N_L(\gamma) - \E\biggl[\sum_{L : 3^L \leq 3^k + m} 3^L N_L(\gamma)\biggr] \leq \O_{\Gamma_2}\bigl(C(a-1)^{-1}T^{-1}(3^k + m)^{1/2}\bigr) \,\,.
\end{equation}
Using that~$\P[\widetilde{\mathrm{B}}_L(u)] \leq e^{-cT^2 3^{aL}}$, the above display gives, for~$T \geq C$,
\begin{equation}
	\label{e.sum.orlicz}
	\sum_{L : 3^L \leq 3^k + m} 3^L N_L(\gamma) \leq C(3^k + m) e^{-cT^2} + \O_{\Gamma_2}\bigl(C(a-1)^{-1}T^{-1}(3^k + m)^{1/2}\bigr) \,\,.
\end{equation}
Combining the bounds for small and large scales, for~$T \geq C(a-1)^{-1}\lambda^{-1}$,
\begin{equation*}
	\P\biggl[\sum_{L=0}^{\infty} 3^L N_L(\gamma) > c (\lambda3^k+m)\biggr] \leq \exp\bigl(-c(a-1)^{2}T^{2}\lambda^2(3^k+m)\bigr) \,\,.
\end{equation*}
The number of paths of length~$3^k + m$ starting in~$\cu_k$ is at most~$3^{dk}3^{d(3^k+m)}$. A union bound over all paths gives
\begin{equation*}
	\P\biggl[\exists\, \gamma : |\gamma| = 3^k+m,\, \sum_{L=0}^{\infty} 3^L N_L(\gamma) > c (\lambda3^k+m) \biggr] \leq \exp\bigl(C(3^k+m) - c(a-1)^{2}T^{2}\lambda^2(3^k+m)\bigr) \,\,.
\end{equation*}
For~$(a-1)T\lambda \geq C$, the exponent is at most~$-c(a-1)^2T^2\lambda^2(3^k+m)/2$. Summing over~$m \geq 0$ yields~\eqref{e.path.bound}.

\smallskip

\emph{Step 3: Maximal diameter bound~\eqref{e.diameter.bound}.}
We introduce, for each~$k \in \N_0$ and~$z \in \Zd$, the event
\begin{equation*}
	\mathcal{C}_k(z)
	\coloneqq
	\Bigl\{
	\text{there exists a path~$\Gamma$ from~$\{z\}$ to~$\Zd \setminus( z+ \cu_k)$ such that~$\bigcap_{y\in\Gamma} \Bigl(\bigcup_{L \in [0,k]} \mathrm{B}_L(y)\Bigr)$  occurs}
	\Bigr\} \,\,.
\end{equation*}
It is immediate from the assumptions that, for every~$k\in\N$ and disjoint subsets~$\mathcal{Z},\mathcal{Z}^\prime \subseteq \Zd$ with~$\dist(\mathcal{Z},\mathcal{Z}^\prime) >  2\cdot 3^k$,
\begin{equation*}
	\mbox{the families\ $\{ \mathcal{C}_k(z) \,:\, z\in \mathcal{Z} \}$ \ and \ $\{ \mathcal{C}_k(z) \,:\, z\in \mathcal{Z}^\prime \}$ \ are independent}
	\,.
\end{equation*}
We claim that, for every~$k,h\in [5,\infty)\cap \N$,
\begin{equation}
	\label{e.paths.ready.for.iteration}
	\sup_{z\in \Zd} \P \bigl[ \mathcal{C}_{k+h} (z) \bigr]
	\leq
	\Bigl(
	3^{d(k+h)} \sup_{z\in \Zd} \P \bigl[ \mathcal{C}_{k}(z)  \bigr] \Bigr)^{3^{h-4}}
	+
	(k+h)3^{d(k+h)}
	\exp\bigl( - T^2 3^{a (k+1) } \bigr)
	\,.
\end{equation}
Suppose that~$\mathcal{C}_{k+h} (z)$ occurs. Let~$\Gamma$ be a path from~$\{z\}$ to~$\Zd\setminus (z+\cu_{k+h})$ such that, for each~$y\in \Gamma$, there exists~$L \in [0,k+h]\cap \N_0$ such that the event~$\mathrm{B}_{L}(y)$ occurs. Then there exist points~$y_1, \ldots, y_M \in \Gamma$ with~$M  \coloneqq  3^{h-4}$ such that~$\dist(y_i,y_j) \geq 3^{k+1}$ for every~$i\neq j$. Indeed, since~$\Gamma$ has length at least~$\frac12 \cdot 3^{k+h}$ and consecutive selected points are spaced by at least~$3^{k+1}$, we can select at least~$3^{h-1}/2 \geq 3^{h-4}$ such points. It follows that, either~$\mathcal{C}_k(y_i)$ occurs for every~$i \in \{ 1,\ldots,M\}$, or else there exists~$z'\in \Zd\cap (z+\cu_{k+h})$ such that~$\mathrm{B}_{k+1}(z') \cup \cdots \cup \mathrm{B}_{k+h}(z')$ occurs.
In summary, what we have shown is that
\begin{align*}
	\mathcal{C}_{k+h} (z)
	&
	\subseteq
	\Bigl\{
	\exists y_1,\ldots, y_M \in\Zd\cap (z+\cu_{k+h})  \,,
	\ \
	\dist(y_i,y_j) \geq 3^{k+1} \indc_{\{i\neq j\}} \,,
	\ \
	\mathcal{C}_{k}(y_i) \ \text{occurs}
	\Bigr\}
	\notag \\ & \qquad \qquad
	\cup
	\Bigl\{
	\exists z' \in \Zd \cap (z+\cu_{k+h}) \,,
	\ \
	\mathrm{B}_{k+1}(z') \cup \cdots \cup \mathrm{B}_{k+h}(z') \ \text{occurs}
	\Bigr\}
	\,.
\end{align*}
The number of finite sequences~$\{ y_1,\ldots , y_M\} \subseteq \Zd \cap  \cu_{k+h}$ with~$\dist(y_i,y_j) \geq 3^{k+1} \indc_{\{i\neq j\}}$ is at most~$3^{d(k+h)M}$; for any such sequence, we have, by independence,
\begin{equation*}
	\P \biggl[ \bigcap_{i=1}^M\mathcal{C}_{k}(y_i)  \biggr]
	=
	\prod_{i=1}^M
	\P \bigl[ \mathcal{C}_{k}(y_i)  \bigr]
	\leq
	\sup_{z\in \Zd} \P \bigl[ \mathcal{C}_{k}(z)  \bigr]^M
	\,.
\end{equation*}
Therefore a union bound yields
\begin{equation*}
	\P \bigl[ \mathcal{C}_{k+h} (z) \bigr]
	\leq
	\Bigl(
	3^{d(k+h)} \sup_{z\in \Zd} \P \bigl[ \mathcal{C}_{k}(z)  \bigr] \Bigr)^M
	+
	3^{d(k+h)}
	\sum_{l=k+1}^{k+h}
	\sup_{z\in\Zd}
	\P \bigl[ \mathrm{B}_l(z)  \bigr]
	\,.
\end{equation*}
Inserting~\eqref{e.badevent.upperbound.appendix} and taking the supremum over~$z$ yields~\eqref{e.paths.ready.for.iteration}.  We next iterate~\eqref{e.paths.ready.for.iteration} to obtain~\eqref{e.diameter.bound}. Fix~$b \in (0,1) \cap (0,a]$.
Select~$h \in \N$ satisfying
\begin{equation}
	h  \coloneqq  \bigl \lceil 8 (1-b)^{-1} \bigr\rceil \qand T^2 \geq  8 d b^{-1} 3^{2h}  \,\,.
	\label{e.h.T.select}
\end{equation}
Assume by induction that, for some~$k_0 \in \N$,
\begin{equation}
	\sup_{z\in\Zd}
	\P \bigl[ \mathcal{C}_{k} (z) \bigr]
	\leq
	\exp\Bigl( -\frac12 T^2 3^{b (k-h)} \Bigr) \,,
	\quad \forall k\in [0,k_0]\cap \N_0
	\,.
	\label{e.Ck.indy.hyp}
\end{equation}
Note that this is valid for~$k_0 = h$. Indeed, by a union bound and~\eqref{e.h.T.select} for any~$k \in [0,h] \cap \N_0$,
\begin{equation*}
	\sup_{z \in \Zd} \P \bigl[ \mathcal{C}_k(z) \bigr]
	\leq
	\sup_{z \in \Zd}
	\P \biggl[ \bigcup_{L\in [0,k]} \mathrm{B}_L(z)  \biggr]
	\leq
	\sum_{L=0}^k
	\exp\bigl( - T^2 3^{a L} \bigr)
	\leq
	(k+1)  \exp\bigl( - T^2 \bigr)
	\leq
	\exp\Bigl(-\frac12 T^2\Bigr)
	\,.
\end{equation*}
By induction, we obtain from~\eqref{e.paths.ready.for.iteration} and~$b\leq a$ that
\begin{align*}
	\lefteqn{
	\sup_{z \in \Zd} \P \bigl[ \mathcal{C}_{k+h}(z) \bigr]
	} \qquad & \notag \\
	& \leq  \exp\Bigl(  \bigl( \log(3) d (k+h) - \frac12 T^2 3^{b (k-h)} \bigr) 3^{h-4}  \Bigr)
	+ (k+h)3^{d(k+h)}
	\exp\bigl( - T^2 3^{a (k+1) } \bigr)
	\\
	&\leq  \exp\Bigl(  - \frac12 T^2 3^{b k}     3^{(1-b) h-4} +  \log(3) d (k+h) 3^{h-4} \Bigr)  +
	\exp\bigl( - T^2 3^{a (k+1) } + 3 d(k+h) \bigr) \\
	&\leq \frac12  \exp\Bigl( -\frac12 T^2 3^{b k} \Bigr)  + \frac12 \exp\Bigl( -\frac12 T^2 3^{a k} \Bigr)
	\leq
	\exp\Bigl( -\frac12 T^2 3^{b k} \Bigr)
	\,,
\end{align*}
proving the induction step.
We therefore have~\eqref{e.Ck.indy.hyp} for all~$k\in\N$. Using this, the definition of~$\mathrm{B}(z)$ in~\eqref{e.badevent.decomp}, the assumption~\eqref{e.badevent.upperbound.appendix},~$b\leq a$, the definition of~$h$ in~\eqref{e.h.T.select} and a union bound, we obtain,  for every~$k\in\N$,
\begin{align*}
	&
	\P
	\Bigl[\mbox{ there exists a path~$\Gamma$ from~$\Zd\cap \cu_k$ to~$\Zd \setminus \cu_{k+1}$ for which~$\bigcap_{y\in\Gamma} \mathrm{B}(y)$  occurs
	} \Bigr]
	\\ & \qquad\qquad
	\leq 3^{d k}  \Bigl( \sup_{z \in \Zd} \P[\mathcal{C}_{k}(z)] + 3^{d} \sum_{j=k+1}^{\infty} \sup_{z \in \Zd} \P\bigl[\mathrm{B}_j(z)\bigr] \Bigr)
	\\ & \qquad\qquad
	\leq
	C 3^{dk}
	\exp\Bigl(- \frac12  T^2  3^{b (k-h) }\Bigr)
	\leq
	\exp\Bigl(- \frac12 \exp\bigl(-C(1-b)^{-1}\bigr)T^2 3^{bk} \Bigr)
	\,.
\end{align*}
This completes the proof of~\eqref{e.diameter.bound}.
\end{proof}

\section{Functional inequalities}
\label{s.functional.inequalities}

This appendix collects the functional inequalities used throughout the paper. We begin by recalling the definitions of the relevant function space seminorms, then present fractional Poincar\'e-Sobolev inequalities, Besov-space duality estimates, and several technical lemmas relating Besov norms of different orders.

\subsection*{Definitions}

The H\"older seminorm~$[ \cdot]_{C^{0,s}(\cu_m)}$ is defined for each~$s \in (0,1]$ by
\begin{equation*}
	[ u ]_{C^{0,s}(\cu_m)}
	\coloneqq
	\sup_{x,y\in \cu_m, \, x\neq y}
	\frac{|u(x) - u(y)| }{|x-y|^s }
	\,.
\end{equation*}
The fractional Sobolev seminorm is defined for each~$s\in (0,1)$ and~$p\in [1,\infty)$ by
\begin{equation}
	\label{e.fractional.Sobolev.seminorm}
	[  u ]_{\underline{W}^{s,p}(\cu_m)}
	\coloneqq
	\biggl(
	\fint_{\cu_m} \int_{\cu_m}
	\frac{|u(x) - u(y)|^p}{|x-y|^{d+sp}}
	\, dx \, dy
	\biggr)^{\! \nf1p}
	\,.
\end{equation}
We work with the volume-normalized Besov-type seminorm
defined on the cube~$\cu_m$ for exponents~$s \in (0,1)$ and~$p,q\in [1,\infty)$ by
\begin{equation}
	[  u ]_{\underline{B}^s_{p,q}(\cu_m)}
	\coloneqq
	\biggl(
	\sum_{n=-\infty}^m
	3^{-nsq}
	\biggl(
	\avsum_{z\in 3^{n-1} \Zd \cap \cu_m, \, z+\cu_n \subseteq \cu_m}
	\bigl\| u - (u)_{z+\cu_n} \bigr\|_{\underline{L}^p(z+\cu_n)}^p
	\biggr)^{\!\nf qp} \biggr)^{\!\nf1q} \,.
	\label{e.Besov.seminorm.0}
\end{equation}
The corresponding norm is~$\|  f \|_{\underline{B}^s_{p,q}(\cu_m)}  \coloneqq  [  f ]_{\underline{B}^s_{p,q}(\cu_m)} + 3^{-sm} |(f)_{\cu_m}|$. When~$q=p$, this seminorm is equivalent to the fractional Sobolev seminorm: there exists a constant~$C(d)<\infty$ such that, for every~$s\in (0,1)$ and~$p\in [1,\infty)$,
\begin{equation}
	C^{-1} [  u ]_{\underline{W}^{s,p}(\cu_m)}
	\leq
	[ u ]_{\underline{B}^s_{p,p}(\cu_m)}
	\leq
	C [  u ]_{\underline{W}^{s,p}(\cu_m)}
	\,.
	\label{e.Wsp.vs.Bspp.intro}
\end{equation}
We also define the negative-order (or ``weak'') Besov seminorm by
\begin{equation}
	\label{e.negative.Besov.def}
	[ f ]_{\Besov{-s}{p}{q}(\cu_m)}
	\coloneqq
	\Biggl(
	\sum_{n=-\infty}^m
	3^{s qn}
	\biggl(
	\avsum_{z\in 3^n\Zd \cap \cu_m}
	\bigl| (f)_{z+\cu_n}\bigr |^p
	\biggr)^{\!\nf qp}
	\Biggr)^{\! \nf1q}\,.
\end{equation}

\subsection*{Fractional Poincar\'e-Sobolev inequality}

We denote the fractional Sobolev conjugate exponents by
\begin{equation*}
	\left\{
	\begin{aligned}
		& p^*(s,d)  \coloneqq  \frac{pd}{d-sp}\,, & &  p \in [1,\nf{d}{s} ) \,, \\
		& p_*(s,d)  \coloneqq  \frac{pd}{d+sp}\,, & & p \in [1,\infty)\,.
	\end{aligned}
	\right.
\end{equation*}
These satisfy the duality relation~$( p^*(s,d) )^\prime = (p^\prime)_*(s,d)$, where~$q^\prime  \coloneqq  q/(q-1)$ denotes the H\"older conjugate exponent.

The fractional Poincar\'e-Sobolev inequality (see~\cite[Theorem~6.7]{DNPaVa}) states that, for each~$s\in (0,1]$ and~$p\in [1, \nf ds )$ there exists~$C(s,p,d)<\infty$ such that, for every~$f \in W^{s,p}(\cu_m)$,
\begin{equation}
	\label{e.fractional.Sobolev.Poincare}
	\| f - (f)_{\cu_m} \|_{\underline{L}^{p^*(s,d)} (\cu_m)}
	\leq
	C 3^{sm}  [f ]_{\underline{W}^{s,p}(\cu_m)}
	\,.
\end{equation}
The fractional Morrey inequality (see~\cite[Theorem~8.2]{DNPaVa}) states that, for each~$s \in (0,1)$ and~$p \in [\nf ds,\infty)$, there exists a constant~$C(s,p,d) < \infty$ such that, for every~$f \in W^{s,p}(\cu_0)$,
\begin{equation}
	\label{e.frac.morrey}
	\| f - (f)_{\cu_m} \|_{L^\infty(\cu_0)}
	\leq
	C3^{sm} [f ]_{\underline{W}^{s,p}(\cu_m)}
	\,.
\end{equation}

We also use the following fractional Ne\v{c}as inequality on cubes: for every~$s\in(0,\nf12]$,~$p\in(1,\infty)$,~$m\in\Z$, and~$u\in H^1_0(\cu_m)$ for which the right side below is finite,
\begin{equation}
	\label{e.fractional.Necas}
	\|u\|_{\underline W^{1-s,p}(\cu_m)}
	\leq C(p,d)s^{-\nf1p}\|\nabla u\|_{\Wminusul{-s}{p}(\cu_m)}
	\,.
\end{equation}
Indeed, after rescaling to~$\cu_0$, this follows by interpolation between the Poincar\'e inequality and the Ne\v{c}as estimate~$\|u\|_{L^p}\leq C\|\nabla u\|_{W^{-1,p}}$; the factor~$s^{-1/p}$ is the endpoint normalization of the Slobodeckij seminorm. We also have~$(\nabla u)_{\cu_m}=0$, so the full and mean-zero negative norms are equivalent. In particular, if~$p>2d$, then~$1-s>d/p$ throughout this range, and the fractional Morrey inequality gives
\begin{equation*}
	3^{-m}\|u\|_{L^\infty(\cu_m)}
	\leq C(p,d)s^{-\nf1p}3^{-sm}\|\nabla u\|_{\Wminusul{-s}{p}(\cu_m)}
	\,.
\end{equation*}

For the Hilbertian negative Besov seminorm, we use the sharper zero-boundary estimate: for every~$\theta\in(0,\nf12]$,~$m\in\Z$, and~$u\in H^1_0(\cu_m)$,
\begin{equation}
	\label{e.negative.Besov.to.L2.zero.bndr}
	3^{-m}\|u\|_{\underline L^2(\cu_m)}
	\leq C(d)3^{-\theta m}
	[\nabla u]_{\Besov{-\theta}{2}{2}(\cu_m)}\,.
\end{equation}
After rescaling to~$\cu_0$, this follows by solving the Dirichlet problem~$-\Delta w=u$, applying Besov duality and Cauchy--Schwarz across the scales, and using the~$L^2$ Calder\'on--Zygmund estimate. The geometric scale sum is bounded uniformly for~$\theta\leq1/2$.

\subsection*{Besov multiplication inequality}

Let~$p,q \in [1,\infty]$. For every~$\zeta, \xi \in [1,\infty]$, we have
\begin{equation}
	\label{e.besov.mult}
	[ f g ]_{\underline{B}_{p,q}^s(\cu_m)}
	\leq
	4 \| f \|_{L^{\zeta'p}(\cu_m)} [ g ]_{\underline{B}_{\zeta p,q}^{s+\nf{d}{\zeta' p}}(\cu_m)} +   2 \|g\|_{L^{\xi}(\cu_m)}
	[ f ]_{\underline{B}_{p,q}^{s + \nf d \xi}(\cu_m)}
	\,.
\end{equation}
To see this, by H\"older's inequality, we have
\begin{align*}
	\lefteqn{
	\| fg - (fg)_{z+\cu_{n}} \|_{\underline{L}^p(z+\cu_{n})}
	} \quad &
	\notag \\ &
	\leq
	2\| f(g - (g)_{z+\cu_n} ) \|_{\underline{L}^p(z+\cu_{n})}
	+  | (g)_{z+\cu_n} | \| f - (f)_{z+\cu_n}  \|_{\underline{L}^p(z+\cu_{n})}
	\notag \\ &
	\leq
	2 \| f \|_{\underline{L}^{\zeta' p}(z+\cu_n)}
	\| g - (g)_{z+\cu_n} \|_{\underline{L}^{\zeta p}(z+\cu_{n})}
	+ 3^{\frac d\xi (m-n)} \| g \|_{\underline{L}^{\xi}(z+\cu_n)} \| f - (f)_{z+\cu_n}  \|_{\underline{L}^p(z+\cu_{n})}
	\notag \\ &
	\leq
	2 \cdot 3^{\frac{d}{\zeta' p}(m-n)}  \| f \|_{\underline{L}^{\zeta' p}(\cu_m)} \| g - (g)_{z+\cu_n} \|_{\underline{L}^{\zeta p}(z+\cu_{n})}
	+
	3^{\frac d\xi (m-n)} \| g \|_{\underline{L}^{\xi}(\cu_m)} \| f - (f)_{z+\cu_n}  \|_{\underline{L}^p(z+\cu_{n})}
	\,.
\end{align*}
Plugging this into~\eqref{e.Besov.seminorm.0}  and using the triangle inequality yields~\eqref{e.besov.mult}.

\subsubsection*{Subadditivity and superadditivity of seminorms}

For every~$s\in (0,1]$,~$p\in [1,\infty)$,~$m,n \in \Z$ with~$n\leq m$, and~$f \in W^{s,p}(\cu_m)$ and~$g\in L^{p}(\cu_m)$, we have
\begin{equation}
	\avsum_{z \in 3^n \Zd \cap \cu_m}  [ f ]_{\underline{W}^{s,p}(z+\cu_n)}^p
	\leq
	[ f ]_{\underline{W}^{s,p}(\cu_m)}^p
	\label{e.Hs.superadditivity}
\end{equation}
and
\begin{equation}
	\| g \|_{\underline{W}^{-s,p}(\cu_m)}
	\leq
	2 d^{\frac{d+sp'}{2p'}} 3^{s(m-n)}
	\biggl(
	\avsum_{z \in 3^n \Zd \cap \cu_m}
	\| g  \|_{\underline{W}^{-s,p}(z+\cu_n)}^p
	\biggr)^{\nf 1p}
	\,.
	\label{e.H.minus.s.subadditivity}
\end{equation}
We compute
\begin{align*}
	\avsum_{z \in 3^n \Zd \cap \cu_m}  [ f ]_{\underline{W}^{s,p}(z+\cu_n)}^p
	&
	=
	\avsum_{z \in 3^n \Zd \cap \cu_m} \fint_{z+\cu_n} \int_{z+\cu_n} \frac{|f(x) - f(y)|^p}{|x-y|^{d+sp}} \, dx \, dy
	\notag \\ &
	\leq
	\frac{1}{|\cu_m|} \sum_{z \in 3^n \Zd \cap \cu_m} \int_{z+\cu_n} \int_{\cu_m} \frac{|f(x) - f(y)|^p}{|x-y|^{d+sp}} \, dx \, dy
	=
	[ f ]_{\underline{W}^{s,p}(\cu_m)}^p \,.
\end{align*}
This is~\eqref{e.Hs.superadditivity}. The inequality~\eqref{e.H.minus.s.subadditivity} follows by duality. Indeed,
\begin{align*}
	\biggl| \fint_{\cu_m} f g  \biggr|
	\leq
	\avsum_{z\in 3^n\Zd\cap \cu_m}
	\biggl|\fint_{z+\cu_n}  fg  \biggr|
	&
	\leq
	\avsum_{z\in 3^n\Zd\cap \cu_m}
	\| f \|_{\underline{W}^{s,p'}(z+\cu_n)}
	\| g \|_{\underline{W}^{-s,p}(z+\cu_n)}
	\notag \\ &
	\leq
	\biggl(
	\avsum_{z\in 3^n\Zd\cap \cu_m}
	\| f \|_{\underline{W}^{s,p'}(z+\cu_n)}^{p'}
	\biggr)^{\!\nf1{p'}}
	\biggl(
	\avsum_{z\in 3^n\Zd\cap \cu_m}
	\| g \|_{\underline{W}^{-s,p}(z+\cu_n)} ^{p}
	\biggr)^{\!\nf1p}
	\notag \\ &
	\leq
	3^{s(m-n)}\| f \|_{\underline{W}^{s,p'}(\cu_m)}
	\biggl(
	\avsum_{z\in 3^n\Zd\cap \cu_m}
	\| g \|_{\underline{W}^{-s,p}(z+\cu_n)} ^{p}
	\biggr)^{\!\nf1p}
	\,.
\end{align*}
Here we used, in the last inequality,~\eqref{e.Hs.superadditivity} and the flat part of the norm~$\| \cdot \|_{\underline{W}^{s,p'}(\cu_m)}$ to obtain
\begin{align*}
	\biggl( \avsum_{z\in 3^n\Zd\cap \cu_m}
	\| f \|_{\underline{W}^{s,p'}(z+\cu_n)} ^{p'} \biggr)^{\nf 1{p'}}
	&
	\leq
	2\biggl(
	\avsum_{z\in 3^n\Zd\cap \cu_m}
	\bigl(
	3^{-sp'n} \|  f \|_{\underline{L}^{p'}(z+\cu_n)}^{p'}
	+
	[ f ]_{\underline{W}^{s,p'}(z+\cu_n)} ^{p'}\bigr)
	\biggr)^{\nf 1{p'}}
	\\ & \notag
	\leq
	2 \bigl( 3^{-sn}  \|  f \|_{\underline{L}^{p'}(\cu_m)} + [ f ]_{\underline{W}^{s,p'}(\cu_m)} \bigr)
	\leq
	2 d^{\frac{d+sp'}{2p'}} 3^{s(m-n)} \| f \|_{\underline{W}^{s,p'}(\cu_m)}
	\,.
\end{align*}
Taking the supremum over~$\| f \|_{\underline{W}^{s,p'}(\cu_m)} \leq 1$ yields~\eqref{e.H.minus.s.subadditivity}.

\subsection*{Duality between positive and negative Besov seminorms}
\label{ss.Besov.duality}

The following duality estimates pairs the negative-order seminorm~\eqref{e.negative.Besov.def} with the positive-order seminorm. For every~$m \in \Z$,~$s > 0$,~$p,q \in (1,\infty)$, and functions~$f,g$ with~$(g)_{\cu_m} = 0$,
\begin{equation}
	\label{e.duality.for.cubes}
	\biggl| \fint_{\cu_m} f g  \biggr|
	\leq
	3^{d+s}
	[ f ]_{\Besov{-s}{p}{q}(\cu_m)} [ g ]_{\underline{B}^s_{p',q'}(\cu_m)}
	\,.
\end{equation}
The proof can be found from~\cite[Lemma A.1]{AK.HC}.

\subsection*{Fractional Sobolev embedding from negative Besov norms}

For every~$s \in (0,1]$,~$p \in [1, (1-s)^{-1} d )$, and~$q \in [2 , q^*]$ with
\begin{equation}
	\label{e.qstar.def}
	q^*  \coloneqq  \frac{dp}{d-(1-s)p}
	\,,
\end{equation}
there exists~$C(d)<\infty$ such that, for every~$u \in H^{1}(\cu_0)$,
\begin{equation}
	\| u - (u)_{\cu_0} \|_{L^{q}(\cu_0)}
	\leq
	C q
	[\nabla u]_{\Besov{-s}{p}{1}(\cu_0)}
	\label{e.fractional.Sobolev.embedding}
\end{equation}
and, for every~$u \in  H_0^1(\cu_0)$,
\begin{equation}
	\| u \|_{L^{q}(\cu_0)}
	\leq
	C  q [\nabla u]_{\Besov{-s}{p}{1}(\cu_0)}
	\,.
	\label{e.fractional.Sobolev.embedding.zero.bndr}
\end{equation}

\begin{proof}[Proofs of~\eqref{e.fractional.Sobolev.embedding} and~\eqref{e.fractional.Sobolev.embedding.zero.bndr}]
For~\eqref{e.fractional.Sobolev.embedding}, assume without loss of generality that~$(u)_{\cu_0} = 0$ and consider the Neumann problem
\begin{equation*}
	\left\{
	\begin{aligned}
		& -\Delta w = |u|^{q-2} u - (|u|^{q-2} u)_{\cu_0}  & &
		\mbox{in } \cu_0\,,
		\\
		&  \mathbf{n}_U \cdot \nabla w  = 0 & &
		\mbox{on } \partial \cu_0\,.
	\end{aligned}
	\right.
\end{equation*}
For~\eqref{e.fractional.Sobolev.embedding.zero.bndr}, we instead consider the Dirichlet problem
\begin{equation*}
	\left\{
	\begin{aligned}
		& -\Delta w = |u|^{q-2} u  & &
		\mbox{in } \cu_0\,,
		\\
		& w = 0 & &
		\mbox{on } \partial \cu_0\,.
	\end{aligned}
	\right.
\end{equation*}
In both cases, the Calder\'on-Zygmund estimates (see, for example,~\cite{Stein}) give us
\begin{equation*}
	\bigl\| \nabla^2 w \bigr\|_{\underline{L}^{q'}(\cu_0)}
	\leq
	C (q'-1)^{-1} \| u \|_{\underline{L}^q(\cu_0)}^{q-1}
	\leq
	C q \| u \|_{\underline{L}^q(\cu_0)}^{q-1}
	\,.
\end{equation*}
By the duality estimate~\eqref{e.duality.for.cubes},
\begin{align*}
	\| u \|_{\underline{L}^q(\cu_0)}^q
	& =
	\fint_{U} \nabla w \cdot \nabla u
	\leq
	C [\nabla u]_{\Besov{-s}{p}{1}(\cu_0)}
	\| \nabla w\|_{\underline{B}_{p',\infty,1}^s(\cu_0)}
	\,.
\end{align*}
Denoting~$Z_n  \coloneqq  \{ z\in 3^{n-1} \Zd \cap \cu_0 : z+\cu_n \subseteq \cu_0 \}$ for~$n \in -\N_0$, we estimate
\begin{align*}
	\| \nabla w\|_{\underline{B}_{p',\infty,1}^s(\cu_0)}^{p'}
	&
	=
	\sup_{n \in\Z\cap (-\infty,0]}
	3^{-p'sn}
	\avsum_{z\in Z_n}
	\bigl\| \nabla w - (\nabla w)_{z+\cu_n} \bigr\|_{\underline{L}^1(z+\cu_n)}^{p'}
	\notag \\ &
	\leq
	C \sup_{n \in\Z\cap (-\infty,0]}
	3^{p'(1-s)n}
	\avsum_{z\in Z_n}
	\bigl\| \nabla^2 w \bigr\|_{\underline{L}^1(z+\cu_n)}^{p'}
	\notag \\ &
	\leq
	C \bigl\| \nabla^2 w \bigr\|_{\underline{L}^{q'}(U)}^{q'}
	\sup_{n \in\Z\cap (-\infty,0]}
	3^{p'(1-s)n}
	\max_{z\in Z_n} \bigl\| \nabla^2 w \bigr\|_{\underline{L}^{q'}(z+\cu_n)}^{p'-q'}
	\notag \\ &
	\leq
	C \bigl\| \nabla^2 w \bigr\|_{\underline{L}^{q'}(U)}^{p'}
	\sup_{n \in\Z \cap (-\infty,0]]}
	3^{p'(1-s)n}
	3^{-d (\frac{p'}{q'}-1) n}
	\notag \\ &
	\leq
	C q^{p'} \| u \|_{\underline{L}^q(U)}^{p'(q-1)}
	\,.
\end{align*}
Combining the previous two displays yields both~\eqref{e.fractional.Sobolev.embedding} and~\eqref{e.fractional.Sobolev.embedding.zero.bndr}. Note that the condition on~$q$ ensures
\begin{equation*}
	p'(1-s)-d\Big(\tfrac{p'}{q'}-1\Big)\ge 0
	\;\Longleftrightarrow\;
	\frac1{q'}\le \frac1{p'}+\frac{1-s}{d}
	\;\Longleftrightarrow\;
	q\le \frac{dp}{d-(1-s)p}
	= q^*
	\,.
	\qedhere
\end{equation*}
\end{proof}
The following estimate, which controls the Besov norm of a product involving a solution and a test function, can be found in~\cite[Lemma A.3]{AK.HC}.

\begin{lemma}
\label{l.divcurl}
There exists~$C(d) < \infty$ such that, for every~$m \in \Z$,~$s \in (0,1)$,~$u \in H^1(\cu_m)$, and~$\varphi \in W^{2,\infty}(\cu_m)$,
\begin{equation}
	\label{e.divcurl.est}
	\left\| (u- (u)_{\cu_m}) \nabla \varphi \right\|_{\underline{B}_{2,\infty}^{s}(\cu_m)}
	\leq
	C
	3^{m}
	\| \nabla \varphi\|_{\underline{W}^{1,\infty}(\cu_m)}
	\left[ \nabla u  \right]_{\Besov{s-1}{2}{1}(\cu_{m})}
	\,.
\end{equation}
\end{lemma}

\subsection*{Relating positive and negative Besov seminorms of a function and its gradient}

The following lemma shows that a positive-order Besov seminorm of a function can be controlled by a negative-order Besov seminorm of its gradient.

\begin{lemma}
\label{l.Besov.grad.to.function}
There exists~$C(d)<\infty$ such that, for every~$m\in \Z$,~$p \in (1,\infty)$,~$s\in(0,1)$, and~$v \in W^{1,p}(U)$,
\begin{equation}
	\label{e.v.in.Bpps.vs.nabla.v.in.B.pp.minus.t}
	[ v ]_{\underline{B}_{p,\infty}^{s}(\cu_m)}
	\leq
	C p s^{-1} [ \nabla v ]_{\underline{B}_{p,p}^{s-1}(\cu_m)}
	\,.
\end{equation}
\end{lemma}

\begin{proof}
Denoting~$Z_{n,m}  \coloneqq  \{ z \in 3^{n-1}\Zd \cap \cu_m \, :\, z + \cu_n \subset \cu_m\}$ and using~\eqref{e.fractional.Sobolev.embedding}, we estimate
\begin{align*}
	[ v ]_{\underline{B}_{p,\infty}^{s}(\cu_m)}^p
	& =
	\sup_{n \leq m} 3^{-sp n}  \avsum_{z \in Z_{n,m}  } \| v - (v)_{z+\cu_n} \|_{\underline{L}^p(z+\cu_n)}^p
	\notag \\ &
	\leq
	(C p)^p
	\sup_{n \leq m}  3^{-sp n}  \avsum_{z \in Z_{n,m}  }
	[\nabla v]_{\Besov{-1}{p}{1}(z+\cu_n)}^p
	\notag \\ &
	=
	(C p)^p
	\sup_{n \leq m}  3^{-sp n}  \avsum_{z \in Z_{n,m}  }  \biggl( \sum_{k=-\infty}^{n} 3^{k} \
	\biggl( \avsum_{z' \in z+ 3^k \Zd \cap \cu_n } \bigl|(\nabla v)_{z'+\cu_k} \bigr|^p \biggr)^{\! \nf 1p}  \biggr)^p
	\notag \\ &
	\leq
	(C p s^{-1})^p
	\sup_{n \leq m}  \sum_{k=-\infty}^{n} 3^{(1-s)pk}
	\avsum_{z \in Z_{n,m} } \avsum_{z' \in z+ 3^k \Zd \cap \cu_n } \bigl|(\nabla v)_{z'+\cu_k} \bigr|^p
	\notag \\ &
	\leq
	(C p s^{-1})^p [ \nabla v ]_{\underline{B}_{p,p}^{s-1}(\cu_m)}^p
	\,.
\end{align*}
This yields~\eqref{e.v.in.Bpps.vs.nabla.v.in.B.pp.minus.t}.
\end{proof}

\section{Concentration inequalities}
\label{s.concentration}

In this appendix, we prove a concentration inequality for indicator functions of rare events, which is needed for the ``approximate independence between scales'' argument in Section~\ref{ss.scale.local}. It can be viewed as a concentration inequality for exceedance counts in a triangular array with weak dependence. While concentration inequalities for~$m$-dependent sequences are classical, the specific structure here---geometric weights within rows combined with column-wise~$r$-dependence---appears to be new. It is motivated by the multiscale structure of the renormalization group iteration in Section~\ref{s.RG.flow}.

\begin{proposition}
\label{p.concentration.for.scales}
Let~$p\in[1,\infty)$,~$r \in \N$ and~$s  \in [\nf1p,1]$. Let~$\{ X_{k,j}\,:\,  k \in \Z, \, j \in \Z \}$ be a collection of nonnegative random variables satisfying the moment bound
\begin{equation}
	\label{e.Xk.p.moment.twosided}
	\E \bigl[ X_{k,j}^{p} \bigr]
	\leq
	1
	\,, \quad \forall k,j \in \Z\,,
\end{equation}
and~$r$-dependence across columns: for all subsets~$I,J \subseteq \Z$ with~$\dist(I,J) \geq r$,
\begin{equation}
	\label{e.independent.columns.twosided}
	\{ X_{k,j}\,:\, k \in \Z\,, \, j \in I \}
	\quad \text{and} \quad
	\{ X_{k,j}\,:\, k \in \Z\,, \, j \in J \}
	\quad \text{are independent.}
\end{equation}
There exists a universal constant~$C > 0$ such that, for every~$m_0 \in \Z$,~$m \in\N_0$ and~$\theta \in (0,1]$,
\begin{equation}
	\label{e.no.bad.scales.twosided}
	\P \biggl[ \frac{1}{m+1} \sum_{k=m_0}^{m_0+m} \indc{\biggl\{ \sum_{j \in \Z} 3^{-s |k-j|} X_{k,j} > 6s^{-1} C^{\nf 1p} \theta^{-\nf 1p} \biggr\}} > \theta \biggr]
	\leq
	\exp \biggl(- \frac{  s p \theta}{16 r} (m+1)  \biggr)
	\,.
\end{equation}
\end{proposition}

\begin{proof}
Without loss of generality, we assume that~$m_0 = 0$ (otherwise consider~$\tilde{X}_{k,j}\coloneqq X_{k+m_0,j}$ instead of~$X_{k,j}$). Let~$\lambda > 0$ be a parameter to be chosen below. We use~$C_1, C_2, \ldots$ to denote positive universal constants. Define the weighted sums
\begin{equation}
	\label{e.Yk.def.twosided}
	Y_k  \coloneqq  \sum_{j \in \Z} 3^{-s |k-j|} X_{k,j}\,, \quad k\in\N_0\,.
\end{equation}

\smallskip

\emph{Step 1.} Reduction to a counting problem. Observe that, for each~$k \in \N_0$ and~$\lambda>0$,
\begin{equation*}
	\sup_{j \in \Z} 3^{-\frac{s}{2}|k-j|} X_{k,j}
	\leq \frac12\lambda
	\quad\implies\quad
	Y_k
	\leq \frac{\lambda}{2}  \sum_{\ell \in \Z} 3^{-\frac{s}{2} |\ell|}
	\leq 2 \lambda  s ^{-1}  \,.
\end{equation*}
Taking the contrapositive of this implication and summing over~$k$, we obtain
\begin{equation}
	\label{e.YktoZk.twosided}
	\sum_{k=0}^{m} \indc \{  Y_k > 2 \lambda s^{-1}  \}
	\leq
	\sum_{k=0}^{m}
	\sum_{j \in \Z}
	\indc{\bigl\{ X_{k,j} > \tfrac{1}{2} \lambda   3^{\frac{s}{2}|k-j|} \bigr\}}
	=
	\sum_{j \in \Z}
	Z_{j}
	\,,
\end{equation}
where we define
\begin{equation}
	\label{e.Zj.def.twosided}
	Z_{j} \coloneqq
	\sum_{k=0}^{m}
	\indc{\bigl\{ X_{k,j} > \tfrac{1}{2} \lambda  3^{\frac{s}{2}|k-j|} \bigr\}}
	\,.
\end{equation}

\smallskip

\emph{Step 2.} Tail bound for~$Z_j$. Fix~$j \in \Z$ and
let~$d_j \coloneqq \dist(j, \{0,\ldots,m\})$ denote the distance from~$j$ to the interval~$\{0,\ldots,m\}$. We claim that
\begin{equation}
	\label{e.Zj.tail.twosided}
	\P \bigl[ Z_j \geq n \bigr]
	\leq
	C_2
	( 2 \lambda^{-1})^p
	3^{-\frac{sp}{2} d_j}
	3^{-\frac{sp}{4} (n-2)_+}
	\,,
\end{equation}
where~$(n-2)_+ = \max\{n-2, 0\}$. To see this, note that for each~$\ell \geq d_j$, there are at most~$2$ values of~$k \in \{0,\ldots,m\}$ with~$|k-j| = \ell$, and there are no such~$k$ for~$\ell < d_j$. If~$Z_j \geq n$, then at least~$\lceil n/2 \rceil$ distinct values of~$\ell = |k-j|$ must occur among the~$k$ contributing to~$Z_j$. Since all such~$\ell$ satisfy~$\ell \geq d_j$, the largest must satisfy
\[
\ell_{\max} \geq d_j + \lceil n/2 \rceil - 1 \geq d_j + \frac{n-2}{2}\,.
\]
Hence, if~$Z_j \geq n$, there exists~$\ell \geq d_j + \frac{(n-2)_+}{2}$ and~$k \in \{0,\ldots,m\}$ with~$|k-j| = \ell$ such that~$X_{k,j} > \frac{\lambda}{2} 3^{\frac{s}{2}\ell}$. By a union bound and Markov's inequality, using~\eqref{e.Xk.p.moment.twosided},
\begin{align}
	\P \bigl[ Z_j \geq n \bigr]
	& \leq
	\sum_{\ell \geq d_j + (n-2)_+/2}
	\sum_{k : |k-j| = \ell}
	\P\bigl[ X_{k,j} > \tfrac12 \lambda 3^{\frac{s}{2}\ell} \bigr]
	\notag \\ &
	\leq
	2\Bigl( \frac{2}{ \lambda}\Bigr)^{\!p}
	\sum_{\ell \geq d_j + (n-2)_+/2}
	3^{-\frac{sp}{2} \ell}
	\notag \\ &
	=
	2\Bigl( \frac{2}{ \lambda}\Bigr)^{\!p}
	\frac{3^{-\frac{sp}{2}(d_j + (n-2)_+/2)}}{1 - 3^{-\frac{sp}{2}}}
	\leq
	C_2 (2\lambda^{-1})^p
	3^{-\frac{sp}{2} d_j}
	3^{-\frac{sp}{4} (n-2)_+}
	\,,
	\label{e.Zj.tail.derivation}
\end{align}
where we can take~$C_2 \coloneqq  \frac{4}{1 - 3^{-1/2}}$ since~$sp \geq 1$. This is~\eqref{e.Zj.tail.twosided}.

\smallskip

\emph{Step 3.} Exponential moment bound for~$Z_j$.  We claim that, for~$C_3 \coloneqq 2(1-3^{-\nf 18})^{-1}C_2$ and~$a \coloneqq \frac18 (\log 3) sp$,
\begin{equation}
	\label{e.Zj.mgf.twosided}
	\E\bigl [e^{a Z_j} \bigr]
	\leq
	1+ C_3 (3 \lambda^{-1})^p 3^{-\frac{sp}{2} d_j}  \,.
\end{equation}
From~\eqref{e.Zj.tail.twosided}, we have, since~$e^a = 3^{\frac18 sp}$ and~$Z_j$ is integer-valued,
\begin{align*}
	\E\bigl [e^{a Z_j} \bigr] - 1
	\leq
	\sum_{n=1}^\infty 3^{\frac{sp}{8} n} \P[Z_j = n]
	&\leq
	C_2 (2\lambda^{-1})^p 3^{-\frac{sp}{2} d_j}
	\sum_{n=1}^\infty 3^{\frac{sp}{8} n} 3^{-\frac{sp}{4}(n-2)_+}
	\\
	&\leq
	2 C_2 (2\lambda^{-1})^p 3^{-\frac{sp}{2} d_j} 3^{\frac{sp}{4}}
	\sum_{n=0}^\infty 3^{-\frac{sp}{8}n}
	\leq
	C_3 (2 \cdot 3^{\nf s4} \lambda^{-1})^p 3^{-\frac{sp}{2} d_j}
\end{align*}
since~$p \geq s^{-1}$.
Hence~\eqref{e.Zj.mgf.twosided} follows since~$2 \cdot 3^{\nf s4} \leq 3$.

\smallskip

\emph{Step 4.}
We show, using independence, that
\begin{equation}
	\label{e.full.mgf.bound}
	\E\Biggl[ \exp\biggl( \frac{a}{r} \sum_{j \in \Z} Z_j \biggr) \Biggr]
	\leq
	\exp \bigl(
	C_6 r^{-1} (3\lambda^{-1})^p (m+1)
	\bigr)
	\,.
\end{equation}
The~$r$-dependence assumption~\eqref{e.independent.columns.twosided} implies that for~$|j - j'| \geq r$, the random variables~$Z_j$ and~$Z_{j'}$ are independent. We decompose the sums into~$r$ subsequences by residue class modulo~$r$. Applying H\"older's inequality and then independence within each residue class,
\begin{align*}
	\E\Biggl[ \exp\biggl(  \frac{a}{r}  \sum_{j \in \Z} Z_j  \biggr) \Biggr]
	&=
	\E\Biggl[ \exp\biggl(  \frac{a}{r}
	\sum_{b=0}^{r-1}
	\sum_{j \in \Z} Z_{jr + b} \biggr) \Biggr]
	\\
	&\leq
	\prod_{b=0}^{r-1}
	\E\Biggl[ \exp\biggl( a \sum_{j \in \Z} Z_{jr + b} \biggr) \Biggr]^{\!\nf 1r}
	=
	\prod_{b=0}^{r-1}
	\prod_{j \in \Z}
	\E\bigl[e^{a Z_{jr + b}}\bigr]^{\nf 1r}
	\,.
\end{align*}
Using~\eqref{e.Zj.mgf.twosided} and~$\log(1+x) \leq x$, this is bounded by
\begin{equation*}
	\prod_{b=0}^{r-1}
	\prod_{j \in \Z}
	\E\bigl[e^{a Z_{jr + b}}\bigr]^{\nf 1r}
	\leq
	\exp \biggl(
	\frac1r
	C_3 (3\lambda^{-1})^p \sum_{j \in \Z} 3^{-\frac{sp}{2} d_j}
	\biggr)
	\,.
\end{equation*}
Now,
\[
\sum_{j \in \Z} 3^{-\frac{sp}{2} d_j}
= (m+1) + 2\sum_{\ell=1}^\infty 3^{-\frac{sp}{2}\ell}
= (m+1) + \frac{2 \cdot 3^{-\frac{sp}{2}}}{1 - 3^{-\frac{sp}{2}}}
\leq (m+1) + C_4
\leq C_5(m+1)\,,
\]
since~$sp \geq 1$ implies the geometric series is bounded by a universal constant, and~$m+1 \geq 1$. Therefore~\eqref{e.full.mgf.bound} follows with~$C_6 = C_3 C_5$.

\smallskip

\emph{Step 5.} The Conclusion. Combining~\eqref{e.YktoZk.twosided} with~\eqref{e.full.mgf.bound} and applying Markov's inequality,
\begin{align}
	\P \biggl[ \frac{1}{m+1} \sum_{k=0}^m \indc{\{ Y_k > 2 \lambda s^{-1}  \}} > \theta \biggr]
	&
	\leq
	\E\Biggl[ \exp\biggl( \frac{a}{r} \sum_{j \in \Z} Z_{j} \biggr) \Biggr]
	\exp \biggl( - \frac{a \theta}{r}(m+1)  \biggr)
	\notag \\[4pt] &
	\leq
	\exp \biggl( - \frac{a}{r} (m+1) \biggl(\theta
	- \frac{C_6}{ a} (3\lambda^{-1})^p
	\biggr) \biggr)
	\,.
	\label{e.propb2.ready.to.pick.lambda.twosided}
\end{align}
Choose~$\lambda \coloneqq 3 (2C_6 / (a\theta))^{\nf 1p}$, so that the inner parenthesis in~\eqref{e.propb2.ready.to.pick.lambda.twosided} equals~$\theta/2$, giving
\[
\P \biggl[ \frac{1}{m+1} \sum_{k=0}^m \indc{\bigl\{ Y_k > 6 s^{-1} (2C_6 a^{-1} \theta^{-1} )^{\nf 1p}   \bigr\}} > \theta \biggr]
\leq
\exp \biggl( - \frac{a \theta}{2r} (m+1) \biggr)
\,.
\]
Recalling~$a \coloneqq (\log 3) 8^{-1} sp$ then gives~\eqref{e.no.bad.scales.twosided}
for a suitable universal constant~$C$, using~$sp \geq 1$.
\end{proof}

\begingroup
\small

\subsubsection*{\bf Acknowledgments}
The authors thank Antti Kupiainen for very helpful comments on an earlier draft of this paper. S.A. and T.K. acknowledge support from the European Research Council (ERC) under the European Union's Horizon Europe research and innovation programme, grant agreement number 101200828. S.A. acknowledges support from NSF grant DMS-2350340. A.B. was supported by NSF grant DMS-2202940. T.K. was supported by the Academy of Finland.

{
\bibliographystyle{alpha}
\bibliography{refs}
}

\endgroup

\end{document}